\documentclass{amsart}
\usepackage[utf8]{inputenc}
\usepackage[T1]{fontenc}
\usepackage{imakeidx}
\usepackage{makeidx}
\makeindex[columns=2, options= -s example_style.ist]
\usepackage{eurosym}
\usepackage{amscd}
\usepackage{amsthm}
\usepackage{extarrows}
\usepackage{latexsym}
\usepackage{mathrsfs}
\usepackage{pst-all,multido,ifthen}
\usepackage{amssymb,amsmath,amsfonts}
\usepackage{tikz}
\usepackage{tikz-cd}
\usepackage[all]{xy}
\usepackage{subfigure}
\usepackage{xcolor}
\usepackage[unicode=true]{hyperref}
\hypersetup{
    colorlinks=true,      % False: boxed links; True: colored links
    citecolor=blue!70!black, % Highly visible dark blue for citations
    linkcolor=red!70!black,  % Dark red for internal links (sections, figures)
    urlcolor=blue!90!black   % Bright dark blue for web URLs
}
\usepackage{stmaryrd}
\newtheorem{theorem}{Theorem}[section]
\newtheorem{corollary}[theorem]{Corollary}
\newtheorem{conjecture}[theorem]{Conjecture}
\newtheorem{criterion}[theorem]{Criterion}
\newtheorem{definition}[theorem]{Definition}
\newtheorem{example}[theorem]{Example}
\newtheorem{fact}[theorem]{Fact}
\newtheorem{lemma}[theorem]{Lemma}
\newtheorem{remark}[theorem]{Remark}
\newtheorem{notation}[theorem]{Notation}
\newtheorem{proposition}[theorem]{Proposition}

\def\Aut{\textup{Aut}}

\def\et{\textup{\'et}}
\def\Alb{\textup{Alb}}

\def\a{\textup{a}}

\def\e{\textup{e}}

\def\f{\textup{f}}
\def\ab{\textup{ab}}

\def\Model{\textup{Model}}
\def\an{\textup{an}}

\def\r{\textup{r}}
\def\sup{\textup{sup}}
\def\BE{\textup{BE}}
\def\MBE{\textup{MBE}}
\def\pMor{\textup{p-Mor}}

\def\C{\textup{C}}
\def\D{\textup{D}}
\def\Div{\textup{Div}}
\def\chr{\textup{char}}

\def\Der{\textup{Der}}

\def\EE{\textup{EE}}
\def\Monomial{\textup{Monomial}}
\def\AffMonomial{\textup{AffMonomial}}

\def\End{\textup{End}}
\def\Res{\textup{Res}}

\def\Ext{\textup{Ext}}

\def\Triv{\textup{Triv}}
\def\Irr{\textup{Irr}}
\def\Frob{\textup{Frob}}
\def\Frac{\textup{Frac}}
\def\proj{\textup{proj}}

\def\GA{\textup{GA}}
\def\TGA{\textup{TGA}}
\def\Gal{\textup{Gal}}
\def\GL{\pmb{\textup{GL}}}
\def\gr{\textup{gr}}

\def\c{\textup{c}}
\def\s{\textup{s}}
\def\z{\textup{z}}

\def\h{\textup{h}}
\def\Hom{\textup{Hom}}
\def\Imm{\textup{Im}}

\def\Intr{\textup{IntRat}}
\def\Nointr{\textup{NorIntRat}}
\def\Smintr{\textup{SmoIntRat}}
\def\Affintr{\textup{AffIntRat}}

\def\ciedim{\textup{ciedim}}
\def\wciedim{\textup{wciedim}}
\def\sciedim{\textup{sciedim}}
\def\edim{\textup{edim}}

\def\Ker{\textup{Ker}}

\def\univ{\textup{univ}}

\def\m{\textup{m}}

\def\N{\textup{N}}

\def\nat{\textup{nat}}

\def\perf{\textup{perf}}
\def\proj{\textup{proj}}
\def\Pic{\textup{Pic}}

\def\QFE{\textup{QFE}}
\def\QF{\textup{QF}}
\def\red{\textup{red}}
\def\Reg{\textup{Reg}}
\def\Sing{\textup{Sing}}

\def\i{\textup{i}}
\def\l{\textup{l}}
\def\q{\textup{q}}
\def\r{\textup{r}}
\def\n{\textup{n}}
\def\gen{\textup{gen}}
\def\w{\textup{w}}
\def\sc{\textup{sc}}
\def\SL{\pmb{\textup{SL}}}

\def\Spec{\textup{Spec\,}}

\def\t{\textup{t}}
\def\tr{\textup{tr}}
\def\Tr{\textup{Tr}}
\def\s{\textup{s}}
\def\sep{\textup{sep}}
\def\t{\textup{t}}
\def\top{\textup{top}}
\def\Inv{\textup{Inv}}
\begin{document}
\title[On endomorphisms of affine spaces and the Jacobian problem]{On endomorphisms of affine spaces and the Jacobian problem}
\author{Alexander Borisov, Ofer Gabber, and Adrian Vasiu}
\date{\today}
\maketitle

\noindent
{\bf ABSTRACT.} Let $p$ be a prime. We provide examples which show that \'etale endomorphisms of affine planes over an algebraically closed field $k$ of characteristic $p$ can have fibers of arbitrary finite cardinal. Let $(l,m)\in\mathbb N\times\mathbb N^{\ast}$. We provide examples of such \'etale endomorphisms whose images have complements of cardinality $l$ and whose geometric degrees are $pm$. Several conjectures are disproved, and in particular we provide an analog over $k$ of Kulikov's counterexample to the `Generalized Jacobian Conjecture' of Bass for surfaces. Surjective (resp.\ surjective and non-surjective) counterexamples to Adjamagbo's Jacobian Conjecture over each $k$ are included in dimension $2$ (resp.\ in all dimensions at least $3$); for instance, if $p=2$, then we show for each $m$ there exist surjective \'etale endomorphisms of the affine spaces over $k$ of dimension at least $3$ of geometric degree $m$. If $e:X\rightarrow X$ is an endomorphism of a variety over an algebraically closed field $K$, then we show that there exists $n\in\mathbb N$ such that $\Imm(e^n)=\Imm(e^{n+1})$ provided either (i) $e$ is quasi-finite or (ii) $\dim(X)\le 2$ and $X\setminus\Imm(e)$ is finite. We provide examples of endomorphisms of affine planes of dimensions at least $3$ over $K$ whose images have complements of cardinality $l$ and for all $n\in\mathbb N$ we have $\Imm(e^n)\neq\Imm(e^{n+1})$. We prove that all affine moduli schemes of \'etale endomorphisms of affine spaces with Jacobian matrices of determinant $1$ over $K$ are connected and classify all those that are smooth. We prove that the classical Jacobian Conjecture over $\mathbb C$ and Adjamagbo's analog of it over $k$ hold for \'etale endomorphisms of affine spaces that are composites $g\circ f$, where $f$ is quasi-finite and a locally closed embedding in codimension $1$ outside a specific finite subset and $g$ is a projection that omits one affine coordinate.

\bigskip\noindent
{\bf KEY WORDS:} affine space, automorphism, complete intersection, endomorphism, \'etale, field, hypersurface, invariant, Jacobian Conjecture, matrix, morphism, parametrization, prime number, quasi-finite, ring, scheme, smooth, surface, surjective, torsor, variety, vector bundle, zeta function, $p$-basis.

\bigskip\noindent
{\bf MSC 2020:} 11M38, 13A35, 13A50, 13B02, 13B05, 13B10, 13B25, 13B40, 14E05, 14E08, 14E20, 14E25, 14F35, 14G17, 14J26, 14J30, 14J40, 14J50, 14J70, 14L24, 14L30, 14R05, 14R10, 14R15, 14R20, 14R25, and 20G40.

\newpage\tableofcontents

\newpage\section{Introduction}\label{S1}

The classical Jacobian Conjecture\index{Jacobian Conjecture} states that if a collection of $2$ complex polynomials in $2$ indeterminates has a non-zero constant Jacobian determinant, then the corresponding polynomial map, to be called a (complex) Kraus--Keller map\index{Kraus--Keller map} (in dimension $2$), has a polynomial inverse. In other words, every locally analytically invertible polynomial self-map of an affine plane over $\mathbb C$ is globally invertible, by a polynomial map. This conjecture goes back to the 1884 and 1939 papers of Kraus and Keller (respectively), see \cite{Kraus} and \cite{Kel} and \cite{RD}, and has received considerable attention. In fact, it is known to be a very tricky problem that has generated a large number of incorrect proofs. 

On its face, the conjecture seems to be solely about characteristic $0$, and, indeed, many attempts to solve it use analytic methods. However, it is also known that if it is false, one can construct a counterexample over the integers that involves more indeterminates than $2$ (see \cite{CvdD}, Subsect.\ 4.3). 

Let $n$ be a positive integer. It is also well-known that if a polynomial map that defines an endomorphism of the complex affine space of dimension $n$ is injective then it must be surjective and, moreover, its inverse must also be a polynomial map. That injectivity implies surjectivity is a particular case of Ax--Grothendieck's Theorem (see \cite{Ax}, Thm., \cite{Gro4}, Prop.\ (10.4.11); see also \cite{Bore1}, Thm.\ 2.3 and Sect.\ 3), whose original proof is by {\it reduction to positive characteristic}. Due to this, the Jacobian Conjecture can be viewed as a strengthening of the following implicit (folklore) conjecture.

\begin{conjecture}[Surjectivity Conjecture\index{Surjectivity Conjecture} over $\mathbb C$]\label{CJ1}
Each complex Kraus-Keller map in dimension $2$ is surjective.\index{Jacobian Conjecture!Surjectivity}
\end{conjecture}

Let $p$ be a prime number. If taken verbatim, the Jacobian Conjecture in characteristic $p$ is false already in dimension $1$, as the polynomial $x+x^p$ has derivative $1$ but the function $x\mapsto x+x^p$ is not injective; however, it is surjective over algebraically closed fields. In fact, in dimension $1$ the Surjectivity Conjecture is clearly true over all algebraically closed fields.

Generalizing the above example, one can consider in characteristic $p$ a natural class of Kraus-Keller maps, that we call {\it basic endomorphisms}:
$$e(x_1,\ldots,x_n)= \bigl(x_1+f_1(x_1^p,\ldots,x_n^p),\ldots,x_n+f_n(x_1^p,\ldots,x_n^p)\bigr).$$
It is natural to ask if over algebraically closed fields they are surjective. 

Finding a preimage by $e$ of the point $(0,\ldots,0)$ and, by a shift, of any other point, leads to a system of equations of a particular kind, that we call a Frobenius system:
\begin{equation}\label{EQ1}
x_i+f_i(x_1^p,\ldots,x_n^p)=0,\ \ \quad i\in\{1,\ldots,n\}.
\end{equation}
Based on the case when the $f_i$s have degrees at most 1 (see \cite{V}, Thm.\ 2.4.1(b)) and other explicit examples, the third author had conjectured that Frobenius systems over algebraically closed fields are always consistent.\footnote{Talk at the 60th
Anniversary of Gerd Faltings, MPI, Bonn, Germany on June 13, 2014.} His persistent effort to prove it was the original impetus and the starting point of our investigation.

Incidentally, if the $p$-th power is applied to the indeterminates of the first terms rather than to the second terms of the left hand side of System (\ref{EQ1}), the system
$$x_i^p+f_i(x_1,\ldots,x_n)=0,\ \ \quad i\in\{1,\ldots,n\}$$
\phantomsection{is consistent if $p$ is large enough (say $p$ is greater than the degree of $f_i$ for each $i\in\{1,\ldots,n\}$). Moreover, if the field is algebraically closed of characteristic $p$ and a dominant map $(x_1,\ldots,x_n)\mapsto \bigl(f_1(x_1,\ldots,x_n),\ldots,f_n(x_1,\ldots,x_n)\bigr)$ is fixed, then the union indexed by $q\in\mathbb N^{\ast}$ of the set of solutions of the system}\label{PH0} 
$$x_i^{p^q}+f_i(x_1,\ldots,x_n)=0,\ \ \quad i\in\{1,\ldots,n\}$$ is Zariski dense in the affine space (e.g., see \cite{BS}, Thm.\ 1.4 for the particular case of an algebraic closure of $\mathbb F_p$ and see \cite{SV}, Thm.\ 1.3 for the general case\footnote{The general case follows from the particular case via a standard specialization argument as for $q\gg 1$ such systems define finite flat endomorphisms of affine spaces.}).
This added to the hope that Frobenius systems are consistent, though the first author considered this problem to be about as hard as the Jacobian Conjecture.

The actual answer was opposite to the initial expectation. Indeed, we have the following particular case of Example \ref{EX7} for $m=1$ and $f(y)=1-y+y^p$.

\begin{example}\normalfont\label{EX1}
Let $i\in\{0,1,\ldots,p-1\}$. The Frobenius system 
$$\left\{ \begin{array}{l}\!\! x\!=\bigl(1\!\!+\!\!x^p(1\!\!-\!\!y\!\!+\!\!y^p)^p\bigr)\, \!\bigl[1\!\!-\!\!x\bigl(1\!\!-\!\!x(1\!\!-\!\!y\!\!+\!\!y^p)\!\!+\!\!x^p(1\!\!-\!\!y\!\!+\!\!y^p)^p\bigr)^{p\!-\!1}\bigl((x(1\!\!-\!\!y\!+\!y^p))^{p\!-\!1}\!\!-\!\!1\bigr)\bigr]^p \\

\!\!y=i+x^p(1-y+y^p)^p\end{array}\right.$$
has no solutions in any field of characteristic $p$.
\end{example}

To show that the assumption that it has a solution leads to a contradiction, let $z:=x(1-y+y^p)$. Hence $y=i+z^p$ by the second equation. As $i^p=i$ in fields of characteristic $p$ and $z^p=y-i$, the first equation can be rewritten as
$$x=(1+z^p)\bigl[1-x^p(1-y+y^p)^{p-1}\bigl((y-i)^{p-1}-1\bigr)\bigr].$$
Multiplying this by $1-y+y^p$, we get $z=(1+z^p)\bigl[1-y+y^p-(y-i)\bigl((y-i)^{p-1}-1\bigr)\bigr]$, which simplifies to $z=1+z^p$. Thus $1-z+z^p=0$. Taking 
this to the power $p$, we get $1-y+y^p=0$. Hence $y=i$ by the second equation. Thus $1-y+y^p=1$, which contradicts $1-y+y^p=0$.

This example implies that, in characteristic $p$, the naive analog of the Surjectivity Conjecture is false, even for basic endomorphisms in dimension 2. In particular, it answers negatively the question (1) of \cite{Lang-J}, paragraph after Rmk.\ 2.11.

An important work on the Jacobian Conjecture in characteristic $p$ was done by Adjamagbo. In particular, he proposed the following conjecture as the natural analog of the Jacobian Conjecture in characteristic $p$. 

\begin{conjecture}[Jacobian Conjecture in characteristic $p$\index{Jacobian Conjecture!in characteristic $p$}, \cite{Ad}, Conj.\ 3.1]\label{CJ2}
Every Kraus--Keller map in dimension $n$ over an algebraically closed field of characteristic $p$ is an automorphism if its geometric degree is not a multiple of $p$. [The geometric degree is the cardinality of the preimage of a generic point or the degree of the corresponding field extension.]
\end{conjecture}

In this paper we disprove Conjecture \ref{CJ2} in all dimensions $n\ge 2$ and over each $k$ (see Section \ref{S28} for counterexamples). Despite this setback, it feels critical to mention the following implicit (folklore) conjecture.

\begin{conjecture}[Weak Jacobian Conjecture in positive characteristic\index{Jacobian Conjecture!Weak in positive characteristic}]\label{CJ2.5}
For each $m\in\mathbb N^{\ast}\setminus\{1\}$ there exists $N(m)\in\mathbb N^{\ast}$ such that for every Kraus--Keller map in dimension $2$ over an algebraically closed field of characteristic $p\ge N(m)$ we have either $\deg(e)=1$ or $\deg(e)>m$.
\end{conjecture}

Stronger forms of this conjecture would provide some formulas for the $N(m)$s, with the stronger form one could  envision being to take $N(m)=m+1$.

While this conjecture seems to be weaker than the Jacobian Conjecture, it is stronger in the following sense: its validity implies the classical Jacobian Conjecture. At this time we cannot prove or disprove it. 

We also disprove another conjecture by Adjamagbo called the `Unirational Jacobian Conjecture $\textup{UJC}(n,p)$' (see \cite{Ad}, Conj.\ 3.4) for all integers $n\ge 2$ and primes $p$ (for $n\in\{2,3\}$ see Theorem \ref{T2}). `Conjecture $\textup{UJC}(n,p)$' was meant to be a characteristic $p$ analog of the `Generalized Jacobian Conjecture' of Bass in characteristic $0$ (see \cite{Ba2}, Rmk.\ after Prop.\ (4.3)) that was disproved by Kulikov in dimension $2$ (see \cite{Kul}, Sect.\ 3; the explicit equations of loc. cit. work over any algebraically closed field of characteristic different from $2$ and $3$). 

The main goal of this paper is to initiate a systematic study of several natural invariants of basic endomorphisms and general Kraus--Keller maps (\'etale endomorphisms) in arbitrary dimension $n$ over algebraically closed fields of characteristic $p$. For the sake of generality and for including also the characteristic $0$ case whenever possible, we often work with larger classes of endomorphisms of affine spaces over arbitrary algebraically closed fields, such as quasi-finite or even only dominant endomorphisms. Example \ref{EX1} above is a particular case of one of several series of interesting \'etale endomorphisms that we were able to construct, sometimes with the help of other smooth affine hypersurfaces that contain affine spaces as open dense subvarieties. Note that while most of our theorems are stated as existential results, our proofs are ultimately constructive, even if it is not always possible to write down explicit formulas.

In the next section we give precise definitions and describe our main results, that are proven in the subsequent sections. Finally, in the last section we propose a number of open problems, questions, and conjectures that naturally arose from our work. Perhaps, the most interesting open problems and questions are the following.

\medskip
{\bf (1)} Can some of our methods (resp.\ examples) be tweaked to prove (resp.\ disprove) Conjecture \ref{CJ2.5} in dimension $2$?

\smallskip
{\bf (2)} For each prime $p$, classify all integers $m\ge p+1$ that are relatively prime to $p$ and for which there exists an \'etale endomorphism of the affine plane over an algebraically closed field of characteristic $p$ of geometric degree $m$.

\smallskip
{\bf (3)} Can some of our examples and methods be used to prove or disprove the Surjectivity Conjecture over $\mathbb C$ in dimension $2$?

\smallskip
{\bf (4)} For each prime $p$, classify all non-decreasing eventually constant sequences $(a_n)_{n\in\mathbb N}$ of natural numbers with the properties that $a_0=0$, the first order difference sequence $(a_{n+1}-a_n)_{n\in\mathbb N}$ is non-increasing, and there exists an \'etale endomorphism $e$ of an affine plane over an algebraically closed field of characteristic $p$ such that for every $n\in\mathbb N$, the complement of the image of $e^n$ has $a_n$ points. 

\section{Notation, definitions, and description of main results}\label{S2}

\phantomsection{Let $K$ be an arbitrary algebraically closed field and let $\chr(K)$ be its characteristic. For a polynomial $f(t)\in K[t]\setminus\{0\}$\index{polynomial}, let $\z(f)\in\mathbb N$ be the number of distinct roots of $f$ in $K$.}\label{PH0.5} 

Recall that $p$ is a prime number. Let $k$ be an algebraically closed field with $\chr(k)=p$. 

\phantomsection{For $(l,m)\in\mathbb N^2$, let}\label{PH1} 
$$\llbracket l,m\rrbracket:=\{i\in\mathbb N|l\le i\le m\}$$
(so $\llbracket l,m\rrbracket=\emptyset$ if $l>m$). 

\phantomsection{Let $R$ be a commutative $\mathbb Z$-algebra. Let $R^{\ast}$ be the multiplicative group of units of $R$. Let $n\in\mathbb N^\ast$. For indeterminates $x_1,\ldots,x_n$, the $R$-algebra $R[x_1,\ldots,x_n]$ is $\mathbb N$-graded by $\deg(1)=0$ and $\deg(x_i)=1$ for each $i\in \llbracket1,n\rrbracket$.}\label{PH2}

\phantomsection{All the standard schemes or group schemes}\label{PH3}
$$\mathbb A^n,\;\mathbb P^n,\;\SL_2,\;\GL_2,\;\mathbb G_a,\;\mathbb G_m,\;\text{etc.}$$ 
\phantomsection{are over $\Spec \mathbb Z$ and we use their pullbacks}\label{PH4}
$$\mathbb A^n_R,\;\mathbb P^n_R,\;\SL_{2,R},\;\GL_{2,R},\;\mathbb G_{\a,R},\;\mathbb G_{\m,R},\;\text{etc.}$$ 
\phantomsection{to $\Spec R$. Thus, with $w$, $x$, $y$, and $z$ as indeterminates, we have}\label{PH5}
$$\mathbb A^n=\Spec(\mathbb Z[x_1,\ldots,x_n]),\;\SL_2=\Spec\bigl(\mathbb Z[w,x,y,z]/(wz-xy-1)\bigr),$$
$$\GL_2=\Spec\bigl(\mathbb Z[w,x,y,z,(wz-xy)^{-1}]\bigr),\;\mathbb G_a=\Spec(\mathbb Z[t]),$$
$$\;\mathbb G_m=\Spec(\mathbb Z[t,t^{-1}]),\;\mathbb A^n_R=\Spec(R[x_1,\ldots,x_n]),\;\text{etc}.$$

\phantomsection{Let $\End_n(R)$ be the monoid of endomorphisms\index{endomorphism} of $\mathbb A^n_R$ over $\Spec R$ under the composition (denoted multiplicatively) and let $\GA_n(R)$ be its subgroup of invertible elements, i.e., of automorphisms of $\mathbb A^n_R$ over $\Spec R$.}\label{PH6} 

If $m\in\mathbb N^\ast$ and $(g_1,\ldots,g_m)\in R[x_1,\ldots,x_n]^m$, by the degree of the $m$-tuple $(g_1,\ldots,g_m)$ or of the system $\mathcal V$ of $m$ polynomial equations in $n$ indeterminates
$$g_i(x_1,\ldots,x_n)=0,\ \ \quad i\in \llbracket1,m\rrbracket$$
\phantomsection{it defines, we mean}\label{PH7} 
$$\pi(g_1,\ldots,g_m)=\pi(\mathcal V):=\max\bigl(\deg(g_i)|i\in \llbracket1,m\rrbracket\bigr)\in\mathbb N\cup\{-\infty\}.\footnote{By convention, the degree of the zero polynomial is $-\infty$.}$$
If $m=n$, by the endomorphism of $\mathbb A^n_R$ over $\Spec R$ defined by $(g_1,\ldots,g_n)$ we mean the endomorphism $\e(g_1,\ldots,g_n)\in\End_n(R)$ defined by the $R$-algebra endomorphism of $R[x_1,\ldots,x_n]$ that maps $x_i$ to $g_i$ for each $i\in \llbracket1,n\rrbracket$.

We use the following formalism of {\it invariants}\index{invariant} for such endomorphisms.

\begin{definition}\label{D1}
Let $M$ be a subset of $\End_n(R)$ stable under left and right translations by elements of $\GA_n(R)$. Let $\mathcal N$ be an arbitrary set. A function $\Im:M\rightarrow\mathcal N$ is called an $\mathcal N$-valued\index{invariant!$\mathcal N$-valued invariant}:

\medskip
{\bf (1)} inner invariant\index{invariant!inner invariant} for $M$ if $\Im(aea^{-1})=\Im(e)$ for each pair $(a,e)\in\GA_n(R)\times M$;

\smallskip
{\bf (2)} left (resp.\ right) invariant\index{invariant!left invariant}\index{invariant!right invariant}  for $M$ if $\Im(ae)=\Im(e)$ (resp.\ $\Im(ea)=\Im(e)$) for each pair $(a,e)\in\GA_n(R)\times M$;

\smallskip
{\bf (3)} invariant for $M$ if $\Im(aeb)=\Im(e)$ for each triple $(a,b,e)\in\GA_n(R)^2\times M$.
\end{definition}

\phantomsection{In practice, $M$ is a set in the chain of inclusions}\label{SETS}
$$\GA_n(R)\subset\cdots\subset\EE_n(R)\subset \QFE_n(R)\subset\QF_n(R)\subset\D_n(R)\subset\End_n(R),$$
\phantomsection{where $\D_n(R)$ is the submonoid of all endomorphisms that are fiberwise dominant, $\QF_n(R)$ is the submonoid of all quasi-finite endomorphisms of $\mathbb A^n_R$, $\QFE_n(R)$ is the submonoid of all quasi-finite endomorphisms of $\mathbb A^n_R$ that are fiberwise generically \'etale, $\EE_n(R)$ is the submonoid of all \'etale endomorphisms (also called Kraus--Keller maps) of $\mathbb A^n_R$, and the dots stand for extra sets to be introduced in the case when $R$ is an $\mathbb F_p$-algebra. Often $\mathcal N$ is a subset of either $\mathbb N\cup\{-\infty,\infty\}$ (ordered as a subset of the extended real line) or $\mathcal P_f(\mathbb N)$ (the set of finite subsets of $\mathbb N$).}\label{PH7.9} 

Though our emphasis is on the case when $n\le 3$ and $R=k$, whenever possible we work with an arbitrary $n$ and $R=K$. 

\begin{definition}\label{D2}
Assume that $R$ is an $\mathbb F_p$-algebra.

\medskip
{\bf (1)} By an {\it $F_n$-system}\index{Frobenius system}\index{Frobenius system!$F_n$ system} (Frobenius system in $n$ indeterminates)\index{Frobenius system!in $n$ indeterminates} $\mathcal S$ over $R$ we mean a system of $n$ polynomial equations in $n$ indeterminates of the form 
$$x_i+f_i(x_1^p,\ldots,x_n^p)=0,\ \ \quad  i\in \llbracket1,n\rrbracket,$$
where $(f_1,\ldots,f_n)\in R[x_1,\ldots,x_n]^n$.
If $S$ is an $\mathbb F_p$-subalgebra of $R$, we say that $\mathcal S$ is defined over $S$ if we have $f_i\in S[x_1,\ldots,x_n]$ for all $i\in \llbracket1,n\rrbracket$. By the endomorphism defined by $\mathcal S$ we mean the endomorphism 
$$e_{\mathcal S}:=e\bigl(x_1+f_1(x_1^p,\ldots,x_n^p),\ldots,x_n+f_n(x_1^p,\ldots,x_n^p)\bigr)\in\End_n(R).$$ 
If $R=k$, let $\s(\mathcal S)\in\mathbb N$ be the number of solutions of $\mathcal S$ in $k^n$. 

\smallskip
{\bf (2)} By a {\it basic endomorphism}\index{endomorphism!basic endomorphism} of $\mathbb A^n_R$ we mean an endomorphism $e\in\End_n(R)$ such that an element of $\GA_n(R)e\GA_n(R)$ is defined by an $F_n$-system over $R$.
\end{definition}

From the very definitions, each basic endomorphism of $\mathbb A^n_R$ is \'etale as, referring to Definition \ref{D2}(1), for the $R[y_1,\ldots,y_n]$-algebra
$$R[x_1,y_1,\ldots,x_n,y_n]/\bigl(y_i-x_i-f_i(x_1^p,\ldots,x_n^p)|i\in \llbracket1,n\rrbracket\bigr)=R[x_1,\ldots,x_n],$$
the $2n$ differentials 
$$dy_i,d\bigl(y_i-x_i-f_i(x_1^p,\ldots,x_n^p)\bigr)=dy_i-dx_i\in\Omega_{R[x_1,y_1,\ldots,x_n,y_n]/R}\;\;\textup{with}\;\; i\in \llbracket1,n\rrbracket$$
\phantomsection{are linearly independent even at all points of $\Spec(R[x_1,y_1,\ldots,x_n,y_n])\cong\mathbb A^{2n}_R$. Thus for the set $\BE_n(R)$ of basic endomorphisms of $\mathbb A^n_R$ we have inclusions}\label{PH8}
$$\GA_n(R)\subset\BE_n(R)\subset\EE_n(R).$$
Clearly, $\BE_n(R)$ is stable under left and right translations by elements of $\GA_n(R)$.

We study the following list of invariants or inner invariants of suitable submonoids of $\End_n(K)$, formed by endomorphisms whose domain and codomain are denoted also as $\mathbb A^n_{K,\s}=\Spec(K_{\s}[x_1,\ldots,x_n])$ and $\mathbb A^n_{K,\t}=\Spec (K_{\t}[x_1,\ldots,x_n])$ (respectively), with the indices $\s$ and $\t$ standing for source and target (respectively).

\medskip
{\bf (I)} \phantomsection{The {\it algebraic degree invariant}\index{algebraic degree invariant}}\label{PH10} 
$$\pi=\pi_{n,K}:\End_n(K)\rightarrow\mathbb N\cup\{-\infty\}$$ 
is defined by the rule: $\pi(e)$ is the minimum of the set of all $\pi(g_1,\ldots,g_n)$ with $(g_1,\ldots,g_n)\in K_{\s}[x_1,\ldots,x_n]^n$ varying subject to the condition that it defines an endomorphism in the double coset $\GA_n(K)e\GA_n(K)$. We similarly define the inner\index{algebraic degree invariant!inner algebraic degree invariant}, left\index{algebraic degree invariant!left algebraic degree invariant}, and right\index{algebraic degree invariant!right algebraic degree invariant} algebraic degree invariants by 
$$\pi_{\i}=\pi_{\textup{inner},n,K}:\End_n(K)\rightarrow\mathbb N\cup\{-\infty\},$$
$$\pi_{\l}=\pi_{\textup{left},n,K}:\End_n(K)\rightarrow\mathbb N\cup\{-\infty\},$$
and
$$\pi_{\r}=\pi_{\textup{right},n,K}:\End_n(K)\rightarrow\mathbb N\cup\{-\infty\}$$ 
(respectively) via the rules: $\pi_{\i}(e)$, $\pi_{\l}(e)$, or $\pi_{\r}(e)$ is the minimum of the set of all $\pi(g_1,\ldots,g_n)$ with $(g_1,\ldots,g_n)\in K_{\s}[x_1,\ldots,x_n]^n$ varying subject to the condition that it defines an endomorphism in the conjugacy class $\{aea^{-1}|a\in\GA_n(K)\}$, in the left coset $\GA_n(K)e$, or in the right coset $e\GA_n(K)$ (respectively).

These types of algebraic degrees measure how efficient the $n$-tuple $\pi(g_1,\ldots,g_n)$ defines suitable cosets of $e$ from the perspective of computing other invariants. 

\smallskip
{\bf (II)} \phantomsection{The {\it geometric degree invariant}\index{geometric degree invariant}}\label{PH11}
$$\deg=\deg_{n,K}:\D_n(K)\rightarrow\mathbb N^\ast$$ 
is a multiplicative invariant defined by the rule: $\deg(e)$ is the separable degree of the finite field extension $K_{\t}(x_1,\ldots,x_n)\rightarrow K_{\s}(x_1,\ldots,x_n)$ induced by $e$; so, if $e\in\QFE_n(K)$, then $\deg(e)$ is the degree of the mentioned field extension.

\phantomsection{For a scheme $W$, let $W_{\red}$ be its associated reduced scheme.}\label{PH11a} 

\smallskip
{\bf (III)} \phantomsection{The {\it non-generic fibers set invariant}\index{non-generic fibers set invariant}}\label{PH12}  
$$\mathcal F=\mathcal F_{n,K}:\QF_n(K)\rightarrow\mathcal P_f(\mathbb N)$$ 
is defined by 
$$\mathcal F(e):=\{s\in \llbracket0,\deg(e)-1\rrbracket|\exists P\in\mathbb A^n_{K,\t}(K)\;\textup{such that}\;\bigl(e^{-1}(P)\bigr)_{\red}\cong\Spec(K^s)\}.$$

We recall that if $e\in\EE_n(K)$, then $e$ is finite iff the set $\mathcal F(e)$ is empty by \cite{Gro5}, Cor.\ (18.2.9). 

\phantomsection{In all that follows by a variety $X$ over $K$ we mean a reduced scheme of finite type over $\Spec K$. The identity automorphism of $X$ is denoted as $1_X$. For a non-empty variety, let $\Irr(X)$ be its set of irreducible components, let $\Sing(X)$ be its singular locus, let $\Reg(X):=X\setminus\Sing(X)$ be its regular locus, and let $\dim(X)$ be its dimension; if $X$ is affine, by its embedding dimension $\edim(X)$ we mean the smallest $n\in\mathbb N$ such that $X$ is isomorphic to a closed subvariety of $\mathbb A^n_K$. If $X$ is regular, let $T_X$ be the tangent bundle over $X$. If $X$ is equidimensional of dimension $1$ (resp.\ $2$), we call $X$ a curve (resp.\ surface) over $K$.}\label{PH12a}% (by convention, $\dim(\emptyset)=-\infty$)

\smallskip
{\bf (IV)} \phantomsection{The {\it iterate inner invariant}\index{iterate invariant}}\label{PH13} 
$$\iota=\iota_{n,K}:\End_n(K)\rightarrow\mathbb N\cup\{\infty\}$$ 
can be defined for each endomorphism $e:\mathcal B\rightarrow\mathcal B$ of a set $\mathcal B$: if there exists $m\in\mathbb N$ such that $\Imm(e^m)=\Imm(e^{m+n})$ for all $n\in\mathbb N$ (equivalently, for $n=1$), then $\iota(e)$ is the smallest such $m$; if no such $m$ exists, then $\iota(e)$ is $\infty$. 

\phantomsection{By the {\it complement} of $e$ we mean the subset $\Gamma_e:=\mathcal B\setminus\Imm(e)$. If $d:X\rightarrow X$ is an endomorphism of a variety $X$ over $K$, then $\iota(d)=\iota\bigl(d(K)\bigr)$ and its complement $\Gamma_d$ is endowed with its reduced subscheme structure if it is locally closed.}\label{PH13a}

Note that for each $s\in\mathbb N^{\ast}$ we have $\iota(e)\le s\iota(e^s)$.

\smallskip\label{PH14}
{\bf (V)} The {\it gap (or missed points) invariant}\index{gap (or missed points) invariant} 
$$\varphi=\varphi_{n,K}:\QF_n(K)\rightarrow\mathbb N\cup\{\infty\}$$ 
can be also defined for each endomorphism $e:\mathcal B\rightarrow\mathcal B$ of a set $\mathcal B$: if $\Gamma_e$ is finite then $\varphi(e)$ is the cardinality of $\Gamma_e$, and if $\varphi_e$ is infinite then $\varphi(e)=\infty$. 

Lemma \ref{L6}(2) recalls the well-known fact that for $e\in QF_2(K)$ we have $\varphi(e)\in\mathbb N$.

\phantomsection{If $\bigl(\iota(e),\varphi(e)\bigr)\in\mathbb N^2$, then $\varphi\bigl(e^{\iota(e)}\bigr)$ is the smallest cardinality of a subset $\mathcal C\subset\mathcal B$ with the property that $e$ induces a surjective endomorphism of $\mathcal B\setminus\mathcal C$.}\label{PH14c}

To define the last invariant we first need the following definition.

\begin{definition}\label{D3}
{\bf (1)} By a Jacobian variety\index{Jacobian variety} over $K$ of dimension $n\ge 2$ we mean an endomorphism $\psi:X\rightarrow X$ of a normal connected variety $X$ over $K$ of dimension $n$ such that $\psi(X)$ is an open subvariety of $X$, the induced morphism $X\rightarrow\psi(X)$ is finite, and the restriction $\psi|\psi(X)$ of $\psi$ to $\psi(X)$ is \'etale. 

\smallskip
{\bf (2)} For a Jacobian variety $\psi:X\rightarrow X$, let $X^{\et}$ be the \'etale locus of $\psi$. 

\medskip\noindent
{\bf (2.a)} We say that $\psi$ is \'etale\index{Jacobian variety!\'etale} if $\psi$ is so, i.e., if $X^{\et}=X$. 

\smallskip\noindent
{\bf (2.b)} We say that $\psi$ is regular\index{Jacobian variety!regular} or affine\index{Jacobian variety!affine} if $X$ is so. 

\smallskip\noindent
{\bf (2.c)} By the class rank of $\psi$ we mean the number $\rho(\psi)$ of irreducible components of $\Gamma_{\psi}$ (i.e., the cardinality of $\Irr_{\Gamma_{\psi}}$). 

\smallskip\noindent
{\bf (2.d)} By the \'etale class rank of $\psi$ we mean the number $\rho_{\et}(\psi)$ of irreducible components of $\Gamma_{\psi}$ whose generic points are in $X^{\et}$.

\smallskip
{\bf (3)} By an isomorphism\index{Jacobian variety!isomorphism between} (resp.\ a quasi-isomorphism\index{Jacobian variety!quasi-isomorphism between}) between two Jacobian varieties $\psi:X\rightarrow X$ and $\psi_1:X_1\rightarrow X_1$ over $K$ we mean an isomorphism $\tilde\imath:X\rightarrow X_1$ (resp.\ a pair $(\tilde\imath_{\s},\tilde\imath_{\t})$ of isomorphisms $\tilde\imath_{\s},\tilde\imath_{\t}:X\rightarrow X_1$) over $K$ such that $\psi_1\circ\tilde\imath=\tilde\imath\circ\psi$ (resp.\ $\psi_1\circ\tilde\imath_{\s}=\tilde\imath_{\t}\circ\psi$ and $\tilde\imath_{\s}\bigl(\Imm(\psi)\bigr)=\Imm(\psi_1)$); we identify $\tilde\imath$ with the quasi-isomorphism $(\tilde\imath,\tilde\imath)$. 
\end{definition}%EXAMPLE

\phantomsection{Let $\mathbb J_n(K)$ (resp.\ $\mathbb J^{\q}_n(K)$) be the set of isomorphism (resp.\ quasi-isomorphism) classes $[\psi]$ of Jacobian varieties $\psi$ over $K$ of dimension $n$.}\label{EXTRA6}

\smallskip
{\bf (VI)} \phantomsection{For $n\ge 2$, the {\it Jacobian variety inner invariant}\index{Jacobian variety!Jacobian variety inner invariant}}\label{PH15} 
$$\Psi=\Psi_{n,K}:\EE_n(K)\rightarrow\mathbb J_n(K)$$
and its associated {\it class rank}\index{class rank} (resp.\ {\it \'etale class rank})\index{class rank!\'etale} {\it invariant}
$$\rho=\rho_{n,K}:\EE_n(K)\rightarrow\mathbb N\;\;\textup{(resp.}\; \rho_{\et}=\rho_{\et,n,K}:\EE_n(K)\rightarrow\mathbb N)$$
are defined as follows. Let $X_e$ be the normalization of $\mathbb A^n_{K,\t}$ in the field of fractions of $\mathbb A^n_{K,\s}$; we have an open embedding $\imath_e:\mathbb A^n_{K,\s}\rightarrow X_e$ by Zariski's Main Theorem. Let $\psi_e$ be the composite of the finite morphism $X_e\rightarrow \mathbb A^n_{K,\t}$ with $\imath_e:\mathbb A^n_{K,\t}=\mathbb A^n_{K,\s}\rightarrow X_e$; so $\rho(e):=\rho(\psi_e)$ is the number of irreducible components (they have dimension $n-1$, see \cite{Gro6}, Exp.\ V, Ex.\ 3.4) of the complement $X_e\setminus\Imm(\imath_e)=\Gamma_{\psi_e}$ and we call it the class rank of $e$ and $\rho_{\et}(e):=\rho_{\et}(\psi_e)$ is the number of such irreducible components whose generic points are contained in the \'etale locus of $\psi_e$ and we call it the \'etale class rank of $e$.\footnote{It is easy to see that the divisor class group $Cl(X_e)$ is isomorphic to $\mathbb Z^{\rho(e)}$.\label{FOOT4}} We call $\psi_e:X_e\rightarrow X_e$ the Jacobian variety associated to $e$\index{Jacobian variety!associated to \'etale endomorphism}. Let 
$$\Psi(e):=[\psi_e].$$

The diagram 
\begin{equation}\label{EQ00}
\xymatrix{
{\mathbb A_{K,s}^n} \ar[r]^{\imath_e} \ar[dr]^{e} & X_e \ar[d]^{\textup{finite}} \ar[dr]^{\psi_e} \\
 & {\mathbb A_{K,t}^n} \ar[r]^{\imath_e} & X_e
}
\end{equation}
summarizes the notation introduced which will be used throughout the paper.

If $\rho(e)>0$ (i.e., $e\neq\psi_e$), then $\iota(\psi_e)=\iota(e)+1$.

\phantomsection{Note that $e$ and $\psi_e$ determine each other and the image of $\Psi$ is in natural bijection to the quotient set of the conjugacy relation on $\EE_n(K)$ defined by the conjugacy action of $\GA_n(K)$ on it. Similarly, the image of the {\it Jacobian variety invariant}}\index{Jacobian variety!Jacobian variety invariant}\label{PH15a}
$$\Psi^{\q}=\Psi^{\q}_{n,K}:\EE_n(K)\rightarrow\mathbb J^{\q}_n(K)$$
defined by the same rule $\Psi^{\q}(e):=[\psi_e]$ is in natural bijection to the set of double cosets $\GA_n(K)\backslash\EE_n(K)/\GA_n(K)$.

As the local rings of $X_e$ of dimension $2$ are Cohen--Macaulay, the finite morphism $X_e\rightarrow\mathbb A^n_{K,\t}$ is flat in codimension at most $2$ by \cite{Gro3}, Prop.\ (6.1.5) or \cite{Ma}, Ch.\ 8, Thm.\ 23.1. Thus, if $n=2$, then $\psi_e$ is flat and we call $\psi_e$ a Jacobian surface.\index{Jacobian variety!Jacobian surface} 

The first goal of the paper is to study extensively the main six (inner) invariants introduced above for $e\in\EE_n(k)$; other invariants that help in or complement the study of the main six ones are introduced when required. 

\phantomsection{For a prime $p$ and $q\in\mathbb N^{\ast}$, let $\mathbb F_{p^q}$ be a finite field with $p^q$ elements.}\label{PH15c}

Our basic results for $n=2$ are three explicit examples (Examples \ref{EX7}, \ref{EX8}, and \ref{EX11}) whose existential essences are combined here as one theorem.

\begin{theorem}\label{T1}
{\bf (1)} For each $l\in \mathbb N$, there exists an $F_2$-system $\mathcal U$ over $k$ with $\s(\mathcal U)=l<\deg(e_{\mathcal U})$ (i.e., there exists $e\in\BE_2(k)$ such that $l\in\mathcal F(e)$).

\smallskip
{\bf (2)} Let the triple $(m,q,s)\in (\mathbb N^{\ast})^2\times\mathbb N$ be such that $p^q\ge s+1$. Then there exists $e\in\EE_2(k)$ defined over $\mathbb F_{p^q}$ and with invariants $\deg(e)=pm$, $\varphi(e)=s$, $\psi_e$ \'etale, $\rho_{\et}(e)=\rho(e)=(pm-1)(s+1)$, and $\pi(e)\le\pi_{\l}(e)=\pi_{\r}(e)=(3s+2)(pm-1)-1$. Moreover, $\mathcal F(e)$ is $\{1\}$ if $s=0$ and $pm=2$, is $\mathcal F(e)=\{0,1\}$ if $s\ge 1$ and $pm=2$, is $\{1,2,pm-1\}$ if $s=0$ and $pm>2$, and is $\{0,1,2,pm-1\}$ if $s\ge 1$ and $pm>2$.

\smallskip
{\bf (3)} Let $(m,q)\in (\mathbb N^\ast)^2$ be such that $p$ divides $mq$. Then there exists $e\in\EE_2(k)$ defined over $\mathbb F_p$ with invariants $\deg(e)=mq$, $\varphi(e)=mq-1$, $\mathcal F(e)=\{0,1,mq-1\}$, $\pi(e)\le\pi_{\l}(e)=\pi_{\r}(e)=mq+\max\bigl((mq-1)q,m\bigr)\ge 2mq$, and $\psi_e$ regular. If $mq>2$ we have $\iota(e)=1$, $\rho_{\et}(e)=mq-1$, $\rho(e)=mq$, and $\psi_e$ is non-\'etale, and if $mq=2$ (so $p=2$) we have $\iota(e)=2$, $\rho_{\et}(e)=\rho(e)=2$, and $\psi_e$ is \'etale. Moreover, if $p$ divides both $m$ and $q$ then $e\in\BE_2(k)$ and if $q=1$ then $\pi_{\i}(e)=\pi_{\l}(e)=\pi_{\r}(e)=2m$.
\end{theorem}

\phantomsection{A key novelty of Theorem \ref{T1}(1) is the case $l=0$. To compute the number of solutions in $k^2$ of an $F_2$-system $x+f^p=0=y+g^p$ over $k$, so $(f,g)\in k[x,y]^2$, one looks at the surjective \'etale (first projection) morphism}\label{PH16b} 
$$\Spec\bigl(k[x,y]/(y+g^p)\bigr)\rightarrow\Spec(k[x])$$ 
and to get examples of inconsistent $F_2$-systems over $k$ one needs to choose $(f,g)$ so that $x+f^p+(y+g^p)\in k[x,y]/(y+g^p)$ is a unit. 

\phantomsection{First, one is led to identify integral curves $C$ on $\Spec (k[x,y]\bigl[\frac{1}{x+y^p}\bigr])$ whose normalizations have open subvarieties that map \'etale surjectively onto the affine line $\Spec (k[x])$ via the natural first projection. There exist many such curves including the ones given by equations $y^l(1+y^{pm})=1$ with $(l,m)\in(\mathbb N\setminus p\mathbb N)\times\mathbb N^{\ast}$ (cf.\ Example \ref{EX28}). The simplest such curve $C$ is defined by the equation $y(x+y^p)=1$ and is isomorphic to $\Spec (k[t,t^{-1}])$ via the finite \'etale morphism}\label{PH16}
$$\Sigma:\mathbb G_{\m,k}\rightarrow\mathbb A^1_k=\Spec(k[x_1])$$
over $\Spec k$ defined at the level of $k$-algebras by the rule: $x_1\mapsto t^{-1}-t^p$. 

Second, one is led to find ways to dominate such a `good' curve $C$ via a morphism $\Spec\bigl(k[x,y]/(y+g^p)\bigr)\rightarrow C$ that commutes with the first projections.

Theorem \ref{T1}(2) gives that for all $m\in\mathbb N^{\ast}$ we have an identity 
\begin{equation}\label{EQ1.1}
\{\varphi(e)|e\in\EE_2(k),\deg(e)=pm,\psi_e\;\textup{is \'etale}\}=\mathbb N.
\end{equation} 
Theorem \ref{T6} proves the identity $\{\varphi(e)|e\in\EE_2(k),\psi_e\;\textup{is non-\'etale}\}=\mathbb N$. Directly from either one of these identities we get (see Section \ref{S18}) the following result.

\begin{corollary}\label{C1}
Let $X$ be a non-empty affine variety over $k$ and let $l:=\edim(X)$. Let $(n,q)\in (\mathbb N^{\ast})^2$ with $n\ge l+2$. Then there exist $e\in\EE_n(k)$ with $\iota(e)\in\mathbb N^{\ast}$ and a sequence $q=q_1\le q_2\le\cdots\le q_{\iota(e)}$ of integers such that for each $i\in \llbracket1,\iota(e)\rrbracket$ the complement $\Gamma_{e^i}$ is a closed subvariety of $\mathbb A^n_{k,\t}$ isomorphic to $q_i$ copies of $X$. 
\end{corollary}

If $X=\mathbb A^l_k$, then the inequality $n\ge l+2$ of Corollary \ref{C1} cannot be improved by Lemma \ref{L6}(2); moreover, if $n=l+2$ then one can exceptionally take $e$ to be a cartesian product $d\times 1_X$ with $d\in\EE_2(k)$ such that $\varphi(d)=q$.

\phantomsection{If $\kappa$ is a field of characteristic $p$ and $e\in\End_n(\kappa)$ is defined by an $n$-tuple $(g_1,\ldots,g_n)\in \kappa[x_1,\ldots,x_n]^n$, then it is well-known that the following three statements are equivalent (see \cite{No}, Thm.\ (6.2); see also \cite{Ni}, Thm.\ 1.4).}\label{PH16a}

\medskip
{\bf (N1)} We have $e\in\EE_n(\kappa)$.

\smallskip
{\bf (N2)} We have $\kappa[x_1,\ldots,x_n]=\kappa[x_1^p,\ldots,x_n^p,g_1,\ldots,g_n]$.

\smallskip
{\bf (N3)} The $n$-tuple $(g_1,\ldots,g_n)$ is a $p$-basis\index{$p$-basis} of $\kappa[x_1,\ldots,x_n]$, i.e., 
$$\Bigl(\prod_{i=1}^n g_i^{m_i}\bigl|m_i\in \llbracket0,p-1\rrbracket\; \forall\; i\in \llbracket1,n\rrbracket\Bigr)$$ 
is a $\kappa[x_1^p,\ldots,x_n^p]$-basis of the $\kappa[x_1^p,\ldots,x_n^p]$-module $\kappa[x_1,\ldots,x_n]$.

\medskip
Based on the equivalence $(N1)\Leftrightarrow (N3)$ we have the following interpretation of the inclusion $\{0,1\}\subset \mathcal F(e)$ of the case $(m,q,s)=(p,1,1)$ of Theorem \ref{T1}(2) or of the case $mq=p$ of Theorem \ref{T1}(3) in terms of $p$-bases (see Subsection \ref{S14.1} for the case $mq=p$ that involves smaller algebraic degrees).

\begin{corollary}\label{C2}
Let $\kappa$ be an arbitrary field of characteristic $p$, let $n\ge 2$ be an integer, and let $\mathcal R:=\kappa[x_1,\ldots,x_n]$. Then the following properties hold.

\medskip
{\bf (1)} There exists a $p$-basis $(g_1,\ldots,g_n)$ of $\mathcal R$ such that the degree of the finite field extension $\kappa(g_1,\ldots,g_n)\rightarrow \kappa(x_1,\ldots,x_n)$ is $p$ and moreover we have an equality $(g_1,\ldots,g_n)=(x_1,\ldots,x_n)$ of ideals of $\mathcal R$.

\smallskip
{\bf (2)} There exists a $p$-basis $(g_1,\ldots,g_n)$ of $\mathcal R$ such that the degree of the finite field extension $\kappa(g_1,\ldots,g_n)\rightarrow \kappa(x_1,\ldots,x_n)$ is $p$ and $\sum_{i=1}^n \mathcal Rg_i=\mathcal R$ (i.e., and the $n$-tuple $(g_1,\ldots,g_n)\in\mathcal R^n$ is unimodular).\end{corollary}

Conjecture \cite{Ad}, Conj.\ 3.4 is not properly formulated but it was meant to be the characteristic $p$ analog of the `Generalized Jacobian Conjecture' of Bass in characteristic $0$ (see \cite{Ba2}, Rmk.\ after Prop.\ (4.3)) which was disproved by Kulikov in dimension $2$ (see \cite{Kul}, Sect.\ 3). The following theorem is a $3$- and $2$-dimensional non-units analogs of the \'etale morphism $\Sigma$ that disproves \cite{Ad}, Conj.\ 3.4 (in all formulations), that is an output of the geometry of $\SL_{2,k}$, that is closely related to \cite{Kraft}, Sect.\ 2, Ex.\ 1, i.e., the hypersurface of $\Spec(\mathbb C[t,x,y,z])$ defined by the equation $x+x^2y+z^2+t^2=0$ equipped with a $\mathbb G_{\m,\mathbb C}$-action, and that refines \cite{Kam}, Thm.\ in the particular cases of simply connected semisimple groups $\SL_{2,\mathbb F_{p^q}}$ with $q\in\mathbb N$ and affine irreducible varieties $\mathbb A^3_k$ and $\mathbb A^2_k$ over $\Spec k$.

\begin{theorem}\label{T2}
The following properties hold for $n\in\{2,3\}$.

\medskip
{\bf (1)} Let $q\in\mathbb N^\ast$. There exists a finite Galois cover
$$\mathcal M_n:\mathbb X_k^n\rightarrow\mathbb A^n_k$$ 
of Galois group $\SL_2(\mathbb F_{p^q})$ between rational affine varieties over $k$ with group of units equal to $k^{\ast}$ which is defined over $\mathbb F_{p^q}$. Hence $\mathcal M_n$ factors through finite \'etale morphisms 
$$\mathbb Y_k^n\rightarrow\mathbb A^n_k\;\;\;\textup{and}\;\;\;\mathbb W^n_k\rightarrow\mathbb A^n_k$$
of degree $p^{2q}-1$ and $p^q+1$ (respectively) between rational affine varieties over $k$ with group of units equal to $k^{\ast}$ which are defined over $\mathbb F_{p^q}$.

\smallskip
{\bf (2)} If $p=2$ and $q=1$, then $\mathcal M_n$ factors through a finite Galois cover $\mathbb V_k^n\rightarrow\mathbb A^n_k$ of degree $2$ defined over $\mathbb F_2$, where $\mathbb V_k^n$ is rational, has group of units equal to $k^{\ast}$, and has a finite Galois cover $\mathbb X^n_k\rightarrow\mathbb V^n_k$ of degree $3$ (hence its prime-to-$2$ fundamental group is non-trivial).

\smallskip
{\bf (3)} If $p=3$ and $q=1$, then $\mathcal M_n$ factors through a finite Galois cover $\mathbb V_k^n\rightarrow\mathbb A^n_k$ of degree $3$ defined over $\mathbb F_3$, where $\mathbb V_k^n$ is rational, has group of units equal to $k^{\ast}$, and has a finite Galois cover $\mathbb X^n_k\rightarrow\mathbb V_k^n$ with Galois group isomorphic to the quaternion group $Q_8$ of order $8$ (hence its prime-to-$3$ fundamental group is non-trivial).

\smallskip
{\bf (4)} In dimension $2$, the morphisms of (1) to (3) are obtained from the analog ones in dimension $3$ via taking quotients through a left $\mathbb G_{\a,k}$-action defined over $\mathbb F_{p^q}$ (i.e., we have $\mathbb X^2_k=\mathbb G_{\a,k}\backslash \mathbb X^3_k$, $\mathcal M_2=\mathbb G_{\a,k}\backslash \mathcal M_3$, etc.).
\end{theorem}

The proof of Theorem \ref{T2} uses Nori's idea (see \cite{Kam}) to take $\mathcal M_3$ as the pullback of Lang $\SL_2(\mathbb F_{p^q})$-torsor $\SL_{2,k}\rightarrow \SL_{2,k}$ via a `good' morphism $\mathbb A^3_k\rightarrow \SL_{2,k}$. The finite Galois cover $\mathbb V_k^n\rightarrow\mathbb A^n_k$ of Theorem \ref{T2}(2) (resp.\ 2.7(3)) is of Artin--Schreier type by either \cite{Mi2}, Thm.\ 1.1 or Lemma \ref{F6}. Moreover, as $\mathbb V^n_k$ has group of units equal to $k^{\ast}$, the existence of a non-trivial cyclic cover of it of degree $3$ (resp.\ $8$) implies by Kummer theory that it is non-factorial.

\phantomsection{The classification of Jacobian varieties $\psi_e$ with $e\in\EE_n(K)$ would implicitly classify all finite covers of $\mathbb A^n_K$ of the form $X_e$ with $e\in\EE_n(K)$. If $e\in\EE_n(K)$ is such that $\psi_e$ is regular and non-\'etale, then $X_e$ is not a complete intersection (see Theorem \ref{T3+}(6)) and it is not clear even how to begin the classification of such $X_e$s.}\label{PH17}\footnote{Each $X_e$ admits a closed embedding into $\mathbb A^{n+d_e}_K$ (resp. $\mathbb A^n_K\times_{\Spec K}\mathbb P^{d_e}_K$) for a smallest $d_e\in\mathbb N^{\ast}$ which can be bounded from above (especially if $X_e$ is regular). But how to describe such an embedding when $d_e\ge 2$ (resp.\ $d_e\ge 1$) in a way that classifications are possible is a fundamental open problem.} Thus, to begin with we introduce the following definition.

\begin{definition}\label{D6-}
Let $e\in\EE_n(K)$.

\medskip
{\bf (1)} By an approximation\index{approximation} of $X_e=\Spec A_e$ we mean any affine variety of the form $Z_e=\Spec B_e$ with $B_e$ a $K_{\t}[x_1,\ldots,x_n]$-subalgebra of $A_e$. 

\smallskip
{\bf (2)} The approximation $Z_e$ of $X_e$ is called flat\index{approximation!flat} if $B_e$ is a free $K_{\t}[x_1,\ldots,x_n]$-module and is called a model\index{approximation!model} of $X_e$ if $\Frac(B_e)=\Frac(A_e)$. 

\smallskip
{\bf (3)} A model $Z_e$ of $X_e$ is called a monomial model\index{approximation!monomial model} if there exists $y\in B_e$ such that $\bigl(y^i|i\in\llbracket0,\deg(e)-1\rrbracket\bigr)$ is a $K_{\t}[x_1,\ldots,x_n]$-basis of the $K_{\t}[x_1,\ldots,x_n]$-module $B_e$.
\end{definition}

Each monomial model is a flat model (approximation). Each model $Z_e=\Spec B_e$ of $X_e$ determines $X_e$ uniquely: $A_e$ is the normalization of $B_e$. But there exist many (flat) models of $X_e$ and they do not determine $e$. So we also consider the morphisms $Z_e\rightarrow\mathbb A^n_{K,\t}$ given by the inclusions $K_{\t}[x_1,\ldots,x_n]\subset B_e$: they determine the class $[\psi_e]\in\mathbb J^{\q}_n(K)$ if $X_e$ has a unique open subvariety $W$ which is isomorphic to $\mathbb A^n_K$ and for which the composite of the open embedding $W\rightarrow X_e$ with $\psi_e$ induces an \'etale morphism $W\rightarrow\mathbb A^n_{K,\t}$. 

So we first consider hypersurfaces in $\mathbb A^{n+1}_K$ from two points of view. 

\phantomsection{First, to study all monomial models of the $X_e$s with $e\in\EE_n(K)$, let}\label{EXT12} 
$$\Monomial_n(K):=\{f\in  K[x_1,\ldots,x_{n+1}]|f\; \textup{is irreducible and monic in}\; x_1\}.$$ 
\phantomsection{For each $f\in\Monomial_n(K)$, let $\mathbb H^n_f$ be the hypersurface in $\mathbb A^{n+1}_K$ which is the zero locus $f=0$ and let $\mathbb H^{n,\n}_f$ be its normalization. Let}\label{EXT13}
$$\AffMonomial_n(K):=\{f\in \Monomial_n(K)|\exists\; \textup{an open embedding}\; \mathbb A^n_K\rightarrow\mathbb H^n_f\}.$$

Second, for complementing the mentioned study and for identifying many $X_e$s with $e\in\EE_n(K)$ such that $\psi_e$ is \'etale, let
$$\Intr_n(K):=\{(f,g)\in K[x_1,\ldots,x_n]^2|g\neq 0,\; \textup{g.c.d.}(f,g)=1\}.$$
For $(f,g)\in\Intr_n(K)$, the zero locus
\begin{equation*}\label{EQ0}
f(x_1,\ldots,x_n)+g(x_1,\ldots,x_n)x_{n+1}=0
\end{equation*}
\phantomsection{is an integral rational hypersurface $\mathbb H^n_{f,g,K}$ in $\mathbb A^{n+1}_K$ but in general it is neither normal nor contains an open subvariety isomorphic to $\mathbb A^n_K$. Thus one is led to consider the following three subsets of $\Intr_n(K)$:}\label{PH18}
$$\Nointr_n(K):=\{(f,g)\in\Intr_n(K)|\mathbb H^n_{f,g,K}\;\textup{is normal}\},$$
$$\Smintr_n(K):=\{(f,g)\in\Intr_n(K)|\mathbb H^n_{f,g,K}\;\textup{is smooth}\},$$
$$\Affintr_n(K):=\{(f,g)\in\Intr_n(K)|\exists\; \textup{an open embedding}\; \mathbb A^n_K\rightarrow\mathbb H^n_{f,g,K}\}.$$

The hypersurfaces $\mathbb H^n_{f,g,K}$ with $(f,g)\in\Intr_n(K)$ such that $g=\prod_{i=1}^n x_i^{m_i}$ with $(m_1,\ldots,m_n)\in\mathbb N^n$ (i.e., $g$ is a primitive monomial) and with $f$ of certain types have been extensively studied in the literature from many points of views, including automorphisms, group actions, invariants, cancellation problems, and isomorphism classes. Some more general types of $g$s are considered in \cite{Du2} and \cite{GG}, e.g., in \cite{GG}, Thm.\ B, $g$ is a monic polynomial in one of the variables.

The main object of interest is the intersection $\Nointr_n(K)\cap\Affintr_n(K)$. One would like to classify (up to isomorphisms) all open embeddings $\mathbb A^n_K\rightarrow\mathbb H^n_{f,g,K}$ with $(f,g)\in\Nointr_n(K)\cap\Affintr_n(K)$ and, when $\chr(K)=p$, all finite \'etale covers $\mathbb H^n_{f,g,K}\rightarrow\mathbb A^n_K$ with $(f,g)\in\Smintr_n(K)\cap\Affintr_n(K)$. 

We note that the intersection $\Nointr_n(K)\cap\Affintr_n(K)$ can be used to study complete intersections in $\mathbb A^n_K$ isomorphic to certain $X_e$s with $e\in\EE_n(K)$ (e.g., see Examples \ref{EX34} and \ref{EX35}). However, it does not suffice to study $\mathbb J_n(K)$ intrinsically via hypersurfaces. 

In what follows we introduce a two polynomial parameters family of finite \'etale covers of $\mathbb A^2_k$. For the sake of generality, we introduce them over $\Spec R$ as follows.

\phantomsection{We say that a monic polynomial $f(t)\in R[t]$ is separable\index{polynomial!separable monic polynomial} if its reduction modulo each maximal ideal of $R$ is so. Let}\label{PH19} 
$$\Theta_R:=\{(f,g)\in (R[t])^2|f\;\textup{is monic separable},g\;\textup{is monic},f(0)\in R^{\ast}\};$$
e.g., for $n\ge 2$ we have $(t^{n-1}-1,1)\in\Theta_{\mathbb Z[\frac{1}{n-1}]}$. For $(f,g)\in\Theta_R$ we consider the affine smooth surface
$$\mathbb S^2_{f,g,R}:=\Spec\bigl(R[x,y,z]/(xf(x)+yg(y)z)\bigr)$$
over $\Spec R$; so $\mathbb S^2_{1,1,\mathbb Z}\cong\mathbb A^2$. For $n\ge 2$ we single out the affine smooth surface
$$\mathbb S^2_n:=\mathbb S^2_{t^{n-1}-1,1,\mathbb Z[\frac{1}{n-1}]}=\Spec\Bigl(\mathbb Z\Bigl[\frac{1}{n-1}\Bigr][x,y,z]/(x^n-x+yz)\Bigr)$$
over $\mathbb Z[\frac{1}{n-1}]$. We call
$$\mathbb S^2:=\mathbb S^2_2=\mathbb S^2_{t-1,1,\mathbb Z}=\Spec\bigl(\mathbb Z[x,y,z]/(x^2-x+yz)\bigr)$$
the arithmetic $2$-dimensional sphere (over $\Spec\mathbb Z$) due to the existence of isomorphisms
$\mathbb S^2_{\mathbb Z[\frac{1}{2}]}\cong \Spec\bigl(\mathbb Z[\frac{1}{2}][x,y,z]/(x^2+y^2-z^2-1)\bigr)$ and, with $i\in\mathbb C$, $i^2=-1$, 
$$\mathbb S^2_{\mathbb Z[i][\frac{1}{2}}\cong \Spec\Bigl(\mathbb Z[i]\Bigl[\frac{1}{2}\Bigr][x,y,z]/(x^2+y^2+z^2-1)\Bigr).$$

Over $K$, groups of automorphisms, group actions, and isomorphism classes related to certain subclasses of the surfaces $\mathbb S^2_{f,g,K}$ have been extensively studied in the literature starting with \cite{GD}, especially when $\chr(K)=0$. For instance, $\mathbb S^2_K$ is isomorphic to $\bigl(\mathbb P^1_K\times_{\Spec K}\mathbb P^1_K\bigr)\setminus\Delta_{\mathbb P^1_K}$, where $\Delta_{\mathbb P^1_K}$ is the image of the diagonal embedding $\mathbb P^1_K\rightarrow\mathbb P^1_K\times_{\Spec K}\mathbb P^1_K$ (e.g., see \cite{BvS}, Lem.\ 2.1; for $K=\mathbb C$ see also \cite{Wr2}, Prop.\ 1.1); the ind-group scheme $\Aut(\mathbb S^2_K)$ has been described as an amalgamated product in \cite{GD}, Sect.\ 10.1 if $\chr(K)\neq 2$ and in \cite{Lam}, Thm.\ 4 in general. 

For $(f,t^m)\in\Theta_K$ with $m\in\mathbb N$, the surfaces $\mathbb S^2_{f,t^m,K}$ are called (often special or classical) Danielewski surfaces in \cite{FM-J}, \cite{Du1}, etc.\ after the author of an unpublished 1989 manuscript which was first quoted in \cite{F} and in which the particular case $f(t)=t-1$ over $\mathbb C$ is studied in connection to the cancellation problem (by $\mathbb A^1_{\mathbb C}$); if $m=0$ (resp.\ if $m\ge 1$), generators of their groups of automorphisms (resp.\ generators of their groups of automorphisms and isomorphism classes) are described in \cite{M-L1}, Thm.\ (resp.\ in \cite{M-L2}, Thms.\ 1 and 2 for $K=\mathbb C$ and in \cite{Cr}, Thm.\ 4.2 and Cor.\ 4.3 applied to $h(x)=-1$). 

\phantomsection{Generalized Danielewski surfaces are zero loci $L(x,y)+y^nz=0$ (hypersurfaces) in $\mathbb A^3_K$ with $L(x,y)\in K[x,y]$ such that $\deg\bigl(L(x,0)\bigr)\ge 2$; over $\mathbb C$ they were introduced (in this generality) and classified in \cite{Pol}, Def.\ 2 and Sects.\ 4 and 5.}\label{PG19a} 

The terminology `generalized' is supported by two facts: (i) for $n\ge 1$ and $g\in K[t]$ such that $\deg(g)<n$ the zero loci $x^2+xg(y)+y^nz=0$ classify surfaces that are non-trivial line bundles over the affine line over $\Spec K$ with the origin doubled (see \cite{Wi}, Sect.\ 2), and (ii) over $\mathbb C$, for $n\ge 2$ and $(f,g)\in K[t]^2$ with $\deg(f)\ge 2$ and $\deg(g)\le n-1$, the surfaces that are the zero loci $f(x)+y^nz=0$ and $f(x)g(y)+y^nz=0$ are isomorphic by \cite{FM-J}, Thm.\ 2. 

There exist many morphisms between these surfaces over $\Spec R$, among which two types are introduced here. 

\phantomsection{First, we have an open embedding}\label{PH20} 
\begin{equation}\label{EQ0a}
\imath_{f,g,R}:\mathbb A^2_R\rightarrow\mathbb S^2_{f,g,R}
\end{equation}
defined on valued points by the rule 
\begin{equation}\label{EQ0b}
(x,y)\mapsto\bigl(xyg(y),y,-xf(xyg(y))\bigr);
\end{equation}
\phantomsection{its complement}\label{PH21}  
\begin{equation}\label{EQ0c}
\mathbb D_{f,g,R}^1:=\mathbb S^2_{f,g,R}\setminus\Imm(\imath_{f,g,R})
\end{equation} 
is the zero locus $f(x)=yg(y)=0$. Over a field $\mathcal K$ of characteristic not dividing $n-1$ that contains the $n-1$-th roots of $1$ and all zeros of $g(y)$ in an algebraic closure of $\mathcal K$, $\mathbb D_{t^{n-1}-1,g,\mathcal K}^1$ is a disjoint union of $(n-1)\z(tg)$ copies of $\mathbb A^1_{\mathcal K}$; e.g., $\mathbb D^1:=\mathbb D^1_{t-1,1,\mathbb Z}\cong \mathbb A^1$ and $\mathbb D^1_{t^{n-1}-1,1,\mathcal K}$ is isomorphic to $n-1$ copies of $\mathbb A^1_{\mathcal K}$. Thus for $n\le 2$ and $\deg(g)\in\mathbb N$ we get the following fact.

\begin{fact}\label{F1}
For $m\in\mathbb N^{\ast}$ (resp.\ $m=0$), there exists a (monic) connected smooth hypersurface in $\mathbb A^3_K$ which is the disjoint union of $\mathbb A^2_K$ and of $m$ copies of $\mathbb A^1_K$ and which is a finite flat cover of $\mathbb A^2_K$ of degree $2$ (resp.\ $1$). 
\end{fact}

The existence of the open embeddings $\imath_{f,g,K}$ implies that we have inclusions
$$\Theta_K\subset\Smintr_2(K)\cap\Affintr_2(K)$$
and
$$\{x_1f(x_1)+x_2g(x_2)x_3|(f,g)\in\Theta_K\}\subset\AffMonomial_2(K).$$

Second, if the pair $(f,g)\in\Theta_R$ is such that we have a product decomposition $f=f_1f_2$ of monic polynomials, then $(f_1,g)\in\Theta_R$ and the rule on valued points $(x,y,z)\mapsto \bigl(x,y,f_2(x)z\bigr)$ defines an open embedding 
$$\imath_{f_1,g;f,R}:\mathbb S^2_{f_1,g,R}\rightarrow\mathbb S^2_{f,g,R}$$ 
whose complement is the zero locus $f_2(x)yg(y)=0$. Note that 
$$\imath_{f_1,g;f,R}\circ\imath_{f_1,g,R}=\imath_{f,g,R}.$$

\phantomsection{Let}\label{PH22}  
$$c_{f,g,R}:\mathbb S_{f,g,R}^2\rightarrow \mathbb A^2_R$$ 
be the finite flat morphism defined by the rule on valued points $(x,y,z)\mapsto (y,z)$; the reduced closed subscheme $c_{f,g,R}(\mathbb D^1_{f,g,R})$ of $\mathbb A^2_R$ is the zero locus $x_1g(x_1)=0$.

\phantomsection{If $R$ is an $\mathbb F_p$-algebra and $m\in\mathbb N^{\ast}$, let}\label{PH23}  
$$\Theta^1_{R,m}:=\{(f,g)\in\Theta_R|f(t)+tf'(t)\in R^{\ast},\deg(f)=pm-1\}.$$
If $(f,g)\in\Theta^1_{k,m}$, then $c_{f,g,k}:\mathbb S_{f,g,k}^2\rightarrow\mathbb A^2_k$ is a finite \'etale cover of degree $pm$ which is Galois if $tf(t)$ is an additive polynomial; its composite with the open embedding $\imath_{f,g,k}:\mathbb A^2_k\rightarrow\mathbb S_{f,g,k}^2$ is a surjective non-finite endomorphism 
$$e_{f,g,k}:=c_{f,g,k}\circ\imath_{f,g,k}\in\EE_2(k)$$ 
with $\deg(e_{f,g,k})=pm$ which is defined over the $\mathbb F_p$-subalgebra of $k$ generated by the coefficients of $f(t)$ and $g(t)$ by the rule $(x,y)\mapsto \bigl(y,-xf(xyg(y))\bigr)$ on valued points and is finite over the complement of the zero locus $x_1g(x_1)=0$ in $\mathbb A^2_{k,\t}$.%GROUP1

For instance, if $(f,g)\in\Theta^1_{k,m}$, then the composite $\imath_{f,g,k}\circ c_{f,g,k}:\mathbb S_{f,g,k}^2\rightarrow\mathbb S_{f,g,k}^2$ is the \'etale Jacobian surface $\psi_{e_{f,g,k}}$ associated to $e_{f,g,k}$; so $\Gamma_{\psi_{e_{f,g,k}}}$ is isomorphic to $(pm-1)\z(tg)$ copies of $\mathbb A^1_k$ and we have $\rho_{\et}(e_{f,g,k})=\rho(e_{f,g,k})=(pm-1)\z(tg)$. This proves Theorem \ref{T1}(2) for $s=0$; the proof for $s>0$ is similar but using a different open embedding $\mathbb A^2_{k,\s}\rightarrow\mathbb S^2_{f,g,k}$ (see Proposition \ref{PR4.5}).

For each $(l,q)\in\mathbb N^2$ we construct (see Corollary \ref{C6}) \'etale endomorphisms of suitable affine open subvarieties of the surfaces $\mathbb S_{f,t,k}^2$ whose images have complements that are disjoint unions of $l$ copies of $\mathbb A^1_k$ and of $q$ points ($q$ copies of $\mathbb A^0_k$). 

If $\chr(K)=0$, then for each $(f,g)\in\Theta_K$ the surface $\mathbb S^2_{f,g,K}$ is simply connected as it has an open subvariety isomorphic to $\mathbb A^2_K$. In contrast, for each $(f,g)\in\Theta_k$ there exists a finite \'etale cover $\mathbb S^2_{f,g,k}\rightarrow\mathbb A^2_k$ and there exists $e\in\EE_2(k)$ such that $\psi_e$ is \'etale with $X_e\cong\mathbb S^2_{f,g,k}$ (see Proposition \ref{PR8}(1) and (2); recall that $\mathbb S^2_{1,g,k}\cong\mathbb A^2_k$). The existence of finite \'etale covers $\mathbb S^2_k\rightarrow\mathbb A^2_k$ shows that Kumar--Murthy--Nori Theorem (see \cite{Wr2}, Thm.\ 1.3) does not hold in characteristic $p$.

If $(f,g)\in\Smintr_n(K)$, then the tangent bundle $T_{\mathbb H^n_{f,g,K}}$ is trivial, hence each open affine non-empty subvariety $W$ of $T_{\mathbb H^n_{f,g,K}}$ is a complete intersection by Theorem \ref{T3}(2). If moreover $n=2$, $(f,g)\in\Theta_K$, and $W$ contains an open subvariety isomorphic to $\mathbb A^2_K$, then there exist examples in which $W\cong\mathbb H^2_{\tilde f,\tilde g,K}$ for a suitable pair $(\tilde f,\tilde g)\in\Smintr_2(K)$ and in addition either $W\cong\mathbb S^2_{\tilde f_1,\tilde g_1,K}$ for a suitable pair $(\tilde f_1,\tilde g_1)\in\Theta_K$ (for $f=f_1f_2$ and $(f,g)\in\Theta_K$, see $\Imm(\imath_{f_1,g;f,K})$) or $W\not\cong \mathbb S^2_{\tilde f_1,\tilde g_1,K}$ for each $(\tilde f_1,\tilde g_1)\in\Theta_K$ (see Example \ref{EX22}).

If $e\in\EE_2(K)$ is non-finite, then $X_e$ (resp.\ the \'etale locus $X_e^{\et}$ of $\psi_e$) has a canonical open affine cover associated to the open embedding $\imath_e:\mathbb A^2_{K,\s}\rightarrow X_e$ formed by connected normal surfaces that are disjoint unions $\Imm(\imath_e)\cup Y$ with $Y\in\Irr(\Gamma_{\psi_e})$ (resp.\ $Y\in\Irr(\Gamma_{\psi_e})$ and $Y\subset X_e^{\et}$) by Definition \ref{D6}(1) (resp.\ Corollary \ref{C2.8}(2)). This leads to the following formalism of sphere-like surfaces.

\begin{definition}\label{D4}
An affine normal connected surface $X$ over $K$ is called a sphere-like surface\index{sphere-like surface} if it has an open subvariety $U$ with the properties that $U\cong\mathbb A^2_K$ and $X\setminus U$ is non-empty irreducible. A sphere-like surface $X$ is called:

\medskip
{\bf (1)} a pseudo-sphere\index{sphere-like surface!pseudo-sphere} if it is not regular;

\smallskip
{\bf (2)} a quasi-sphere\index{sphere-like surface!quasi-sphere} if it is regular but its tangent bundle is non-trivial;

\smallskip
{\bf (3)} an almost sphere\index{sphere-like surface!almost sphere} if it is regular with trivial tangent bundle and $X\not\cong\mathbb S^2_K$;

\smallskip
{\bf (4)} a sphere\index{sphere-like surface!sphere} if it is isomorphic to $\mathbb S^2_K$.
\end{definition}

In this paper we develop the theory of sphere-like surfaces only as far as needed for computing invariants and for relating to prior works. 

A morphism $\Sigma:W\rightarrow X$ of varieties over $K$ is called an $\mathbb A^1$-fibration\index{$\mathbb A^1$-fibration} if for a general point $x\in X$, the fiber $\Sigma^{-1}(x)$ is an affine line over the spectrum of the residue field of $x$. Each sphere-like surface over $K$ is an $\mathbb A^1$-fibration over $\mathbb P^1_K$ and a disjoint union of an open subvariety isomorphic to $\mathbb A^2_K$ and a closed subvariety isomorphic to $\mathbb A^1_K$ by Theorem \ref{T10}(1) and (5); if $\chr(K)=0$ then the first property and, in the case when $X$ is regular, the second property were known before (e.g., see \cite{Mi3}, Sect.\ 2 and Lem.\ 1) but the proofs we include are different. 

If $e\in\EE_2(K)$, then the canonical open affine cover of $X_e$ associated to $\imath_e$ consists precisely of $|\Sing(X_e)|$ pseudo-spheres, of $\rho(e)-\rho_{\et}(e)-|\Sing(X_e)|$ quasi-spheres, and of $\rho_{\et}(e)$ almost spheres or spheres (see Corollary \ref{C2.9}); in particular, $\psi_e$ is \'etale iff the canonical open affine cover of $X_e$ associated to $\imath_e$ consists only of almost spheres and spheres (cf.\ Corollary \ref{C2.8}(2)). 

For examples of pseudo-spheres see Example \ref{EX10} and Proposition \ref{PR4.1}. For examples of quasi-spheres see each $\Imm(\imath_e)\cup Y$ with $e\in\EE_2(k)$ such that $X_e$ regular and non-\'etale and $Y\in\Irr(X_e\setminus X_e^{\et})$. If $m\ge 1$, then $\mathbb S^2_{t-1,t^m,K}$ is an almost sphere (see Example \ref{EX20-}). 

The geometries of $\SL_2$, $\mathbb S^2$, and $\mathbb A^2$ are interrelated. For instance, $\mathbb S^2$ is a quotient of $\SL_2$ and a line bundle over $\mathbb P^1$ and we have an isomorphism $\mathbb A^2\rightarrow\mathbb S^2\setminus \mathbb D^1$ induced by $\imath_{t-1,1,\mathbb Z}:\mathbb A^2\rightarrow \mathbb S^2$. The interrelations are explained and used in Section \ref{S25}.
 
Sections \ref{S3} to \ref{S10} contain different set-theoretical, algebraic, complete intersection and geometric preliminaries on iterates of functions, inner, left, and right algebraic degrees and cosets of endomorphisms $e\in\End_n(K)$, the submonoid of $\EE_2(k)$ formed by $p$-morphisms introduced in \cite{Lang-J}, Def.\ 2.2, normalizations of finite flat extensions of $k[[x_1,x_2]]$, complete intersection embedding dimensions (see Definition \ref{D5}(3)), and rational singularities. Though $\EE_1(k)=\BE_1(k)$, the set $\EE_2(k)\setminus\BE_2(k)$ is very large as, for instance, it contains all endomorphisms constructed in \cite{No}, Ex.\ (6.10) with $(a,b)\in (k^{\ast})^2$ (see also \cite{Lang-J}, Ex.\ 4.1) and generalized in \cite{Lang-J}, Props.\ 4.5 and 4.8 and Rmk.\ 4.9 based on Proposition \ref{PR1+}(3).

Section \ref{S11} presents properties of the sphere-like surfaces and the open embeddings $X_e^{\et}\subset X_e$. 

Section \ref{S12} proves Theorem \ref{T1}(1) (see Example \ref{EX7}). 

Section \ref{S13} proves Theorem \ref{T1}(2) (see Example \ref{EX8}). 

Section \ref{S14} proves Theorem \ref{T1}(3) (see Example \ref{EX11}) and Corollary \ref{C2} (see Subsection \ref{S14.1}). 

Section \ref{S15} constructs finite endomorphisms and `good' extensions of endomorphisms of affine spaces; they are used extensively in the subsequent sections in order to compute geometric degrees, extend properties in dimension $n$ to all dimensions greater than $n$, and study moduli spaces of endomorphisms. 

Properties of iterates are gathered in Section \ref{S16}. To describe them, let $e:X\rightarrow X$ be an endomorphism of a variety over $K$. Theorem \ref{T4}(1) proves that if $e$ is quasi-finite then $\iota(e)\in\mathbb N$, again by {\it reduction to positive characteristic}. Previously, this was known only when $K=\mathbb C$ (hence when $\chr(K)=0$), $X$ is affine, and $e$ is universally open with $\Gamma_e$ finite (see \cite{PNCG}, Thm.\ 1).\footnote{Recall that if $K=\mathbb C$ and $e^{\an}:X^{\an}\rightarrow X^{\an}$ is the holomorphic endomorphism associated to $e$, then $e$ is universally open iff $e^{\an}$ is open by \cite{MB}, Thm.\ 1.1.} If $\Gamma_e$ is a finite set and $\dim(X)\le 2$ (resp.\ and $e$ is defined over a finite field), then $\iota(e)\in\mathbb N$ by Proposition \ref{PR6} (resp.\ Corollary \ref{C3}). If $n\ge 3$ and $K$ is not an algebraic closure of a finite field, then for each $l\in\mathbb N^{\ast}$ there exists $d\in\End_n(K)$ such that $\varphi(d)=l$ and $\iota(d)=\infty$ by Proposition \ref{PR5}; if $n=3$, then $\deg(d)$ can be an arbitrary integer greater or equal to $3$. Corollary \ref{C3} bounds $\iota(e)$ when $\Gamma_e$ is a finite set, $\chr(K)=p$, and $e$ is defined over a finite field.

Section \ref{S17} proves Theorem \ref{T6} and for a given $r\in\mathbb N^{\ast}$ provides examples (see Corollary \ref{C5}) of endomorphisms $e$ in $\QF_2(K)$ or $\EE_2(k)$ with $\varphi(e)=r$ and $\iota(e)$ arbitrarily large.

Section \ref{S18} proves Corollary \ref{C1} based on a general lemma (see Lemma \ref{L13}) that gives large algebraic and geometric degrees and includes Example \ref{EX18} that illustrates how one can work around to control these degrees. 

Section \ref{S19} proves Theorem \ref{T2}. 

Section \ref{S20} provides examples of finite \'etale covers of $\mathbb A^2_k$ and applies them to show that for each $r\in\mathbb N$ and $m\in 2\mathbb N^{\ast}$ there exists $e\in\EE_2(k)$ such that $\rho(e)=\rho_{\et}(e)=r$ and $\deg(e)=pm$ (see Corollary \ref{C8}(1)). Example \ref{EX19} shows that the Nollet--Xavier Conjecture over $\mathbb C$ (see \cite{Je}, Sect.\ 1 and \cite{NX}, Quest.\ 6) becomes false over $k$ even for $n=2$: for each $l\in\mathbb N^{\ast}$, there exists $e\in\EE_2(k)$ such that the non-proper locus of $e$ is isomorphic to $l$ copies of $\mathbb A^1_k$ (so for $l=1$ it is smooth and integral). 

Section \ref{S21} shows that the differential equation of \cite{Lang-J}, Prop.\ 4.8 generalizes to all $n\ge 2$ producing many examples of \'etale endomorphisms in $\EE_n(k)$ whose invariants are studied in order to get applications and results that are in contrast to the analogous results in characteristic $0$ (see Remarks \ref{R17}(1) and \ref{R19} that pertain to prior works \cite{Gw}, \cite{YD}, and \cite{Zh}) or for $n>2$ are in contrast to results for $n=2$ (see Example \ref{EX23}). In particular we show that for each $n\ge 3$ there exist plenty of $e\in\EE_n(k)$ for which $\psi_e$ is non-regular (see Example \ref{EX25}).

Section \ref{S22} introduces different types of endomorphisms $e\in\EE_n(K)$ that are defined via decompositions of $e$ as composites of suitable quasi-finite endomorphisms $\mathbb A^n_K\rightarrow\mathbb A^{m+n}_K$ (such as locally closed embeddings) with projections $\mathbb A^{m+n}_K\rightarrow\mathbb A^n_K$ on the last $n$ coordinates. 

Section \ref{S23} gathers additional examples that are required in the subsequent sections or are of interest in their own. 

Section \ref{S24} proves the surjectivity of each basic endomorphism $e\in\BE_2(k)$ defined by a pair $(f_1,f_2)\in k[x_1,x_2]^p$ with $\pi(f_1,f_2)\le 2p$. 

Section \ref{S25} proves variants of Theorem \ref{T1}(2) that often compute the non-generic fibers set invariants. 

Section \ref{S26} studies moduli affine schemes of different types of \'etale endomorphisms or of automorphisms of $\mathbb A^n$ that were introduced in earlier works such as \cite{No}, Sect.\ 7 and \cite{BCW2}, Sect.\ I, Subsect.\ 1. 

For instance, we show that for each $q\in\mathbb N^\ast$ the moduli affine scheme $SEE_{n,q,K}$ over $\Spec K$ that parametrizes \'etale endomorphisms defined by $n$-tuples $(f_1,\ldots,f_n)$ in $R[x_1,\ldots,x_n]^n$ with $\pi(f_1,\ldots,f_n)\le q$ and with Jacobian determinants $1$ are connected (see Corollary \ref{C9}) and are smooth iff either $n=q=\chr(K)=2$ or $\min(n,q)=1$ (see Theorem \ref{T7}). 

Conjecture \cite{MR}, Conj.\ 3.4 (resp.\ 5.3(2)) is disproved in Example \ref{EX39} (resp.\ Remark \ref{R27}). 

Section \ref{S27} begins the study of a conjecture of Adjamagbo (see \cite{Ad}, Conj.\ 3.1 or $\textup{CJC}(n,p)$; see also Conjecture \ref{CJ2}) over subrings of our fields $k$ (i.e., over integral domains of characteristic $p$) but it suffices to prove it over the fields $k$ and it suffices to disprove in all cases over $\mathbb F_p$ itself. With our notation, this conjecture over $k$ can be reformulated as follows.

\begin{conjecture}[Adjamagbo 1995]\label{CJ3}
If $e\in\EE_n(k)$, then $\deg(e)\in\{1\}\cup p\mathbb N^\ast$. 
\end{conjecture}

Corollary \ref{C10} proves the following characteristic $p$ analogs of classical results of Campbell over $\mathbb C$ (see \cite{C}, Thm.), Razar (see \cite{Raz}, Thm.\ 2) in characteristic $0$, and Wright (see \cite{Wr1}, Thms.\ 3.2, 3.3, and 3.7) in characteristic $0$: if $e\in\EE_n(k)$ is such that either it induces a Galois field extension $k_{\t}(x_1,\ldots,x_n)\rightarrow k_{\s}(x_1,\ldots,x_n)$ of degree prime to $p$ or it is finite with $\deg(e)<p$, then $e$ is an isomorphism. 

Proposition \ref{PR16} proves that if $\chr(K)=0$ (resp.\ $\chr(K)=p)$ and $e\in\EE_n(K)$ factors as the composite of a quasi-finite morphism $\mathbb A^n_K\rightarrow\mathbb A^{n+1}_K$ which is a locally closed embedding at points of codimension $1$ outside a specific finite subset with a projection $\mathbb A^{n+1}_K\rightarrow\mathbb A^n_K$ on the last $n$ coordinates, then the Jacobian Conjecture (resp.\ Conjecture \ref{CJ3}) holds for $e$. This result generalizes \cite{Fo1}, Thm.\ 1 which assumes that $\chr(K)=0$ and that the morphism $\mathbb A^n_K\rightarrow\mathbb A^{n+1}_K$ is actually a closed embedding.

If $\chr(K)=0$, then there exists no $e\in\EE_2(K)$ with $\deg(e)\in\{2,3,4,5\}$ by \cite{\.Z}, Thm.\ 6.12. In \cite{Mos2}, Thms. 3.2 and 3.3 it is claimed that there exists no $e\in\EE_2(\mathbb C)$ with $\deg(e)$ a prime number. We do not know other general works on geometric degrees at least $3$ for $n=2$. 

Section \ref{S28} provides three general families of endomorphisms $e\in\EE_3(K)$ that disprove Conjecture \ref{CJ3} in dimension $3$ and hence in any dimension $n\ge 3$, including surjective and non-surjective families whose complements are complete intersections in $\mathbb A^3_k$ isomorphic to either $\mathbb A^1_k\setminus\{0\}$ or a finite disjoint union of curves isomorphic to $\mathbb A^1_k$ (see Theorems \ref{T6.3}, \ref{T6.4}, and \ref{T6.5}). 

In particular, for $p=2$, Corollary \ref{C8.2}(1) (resp.\ \ref{C8.2}(2)) shows that for each $j\in\mathbb N^{\ast}$ (resp.\ $j\in\mathbb N^{\ast}\setminus\{1\}$ not divisible by a prime congruent to $1$ modulo $4$) and every $n\in\mathbb N^{\ast}\setminus\{1,2\}$, there exists a surjective (resp.\ non-surjective) $e\in\EE_n(k)$ with $\deg(e)=j$.

The families are obtained via a method\index{method} in three steps, which we call the {\it linear-tame-perturbation method}, shortly the LTP method.\index{method!linear-tame-perturbation}

First we start with a simple \'etale endomorphism $d_0\in\EE_2(k)$ defined by a rule of the form $(x,y)\mapsto \bigl(f_1(x)+f_2(x)y,h_1(x,y)\bigr)$ with $(f_1,f_2)\in k[t]^2$ and $h_1\in k[t,w]$ such that $\deg(f_2)\ge 1$. 

Second we add an extra variable by considering $d_1:=d_0\times 1_{\mathbb A^1_k}\in\EE_3(k)$; it is defined by the rule $(x,y,z)\mapsto \bigl(f_1(x)+f_2(x)y,h_1(x,y),z\bigr)$ and we consider a composite $d_2=\e(g_1,g_2,g_3):=ad_1\in\EE_3(k)$, with $a\in\GA_3(k)$ a special tame automorphism defined by a simple rule $(x,y,z)\mapsto \bigl(x,y+f_3(x)z,z\bigr)$ with $f_3\in k[t]$. Here $(g_1,g_2,g_3)\in k_{\s}[x_1,x_2,x_3]$; for instance, $g_1(x_1,x_2,x_3)=f_1(x_1)+f_2(x_2)x_2$.

Third we consider polynomials $h_2(t,w)\in k[t,w]$ and perturbations 
$$d_3=\e\bigl(g_1(x_1,x_2,x_3),g_2(x_1,x_2,x_3),g_3(x_1,x_2,x_3)+h_2(x_1,x_3)^p\bigr)$$ 
of $d_2$ in the sense of Definition \ref{D4.1}(1) which are hence automatically in $\EE_3(k)$ and which involve monomials $h_2$ in the (first and third) variables fixed by $a$.

The resulting field extensions $k_{\t}(x_1,x_2,x_3)\rightarrow k_{\s}(x_1,x_2,x_3)=k_{\s}(x,y,z)$ obtained by considering the systems of equations
$$g_1(x,y,z)-x_1=g_2(x,y,z)-x_2,g_3(x,y,z)+h_2(x,z)^p-x_3=0$$
are generated by $x$ as one can easily express initially $y$ and then $z$ in terms of $x_1$, $x_2$, $x_3$, and $x$. One choses $h_2$ such that the minimal polynomial $L$ in $k_{\t}(x_1,x_2,x_3)[w]$ of $x$ over $k_{\t}(x_1,x_2,x_3)$ has degree at least $2$ congruent to $1$ or $-1$ modulo $p$ and is a quotient $\frac{L^{\textup{theroretical}}}{w^i}$, where $L^{\textup{theroretical}}$ is the `theoretical minimal polynomial' obtained from equations and of degree in $p\mathbb N^{\ast}$ and $i\in\{1\}\cup (p\mathbb N^{\ast}-1)$.\footnote{If this could be worked out with $h_2(x_1,x_2)^p$ instead of $h_1(x_1,x_3)^p$, then one would obtain counterexamples in dimension $2$ using the LTP method as well.}

Section \ref{S28} also provides one general family of endomorphisms $e\in\EE_2(K)$ that disproves Conjecture \ref{CJ3} in dimension $2$ (see Theorem \ref{T6.6} which only studies the geometric degrees and the surjectivity and hence the gap invariant $\varphi$ but none of the other invariants). 

In particular, for $p=2$, Corollary \ref{C8.3} shows that for each $j\in\mathbb N^{\ast}$ there exists $e\in\EE_2(k)$ with $\deg(e)=j$.

The method used is similar to the one above: it involves the analogs of the mentioned steps adapted to dimension $2$; we call it the {\it substitution-\'etale tame-perturbation method}, shortly the SETP method.\index{method!substitution-\'etale tame-perturbation}

While Theorem \ref{T6.6} provides many geometric degrees not divisible by $p$, it falls short from getting any substantial improvement of Equation (\ref{EQ1.1}). For instance, all the gap invariants computed in Theorem \ref{T6.6}(2.b) are congruent to $-1$ modulo $p$.

From Corollaries \ref{C8.2} and \ref{C8.3} we get directly the following consequence.

\begin{theorem}\label{T2.1}
Suppose that $p=2$. Then for each $(n,m)\in (\mathbb N^{\ast}\setminus\{1\})\times\mathbb N^{\ast}$ there exists $e\in\EE_n(k)$ with $\deg(e)=m$ and therefore there exists a smallest $q_{2,n,m}\in\mathbb N^{\ast}$ such that the set $\{d\in\EE_n(k)|\deg(d)=m,\pi(d)=q_{2,n,m}\}$ is non-empty.
\end{theorem}

Section \ref{S29} presents some general properties on extensions of vector bundles on affine open subvarieties of the affine spaces $\mathbb A^n_K$ that are required in the next section.

Section \ref{S30} studies endomorphisms $e\in\EE_n(K)$ with $\deg(e)=3$ for all $K$s, including estimates of complete intersection embedding  dimensions and partial (geometric or via equations) classifications.

For the first general and explicit applications of Theorem \ref{T10} see the link \url{https://people.math.binghamton.edu/adrian/AI.pdf} and Remark \ref{R27.9}(1) to (3)  when $n\ge 3$ and see Corollary \ref{C26} when $n=2$. In particular, for $\chr(K)=0$ (resp.\ $\chr(K)\ge 5$) Remark \ref{R27.9}(1) could be used to classify using computer programming for each given triple $(n,q,l)\in (\mathbb N^{\ast}\setminus\{1,2\})^3$ all \'etale endomorphisms $e=\e(f_1,\ldots,f_n)\in\EE_n(K)$ such that $\deg(e)=3$, $\pi(f_1,\ldots,f_n,h_0,g_1,g_2)\le q$ (with $h_0$, $g_1$, and $g_2$ as in Remark \ref{R27.9}(1)), and the coefficients of these $n+3$ polynomials are in $\llbracket-l,l\rrbracket$ (resp.\ in the finite subfield of $K$ with $\chr(K)^l$ elements).

If $\chr(K)\ge 5$, then the Conjecture on the Classification of Lines in $\mathbb A^2_K$ (see the equivalent conjectures in \cite{Moh2}, Conj.\ 1 and \cite{Da}, Conj.\ 2) implies that there exists no $e\in\EE_2(K)$ with $\deg(e)=3$ (see Corollary \ref{C26}(1) and (3)). Based on this, one could say that Conjecture \ref{CJ2.5} passed the geometric degree $3$ test, which was the original intent of Section \ref{S30}.

Section \ref{S31} complements Section \ref{S30} via examples for $p\in\{2,3\}$ and $n=2$ that emphasize special features and that point out that many hypotheses used in Section \ref{S30} are required.

Section \ref{S32} is a first appendix that introduces an ad-hoc notion of `weakly separated scheme' that is used to give alternative proofs to parts of Section \ref{S16} that involve iterations of non-separated endomorphisms of varieties.

Section \ref{S33} is a second appendix that refines and corrects Nori--Kambayashi's work on pullbacks of Lang torsors (see \cite{Kam}) based on the more recent work \cite{CGP} on split unipotent groups.

Section \ref{S34} is a third appendix that lists open problems, questions, and conjectures with a few examples, remarks, and comments. 

Theorem \ref{T1} neither is implied by nor implies Conjecture \ref{CJ2.5} in dimension $2$ (see also \cite{vdE}, Part II, Ch.\ 10, Sect.\ 10.3, Subsect.\ (10.3.16), Conj.\ $\textup{JC}(\mathbb F_p,\,n,\,p)$).%Examples \ref{EX4}

In characteristic $0$ there exist \'etale endomorphisms of surfaces that are not automorphisms (e.g., see \cite{Mi4}, Sect.\ 4, Exs.\ 1 and 3 and \cite{DP}, Thm.\ B) but it seems that all of them are surjective (and in most cases one even expects that they are finite, see \cite{Mi6}, Sect.\ 4). 

Replacing \'etaleness by the weaker quasi-finiteness, the following general example improves \cite{vdE}, Part III, Ch.\ D, Sect.\ D.2, Ex.\ D.2.3.

\begin{example}\normalfont\label{EX2}
Let $q\in\mathbb N^\ast$. Let $f(t)\in K[t]$ be a polynomial with $f(0)\neq 0$ and $\z(f)=q$. Let $e\in\End_2(K)$ be given on valued points by the rule
$$(x,y)\mapsto\bigl(xf(xy)+y,xy\bigr).$$ 
Then we have $e\in\QFE_2(K)$ with $\varphi(e)=q$, $\deg(e)=2$, and $\mathcal F(e)=\{0,1\}$. In particular, the field extension $K_{\t}(x_1,x_2)\rightarrow K_{\s}(x_1,x_2)$ induced by $e$ is isomorphic to $K_{\t}(x_1,x_2)\rightarrow K_{\t}(x_1,x_2)[z]/\bigl(z^2-\frac{x_1}{f(x_2)}z+\frac{x_2}{f(x_2)}\bigr)$. 

For $(\alpha,\beta)\in K^2$, the system of equations 
$$xf(xy)+y-\alpha=0=xy-\beta$$ in the indeterminates $x$ and $y$ has in $K^2$ no solution if $\alpha=0=f(\beta)$, one solution if $\alpha\neq 0=f(\beta)$ or $f(\beta)\neq 0=\alpha^2-4\beta f(\beta)$, and two solutions otherwise, i.e., if $f(\beta)[\alpha^2-4\beta f(\beta)]\neq 0$.\footnote{We have $\pi_{\r}(e)=\pi_{\l}(e)=2\deg(f)+1$ by Proposition \ref{PR1}(1) and (2). Hence for $\deg(f)=q$ we get that $\pi_{\r}(e)=\pi_{\l}(e)=2q+1$ and $\varphi(e)=q$ while Theorem \ref{T1}(3) (resp.\ Theorem \ref{T1}(2)) "applied to $(m,q)$ equal to $(q+1,1)$ (resp.\ to $m=1$) provides for $q\in p\mathbb N^{\ast}-1$ (resp.\ $q\in\mathbb N$) endomorphisms $e\in\EE_2(k)$ with $\pi_{\r}(e)=\pi_{\l}(e)=2q+2$ (resp.\ $\pi_{\r}(e)=\pi_{\l}(e)=3(p-1)q+2p-3$) and $\varphi(e)=q$.}
\end{example}

For $e\in\QFE_2(K)$ we have $\deg(e)=1$ iff $e\in\GA_2(K)$ (see Lemma \ref{L6}(4)). From this and the non-existence of $e\in\EE_2(K)$ with $\deg(e)=2$ and $\chr(K)\neq 2$, we get that for $\varphi(e)>0$, the equalities $\deg(e)=2$ in characteristic different from $2$ of Example \ref{EX2} and $\deg(e)=2$ in characteristic $2$ of Theorem \ref{T1}(2) or (3) are optimal in the sense that they give the smallest possible geometric degrees. We also note that if $\chr(K)\neq 2$, then there exists no $e\in \EE_2(K)\setminus\GA_2(K)$ with $\pi(e)=2$ (see \cite{Wa2}, Thm.\ 62; see also \cite{BCW2}, Thm.\ (2.4) and Cor.\ (1.4)).

\section{Set-theoretical preliminaries}\label{S3}

We set the notation and prove basic properties of iterates of endomorphisms of sets. 

For a set $\mathcal B$, let 
$$|\mathcal B|\in\mathbb N\cup\{\infty\}$$ 
be its cardinality if it is finite and be $\infty$ if it is infinite.

\begin{notation}\normalfont\label{NOT1}
Let $e:\mathcal B\rightarrow\mathcal B$ be an endomorphism of a set $\mathcal B$.

\medskip
{\bf (1)} Recall that $\Gamma_e:=\mathcal B\setminus\Imm(e)$ is the complement of $e$ and $\varphi(e)=|\Gamma_e|$ is the gap invariant of $e$. See Section \ref{S2} for $\iota(e)\in\mathbb N\cup\{\infty\}$. 

\smallskip
{\bf (2)} For $l\in\mathbb N$, let $e^l$ be the $l$-th iterate of $e$ (so $e^0:=1_{\mathcal B}$ and $e^{l+1}:=e^l\circ e$). Let
$$\quad\quad\quad\quad\quad\quad\quad\quad\mathcal E_l(e):=\Imm(e)\setminus\Imm(e^{l+1})=\{P\in\Imm(e)|e^{-1}(P)\subset\Gamma_{e^l}\}\subset\mathcal B\,$$
and
$$\quad\mathcal E_l^1(e):=\{P\in\mathcal E_l(e)||e^{-1}(P)|=1\}\subset\mathcal B.\quad\quad$$ 

{\bf (3)} Let 
$$\quad\quad\quad\quad\quad\quad\varsigma_e:=\min(|e^{-1}(P)||P\in\Imm(e)\in\mathbb N^{\ast}\cup\{\infty\}\quad\quad$$ 
$$\Gamma_e^+:=\cup_{i\in\mathbb N} e^i(\Gamma_e)\subset\mathcal B,\quad\quad\quad\;\;\,\quad\quad$$
$$\mathcal E^1(e):=\cup_{l\ge 0} \mathcal E_l^1(e)\subset\mathcal B,\quad\quad\quad\quad\;\;\;\quad\quad$$
and
$$\varrho_e:=|\mathcal E^1(e)|\in\mathbb N\cup\{\infty\}.\quad\quad\quad\quad$$ 

{\bf (4)} If $\varrho_e\in\mathbb N$, let $\xi_e\in\mathbb N$ be the smallest such that $\mathcal E^1(e)=\mathcal E_{\xi_e}^1(e)$. 
\end{notation}

\begin{lemma}\label{L1}
Assume $e:\mathcal B\rightarrow\mathcal B$ is such that $\Gamma_e$ is finite. Then the following properties hold.

\medskip
{\bf (1)} The sequence $\bigl(\varphi(e^{l+1})-\varphi(e^l)\bigr)_{l\ge 0}$ is a non-increasing sequence in $\mathbb N$. In particular, $\varphi(e^l)\in\mathbb N$ for all $l\in\mathbb N$.

\smallskip
{\bf (2)} For each $l\in\mathbb N$ we have an inequality $\varphi(e^l)\le l\varphi(e)$.

\smallskip
{\bf (3)} For each $l\in\mathbb N$ we have an identity $|\mathcal E_l(e)|=\varphi(e^{l+1})-\varphi(e)$.

\smallskip
{\bf (4)} For each $l\in\mathbb N$ we have inclusions $\Imm(e^{l})\setminus\Imm(e^{l+1})\subset e^l(\Gamma_e)$. In particular, $\mathcal E^1_l(e)\subset \mathcal E_l(e)\subset\cup_{i=1}^l e^i(\Gamma_e)\subset\Gamma_e^+$.

\smallskip
{\bf (5)} For all pairs $(l,n)\in\mathbb N^2$ with $l\le n$ we have inclusions
$$\Gamma_{e^{n+1}}\subset \Gamma_{e^l}\cup [\cup_{i=l}^{n} e^i(\Gamma_e)]\subset\cup_{i=0}^n e^i(\Gamma_e).$$

\smallskip
{\bf (6)} If $\Gamma_e^+$ is a finite set, then $\iota(e)\le |\Gamma_e^+\setminus\Gamma_e|$.

\smallskip
{\bf (7)} Assume $\varsigma_e\ge 2$. Then for each $l\in\mathbb N$ we have inequalities
$$\varphi(e^{l+1})-\varphi(e)\le\big\lfloor\frac{\varphi(e^l)}{\varsigma_e}\big\rfloor\le \frac{\varphi(e^l)}{\varsigma_e}.$$ 
Moreover, $\varphi(e^{l+1})\le\lfloor\frac{\varsigma_e\varphi(e)}{\varsigma_e-1}\rfloor$ and if $\varphi(e)>0$ we have $\varphi(e^{l+1})<\frac{\varsigma_e\varphi(e)}{\varsigma_e-1}$. In particular, the non-decreasing sequence $\bigl(\varphi(e^i)\bigr)_{i\ge 0}$ has at most $\lfloor\frac{\varsigma_e\varphi(e)}{\varsigma_e-1}\rfloor$ distinct terms and we have $\iota(e)\le \lfloor\frac{\varsigma_e\varphi(e)}{\varsigma_e-1}\rfloor$.

\smallskip
{\bf (8)} If $\varsigma_e=1$, then for each $l\in\mathbb N$ we have inequalities 
$$|\mathcal E_l^1(e)|\ge 2\varphi(e^{l+1})-\varphi(e^l)-2\varphi(e)\ge\varphi(e^l)-2\varphi(e).$$

{\bf (9)} Assume $\iota(e)=\infty$. Then $\varsigma_{e^s}=1$ for each $s\in\mathbb N^{\ast}$.

\smallskip
{\bf (10)} We have $\iota(e)=\infty$ iff the subset $\mathcal E^1(e)$ of $\Gamma_e^+$ is infinite.

\smallskip
{\bf (11)} We have an identity $\mathcal E^1(e)=\{P\in\Gamma_e^+||e^{-1}(P)|=1|\}\setminus\cap_{l\in\mathbb N} \Imm(e^l)$.

\smallskip
{\bf (12)} Assume $e$ is not surjective, $\iota(e)\in\mathbb N$ and $\varsigma_e=1$. Then $\varphi(e^l)\le 2\varphi(e)+\varrho_e-1$ for each $l\in\mathbb N$. Moreover, $\xi_e+1\le\iota(e)\le \varphi(e)+\varrho_e\le \varphi(e)(\xi_e+1)$.\end{lemma}

\begin{proof}
As $e$ maps $\Imm(e^n)$ onto $\Imm(e^{n+1})$ for $n\in\mathbb N$, it also maps $\Imm(e^l)\setminus\Imm(e^{l+1})$ onto $\Imm(e^{l+1})\setminus\Imm(e^{l+2})$. From this part (1) follows. 

As $\varphi(e^0)=0$, we have
$$\varphi(e^l)=\sum_{i=1}^l \bigl(\varphi(e^{i})-\varphi(e^{i-1})\bigr)\le \sum_{i=1}^l \bigl(\varphi(e)-\varphi(e^0)\bigr)=l\varphi(e),$$
where the inequality follows from part (1); so part (2) holds. 

Part (3) is clear.

\phantomsection{To prove part (4), let $Q\in \Imm(e^{l})\setminus\Imm(e^{l+1})$. Let $P\in\mathcal B$ be such that $Q=e^l(P)$. If $P\notin\Gamma_e$, then there exists $P^{\prime}\in\mathcal B$ such that $P=e(P^{\prime})$; so $Q=e^{l+1}(P^{\prime})\in\Imm(e^{l+1})$, a contradiction. Thus $P\in\Gamma_e$, hence $Q\in e^l(\Gamma_e)$, and part (4) holds.}\label{PH23d} 

Part (5) is a standard application of part (4). 

To prove part (6), we can assume that $\Gamma_e\neq\emptyset$; hence $\Gamma_e\neq\Gamma_e^+$. For $n\in\mathbb N^{\ast}$ we have inclusions $e^n(\Gamma_e)\subset e(\Gamma_e^+)\subset\Gamma_e^+\setminus\Gamma_e$; hence for $s\in\mathbb N$ we also have $e^{n+s}(\Gamma_e)\subset e^s(\Gamma_e^+\setminus\Gamma_e)$. Let $m\in \llbracket1,|\Gamma_e^+\setminus\Gamma_e|\rrbracket$ be the smallest integer such that $e^{m-1}(\Gamma_e^+\setminus\Gamma_e)=e^{m-1+l}(\Gamma_e^+\setminus\Gamma_e)$ for all $l\in\mathbb N$. As for $(s,l)\in\mathbb N^2$ we have $e^{m+s}(\Gamma_e)\subset e^{m-1}(\Gamma_e^+\setminus\Gamma_e)=e^{m+l}(\Gamma_e^+\setminus\Gamma_e)$, we get that $\Gamma_{e^{m+l}}\cap e^{m+s}(\Gamma_e)=\emptyset$. Thus, by taking $(l,n)$ in part (5) to be $(m,m+l-1)$ with $l\in\mathbb N^{\ast}$, it follows that we have an inclusion $\Gamma_{e^{m+l}}\subset\Gamma_{e^{m}}$. As $\Imm(e^{m+l})\subset\Imm(e^m)$ holds in general, we incur that $\Imm(e^{m+l})=\Imm(e^m)$ for all $l\in\mathbb N^{\ast}$ and hence for all $l\in\mathbb N$, so $\iota(e)\le m\le|\Gamma_e^+\setminus\Gamma_e|$.

For part (7) we note that we have a disjoint union decomposition $\sqcup_{P\in\mathcal E_l(e)} e^{-1}(P)$ of a subset of $\Gamma_{e^l}$ into sets that have at least $\varsigma_e$ elements, so $\varsigma_e|\mathcal E_l(e)|\le\varphi(e^l)$. From this and part (3) we get that $\varphi(e^{l+1})-\varphi(e)\le\frac{\varphi(e^l)}{\varsigma_e}$. So the first statement of part (7) holds and by induction on $l\in\mathbb N$ we get that $\varphi(e^{l+1})\le\varphi(e)(\sum_{i=0}^l\varsigma_e^{-i})$, from which the second statement of part (7) follows.

For part (8), let $\mathcal E_l^{\ge 2}(e):=\mathcal E_l(e)\setminus\mathcal E_l^1(e)$. We have relations
$$\varphi(e^{l+1})-\varphi(e)=|\mathcal E_l^1(e)|+|\mathcal E_l^{\ge 2}(e)|\le |\mathcal E_l^1(e)|+2|\mathcal E_l^{\ge 2}(e)|\le\varphi(e^l)$$ by part (3) and the mentioned disjoint union decomposition. Therefore we have $|\mathcal E_l^{\ge 2}(e)|\le\varphi(e)+\varphi(e^l)-\varphi(e^{l+1})$ and thus 
$$|\mathcal E_l^1(e)|=\varphi(e^{l+1})-\varphi(e)-|\mathcal E_l^{\ge 2}(e)|\ge 2\varphi(e^{l+1})-\varphi(e^l)-2\varphi(e).$$
So part (8) holds. 

If $\iota(e)=\infty$, then $\iota(e^s)=\infty$, and hence from part (7) applied to $e^s$ we get that we must have $\varsigma_{e^s}<2$. So $\varsigma_{e^s}=1$ and thus part (9) holds.

As for each $l\in\mathbb N$ we have inequalities $\varphi(e^l)-2\varphi(e)\le|\mathcal E_l^1(e)|\le \varphi(e^{l+1})-\varphi(e)$ by parts (3) and (8) and an inclusion $\mathcal E_l^1(e)\subset\mathcal E_{l+1}^1(e)$, $\mathcal E^1(e)$ is an infinite set iff the non-decreasing sequence $\varphi(e^l)_{l\in\mathbb N}$ diverges to $\infty$ and hence iff $\iota(e)=\infty$. So part (10) holds.

The `$\subset$' of part (11) follows from definitions. For $Q\in \Gamma_e^+\setminus\cap_{l\in\mathbb N} \Imm(e^l)$ with $|e^{-1}(Q)|=1$, let $l\in\mathbb N$ be the smallest such that $Q\notin\Imm(e^{l+1})$. As $|e^{-1}(Q)|=1$, we get that $Q\in\mathcal E_l^1(e)\subset\mathcal E^1(e)$ and the `$\supset$' of part (11) also holds.

For part (12), by taking $l$ in part (8) to be $\iota(e)-1$, we get $\varphi\bigl(e^{\iota(e)}\bigr)-\varphi\bigl(e^{\iota(e)-1}\bigr)\ge 1$ and hence $\varphi\bigl(e^{\iota(e)}\bigr)\le 2\varphi(e)+|\mathcal E_{\iota(e)-1}^1(e)|-1\le 2\varphi(e)+\varrho_e-1$. As for each $l\in\mathbb N$ we have an inequality $\varphi(e^l)\le\varphi\bigl(e^{\iota(e)}\bigr)$, the first statement of part (12) holds. 

As $\varphi(e)<\varphi(e^2)<\cdots<\varphi\bigl(e^{\iota(e)}\bigr)$ we get that $\iota(e)-1\le\varphi\bigl(e^{\iota(e)}\bigr)-\varphi(e)$. As $\varphi\bigl(e^{\iota(e)}\bigr)-\varphi(e)\le \varphi(e)+|\mathcal E_{\iota(e)-1}^1(e)|-1$ by the prior paragraph, we get that $\iota(e)-1\le \varphi(e)+|\mathcal E_{\iota(e)-1}^1(e)|-1$. Thus $\iota(e)\le\varphi(e)+|\mathcal E_{\iota(e)-1}^1(e)|\le\varphi(e)+\varrho_e$. 

To prove that $\xi_e+1\le\iota(e)$ we can assume that $\xi_e>0$ and hence the set $\mathcal E_{\xi_e}^1(e)\setminus \mathcal E_{\xi_e-1}^1(e)$ makes sense and it is non-empty, which implies that the inclusion $\Imm(e^{\xi_e+1})\subset\Imm(e^{\xi_e})$ is strict. Thus $\xi_e+1\le\iota(e)$. By parts (2) and (3), we estimate $\varrho_e=|\mathcal E_{\xi_e}^1(e)|\le |\mathcal E_{\xi_e}(e)|=\varphi(e^{\xi_e+1})-\varphi(e)\le\xi_e\varphi(e)$. 

The second statement of part (12) follows from the last two paragraphs.\end{proof}

\section{Polynomial preliminaries}\label{S4}

In this section we study affine coordinates in relation to algebraic degrees and (left, right, or double) cosets of $e$. 

\begin{notation}\normalfont\label{NOT2}
{\bf (1)} If $e\in\End_n(K)$ is defined by $(g_1,\ldots,g_n)\in K_{\s}[x_1,\ldots,x_n]^n$, let 
$$e^{\#}: K_{\t}[x_1,\ldots,x_n]\rightarrow K_{\s}[x_1,\ldots,x_n]$$ 
be the $K$-algebra endomorphisms such that we have $e^{\#}(x_i)=\Frac(e^{\#})(x_i)=g_i$ for each $i\in \llbracket1,n\rrbracket$; so $e=\Spec(e^{\#})$. If moreover $e\in\D_n(K)$, then $e^{\#}$ is injective and it induces a field extension $\Frac(e^{\#}): K_{\t}(x_1,\ldots,x_n)\rightarrow K_{\s}(x_1,\ldots,x_n)$.

\smallskip
{\bf (2)} The elements of the set
$$\mathcal A_n(K):=\{\{y_1,\ldots,y_n\}\subset K[x_1,\ldots,x_n]|K[y_1,\ldots,y_n]=K[x_1,\ldots,x_n]\}$$ are called affine coordinates\index{affine coordinates} of $K[x_1,\ldots,x_n]$.
\end{notation}

The following lemma computes (left or right) algebraic degrees.

\begin{lemma}\label{L2}
Let $e\in\End_n(K)$ be defined by $(g_1,\ldots,g_n)\in K_{\s}[x_1,\ldots,x_n]^n$. Then $\pi(e)$ is the smallest possible value one gets for
$$\max\bigl(\deg(h_i(g_1(y_1,\ldots,y_n),\ldots,g_n(y_1,\ldots,y_n))|i\in \llbracket1,n\rrbracket\bigr)$$ 
when $(\{h_1,\ldots,h_n\},\{y_1,\ldots,y_n\})$ varies through all elements of $\mathcal A_n(K)^2$. Similarly,
$$\pi_{\l}(e)=\min\bigl(\max(\deg(h_i(g_1(\underline{x}),\ldots,g_n(\underline{x}))|i\in \llbracket1,n\rrbracket)|(h_1,\ldots,h_n)\in\mathcal A_n(K)\bigr),$$
where $\underline{x}:=x_1,\ldots,x_n$, and 
$$\pi_{\r}(e)=\min\bigl(\max(\deg(g_1(y_1,\ldots,y_n),\ldots,g_n(y_1,\ldots,y_n)))|(y_1,\ldots,y_n)\in\mathcal A_n(K)\bigr).$$
\end{lemma}

\begin{proof} Let $(a,b)\in\GA_n(K)^2$. For $i\in \llbracket1,n\rrbracket$ let $h_i:=a^{\#}(x_i)$ and $y_i:=b^{\#}(x_i)$. The rule $(a,b)\mapsto (\{h_1,\ldots,h_n\},\{y_1,\ldots,y_n\})$ defines a bijection $\GA_n(K)^2\rightarrow\mathcal A_n(K)^2$ and we compute
$$(aeb)^{\#}(x_i)=(b^{\#}\circ e^{\#}\circ a^{\#})(x_i)=h_i\bigl(g_1(y_1,\cdots,y_n),\ldots,g_n(y_1,\cdots,y_n)\bigr).$$ Thus the first part of the lemma follows from the definition of $\pi(e)$. 

The second part of the lemma for $\pi_{\l}(e)$ (resp.\ $\pi_{\r}(e)$) follows from its definition and the first part of the lemma applied only with $b$ (resp.\ $a$) being $1_{\mathbb A^2_K}$, i.e., with $y_i=x_i$ (resp.\ $h_i=x_i$) for each $i\in \llbracket1,n\rrbracket$. 
\end{proof}

\begin{proposition}\label{PR1}
Let $e\in\End_n(K)$ be defined by an $n$-tuple $(g_1,\ldots,g_n)$ in $K_{\s}[x_1,\ldots,x_n]^n$ such that for each $i\in \llbracket1,n\rrbracket$ we can decompose $g_i=\sum_{j=1}^{d_i} \mu_{i,j}$ as a sum of non-zero monomials with $d_i\in\mathbb N^{\ast}$. We consider the following right and left monomial conditions.\index{monomial condition}

\medskip
{\bf ($\triangleright$)} If $(i,j)\in \llbracket1,n\rrbracket\times\llbracket1,d_i-1\rrbracket$, then $\mu_{i,j}$ divides $\mu_{i,d_i}$ and $\deg(\mu_{i,j})<\deg(\mu_{i,d_i})$.\index{monomial condition!($\triangleright$)}

\smallskip
{\bf ($\triangleleft$)} The primitive monomials proportional to the $\mu_{i,d_i}$s ($i\in \llbracket1,n\rrbracket$) are multiplicatively independent.\index{monomial condition!($\triangleleft$)}

\medskip
Let $\Pi:=\max(\deg(\mu_{i,d_i})|i\in \llbracket1,n\rrbracket)$. Then the following properties hold.

\medskip
{\bf (1)} If Condition ($\triangleright$) holds, then $\pi_{\r}(e)=\Pi$.

\smallskip
{\bf (2)} If Conditions ($\triangleright$) and ($\triangleleft$) hold, then $\pi_{\l}(e)=\Pi$.

\smallskip
{\bf (3)} For $i\in \llbracket1,n\rrbracket$ let $s_i\in\mathbb N$ be the greatest such that $x_i^{s_i}|\mu_{i,d_i}$ and let $\epsilon_i\in\{0,1\}$ be such that $\epsilon_i=0$ iff $\mu_{i,d_i}\in K^{\ast}x_i^{s_i}$. If Condition ($\triangleright$) holds, then we have inequalities 
$$\max\bigl(\min(\deg(\mu_{i,d_i})|i\in \llbracket1,n\rrbracket),\max(s_1+\epsilon_1,\ldots,s_n+\epsilon_n)\bigr)\le\pi_{\i}(e)\le\Pi.$$

{\bf (4)} Assume $n=2$, Condition ($\triangleright$) holds, and (to fix the ideas) $\Pi=\deg(\mu_{2,d_2})$. If either $\mu_{2,d_2}\in K^{\ast}x_2^{\deg(\mu_{2,d_2})}$ or $x_1x_2|\mu_{1,d_1}$ and $\deg(\mu_{1,d_1})=\Pi-1$, then $\pi_{\i}(e)=\Pi$.
\end{proposition}

\begin{proof}
\phantomsection{The `$\le\Pi$' is clear for parts (1) to (3). To prove the `$\ge$' for part (1), based on Lemma \ref{L2} it suffices to show that for each $\{y_1,\ldots,y_n\}\in\mathcal A_n(K)$ and every $i\in \llbracket1,n\rrbracket$ we have $\deg\bigl(g_i(y_1,\ldots,y_n)\bigr)\ge\deg(\mu_{i,d_i})$. As Condition ($\triangleright$) holds, we have an identity $\deg\bigl(g_i(y_1,\ldots,y_n)\bigr)=\deg\bigl(\mu_{i,d_i}(y_1,\ldots,y_n)\bigr)$. As $\deg(y_j)\ge 1$ for each $j\in \llbracket1,n\rrbracket$ and $\mu_{i,d_i}$ is a monomial, $\deg\bigl(\mu_{i,d_i}(y_1,\ldots,y_n)\bigr)\ge\deg(\mu_{i,d_i})$, and hence $\deg\bigl(g_i(y_1,\ldots,y_n)\bigr)\ge\deg(\mu_{i,d_i})$. So part (1) holds.}\label{PH23a}

To prove the `$\ge$' for part (2), based on Lemma \ref{L2} it suffices to show that for each $\{h_1,\ldots,h_n\}\in\mathcal A_n(K)$ and every $i\in \llbracket1,n\rrbracket$, there exists $j\in \llbracket1,n\rrbracket$ such that for $L_j=L_j(x_1,\ldots,x_n):=h_j\bigl(g_1(x_1,\ldots,x_n),\ldots,g_n(x_1,\ldots,x_n)\bigr)$, we have an inequality $\deg(L_j)\ge\deg(\mu_{i,d_i})$. We choose $j$ such that $\frac{\partial h_j}{\partial x_i}(0)\neq 0$. Let $$L_j^{\top}=L_j^{\top}(x_1,\ldots,x_n):=h_j\bigl(\mu_{1,d_1}(x_1,\ldots,x_n),\ldots,\mu_{n,d_n}(x_1,\ldots,x_n)\bigr).$$ Condition ($\triangleright$) implies that each monomial of $L_j-L_j^{\top}$ divides a suitable monomial of $L_j^{\top}$ in such a way that the quotient has positive degree, and Condition ($\triangleleft$) implies that $\deg(L_j^{\top})$ is the maximum of all degrees 
$$\deg\bigl(h_{j,s}(\mu_{1,d_1}(x_1,\ldots,x_n),\ldots,\mu_{n,d_n}(x_1,\ldots,x_n)\bigr),$$ 
where $h_{j,s}$s are the non-zero monomials of $h_j$, attained at a unique index $s$. Thus $\deg(L_j-L_j^{\top})<\deg(L_j^{\top})$, hence $\deg(L_j)=\deg(L_j^{\top})$. As $\frac{\partial h_j}{\partial x_i}(0)\neq 0$, one of the $h_{j,s}$s is in $K^{\ast}x_i$ and hence $\deg(L_j^{\top})\ge\deg(\mu_{i,d_i})$. We conclude that $\deg\bigl(L_j(x_1,\ldots,x_n)\bigr)\ge\deg(\mu_{i,d_i})$. So part (2) holds.

For parts (3) and (4), let $a\in\GA_n(K)$ and $d:=aea^{-1}\in\End_n(K)$. For $i\in \llbracket1,n\rrbracket$ let $y_i:=(a^{-1})^{\#}(x_i)$. As $d^{\#}(y_i)=[(a^{-1})^{\#}\circ e^{\#}\circ a^{\#}](y_i)=g_i(y_1,\ldots,y_n)$, from Condition ($\triangleright$) it follows that 
\begin{equation}\label{EQ1+}
\deg\bigl(d^{\#}(y_i)\bigr)=\deg\bigl(\mu_{i,d_i}(y_i,\ldots,y_n)\bigr).
\end{equation}
Let $q\in \llbracket1,n\rrbracket$ be such that $\deg(y_q)=\min(\deg(y_i)|i\in \llbracket1,n\rrbracket)$. For $i\in \llbracket1,n\rrbracket$ let $f_i:=d^{\#}(x_i)$. Let $l:=\pi(f_1,\ldots,f_n)$. As $\deg\bigl(\mu_{q,d_q}(y_1,\ldots,y_n)\bigr)\ge\deg(\mu_{q,d_q})\deg(y_q)$ and $\deg\bigl(d^{\#}(y_q)\bigr)\le l\deg(y_q)$, from Equation (\ref{EQ1+}) applied to $i=q$ it follows that $l\deg(y_q)\ge\deg(\mu_{q,d_q})\deg(y_q)$, hence $l\ge\deg\bigl(\mu_{q,d_q})\ge\min(\deg(\mu_{i,d_i})|i\in \llbracket1,n\rrbracket\bigr)$. Thus $\min\bigl(\deg(\mu_{i,d_i})|i\in \llbracket1,n\rrbracket\bigr)\le\pi_{\i}(e)$. For $i\in \llbracket1,n\rrbracket$, as $\deg\bigl(d^{\#}(y_i)\bigr)\le l\deg(y_i)$ and $s_i\deg(y_i)+\epsilon_i\le \deg\bigl(\mu_{i,d_i}(y_i,\ldots,y_n)\bigr)$, we have $s_i+\epsilon_i=s_i+\lceil\frac{\epsilon_i}{\deg(y_i)}\rceil\le l$ by Equation (\ref{EQ1+}). So part (3) holds. 

For part (4), based on part (3) it suffices to show that $l\ge \deg(\mu_{2,d_2})$. To check this we can assume that $q=1$ and $\deg(y_1)<\deg(y_2)$. If $\mu_{2,d_2}\in K^{\ast}x_2^{\deg(\mu_{2,d_2})}$, then with the notation of the last paragraph we have $\epsilon_2=0$ and $s_2=\deg(\mu_{2,d_2})$. From this and the inequality $l\ge s_2+\epsilon_2$ of part (3) we get that $l\ge \deg(\mu_{2,d_2})$. If $x_1x_2$ divides $\mu_{1,d_1}$ and $\deg(\mu_{1,d_1})=\Pi-1$, then $\deg\bigl(\mu_{1,d_1}(y_1,y_2)\bigr)>\deg(\mu_{1,d_1})\deg(y_1)$ and thus $l>\deg(\mu_{1,d_1})$. Hence $l\ge\deg(\mu_{1,d_1})+1=\Pi$. So part (4) holds.\end{proof}

\begin{example}\normalfont\label{EX3}
Let $(i,j,l)\in (\mathbb N^{\ast})^3$ with $l\ge 2$. Let $\mu:=x_1^ix_2^j$ and
$$e:=\e(x_1+\mu^l,x_2+\mu)\in\End_2(K).$$ 
If $\chr(K)=p$ and $\{il,j\}\subset p\mathbb N^{\ast}$ (resp.\ $\{i,j\}\subset p\mathbb N^{\ast}$), then $e\in\EE_2(K)$ (resp.\ $e\in\BE_2(K)$). We have identities $\pi_{\r}(e)=l\deg(\mu)=l(i+j)$ by Proposition \ref{PR1}(1). We have $\pi_{\l}(e)\le (l-1)\deg(\mu)+1<\pi_{\r}(e)$ as $\deg(\mu)=i+j\ge 2$ and for $a:=\e(x_1-x_2^l,x_2)\in\GA_2(K)$, we have $ae=\e(x_1-\sum_{m=1}^l \binom{l}{m}x_2^m\mu^{l-m},x_2+\mu)$. If $\chr(K)=p$ and $l=p$, then we have $\pi_{\l}(e)\le\pi(ae)=\max(p,i+j)$. 

The primitive monomials $x_2\mu^{l-1}$ and $\mu$ are multiplicatively independent and $x_2^m\mu^{l-m}$ divides $x_2\mu^{l-1}$ for each $m\in\llbracket2,l\rrbracket$. Thus, if $p\nmid l$, then Conditions ($\triangleright$) and ($\triangleleft$) of Proposition \ref{PR1} hold for $ae$; hence $\pi_{\l}(e)=\pi_{\l}(ae)=(l-1)\deg(\mu)+1<\pi_{\r}(e)$ by Proposition \ref{PR1}(2). The last two sentences imply that Proposition \ref{PR1}(2) does not hold without its Condition ($\triangleleft$). We have $\max(i+j,il+1)\le\pi_{\i}(e)$ by Proposition \ref{PR1}(3). If $\chr(K)=p$, $\{i,j\}\subset p\mathbb N^{\ast}$ and $p\nmid l$, then $e\in\BE_2(K)$ and $p\nmid \pi_{\l}(e)$.

Assume $l=2$. We have $ae=\e(x_1-x_2^2-2x_2\mu,x_2+\mu)$, $a^{-1}=\e(x_1+x_2^2,x_2)$, and $aea^{-1}=\e\bigl(x_1-x_2^2-2x_2(x_1+x_2^2)^ix_2^j,x_2+(x_1+x_2^2)^ix_2^j)\bigr)$. Therefore $\pi_{\i}(e)\le 2i+j$ if $\chr(K)=2$ and $\pi_{\i}(e)\le 2i+j+1$ if $\chr(K)\neq 2$. If $j=1$, then for $\chr(K)\neq 2$ we get the following relations 
$$\pi_{\l}(e)=\deg(\mu)+1=i+2\le 2i+1\le\pi_{\i}(e)\le 2i+2=\pi_{\r}(e)$$ 
and for $\chr(K)=2$ we get the following relations
$$\pi_{\l}(e)\le\deg(\mu)=i+1<2i+1=\pi_{\i}(e)< 2i+2=\pi_{\r}(e).$$ 
 Similarly, if $2\le j<i$, then for $\chr(K)\neq 2$ we get the following relations
$$\pi_{\l}(e)=i+j+1<2i+1\le \pi_{\i}(e)\le 2i+j+1< 2i+2j=\pi_{\r}(e)$$
and for $\chr(K)=2$ we get the following relations
$$\pi_{\l}(e)\le i+j<2i+1\le \pi_{\i}(e)\le 2i+j< 2i+2j=\pi_{\r}(e).$$\end{example}

\section{Formal power series preliminaries}\label{S5}

We present regularity criteria for local complete normal spectra that are finite flat over $\Spec (k[[x,y]])$ and that are obtained by taking normalizations. 

\begin{lemma}\label{L3} Let $q\in\mathbb N^{\ast}$, $R:=k[[x,y]]$ and $\alpha\in R^{\ast}$. For $m\in\mathbb N^{\ast}$ we consider the $R$-algebra $R_m:=k[[x^{\frac{1}{m}},y]]$ (so $R_1=R$) and for each $i\in\llbracket0,q-1\rrbracket$ we consider a triple $(l_i,m_i,u_i)\in\mathbb N^2\times (\{0\}\cup R^{\ast})$. We assume that $u_0\in R^{\ast}$. Let 
$$h(z):=z^{pq}+\alpha z^{pq-1}+\sum_{i=0}^{q-1} u_ix^{l_i}y^{m_i}z^{pi}\in R[z]\;\;\;\textup{and}\;\;\;A:=R[z]/\bigl(h(z)\bigr).$$ 
Then the following properties hold for $A$ and its normalization $A^{\n}$.

\medskip
{\bf (1)} The $R$-algebra $A$ is generically \'etale.

\smallskip
{\bf (2)} If $p=2$ and $q=1$, then the $R$-algebra $A$ is \'etale (so $A=A^{\n}\cong R^2$).

\smallskip
{\bf (3)} If $pq>2$, then the non-\'etale locus of $\, \Spec A\rightarrow\Spec R$ is the zero locus $z=x^{l_0}y^{m_0}=0$ in $\Spec A$.

\smallskip
{\bf (4)} Suppose that $pq-1$ divides $m_0$ and that for each $i\in\llbracket0,q-1\rrbracket$ with $u_i\neq 0$ the inequalities $m_i\ge \frac{pq-1-pi}{pq-1}m_0$ and $l_i\ge\frac{pq-1-pi}{pq-1}l_0$ hold. Then the $R$-algebra $A^{\n}$ is regular. If moreover $pq-1$ and $l_0$ are relatively prime and $l_i>\frac{pq-1-pi}{pq-1}l_0$ for each $i\in\llbracket1,q-1\rrbracket$ with $u_i\neq 0$, then the $R$-algebra $A^{\n}$ is isomorphic to $R\times R_{pq-1}$.

\smallskip
{\bf (5)} If $(l_0,m_0)\in\{(0,1),(1,0)\}$, then $A$ is regular.

\smallskip
{\bf (6)} If $(l_0,m_0)=(1,1)$, then $A=A^{\n}$ and the singular locus of $\Spec A$ is only one point defined by $x=y=z=0$. 
\end{lemma}

\begin{proof}
Parts (1) to (3) follow directly from the fact that $\frac{\partial h}{\partial z}=-\alpha z^{pq-2}$. 

By substituting $z:=\alpha z^\prime$ in $\alpha^{-pq}h(z)$ we can assume that $\alpha=1$. 

For part (4), let $d_0$ be the greatest common divisor of $pq-1$ and $l_0$. Writing $m_0=s_0(pq-1)$, $pq-1=n_0d_0$, and $l_0=d_0r_0$ with $(s_0,n_0,r_0)\in\mathbb N^3$ and $n_0$ relatively prime to $r_0$, for $t:=u_0^{\frac{1}{pq-1}}x^{\frac{r_0}{n_0}}y^{s_0}\in R_{n_0}$ we have $t^{pq-1}=u_0x^{l_0}y^{m_0}$. 

Substituting $w:=tz^{-1}$ in $tz^{-pq}h(z)$ we get that the normalization $(A\otimes_R R_{n_0})^{\n}$ of the $R_{n_0}$-algebra $A\otimes_R R_{n_0}$ is isomorphic as an $R_{n_0}$-algebra to the finite \'etale $R_{n_0}$-algebra $R_{n_0}[w]/\bigl(g(w)\bigr)$ with
$$g(w):=w^{pq}+\Bigl(\sum_{i=1}^{q-1} u_ix^{l_i}y^{m_i}t^{pi-pq+1}w^{p(q-i)}\Bigr)+w+t\in R_{n_0}[w].$$ 
This makes sense as $u_ix^{l_i}y^{m_i}t^{pi-pq+1}\in R_{n_0}$ for $i\in\llbracket0,q-i\rrbracket$; more precisely, we have $u_ix^{l_i}y^{m_i}t^{pi-pq+1}\in\{0\}\cup R^{\ast}x^{r_i}y^{s_i}$, where for $i\in\llbracket0,q-i\rrbracket$ with $u_i\neq 0$ we have $r_i:=l_i+\frac{pil_0}{pq-1}-l_0\in\frac{1}{n_0}\mathbb N$ and $s_i:=m_i+pis_0-m_0\in\mathbb N$ by hypotheses. Therefore $(A\otimes_R R_{n_0})^{\n}\cong R_{n_0}^{pq}$. From Galois theory we get that the only intermediate fields between $\Frac(R)$ and $\Frac(R_{n_0})$ are $\Frac(R_n)$ with $n$ a divisor of $n_0$. Hence each local direct factor of $A^{\n}$ is isomorphic to such an $R_n$ and therefore $A^{\n}$ is regular. 

We are left to prove that if $d_0=1$, so $n_0=pq-1$, and $r_i>0$ if $u_i\neq 0$ and $i>0$, then $A^{\n}\cong R\times R_{pq-1}$. Based on part (2) we can assume $pq>2$ and it suffices to show that there exist $R$-algebra homomorphisms $\sigma_1:A\rightarrow R$ and $\sigma_2:A\rightarrow R_{pq-1}$ which extend to surjective $R$-algebra homomorphisms $\sigma_1^{\n}:A^{\n}\rightarrow R$ and $\sigma_2^{\n}:A^{\n}\rightarrow R_{pq-1}$.

There exists a unique $R$-algebra retraction $\sigma_1:A\rightarrow R$ that maps $z+\bigl(h(z)\bigr)$ to $-1+(x,y)$ by Hensel's Lemma (e.g., see \cite{Ei}, Part I, Thm.\ 7.3); as $\sigma_1$ is surjective, so is $\sigma_1^{\n}$. Hensel's Lemma also gives the existence of an element $v\in R_{pq-1}^{\ast}$ with 
$tv^{pq}+v^{pq-1}+(\sum_{i=1}^{q-1} u_ix^{r_i}y^{s_i}v^{pi})+1=0$ (recall $r_i>0$ if $u_i\neq 0$ and $i>0$). 

By multiplying first the last equation by $t^{pq-1}=u_0x_0^{l_0}y_0^{m_0}$ and then substituting $z_0:=tv\in R_{pq-1}$, we get $z_0^{pq}+z_0^{pq-1}+\sum_{i=0}^{q-1} u_ix^{l_i}y^{m_i}z_0^{pi}=0$. Hence there exists an $R$-algebra homomorphism $\sigma_2:A\rightarrow R_{pq-1}$ that maps $z+\bigl(h(z)\bigr)$ to $z_0$. As $l_0$ is relatively prime to $pq-1$, for a divisor $n$ of $n_0=pq-1$ we have $z_0\in\Frac(R_n)$ iff $n=pq-1$. Thus, as $z_0\in\Imm(\sigma_2)$, $\sigma_2^{\n}\otimes 1_{\Frac(R)}$ is surjective and this implies that $\sigma_2^{\n}$ is surjective.

Part (5) follows from Nagata's Jacobian Criterion (see \cite{Gro2}, Prop.\ 22.7.2) as $h$ and its partial derivatives generate $R[z]$. 

For part (6), $h$ and its partial derivatives generate the ideal $(x,y,z)$ of $R[z]$, so $\{(0,0,0)\}$ is the singular locus of $\Spec A$ by loc.\ cit. Due to this and the fact that $A$ is a free module of finite rank over the regular local ring $R$ of dimension $2$, properties $(R_1)$ and $(S_2)$ hold for $A$; so $A$ is normal by Serre's Criterion (see \cite{Gro3}, Thm.\ (5.8.6) or \cite{Ma}, Ch.\ 8, Thm.\ 23.8).\end{proof}

\begin{lemma}\label{L4} Let $S:=k[[x]]$ and $R:=S[[y]]=k[[x,y]]$. For $q\in\mathbb N^{\ast}$ and $(f_1,\ldots,f_q)\in S^q$ with $f_q\neq 0$, let $h(z):=z+\sum_{i=1}^q f_i(x)z^{pi}\in S[z]$. Assume that either $q$ is a power of $p$ and $h(z)$ is additive (i.e., and $f_i=0$ if $i\in \llbracket1,q\rrbracket$ is not a power of $p$) or $q=2$. Then the normalization $A$ of $R$ in the \'etale $k((x,y))$-algebra $k((x,y))[z]/\bigl(y+h(z)\bigr)$ is regular.
\end{lemma}

\begin{proof}
If $f_q(x)\in S^{\ast}$, then $A\cong R[z]/\bigl(y+h(z)\bigr)$ is a finite \'etale $R$-algebra isomorphic to $R^{pq}$. So we can assume that $f_q(0)=0$. 

Assume $h$ is additive. Let $\prod_{i=1}^s S_i$ be the product decomposition into complete discrete valuation rings of the normalization of $S$ in $\Frac(S)[z]/\bigl(h(z)\bigr)$. For each $i\in \llbracket1,s\rrbracket$ let $z_i\in\Frac(S_i)$ be the image of $z+\bigl(h(z)\bigr)$ in $\Frac(S_i)$; it generates the separable field extension $\Frac(S)\rightarrow\Frac(S_i)$ and we have $h(z_i)=0$. Let $z_0\in R$ be the only solution of the equation $y+h(z_0)=0$ with the property that $z_0\in -y+Ry^p$; it exists by Hensel's Lemma. As $h$ is additive, $y+h(z_0+z_i)=y+h(z_0)+h(z_i)=0$ for each $i\in \llbracket1,s\rrbracket$. Hence for each $i\in \llbracket1,s\rrbracket$ there exists an $R$-algebra homomorphism $\sigma_i:A\rightarrow S_i[[y]]$ under which, at the level of rings of fractions, $z+\bigl(y+h(z)\bigr)$ is mapped to $z_0+z_i$. As $\Frac(S_i)=\Frac(R)(z_i)$ and $\sigma_i$ factors through a local integral domain direct factor of $A$, $\sigma_i$ is surjective. The resulting $R$-algebra homomorphism $(\sigma_1,\ldots,\sigma_s):A\rightarrow\prod_{i=1}^s S_i[[y]]$ is an isomorphism as it is so after tensoring with $\Frac(S)$ over $S$ (one can see this by working modulo $y$). 

Assume $q=2$. As $h(z)$ is additive if $p=2$, we can assume that $p\ge 3$. This case is similar to the previous one except that under $\sigma_i$, $z+\bigl(y+h(z)\bigr)$ maps to $z_0+z_i+t_i$ for some $t_i\in S_i[[y]]y^p$. We only show the existence of $t_i$. As we have $y+z_0+f_1(x)z_0^p+f_2(x)z_0^{2p}=0=z_i+f_1(x)z_i^p+f_2(x)z_i^{2p}$, we get 
$$y+z_0+z_i+t_i+f_1(x)(z_0+z_i+t_i)^p+f_2(x)(z_0+z_i+t_i)^{2p}$$
$$=2f_2(x)z_0^pz_i^p+t_i+[f_1(x)+2f_2(x)z_0^p+2f_2(x)z_i^p]t_i^p+f_2(x)t_i^{2p}.$$
As $z_0^p\in Ry^p$, the existence of $t_i$ follows from Hensel's Lemma once we check that $f_2(x)z_i^p\in S_i$. This is clear if $z_i=0$, hence we can assume that the index $i\in \llbracket1,s\rrbracket$ is such that $1+f_1(x)z_i^{p-1}+f_2(x)z_i^{2p-1}=0$. Thus $-1=z_i^{p-1}[f_1(x)+f_2(x)z_i^p]$, hence $z_i^{-p+1}\in S_i$ and $f_2(x)z_i^p=-z_i^{-p+1}-f_1(x)\in S_i$.\end{proof}

\begin{lemma}\label{L5} Let $R:=k[[x,y]]$. For $q\in\mathbb N^{\ast}$, $(f_0,\ldots,f_{q-1})\in R^q$, and $u\in R^{\ast}$, let $h(z):=z+(\sum_{i=0}^{q-1} f_i(x,y)z^{pi})+uxz^{pq}\in R[z]$. Then the normalization $A$ of $R$ in the \'etale $k((x,y))$-algebra $k((x,y))[z]/\bigl(h(z)\bigr)$ is regular.
\end{lemma}

\begin{proof}
If $f_0(0,0)\in k^{\ast}$, then using the substitution $w:=z^{-1}$ in $f_0^{-1}z^{-pq}h(z)$, it follows from Lemma \ref{L3}(5) that $A$ is regular. So we can assume that $f_0(0,0)=0$. 

Hensel's Lemma implies that there exists a unique $z_0\in -f_0+Rx^p+Ry^p$ such that $h(z_0)=0$. Using the substitution $u:=z_0+z$, we can assume that $f_0=0$. It suffices to show that the normalization of $R$ in the \'etale $k((x,y))$-algebra $k((x,y))[z]/\bigl(g(z)\bigr)$ is regular, where $g(z):=1+(\sum_{i=1}^{q-1} f_i(x,y)z^{pi-1})+uxz^{pq-1}$. Using the substitution $w:=z^{-1}$ in $z^{-pq+1}g(z)$, it suffices to show that for 
$$L(w):=w^{pq-1}+\bigl(\sum_{i=1}^{q-1}f_i(x,y)w^{p(q-i)}\bigr)+ux\in R[w],$$ 
\phantomsection{the $R$-algebra $B:=R[w]/\bigl(L(w)\bigr)$ is regular. As $L_x$ is a unit of $R[w]/(L_w)$, $L$ and its partial derivatives $L_x$ and $L_w$ generate $R[w]$ and hence $B$ is regular by Nagata's Jacobian Criterion (see \cite{Gro2}, Prop.\ 22.7.2).}\label{PH23b}\end{proof}

\section{Preliminaries on endomorphisms of affine spaces}\label{S6}

In this section we include properties of several classes of endomorphisms of affine spaces over $K$ which are of theoretical nature or pertain to new invariants. In particular, for $n\ge 2$ the study of basic endomorphisms $e\in\BE_n(k)$ relates closely to the study of $p$-morphisms in $\EE_n(k)$ which for $n=2$ were introduced in \cite{Lang-J}.

We recall the following well-known property.

\begin{lemma}\label{F1.0}
Let $d\in\End_n(K)$. Then $d\in\QF_n(K)$ iff $d$ is flat.
\end{lemma}

\begin{proof}
Assume $d$ is flat. Then for each $P\in \mathbb A^n_{K,\s}(K)=K^n$, the fiber $d^{-1}\bigl(d(P)\bigr)$ has dimension $0$ at $P$ by \cite{Ma}, Ch.\ 5, Thms.\ 15.1, and hence $d\in\QF_n(K)$.

If $d\in\QF_n(K)$, then $d$ is flat at each $K$-valued point of $d$ by \cite{Ma}, Ch.\ 5, Thms.\  23.1, and thus $d$ is flat.
\end{proof} 

We recall the argument for the inclusion $\Imm(\varphi_{2,K})\subset\mathbb N$ in the following general form, some part of which is only a variation of \cite{Raz}, Lem.\ 1.

\begin{lemma}\label{L6}
Let $d\in\QF_n(K)$. Then the following properties hold.

\medskip
{\bf (1)} The endomorphism $d$ is universally open.

\smallskip
{\bf (2)} The complement $\Gamma_d$ is closed and either empty or of codimension at least $2$. In particular, it is empty if $n=1$ and it is finite if $n=2$. 

\smallskip
{\bf (3)} Viewing $d^{\#}$ as an inclusion, inside $K_{\s}(x_1,\ldots,x_n)$ we have 
$$K_{\t}[x_1,\ldots,x_n]=K_{\s}[x_1,\ldots,x_n]\cap K_{\t}(x_1,\ldots,x_n).$$

\smallskip
{\bf (4)} We have $d\in\GA_n(K)$ iff $d$ is birational (i.e., $\deg(d)=1$ and, if $\chr(K)=p$, $d$ is generically \'etale).
\end{lemma}

\begin{proof} 
Part (1) is a particular case of \cite{Gro4}, Cor.\ (14.4.3).

For part (2), we note that $\Imm(d)$ is a non-empty open in $\mathbb A^n_{K,\t}$ by part (1). We show that the assumption that the reduced closed subvariety $\Gamma_d$ has an irreducible component of dimension $n-1$ leads to a contradiction. We incur the existence of a prime element $h\in K_{\t}[x_1,\ldots,x_n]$ such that the zero locus $h=0$ does not intersect $\Imm(d)$. This implies that $d^{\#}(h)\in K_{\s}[x_1,\ldots,x_n]^{\ast}=K^{\ast}=d^{\#}(K^{\ast})$. Thus, as $d^{\#}$ is injective, $h\in K^{\ast}$, a contradiction to $h$ being a prime element. So part (2) holds.

\phantomsection{For part (3), the `$\subset$' part is clear and to show the `$\supset$' part let $O_{\t}$ be a local ring of $\mathbb A^n_{k,\t}$ which is a discrete valuation ring. From part (2) we get that there exists a local ring $O_{\s}$ of $\mathbb A^n_{k,\s}$ which dominates $O_{\t}$; so $O_{\s}$ is a discrete valuation ring. We get that $K_{\s}[x_1,\ldots,x_n]\cap K_{\t}(x_1,\ldots,x_n)$ is contained in $O_{\s}\cap K_{\t}(x_1,\ldots,x_n)=O_{\t}$. Thus, as $K_{\t}[x_1,\ldots,x_n]$ is the intersection of all such local rings $O_{\t}$ of $\mathbb A^n_{k,\t}$ (see \cite{Ma}, Ch.\ 4, Thm.\ 11.5), we get that the `$\supset$' part also holds. Hence part (3) holds.}\label{PH24} 

For part (4), the `only if' part is clear and the `if' part follows from part (3) as $d$ birational implies $K_{\s}[x_1,\ldots,x_n]\cap K_{\t}(x_1,\ldots,x_n)=K_{\s}[x_1,\ldots,x_n]$.\end{proof} 

\begin{remark}\normalfont\label{R1}
Lemma \ref{L6}(2) implies that the rule $(d,e)\mapsto de$ defines a composite map $[\D_2(K)\setminus\QF_2]\times\QF_2(K)\rightarrow\D_2(K)\setminus\QF_2$. However, for $n\ge 3$ there exists $(d,e)\in [\D_n(K)\setminus\QF_n(K)]\times\QF_n(K)$ such that $de\in\QF_n(K)$. To exemplify this we can assume that $n=3$. Let $d$ be defined by the rule $(x,y,z)\mapsto (xy,xz+y,z)$. The only infinite fiber of $d$ is $d^{-1}(0,0,0)=\mathbb A^1_K\times_{\Spec K} \{(0,0)\}$. Let $e$ be such that its image is contained in $\mathbb A^3_{K,\t}\setminus [\mathbb A^1_K\times_{\Spec K} \{(0,0)\}]$, say $e$ is $1_{\mathbb A^1_K}$ times the translation by $(0,1)$ of an endomorphism in $QF_2(K)$ as in Example \ref{EX2} for $f(x)=x+1$ and $q=1$; concretely, $e$ is defined by the rule $(x,y,z)\mapsto (x,y^2z+y+z,yz+1)$. Clearly, $(e,de)\in\QF_3(K)^2$.
\end{remark}

\phantomsection{Let $\Der_n(K)$ denote the free $K_{\s}[x_1,\ldots,x_n]$-module of $K$-linear derivations of $K_{\s}[x_1,\ldots,x_n]$. Let $\{y_1,\ldots,y_n\}\in\mathcal A_n(K)$. For each $i\in \llbracket1,n\rrbracket$ let $D_{y_i}\in\Der_n(K)$ be such that for every $j\in \llbracket1,n\rrbracket$, $D_{y_i}(y_j)$ is $1$ if $j=i$ and is $0$ otherwise. Then $(D_{y_1},\ldots,D_{y_n})$ is a $K_{\s}[x_1,\ldots,x_n]$-basis of the $K_{\s}[x_1,\ldots,x_n]$-module $\Der_n(K)$.}\label{PH25}

The commutator bracket $[,]$ makes $\Der_n(K)$ a Lie algebra over $K_{\s}[x_1^p,\ldots,x_n^p]$ and we have $[D_{y_i},D_{y_j}]=0$ for all $i,j\in \llbracket1,n\rrbracket$. 

\begin{lemma}\label{L7}
Let $e=\e(g_1,\ldots,g_n)\in\D_n(k)$ and $d=\e(h_1,\ldots,h_n)\in\D_n(k)$. Then the following properties hold.

\medskip
{\bf (1)} We have $\GA_n(K)e=\GA_n(K)d$ iff $K_{\s}[g_1,\ldots,g_n]=K_{\s}[h_1,\ldots,h_n]$.

\smallskip
{\bf (2)} We have an identity $\GA_n(K)e\GA_n(K)=\GA_n(K)d\GA_n(K)$ iff there exists  $\{y_1,\ldots,y_n\}\in\mathbb A_n(K)$ such that the  identification $K_{\s}[x_1,\ldots,x_n]=K_{\s}[y_1,\ldots,y_n]$ restricts to an identification 
\begin{equation}\label{EQ22}
K_{\s}[g_i(x_1,\ldots,x_n)|i\in \llbracket1,n\rrbracket]=K_{\s}[h_i(y_1,\ldots,y_n)|i\in \llbracket1,n\rrbracket].
\end{equation}

\smallskip
{\bf (3)} If there exists $c\in\End_n(K)$ such that $ce$ is \'etale, then both $c$ and $e$ are \'etale.

{\bf (4)} Suppose that $(e,d)\in\EE_n(K)^2$ and $\deg(e)=\deg(d)$. Then part (2) holds with the equality in Equation (\ref{EQ22}) replaced by either one of the two possible inclusions.
\end{lemma}

\begin{proof} For $a\in\GA_n(K)$ we have $ae=d$ iff $d^{\#}=e^{\#}\circ a^{\#}$. Thus, such an $a$ exists iff $\Imm(e^{\#})=\Imm(d^{\#})$ and part (1) follows. 

For $(a,b)\in\GA_n(K)^2$ we have $aeb=d$ iff $bae=bdb^{-1}$. As we have an identity $(bdb^{-1})^{\#}=(b^{-1})^{\#}\circ d^{\#}\circ b^{\#}$, for $i\in \llbracket1,n\rrbracket$ we define $y_i:=(b^{\#})^{-1}(x_i)=(b^{-1})^{\#}(x_i)$. Thus $\{y_1,\ldots,y_n\}\in\mathcal A_n(K)$ is defined by $b^{-1}$ and determines $b^{-1}$ (or $b$) uniquely. So, as $\Imm\bigl((bdb^{-1})^{\#}\bigr)=K_{\s}[h_1(y_1,\ldots,y_n),\ldots,h_n(y_1,\ldots,y_n)]\subset K_{\s}[y_1,\ldots,y_n]$, part (2) follows from the last paragraph applied to $(a,e,d)=(ba,e,bdb^{-1})$.

For part (3), clearly $e$ is \'etale and the \'etale locus of $c$ contains $\Imm(e)$. Thus the non-\'etale locus of $c$ is either empty or its codimension in $\mathbb A^n_{K,\t}$ is at least 2 by Lemma \ref{L6}(2). Thus, being the zero locus of the Jacobian determinant, the non-\'etale locus of $c$ is empty, i.e., $c$ is \'etale and part (3) holds.

For part (4) it suffices to show that if $K_{\s}[g_1(x_1,\ldots,x_n),\ldots,g_n(x_1,\ldots,x_n)]$ contains $K_{\s}[h_1(y_1,\ldots,y_n),\ldots,h_n(y_1,\ldots,y_n)]$, then Equation (\ref{EQ22}) holds. Due to this inclusion we have a factorization $d=ce$, where the endomorphism $c:\mathbb A^n_{K,\t}\rightarrow \mathbb A^n_{K,\t}$ is birational (as $\deg(e)=\deg(d)$) and \'etale by part (3); from Lemma \ref{L6}(4) we get that $c$ is an isomorphism. Thus Equation (\ref{EQ22}) holds by part (2).\end{proof}

\begin{definition}\label{D4.1}
Let $e=\e(g_1,\ldots,g_n)\in\End_n(k)$ and $d=\e(h_1,\ldots,h_n)\in\End_n(k)$. Let $\mathfrak V_n(k)$ be the set of $k$-vector subspaces of $k_{\s}[x_1,\ldots,x_n]$ that contain $k_{\s}[x_1^p,\ldots,x_n^p]$.

\medskip
{\bf (1)} We say that $e$ is a perturbation\index{perturbation} of $d$ if $f_i-g_i\in k_{\s}[x_1^p,\ldots,x_n^p]$ for all $i\in \llbracket1,n\rrbracket$.

\smallskip
{\bf (2)} The perturbed image left invariant\index{perturbation!perturbed image left invariant} 
$$\Re=\Re_n(k):\End_n(k)\rightarrow\mathfrak V_n(k)$$ 
is defined by the rule $\Re(e):=\Imm(e^{\#})+k_{\s}[x_1^p,\ldots,x_n^p]=k_{\s}[g_1,\ldots,g_n]+k_{\s}[x_1^p,\ldots,x_n^p]$.

\smallskip
{\bf (3)} The basic number invariant\index{basic number invariant} 
$$\mathfrak B=\mathfrak B_n(k):\End_n(k)\rightarrow \llbracket0,n\rrbracket$$
is defined by the rule: $\mathfrak B(e)$ is the largest integer in $\llbracket0,n\rrbracket$ such that there exists $\{y_1,\ldots,y_n\}\in\mathcal A_n(K)$ with $\{y_1,\ldots,y_{\mathfrak B(e)}\}\subset\Re(e)$.

\smallskip
{\bf (4)} Suppose that $(d,e)\in\EE_n(k)$. We say that $e$ is a multiplicative perturbation\index{multiplicative perturbation} of $d$ if there exists $i\in\llbracket1,n\rrbracket$ and a pair $(\tilde g_i,h_i)\in k_{\s}[x_1,\ldots,x_n]^2$ such that $g_i=f_i+\tilde g_ih_i^p$, $\e(f_1,\ldots,f_{i-1},f_i+\tilde g_i,f_{i+1},\ldots,f_n)\in\EE_n(k)$, and $\tilde g_i\notin k_{\s}[x_1^p,\ldots,x_n^p]$.
\end{definition} 

\begin{remark}\normalfont\label{R1.5}
{\bf (1)} We view $\Re(e)$ as a left $k_{\t}[x_1^p,\ldots,x_n^p]$-module via $e^{\#}$. 

\smallskip
{\bf (2)} The relation $\{(e,d)\in\End_n(k)^2|e\;\textup{is a perturbation of}\; d\}$ on $\End_n(k)$ is an equivalence relation. 

\smallskip
{\bf (3)} If $e$ is a perturbation of $d$, then the Jacobian matrices of $e$ and $d$ are equal, hence $e\in\EE_n(k)$ iff $d\in\EE_n(k)$.

\smallskip
{\bf (4)} The $k$-subalgebra of $k_{\s}[x_1,\ldots,x_n]$ generated by $\Re(e)$ is $k_{\s}[x_1,\ldots,x_n]$ itself iff $e\in\EE_n(k)$ by the equivalence of $(N1)\Leftrightarrow (N2)$ of Section \ref{S2}. Thus, if $\mathfrak B(e)=n$, then $e\in\EE_n(k)$.
\end{remark}

\begin{definition}\label{D4.2}
Let $e\in\End_n(K)$ be generically \'etale. We say that $e$ is a quasi-automorphism\index{quasi-automorphism} if any one of the following two equivalent statements holds.

\medskip
{\bf (1)} There exists $a\in\GA_n(K)$ such that $D_{x_i}\bigl(\Imm(a^{\#}\circ e^{\#})\bigr)\subset \Imm(a^{\#}\circ e^{\#})$ for each $i\in \llbracket1,n\rrbracket$.

\smallskip
{\bf (2)} There exists $\{y_1,\ldots,y_n\}\in\mathcal A_n(K)$ such that $D_{y_i}\bigl(\Imm(e^{\#})\bigr)\subset \Imm(e^{\#})$ for each $i\in \llbracket1,n\rrbracket$.\end{definition}

The equivalence $(1)\Leftrightarrow (2)$ is easily checked via $y_i=a^{\#}(x_i)$ for $i\in \llbracket1,n\rrbracket$. Clearly, each automorphism in $\End_n(K)$ is a quasi-automorphism. 

\begin{lemma}\label{F1.1}
Suppose that $\chr(K)=0$. Then $e$ is an automorphism iff it is a quasi-automorphism.
\end{lemma}

\begin{proof}
The `only if' part is clear. For the `if' part it suffices to show that if $(f_1,\ldots,f_n)\in K_{\s}[x_1,\ldots,x_n]^n$ defines $e$ and $D_{x_i}(f_j)\in K_{\s}[f_1,\ldots,f_n]$ for each $(i,j)\in \llbracket1,m\rrbracket^2$, then $K_{\s}[f_1,\ldots,f_n]=K_{\s}[x_1,\ldots,x_n]$, i.e., $x_i\in K_{\s}[f_1,\ldots,f_n]$ for each $i\in \llbracket1,m\rrbracket$. To check this we can assume that $i=1$. As $e$ is generically \'etale, we have $e\in\D_n(K)$. So there exists $j\in \llbracket1,m\rrbracket$ such that $f_j\notin K_{\s}[x_2,\ldots,x_n]$. Among the primitive monomials $\prod_{i=1}^n x_i^{l_i}$ with $(l_1,\ldots,l_n)\in\mathbb N^n$ that show up in $f_j$ with a non-zero coefficient, we choose the one which is maximal with respect to the lexicographic order on $\mathbb N^n$; so $l_1\ge 1$. As $(D_n^{l_n}\circ D_{n-1}^{l_{n-1}}\circ D_2^{l_2}\circ D_1^{l_1-1})(f_j)\in K^{\ast}x_1$ by the maximality property, we have $K^{\ast}x_1\cap K_{\s}[f_1,\ldots,f_n]\neq\emptyset$. So $x_1\in K_{\s}[x_2,\ldots,x_n]$ and the `if' part holds.\end{proof}

In contrast to Lemma \ref{F1.1}, for $n\ge 2$ there exist quasi-automorphisms in $\End_n(k)$ that are finite but not \'etale, e.g., $\e(x_1x_2^p+x_1^p,x_2,\ldots,x_n)$. Moreover, next we show that for each $n\in\mathbb N^{\ast}$, every $e\in\BE_n(k)$ is a quasi-automorphism. 

\begin{proposition}\label{L8.1}
For $e\in\End_n(k)$ we consider the following statements.

\medskip
{\bf (1)} We have $e\in\BE_n(k)$.

\smallskip
{\bf (2)} We have $\mathfrak B(e)=n$. 

\smallskip
{\bf (3)} There exists $c\in\EE_n(k)$ such that $ce\in\BE_n(k)$.

\smallskip
{\bf (4)} The endomorphism $e$ is \'etale and a quasi-automorphism. 

\smallskip
{\bf (5)} The endomorphism $e$ is a quasi-automorphism. 

\smallskip
{\bf (6)} We have $\mathfrak B(e)\ge 1$

\medskip
Then the following implications $(1)\Rightarrow (2)\Leftrightarrow (3)\Rightarrow (4)\Rightarrow (5)\Rightarrow (6)$ hold.\end{proposition}

\begin{proof}
We write $e=\e(g_1,\ldots,g_n)$ with $(g_1,\ldots,g_n)\in k_{\s}[x_1,\ldots,x_n]^n$. 

Clearly, $(1)\Rightarrow (3)$ and $(4)\Rightarrow (5)$. 

To prove that $(2)\Rightarrow (3)$, let $\{y_1,\ldots,y_n\}\in\mathcal A_n(k)$ be such that we have an inclusion $\{y_1,\ldots,y_n\}\subset\Re(e)$. For $i\in \llbracket1,n\rrbracket$ let $v_i\in k_{\s}[x_1,\ldots,x_n]$ be such that $y_i-v_i^p\in\Imm(e^{\#})$. If $d:=\e(y_1-v_1^p,\ldots,y_n-v_n^p)\in\End_n(k)$, then $d\in\BE_n(k)$ and we have a factorization $d=ce$ with $c\in\End_n(K)$ such that $c^{\#}$ is defined by the inclusion $\Imm(d^{\#})\subset\Imm(e^{\#})$. As $(c,e)\in\EE_n(k)^2$ by Lemma \ref{L7}(3), $(2)\Rightarrow (3)$. 

To prove that $(3)\Rightarrow (2)$, let $\{z_1,\ldots,z_n\}\in\mathcal A_n(k)$ be such that we have an inclusion $\{z_1,\ldots,z_n\}\subset\Re(ce)$ by the implication $(1)\Rightarrow (2)$ applied to $ce$. As $\Re(ce)\subset\Re(e)$, we get that $\{z_1,\ldots,z_n\}\subset\Re(e)$ and thus $(3)\Rightarrow (2)$.

To prove that $(3)\Rightarrow (4)$, let the pair $(a,b)\in\GA_n(k)^2$ be such that we have $aceb=\e(x_1+f_1^p,\ldots,x_n+f_n^p)$ with $(f_1,\ldots,f_n)\in k_{\s}[x_1,\ldots,x_n]^n$. By replacing $(c,e)$ with $(ac,eb)$ we can assume that $a=b=1_{\mathbb A^n_k}$; so $ce=\e(x_1+f_1^p,\ldots,x_n+f_n^p)$ is in $\BE_n(k)$. So $c$ and $e$ are \'etale by Lemma \ref{L7}(3). 

As the Jacobian matrix of $\e(x_1+f_1^p,\ldots,x_n+f_n^p)$ is the identity matrix, from the chain rule for Jacobian matrices it follows that for $(\alpha_1,\ldots,\alpha_n)\in k^n$, the inverse of the Jacobian matrix $J(g_1,\ldots,g_n)(\alpha_1,\ldots,\alpha_n)\in\GL_n(k)$ depends only on $e(\alpha_1,\ldots,\alpha_n)$. Thus $J(g_1,\ldots,g_n)(\alpha_1,\ldots,\alpha_n)$ itself depends only on $e(\alpha_1,\ldots,\alpha_n)$. From this and standard Galois descent with respect to the Galois extension of $k_{\t}(x_1,\ldots,x_n)$ generated by $\Frac(e^{\#})$, it follows that $J(g_1,\ldots,g_n)$ has entries in $\Imm\bigl(\Frac(e^{\#})\bigr)$. As the entries are also in $k_{\s}[x_1,\ldots,x_n]$, from Lemma \ref{L6}(3) we get that they are in $\Imm(e^{\#})=k_{\s}[g_1,\ldots,g_n]$. Hence $(3)\Rightarrow (4)$.

To prove that $(5)\Rightarrow (6)$, it suffices to show that if $e$ is generically \'etale and such that for $A:=k_{\s}[g_1,\ldots,g_n]$ we have $D_{x_i}(A)\subset A$ for each $i\in\llbracket1,n\rrbracket$, then there exists $(\alpha_1,\ldots,\alpha_n)\in k^n\setminus\{(0,\ldots,0)\}$ such that $\sum_{i=1}^n \alpha_ix_i \in\Re(e)$. 

Let $f\in A\setminus k_{\s}$ be such that 
\begin{equation}\label{EQ39.7}
\deg(f)=\min\bigl(\deg(h)|h\in A\setminus k_{\s}\bigr). 
\end{equation} 
As $(x_1,\ldots,x_n)$ is a $p$-basis of $k_{\s}[x_1,\ldots,x_n]$, we can write uniquely 
$$f=\sum_{(i_1,\ldots,i_n)\in \llbracket0,p-1\rrbracket^n} h^p_{i_1,\ldots,i_n}\prod_{j=1}^n x_j^{i_j}$$ 
with each $h_{i_1,\ldots,i_n}\in K_{\s}[x_1,\ldots,x_n]$. Clearly, we have
\begin{equation}\label{EQ39.8}
\deg(f)=\max\bigl(\{p\deg(h_{i_1,\ldots,h_n})+\sum_{j=1}^n i_j|(i_1,\ldots,i_n)\in \llbracket0,p-1\rrbracket^n\}\bigr).
\end{equation} 
If $f=h_{0,\ldots,0}^p$, then, as $e$ is generically \'etale, we have $h_{0,\ldots,0}\in A$. The inequality $\deg(h_{0,\ldots,0})<\deg(f)$ contradicts Equation (\ref{EQ39.7}). So there exists an $n$-tuple $(j_1,\ldots,j_n)\in \llbracket0,p-1\rrbracket^n\setminus\{(0,\ldots,0)\}$ such that $h_{j_1,\ldots,j_n}\neq 0$. As $D_{x_i}(A)\subset A$ for each $i\in \llbracket1,n\rrbracket$, Equation (\ref{EQ39.7}) and Equation (\ref{EQ39.8}) applied to $f$ and its partial derivatives give that $h_{i_1,\ldots,i_n}=0$ if $\sum_{j=1}^n i_j\ge 2$. Hence by denoting $h_0:=h_{0,\ldots,0}$ and $h_j:=h_{0,\ldots,0,1,0,\ldots,0}$ with $1$ as the $j$-th entry, we have $f=h_0^p+\sum_{i=1}^n x_ih_i^p$ and there exists $i\in \llbracket1,n\rrbracket$ such that $h_i\neq 0$. If $h_i\neq 0$, then $\deg(f)\ge 1+p\deg(h_i)$ by Equation (\ref{EQ39.8}). From this, the fact that $D_{x_i}(f)=h_i^p\in A$, and Equation (\ref{EQ39.7}) it follows that $h_i\in k_{\s}$ for each $i\in \llbracket1,n\rrbracket$. So by taking $\alpha_i:=h_i^p$ for each $i\in \llbracket1,n\rrbracket$ we have $(\alpha_1,\ldots,\alpha_n)\in k^n\setminus\{(0,\ldots,0)\}$ and $\sum_{i=1}^n \alpha_ix_i =f-h_0^p\in\Re(e)$. So $(5)\Rightarrow (6)$.
\end{proof} 

To refine slightly Proposition \ref{L8.1} for $n=p=2$ we first prove a general lemma.

\begin{lemma}\label{L8.9}
The following properties hold for $e=\e(g_1,\ldots,g_n)\in\QFE_n(s)$.

\medskip
{\bf (1)} If $f\in k_{\s}[x_1,\ldots,x_n]$ is such that $f^p\in k_{\s}[g_1,\ldots,g_n]$, then $f\in k_{\s}[g_1,\ldots,g_n]$.

\smallskip
{\bf (2)} Assume $e\in\EE_n(k)$. If $(f,h)\in k_{\s}[x_1,\ldots,x_n]\times k_{\s}[g_1,\ldots,g_n]$ is such that $fh^p\in\Re(e)$, then $f\in\Re(e)$.
\end{lemma}

\begin{proof}
For part (1), as the field extension $k_{\s}(g_1,\ldots,g_n)\rightarrow k_{\s}(x_1,\ldots,x_n)$ is separable, we have $f\in k_{\s}(g_1,\ldots,g_n)$. Thus $f\in k_{\s}(g_1,\ldots,g_n)\cap k_{\s}[x_1,\ldots,x_n]$, hence $f\in k_{\s}[g_1,\ldots,g_n]$ by Lemma \ref{L6}(3). So part (1) holds.

For part (2), we consider a pair $(f_0,g)\in k_{\s}[x_1,\ldots,x_n]\times k_{\s}[g_1,\ldots,g_n]$ such that $fh^p=f_0^p+g$. There exist unique polynomials $h_{i_1,\ldots,i_n}\in k_{\s}[g_1,\ldots,g_n]$ indexed by $n$-tuples $(i_1,\ldots,i_n)\in \llbracket0,p-1\rrbracket^n$ such that $g=\sum_{(i_1,\ldots,i_n)\in \llbracket0,p-1\rrbracket^n} h_{i_1,\ldots,i_n}^p\prod_{j=1}^n g_j^{i_j}$. Let $f_1:=f_0+h_{0,\ldots,0}\in k_{\s}[x_1,\ldots,x_n]$. We get the identity
\begin{equation}\label{EQ39.9}
f=\left(\frac{f_1}{h}\right)^p+\sum_{(i_1,\ldots,i_n)\in \llbracket0,p-1\rrbracket^n\setminus\{(0,\ldots,0)\}} \left(\frac{h_{i_1,\ldots,i_n}}{h}\right)^p\prod_{j=1}^n g_j^{i_j}.
\end{equation}

From the equivalence $(N1)\Leftrightarrow (N3)$ (see Section \ref{S2}) we get a $k_{\s}(x_1^p,\ldots,x_n^p)$-basis $(\prod_{j=1}^n g_j^{i_j}|(i_1,\ldots,i_n)\in \llbracket0,p-1\rrbracket^n)$ of $k_{\s}(x_1,\ldots,x_n)$. 

From the last two sentences and $(N3)$ we get that $\frac{f_1}{h}\in k_{\s}[x_1,\ldots,x_n]$ and $\frac{h_{i_1,\ldots,i_n}}{h}\in k_{\s}[x_1,\ldots,x_n]$ for each $(i_1,\ldots,i_n)\in \llbracket0,p-1\rrbracket^n\setminus\{(0,\ldots,0)\}$. Hence $\frac{h_{i_1,\ldots,i_n}}{h}\in k_{\s}[x_1,\ldots,x_n]\cap k_{\s}(g_1,\ldots,g_n)$ and thus $\frac{h_{i_1,\ldots,i_n}}{h}\in k_{\s}[g_1,\ldots,g_n]$ by Lemma \ref{L6}(3) for each $(i_1,\ldots,i_n)\in \llbracket0,p-1\rrbracket^n\setminus\{(0,\ldots,0)\}$. From this and Equation (\ref{EQ39.9}) we get that $f\in k_{\s}[x_1,\ldots,x_n]+k_{\s}[g_1,\ldots,g_n]=\Re(e)$. So part (2) holds.
\end{proof} 

\begin{corollary}\label{C2.3}
Assume $n=p=2$. Then the statements (2) to (4) of Proposition \ref{L8.1} are equivalent.
\end{corollary}

\begin{proof} 
Based on Proposition \ref{L8.1} it suffices to show that for a quasi-automorphism $e=\e(g_1,g_2)\in\EE_2(k)$ we have $\mathfrak B(e)=2$. To ease notation we assume that for $A:=k_{s}[g_1,g_2]$ we have $D_{x_i}(A)\subset A$ for each $i\in\{1,2\}$. For $i\in\{1,2\}$ and $j\in\{0,1,2,3\}$, let $h_{ij}\in k_{s}[x_1,x_2]$ be such that $g_i=h_{i0}^2+h_{i1}^2x_1+h_{i2}^2x_2+h_{i3}^2x_1x_2$. 

We consider two cases as follows. 

{\bf Case 1: $h_{13}$ and $h_{23}$ are not both $0$.} To fix the ideas we can assume that $h_{13}\neq 0$. We have $D_{x_1}(g_1)=h_{11}^2+h_{13}^2x_2\in A$ and $D_{x_2}(g_1)=h_{12}^2+h_{13}^2x_1\in A$. Hence for $i\in\{1,2\}$ we have $h_{13}^2x_i\in A+k_{s}[x_1^2,x_2^2]=\Re(e)$ and from Lemma \ref{L8.9}(2) we get that $x_i\in\Re(e)$. Thus $\mathfrak B(e)=2$.

{\bf Case 2: $h_{13}=h_{23}=0$.} For $(i,j)\in\{1,2\}^2$ we have $D_{x_j}(g_i)=h_{ij}^2\in A$ and thus $h_{ij}\in A$ by Lemma \ref{L8.9}(1). As the determinant of the Jacobian matrix $J(g_1,g_2)=\begin{bmatrix} 
h_{11}^2 & h_{12}^2 \\
h_{21}^2 & h_{22}^2 \\ 
\end{bmatrix}$ is $\alpha\in k^{\ast}$, we have $J(g_1,g_2)^{-1}=\alpha^{-1}\begin{bmatrix} 
h_{22}^2 & h_{12}^2 \\
h_{21}^2 & h_{11}^2 \\ 
\end{bmatrix}$. From this and the matrix identity $\begin{bmatrix} 
g_1\\
g_2\\ 
\end{bmatrix}=\begin{bmatrix} 
h_{10}^2\\
h_{20}^2\\ 
\end{bmatrix}+J(g_1,g_2)\begin{bmatrix} 
x_1\\
x_2\\ 
\end{bmatrix}$, we get another matrix identity $\begin{bmatrix} 
x_1\\
x_2\\ 
\end{bmatrix}=\alpha^{-1}\begin{bmatrix} 
h_{22}^2 & h_{12}^2 \\
h_{21}^2 & h_{11}^2 \\ 
\end{bmatrix}\begin{bmatrix} 
g_1+h_{10}^2\\
g_2+h_{20}^2\\ 
\end{bmatrix}$. From this and the fact that $\{g_1,g_2,h_{11},h_{12},h_{21},h_{22}\}\subset A$ we get that $(x_1,x_2)\in\Re(e)^2$. Hence  $\mathfrak B(e)=2$.
\end{proof} 

\begin{remark}\normalfont\label{R2}
If $f(x_1,\ldots,x_{n-1})\in k_{\s}[x_1,\ldots,x_{n-1}]\setminus k_{\s}[x_1^p,\ldots,x_{n-1}^p]$, then the $n$-tuple $\bigl(x_1,\ldots,x_{n-1},x_n+f(x_1,\ldots,x_{n-1})x_n^{pm}\bigr)\in k_{\s}[x_1,\ldots,x_n]^n$ defines $e\in\EE_n(k)$ with $x_n\notin\Re(e)$. To show this, by substituting $x_i=x^{q_i}$ for $i\in \llbracket1,n-1\rrbracket$ with $q_i\in\{1\}\cup p\mathbb N^\ast$ such that $f(x^{q_1},\ldots,x^{q_{n-1}})\in k_{\s}[x]\setminus k_{\s}[x^p]$, we can assume that $n=2$. For $n=2$ we show that $\left(k_{\s}[x_2]\setminus k_{\s}[x_2^p]\right)\cap \Re(e)=\emptyset$. It suffices to show that the assumption that there exists $(g,h,F)\in (k_{\s}[x,y])^2\times \left(k_{\s}[x_2]\setminus k_{\s}[x_2^p]\right)$ such that
\begin{equation}\label{EQ40}
g\bigl(x_1,x_2+f(x_1)x_2^{pm}\bigr)=F(x_2)+h(x_1^p,x_2^p)
\end{equation}
leads to a contradiction. Applying $D_{x_1}$ to this identity, we get that 
$$g_x\bigl(x_1,x_2+f(x_1)x_2^{pm}\bigr)+g_y\bigl(x_1,x_2+f(x_1)x_2^{pm}\bigr)f'(x_1)x_2^{pm}=0.$$ 
If $g_y=0$, then also $g_x=0$; so $g\in k_{\s}[x^p,y^p]$ and the left-hand side of Equation (\ref{EQ40}) is in $k_{\s}[x_1^p,x_2^p]$ and cannot be equal to $F(x_2)+h(x_1^p,x_2^p)$. Thus $g_y\neq 0$ and
$$f'(x_1)x_2^{pm}=\frac{-g_x\bigl(x_1,x_2+f(x_1)x_2^{pm}\bigr)}{g_y\bigl(x_1,x_2+f(x_1)x_2^{pm}\bigr)}$$ 
belongs to the intersection of the field of fractions of $k_{\s}[x_1,x_2+f(x_1)x_2^{pm}]$ with $k_{\s}[x_1,x_2]$ and hence (see Lemma \ref{L6}(3)) it belongs to $k_{\s}[x_1,x_2+f(x_1)x_2^{pm}]$. The relation $f(x_1)\in k_{\s}[x_1]\setminus k_{\s}[x_1^p]$ implies that $0\le \deg\bigl(f'(x_1)\bigr)<\deg\bigl(f(x_1)\bigr)$, and it is easy to see that $f'(x_1)x_2^{pm}\notin k_{\s}[x_1,x_2+f(x_1)x_2^{pm}]$, contradiction.
\end{remark}

The following notions were either explicitly or in essence introduced in \cite{Lang-J}, Defs.\ 2.2 and 2.15 for $n=2$.

\begin{definition}\label{D4.3}
Assume $n\ge 2$. For $f\in k[x_1,\ldots,x_n]$, let $e_f\in\EE_n(k)$ be defined by the rule $(x_1,\ldots,x_n)\mapsto\bigl(x_1,\ldots,x_n+f(x_1,\ldots,x_{n-1},x_n^p)\bigr)$. Let $\pMor_n(k)$ be the submonoid of $\EE_n(k)$ generated by $\GA_n(k)\cup \{e_f\in\EE_n(k)|f\in k[x_1,\ldots,x_n]\}$. 

\medskip
{\bf (1)} The elements of $\pMor_n(k)$ are called $p$-morphisms\index{$p$-morphism} and the elements in $\GA_n(k)\{e_f\in\EE_n(k)|f\in k[x_1,\ldots,x_n]\}\GA_n(k)$ are called elementary $p$-morphisms\index{$p$-morphism!elementary}.

\smallskip
{\bf (2)} Let $f\in k[x_1,\ldots,x_n]\setminus k[x_1^p,\ldots,x_n^p]$. We decompose $f=\sum_{i=0}^{\deg(f)} f_i$, where each $f_i$ is homogeneous of degree $i$.\footnote{By convention, the homogeneous polynomials of degree $i\in\mathbb N$ contain the $0$ polynomial.} By the modulo $p$-th powers degree\index{modulo $p$-th powers degree} of $f$, to be denoted as $\deg_p(f)$, we mean $m\in \llbracket1,\deg(f)\rrbracket$ such that $f_m\notin k[x_1^p,\ldots,x_n^p]$ but $f_i\in k[x_1^p,\ldots,x_n^p]$ for each $i\in\llbracket m+1,\deg(e)\rrbracket$.

\smallskip
{\bf (3)} Assume $f\notin k[x_1^p,\ldots,x_n^p]$. By a standard writing modulo $p$-th powers\index{standard writing} of $F$, shortly a standard writing of $f$, we mean a decomposition $f=g^p+\sum_{i=1}^{\deg_p(f)} h_i$ with $g\in k[x_1,\ldots,x_n]$ and each $h_i\in k[x_1,\ldots,x_n]$ homogeneous of degree $i$. A standard writing $f=g^p+\sum_{i=1}^{\deg_p(f)} h_i$ of $f$ is called reduced\index{standard writing!reduced} if either $p\nmid \deg_p(f)$ or $p\mid \deg_p(f)$ and each primitive monomial in $k[x_1^p,\ldots,x_n^p]$ of degree $\deg_p(f)$ has coefficient $0$ in $h_{\deg_p(f)}$.

\smallskip
{\bf (4)} \phantomsection{Let $r\in\mathbb N^{\ast}$. Assume $f\notin k[x_1^p,\ldots,x_n^p]$. We say that $f$ has $r$ points at infinity modulo $p$\index{points at infinity modulo $p$}, and we write $\infty_p(f)=r$, if there exists one (equivalently, if in each) reduced standard writing $f=g^p+\sum_{i=1}^{\deg_p(f)} h_i$ of $f$ such that $h_{\deg_p(f)}$ is a product $h^p\prod_{i=1}^r \ell_i^{s_i}$ with $h\in k[x_1,\ldots,x_n]\setminus\{0\}$, $(s_1,\ldots,s_r)\in \llbracket1,p-1\rrbracket^r$, and $\ell_1,\ldots,\ell_s$ non-associated irreducible forms (if $n=2$, they are linear forms).}\label{EXT9}

\smallskip
{\bf (5)} Let $e=\e(g_1,\ldots,g_n)\in\End_n(k)$ and $r\in\mathbb N^{\ast}$. We say that $e$ or $(g_1,\ldots,g_n)$ has $r$ points at infinity modulo $p$\index{endomorphism!has points at infinity modulo $p$} if for each $i\in \llbracket1,n\rrbracket$ we have $g_i\notin k_{\s}[x_1^p,\ldots,x_n^p]$ and $\infty_p(g_i)=r$.
\end{definition}

\begin{remark}\normalfont\label{R3}
Let $(g_1,g_2)\in k[x_1,x_2]^2$ and $f\in k[x_1,\ldots,x_n]\setminus k[x_1^p,\ldots,x_n^p]$. 

\medskip
{\bf (1)} If $\e(g_1,g_2)\in\EE_2(k)$, then $\infty_p(g_1)=1$ iff $\infty_p(g_2)=1$ by \cite{Lang-J}, Lem.\ 2.17. Thus $e(g_1,g_2)$ has $1$ point at infinity modulo $p$ iff $1\in\{\infty_p(g_1),\infty_p(g_2)\}$.%and each $d\in\EE_2(k)$ has $r$ point at infinity modulo $p$ with $r\{1,2\}$. 

\smallskip
{\bf (2)} If $n=2$ and $\infty_p(f)=1$, then $p\nmid\deg_p(f)$. 

\smallskip
{\bf (3)} If $p\nmid \deg_p(f)$ (resp.\ $p\mid \deg_p(f)$), then in each standard (resp.\ reduced standard) writing $f=g^p+\sum_{i=1}^{\deg_p(f)} h_i$ of $f$, the form $h_{\deg_p(f)}$ is uniquely determined.

\smallskip
{\bf (4)} For each $g\in k[x_1,\ldots,x_n]$, we have $\infty_p(f)=\infty_p(f+g^p)$. 

\smallskip
{\bf (5)} If $n=2$ and $p\nmid\deg_p(f)$, then $\infty_p(f)$ is as in \cite{Lang-J}, Def.\ 2.15. But if $p\mid\deg_p(f)$, then $\infty_p(f)$ is in general different from the one defined in loc.\ cit.: the change is made so that part (4) holds. If $p=2$ and $f:=x_1x_2(x_1+x_2)(x_1+\alpha x_2)$  with $\alpha\in k\setminus\{0,1\}$, then $f=x_1^3x_2+(\alpha+1)x_1^2x_2^2+\alpha x_1x_2^3$, $\deg(f)=\deg_p(f)=4$, in the reduced standard writings of $f$ we have $h_4=x_1^3x_2+\alpha x_1x_2^3=x_1x_2(x_1+\sqrt{\alpha}x_2)^2$, hence $\infty_2(f)=2$ while loc.\ cit.\ gives the value $4$.

\smallskip
{\bf (6)} Let $e\in\pMor_n(k)$. If $e$ is not an elementary $p$-morphism, then $p^2\mid\deg(e)$. If $e$ is an elementary $p$-morphism, then we have $\mathcal F(e)\subset\{1\}\cup p\mathbb N^{\ast}$. 
\end{remark}

%\begin{example}\normalfont\label{EX3.1}
%Let $m\in\mathbb N^{\ast}$. Let $e\in\EE_2(k)$ be defined by a pair $\bigl(x_1,x_2+f(x_1)x_2^{pm}\bigr)\in B^2$ with $f\in k_{\s}[x_1]\setminus k_{\s}[x_1^p]$. We show that the assumption that $e\in\BE_2(k)$ leads to a contradiction. This assumption implies that there exists $a=\e(y_1,y_2)\in\GA_2(k)$ such that for $d:=ea=\e\bigl(y_1,y_2+f(y_1)y_2\bigr)$, $D_{x_1}$ and $D_{x_2}$ are $k$-linear derivations of the $k$-subalgebra $A:=k_{\s}[y_1,y_2+f(y_1)y^2]$ of $B:=k_{\s}[x_1,x_2]$. For $(i,j)\in \{1,2\}^2$, let $g_{i,j}\in A$ be such that $D_{y_1}:=g_{11}D_{x_1}+g_{12}D_{x_2}$ and $D_{y_2}:=g_{21}D_{x_1}+g_{22}D_{x_2}$. As $D_{y_1}(y_1)=1$, we get that $g_{12}D_{x_2}(y_1)=1\end{example}

\begin{lemma}\label{L9}
Let $f\in k[x,y]\setminus k[x^p,y^p]$. Let $l:=\max\bigl(\deg(f_x),\deg(f_y)\bigr)\in\mathbb N$ and $m:=\infty_p(f)$. For $z\in\{x,y\}$ let $h_z\in k[x,y]$ be the homogeneous component of $f_z$ of degree $l$. Then $\infty_p(f)=1$ iff $l=m-1$ and one of the following partial derivatives conditions holds.\index{partial derivative condition}

\medskip
{\bf ($\partial_1$)} We have $h_xh_y\neq 0$ and $kh_x=kh_y=k\ell^sh_0^p$ with $h_0\in k[x,y]$, $s\in \llbracket0,p-2\rrbracket$, and $\ell\in k[x,y]$ a linear form.\index{partial derivative condition!($\partial_1$)}

\medskip
{\bf ($\partial_2$)} We have $h_xh_y=0$ and the non-zero factor $h_z$ with $z\in\{x,y\}$ is of the form $z^sh_0^p$ with $h_0\in k[x,y]$ and $s\in \llbracket0,p-2\rrbracket$.\index{partial derivative condition!($\partial_2$)}
\end{lemma}

\begin{proof}
 We consider a standard writing $f=g^p+\sum_{i=1}^m h_i$ of $f$; so each $h_i\in k[x,y]$ is homogeneous of degree $i$ and $h_m\notin k[x^p,y^p]$. 
 
Assume $\infty_p(f)=1$. Thus $h_m=\ell^qh^p$ with $h\in k[x,y]\setminus\{0\}$, $q\in \llbracket1,p-1\rrbracket$ and $\ell=\alpha x+\beta y$ with $(\alpha,\beta)\in k^2\setminus \{(0,0)\}$. We have $m=q+p\deg(h)$. Clearly, $l=m-1$. If $\alpha\neq 0$ (resp.\ $\alpha=0$), then $h_x=q\alpha\ell^{q-1}h^p$ (resp.\ $\beta\neq 0$, $h_y=q\beta\ell^{q-1}h^p$, and $h_x$=0). Similarly, if $\beta\neq 0$ (resp.\ $\beta=0$), then $h_y=q\beta\ell^{q-1}h^p$ (resp.\ $\alpha\neq 0$, $h_x=q\alpha\ell^{q-1}h^p$, and $h_y=0$). Thus for $\alpha\beta\neq 0$ (resp.\ $\alpha\beta=0$), Condition ($\partial_1$) (resp.\ ($\partial_2$)) holds with $s:=q-1$.

Assume now that $l=m-1$ and one of the two Conditions ($\partial_1$) and ($\partial_2$) holds. Let $r:=\infty_p(f)\in\mathbb N^{\ast}$. We assume that $f=g^p+\sum_{i=1}^m h_i$ is a reduced standard writing of $f$. We have a product decomposition $h_m=h^p\prod_{i=1}^r \ell_i^{s_i}$ with $h\in k[x,y]\setminus\{0\}$, $(s_1,\ldots,s_r)\in \llbracket1,p-1\rrbracket^r$, and the $\ell_i$s are non-proportional non-zero linear forms. Thus $l+1=m=p\deg(h)+\sum_{i=1}^r s_i$. Let $F:=\prod_{i=1}^r \ell_i^{s_i}$. We have $f_x=F_xh^p+L_x$ and $f_y=F_yh^p+L_y$, where for $z\in\{x,y\}$ we have $L_z\in k[x,y]$ with $\deg(L_z)\le m-2$. 

Assume Condition ($\partial_2$) holds. To fix the ideas we can assume that $z=x$. Thus $h_x=x^sh_0^p$ and $h_y=0$. So $h_x=x^sh_0^p=F_xh^p$; as $s\in \llbracket0,p-2\rrbracket$, it follows that $h|h_0$. Let $g_0\in k[x,y]$ be such that $h_0=hg_0$; so $h$, $h_0$, and $g_0$ are homogeneous polynomials. Hence $F_x=x^sg_0^p$. Therefore $F=\frac{x^{s+1}g_0^p}{s+1}+g_1(x^p,y)$ with $g_1\in k[x,y]$ such that $g_1(x^p,y)$ is homogeneous of degree $s+1+p\deg(g_0)\le p-1+p\deg(g_0)$. 

By reasons of degree we get that $g_1=y^{s+1}g_2$ with $g_2\in k[x,y]$ a polynomial such that $g_2(x^p,y)$ is homogeneous of degree $p\deg(g_0)$. This implies that there exists a homogeneous polynomial $g_3\in k[x,y]$ such that $g_2(x^p,y)=g_3(x,y)^p$; if $g_3\neq 0$, then $\deg(g_3)=\deg(g_0)$. Hence we have $F=\frac{x^{s+1}g_0^p}{s+1}+y^{s+1}g_3(x,y)^p$. As $h_y=0$ and $F_y=(s+1)y^sg_3(x,y)^p$, we get that $\deg(F_yh^p)<s+p\deg(h_0)=l=m-1$. Thus $\deg(g_3)<\deg(g_0)$ and hence $g_3=0$. So $g_1=0$. Thus $F=\frac{x^{s+1}}{s+1}$, hence $r=1$.

Assume Condition ($\partial_1$) holds. We write $\ell=\alpha x+\beta y$ with $(\alpha,\beta)\in k^2\setminus\{(0,0)\}$. To fix the ideas we can assume that $\alpha\neq 0$. By replacing $h_0$ with $\sqrt[p]{\alpha}h_0$ we can assume that $\alpha=1$. As in the previous paragraph we argue that there exists $g_3\in k[x,y]$ homogeneous such that $F=\frac{\ell^{s+1}g_0^p}{s+1}+y^{s+1}g_3(x,y)^p$. As $F_xh^p=kf_x=kf_y=F_yh^p$ we get that $ky^sg_3(x,y)^p\subset k\ell^sh_0^p$, thus $g_3=0$, hence $\infty_p(f)=1$.\end{proof}%Thus $\deg(F_x)=\deg(F_y)=\deg(F)-1$ and $kF_x=kF_y$. Let $q:=\deg(F)$ and $\delta\in k^{\ast}$ be such that $f_x=\delta F_y$. We write $F=\sum_{i=0}^q \alpha_ix^iy^{q-i}$ with $(\alpha_0,\ldots,\alpha_q)\in k^{q+1}$. As $\deg(F_y)=\deg(F)-1$, we have $p\nmid q$ and $\alpha_0\alpha_q\neq 0$. 

\begin{proposition}\label{PR1+}
Let $e\in\EE_2(k)$. Then the following properties hold.

\medskip
{\bf (1)} If $e$ has $1$ point at infinity modulo $p$, then for each $a\in\GA_2(k)$ the endomorphism $ae$ has $1$ point at infinity modulo $p$. 

\smallskip
{\bf (2) (Lang)} If $e\in\pMor_2(k)$, then $e$ has $1$ point at infinity modulo $p$. 

\smallskip
{\bf (3)} Each $e\in\BE_2(k)$ has $1$ point at infinity modulo $p$.
\end{proposition}

\begin{proof}
For part (1), let $(f_1,f_2)\in k_{\s}[x_1,x_2]$ be such that $e=\e(f_1,f_2)$; by definition we have $\infty_p(f_1)=\infty_p(f_2)=1$. Based on the amalgamated structure of $\GA_2(k)$ (see \cite{vdK}, Thms.\ 1 and 2; cf.\ also \cite{Kraft}, Sect.\ 5, Thm.\ 7), to prove part (1) we can assume that $a$ is defined by a pair $(F_1,F_2)\in k[x_1,x_2]$ such that either $\pi(F_1)=\pi(F_2)=1$ or $F_1=x_1$. If $F_1=x_1$, then $ae$, being defined by a pair $(f_1,f_3)$ with $f_3\in k_{\s}[x_1,x_2]$, has $1$ point at infinity modulo $p$. So we can assume that $\pi(F_1)=\pi(F_2)=1$.

Hence there exists $\begin{bmatrix} 
\gamma_{11} & \gamma_{12} \\
\gamma_{21} & \gamma_{22} \\ 
\end{bmatrix}\in\GL_2(k)$ and $(\gamma_1,\gamma_2)\in k^2$ such that we have $ae=\e(\gamma_{11} f_1+\gamma_{12}f_2+\gamma_1,\gamma_{21}f_1+\gamma_{22}f_2+\gamma_2)$. For each $i\in\{1,2\}$, let $m_i:=\deg_p(f_i)$ and a standard writing $f_i=g_i^p+\sum_{j=1}^{m_i} h_{i,j}$ of $f_i$ with $g_i\in k_{\s}[x_1,x_2]$ and each $h_{i,j}\in k_{\s}[x_1,x_2]$ a homogeneous polynomial of degree $j$. For each $i\in\{1,2\}$, we write $h_{i,m_i}=H_i^p\ell_i^{s_i}$ with $H_i\in k_{\s}[x_1,x_2]\setminus\{0\}$, $s_i\in \llbracket1,p-1\rrbracket$, and $\ell_i=\alpha_i x_1+\beta_i x_2$ for a pair $(\alpha_i,\beta_i)\in k^2\setminus\{(0,0)\}$. 

To show that $1\in\{\infty_p(\gamma_{11} f_1+\gamma_{12}f_2+\gamma_1),\infty_p(\gamma_{21} f_1+\gamma_{22}f_2+\gamma_2)\}$ we can assume that $\gamma_1=\gamma_2=0$ and we consider two cases as follows.

{\bf Case 1: $m_1\neq m_2$.} To fix the ideas we can assume that $m_1>m_2$. Let $i\in\{1,2\}$ be such that $\gamma_{i,1}\neq 0$. We have $\gamma_{i1} f_1+\gamma_{i2}f_2=\gamma_{i1}g_1^p+\gamma_{i2}g_2^p+\sum_{j=1}^{m_1-1} h_{i,j,a}+\gamma_{i1}H_1^p\ell_1^{s_1}$ with each $h_{i,j,a}\in k_{\s}[x_1,x_2]$ homogeneous of degree $j$. Thus $\infty_p(\gamma_{i1} f_1+\gamma_{i2}f_2)=1$.

{\bf Case 2: $m_1=m_2$.} The case $m_1=1$ is trivial, so we assume that $m_1\ge 2$. As $m_i\equiv s_i\;(\textup{mod}\; p)$, we have $s_1=s_2$, thus $\deg(H_1)=\deg(H_2)$. As $p\nmid m_1m_2$, from \cite{Lang-J}, Prop.\ 2.16 we get that $\ell_1$ and $\ell_2$ are proportional. Thus we can assume that $\ell_1=\ell_2$. As $\gamma_{11}\gamma_{22}-\gamma_{12}\gamma_{21}\in k^{\ast}$, either $\gamma_{11}H_1^p+\gamma_{12}H^2_2$ or $\gamma_{21}H_1^p+\gamma_{22}H^2_2$ is $H_a^p$ with $H_a$ a non-zero homogeneous polynomial of degree $\deg(H_1)$. Therefore $1\in\{\infty_p(\gamma_{11} f_1+\gamma_{12}f_2),\infty_p(\gamma_{21} f_1+\gamma_{22}f_2)\}$ as one can easily check based either on an argument similar to the one of Case 1 or on Lemma \ref{L9}. 

As $1\in\{\infty_p(\gamma_{11} f_1+\gamma_{12}f_2),\infty_p(\gamma_{21} f_1+\gamma_{22}f_2)\}$ in both cases, part (1) holds by Remark \ref{R3}(1). 

See \cite{Lang-J}, Thm.\ 2.18 for part (2).\footnote{The proof of part (1) follows closely the proof of \cite{Lang-J}, Thm.\ 2.18.}

To prove part (3) we write $e=bdc$ with $(b,c)\in\GA_2(k)^2$ and $d\in\EE_2(k)$ defined by a pair $\bigl(x_1+g_1(x_1^p,x_2^p),x_2+g_2(x_1^p,x_2^p)\bigr)$ with $(g_1,g_2)\in k_{\s}[x_1,x_2]^2$. Based on part (1) we can assume that $b=1_{\mathbb A^2_k}$, thus $e=dc$. Let $(y_1,y_2)\in\mathcal A_2(k)$ be such that $c=\e(y_1,y_2)$, then $dc=\e\bigl(y_1+f_1(y_1^p,y_2^p),y_2+f_2(y_1^p,y_2^p)\bigr)$. As $c$ has $1$ point at infinity modulo $p$ by part (2), we have $\infty_p(y_1)=\infty_p(y_2)=1$ and therefore $dc$ has $1$ point at infinity modulo $p$. Thus part (c) holds.\end{proof}

For $p>2$, see Example \ref{EX22.9} for pairs $(e,a)\in\EE_2(k)\times\GA_2(k)$ such that $e$ and $ea$ have $1$ point and $2$ points (respectively) at infinity modulo $p$.

\section{On complete intersections}\label{S7}

For a scheme $W$, let $\mathcal O_W$ be its structure ring sheaf. Let $\mathcal O(W):=\mathcal O_W(W)$. 

\phantomsection{For a variety $W$ over $K$, the partial order on the set}\label{PH100} 
$$\Triv(W):=\{Z\subset W|Z\;\textup{is affine open regular}\;\textup{and}\;T_Z\;\textup{is trivial}\}$$
is the inclusion. After recalling two definitions, we introduce four more.

\begin{definition}\label{D5}
Let $X=\Spec R$ be a non-empty equidimensional affine variety over $K$.

\medskip
{\bf (1)} Suppose that $\chr(K)=p$. By a differential basis\index{differential basis} of $R$ we mean a subset $\{f_1,\ldots,f_{\dim(X)}\}$ of $R$ such that $(df_1,\ldots,df_{\dim(X)})$ is an $R$-basis of the $R$-module $\Omega_{R}=\Omega_{R/K}$ of (relative) differentials.

\smallskip
{\bf (2)} A closed embedding $X\rightarrow\mathbb A^n_K$, i.e., an identification $R=K[x_1,\ldots,x_n]/I$ of $K$--algebras, is called a complete intersection\index{complete intersection} if the ideal $I$ is generated by $n-\dim(X)$ elements.

\smallskip
{\bf (3)} The complete intersection embedding dimension\index{complete intersection!embedding dimension} of $X$ is $\ciedim(X)\in\mathbb N\cup\{\infty\}$ defined as follows. If there exists a complete intersection closed embedding $X\rightarrow\mathbb A^N_K$, then $\ciedim(X)$ is the smallest possible value for $N$. If no such $N$ exists, then $\ciedim(X):=\infty$. We say that $X$ is a complete intersection if $\ciedim(X)\in\mathbb N$.

\smallskip
{\bf (4)} Suppose that $X$ is regular. By the weak complete intersection embedding dimension\index{complete intersection!weak embedding dimension} of $X$ we mean $\wciedim(X)\in\mathbb N\cup\{\infty\}$ defined as follows. If there exists a closed embedding $\imath:X\rightarrow\mathbb A^N_K$ whose normal bundle $N_{\imath}$ is trivial, then $\wciedim(X)$ is the smallest possible value for $N$. If no such $N$ exists, then $\wciedim(X):=\infty$.

\smallskip
{\bf (5)} Suppose that $\ciedim(X)\in\mathbb N$. By the strong complete intersection embedding dimension\index{complete intersection!strong embedding dimension} of $X$ we mean the smallest $\sciedim(X)\in\mathbb N$ such that the inequality $\ciedim(X)\le\sciedim(X)$ holds and each closed embedding $X\rightarrow\mathbb A^{\sciedim(X)+m}_K$ with $m\in\mathbb N$ is a complete intersection.

\smallskip
{\bf (6)} Suppose that $X$ is normal connected and contains an open subvariety isomorphic to $\mathbb A^{\dim(X)}_K$. Let $\rho_{\et}(X)\in\mathbb N$ be the largest $n$ such that there exists $U\in\Triv(X)$ and an open subvariety $V$ of $U$ for which we have $V\cong\mathbb A^{\dim(X)}_K$ and $|\Irr(U\setminus V)|=n$.
\end{definition}

The existence of $\sciedim(X)$ is implied by Theorem \ref{T3}(3.a) below. 

We first group together several results in the literature which also give the existence of $\sciedim(X)$ and then apply them in the context of $X_e$ with $e\in\EE_n(K)$.

\begin{theorem}\label{T3}
Let $X$ be a non-empty equidimensional affine regular variety over $K$. Then the following properties hold.

\medskip
{\bf (1)} We have $\edim(X)\le 2\dim(X)+1$. 

\smallskip
{\bf (2)} We have $\ciedim(X)\in\mathbb N$ iff $T_X$ is trivial.

\smallskip
{\bf (3)} If $T_X$ is trivial and $X$ is integral, then the following properties hold.

\medskip\noindent
{\bf (3.a)} We have inequalities
$$\edim(X)\le\wciedim(X)\le\ciedim(X)\le\sciedim(X)\le 2\dim(X)+2$$ 
with $\ciedim(X)-\wciedim(X)\in\{0,1\}$.

\smallskip\noindent
{\bf (3.b)} If $\edim(X)=2\dim(X)$ and $\dim(X)\ge 1$, then we have an inequality $\ciedim(X)\le 2\dim(X)+1$ if $\dim(X)\ge 3$ and an identity $\ciedim(X)=2\dim(X)$ if $\dim(X)\in\{1,2\}$. 

\smallskip\noindent
{\bf (3.c)} If $\edim(X)=2\dim(X)-l$ with $l\in\llbracket1,\dim(X)-1\rrbracket$ and if for $\dim(X)\ge 3+l$ the factorial number $\bigl(\dim(X)-1\bigr)!$ is invertible in $K$, then we have an inequality $\ciedim(X)\le 2\dim(X)$ if $\dim(X)\ge 3+l$ and an identity $\ciedim(X)=2\dim(X)-l$ if $\dim(X)\in \{1+l,2+l\}$. 

\medskip
{\bf (4)} Let $R$ be a finitely generated $K$-algebra with equidimensional spectrum. If $\chr(K)=p$, then $R$ has a $p$-basis iff $R$ is regular and has a differential basis.\end{theorem}

\begin{proof}
Part (1) is well-known (e.g., see \cite{Sr}, Cor.\ 1). 

Clearly, $\edim(X)\le\wciedim(X)$. 

In this paragraph we assume that $\ciedim(X)\in\mathbb N$. For each $m\in\mathbb N$ there exists a complete intersection closed embedding $X\rightarrow\mathbb A^{\ciedim(X)+m}_K$. If the inequality holds $\ciedim(X)+m\ge 2\dim(X)+2$, then all closed embeddings $X\rightarrow\mathbb A^{\ciedim(X)+m}_K$ form one orbit under the natural right action of $\GA_{\ciedim(X)+m}(K)$ by \cite{Sr}, Thm.\ 2, and hence are complete intersections. Thus $\sciedim(X)\le 2\dim(X)+2$.\footnote{If $X$ is integral, then $\sciedim(X)\le 2\dim(X)+2$ by \cite{Kum1}, Thm.2) or \cite{Kum2}, Cor.\ 4.} Clearly, $\wciedim(X)\le\ciedim(X)$. 

For a closed embedding $\jmath:X\rightarrow\mathbb A^{\edim(X)+\epsilon}_K$ with $\epsilon\in\mathbb N$ we have a normal bundle $N_{\jmath}$ over $X$ defined by a short exact sequence $0\to T_X\rightarrow\jmath^*(T_{\mathbb A^{\edim(X)+\epsilon}_K})\rightarrow N_{\jmath}\rightarrow 0$ of vector bundles over $X$. If $\jmath$ is a complete intersection, then the dual of $N_{\jmath}$, called the conormal bundle of $\jmath$, is trivial and thus $N_{\jmath}$ is trivial; thus $T_{X}$ is stable free and hence trivial by \cite{Su}, Thm.\ 1. Hence the `only if' of part (2) holds.

To prove the `if' of part (2) and part (3), we can assume that we have $\dim(X)>0$, $X\not\cong\mathbb A_K^{\dim(X)}$, and $T_X$ is trivial. We write $\edim(X)=2\dim(X)+1-i$ with $i\in\llbracket0,\dim(X)\rrbracket$ by part (1). As $N_{\jmath}$ has rank $\dim(X)+1-i+\epsilon\ge 1$ and is stable free by the above short exact sequence, it follows that it is trivial if $\epsilon\ge i$ by \cite{Ba1}, Thm.\ 9.3, if $\epsilon=i-1$ by \cite{Su}, Thm.\ 1, and if $\epsilon=i-2$ and for $\dim(X)\ge 4$ the factorial $\bigl(\dim(X)-1\bigr)!$ is invertible in $K$ by \cite{FRS}, Thm. 

Let $\n(X):=\wciedim(X)-\edim(X)$. From the prior paragraph it follows that for $i\in\{0,1\}$ we have $\n(X)=0$ and for $i\ge 2$ we have $\n(X)\le i-1$; moreover, if $i\ge 2$ and either $\dim(X)\in\{2,3\}$ or $\dim(X)\ge 4$ and $\bigl(\dim(X)-1\bigr)!$ is invertible in $K$, then $\n(X)\le i-2$. If $i\ge 2$, whenever we are in the context of part (3.c) we have $l=i-1$.

Let $\iota:X\rightarrow\mathbb A^{\edim(X)+\n(X)}_K$ be a closed embedding with trivial normal bundle $N_{\imath}$. Clearly, $\ciedim(X)\ge\edim(X)+\n(X)$. As the dual of $N_{\imath}$ is also trivial, $X$ is a connected component of a complete intersection $X_+=X_+(\imath)$ in $\mathbb A^{\edim(X)+\n(X)}_K$. 

If $\dim(X)-i+\n(X)=0$, then $X$ is a hypersurface in $\mathbb A^{\dim(X)+1}_K=\mathbb A^{\edim(X)+\n(X)}_K$ and thus we can assume that $X_+=X$. If $\dim(X)-i+\n(X)=1$, then $X$ is a hypersurface in $\mathbb A^{\dim(X)+2}_K=\mathbb A^{\edim(X)+\n(X)}_K$ and by \cite{Kum1}, Thm.1) we can assume that $X_+=X$. Thus for $\dim(X)-i+\n(X)\in\{0,1\}$ we have $X_+=X$ and $\ciedim(X)\le \edim(X)+\n(X)$, hence $\ciedim(X)=\edim(X)+\n(X)$. Similarly, if $\dim(X)-i+\n(X)\ge 2$ and $X_+=X$, then $\ciedim(X)=\edim(X)+\n(X)$.

Assume now that $\dim(X)-i+\n(X)\ge 2$ and $W=W(\imath):=X_+\setminus X$ is non-empty for each embedding $\imath$ (i.e., and $\ciedim(X)\ge \edim(X)+\n(X)+1$). For $(f_1,\ldots,f_{\dim(X)-i+\n(X)+1},g)\in K[x_1,\ldots,x_{\edim(X)+\n(X)+1}]^{\dim(X)-i+\n(X)+2}$ with\break $g(W)=\{0\}$, $g(X)=\{1\}$ and the ideal $(f_1,\ldots,f_{\dim(X)-i+\n(X)+1})$ defines $X_+$, 
$$\Spec\bigl(K[x_1,\ldots,x_{\edim(X)+\n(X)+1}]/(f_1,\ldots,f_{\dim(X)-i+\n(X)+1},1-gx_{\edim(X)+\n(X)+1})\bigr)$$ 
is isomorphic to $X$, hence $\ciedim(X)= \edim(X)+\n(X)+1$. Therefore if we have $\edim(X)=2\dim(X)-l$ with $l\in\llbracket1,\dim(X)-1\rrbracket$, then from $\dim(X)\ge i+2+\n(X)$ we get that $\dim(X)\ge l+3+\n(X)\ge l+3$; if moreover $\bigl(\dim(X)-1\bigr)!$ is invertible in $K$, then $\n(X)\le i-2=l-1$ and hence $\ciedim(X)\le\edim(X)+l=2\dim(X)$.

The last two paragraphs and the above relations on $\n(X)$ imply that the `if' of part (2) and part (3) holds.

For part (4), if $R$ has a $p$-basis then $R$ is regular by \cite{Kun}, Cor.\ 2.7 applied to the local rings of $R$. Thus to prove part (4) we can assume that $R$ is regular and, due to the equidimensional hypothesis, an integral domain. Thus part (4) holds by either \cite{KN}, Thm.\ or \cite{Aba}, Prop.\ 3.2 and Cor.\ 4.2.\end{proof}

\begin{theorem}\label{T3+}
Let $e\in\EE_n(K)$ with $n\ge 2$. Let $Y=\Spec R$ be an open affine subvariety of $X_e$ that contains $\Imm(\imath_e)$.

\medskip
{\bf (1)} The following statements are equivalent.

\medskip\noindent
{\bf (1.a)} We have $Y\subset X_e^{\et}$.

\smallskip\noindent
{\bf (1.b)} The generic points of the irreducible components of $Y\setminus\Imm(\imath_e)$ belong to $X_e^{\et}$.

\smallskip\noindent
{\bf (1.c)} We have $Y=\Reg(Y)$ and $T_Y$ is trivial (i.e., $Y\in\Triv(X_e)$).

\smallskip\noindent
{\bf (1.d)} The tangent bundle $T_{\Reg(Y)}$ over $\Reg(Y)$ is trivial.

\smallskip\noindent
{\bf (1.e)} We have $Y=\Reg(Y)$ and $\ciedim(Y)\in\mathbb N$.

\medskip
{\bf (2)} If $\chr(K)=p$, then statements (1.a) to (1.e) of part (1) are also equivalent to either one of the following extra statements.

\medskip\noindent
{\bf (2.a)} The ring $R$ has a $p$-basis.

\smallskip\noindent
{\bf (2.b)} The ring (or the $K$-algebra) $R$ has a differential basis.

\medskip
{\bf (3)} The open regular subset $X_e^{\et}$ of $X_e$ is a maximal element of $\Triv(X_e)$.

\smallskip
{\bf (4)} Assume $n=2$, $\psi_e$ is non-finite \'etale, and $\Imm(\imath_e)\neq Y$. If $\edim(Y)\in\{3,4\}$, then $\ciedim(Y)=\edim(Y)$. If $\edim(Y)=5$, then $\ciedim(Y)\in\{5,6\}$.

\smallskip
{\bf (5)} We have the relations $\rho_{\et}(e)=\rho_{\et}(X_e^{\et})\le\rho_{\et}(X_e)$.

\smallskip
{\bf (6)} If $\psi_e$ is regular but non-\'etale, then $\ciedim(X_e)=\infty$.
\end{theorem}

\begin{proof}
\phantomsection{To prove part (1), let $Z:=\Reg(Y)$. As $Y$ is normal, we have an inequality $\dim(Y\setminus Z)\le \dim(Y)-2$. Let $\psi_Y:Y\rightarrow\mathbb A^n_{K,\t}$ and $\psi_Z:Z\rightarrow\mathbb A^n_{K,\t}$ be defined by $\psi_e$ via restrictions of the source and target; we have a functorial morphism $T(\psi_Z):T_Z\rightarrow\psi_Z^*(\mathbb A^n_{K,\t})$ between bundles over $Z$ whose target is trivial. Let $Y^{\et}:=Y\cap X_e^{\et}$; it is the \'etale locus of $\psi_Y$. See \cite{Stacks1} or \cite{GZ}, Thm.\ 1 for the affineness of $X_e^{\et}$ and hence also of $Y^{\et}$. Thus the complement $Y\setminus Y^{\et}$ is either empty or of pure codimension $1$ (see \cite{Gro6}, Exp.\ V, Ex.\ 3.4). Hence $Y=Y^{\et}$ iff $\dim(Y\setminus Y^{\et})\le \dim(Y)-2$.}\label{PH101}

Clearly, $(1.a)\Rightarrow (1.b)$ and $(1.c)\Rightarrow (1.d)$. Moreover, $(1.c)\Leftrightarrow (1.e)$ by Theorem \ref{T3}(2). 

If $(1.b)$ holds, then $\dim(Y\setminus Y^{\et})\le\dim(Y)-2$ and hence $Y=Y^{\et}$; thus $Z=Y$, $\psi_Y$ is \'etale and therefore $T(\psi_Y)$ is an isomorphism. Hence $(1.b)\Rightarrow (1.c)$. 

If $(1.d)$ holds, then the morphism $T(\psi_Z)$ is between trivial bundles of dimension $n$ and an isomorphism over $\Imm(\imath_e)$. Thus its determinant is a global function on $Z$ that maps to $\mathcal O\bigl(\Imm(\imath_e)\bigr)=K^{\ast}$ and hence $T(\psi_Z)$ is an isomorphism, which implies that $\psi_Z$ is \'etale. Hence $Z\subset Y^{\et}$. Thus $\dim(Y\setminus Y^{\et})\le \dim(Y)-2$ and hence $Y=Y^{\et}$. Thus $(1.d)\Rightarrow (1.a)$. 

We conclude that part (1) holds.

For part (2) we assume that $\chr(K)=p$. Clearly, $(1.a)\Rightarrow (2.a)$.

We have $(2.a)\Leftrightarrow (2.b)$ by Theorem \ref{T3}(4). 

We begin to show that $(2.a)\Rightarrow (1.a)$ by considering the $K$-algebra homomorphism $K_{\t}[x_1,\ldots,x_n]\rightarrow R$ that defines $\psi_Y$ and we view it as an inclusion. We also consider the $R^p$-subalgebra 
$$S:=\sum_{(i_1,\ldots,i_n)\in \llbracket0,p-1\rrbracket^n} R^p\prod_{j=1}^n x_j^{i_j}$$ 
of $R$. From the equivalence $(N1)\Leftrightarrow (N3)$ of Section \ref{S2} applied to $e$ it follows that in fact we have a direct sum decomposition $S=\oplus_{(i_1,\ldots,i_n)\in \llbracket0,p-1\rrbracket^n} R^p\prod_{j=1}^n x_j^{i_j}$ of $R^p$-modules. The assumption that $R$ has a $p$-basis implies that the $R^p$-module $R$ is free of rank $p^n$. The inclusion $S\rightarrow R$ is an $R^p$-linear transformation between free modules of rank $p^n$; its determinant $u\in R^p$ is uniquely determined up to multiplication by an element in $R^{\ast}=K^{\ast}$. The mentioned equivalence also gives that $u\in (R^p)^{\ast}=K^*$. Thus $S=R$, hence $R$ has a $p$-basis formed by elements of $K_{\t}[x_1,\ldots,x_n]$. This implies that $\psi_Y$ is unramified and thus, being also flat in codimension at most $2$ (as so is $\psi_e$), it is \'etale in codimension at most $2$, i.e., $\dim(Y\setminus Y^{\et})\le\dim(Y)-2$. Thus $Y=Y^{\et}$ and hence $(2.a)\Rightarrow (1.a)$. We conclude that part (2) holds.

As $X_e^{\et}$ is an affine variety, we have $X_e^{\et}\in\Triv(X_e)$ by part (1) applied to $Y:=X_e^{\et}$. If $W\in\Triv(X_e)$ contains $X_e^{\et}$, then $W\subset X_e^{\et}$ by part (1) applied to $Y:=W$, hence $W=X_e^{\et}$. Thus part (3) holds. 

Part (4) follows from Theorem \ref{T3}(2.b) and (3.c) and part (1). 

Based on part (3), the inequality $\rho_{\et}(X_e^{\et})\le\rho_{\et}(X_e)$ follows from definitions. For an open subvariety $U$ of $X_e^{\et}$ with $U\cong\mathbb A^n_K$, $Cl(X_e^{\et})\cong\mathbb Z^{|\Irr(X_e^{\et}\setminus U)|}$. As $|\Irr(X_e^{\et}\setminus U)|$ does not depend on $U$, it is $\rho_{\et}(X_e^{\et})$. When $U=\Imm(\imath_e)$ we get that $Cl(X_e^{\et})\cong \mathbb Z^{\rho_{\et}(e)}$ and $\rho_{\et}(X_e)=\rho_{\et}(e)$. So part (5) holds. 

Part (6) follows from Theorem \ref{T3}(2) and part (3).\end{proof}

\section{On models of $\mathbb P^1_{\mathcal K}$}\label{S8}

Let $D$ be a discrete valuation ring.  Let $\mathcal K:=\Frac(D)$. In this section we study normal projective flat models $X$ of $\mathbb P^1_{\mathcal K}$ over $\Spec D$; so $X_{\mathcal K}=\mathbb P^1_{\mathcal K}$ and $X$ is integral. In practice, $D$ is (the completion) of a local ring $K[x]_{(x-\alpha_0)}$ with $\alpha_0\in K$ and our $X$s are obtained naturally from suitable subschemes of some $X_e$ with $e\in\EE_2(K)$, but for the sake of generality, we work over more general settings that involve contractions of curves and Stein factors. 

Let $D^{\h}$ be the henselization of $D$ and $\mathcal K^{\h}:=\Frac(D^{\h})$. Let $\kappa$ be the residue field of $D$ (hence also of $D^{\h})$. 

The  normal projective flat models $X$ of $\mathbb P^1_{\mathcal K}$ correspond to finite non-empty sets of divisorial valuations of the field of fractions $\mathcal K(t)$ of $\mathbb P^1_{\mathcal K}$, which are (equivalent classes of) discrete valuations coming from irreducible curves in the special fiber $X_{\kappa}$ of $X$. The terminology `divisorial valuations' is more used and established and pertains to normal models $X$ which are known to be of finite type over $\Spec D$ (see \cite{Em}, \cite{Pi}, \cite{Sc}, \cite{HS}, \cite{HKP}, \cite{GP}, etc.). In more general contexts (without finiteness of the normalization) one can speak about ``divisorial valuations" which are closely related to ``Zariski prime divisors"\index{Zariski prime divisor} in the sense of \cite{Ga1}, p.\ 108.  

We formalize divisorial valuations as follows. 

\begin{definition}\label{D4.9}
Let $W$ be a separated noetherian integral scheme with field of fractions $\mathcal K$ and let $\mathcal K_1$ be a finitely generated field extension of $\mathcal K$. Then a discrete valuation of $\mathcal K_1$ with valuation ring $\mathcal O_1$ is said to be a divisorial valuation\index{divisorial valuation} of $\mathcal K_1/W$ iff the morphism $\Spec\mathcal K_1\rightarrow\Spec\mathcal K$ extends to a morphism $\Spec\mathcal O_1\rightarrow W$ and the $W$-scheme $\Spec\mathcal O_1$ is the normalization of an integral $W$-scheme essentially of finite type (equivalently, iff there exists a finite type model $W_1\rightarrow W$ of $\Spec\mathcal K_1\rightarrow\Spec\mathcal K$ such that $\mathcal O_1$ corresponds to a Zariski prime divisor of $W_1$ in the sense of \cite{Ga1}, p.\ 108).
\end{definition}

\begin{remark}\normalfont\label{R3.9} 
{\bf (1)} If $W$ is universally Japanese, then in Definition \ref{D4.9} we can omit the words ``the normalization of a $W$-scheme." 

\smallskip
{\bf (2)} Suppose that $W$ is universally catenary. Then the equivalent conditions of Definition \ref{D4.9} are also equivalent to the condition that $\mathcal O_1$ dominates a local ring $\mathcal O_{W,x}$ of $W$ at a point $x\in W$ and the transcendence degree of the residue field of $\mathcal O_1$ over the residue field of $x$ is equal to $\textup{trdeg}_{\mathcal K} \mathcal K_1+\dim(\mathcal O_{W,x})-1$. 

\smallskip
{\bf (3)} If $W=\Spec A$ is affine, then the set of divisorial valuations with respect to $A$ defined in \cite{HS}, Def.\ 9.3.1 is a subset of the set of divisorial valuations of $\mathcal K_1/W$ defined in Definition \ref{D4.9}. The two sets coincide iff $A$ is universally catenary.

\smallskip
{\bf (4)} If $\mathcal K$ is algebraically closed, the divisorial valuations of $\mathcal K_1/\mathcal K$ are called Zariski prime divisors in \cite{Pop}.
\end{remark}

In what follows, all divisorial valuations are of $\mathcal K(t)/\Spec D$ and dominate $D$. In simpler words, they are discrete valuation rings of the field of fractions $\mathcal K(t)$ of $\mathbb P^1_D$ that dominate $D$ and their residue fields are transcendental over $\kappa$. It is known that the normalizations of the finite type models of $\mathcal K(t)/D$ are finite (e.g., by \cite{GLL}, Prop.\ 7.6(a)).

\begin{remark}\normalfont\label{R4}
{\bf (1)} Let $A$ be a henselian local ring of residue field $\kappa$. Let $Y$ be a proper algebraic space over $\Spec A$ of relative dimension $\le 1$. The special fiber $Y_{\kappa}$ of $Y$ is a proper algebraic space of dimension $\le 1$ over $\Spec\kappa$ and hence it is a proper scheme over $\Spec(\kappa)$ by either \cite{Stacks6} or \cite{Kn}, Ch.\ 5, Thm.\ 4.9 and Ch.\ 2, Cor.\ 6.16. Being proper and of dimension $\le 1$, $Y_{\kappa}$ is projective over $\Spec(\kappa)$. If $Y$ is a scheme, then $Y$ is projective by \cite{Gro5}, Cor.\ (21.9.12) and Rmk.\ (21.9.13).\footnote{\label{FOOT6}In general, one can lift an ample line bundle on $Y_{\kappa}$ to an ample line bundle on $Y$ as in the proof of the proper base change theorem of \'etale cohomology (see \cite{SGA4-3}, Exp.\ XVIII, Cor.\ 3.2), hence $Y$ is projective over $\Spec A$ and in particular it is a scheme.}%The reductions of $Y$ modulo powers of the maximal ideal of $D$ are also schemes by \cite{Kn}, Ch.\ 3, Cor.\ 3.6 and thus projective schemes by \cite{Gro5}, Cor.\ (21.9.12) and Rmk.\ 21.9.13. 

\smallskip
{\bf (2)} If $Z$ is a proper algebraic space over $\Spec D$ with $Z_{\mathcal K}=\mathbb P^1_{\mathcal K}$ and special fiber of dimension at most $1$ (e.g., if $Z_{\mathcal K}$ is Zariski dense in $Z$) and $Z_{D^{\h}}:=Z\times_{\Spec D} \Spec D^{\h}$ is a scheme, then we have an ample line bundle over $Z_{D^{\h}}$ by part (1); as the pullback homomorphism $\Pic(\mathbb P^1_{\mathcal K})\rightarrow\Pic(\mathbb P^1_{\mathcal K^{\h}})$ is an isomorphism, it descends to $Z$. The descended line bundle is still ample and thus $Z$ is projective.\footnote{Based on Footnote \ref{FOOT6}, the same holds without assuming that $Z_{D^{\h}}$ is a scheme.}

\smallskip
{\bf (3)} Following \cite{BLR}, Ch.\ 6, Sect.\ 6.7, assume that $\kappa$ is algebraically closed but not an algebraic closure of a finite field and $D=\kappa[t]_{(t)}$. Let $(f,g)\in\kappa[x,y,z]^2$ be a pair of non-zero non-proportional homogeneous polynomials of degree $3$. We assume that the zero locus $f=0$ is smooth and has a $\kappa$-valued inflection point $P$. The zero locus $f+tg=0$ is a cubic $X$ in $\mathbb P^2_D$ and we endow the special fiber $X_{\kappa}$ of $X$ with the group law whose identity element is $P$. Let $W$ be the blow up of $X$ at a $\kappa$-valued point of infinite order modulo the subgroup generated by the zero locus $f=g=0$. The blow down $f:W\rightarrow Z$ of the proper transform of $X_{\kappa}$ in the category of algebraic spaces (see \cite{Em}, Cor.\ or \cite{Pi}, Cor.\ 4 for the pullback of $f$ to $\Spec D^{\h}$ which is a morphism of projective schemes), gives a normal proper flat algebraic space model $Z$ of the generic fiber of $W$ which is not a scheme.\end{remark}

Recall the following definition introduced in \cite{Kl}, Sect.\ 5, Def.\ 2. 

\begin{definition}\label{D5a}
Let $R$ be a local ring. Let $W$ be a proper scheme over $\Spec R$ of relative dimension $\le 1$. By a Stein factor\index{Stein factor} of $W$ we mean a proper morphism $f:W\rightarrow V$ of schemes over $\Spec R$ such that the homomorphism $\mathcal O_V\rightarrow f_*(\mathcal O_W)$ is an isomorphism. An isomorphism between two Stein factors $f_1:W\rightarrow V_1$ and $f_2:W\rightarrow V_2$ of $W$ is an isomorphism $g:V_1\rightarrow V_2$ over $\Spec R$ such that $g\circ f_1=f_2$. 
\end{definition}

\begin{proposition}\label{L9a}
Let $R$ be a local ring of residue field $\kappa$. Let $W$ be a proper scheme over $\Spec R$ of relative dimension $1$. Let $\mathcal B$ be the set of $1$-dimensional irreducible components of the special fiber $W_{\kappa}$ of $W$. Then the following properties hold.

\medskip
{\bf (1)} Each Stein factor $f:W\rightarrow V$ of $W$ is uniquely determined up to isomorphism by the subset $\mathcal C_f:=\{C\in \mathcal B|f(C)\;\textup{is one point}\}$ of $\mathcal B$. 

\smallskip
{\bf (2)} Suppose that $R=D$ is a discrete valuation ring and $W$ is flat over $\Spec D$ with an integral generic fiber. Then the following properties hold.

\medskip\noindent
{\bf (2.a)} For each subset $\mathcal C$ of $\mathcal B$ that is distinct from $\mathcal B$ there exists a Stein factor $f^{\h}_{\mathcal C}:W_{D^{\h}}\rightarrow V_{D^{\h}}$ of $W_{D^{\h}}$ such that $\mathcal C_{f^{\h}_{\mathcal C}}=\mathcal C$ and its generic fiber is an isomorphism and thus $f^{\h}_{\mathcal C}$ descends to a morphism $f_{\mathcal C}:W\rightarrow V$ of algebraic spaces over $\Spec D$.\footnote{The Stein factor $f_{\mathcal B}:W\rightarrow V$ that contracts $W_{\kappa}$ exists but $V=\Spec \mathcal O(W)$ is finite flat over $\Spec D$, so the generic fiber of $f_{\mathcal B}$ is not an isomorphism.}

\smallskip\noindent
{\bf (2.b)} If $W$ is normal, then $V_{D^{\h}}$ and $V$ are normal.

\smallskip\noindent
{\bf (2.c)} If $W_{\mathcal K}=\mathbb P^1_{\mathcal K}$, then $f_{\mathcal C}:W\rightarrow V$ is projective. 

\smallskip\noindent
{\bf (2.d)} Let $Q\subset Z$ be a closed subset of pure codimension $1$ with a non-empty generic fiber. The morphism $Z\setminus Q\rightarrow\Spec\mathcal O(Z\setminus Q)$ is proper and is an isomorphism iff $Q$ meets all irreducible components of the special fiber. In particular, if $Z_{\kappa}$ is irreducible, then $Z\setminus Q\rightarrow\Spec\mathcal O(Z\setminus Q)$ is an isomorphism.
\end{proposition}

\begin{proof}
Part (1) is stated in \cite{Sc}, Sect.\ 2 when $R$ is noetherian by just quoting \cite{Gro1}, Thm.\ (5.4.1) but the last reference does not suffice. If $R$ is a discrete valuation ring and $W$ is flat over $\Spec R$, then part (1) holds by \cite{Pi}, Prop.\ 4. 

In general, let $f_1:W\rightarrow V_1$ and $f_2:W\rightarrow V_2$ be two Stein factors of $W$ such that $\mathcal C_{f_1}=\mathcal C_{f_2}$. Let $V_3$ be the schematic image of the product morphism $f_1\times f_2:X\rightarrow V_1\times_{\Spec R} V_2$. For $i\in\{1,2\}$, the $i$-th projection $\pi_i:V_3\rightarrow V_i$ is finite on special fibers over $\Spec\kappa$, hence it has finite fibers by Chevalley's Semi-continuity Theorem (see \cite{Gro4}, Thm.\ (13.1.3)). Thus $\pi_i$ is finite; by the Stein property it is an isomorphism. Hence $g:=\pi_2\circ\pi_1^{-1}:V_1\rightarrow V_2$ is an isomorphism between $f_1:W\rightarrow V_1$ and $f_2:W\rightarrow V_2$. Thus part (1) holds. 

Part (2.a) holds by \cite{Em}, Cor.\ or \cite{Pi}, Cor.\ 4. 

Part (2.b) holds as $V_{D^{\h}}$ is normal by \cite{Pi}, Prop.\ 3. 

Part (2.c) follows from part (2) and Remark \ref{R4}(2) applied to the algebraic space $V$ over $\Spec D$ and the scheme $V_{D^{\h}}$ over $\Spec D^{\h}$. 

Part (2.d) follows from the arguments in the proof of \cite{Sc}, Thm.\ 2.2.\end{proof}

\begin{corollary}\label{C2.4}
Let $\Model(\mathbb P^1_{\mathcal K})$ be the set of isomorphism classes of normal proper (equivalently projective) flat models of $\mathbb P^1_{\mathcal K}$. Let $\mathcal P_f^0\bigl(\Div(\mathcal K(t)/\Spec D)\bigr)$ be the set of non-empty finite sets of divisorial valuations of $\mathcal K(t)/\Spec D$ that dominate $D$. Then there exists a bijection
$$\Psi_{\mathcal K}:\Model(\mathbb P^1_{\mathcal K})\rightarrow \mathcal P_f^0\bigl(\Div(\mathcal K(t)/\Spec D\bigr)$$ 
defined by the rule $\Psi_{\mathcal K}([X])=\Div_X$, where for a normal proper flat model $X$ of $X_{\mathcal K}=\mathbb P^1_{\mathcal K}$ we consider the set $\mathcal B_X$ of $1$-dimensional irreducible components of $X_{\kappa}$ and we associate to it the set of divisorial valuations
$$\Div_X:=\{\mathcal O_{X,\eta_C}|C\in\mathcal B_X,\eta_C\;\textup{is the generic point of C}\}.$$
\end{corollary}

\begin{proof}
Let $X'$ be a normal proper flat model of $X_{\mathcal K}=\mathbb P^1_{\mathcal K}$ with $\Div_{X'}=\Div_X$. Let $X^+$ be the normalization of the Zariski closure in $X\times_{\Spec D} X'$ of the diagonal embedding of $\mathbb P^1_{\mathcal K}\rightarrow (X\times_{\Spec D} X')\times_{\Spec D} \Spec\mathcal K=\mathbb P^1_{\mathcal K}\times_{\Spec\mathcal K} \mathbb P^1_{\mathcal K}$. The two projections $X^+\rightarrow X$ and $X^+\rightarrow X'$ are Stein factors of $X^+$. Thus, as we have $\Div_X=\Div_{X'}$, these two Stein factors of $X^+$ are isomorphic by Proposition \ref{L9a}(1). Hence there exists an isomorphism $X\rightarrow X'$ that extends $1_{\mathbb P^1_{\mathcal K}}$. Thus $\Psi_{\mathcal K}$ is injective.

Given $\Div\in\mathcal P_f^0\bigl(\Div(\mathcal K(t)\Spec D\bigr)$, let $Y$ be an affine normal flat scheme of finite type over $\Spec D$ such that $Y\times_{\Spec D}\Spec\mathcal K$ is an open subvariety of $\mathbb P^1_{\mathcal K}$ and the set of local rings of $Y$ at generic points of irreducible components of its special fibre $Y_{\kappa}$ is $\Div$. Let $Y^+$ be the quasi-projective variety obtained by gluing $Y$ and $\mathbb P^1_{\mathcal Y}$ and let $X$ be a projective normal scheme over $\Spec D$ that contains $Y^+$ as an open Zariski dense subscheme; so $X$ is a normal proper flat model of $\mathbb P^1_{\mathcal K}$ with $\Div\subset\Div_X$. Let $f:X\rightarrow V$ be a Stein factor such that $V$ is a normal projective flat model of $\mathbb P^1_K$ with $\Psi_{\mathcal K}(V)=\Div_V=\Div$ by Proposition \ref{L9a}(2.a) to (2.c). Thus $\Psi_{\mathcal K}$ is surjective.\end{proof}

\begin{example}\normalfont\label{EX5-}
We take $K$ to be an algebraic closure of $\mathcal K$. Let $\widehat{D}$ be the completion of $D$. Let $\widehat{\mathcal K}:=\Frac(\widehat{D})$. Let $\overline{\widehat{\mathcal K}}$ be an algebraic closure of $\widehat{\mathcal K}$. We fix a divisorial valuation $O$ of $\mathcal K(t)$. Let $\mathfrak w:\mathcal K(t)\rightarrow \Gamma_{\mathcal K(t)}\cup\{\infty\}$ be the discrete valuation that defines it; so $\mathfrak w$ is surjective, its valuation ring is $O$, and we have an isomorphism $\Gamma_{\mathcal K(t)}\cong \mathbb Z$ of totally ordered groups. Let $\overline{\mathfrak w}: K(t)\rightarrow \Gamma_{K(t)}\cup\{\infty\}$ be a valuation that extends $\mathfrak w$; so $\overline{\mathfrak w}$ is surjective and $\Gamma_{K(t)}$ is a totally ordered group. Let $O_{\overline{\mathfrak w}}$ be the valuation ring of $\overline{\mathfrak w}$. Let $\overline{D}:=O_{\overline{\mathfrak w}}\cap K$ and let $\overline{\mathfrak v}:K\rightarrow\Gamma_K\cup\{\infty\}$ and $\mathfrak v:\mathcal K\rightarrow\Gamma_{\mathcal K}\cup\{\infty\}$ be induced by $\overline{\mathfrak w}$; so $\overline{\mathfrak v}$ and $\mathfrak v$ are surjective and $\Gamma_K$ and $\Gamma_{\mathcal K}$ are totally ordered subgroups of $\Gamma_{K(t)}$. Let $\overline{D}$ and $D$ be the valuation rings of $\overline{\mathfrak v}$ and $\mathfrak v$ (respectively). We denote also by $\overline{\mathfrak v}$ the valuation of the completion  $\widehat{K}$ of $K$ with respect to $\overline{\mathfrak v}$. We naturally identify $\overline{\widehat{\mathcal K}}$ with a subfield of $\widehat{K}$. We have 
$$\Gamma_{\mathcal K}\leqslant\Gamma_{\mathcal K(t)}\leqslant\Gamma_{K(t)}=\Gamma_K=\Gamma_{\mathcal K}\otimes_{\mathbb Z} \mathbb Q=\Gamma_{\mathcal K(t)}\otimes_{\mathbb Z} \mathbb Q$$ 
(see \cite{Epp}, Thm.\ (2.0) and \cite{Kuh} for the first equality).

Let $k$ be the residue field of $\overline{D}$; it is an algebraic closure of $\kappa$. The residue field $k(\overline{\mathfrak w})$ of $O_{\overline{\mathfrak w}}$ has transcendental degree $1$ over $k$, i.e., $\overline{\mathfrak w}$ is a residual transcendental extension of $\overline{\mathfrak v}$ in the sense of \cite{APZ1}, Sect.\ 1; such extensions were classified in \cite{AP}, \cite{APZ1}, and \cite{APZ2}. In particular, there exists $(\alpha_0,r)\in K\times \overline{\mathfrak v}(K^{\ast})$ such that 
$$\mathfrak w\Bigl(\sum_{i=0}^s \delta_i(t-\alpha_0)^i\Bigr)=\inf\{\overline{\mathfrak v}(\delta_i)+ir|i\in\llbracket0,s\rrbracket\}\in\Gamma_{K(t)}\cup\{\infty\}$$ for all $(s,\delta_0,\ldots,\delta_s)\in\mathbb N\times K^{s+1}$ (see \cite{AP}, Prop.\ 2), equivalently the norm on $K(t)$ is the $\sup$ norm on the disk $\mathcal B_K(\alpha_0,r):=\{x\in K|\overline{\mathfrak v}(x-\alpha_0)\ge r\}$ or on the (rigid analytic) disk $\mathcal B_{\widehat{K}}(\alpha_0,r):=\{x\in \widehat{K}|\overline{\mathfrak v}(x-\alpha_0)\ge r\}$.

The divisorial valuations of $\mathcal K(t)$ are in bijection to those of $\mathcal K^{\h}(t)$, which correspond to $\Gal(K/\mathcal K^{\h})$-orbits of the action of $\Gal(K/\mathcal K^{\h})$ on the set of such disks $\mathcal B_K(\alpha_0,r)$ (see \cite{APZ2}, Thm.\ 2.2; cf.\ \cite{AP}, Prop.\ 3). 

Let $l\in\mathbb N$ be the smallest such that there exists $\beta_0\in\overline{\widehat{\mathcal K}}$ of degree $l$ over $\widehat{\mathcal K}$ with $\mathcal B_{\widehat{K}}(\beta_0,r)=\mathcal B_{\widehat{K}}(\alpha_0,r)$. We approximate $\beta_0$ by separable elements (if $\chr(\mathcal K)=p$, approximate a solution of $X^{p^q}=\gamma$ by a solutions of $X^{p^q}+\epsilon X=\gamma$ with $\epsilon$ near $0$). Then by Krasner's Lemma we can assume that $\beta_0\in K$ is separable over $\mathcal K$ of degree $l$; this was first proved in \cite{APZ2}, Thm.\ 3.1. 

If $\beta\in K\cap\mathcal B_{\widehat{K}}(\alpha_0,r)$ is algebraic over $\mathcal K$ of degree $l$, then there exists a unique extension of $\mathfrak v$ to $\mathcal K(\beta)$, i.e., the local degree (i.e., the degree after the completion) is still $l$; hence the valuation ring of $\mathcal K(\beta)$ is a finite $D$-module.
\end{example} 

\begin{theorem}\label{T11}
We consider a divisorial valuation $O$ of $\mathcal K(t)/\Spec D$ such that $\{O\}\in\mathcal P_f^0\bigl(\Div(\mathcal K)/\Spec D\bigr)$. Let $X=X_O$ be a flat normal model of $X_{\mathcal K}=\mathbb P^1_{\mathcal K}$ such that $\Psi_{\mathcal K}([X])=\{O\}$ (see Corollary \ref{C2.4}). Then the following properties hold.

\medskip
{\bf (1)} If $Q$ is the closure of $(1:0)\in\mathbb P^1_{\mathcal K}(\mathcal K)$, then $X\setminus Q$ is affine and it is $\Spec (\mathcal K[t]\cap O)$. 

\smallskip
{\bf (2)} The flat model $\Spec (\mathcal K[t]\cap O)$ of $\mathbb A^1_{\mathcal K}$ has at most one singular point.

\smallskip
{\bf (3)} The flat normal model $X$ of $\mathbb P^1_{\mathcal K}$ has at most two singular points.
\end{theorem}

\begin{proof} Part (1) follows from Proposition \ref{L9a}(2.d). 

For parts (2) and (3), we use the notation of Example \ref{EX5-}, with the finite field extension $\mathcal K\rightarrow \mathcal K(\alpha_0)$ separable of degree $l$. Let $L\in\mathcal K[t]$ be the monic minimal polynomial of $\alpha_0$ over $\mathcal K$. Let $\mathcal K':=\mathcal K(\alpha_0)$; its valuation ring $D'$ is a free $D$-module of rank $l$. Let $\kappa'$ be the residue field of $D'$ and let $i_{D'}$ be the ramification index of the extension $D\rightarrow D'$. We have 
$$\Gamma_{\mathcal K}\subset \Gamma_{\mathcal K'}\subset \Gamma_{\mathcal K'}+\mathbb Zr\subset \Gamma_{\mathcal K}\otimes_{\mathbb Z} \mathbb Q$$ and inclusions $\mathcal K^{\h}\subset\mathcal K''\subset (\mathcal K')^{\h}$, where $\Gal(K/\mathcal K'')$ is the subgroup of $\Gal(K/\mathcal K^{\h})$ that fixes $\mathcal B_K(\alpha_0,r)$.\footnote{The field extension $\mathcal K''\rightarrow (\mathcal K')^{\h}$ is an extension with trivial tame part. In particular, if it is non-trivial, then $\chr(\kappa)=p$ and it has degree a power of $p$.}

For the disk $\mathcal B_{\Frac(\widehat{D'})}(\alpha_0,r)$ defined similarly, where $\widehat{D'}$ is the completion of $D'$, the valuation group is $\Gamma_{\mathcal K'}+\mathbb Zr$. We have $\overline{\mathfrak w}(L)\in \Gamma_{\mathcal K'}+\mathbb Zr$; let $m$ be the order of $\mathfrak{w}(L)$ modulo $\Gamma_{\mathcal K'}$. Thus $\overline{\mathfrak w}(L^m)=\overline{\mathfrak v}(\delta)$ for some $\delta\in (\mathcal K')^{\ast}$. We have a $\mathcal K$-linear natural isomorphism $\oplus_{j=0}^{l-1}\mathcal Kt^j\cong \mathcal K'$. If $0\neq\gamma\in\mathcal K'$ is represented by $g_{\gamma}(t)\in\oplus_{j=0}^{l-1}\mathcal Kt^j$, then $g_{\gamma}$ has no zeros in the disk so it must be a constant $\gamma_0\in\mathcal K$ plus something with greater value of $\overline{\mathfrak v}$. Thus $\overline{\mathfrak w}(\frac{L^m}{g_{\delta}})=0$. We have to compare the disk $\mathcal B_{\mathcal K(t)}(\alpha_0,r)$ and its ring of polynomial functions $\mathcal K[t]\cap O$ to the disk $\mathcal B_{\widehat{\mathcal K'}}(0,\overline{\mathfrak w}(L))$. Note that $\frac{L^m}{g_{\delta}}\in O_{\overline{\mathfrak w}}$ defines an element $\zeta\in k(\overline{\mathfrak w})$ which is transcendental over $\kappa'$ (see \cite{APZ1}, Thm.\ 2.1b)).

Let $O'$ be the discrete valuation ring of the `standard' valuation $\mathfrak w'$ of $\mathcal K'(t)$ with the properties: (i) the restrictions of $\mathfrak w'$ and $\overline{\mathfrak w}$ to $\mathcal K'$ coincide, and (ii) $\mathfrak w'(t)=\overline{\mathfrak w}(L)$.

For $h\in\mathcal K[t]$ we consider its $L$-expansion $h=\sum_{j=0}^{n_h} h_jL^j$; so $n_h\in\mathbb N$ is the smallest such that $\deg(h_j)<l$ for each $j\in\llbracket0,n_h\rrbracket$. By sending each $h_j$ to its evaluation $h_j(\alpha_0)\in\mathcal K'$ and $L$ to $t$ we get a $\mathcal K$-linear isomorphism $\mathcal K[t]\rightarrow\mathcal K'[t]$, which preserves the $\sup$ norm on the relevant disks. It is not a $\mathcal K$-algebra isomorphism in general but we have equality of products modulo elements with smaller $\sup$ norm. Thus, this gives an isomorphism of associated graded rings
$$\gr(\mathcal K[t]\cap O)\rightarrow\gr(\mathcal K'[t]\cap O'),$$ the filtrations being $\bigl(\mathcal (K[t]L^j)\cap O\bigr)_{j\in\mathbb N}$ and $\bigl(\mathcal (K'[t]t^j)\cap O'\bigr)_{j\in\mathbb N}$. This implies that the residue field of $O'$ is $\kappa'(\zeta)$ and that the value group $\Gamma_{K'(t)}$ is $\Gamma_{\mathcal K'}+\mathbb Z\overline{\mathfrak w}(L)$.

We consider two cases as follows.

{\bf Case 1: $\overline{\mathfrak w}\in\Gamma_{\mathcal K'}$ (i.e., $m=1$).} Then $O'\cap\mathcal K'[t]\cong D'[t/\delta]$ is regular. From this and the information on graded rings we get that $\mathcal K[t]\cap O$ is regular by \cite{Ma}, Ch.\ 6, Thm.\ 17.10.

{\bf Case 2: $\overline{\mathfrak w}\notin\Gamma_{\mathcal K'}$ (i.e., $m>1$).} Then we can still verify the regularity after inverting $\frac{L^m}{g_{\delta}}$. This leaves one possible singular point in the special fiber $\Spec\kappa'[\zeta]$ of $\Spec (\mathcal K[t]\cap O)$ which is the zero locus $\zeta=0$. 

Part (2) follows from the two cases above. 

Part (3) follows from part (1) and the union $X=Q\cup [\Spec (\mathcal K[t]\cap O)]$.\end{proof}

\begin{remark}\normalfont\label{R5} We refer to the proof of Theorem \ref{T11}. Let $C:=(X_{\kappa})_{\red}$. Let $(C\cap Q)_{\red}$ define the closed point $P$ of $C$ at infinity.

\medskip
{\bf (1)} The special fiber $C\setminus\{P\}$ of $\Spec (\mathcal K[t]\cap O)$ is $\Spec\kappa'[\zeta]\cong\mathbb A^1_{\kappa'}$.

\smallskip
{\bf (2)} In Case 2, the closed point $\zeta=0$ of $\Spec\kappa'[\zeta]$ lifts to a $D'$-valued point of $\Spec (\mathcal K[t]\cap O)$ given by evaluation at $\alpha_0$; in fact we have a $D$-algebra surjection $\mathcal K[t]\cap O\rightarrow D'$ as a uniformizer $v'$ of $D'$ lifts to an element given by a polynomial $h_{v'}\in\mathcal K[t]$ of degree $<l$. It follows that the ideal of the reduced special fiber is not principal at this point (as otherwise it would be generated by $h_{v'}$ which does not have the right valuation at the generic point of the special fiber), so the point is singular. 

\smallskip
{\bf (3)} From (2) and Case 1 we get that $\Spec (\mathcal K[t]\cap O)$ is singular iff $m>1$.

\smallskip
{\bf (4)} The reduced special fiber $C$ of $X=X_O$ is the pinching pushout 
\[\xymatrix{
\Spec\kappa' \ar[r] \ar[d] & \Spec\kappa \ar[d] \\
\mathbb P^1_{\kappa} \ar[r] & C.
}\]

{\bf (5)} The ideal of $C$ at $P$ is not principal if $i_{D'}m>1$ and so we have a singularity. It is principal generated by a uniformizer $v$ of $D$ if $i_{D'}=m=1$, but if $\kappa'\neq\kappa$ we still have a singularity as the zero locus $v=0$ is singular and $v$ is not in the square of the maximal ideal. 

\smallskip
{\bf (6)} Suppose that the field $\kappa$ is algebraically closed. Then $\kappa'=\kappa$ and from remarks (3) and (5) we get that $|\Sing(X)|=2$ iff $m>2$ and iff $|\Spec (\mathcal K[t]\cap O)|=1$ and that $|\Sing(X)|=1$ iff $i_{D'}\ge 2$.

\smallskip
{\bf (7)} If $\mathcal K'=\mathcal K$ and $m=1$ in the proof of Theorem \ref{T11}, then $X=\mathbb P^1_D$. 
\end{remark}

\begin{corollary}\label{C2.5}
Let $R:=k[[x]]$. Let $R\rightarrow S$ be an extension of discrete valuation rings which is generically finite Galois of Galois group $G$. If the action of $G$ on $S$ extends to an action of $G$ on $S[t]$, then $S[t]^G$ has at most one singular point.
\end{corollary}

\begin{proof}
We take $D:=R$, $\mathcal K=k((x))$ and $O:=\mathcal O_{\Spec S[t]^G,\eta}$, where $\eta$ is the generic point of $\Spec (S[t]^G/S[t]^Gx)_{\red}$. This makes sense as $\Spec S[t]^G[\frac{1}{x}]$, being a form of $\mathbb A^1_{k((x))}$, is isomorphic to $\mathbb A^1_{k((x))}$ and thus we identify it with $\Spec\bigl(k((x))[t]\bigr)$. With the notation of Theorem \ref{T11} and its proof, we have $\kappa=\kappa'=k$ and $\mathcal K[t]\cap O=S[t]^G$. So the corollary holds by Theorem \ref{T11}(2).\end{proof}

\section{On rational singularities}\label{S9}

For (quasi-)excellent rings and schemes we refer to \cite{Gro3}, Subsect.\ 7.8 and \cite{T}, Exp.\ 1. We recall the following general result.

\begin{lemma}\label{L9+} Let $W$ be a noetherian scheme. The following two properties hold.

\medskip
{\bf (1)} If $W$ is normal, then it is quasi-excellent iff it is excellent. 

\smallskip
{\bf (2)} Let $f:V\rightarrow W$ be a finite surjective morphism. Then $V$ is quasi-excellent iff $W$ is quasi-excellent.
\end{lemma}

\begin{proof}
The `if' of part (1) is clear. For the `only if' of part (1) it suffices show that a quasi-excellent normal noetherian ring $R$ is universally catenary. To check this we can assume that $R$ is local; its formal fibers are geometrically regular and thus also geometrically normal and from \cite{Stacks5} we get that $R$ is universally catenary. 

To prove part (2) we can assume that $W$ is affine, hence $V$ is also affine, and this case was proved in \cite{Gre}, Thm.\ 3.1i).\end{proof}

For rational singularities see \cite{Lip1}, Def.\ 1.1.  Part (2) of the following result in characteristic $0$ is a particular case of either \cite{Bout}, Thm.\ or \cite{Kem}, Thm.\ 7. 

\begin{proposition}\label{PR25}
Let $G$ be a finite group. Let $W$ be a noetherian scheme over $\Spec\mathbb Z[\frac{1}{|G|}]$. We consider an action of $G$ on $W$ with the property that for each $x\in W$ the orbit $G\cdot x$ is contained in an affine open subscheme of $W$. Then the following properties hold.

\medskip
{\bf (1)} The universal geometric quotient $W\rightarrow W/G$ in the sense of \cite{MFK}, Ch.\ 0, Defs.\ 0.6 and 0.7 exists and it is an affine morphism. 

\smallskip
{\bf (2)} If $\dim(W)=2$ and $W$ is normal, excellent and has at most rational singularities, then $X/G$ has at most rational singularities.
\end{proposition}

\begin{proof}
We can assume that $W=\Spec R$ is affine. Taking $W/G:=\Spec R^G$, clearly $W\rightarrow W/G$ is surjective and $G$-invariant. If $W$ is normal, so is $W/G$. The surjective morphism $W\rightarrow W/G$ is finite by \cite{T}, Exp.\ IV, Prop.\ 2.2.3(ii)(a). 

The commutative $\mathbb Z$-algebra $R^G$ is noetherian due to the existence of the Raynolds operator $R\rightarrow R^G$, i.e., of the $R^G$-linear retraction that maps $x\in R$ to the average sum $\frac{1}{|G|}\sum_{g\in G} gx\in R^G$ (e.g., see \cite{T}, Exp.\ IV, Prop.\ 2.2.3(i)(a)). The Reynolds operator implies that for each $R^G$-algebra $S$ we have $(R\otimes_{R^G} S)^G=S$. If $S$ is an algebraically closed field, then $G$ acts transitively on $\Spec R\otimes_{R^G} S$ as the restriction function $\Hom(R,S)/G\rightarrow\Hom(R^G,S)$ is a bijection by \cite{Bour}, Ch.\ V, Subsect.\ 2.1, Cor.\ 4 after Thm.\ 1 and Subsect.\ 2.2, Cor.\ after Thm.\ 2. So $W\rightarrow W/G$ is a universal geometric quotient and part (1) holds.

For part (2), $R^G$ is excellent by Lemma \ref{L9+}(1) and (2). Being also normal, we can strictly henselize by \cite{Lip1}, Prop.\ 16.5 and \cite{Gre}, Cor.\ 5.6, i.e., we can assume that $W=\Spec R$ and $W/G=\Spec R^G$ are affine normal strictly henselian local of dimension $2$ and that $G$ acts trivially on the residue field of $R$. 

\phantomsection{Let $\Sigma:V\rightarrow W/G$ be a resolution of singularities by \cite{Lip2}, Thm. We can desingularize $Z:=(V\times_{W/G} W)^\prime$ (the $\prime$ means that we take the schematic closure of the punctured spectrum of $R$, equivalently the reduced scheme) by normalizations and blowing up of closed points. This is $G$-equivariant. We consider a $G$-equivariant desingularization $Z_1\rightarrow W$ that maps to $V$. 

For the resulting morphism $\phi:Z_1/G\rightarrow V$, the homomorphism of ringed sheaves $\mathcal O_V\rightarrow \phi_*(\mathcal O_{Z_1/G})$ is an isomorphism. Then a Čech computation gives $H^1(Z_1/G,\mathcal O_{Z_1/G})=H^1(Z_1,\mathcal O_{Z_1})^G=0$ and $H^1(V,\mathcal O_V)\to H^1(Z_1/G,\mathcal O_{Z_1/G})$ is injective by the Leray spectral sequence. Thus $H^1(V,\mathcal O_V)=0$, i.e., $R/G$ has a rational singularity.}\label{PH50}\end{proof}

\begin{proposition}\label{PR25+}
Let $D$ be a discrete valuation ring. Let $R$ be a $D$-subalgebra of $\Frac(D)[t]$ which is of finite type over $D$. Then the normalization $R^{\n}$ of $R$ has only rational singularities.
\end{proposition}

\begin{proof}
We can assume that $D\neq R\neq\Frac(D)$; so $\Spec R^{\n}\otimes_D \Frac(D)\cong\mathbb A^1_{\Frac(D)}$. Hence the $R$-module $R^{\n}$ is finite by \cite{GLL}, Prop.\ 7.6(a). We compactify $\Spec R^{\n}$ to a normal integral projective scheme $X$ over $\Spec D$ with general fiber $\mathbb P^1_{\Frac(D)}$. The rational map $\mathbb P^1_D\dashrightarrow X$ becomes regular  $\begin{tikzcd}[row sep=tiny]
W \arrow[rd] \arrow[rr, "\phi"] & & \mathbb P^1_D \arrow[ld, dashrightarrow] \\
& X
\end{tikzcd}$ after blowing up finitely many points.  As $X\rightarrow\Spec D$ is cohomologically flat by (the implication $(ii)\Rightarrow (iv)$ in) \cite{Ray}, Thm.\ 8.2.1, we have $H^1(X,\mathcal O_X)=0$. As the blowing up does not change the cohomology of the structure sheaf, $H^1(W,\mathcal O_W)=0$ and $R^1 \phi_*(\mathcal O_W)$ is a skyscrapper sheaf supported at finitely many closed points. To show a singular point $P\in X$ is rational we have to show that $\bigl(R^1 \phi_*(\mathcal O_W)\bigr)_P=0$. As $H^0\bigl(X,R^1 \phi_*(\mathcal O_W)\bigr)$ surjects to this, it suffices to show that $H^0\bigl(X,R^1 \phi_*(\mathcal O_W)\bigr)=0$. The differential $H^0\bigl(X,R^1 \phi_*(\mathcal O_W)\bigr)\rightarrow H^2\bigl(X,\phi_*(\mathcal O_W)\bigr)$ vanishes by a general cohomological dimension fact on proper schemes in terms of the dimension of fibers. Thus $H^0\bigl(X,R^1 \phi_*(\mathcal O_W)\bigr)=0$ as $H^1(W,\mathcal O_W)=0$, hence $\bigl(R^1 \phi_*(\mathcal O_W)\bigr)_P=0$.\end{proof}

\section{Geometric preliminaries}\label{S10}

For $(n,m)\in (\mathbb N^{\ast})^2$ with $n\ge m$, a smooth morphism $\mathbb A^n_R\rightarrow\mathbb A^m_R$ is called a projection if it is isomorphic to the one given on valued points by the rule $(x_1,\ldots,x_n)\mapsto (x_1,\ldots,x_m)$. We list the geometric notation for $e\in\EE_n(K)$.

\begin{notation}\normalfont\label{NOT3}
For $e\in\EE_n(K)$, let the Jacobian variety $\psi_e:X_e\rightarrow X_e$ and the open embedding $\imath_e:\mathbb A^n_{K,\s}\rightarrow X_e$ be as in Diagram (\ref{EQ00}) of Section \ref{S1} and let the finite subset $\mathcal F(e)$ of $\mathbb N$ and the numerical invariants $\rho(e)$, $\rho_{\et}(e)$, and $\deg(e)$ be as in Section \ref{S2}. We write $X_e=\Spec A_e$. Let $N_e:=\psi_e(\Gamma_{\psi_e})$ be the non-proper locus of $e$. Let $\nu(e):=|\Irr(N_e)|\in\mathbb N$. Let $X_e^{\et}$ be the open subvariety of $X_e$ which is the \'etale locus of $\psi_e$. Let $N_e^{\textup{n}-\et}:=\psi_e(X_e\setminus X_e^{\et})$. Let $\nu_{\textup{n}-\et}(e):=|\Irr(N_e^{\textup{n}-\et})|\in\mathbb N$.
\end{notation}

We have $N_e=\emptyset$ iff $\Gamma_{\psi_e}=\emptyset$ iff $\imath_e$ is an isomorphism (i.e., $X_e=\mathbb A^n_{K,s}$) and iff $e$ is finite. If $N_e\neq\emptyset$, then as $\psi_e$ is quasi-finite and $\Gamma_{\psi_e}=X_e\setminus\Imm(\imath_e)$ is of pure codimension $1$ (see \cite{Gro6}, Exp.\ V, Ex.\ 3.4), $N_e=\psi_e(\Gamma_{\psi_e})$ is equidimensional of dimension $n-1$. In general, the inequalities 
$$\nu(e)\le\rho(e)\le\nu(e)[\deg(e)-1]$$ 
hold and $N_e$ is the smallest (hence reduced) closed subscheme of $\mathbb A^n_{K,\t}$ such that $e$ is finite \'etale over $\mathbb A^n_{K,\t}\setminus N_e$. Similarly, $N_e^{\textup{n}-\et}$ is equidimensional of dimension $n-1$, the inequalities 
$$\nu(e)_{\textup{n}-\et}\le\rho(e)-\rho_{\et}(e)\le\nu_{\textup{n}-\et}(e)[\deg(e)-1]$$ 
hold and $N_e^{\textup{n}-\et}$ is the smallest closed subvariety of $\mathbb A^n_{K,\t}$ such that $\psi_e$ is finite \'etale over $\mathbb A^n_{K,\t}\setminus N_e^{\textup{n}-\et}$. Moreover, as $X_e\setminus X_e^{\et}\subset X_e\setminus\Imm(\imath_e)$, we have $N_e^{\textup{n}-\et}\subset N_e$ and $\nu_{\textup{n}-\et}(e)\le\nu(e)$. 

\begin{remark}\normalfont\label{R6}
Let $e\in\QFE_n(K)$. For $l\in\mathcal F(e)\cup\{\deg(e)\}$ we consider the reduced locally closed subscheme $N_l(e)$ of $\mathbb A^n_{K,\t}$ defined by 
$$N_l(e)(K):=\{P\in\mathbb A_{K,\t}^n(K)||e^{-1}(P)|=l\}.$$
Assume now that $e\in\EE_n(K)$; we have $N_e=\cup_{l\in\mathcal F(e)} N_l(e)$. As $N_0(e)$ is either empty or of codimension at least $2$ (see Lemma \ref{L6}(2)) while $N_e$ is equidimensional of dimension $n-1$, it follows that if $e$ is non-finite (i.e., if $\mathcal F(e)\neq\emptyset$), then $\mathcal F(e)\cap\mathbb N^{\ast}$ is non-empty and hence 
$$\varsigma_e=\min(\mathcal F(e)\cap\mathbb N^{\ast})\in\mathbb N^{\ast}$$ 
(see Notation \ref{NOT1} for $\varsigma_e)$. The invariant 
$$\upsilon(e):=\max\bigl(\mathcal F(e)\cap\mathbb N^{\ast}\bigr)\in\mathbb N^{\ast}$$ 
helps in the study of $N_e$.
\end{remark}

\begin{example}\normalfont\label{EX5}
Let $e\in\EE_n(K)$ be non-finite. For $l\in\mathcal F(e)$, as the Zariski closure of $N_l(e)$ in $\mathbb A^n_{K,\t}$ (equivalently, in $N_e$) is contained in $\cup_{m\in\mathcal F(e),m\le l} N_m(e)$ by \cite{Gro5}, Prop.\ (18.2.8), we have $\Irr\bigl(N_{\upsilon(e)}(e)\bigr)\subset\Irr(N_e)$, hence $N_{\upsilon(e)}(e)$ is equidimensional of dimension $n-1$. If $\upsilon(e)=\deg(e)-1$, then $\psi_e$ is \'etale above the generic points of irreducible components of $N_{\upsilon(e)}(e)$ and thus 
$$\rho_{\et}(e)\ge |\Irr\bigl(N_{\upsilon(e)}(e)\bigr)|\ge 1;$$ 
if moreover, $\mathcal F(e)=\{\deg(e)-1\}$, then $X_e$ is a finite \'etale cover of $\mathbb A^n_{K,\t}$ by the purity of the branch locus (see \cite{Gro6}, Exp.\ X, Thm.\ 3.4 (i)) and therefore we have $\nu(e)=\rho_{\et}(e)=\rho(e)>0$.\end{example}

\begin{example}\normalfont\label{EX5+}
Let $e\in\EE_n(K)$. We define $W_e^{0+}:=\mathbb A^n_{K,\s}\times_{e,\mathbb A^n_{k,\t},\psi_e} X_e$. The second projection $W_e^{0+}\rightarrow X_e$ is \'etale; moreover, it is surjective iff $e$ is so. The first projection $f^{0+}_e:W_e^{0+}\rightarrow\mathbb A^n_{K,\s}$ is a finite surjective generically \'etale morphism between normal schemes which is flat in codimension at most $2$ and has a section $1_{\mathbb A^n_{K,\s}}\times\imath_e:\mathbb A^n_{K,\s}\rightarrow W_e^+$. Let $W_e^0:=W_e^{0+}\setminus\Imm(1_{\mathbb A^n_{K,\s}}\times\imath_e)$; it is a clopen subscheme of $W_e^{0+}$ which is empty iff $\deg(e)=1$ and iff $e\in\GA_n(k)$ by Lemma \ref{L6}(4). 

Suppose that $e\in\EE_n(K)\setminus\GA_n(K)$. Then $W_e^0$ is equipped with a finite surjective morphism $f_e^0:W_e^0\rightarrow\mathbb A^n_{K,\s}$
which is generically \'etale of degree $\deg(e)-1$ and flat in codimension at most $2$. Note that $\psi_e$ is regular (resp.\ \'etale) above $\Imm(e)$ iff $W_e^0$ is regular (resp.\ iff $f_e$ is \'etale). Let $W_e$ be the normalization of $X_e$ in the ring of fractions of $W_e^0$. So $W_e$ is a normal affine variety over $K$ which contains $W_e^0$ as an open subvariety and is equipped with a finite surjective morphism
$$f_e:W_e\rightarrow X_e$$
that extends $f_e^0$.\end{example}

Recall the following well-known lemma.

\begin{lemma}\label{L10}
Let $f:X\rightarrow Y$ be a dominant morphism of irreducible separated curves over $K$. For $W\in\{X,Y\}$, let $W^{\n}$ be the normalization of $W$, and let $W^{\c}$ be the complete smooth connected curve having $W^{\n}$ as an open subvariety. Then the following properties hold.

\medskip
{\bf (1)} If $f$ is surjective, then $|X^{\c}\setminus X^{\n}|\ge |Y^{\c}\setminus Y^{\n}|$. 

\smallskip
{\bf (2)} If $X$ is rational, then $Y$ is rational.

{\bf (3)} If $|X^{\c}\setminus X^{\n}|=|Y^{\c}\setminus Y^{\n}|$, then $f$ is finite surjective.
\end{lemma}

\begin{proof}
We can assume that $X$ and $Y$ are normal. The morphism $f$ extends to a finite surjective morphism $f^{\c}:X^{\c}\rightarrow Y^{\c}$. Hence $f^{\c}(X^{\c}\setminus X)$ has at least $|Y^{\c}\setminus Y|$ elements and part (1) holds. 

Part (2) follows from Lüroth's Theorem. 

For part (3), the inequalities  $|Y^{\c}\setminus Y|\le |(f^{\c})^{-1}(Y^{\c}\setminus Y)|\le |X^{\c}\setminus X|$ must be equalities, so the inclusion $(f^{\c})^{-1}(Y^{\c}\setminus Y)\subset X^{\c}\setminus X$ is a bijection. Therefore $(f^{\c})^{-1}(Y)=X$ and the finite surjective morphism $(f^{\c})^{-1}(Y)\rightarrow Y$ is $f$.\end{proof}

\begin{example}\normalfont\label{EX5.5}
{\bf (1)} The hypersurface 
$$X_1:=\Spec\bigl(k[x,y,z]/(x^{p+1}+xy+z^n)\bigr)$$ 
in $\mathbb A^3_k$ is regular for $n=1$ and normal with $\Sing(X_1)=\{(0,0,0)\}$ if $n\ge 2$. The non-\'etale locus $Y$ of the finite flat morphism $\phi_1:X_1\rightarrow\mathbb A^2_k$ defined by the projection on the last two coordinates is the zero locus $x^p+y=0$. The induced morphism $Y\to\psi_1(Y)$ is radicial and $\psi_1^{-1}\bigl(\psi_1(Y)\bigr)$ has two irreducible components Y and $W$ with $(Y\cap W)_{\red}=\{(0,0,0)\}$. If $n=1$, then $X_1\cong\mathbb A^2_k$. If $n\ge 2$, then the rule $(x,y)\mapsto (x,-x^{n-1}y^n-x^p,xy)$ defines a birational morphism $\mathbb A^2_k\rightarrow X_1$ which is not an open embedding. 

\smallskip
{\bf (2)} If $p=2$ let $l=3$ and if $p\ge 3$ let $l=2$. For the normal hypersurface 
$$X_2:=\Spec\bigl(k[x,y,z]/(x^{pl}+xy-z^l)\bigr)$$ 
in $\mathbb A^3_k$ we have $\Sing(X)=\{(0,0,0)\}$. The non-\'etale locus of the finite flat morphism $\phi_2:X_2\rightarrow\mathbb A^2_k$ defined by the projection on the last two coordinates is the zero locus $y=0$ and has two irreducible components $Y_1$ and $Y_2$ with $(Y_1\cap Y_2)_{\red}=\Sing(X)$ and with $Y_1$ as the zero locus $x^p-z=0$. The induced morphism $Y_1\rightarrow \phi_2(Y_1)$ is radicial and the rule $(x,y)\mapsto (x,x^{l-1}y^l-x^{pl-1},xy)$ defines a morphism $\mathbb A^2_k\rightarrow X_2$ which is birational but not an open embedding.

\smallskip
{\bf (3)} The hypersurface 
$$X_3:=\Spec\bigl(k[x,y,z]/(x^p+xy+yz)\bigr)$$ 
in $\mathbb A^3_k$ is normal with $\Sing(X_3)=\{(0,0,0)\}$. The non-\'etale locus $W$ of the finite flat morphism $\phi_3:X_3\rightarrow\mathbb A^2_k$ defined by the projection on the last two coordinates is the zero locus $y=0$, thus the induced morphism $W\rightarrow \phi_3(W)$ is an isomorphism.

\smallskip
{\bf (4)} Let $l\in\mathbb N^{\ast}\setminus\{1\}$. Let $A=K\oplus t^lK[t]$ be the $K$-subalgebra of $K[t]$ generated by $\{t^i|i\in\llbracket l,2l-1\rrbracket\}$. For the affine curve $C:=\Spec A$ we have $\edim(C)=l$. The morphism $C_1:=\Spec \bigl(A[x]/(x^l-t^l)\bigr)\rightarrow C$ is finite flat of degree $l$ and $C_1$ has a closed subscheme isomorphic to $\mathbb A^1_K$.\end{example}

Example \ref{EX5.5}(1) or (2) implies that the result \cite{Je}, Thm.\ 2.1 over $\mathbb C$ does not hold in characteristic $p$ even when $n=2$.

Examples \ref{EX5.5}(1) to (4) imply that the following regularity criterion, which is only a variation of Abhyankar's Lemma, is the maximum one can get using local arguments.

\begin{criterion}\label{CRI1}
Suppose that $n\ge 2$. Let $e\in\EE_n(K)$. For $Z\in\Irr(X_e\setminus X_e^{\et})$ let $m_Z$ be the multiplicity of $Z$ in $\psi_e^{-1}\bigl(\psi_e(Z)\bigr)$. Let 
$$\Irr_{\textup{t}}(X_e\setminus X_e^{\et}):=\{Z\in \Irr(X_e\setminus X_e^{\et})|\chr(K)\nmid m_Z,\; Z\rightarrow\psi_e(Z)\;\textup{is generically \'etale}\}.$$ Let $Y\in \Irr_{\textup{t}}(X_e\setminus X_e^{\et})$. Then the following properties hold.

\medskip
{\bf (1)} Let $Q\in Y(K)$ be a point such that we have $\psi_e(Q)\in\Reg\bigl(\psi_e(Y)\bigr)$. We assume that for each $Z\in \Irr(X_e\setminus X_e^{\et})$ with $Z\notin\Irr_{\textup{t}}(X_e\setminus X_e^{\et})$ or $\psi_e(Z)\neq Y$ we have $Q\notin Z(K)$.\footnote{If $n=2$, the second condition always holds by Corollary \ref{C2.6}(1) below.} Then $Q\in\Reg(X_e)$.

\smallskip
{\bf (2)} If $\Irr(X_e\setminus X_e^{\et})=\Irr_{\textup{t}}(X_e\setminus X_e^{\et})$, then $$\Sing(X_e)\subset \bigl(\psi_e^{-1}(\cup_{W\in\Irr(N_e^{\textup{n}-\et})} \Sing(W))\bigr)\cap (X_e\setminus X_e^{\et}).$$
\end{criterion}

\begin{proof}
Let $S$ and $R$ be the completions of the local rings $\mathcal O_{X_e,Q}$ and $\mathcal O_{\mathbb A^n_{K,\t},\psi_e(Q)}$ (respectively); the $R$-algebra $S$ is integral non-\'etale and normal. 

For part (1), as $\psi_e(Q)\in\Reg\bigl(\psi_e(Y)\bigr)$, we can identify $R=K[[x_1,\ldots,x_n]]$ in such a way that $x_1=0$ defines the inverse image of $\psi_e(Y)$ to $\Spec R$. The hypotheses imply that the morphism $\Spec A\rightarrow\Spec R$ is tamely ramified in codimension $1$ and becomes \;etale after inverting $x_1$. From Abhyankar's Lemma for regular divisors (see \cite{Stacks2}) we get that $S\cong R[t]/(t^{m_Y}-x_1)$, hence $S$ is regular. So $Q\in\Reg(X_e)$ and part (1) holds. 

Part (2) follows from part (1). 

If $\psi_e(Y)\cong\mathbb A^1_K$, then $Y(K)\subset\Reg(X_e)$ by part (1), hence part (3) holds.\end{proof}

If $\chr(K)=0$ and $e\in\EE_n(K)$, then $\Gamma_e\subset\Sing(N_e)$ by \cite{Je}, Cor.\ 1.2. We have the following version of loc.\ cit.\ which holds over all $K$ and which is weaker (cf.\ Example \ref{EX5.5}) as it uses in essence local arguments.

\begin{proposition}\label{PR4}
Let $e\in\EE_n(K)$. Then the following properties hold.

\medskip
{\bf (1)} We have $\Gamma_e\subset\Sing(N_e)\cup [\Reg(N_e)\cap\psi_e(X_e\setminus X_e^{\et})]$.

\smallskip
{\bf (2)} If $\psi_e$ is \'etale, then $\Gamma_e\subset\Sing(N_e)$.
\end{proposition}

\begin{proof}
As $(\mathbb A^n_{K,\t}\setminus N_e)\subset \Imm(e)$, we have $\Gamma_e\subset N_e$. We show that the assumption that there exists $P\in\Gamma_e\cap\Reg(N_e)$ such that $\psi_e^{-1}(P)\subset X_e^{\et}$ leads to a contradiction. Let $Y\in\Irr(N_e)$ be the unique irreducible component such that $P\in Y$. Let $(l,s)\in (\mathbb N^{\ast})^2$ be such that $\Irr\bigl(\psi_e^{-1}(Y)\bigr)=\{Y_i|i\in\llbracket1,l\rrbracket\}\cup \{Z_j|j\in \llbracket1,s\rrbracket\}$ with the generic point of each $Y_i$ belonging to $\mathbb A^n_{K,\s}$ and every $Z_j\in\Irr(\Gamma_{\psi_e})$. For each $i\in\llbracket1,l\rrbracket$, as $Y_i$ maps onto $Y$, hence there exist $j_i\in \llbracket1,s\rrbracket$ and $P_{j_i}\in (Y_i\cap Z_{j_i})(K)$ such $\psi_e(P_{j_i})=P$. As $P_{j_i}\in\Sing\bigl(\psi_e^{-1}(Y)\bigr)\cap X_e^{\et}$ maps to $P\in\Reg(Y)$, we reached a contradiction; so part (1) holds.\footnote{If $P\in\Gamma_e\cap\Reg(N_e)\cap\psi_e(X_e\setminus X_e^{\et})$, then by reasons of multiplicities each $W\in\Irr\bigl(\psi_e^{-1}(Y)\bigr)$ passes through a unique point in $\psi_e^{-1}(P)$.} 

Clearly, $(1)\Rightarrow (2)$.\end{proof}

For $e\in\EE_2(k)$ with $\Gamma_e\subsetneq\Sing(N_e)$ (resp.\ $\Gamma_e=\Sing(N_e)$) see Example \ref{EX8} with $\deg(f)=1$ (resp.\ Example \ref{EX8} with $\deg(f)\ge 2$ and Example \ref{EX11}).

%The following definition aims to represent a beginning for the classification of curves in $\mathbb A^2_k$ isomorphic to $\mathbb A^1_k$.

%\begin{definition}\label{D5b}
%Let $\mathbb A^{1;2}(k)$ be the set of curves in $\mathbb A^2_k$ isomorphic to $\mathbb A^1_k$. 

%\medskip
%{\bf (1)} Let $\mathfrak R:\mathbb A^{1;2}(k)\rightarrow\mathbb N\cup\infty}$ be the rectifying invariant\index{rectifying invariant} defined as follows. If there exists $j\in\mathbb N$
%\end{definition}

\section{On sphere-like surfaces}\label{S11}

The next lemma is also well-known (for instance, see \cite{Mi3}, Sect.\ 2 for the case $\chr(K)=0$), but due to its importance in the paper we include a self-contained proof of it.

\begin{lemma}\label{L11}
Let $X$ be a connected affine normal surface over $K$ with an open embedding $\imath:\mathbb A^2_K\rightarrow X$. For a projection $\pi:\mathbb A^2_K\rightarrow\mathbb A^1_K$, let $\chi(\pi):\mathbb A^2_K\rightarrow\mathbb P^1_K$ be its composite with the open embedding $\mathbb A^1_K\subset\mathbb P^1_K$. Then the following properties hold.

\medskip
{\bf (1)} There exists a unique morphism $\chi(\pi)_X:X\rightarrow\mathbb P^1_K$ such that $\chi(\pi)_X\circ\imath=\chi(\pi)$, i.e., we have a commutative diagram

\[\xymatrix{
\mathbb A_K^2 \ar[r]^{\imath} \ar[d]^{\pi} \ar[dr]^{\chi(\pi)} & X \ar[d]^{\chi(\pi)_X} \\
\mathbb A_K^1\ar[r]^{\textup{open}} & \mathbb P_K^1.
}\]

{\bf (2)} For each $Y\in\Irr\bigl(X\setminus\Imm(\imath)\bigr)$, the restriction of $\chi(\pi)_X$ to $Y$ is constant.
\end{lemma}

\begin{proof}
We consider a second projection $\pi^{\perp}:\mathbb A^2_K\rightarrow\mathbb A^1_K$ such that the morphism $\pi\times\pi^{\perp}:\mathbb A^2_K\rightarrow \mathbb A^1_K\times_{\Spec K} \mathbb A^1_K$ is an isomorphism. Let 
$$\phi=\chi(\pi)\times\chi(\pi^{\perp}):\mathbb A^2_K\to\mathbb P^1_K\times_{\Spec K}\mathbb P^1_K;$$
equivalently, $\phi$ is the composite of $\pi\times\pi^{\perp}$ with the product open embedding $\mathbb A^1_K\times_{\Spec K}\mathbb A^1_K\subset\mathbb P^1_K\times_{\Spec K}\mathbb P^1_K$. It suffices to show that $\phi$ extends uniquely to a morphism $\phi_X:X\rightarrow\mathbb P^1_K\times_{\Spec K}\mathbb P^1_K$ and that for each $Y\in\Irr\bigl(X\setminus\Imm(\imath)\bigr)$, the restriction $\phi_X|Y:Y\rightarrow\mathbb P^1_K\times_{\Spec K}\mathbb P^1_K$ of $\phi$ to it is constant. 

As $\mathbb P^1_K\times_{\Spec K}\mathbb P^1_K$ is projective, there exists an open subvariety $X_0$ of $X$ such that $X\setminus X_0$ is finite and we have a morphism $\phi_{X_0}:X_0\rightarrow \mathbb P^1_K\times_{\Spec K}\mathbb P^1_K$ that extends $\phi$. 

Let $W$ be the normalization of the Zariski closure of the diagonal embedding $X_0\rightarrow X\times_{\Spec K} \mathbb P^1_K\times_{\Spec K} \mathbb P^1_K$ in $X\times_{\Spec K} \mathbb P^1_K\times_{\Spec K} \mathbb P^1_K$. We have two projections $\Sigma_1:W\rightarrow X$ and $\Sigma_2:W\rightarrow\mathbb P^1_K\times_{\Spec K} \mathbb P^1_K$, with $\Sigma_1$ projective and birational as $\Sigma_{1,X_0}:W\times_X X_0\rightarrow X_0$ is an isomorphism and with $\Sigma_2$ affine. Therefore $\mathcal O(W)=\mathcal O(X)=\mathcal O(X_0)$. Hence the fibers of $\Sigma_1$ are connected by properties of Stein Factorization (see \cite{Stacks3} or \cite{H2}, Ch.\ III, Cor.\ 11.5 and its proof). 

We show that the assumption that $\Sigma_1$ is not an isomorphism leads to a contradiction. This assumption implies that, by enlarging $X_0$, we can assume that $X_0\neq X$ and that for each $x\in X\setminus X_0$, the fiber $\Sigma_1^{-1}(x)$ is a connected curve. As $W\rightarrow X\times_{\Spec K} \mathbb P^1_K\times_{\Spec K} \mathbb P^1_K$ is quasi-finite, the restriction of the second projection $W\rightarrow \mathbb P^1_K\times_{\Spec K} \mathbb P^1_K$ to each irreducible component of $\Sigma_1^{-1}(x)$ is non-constant. Let $C$ be the Zariski closure of $X_0\setminus\Imm(\imath)$ in $W$. The open embedding $W\setminus C\rightarrow W$ is affine by Nagata's Theorem (see \cite{H1}, Cor.\ 3.3). Thus the resulting morphism $W\setminus C\rightarrow \mathbb P^1_K\times_{\Spec K} \mathbb P^1_K$ is affine and quasi-finite, hence an open embedding (by Zariski's Main Theorem) which contains $\mathbb A^1_K\times_{\Spec K} \mathbb A^1_K$ in such a way that the complement $[W\setminus C]\setminus \mathbb A^1_K\times_{\Spec K} \mathbb A^1_K$ is a non-empty divisor. Therefore $W\setminus C$ is either $\mathbb P^1_K\times_{\Spec K} \mathbb P^1_K$ or isomorphic to $\mathbb P^1_K\times_{\Spec K} \mathbb A^1_K$. Thus the $K$-algebra $\mathcal O(W\setminus C)$ is either $K$ or isomorphic to $K[t]$ which contradicts the fact that it has a $K$-subalgebra isomorphic to $\mathcal O(X)$. 

So $\Sigma_1$ is an isomorphism, hence we can take $\phi_X:=\Sigma_2\circ\Sigma_1^{-1}:X\rightarrow \mathbb P^1_K\times_{\Spec K} \mathbb P^1_K$. The uniqueness of $\phi_X$ and $\chi(\pi)_X$ holds as their targets are separated.

We show that the assumption that there exists $Y\in\Irr\bigl(X\setminus\Imm(\imath)\bigr)$ such that $\phi_X|Y$ is non-constant leads to a contradiction. By replacing $X$ with its open subvariety $X_Y:=X\setminus (\cup_{Z\in\Irr(X\setminus\Imm(\imath),Z\neq Y} Z)$ which is affine by Nagata's Theorem, we can assume that $\Irr\bigl(X\setminus\Imm(\imath)\bigr)=\{Y\}$. As $\phi_X|Y$ is non-constant, $\phi_X$ is quasi-finite, hence an open embedding by Zariski's Main Theorem. This implies that $X$ is isomorphic to $\mathbb A^1_K\times_{\Spec K}\mathbb P^1_K$ minus a finite number of points, hence the $K$-algebra $\mathcal O(X)$ is isomorphic to $K[t]$, a contradiction.\end{proof}

\begin{theorem}\label{T10}
Let $X$ be a connected affine normal surface over $K$ with an open embedding $\imath:\mathbb A^2_K\rightarrow X$. Assume that $Y:=X\setminus\Imm(\imath)$ is irreducible (i.e., $X$ is a sphere-like surface). Then the following properties hold.

\medskip
{\bf (1)} There exists an $\mathbb A^1$-fibration $\Sigma:X\rightarrow\mathbb P^1_K$ such that $\Imm(\imath)=\Sigma^{-1}(\mathbb A^1_K)$ (hence $Y=\bigl(\Sigma^{-1}(1:0)\bigr)_{\red}$) and the induced morphism $\Sigma^{-1}(\mathbb A^1_K)\rightarrow \mathbb A^1_K$ is a line bundle. 

\smallskip
{\bf (2)} Let $q(\Sigma)\in\mathbb N^{\ast}$ be the multiplicity of $Y$ in $\Sigma^{-1}\bigl((1:0)\bigr)$. We have $q(\Sigma)=1$ iff $\Sigma$ is a line bundle.

\smallskip
{\bf (3)} Let $\Sigma:X\rightarrow\mathbb P^1_K$ be an $\mathbb A^1$-fibration as in part (1) and let $q:=q(\Sigma)$. Then there exists a finite flat morphism $c:C\rightarrow\mathbb P^1_K$ with $C$ a smooth connected complete curve over $K$ such that the following properties hold.

\medskip\noindent
{\bf (3.a)} The normalization of $X\times_{\Sigma,\mathbb P^1_K,c} C$ has an affine open cover formed by line bundles over $C$.

\smallskip\noindent
{\bf (3.b)} We can assume (as needed) that either $|c^{-1}(1:0)|=1$ or $C=\mathbb P^1_k$ and the field extension $\Frac(c):K(t)\rightarrow K(t)$ induced by $c$ at the level of fields of fractions is separable (or Galois). 

\smallskip\noindent
{\bf (3.c)} If $\chr(K)\nmid q$, then we can take $c:\mathbb P^1_{K,\s}\rightarrow\mathbb P^1_{K,\t}$ to be defined by the rule $(x:y)\mapsto (x^q:y^q)$ (hence $C=\mathbb P^1_K$, $|c^{-1}(1:0)|=1$, and $\Frac(c)$ is cyclic). 

\medskip
{\bf (4)} The curve $Y$ is rational with one point at infinity.

\smallskip
{\bf (5)} We have $Y\cong\mathbb A^1_K$.

\smallskip
{\bf (6)} We have $|\Sing(X)|\le 1$.

\smallskip
{\bf (7)} If $\chr(K)\nmid q$ and $|\Sing(X)|=1$, then the singularity of $X$ is a quotient singularity by a subgroup of the cyclic group $\mathbb Z/q\mathbb Z$. 

\smallskip
{\bf (8)} If $\chr(K)|q$ and $|Sing(X)|=1$, then the singularity of $X$ is a quotient singularity by a finite group.

\smallskip
{\bf (9)} Each singularity of $X$ is rational. 

\smallskip
{\bf (10)} If $\chr(K)=0$ and $|\Sing(X)|=1$, then the singularity of $X$ is log-terminal.\end{theorem}

\begin{proof}
We identify $\mathbb A^1_K=\Spec K[t]$ and $K(t)$ with the field of fractions of $\mathbb A^1_K$ or $\mathbb P^1_K$.

For $i\in\{1,2\}$ let $\pi_i:\mathbb A^2_K\rightarrow\mathbb A^1_K$ be a projection such that the morphism $\pi_1\times\pi_2:\mathbb A^2_K\rightarrow \mathbb A^1_K\times_{\Spec K} \mathbb A^1_K$ is an isomorphism. 

Let $\phi:\mathbb A^2_K\to\mathbb P^1_K\times_{\Spec K}\mathbb P^1_K$ be the composite of $\pi_1\times\pi_2$ with the product open embedding $\mathbb A^1_K\times_{\Spec K}\mathbb A^1_K\subset\mathbb P^1_K\times_{\Spec K}\mathbb P^1_K$ and let $\phi_X:X\rightarrow\mathbb P^1_K\times_{\Spec K}\mathbb P^1_K$ be its unique extension to $X$ (see Lemma \ref{L11}(1)). The image $\phi_X(Y)$ is one point $(P_1,P_2)\in (\mathbb P^1_K\times_{\Spec K}\mathbb P^1_K)(K)$ by Lemma \ref{L11}(2). We write $P_i=(\alpha_i:\beta_i)$ with $(\alpha_i,\beta_i)\in K^2\setminus\{(0,0)\}$. 

As $X$ is separated, its local ring $\mathcal O_{X,Q}$ at a point $Q\in Y(K)$ cannot dominate its local ring $\mathcal O_{X,\imath(\gamma_1,\gamma_2)}$ with $(\gamma_1,\gamma_2)\in K^2$, hence $\beta_1\beta_2=0$. If $i\in\{1,2\}$ is such that $\beta_i=0$, then the composite of $\phi_X$ with the $i$-th projection $\mathbb P^1_K\times_{\Spec K}\mathbb P^1_K\rightarrow\mathbb P^1_K$ is a flat surjective morphism $\Sigma:X\rightarrow\mathbb P^1_K$ whose generic fiber is isomorphic to $\mathbb A^1_{K(t)}$. For $(\alpha,\beta)\in K^2\setminus\{(0,0)\}$, the fiber $\Sigma^{-1}\bigl((\alpha:\beta)\bigr)$ is isomorphic to $\mathbb A^1_K$ if $\beta\neq 0$ and its reduced variety is $Y$ if $\beta=0$. So $\Sigma$ is an $\mathbb A^1$-fibration on $X$. As the induced morphism $\Sigma^{-1}(\mathbb A^1_K)\rightarrow \mathbb A^1_K$ is isomorphic to a projection $\mathbb A^2_K\rightarrow\mathbb A^1_K$, it is a line bundle. So part (1) holds.

The `if' of part (2) is clear. If $q=1$, i.e., if $\Sigma^{-1}\bigl((1:0)\bigr)$ is reduced (equivalently integral), then $\Sigma$ is an $\mathbb A^1$-line bundle over $\mathbb P^1_K$ by \cite{KM}, Thm.\ 1 and thus a line bundle by \cite{BCW1}, Thm.\ 4.4. So the `only if' of part (2) holds, $X$ is regular and $Y\cong\mathbb A^1_K$ which also follow directly from \cite{KM}, Lem.\ 1.3. So part (2) holds and parts (3) to (10) hold if $q=1$.

To prove (3) to (10) in general we assume $q>1$. 

For a finite field extension $\mathcal K$ of $K(t)$, let $C:=C(\mathcal K)$ be the normalization of $\mathbb P^1_K$ in $\mathcal K$; we have a finite flat morphism $c:C\rightarrow\mathbb P^1_K$. Let $X_C^{\n}$ be the normalization of $X\times_{\Sigma,\mathbb P^1_K,c} C$. We have two projections: a finite surjective morphism $\Sigma_1:X_C^{\n}\to X$ and an $\mathbb A^1$-fibration $\Sigma_2:X_C^{\n}\rightarrow C$. 

From \cite{Epp}, Thm.\ (2.0) and \cite{Kuh} we get that we can choose $\mathcal K$ such that the fiber $c^{-1}\bigl((1:0)\bigr)$ is one point $P$, i.e., $c^{-1}\bigl((1:0)\bigr)_{\red}\cong\Spec K$, with $\Sigma_2^{-1}(P)$ reduced (hence a curve over $K$). Let $s:=|\Irr\bigl(\Sigma_2^{-1}(P)\bigr)|\in\llbracket1,[\mathcal K:K(t)]\rrbracket$; we write $\Irr\bigl(\Sigma_2^{-1}(P)\bigr)=\{Y_1,\ldots,Y_s\}$. For $i\in \llbracket1,s\rrbracket$, let $W_i:=\cup_{Z\in \Irr\bigl(\Sigma_2^{-1}(P)\bigr)\setminus\{Y_i\}} Z$, and let $X_i:=X_C^{\n}\setminus W_i$. So $X_i$ is affine by Nagata's Theorem and the fibers of the restriction $\Sigma_2|X_i:X_i\rightarrow C$ are integral. Thus as above we argue that $\Sigma_2|X_i$ is an $\mathbb A^1$-line bundle over $C$. As $(\Sigma_2|X_i)^{-1}(P)=Y_i\setminus (Y_i\cap W_i)\cong\mathbb A^1_K$ and $Y_i$ is affine, we have $Y_i\cap W_i=\emptyset$ and $Y_i\cong\mathbb A^1_K$. Moreover, $X_C^{\n}=\cup_{i=1}^s X_i$.

\phantomsection{From \cite{Stacks4} we get that we can choose $\mathcal K$ to be a separable extension of $K(t)$ such that for each $P\in c^{-1}\bigl((1:0)\bigr)(K)$ the fiber $\Sigma_2^{-1}(P)$ is reduced. As above we argue that $X_C^{\n}$ has an affine open cover formed by line bundles over $C$. Concretely, if for $l:=|c^{-1}\bigl((1:0)\bigr)(K)|\in\mathbb N^{\ast}$ we write $c^{-1}\bigl((1:0)\bigr)(K)=\{P_1,\ldots,P_l\}$ and if for each $j\in\llbracket1,l\rrbracket$ we write $\Irr\bigl(\Sigma_2^{-1}(P_j)\bigr)=\{Y_{1,j},\ldots,Y_{s_j,j}\}$ with $s_j\in\llbracket1,[\mathcal K:K(t)]\rrbracket$, then for $\lambda:=(i_1,\ldots,i_l)\in \prod_{j=1}^l \llbracket1,s_j\rrbracket$, let $W_{i_j,j}:=\cup_{Z\in \Irr\bigl(\Sigma_2^{-1}(P_j)\bigr)\setminus\{Y_{i_j,j}\}} Z$, and it follows that $X_\lambda:=X_C^{\n}\setminus\cup_{j=1}^l W_{i_j,j}$ is a line bundle over $C$, $Y_{i_j,j}\cap W_{i_j,j}=\emptyset$, and $Y_{i_j,j}\cong\mathbb A^1_K$. We have $X=\cup_{\lambda\in \prod_{j=1}^l \llbracket1,s_j\rrbracket} X_{\lambda}$. By enlarging $\mathcal K$, we can assume it is a Galois extension of $K(t)$ of finite Galois group $G$; clearly, $X=(X_C^{\n})/G$.}\label{PH51}

From the last two paragraphs we get that parts (3.a) and (3.b) hold. Note that part (3.a) implies that $\Sigma_2$ is smooth.  

In this paragraph we assume that $\chr(K)\nmid q$. We take $\mathcal K=K(t^\frac{1}{q})$; the field extension $K(t)\rightarrow\mathcal K$ is cyclic (thus Galois) and the fiber $c^{-1}\bigl((1:0)\bigr)$ is one point $P$ with $\Sigma_2^{-1}(P)$ reduced. Clearly, $c:C=\mathbb P^1_{K,\s}\rightarrow \mathbb P^1_{K,\t}$ is defined by the rule $(x:y)\mapsto (x^q:y^q)$ and $s\mid q$. Thus part (3.c) holds. So part (3) holds.

If $c^{-1}\bigl((1:0)\bigr)$ is one point (resp.\ if $\mathcal K$ is a separable extension of $K(t)$), then as the morphism $\mathbb A^1_K\cong Y_i\rightarrow Y$ (resp.\ $\mathbb A^1_K\cong Y_{i_j,j}\rightarrow Y$) induced by $\Sigma_1$ is finite surjective, part (4) holds by Lemma \ref{L10}(1) and (2).

To prove that parts (5) to (7) hold when $\chr(K)\nmid q$, we choose $c$ as in part (3.c), so $c$ is a finite Galois cover over $\mathbb P^1_{K,\t}\setminus\{(1:0),(0:1)\}$ of Galois group $G\cong\mathbb Z/q\mathbb Z$. Hence $\Sigma_1:X_C^{\n}\rightarrow X$ is a finite Galois cover over the complement of $\Sigma^{-1}(0:1)\cup\mathcal B$ in $X$, where $\mathcal B$ is a finite subset of $Y(K)$. The group $G$ acts on $\Sigma_1$ and we have $X=(X_C^{\n})^G$. As $\chr(K)\nmid q$, all $G$-modules over $K$ are semisimple and the functor $\star^G$ on $G$-modules $\star$ over $K$ is exact. Thus $Y$ is the quotient of the action of $G$ on either $\Sigma_1^{-1}(Y)$ or (as $Y$ is reduced) $\Sigma_1^{-1}(Y)_{\red}=\sqcup_{i=1}^s Y_i$. Hence $Y=(Y_i)^{G_i}$ with $G_i$ a cyclic group of order $\frac{q}{s}$ which is a subgroup of $\Aut(Y_i)$. Thus $Y$ is normal and hence regular; from this and part (4) we get that $Y\cong\mathbb A^1_K$. Moreover, there exists a unique point $P_i\in Y_i(K)$ such that the stabilizer of $P_i$ in $G_i$ is non-trivial, and in fact this stabilizer is $G_i$ itself. So for $Q_i\in Y_i(K)\setminus\{P_i\}$, the completions of the local rings $\mathcal O_{X,\Sigma_1(Q_i)}$ and $\mathcal O_{X_C^{\n},Q_i}$ are naturally identified; thus $\Sigma_1(Q_i)\notin\Sing(X)$. As $Y_i\rightarrow Y$ is surjective, it follows that $\Sing(X)\subset\{\Sigma_1(P_i)\}$ and hence $|\Sing(X)|\le 1$. So parts (5) and (6) if $\chr(K)\nmid q$ and part (7) hold. Note that we can take $\mathcal B=\Sing(X)\subset Y(K)$ and, if $X$ is singular, that $\Sigma_1(P_i)$ does not depend on $i\in \llbracket1,s\rrbracket$.

To prove that parts (5), (6), and (8) hold we could assume that $\chr(K)|q$ but we include arguments which actually work in general. 

For part (5), let $R\cong K[[x]]$ be the completion of $\mathcal O_{\mathbb P^1_K,(1:0)}$ and let 
$$\Sigma_R:X\times_{\Sigma,\mathbb P^1_K} \Spec R\rightarrow\Spec R$$ 
be the pullback of $\Sigma$ via the natural morphism $\Spec R\rightarrow\mathbb P^1_K$. We consider a projective morphism $\Sigma_R^+:X_R^+\rightarrow\Spec R$ that extends $\Sigma_R$ and with $X_R^+$ an integral normal scheme that contains $X_R$ as an open subscheme. The $K$-algebra of global functions on $X_R$ and hence also on $X_R^+$ is $R$. From this and the fact that $\Sigma_R^+$ is flat, we get that the $K$-algebra of global functions on the special fiber $X_K^+$ of $\Sigma_R^+$ (over $\Spec K$) is $K$. As the generic fiber of $\Sigma_R^+$ (over $\Spec(\Frac R)$) is isomorphic to $\mathbb A^1_{\Frac(R)}$, it has a divisor of degree $1$. So $\Sigma_R^+$ is cohomologically flat by (the implication $(ii)\Rightarrow (iv)$ in) \cite{Ray}, Thm.\ 8.2.1. Thus $H^1(X_R^+,\mathcal O_{X_R^+})=0$, hence for the curve $Y^+\in\Irr(X_K^+)$ that contains $Y$ we have $H^1(Y^+,\mathcal O_{Y^+})=0$. Thus $Y^+\cong\mathbb P^1_K$. From this and part (4) we get that part (5) holds.\footnote{Part (5) is also a particular case of Remark \ref{R5}(1).}

For parts (6) and (8), we chose $c$ such that $\mathcal K$ is a Galois extension of $K(t)$ and we use the above notation for the curves $Y_{i_j,j}$ and the Galois group $G$. Fixing a $j\in\llbracket1,l\rrbracket$ and an $i_j\in\llbracket1,s_j\rrbracket$, let $H_{i_j,j}$ be the subgroup of $G$ that stabilizes $Y_{i_j,j}$. Let $S$ be the completion of the local ring $\mathcal O_{C,P_j}$. If $i_j$ is the $j$-th component of the $l$-tuple $\lambda=(i_1,\ldots,i_l)$, then $X_\lambda\times_C \Spec S\cong\mathbb A^1_S$, $H_{i_j,j}$ is a subgroup of the Galois group $\Gal\bigl(\Frac(S)/\Frac(R)\bigr)$ and the morphism 
$$(X_\lambda\times_C \Spec S)/H_{i_j,j}\rightarrow \widehat{X}:=X\times_{\Sigma,\mathbb P^1_K} \Spec R$$ 
induces isomorphisms at the level of completions of local rings of residue fields $K$. Thus $|\Sing(\widehat{X})|$ is the number of singular points of  $(X_\lambda\times_C \Spec S)/H_{i_j,j}$. From this and Corollary \ref{C2.5} applied to $R$ equal to $S^{H_{i_j,j}}$ we get that $|\Sing(\widehat{X})|\le 1$. Thus part (6) holds. As $X=(X_C^{\n})/G$, part (8) holds. 

Part (9) follows from Proposition \ref{PR25}(3) (resp.\ \ref{PR25+}) and part (7) if $\chr(K)\nmid q$ (resp.\ part (8) in general). 

Part (10) follows from part (7) and \cite{Kaw}, Prop.\ 1.7 and the paragraph after it.\end{proof}

\begin{remark}\normalfont\label{R7}
{\bf (1)} For a triple $(l_1,l_2,l_3)\in (\mathbb N^{\ast})^3$ with $g.c.d.(l_1,l_2)=1$ and $\min(l_1,l_2)\ge 2$ the flat integral $k$-algebra homomorphism 
$$k[t^{l_1},t^{l_2}]\rightarrow k[x,t^{l_1},t^{l_2}]/(x-t)^{p^{l_3}}$$ 
has a non-regular integral domain source and a non-reduced target whose reduction is regular. If $p\nmid n$, for the non-regular integral domain $R:=k[x,y]/(x^n-y^p)$, and the flat $R$-algebra $S:=A[z]/(z^p-x)$ whose reduction is isomorphic to $k[x]$, the rule $t\mapsto z^n-y$ defines an $R$-algebra homomorphism $A\otimes_k k[t]/(t^p)\rightarrow S$ that makes $S$ a flat $k[t]/(t^p)$-algebra. These examples explain why the proof of Theorem \ref{T10}(5) is not worked out using only schemes over $Y$.

\smallskip
{\bf (2)} Theorem \ref{T10}(1), (4), and (6) in general and Theorem (5) in the case when $X$ is regular are proved in \cite{Mi3}, Sect.\ 2 and Lem.\ 1 provided $\chr(K)=0$.

\smallskip
{\bf (3)} If $\chr(K)=0$ and $|\Sing(X)|=1$, then the fact that the singularity of $X$ is rational was first proved in \cite{Mi3}, Thm.\ 1(1) and it also follows from the more general result \cite{FZ}, Thm.\ 0.3 or 1.12 which applies as $X\setminus Y=\Imm(\imath)$ is a union of closed $\mathbb A^1_K$ curves in $X$ that do not pass through the singular point of $X$ in $Y$.

\smallskip
{\bf (4)} If in the proof of Theorem \ref{T10}(6) we have $|H_{i_j,j}|=\chr(K)\mid q$, then either $(X_\lambda\times_C \Spec S)/H_{i_j,j}\rightarrow \Spec S^{H_{i_j,j}}$ or $X_\lambda\times_C \Spec S\rightarrow (X_\lambda\times_C \Spec S)/H_{i_j,j}$ is a smooth morphism, hence $\widehat{X}$ is regular; so $X$ is also regular.\end{remark}

\begin{corollary}\label{C2.6}
Let $X$ be a connected affine normal surface over $K$ with an open embedding $\imath:\mathbb A^2_K\rightarrow X$ which is not surjective. Then the following properties hold.

\medskip
{\bf (1)} Each connected component of $X\setminus\Imm(\imath)$ is irreducible and isomorphic to $\mathbb A^1_K$. 

\smallskip
{\bf (2)} Let $Y\in\Irr\bigl(X\setminus\Imm(\imath)\bigr)$. The union $\mathbb S^2_Y:=\Imm(\imath)\cup Y$ is an open affine subvariety of $X$ and we have $|Y\cap\Sing(X)|\le 1$.\end{corollary}

\begin{proof}
Let $Y\in\Irr\bigl(X\setminus\Imm(\imath)\bigr)$ and $W:=\cup_{Z\in \Irr\bigl(X\setminus\Imm(\imath)\bigr)\setminus\{Y\}} Z$. Then $X\setminus W$ is affine by Nagata's Theorem and the union of $\Imm(\imath)$ and $Y\setminus (Y\cap W)$. If $Y\cap W\neq\emptyset$, then $Y\setminus (Y\cap W)$ is an irreducible curve with at least $2$ points at infinity, a contradiction to Theorem \ref{T10}(4) applied to $X\setminus W$. Thus $Y\cap W=\emptyset$. So the connected components of $X\setminus\Imm(\imath)$ are irreducible. As $Y=Y\setminus (Y\cap W)$, $\mathbb S^2_Y=X\setminus W$ is affine. From these and Theorem \ref{T10}(5) and (6) we get that corollary holds.\end{proof}

\begin{corollary}\label{C2.7}
Let $X$ be a connected affine normal surface over $K$. Then $X$ is a sphere-like surface iff there exists a surjective $\mathbb A^1$-fibration $\Sigma:X\rightarrow\mathbb P^1_K$ with all fibers irreducible and with at most one non-reduced fiber.
\end{corollary}

\begin{proof}
The `only if' part holds by Theorem \ref{T10}(1). 

For the `if' part assume that all fibers of $\Sigma:X\rightarrow\mathbb P^1_K$ except $\Sigma^{-1}(1:0)$ are integral. Thus $\Sigma^{-1}(\mathbb A^1_K)\rightarrow\mathbb A^1_K$ is an $\mathbb A^1_K$-bundle by \cite{KM}, Thm.\ 1 and hence a line bundle by \cite{BCW1}, Thm.\ 4.4. So $\Sigma^{-1}(\mathbb A^1_K)\cong \mathbb A^2_K$ and the `if' part holds as $\Irr(X\setminus\Sigma^{-1}(\mathbb A^1_K))=\{\bigl(\Sigma^{-1}(1:0)\bigr)_{\red}\}$ has cardinality $1$.\end{proof}

\begin{definition}\label{D6}
{\bf (1)} Let $X$ be a connected affine normal surface over $K$ with an open embedding $\imath:\mathbb A^2_K\rightarrow X$ which is not surjective. We call the open affine cover $X=\cup_{Y\in\Irr\bigl(X\setminus\Imm(\imath)\bigr)} \mathbb S^2_Y$ as the canonical open affine cover\index{canonical open affine cover} of $X$ by sphere-like surfaces associated to $\imath$, where $\mathbb S^2_Y:=\Imm(\imath)\cup Y$.

\smallskip
{\bf (2)} Let $X$ be a sphere-like surface over $K$ for which there exists an $\mathbb A^1$-fibration $\Sigma:X\rightarrow\mathbb P^1_K$ such that $\Sigma^{-1}(\mathbb A^1_K)\rightarrow \mathbb A^1_K$ is a line bundle and the fiber $\Sigma^{-1}(1:0)$ is irreducible but non-reduced (in other words, $\Sigma$ is as in Theorem \ref{T10}(1) and (2) with $q=q(\Sigma)\ge 2$). For each such $\Sigma$ and every $P\in\mathbb P^1_K(K)\setminus\{(1:0)\}$, $X_{\Sigma,P}:=X\setminus\Sigma^{-1}(P)$ is an affine pseudo-plane over $K$ in the sense of \cite{MM1}, Def.\ 2.1\footnote{Based on Theorem \ref{T10}(5), the hypotheses $F\cong\mathbb A^1$ in \cite{MM1}, Def.\ 2.1(2) is superfluous.} and hence we call it a standard affine pseudo-plane\index{standard aﬀine pseudo-plane} of $X$ (of type $q$).
\end{definition}

\begin{remark}\normalfont\label{R8}
{\bf (1)} Let $e\in\EE_n(K)$. For each approximation $Z_e=\Spec B_e$ of $X_e$, the inclusion $B_e\rightarrow A_e$ defines a finite surjective morphism $X_e\rightarrow Z_e$ and we have an inequality $\rho(e)\ge|\Irr(Z_e\setminus V_e)|$, where $V_e$ is the image of $\Imm(\imath_e)$ in $Z_e$. Note that $X_e$ is the normalization of each model of it.

\smallskip
{\bf (2)} Let $e\in\EE_2(K)$. For two projections $\pi_1,\pi_2:\mathbb A^2_{K,\s}\rightarrow\mathbb A^1_K$, let $Z_e(\pi_1,\pi_2)$ be the Zariski closure in $\mathbb P^1_K\times_{\Spec K} \mathbb P^1_K\times_{\Spec K}\mathbb A^2_{K,\t}$ of the image of the locally closed embedding $\chi(\pi_1)\times\chi(\pi_2)\times e:\mathbb A^2_{K,\s}\rightarrow\mathbb P^1_K\times_{\Spec K} \mathbb P^1_K\times_{\Spec K}\mathbb A^2_{K,\t}$. Then $\chi(\pi_1)\times\chi(\pi_2)\times e$ extends to a finite surjective morphism $X_e\rightarrow Z_e(\pi_1,\pi_2)$ defined by $\chi_e(\pi_1)\times\chi_e(\pi_2)\times\psi_e$. Note that $Z_e(\pi_1,\pi_2)$ is an approximation of $X_e$. Moreover, if $\pi_1\times\pi_2:\mathbb A^2_{K,\s}\rightarrow\mathbb A^1_K\times_{\Spec K} \mathbb A^1_K$ is an isomorphism, then $Z_e(\pi_1,\pi_2)$ is a model of $X_e$.

\smallskip
{\bf (3)} In Definition \ref{D6}(2), the $\mathbb A^1$-fibration $\Sigma_P:X_{\Sigma,P}\rightarrow\mathbb P^1_K\setminus\{P\}\cong \mathbb A^1_K$ induced by $\Sigma$ has only one non-reduced fiber $\Sigma^{-1}(1:0)=\Sigma_P^{-1}(1:0)$ of multiplicity $q$ and has all fibers over $K$-valued points such that their reduced varieties are isomorphic to $\mathbb A^1_K$ by Theorem \ref{T10}(5).\end{remark} 

For $e\in\EE_2(K)$, the \'etale loci $X_e^{\et}$ can be described as follows.

\begin{corollary}\label{C2.8}
For a non-finite $e\in\EE_2(K)$ let $X_e=\cup_{Y\in\Irr(\Gamma_{\psi_e})} \mathbb S^2_Y$ be the canonical open affine cover associated to $\imath_e$. We write $\mathbb S^2_Y=\Spec A_Y$. Then the following properties hold.

\medskip
{\bf (1)} The following three statements are equivalent.

\medskip\noindent
{\bf (1.a)} The open subvariety $\mathbb S^2_Y$ is contained in the \'etale locus $X_e^{\et}$ of $\psi_e$.

\smallskip\noindent
{\bf (1.b)} The generic point $\eta_Y$ of $Y$ belongs to $X_e^{\et}$.

\smallskip\noindent
{\bf (1.c)} The tangent bundle $T_{\mathbb S^2_Y}$ over $\mathbb S^2_Y$ is trivial.

\medskip
Moreover, if $\chr(K)=p$, then statements (1.a) to (1.c) are also equivalent to either one of the following extra statements.

\medskip\noindent
{\bf (1.d)} The ring $A_Y$ has a $p$-basis.

\smallskip\noindent
{\bf (1.e)} The ring (or the $K$-algebra) $A_Y$ has a differential basis.

\medskip
{\bf (2)} We have an identity $X_e^{\et}=\Imm(\imath_e)\bigcup_{Y\in\Irr(\Gamma_{\psi_e}),T_{\mathbb S^2_Y}\;\textup{is a trivial bundle}}\; \mathbb S^2_Y$. If $\chr(K)=p$ then the identity $X_e^{\et}=\Imm(\imath_e)\bigcup_{Y\in\Irr(\Gamma_{\psi_e}),A_Y\;\textup{has a}\; p\textup{-basis}}\; \mathbb S^2_Y$ holds.

\smallskip
{\bf (3)} The open subvariety $X_e^{\et}$ of $X_e$ is affine regular, contains $\Imm(\imath_e)=\Imm(\psi_e)$, and the complement $X_e^{\et}\setminus\Imm(\imath_e)$ is a disjoint union of $\rho_{\et}(e)$ irreducible components of $\Gamma_{\psi_e}$ isomorphic to $\mathbb A^1_K$.

\smallskip
{\bf (4)} The complement $X_e\setminus\Imm(\imath_e)$ is a disjoint union of $\rho(e)$ irreducible components of $\Gamma_{\psi_e}$ isomorphic to $\mathbb A^1_K$.

\smallskip
{\bf (5)} If $Y\subset X_e^{\et}$ then $\mathbb S^2_Y$ is regular, and if $Y\subset X_e\setminus X_e^{\et}$ then $|\Sing(\mathbb S^2_Y)|\le 1$.

\smallskip
{\bf (6)} We have $|\Sing(X_e)|\le\rho(e)-\rho_{\et}(e)$.  

\smallskip
{\bf (7)} If $Y\in \Irr_{\textup{t}}(X_e\setminus X_e^{\et})$ is such that $\psi_e(Y)\cong\mathbb A^1_K$, then $Y\subset\Reg(X_e)$.
\end{corollary}

\begin{proof}
Part (1) is a particular case of Theorem \ref{T3+}(1) and (2). 

Part (2) follows from part (1). 

For part (3), $X_e^{\et}$ is regular as $\Imm(\psi_e)=\mathbb A^2_{K,\s}$ and is affine by Nagata's Theorem as $X_e\setminus X_e^{\et}$ is a union of irreducible divisors of $X_e$ by part (2) and Corollary \ref{C2.6}(1). The last statement of part (3) follows from part (2), Corollary \ref{C2.6}(1), and the definition of $\rho_{\et}(e)$. So part (3) holds. 

Similarly, part (4) follows from Corollary \ref{C2.6}(1) and the definition of $\rho(e)$. 

The first part of (5) is clear and the second part of (5) follows from Corollary \ref{C2.6}(2). 

Part (6) follows from part (5) as the number of $Y\in\Irr(\Gamma_{\psi_e})$ with $A_Y$ not having a $p$-basis is $\rho(e)-\rho_{\et}(e)$.

For part (7), for each $Q\in Y(K)$ we have $\psi_e(Q)\in\Reg\bigl(\psi_e(Y)\bigr)$ as $\psi_e(Y)\cong\mathbb A^1_K$. Based on this and part (4) we get that the hypotheses of Criterion \ref{CRI1}(1) hold for each such point $Q$. So $Y\subset\Reg(X_e)$ by Criterion \ref{CRI1}(1).\end{proof}

Corollary \ref{C2.8}(2) justifies the following definition.

\begin{definition}\label{D7}
{\bf (1)} Let $e\in\EE_2(K)$ be non-finite. If $X_e^{\et}\neq\Imm(\imath_e)$, then we refer to the open affine cover $X_e^{\et}=\cup_{Y\in\Irr(\Gamma_{\psi_e}),T_{\mathbb S^2_Y}\;\textup{is a trivial bundle}}\; \mathbb S^2_Y$ as the canonical open affine cover of $X_e^{\et}$ by sphere-like surfaces associated to $\imath_e$.

\smallskip
{\bf (2)} The {\it sphere class rank invariant}\index{class rank!sphere} and the {\it almost sphere class rank invariant}\index{class rank!almost sphere}
$$\rho_1=\rho_{1,K,2}:\EE_2(K)\rightarrow\mathbb N\;\;\;\textup{and}\;\;\; \rho_2=\rho_{2,K,2}:\EE_2(K)\rightarrow\mathbb N$$
(respectively) are defined for $e\in\EE_2(K)$ by the following three rules. 

\medskip\noindent
{\bf (2.a)} If $e$ is finite then $\rho_1(e):=0$.

\smallskip\noindent
{\bf (2.b)} If $e$ is non-finite then $\rho_1(e):=|\{Y\in\Irr(\Gamma_{\psi_e})|\mathbb S^2_Y\cong\mathbb S^2_k\}|$.

\smallskip\noindent
{\bf (2.c)} We define $\rho_2(e):=\rho_{\et}(e)\setminus\rho_1(e)$.

\medskip
{\bf (3)} The {\it quasi-sphere class rank invariant}\index{class rank!quasi-sphere} and the {\it pseudo-sphere class rank invariant}\index{class rank!pseudo-sphere}
$$\rho_3=\rho_{3,K,2}:\EE_2(K)\rightarrow\mathbb N\;\;\;\textup{and}\;\;\; \rho_4=\rho_{4,K,2}:\EE_2(K)\rightarrow\mathbb N$$
(respectively) are defined by the two following rules.

\medskip\noindent
{\bf (3.a)} We define $\rho_3(e):=\rho(e)-\rho_{\et}(e)-|\Sing(X_e)|$.

\smallskip\noindent
{\bf (3.b)} We define  $\rho_4(e):=|\Sing(X_e)|$.
\end{definition}

From definitions and Corollaries \ref{C2.8}(2) and (4) and \ref{C2.6} we get directly the following result.

\begin{corollary}\label{C2.9}
Let $e\in\EE_2(K)$ be non-finite. Then the canonical open affine cover of $X_e$ by sphere-like surfaces associated to $\imath_e$ consists of exactly $\rho_1(e)$ spheres, $\rho_2(e)$ almost spheres, $\rho_3(e)$ quasi-spheres, and $\rho_4(e)$ pseudo-spheres.
\end{corollary}

\begin{remark}\normalfont\label{R8.3}
For $e\in\EE_2(k)$, all cases of `strictness-equality' in the general two inequalities $0\le\rho_{\et}(e)\le\rho(e)$ occur with the Jacobian surface $\psi_e$ regular. 

\medskip
{\bf (1)} For $0<\rho_{\et}(e)=\rho(e)$ see $e_{f,g,k}$ of Section \ref{S2} with $(f,g)\in\Theta_{k,m}^1$, Example \ref{EX19}, Example \ref{EX23} for $n=2$, and the case $p=3$ of Examples \ref{EX29} and \ref{EX30}.

\smallskip
{\bf (2)} For $0=\rho_{\et}(e)<\rho(e)$ see Remark \ref{R22} and the case $p\ge 5$ of Examples \ref{EX29} and \ref{EX30}.

\smallskip
{\bf (3)} For $0<\rho_{\et}(e)<\rho(e)$ see end of Example \ref{EX11} for $mq>2$.

\smallskip
{\bf (4)} For $0=\rho_{\et}(e)=\rho(e)$ see Remark \ref{R10}.
\end{remark}

\section{Proof of Theorem \ref{T1}(1)}\label{S12}

We construct two families of examples. The easier one, Example \ref{EX6}, works for almost all $l$, and the harder one, Example \ref{EX7}, works for all $l$.

\begin{example}\normalfont\label{EX6}
Let $m\in \mathbb N^\ast$ and $s\in\llbracket0,pm\rrbracket$. Let $\wp(m,s):=p^2m-s(p-1)$. Let $g(x)\in k[x]$ be a polynomial of degree $m$ and let $\gamma_1,\ldots,\gamma_{pm}\in k$ be the $pm$ distinct zeros of the equation $x+g(x)^p=0$. Let $f(x)\in k[x]$ be such that for $i\in\llbracket1,pm\rrbracket$ we have $f(\gamma_i)=0$ iff $i\in \llbracket1,s\rrbracket$; if $s=pm$ we assume that $f\neq 0$. Thus $\deg(f)\ge s$. Then the $F_2$-system $\mathcal W$ over $k$ defined by
$$x+g(x)^p=0=y+y^pf(x)^p$$ 
has $\s(\mathcal W)=s+p(pm-s)=\wp(m,s)$ solutions computed as follows: for $i\in \llbracket1,s\rrbracket$ we have the solution $(\gamma_i,0)$ and for $i\in\llbracket1,pm\rrbracket\setminus \llbracket1,s\rrbracket$ there exist precisely $p$ solutions with the first coordinate equal to $\gamma_i$. From Proposition \ref{PR1}(1) to (3) we get the first two inequalities and the first two equalities in
\begin{equation}\label{EQ2}
pm\le\pi_{\i}(e)\le\pi_{\l}(e)=\pi_{\r}(e)=\pi(\mathcal W)= p\max\bigl(m,\deg(f)+1\bigr)\ge p\max(m,s+1).
\end{equation}
\end{example}
E.g., if $s\le\deg(f)\le m-1$, then we have $\pi_{\i}(e)=\pi_{\l}(e)=\pi_{\r}(e)=\pi(\mathcal W)=pm$.

\begin{lemma}\label{L12}
For $\mathcal D:=\{\wp(m,s)|m\in \mathbb N^\ast, s\in\llbracket0,pm\rrbracket\}$ we have inclusions
$$\llbracket0,p-1\rrbracket\subset\mathbb N\setminus\mathcal D=\{0\}\sqcup\{pq-r|q,r\in \llbracket1,p-1\rrbracket,q\le r\}\subset\llbracket0,(p-1)^2\rrbracket.$$ 
\end{lemma}

\begin{proof} Let $r:=pm-s\in\llbracket0,pm\rrbracket$. We have $\wp (m,s)=pm+(p-1)r$ with $m\in \mathbb N^\ast$ and $r\in\llbracket0,pm\rrbracket$. Thus the set $\mathcal D$ is equal to 
$$\{pm+(p-1)r|m\in\mathbb N^\ast, r\in\llbracket0,pm\rrbracket\}=\{pm+(p-1)r|(m,r)\in\mathbb N^\ast\times \llbracket0,p-1\rrbracket\}.$$
Denoting $q:=m+r$, we have $\wp(m,s)=pq-r$ with $q\ge r+1$. Therefore $\mathcal D=\{pq-r|(q,r)\in\mathbb N^\ast\times \llbracket0,p-1\rrbracket, q\ge r+1\}$. For $r=0$ we get $0\notin\mathcal D\supset p\mathbb N^\ast$ and for $r\in \llbracket1,p-1\rrbracket$ the smallest and the largest numbers in $\mathbb N\setminus\mathcal D$ congruent to $-r$ modulo $p$ are $p-r$ and $pr-r$ (respectively), which for $r=p-1$ are $1$ and $(p-1)^2$ (respectively). The lemma follows.
\end{proof}

\begin{example}\normalfont\label{EX7}
Let $(m,q)\in (\mathbb N^\ast)^2$. Let $f(y)\in k[y]$ be such that $f(0)\neq 0$ and the equation $f(y)=0$ has $q$ distinct solutions in $k$, to be denoted as $\gamma_1,\ldots,\gamma_q$; so $\deg(f)\ge q$. For $i\in \llbracket1,q\rrbracket$, let $\varepsilon_i:=\sqrt[p]{\gamma_i}\in k$. We consider the $k$-algebra
$$A:=k[x,y]/\bigl(y-x^pf(y)^p\bigr).$$
As $f(0)\neq 0$, the zero locus $f(y)=y-x^pf(y)^p=0$ is empty and hence we have $u:=f(y)+\bigl(y-x^pf(y)^p\bigr)\in A^{\ast}$. Moreover, $u$ is a $p$-th power of an element of $A$ as 
$y+\bigl(y-x^pf(y)^p\bigr)=\bigl(xf(y)+\bigl((y-x^pf(y)^p)\bigr)^p\in A$. So the inverse of $u$ is a $p$-th power. Thus there exists a polynomial $g\in k[x,y]$ such that $f(y)\bigl(g(x,y)\bigr)^p-1$ is divisible by $y-x^pf(y)^p$. If $z:=\sqrt[p]{y}$, so $y=z^p$, we get that we have an isomorphism
$$A\cong k[z]\Bigl[\frac{1}{\prod_{i=1}^q (z-\varepsilon_i)}\Bigr]$$ 
induced by the $k$-algebra epimorphism 
$$\theta:k[x,y]\rightarrow k[z,f(z^p)^{-1}]$$
that maps $y$ to $z^p$ and $x$ to $zf(z^p)^{-1}$; note that $\theta\bigl(f(y)x\bigr)=z$, $\theta\bigl(y-x^pf(y)^p\bigr)=0$, and $\theta\bigl(g(x,y)^p\bigr)=f(z^p)^{-1}$.

We now study the solution set $\mathcal J_0\subset k^2$ of the $F_2$-system $\mathcal U_0$ over $k$ which is defined by
$$x-\bigl(g(x,y)\bigr)^p-x^{pm}f(y)^{pm}\bigl(g(x,y)\bigr)^p=y-x^pf(y)^p=0$$
and which satisfies
\begin{equation}\label{EQ3}
\pi(e_{\mathcal U_0})\le\pi(\mathcal U_0)= p[m+m\deg(f)+\deg(g)].
\end{equation}
As $\theta\bigl(x-(g(x,y))^p-x^{pm}f(y)^{pm}(g(x,y))^p\bigr)=f(z^p)^{-1}(-z^{pm}+z-1)$, we have 
$$\mathcal J_0=\Bigl\{\Bigl(\frac{z}{f(z^p)},z^p\Bigr)|z\in k\setminus\{\varepsilon_i|i\in \llbracket1,q\rrbracket\},\; z^{pm}-z+1=0\Bigr\}.$$
As the equation $z^{pm}-z+1=0$ has exactly $pm$ solutions $\delta_1,\ldots,\delta_{pm}$ in $k$, we have $\s(\mathcal U_0)=|\mathcal J_0|=pm-r$ with
$$r:=|\{\varepsilon_i|i\in \llbracket1,q\rrbracket\}\cap \{\delta_i|i\in\llbracket1,pm\rrbracket\}|\in\llbracket0,\min(q,pm)\rrbracket.$$

Let $\alpha\in k^{\ast}$; if $\deg(f)=m$, we assume that $-\alpha f(y)$ is not monic. A similar argument shows that the solution set of the $F_2$-system $\mathcal U_{\alpha}$ 
$$x-\bigl(g(x,y)\bigr)^p-x^{pm}f(y)^{pm}\bigl(g(x,y)\bigr)^p-\alpha=y-x^pf(y)^p=0$$
is $\mathcal J_{\alpha}:=\bigl\{\bigl(\frac{z}{f(z^p)},z^p\bigr)|z\in k\setminus\{\varepsilon_i|i\in \llbracket1,q\rrbracket\},\; z^{pm}-z+1+\alpha f(z^p)=0\bigr\}.$
We have $\deg\bigl(z^{pm}-z+1+\alpha f(z^p)\bigr)=p\max\bigl(m,\deg(f)\bigr)$. Hence, as $\varepsilon_i$ with $i\in \llbracket1,q\rrbracket$ is a zero of $z^{pm}-z+1$ iff it is a zero of $z^{pm}-z+1+\alpha f(z^p)$ and as the equation $z^{pm}-z+1+\alpha f(z^p)=0$ has exactly $p\max\bigl(m,\deg(f)\bigr)$ solutions in $k$, we get that
$$\s(\mathcal U_{\alpha})=|\mathcal J_{\alpha}|=p\max\bigl(m,\deg(f)\bigr)-r.$$ 

Let now $l\in\mathbb N$. We take $m\in\mathbb N^\ast$, $m\ge\lceil\frac{l}{p}\rceil$, $r=pm-l$, and we choose $f(y)$ so that $\{\varepsilon_i|i\in \llbracket1,q\rrbracket\}\cap \{\delta_i|i\in\llbracket1,pm\rrbracket\}$ has precisely $r$ elements; so $\deg(f)\ge q\ge r$. We get that $\s(\mathcal U_0)=l$. If $\deg(f)>m$, then we have relations 
$$\s(\mathcal U_0)=l=pm-r<p\deg(f)-r=\s(\mathcal U_{\alpha})\le\deg(e_{\mathcal U_0})=\deg(e_{\mathcal U_{\alpha}});$$
see \cite{Gro5}, Prop.\ (18.2.8) for the last (i.e., non-strict) inequality. 

Assume $l=0$. We take $m=1$, $q=r=p$, and $f(y)=y^p-y+1\in\mathbb F_p[y]$; note that the set of zeros of $f$ is $\{\delta_i|i\in\llbracket1,p\rrbracket\}=\{\delta_i^{\frac{1}{p}}|i\in\llbracket1,p\rrbracket\}$.

Assume $l=1$ and $p>2$. We choose $m$ such that $p|2^m-1$, i.e., $2$ is a zero of $z^{pm}-z+1$, say $2=\delta_1$. For $f(y)=\prod_{i=2}^{pm} (y-\delta_i)=\frac{y^{pm}-y+1}{y-2}\in\mathbb F_p[y]$ we have $q=r=pm-1$.

If $f(y)\in\mathbb F_p[y]$, then we can choose $g(x,y)\in\mathbb F_p[x,y]$, so $\mathcal U_0$ is defined over $\mathbb F_p$.
\end{example}

Theorem \ref{T1}(1) follows from Example \ref{EX7}.

For $l\in\mathcal D$ (see Lemma \ref{L12}), Equation (\ref{EQ2}) gives simpler and smaller upper bounds for algebraic degrees than Relations (\ref{EQ3}) applied to an $f$ with $\s(\mathcal U_0)=l$.

\begin{remark}\normalfont\label{R9}
With $m$, $q$, and $f(y)$ as in Example \ref{EX7}, we can write everything explicitly, which is how Example \ref{EX1} was obtained. Let $f_1(y)\in k[y]$ be the unique polynomial such that $f(y)=f(0)+yf_1(y)$. Working modulo $\bigl(y-x^pf(y)^p\bigr)$ we get
$$1=f(0)^{-1}\bigl(f(y)-yf_1(y)\bigr)\equiv f(0)^{-1}\bigl(f(y)-x^pf(y)^pf_1(y)\bigr)$$
$$= f(0)^{-1}f(y)\bigl(1-x^p f(y)^{p-1}f_1(y)\bigr)\equiv f(0)^{-1}f(y)\bigl(1-x^p f(x^pf(y)^p)^{p-1}f_1(x^pf(y)^p)\bigr),$$
i.e., we can take $g(x,y)^p$ to be $f(0)^{-1}\bigl(1-x^pf(x^pf(y)^p)^{p-1}f_1(x^pf(y)^p)\bigr).$ Thus
$$(x,y)\mapsto\bigl(x-f(0)^{-1}(1+x^{pm}f(y)^{pm})(1-f_1(x^pf(y)^p)x^pf(x^pf(y)^p)^{p-1}),y-x^pf(y)^p\bigr)$$
is the rule that defines $e_{\mathcal U_0}$.

If $f(y)=1-y+y^p$, then $q=p$, $f(0)=1$, and $f_1(y)=-1+y^{p-1}$; if $m=1$, then $\mathcal U_0$ is the $F_2$-system of Example \ref{EX1}. With $m\in\mathbb N^{\ast}$ arbitrary, let $e=:e_{\mathcal U_0}$ (it depends on $m$ and $f$). To compute $\deg(e)$ we substitute $z:=xf(y)$. The system 
$$x-\bigl(1+x^{pm}f(y)^{pm}\bigr)\bigl(1-x^pf(x^pf(y)^p)^{p-1}f_1(x^pf(y)^p)\bigr)-x_1=0=y-x^pf(y)^p-x_2,$$ 
gives 
$$y=z^p+x_2\;\;\;\textup{and}\;\;\;x=\frac{z}{z^{p^2}-z^p+1+h_2}$$ 
with $h_2:=x_2^p-x_2\in k[x_1,x_2]$. Substituting these expressions in $f(y)^p$ times the first equation of the system and $zf_1(z)=f(z)-1$ we get that
$$zf(z^p+x_2)^{p-1}-x_1f(z^p+x_2)^p=(z^{pm}+1)[f(z^p+x_2)^p-z^pf(z^p)^{p-1}f_1(z^p)]$$
$$=(z^{pm}+1)[f(z^p+x_2)^p-f(z^p)^p+f(z^p)^{p-1}]=(z^{pm}+1)[f(z^p)^{p-1}+h_2^p].$$
The left (resp.\ right) hand side is a polynomial in $z$ of degree $p^3$ (resp.\ $p^3-p^2+pm$) with leading coefficient $-x_1$ (resp.\ $1$). Thus
$$\deg(e)=\max(p^3-p^2+pm,p^3)$$ 
and the field extension $\Frac(e^{\#})$ is isomorphic to $k_{\t}(x_1,x_2)\rightarrow k_{\t}(x_1,x_2)[z]/\bigl(h(z)\bigr)$, where 
$$h(z):=(z^{pm}+1)[(z^{p^2}-z^p+1)^{p-1}+h_2^p]+x_1(z^{p^2}-z^p+1+h_2)^p-z(z^{p^2}-z^p+1+h_2)^{p-1}$$
is in $k[x_1,x_2][z]$. If $m>p$, then $\Spec\bigl(k_{\t}(x_1,x_2)[z]/(h(z))\bigr)$ is a monomial model of $X_e$.\end{remark}

\section{Proof of Theorem \ref{T1}(2)}\label{S13}

For the proof we need a new type of open embeddings of $\mathbb A^2_K$ into the surfaces $\mathbb S^2_{f,g,K}=\Spec\Bigl(K[x,y,z]/\bigl(xf(x)+yg(y)z\bigr)\Bigr)$ with $(f,g)\in\Theta_K$ of Section \ref{S2}.

\begin{proposition}\label{PR4.5}
Let $(f,g)\in\Theta_K$ with $\deg(f)\ge 1$ and $g(0)\neq 0$. Let $\delta\in K^{\ast}$ be such that $f\bigl(\delta g(0)\bigr)=0$. Let $h=h(x,y):=\delta+xy\in K[x,y]$. Then the rule on valued points
$$(x,y)\mapsto \left(hg(yh),yh,\frac{-f\bigl(hg(yh)\bigr)}{y}\right)$$
defines an open embedding 
$$\jmath_{f,g,K;\delta}:\mathbb A^2_K\rightarrow\mathbb S^2_{f,g,K}$$ 
with the property that $\mathbb S^2_{f,g,K}\setminus\Imm(\jmath_{f,g,K;\delta})$ is the disjoint union of $\deg(f)\z(tg)$ curves isomorphic to $\mathbb A^1_K$.
\end{proposition}

\begin{proof}
As $hg(yh)\in \delta g(0)+yK[x,y]$ and $f\bigl(\delta g(0)\bigr)=0$, $y$ divides $-f\bigl(hg(yh)\bigr)$ and hence the rule is well defined. By computing $-f\bigl(hg(yh)\bigr)$ modulo $y^2$ we get
\begin{equation}\label{EQ3+}
\frac{-f\bigl(hg(yh)\bigr)}{y}\in -f'\bigl(\delta g(0)\bigr)g(0)x-\delta^2f'\bigl(\delta g(0)\bigr)g'(0)+yK[x,y].
\end{equation} 

As $y$ (resp.\ $h$) are rational functions in the first and third (resp.\ and second) coordinates of the rule and $x$ is a rational function in $y$ and $yh$, $\jmath_{f,g,K;\delta}$ is birational. To show that it is an open embedding and to identify its complement, for each $(\alpha,\beta,\gamma)\in\mathbb S^2_{f,g,K}(K)$, i.e., for $(\alpha,\beta,\gamma)\in K^3$ with $\alpha f(\alpha)+\beta g(\beta)\gamma=0$, we solve in $K^2$ the system of three equations
$$hg(yh)-\alpha=yh-\beta=\frac{-f\bigl(hg(yh)\bigr)}{y}-\gamma=0$$
in the indeterminates $x$ and $y$ by considering three cases as follows.

{\bf Case 1: $\alpha=0$ (so $\beta g(\beta)\gamma=0$).} As $h(x,0)g(0)=\delta g(0)\neq 0=\alpha$, we have $y\neq 0$. So $\frac{-f\bigl(hg(yh)\bigr)}{y}=\frac{-f(0)}{y}\in K^{\ast}y^{-1}$. So the system is inconsistent if $\gamma=0$. We now assume that $\gamma\neq 0$. Thus $y=\frac{-f(0)}{\gamma}$ and, as $y(\delta+xy)=\beta$, we get that $x=\frac{\gamma[\beta\gamma+\delta f(0)]}{f(0)^2}$. As $hg(yh)=\frac{-\beta g(\beta)\gamma}{f(0)}=0=\alpha$, $\bigl(\frac{\gamma[\beta\gamma+\delta f(0)]}{f(0)^2},\frac{-f(0)}{\gamma}\bigr)$ is the only solution when $\gamma\neq 0$.

{\bf Case 2: $\alpha\neq 0$ and $g(\beta)=0$ (hence $f(\alpha)=0$).} As $\alpha\neq 0=hg(\beta)$, the system is inconsistent. 
 
{\bf Case 3: $\alpha g(\beta)\neq 0$.} By replacing $yh=\beta$ in $hg(yh)=\alpha$ we get that $h=\frac{\alpha}{g(\beta)}$, so $h\neq 0$, $y=\frac{\beta g(\beta)}{\alpha}$, and $xy=h-\delta=\frac{\alpha-\delta g(\beta)}{g(\beta)}$. We consider two subcases as follows.

{\bf Subcase 3.1: $\beta=0$ (hence $f(\alpha)=0$).} The second equation gives $y=0$. Thus, if $\alpha$ is a zero of $f$ different from $h(x,0)g(0)=\delta g(0)$, then the system is inconsistent. If $\alpha=\delta g(0)$, then the first equation and the equation $xy=\frac{\alpha-\delta g(\beta)}{g(\beta)}=0$ hold. As $y=0$, the third equation becomes $-f'\bigl(\delta g(0)\bigr)g(0)x-\delta^2f'\bigl(\delta g(0)\bigr)g'(0)=\gamma$ by Relation (\ref{EQ3+}), so the only solution is $\left(\frac{-\gamma-\delta^2f'\bigl(\delta g(0)\bigr)g'(0)}{f'\bigl(\delta g(0)\bigr)g(0)},0\right)$.%$\frac{-f\bigl(hg(yh)\bigr)}{y}-\gamma=0$ 

{\bf Subcase 3.2: $\beta\neq 0$.} Thus $y\neq 0$ and therefore $x=\frac{\alpha[\alpha-\delta g(\beta)]}{\beta [g(\beta)]^2}$. The pair $\left(\frac{\alpha[\alpha-\delta g(\beta)]}{\beta [g(\beta)]^2},\frac{\beta g(\beta)}{\alpha}\right)$ is indeed the only solution as $\alpha f(\alpha)+\beta g(\beta)\gamma=0$ implies that the third equation $\frac{-f\bigl(hg(yh)\bigr)}{y}-\gamma=0$ also holds.

These cases imply that the birational morphism $\jmath_{f,g,K;\delta}$ is also radicial, hence it is an open embedding. 

These cases also imply that the complement $\mathbb S^2_{f,g,K}\setminus\Imm(\jmath_{f,g,K;\delta})$ is the disjoint union of curves isomorphic to $\mathbb A^1_K$ and given by equations $x=z=0$ (see Case 1) and equations $x-\alpha=y-\beta=0$ with $(\alpha,\beta)\in K^2$ satisfying either $f(\alpha)=g(\beta)=0$ (see Case 2) or $f(\alpha)=\beta=0$ with $\alpha\neq\delta g(0)$ (see Subcase 3.1). The number of these curves is $1+\deg(f)\z(g)+[\deg(f)-1]=\deg(f)[\z(g)+1]=\deg(f)\z(tg)$, where the last identity holds as $g(0)\neq 0$.\end{proof}

The following concrete application of Proposition \ref{PR4.5} proves Theorem \ref{T1}(2).

\begin{example}\normalfont\label{EX8}
Let $(m,s)\in \mathbb N^{\ast}\times\mathbb N$. Let $\alpha_0:=0\in k$ and $\alpha_1:=1\in k$. If $s\ge 2$, let $a_2,\ldots,a_s$ be distinct elements of $k^{\ast}\setminus\{1\}$. If $s=0$ let $g(t):=1$ and if $s\ge 1$ let $(j_1,\ldots,j_s)\in (\mathbb N^{\ast})^{s}$ and $g(t):=\prod_{i=1}^s (t-\alpha_i)^{j_i}\in k[t]$. So $g$ is monic with $g(0)\neq 0$, $\z(g)=s$ and $\deg(g)=\sum_{i=1}^s j_i\ge s$. Let $(f,g)\in\Theta^1_{k,m}$ (so $\deg(f)=pm-1\ge 1$); e.g., we can take $f(t)=t^{pm-1}-1$. The smooth surface $\mathbb S^2_{f,g,k}$ is equipped with the finite \'etale projection $c_{f,g,k}:\mathbb S^2_{f,g,k}\rightarrow\mathbb A^2_{k,\t}$ on the last two coordinates (see Section \ref{S2}). Let $\delta\in k^{\ast}$ be such that $f\bigl(\delta g(0)\bigr)=0$. For $h=h(x,y):=\delta+xy\in k_{\s}[x,y]$, the degree of $\frac{-f\bigl(hg(yh)\bigr)}{y}\in k_{\s}[x,y]$ is $[3\deg(g)+2](pm-1)-1\ge (3s+2)(pm-1)-1$ with equality iff for $s\ge 1$ we have $j_1=\cdots=j_s=1$. Let (see Proposition \ref{PR4.5})
$$e:=c_{f,g,k}\circ\jmath_{f,g,k;\delta}\in\End_2(k).$$ 

As $e$ is the composite of an open embedding with a finite \'etale morphism of degree $pm$, we have $e\in\EE_2(k)$ and $\deg(e)=pm$. Clearly, $X_e=\mathbb S^2_{f,g,k}$ is a finite \'etale cover of $\mathbb A^2_{k,\t}$, so $\psi_e$ is \'etale.

The last paragraph of the proof of Proposition \ref{PR4.5} gives
$$\Gamma_{\psi_e}(k)=\{(0,\beta,0)|\beta\in k\}\cup\{(\alpha,\beta,\gamma)\in k^3|f(\alpha)=\beta g(\beta)=0,(\alpha,\beta)\neq \bigl(\delta g(0),0\bigr)\},$$ 
$$\Gamma_e=\{(\beta,\gamma)\in k^2|\nexists\;\alpha\in k\; \textup{with}\;(\alpha,\beta,\gamma)\in [\mathbb S^2_{f,g,k}\cap\Imm(\jmath_{f,g,k;\delta})](k)\}$$
$$=\{(\beta,0)\in k^2|g(\beta)=0\}.$$ Thus $\rho_{\et}(e)=\rho(e)=\deg(f)\z(tg)=(pm-1)(s+1)$ and $\varphi(e)=s$. 

For $(\beta,\gamma)\in k^2\setminus\Gamma_e$, from the three cases of the proof of Proposition \ref{PR4.5} we get that $|e^{-1}(\beta,\gamma)|$ is $1$ if either $g(\beta)=0\neq\gamma$ or $\beta=\gamma=0$, is $2$ if $\beta=0\neq \gamma$, is $pm-1$ if $\beta g(\beta)\neq 0=\gamma$, and is $pm$ if $\gamma\beta g(\beta)\neq 0$; so for $pm>2$ we have $\mathcal F(e)=\{1,2,pm-1\}$ if $s=0$ and $\mathcal F(e)=\{0,1,2,pm-1\}$ if $s\ge 1$ and for $pm=2$ we have $\mathcal F(e)=\{1\}$ if $s=0$ and $\mathcal F(e)=\{0,1\}$ if $s\ge 1$. If $pm=2$, i.e., $\deg(f)=1$, then $N_e=c_{f,g,k}(\Gamma_{\psi_e})$ is the zero locus $g(x_1)x_2=0$, hence $\nu(e)=s+1$ and $\Gamma_e=\Sing(N_e)$. If $pm\ge 3$, i.e., $\deg(f)\ge 2$, then $N_e=c_{f,g,k}(\Gamma_{\psi_e})$ is the zero locus $x_1g(x_1)x_2=0$, hence $\nu(e)=s+2$ and $\Gamma_e\subsetneq\Sing(N_e)$ as $|\Sing(N_e)|=s+1$. As Conditions ($\triangleright$) and ($\triangleleft$) of Proposition \ref{PR1} hold for the pair $\Bigl(yh,\frac{-f\bigl(hg(yh)\bigr)}{y}\Bigr)$ that defines $e$, Proposition \ref{PR1}(1) and (2) gives
$$\pi(e)\le\pi_{\l}(e)=\pi_{\r}(e)=(3\deg(g)+2)(pm-1)-1\ge (3s+2)(pm-1)-1,$$ 
and the last inequality is an equality iff for $s\ge 1$ we have $j_1=\cdots=j_s=1$. 

Let $l\in\mathbb N^{\ast}$ be such that $p^l\ge s+1$. If we take $f(t)=t^{pm-1}-1$, $\delta=g(0)^{-1}$, and $(\alpha_2,\ldots,\alpha_s)\in (\mathbb F_{p^l}\setminus\{0,1\})^{s-1}$, then $e$ is defined over $\mathbb F_{p^l}$. Thus Theorem \ref{T1}(2) holds.
\end{example}

\begin{proposition}\label{F2}
Let $(f,g)\in\Theta_K$ and $\delta\in K^{\ast}$ be such that $\deg(f)\ge 2$, $g(0)\neq 0$, and $f\bigl(\delta g(0)\bigr)=0$. Then the two open embeddings $\imath_{f,g,K}:\mathbb A^2_K\rightarrow \mathbb S^2_{f,g,K}$ and $\jmath_{f,g,K;\delta}:\mathbb A^2_K\rightarrow \mathbb S^2_{f,g,K}$ are not isomorphic and we have an identity 
$$\Imm(\imath_{f,g,K})\setminus Y_1=\Imm(\jmath_{f,g,K;\delta})\setminus Y_2,$$
with $Y_1$ and $Y_2$ as zero loci $x=z=0$ and $x-\delta g(0)=y=0$ (respectively) in $\mathbb S^2_{f,g,K}$.
\end{proposition}

\begin{proof}
The identity follows from the descriptions of the complements of the images of the two open embeddings (see Displays (\ref{EQ0a}) and (\ref{EQ0b}) and the proof of Proposition \ref{PR4.5}). There exists a unique monic polynomial $f_1(t)\in K[t]$ such that we have an identity $f(t)=\bigl(t-\delta g(0)\bigr)f_1(t)$; so $\deg(f_1)=\deg(f)-1\ge 1$. 

Recall that $h=\delta+xy$. We write $g(t)=g(0)+\sum_{i=1}^{\deg(g)} \alpha_it^i$ with $\alpha_{\deg(g)}=1$ and, if $\deg(g)\ge 2$, $(\alpha_1,\ldots,\alpha_{\deg(g)-1})\in K^{\deg(g)-1}$. Thus 

\begin{equation}\label{EQ3++}
\begin{split}
h_1(x,y):=\frac{h(x,y)g\bigl(yh(x,y)\bigr)-\delta g(0)}{y}= x\bigl[g\bigl(yh(x,y)\bigr)]+\delta\sum_{i=1}^{\deg(g)} \alpha_i y^{i-1}h(x,y)^i\\
= g(0)x+\sum_{i=1}^{\deg(g)} \alpha_i [yh(x,y)]^{i-1}h(x,y)^2\in xg(0)+K[h,yh]\subset K[x,y].\;\;
\end{split}
\end{equation}

Also recall that the open embeddings $\imath_{f,g,K}$ and $\jmath_{f,g,K;\delta}$ are given by the rules $(x,y)\mapsto\bigl(xyg(y),y,-xf(xyg(y))\bigr)$ and $(x,y)\mapsto\bigl(hg(yh),yh,-h_1(x,y)f_1(hg(yh))\bigr)$ (respectively). To show that these two open embeddings are not isomorphic it suffices to show that exists no $a\in\GA_2(K)$ such that $a^{\#}(B_1)=B_2$, where $B_1$ and $B_2$ are the $K$-subalgebras $K_{\s}[x_1x_2g(x_2),x_2,x_1f\bigl(x_1x_2g(x_2)\bigr)]$ and 
$$K_{\s}\bigl[h(x_1,x_2)g\bigl(x_2h(x_1,x_2)\bigr),x_2h(x_1,x_2),h_1(x_1,x_2)f_1\bigl(h(x_1,x_2)g(x_2h(x_1,x_2))\bigr)\bigr]$$ 
(respectively) of $K_{\s}[x_1,x_2]$  that define $\imath_{f,g,K}$ and $\jmath_{f,g,K;\delta}$ (respectively). 

That no such $a$ exists follows from the fact that there exists an element $y_1\in B_1$ that can be extended to $\{y_1,y_2\}\in\mathcal A_2(K)$ (e.g., $y_1=x_2$) but no such element exists in $B_2$. To check this for $B_2$, let $B_3$ be the $B_2$-subalgebra of $K_{\s}[x_1,x_2] $ generated by $h$. Based on the definition of $B_2$ and Equation (\ref{EQ3++}), we have 
$$B_3=K_{\s}[h,x_2h,g(0)x_1f_1\bigl(hg(x_2h)\bigr)]=[x_1x_2,x_2h,g(0)x_1f_1\bigl(hg(x_2h)\bigr)].$$ 
As for each $i\in\mathbb N^{\ast}$, we have $g(0)x_1h(x_2h)^i=g(0)(x_1x_2)h^2(x_2h)^{i-1}\in B_3$, we get that $B_3$ does not depend on $g$ and thus we can assume that $g=1$. Therefore $B_2=B_3=K_{\s}[\delta+x_1x_2,x_2(\delta+x_1x_2),-x_1f_1(\delta+x_1x_2)]$. 

Note that $(1,\bigl(x_1f_1(\delta+x_1x_2)\bigr)^i,\bigl(x_2(1+x_1x_2)\bigr)^i|i\in\mathbb N^{\ast})$ is a $K[\delta+x_1x_2]$-basis of the $K[\delta+x_1x_2]$-module $B_2$. Therefore, as $\deg(f_1)\ge 1$, each irreducible polynomial $z_1\in B_2\subset K_{\s}[x_1,x_2]$  is a non-monic polynomial in either $x_1$ or $x_2$. Thus the zero locus $z_1=0$ is a closed curve $C_1$ of $\mathbb A^2_{K,\s}$ which is not isomorphic to $\mathbb A^1_K$ as the first or second projection $C_1\rightarrow\mathbb A^1_K$ is a quasi-finite non-finite morphism. So there exists no $z_2\in K_{\s}[x_1,x_2] $ such that $\{z_1,z_2\}\in\mathcal A_2(K)$, hence indeed $B_1\not\cong B_2$.\end{proof}

\begin{example}\normalfont\label{EX10}
The normal and general analogs of $\Theta_R$ over a field $\mathcal K$ are 
$$\Theta_{\mathcal K}^{\n}:=(\mathbb N^{\ast}\setminus\{1\})\times\{1\}\times\{(f,g)\in \mathcal K[t]^2|f\;\textup{are}\; g\;\textup{are monic separable},f(0)g(0)\neq 0\}$$
and $\Theta_{\mathcal K}^{\gen}:=(\mathbb N^{\ast}\setminus\{1\})\times\mathbb N^{\ast}\times \{(f,g)\in\Theta_{\mathcal K}|g(0)\neq 0\}\supset \Theta_{\mathcal K}^{\n}$. For $(q,l,f,g)\in\Theta_K^{\gen}$ the affine surface 
$$\mathbb S^2_{q,l,f,g,K}:=\Spec\bigl(K[x,y,z]/(x^qf(x)+y^lg(y)z)\bigr)$$ 
is normal iff $(q,l,f,g)\in\Theta_{\mathcal K}^{\n}$; similarly, $(0,0,0)$ is a normal point of it iff $l=1$ and iff $(0,0,0)$ is a normal singular point of it. If $(q,1,f,g)\in\Theta_K^{\n}$, then we have $\Sing(\mathbb S^2_{q,1,f,g,K})=\{(0,\beta,0)|\beta\in K,g(\beta)=0\}$. The open embeddings $\imath_{f,g,K}$ or $\jmath_{f,g,K;\delta}$ have analogs for the general context in many cases. To exemplify this, let $(q,1,f,g)\in\Theta_K^{\gen}$ with $\deg(f)\ge 1$ and $g=g_1^q$ with $g_1(t)\in K[t]$ monic. Let $\delta\in K^{\ast}$ and $f_1(t)\in K[t]$ be such that $f(t)=\bigl(t-\delta g_1(0)\bigr)f_1(t)$. Let 
$$\jmath_{q,1,f,g,K;\delta}:\mathbb A^2_K\rightarrow \mathbb S^2_{q,1,f,g,K}$$ 
be given by the rule $(x,y)\mapsto \left(hg_1(yh^q),yh^q,\frac{-f(hg_1(yh^q))}{y}\right)$, where $h:=\delta+xy$. The proof of Proposition \ref{PR4.5} applies to give that $\jmath_{q,1,f,g,K;\delta}$ is a birational quasi-finite morphism and $\mathbb S^2_{q,1,f,g,K}\setminus\Imm(\jmath_{q,f,g,1,K;\delta})$ is the disjoint union of $\deg(f)\z(tg_1)$ closed subvarieties isomorphic to $\mathbb A^1_K$, the zero locus $x=z=0$ being the one that contains the normal singular point $(0,0,0)$. The birational morphism $\jmath_{q,1,f,g,K;\delta}$ is an open embedding (i.e., is injective) iff either $(q,1,f,g)\in\Theta_K^{\n}$ (i.e., $g=g_1=1$) or $q=\chr(K)^m$ with $m\in\mathbb N^{\ast}$. If $f_1(t)=g(t)=1$, so $(q,1,f,g)\in\Theta_K^{\n}$, then the rule $(x,y)\mapsto (h,y,-xh^q)$ defines an open embedding $\imath_{q,1,t-\delta,1,K}:\mathbb A^2_K\rightarrow\mathbb S^2_{q,1,t-\delta,1,K}$ and we have $\imath_{q,1,t-\delta,1,K}\circ\Inv_{\mathbb A^2_K}=\Inv_{\mathbb S^2_{q,1,t-\delta,1,K}}\circ\jmath_{q,1,t-\delta,1,K;\delta}$, where $\Inv_{\mathbb A^2_K}$ and $\Inv_{\mathbb S^2_{q,1,t-\delta,1,K}}$ are the involutions of $\mathbb A^2_K$ and $\mathbb S^2_{q,1,t-\delta,1,K}$ (respectively) defined by the rules $(x,y)\mapsto (y,x)$ and $(x,y,z)\mapsto (x,-z,-y)$ (respectively).\end{example}

\begin{example}\normalfont\label{EX10.5}
Let $q\in\mathbb N^{\ast}$. Let $\Sigma:\mathbb S^2_{q,1,t-1,1,K}\rightarrow\mathbb P^1_K$ be the morphism that maps $(\alpha,\beta,\gamma)\in \mathbb S^2_{q,1,t-1,1,K}(K)$ to either $\bigl(1-\alpha:\beta\bigr)$ or $(\gamma:\alpha^q)$ (whichever is defined; if both are defined, they are equal). The multiplicity of $\bigl(\Sigma^{-1}(1:0)\bigr)_{\red}$ in $\Sigma^{-1}(1:0)$ is $q$ and $\Sigma$ is a surjective $\mathbb A^1$-fibration as in Corollary \ref{C2.7}. Defining $\jmath_{1,1,f,g,K;\delta}$ by the same rule $(x,y)\mapsto \left(1+xy,y(1+xy),x\right)$ obtained by allowing $q=1$, then the composite $\Sigma\circ \jmath_{q,1,t-1,1,K;\delta}:\mathbb A^2_K\rightarrow\mathbb P^1_K$ is given by the rule $(x,y)\mapsto \bigl((x:1+xy)^q\bigr)$. Thus, if $q\ge 2$, then $\mathbb S^2_{q,1,t-1,1,K}$ is a sphere-like surface with $|\Sing(\mathbb S^2_{q,1,t-1,1,K})|=1$ by Example \ref{EX10}, we have $(q,1,t-1,1)\in\Theta_K^{\n}$, and $\mathbb S^2_{q,1,t-1,1,K}$ is not a monomial model of any $X_e$ with $e\in\EE_2(K)$ such that $\psi_e$ is regular.
\end{example}
 
\begin{proposition}\label{PR4.1}
Let $m\in\mathbb N$. Let $B_m$ be the $k$-subalgebra of $k_{\s}[x_1,x_2]$ generated by $x_1$ and $x_2^i(1+x_1x_2)^{i+1}$ with $i\in \llbracket0,m\rrbracket$. Then the following properties hold.

\medskip
{\bf (1)} There exists no pair $(g_1,g_2)\in B_m^2$ such that $\e(g_1,g_2)\in\EE_2(k)$.

\smallskip
{\bf (2)} Let $(q,1,f,g_1^q)\in\Theta_k^{\gen}$ with $\deg(f)\ge 1$ and $g_1(t)\in k[t]$ monic. Let $\delta\in k^{\ast}$. Then there exists no $e\in\EE_2(k)$ such that $\iota_e$ factors through the normalization of $\jmath_{q,1,f,g_1^q,k;\delta}$.
\end{proposition}
 
\begin{proof}
Let $h:=1+x_1x_2\in k_{\s}[x_1,x_2]$. We have a chain of strict inclusions 
$$k_{\s}[h]\subsetneq B_0\subsetneq B_1\subsetneq\cdots\subsetneq B_m\subsetneq\cdots.$$ 
We consider the $k$-linear derivations $h^2D_{x_1}$ and $D_{x_2}$ of $k_{\s}[x_1,x_2]$ and its $k$-subalgebra $B_{\infty}:=\cup_{i\ge 0} B_i$. Let $D:=h^2D_{x_1}+D_{x_2}$. One computes that $D(x_1)=h^2$, $D(x_2)=1$, and $D(h)=x_1+x_2h^2$. For $i\in\mathbb N^{\ast}\setminus\{1\}$, as $x_1x_2=h-1$, we compute
\begin{equation}\label{EQ3.8}
\begin{split}
D\bigl(x_2^{i-1}h^i\bigr)&=(i-1)x_2^{i-2}h^i+ix_2^{i-1}(x_1+x_2h^2)h^{i-1}\\
&=(i-1)x_2^{i-2}h^i+ix_2^{i-2}h^{i-1}(h-1)+ix_2^ih^{i+1}\\
&=ix_2^ih^{i+1}+(2i-1)x_2^{i-2}h^i-ix_2^{i-2}h^{i-1}.
\end{split}
\end{equation}
It follows that $D$ is a $k$-linear derivation of $B_m$ iff $p|m+1$.

For $m\ge 1$ we consider the normal affine surface 
$$\mathbb T^2_m:=\Spec\bigl(k[y_0,\ldots,y_{m+1}]/(y_1^3-y_1^2+y_0y_2,y_i^2-y_{i-1}y_{i+1}|i\in\llbracket2,m\rrbracket)\bigr);$$ 
it is a complete intersection in $\mathbb A^{m+2}_k$ with $\Sing(\mathbb T^2_m)=\{(0,\ldots,0)\}$. The rule $(x,y)\mapsto (-x,h,yh^2,y^2h^3,\ldots,y^mh^{m+1})$ defines a morphism 
$$\imath_{m,k}:\mathbb A^2_k\rightarrow\mathbb T^2_m$$ 
which is an open embedding as already the rule $(-x,h,yh^2)$ defines an open embedding by Example \ref{EX10} applied to $(z,x,y)$ instead of $(x,y,z)$. 

For a point $P=(\alpha_0,\ldots,\alpha_{m+1})\in\mathbb T^2_m(k)$, if $\alpha_0=\alpha_1=0$, then we have $\alpha_i=0$ for all $i\in\llbracket0,m\rrbracket$; hence $P$ is the singular point $(0,\ldots,0)$ iff $\alpha_0=\alpha_1=\alpha_{m+1}=0$. 

To prove part (1) we can assume that $p|m+1$; so $D$ is a $k$-linear derivation of $B_m$. Let $T_D$ be the global section of the tangent bundle $T_{\Reg(\mathbb T^2_m)}$ defined by $D$. 

We show that the assumption that $T_D$ is zero at a $k$-valued point $(\beta_0,\ldots,\beta_{m+1})$ of $\Reg(\mathbb T^2_m)(k)$ leads to a contradiction. Based on the formulas for $D(x_1)$ and $D(h)$, this assumption implies that $\beta_1^2=\beta_0+\beta_2=0$ and thus, as we also have $\beta_1^3-\beta_1^2=\beta_0\beta_2$, we get that $\beta_0=\beta_1=\beta_2=0$ and therefore $m\ge 2$ and we have $\beta_i=0$ for all $i\in\llbracket0,m\rrbracket$. From Equation (\ref{EQ3.8}) applied to $i=m\ge 2$, as $\beta_{m-2}=\beta_{m-1}=0$, we get that $m\beta_{m+1}=0.$ As $p|m+1$, it follows that $\beta_{m+1}=0$, so $(\beta_0,\ldots,\beta_{m+1})$ is the singular point, a contradiction. 

Thus $T_D$ is nowhere zero and it generates a direct summand of $T_{\Reg(\mathbb T^2_m)}$ which is free of rank $1$. Based on this, as moreover the conormal bundle of the complete intersection closed embedding $\mathbb T^2_m\rightarrow\mathbb A^{m+2}_k$ is trivial and the tangent bundle $T_{\mathbb A^{m+2}_k}$ is trivial, it follows that $T_{\Reg(\mathbb T^2_m)}$ itself is trivial. 

We show that the assumption that there exists a pair $(g_1,g_2)\in B_m^2$ such that $\e(g_1,g_2)\in\EE_2(k)$ leads to a contradiction. This assumption implies that there exists a morphism $\psi:\mathbb T^2_m\rightarrow\mathbb A^2_k$ whose restriction to $\Imm(\imath_{m,k})$ is \'etale. Denoting $\psi_0:=\psi|\Reg(\mathbb T^2_m)$, as in the proof of the implication $(1.d)\Rightarrow (1.a)$ of Theorem \ref{T3+}(1) we argue that the functorial morphism $T_{\mathbb T^2_m}\rightarrow \psi_0^*(T_{\mathbb A^2_k})$ is an isomorphism. Hence $\psi_0$ is \'etale. By the purity of the branch locus (see \cite{Gro6}, Exp.\ X, Thm.\ 3.4 (i)) we get that $\psi$ itself is \'etale, which contradicts the identity $|\Sing(\mathbb T^2_m)|=1$.

Part (2) follows from part (1) applied to $m=2$ once we remark that up to isomorphism we can assume that $\delta=1$ and that the $k$-subalgebra $B_{q,1,f,g_1^q,k}$ of $k[x,y]$ generated by $hg_1(yh^q)$, $yh^q$ and $\frac{-f(hg_1^q(yh^q))}{y}$ is contained in $B_{2,1,t-1,1,k}$ which, under the identification $(x,y)=(x_1,x_2)$, is $B_2$.
\end{proof}

\section{Proof of Theorem \ref{T1}(3)}\label{S14}

Theorem \ref{T1}(3) follows from the following example.

\begin{example}\normalfont\label{EX11}
Let $(m,q)\in (\mathbb N^\ast)^2$ be such that $p$ divides $mq$ and let $e\in\EE_2(k)$ be defined over $\mathbb F_p$ by the rule on valued points
\begin{equation}\label{EQ3.9}
(x,y)\mapsto\left(x+(1-x^my^{mq-1})^q,y-x^my^{mq}\right).
\end{equation}
If $p\mid q$ (resp.\ $p\mid m$) , then $e$ is a perturbation of the elementary $p$-morphism $\e(x,y-x^my^{mq})$ (resp.\ $\e\bigl(x+(1-x^my^{mq-1})^q,y\bigr)$).

For $(\alpha,\beta)\in k^2$ let $\s(\alpha,\beta)\in\mathbb N$ be such that $e^{-1}(\alpha,\beta)\cong\Spec\bigl(k^{\s(\alpha,\beta)}\bigr)$. So $\s(\alpha,\beta)$ is the number of solutions in $k^2$ of the system $\mathcal J_{\alpha,\beta}$ of equations
\begin{equation}\label{EQ4}
x+(1-x^my^{mq-1})^q-\alpha=0=y-x^my^{mq}-\beta.
\end{equation}
Let $z:=1-x^my^{mq-1}$. Then $\mathcal J_{\alpha,\beta}$ can be rewritten as
\begin{equation}\label{EQ5}
x+z^q=\alpha\;\;\;\textup{and}\;\;\; yz=\beta.
\end{equation}

\phantomsection{To compute every $\s(\alpha,\beta)$, we first consider the case when $\beta\neq 0$. Hence $y\neq 0$ and $z=\beta y^{-1}$. From this and the equation $x+z^q=\alpha$ we get that $x=\alpha-\beta^qy^{-q}$. Thus $xy^q=\alpha y^q-\beta^q$ and $(xy^q)^m=(\alpha y^q-\beta^q)^m$.
Hence $\mathcal J_{\alpha,\beta}$ is equivalent to the system of equations}\label{EXTRA5}
$$x=\alpha-\beta^qy^{-q}\;\;\;\textup{and}\;\;\; y-\beta=(\alpha y^q-\beta^q)^m.$$
If $\alpha\neq 0$, then $p\mid mq$ implies that the equation $y-\beta=(\alpha y^q-\beta^q)^m$ in $y$ has $mq$ distinct solutions in $k$, $0$ being a solution of it iff $\beta$ belongs to the set $\mathbb I_{p,m,q}^{\ast}$, where
$$\mathbb I_{p,m,q}:=\{x\in k|x=(-1)^{m-1}x^{mq}\}\;\;\;\textup{and}\;\;\;\mathbb I_{p,m,q}^{\ast}:=\mathbb I_{p,m,q}\cap k^{\ast}=\mathbb I_{p,m,q}\setminus\{0\}.$$ 
Therefore for $\beta\notin\mathbb I_{p,m,q}$ and $\alpha\neq 0$, we have $\s(\alpha,\beta)=mq$, which implies that $\deg(e)=mq$, and for $\beta\in\mathbb I_{p,m,q}^{\ast}$ and $\alpha\neq 0$ we have $\s(\alpha,\beta)=mq-1$. Similarly, if $\beta\notin\mathbb I_{p,m,q}$, we have $\s(0,\beta)=1$, the only solution being
$$\bigl(-[1+(-1)^m\beta^{mq-1}]^{-q},\beta+(-1)^m\beta^{mq}\bigr),$$ and if $\beta\in\mathbb I_{p,m,q}^{\ast}$ we have $\s(0,\beta)=0$. 

We now consider the case $\beta=0$. Thus $yz=0$. If $y=0$, then $z=1$ and we get one solution $(\alpha-1,0)$ of $\mathcal J_{\alpha,0}$ (see System (\ref{EQ4})). If $z=0$, then $x=\alpha$ (see System (\ref{EQ5})) and hence $z=0=1-\alpha^my^{mq-1}$ and we get no solutions if $\alpha=0$ and $mq-1$ solutions of the form $(\alpha,\sqrt[mq-1]{\alpha^{-m}})$ if $\alpha\neq 0$, all of them being different from the prior solution $(\alpha-1,0)$. We conclude that $\s(0,0)=1$ and $\s(\alpha,0)=mq$ if $\alpha\neq 0$.

From the last two paragraphs we get first that $\mathcal F(e)=\{0,1,mq-1\}$ and second that 
\begin{equation}\label{EQ5.1}
\Gamma_e=\{(0,\beta)|\beta\in\mathbb I_{p,m,q}^{\ast}\},
\end{equation} 
hence $\varphi(e)=mq-1$. Moreover, $N_e$ is the zero locus $x_1\prod_{\beta\in\mathbb I_{p,m,q}^{\ast}} (x_2-\beta)=0$, hence $\nu(e)=1+|\mathbb I_{p,m,q}^{\ast}|=mq$ and $\Gamma_e=\Sing(N_e)$.

For $\beta\in\mathbb I_{p,m,q}^{\ast}$ we have $e(0,\beta)=(1,\beta)$ and $\s(1,\beta)=mq-1$, and therefore $|e^{-1}(1,\beta)\setminus\Gamma_e|=mq-2$. Thus, for $mq>2$ we have $|\mathcal E_2(e)|=0$ and from Lemma \ref{L1}(3) we get that $\Imm(e)=\Imm(e^2)$; so $\iota(e)=1$. Similarly, if $mq=2$ (so $p=2$), then $\mathbb I_{2,m,q}=\mathbb F_2$, $\Imm(e)\setminus\Imm(e^2)=\{(1,1)\}$ and thus $\varphi(e^2)=2$; as $e(1,1)=(1,0)$ and $\s(1,0)=2$, we get $\mathcal E_2(e)=\mathcal E_1(e)$, so $\varphi(e^3)=\varphi(e^2)$ by Lemma \ref{L1}(3) and $\iota(e)=2$. 

Clearly, if $p$ divides both $m$ and $q$, then $e\in\BE_2(k)$. 

Moreover, if $q=1$ we have $\pi_{\i}(e)=\pi_{\l}(e)=\pi_{\r}(e)=2m$ and if $q\ge 2$ we have $mq+m\le\pi_{\i}(e)\le\pi_{\l}(e)=\pi_{r}(e)=mq+(mq-1)q$ by Proposition \ref{PR1}(1) to (4). For $r:=mq\in p\mathbb N^\ast$ fixed, the smallest possible value of $\pi_{\l}(e)=\pi_{\r}(e)$ is $2r=2mq$ obtained only for $(m,q)=(r,1)$.

Recall that $X_e=\Spec A_e$; we have inclusions $k_{\t}[x_1,x_2]\subset A_e\subset k[x,y]$ defined by $e$, with $x_1=x+z^q$, $x_2=yz$, and $z=1-x^my^{mq-1}$. So here we identify $\mathbb A^2_{k,\s}=\Spec k[x,y]$. It follows that $e$ is the composite of birational quasi-finite morphisms $\mathbb A^2_{k,\s}\rightarrow X\rightarrow Y$ with a finite flat morphism $Y\rightarrow\mathbb A^2_{k,\t}$ of degree $mq$, where $X:=\Spec\bigl(k_{\t}[x_1,x_2][z]/(f(z))\bigr)$ with
$$f(z):=z^{mq}-z^{mq-1}+(x_1-z^q)^mx_2^{mq-1}$$ 
and where $Y:=\Spec\bigl(k_{\t}[x_1,x_2][w]/(g(w)\bigr)$ is a monomial model of $X_e$ with
$$g(w):=w^{mq}-w^{mq-1}+\sum_{i=1}^m (-1)^{m-i}\binom{m}{i}x_1^ix_2^{mq-1}\bigl(1+(-1)^mx_2^{mq-1}\bigr)^{qi-1}w^{q(m-i)}$$ in $k_{\t}[x_1,x_2][w]$ obtained from $\bigl(1+(-1)^mx_2^{mq-1}\bigr)^{mq-1}f(z)$ through the substitution $w:=z\bigl(1+(-1)^mx_2^{mq-1}\bigr)$. Thus $w\in A_e$. If $mq>2$, then for $\alpha\in k^{\ast}$ the $mq-1$ points of $\mathbb A^2_{k,\s}$ of the form $(\alpha,\sqrt[mq-1]{\alpha^{-m}})$ map to the same point of $Y$ and therefore $Y$ is not normal by Zariski's Main Theorem, i.e., $Y\neq X_e$. If $mq=2$ (so $p=2$) then $\mathbb I_{p,m,q}$ is $\mathbb F_2$ and $g(w)$ is $w^2+w+x_1^2x_2(x_2+1)$ if $(m,q)=(2,1)$ and is $w^2+w+x_1x_2(x_2+1)$ if $(m,q)=(1,2)$; hence $Y=X_e\rightarrow\mathbb A^2_{k,\t}$ is a finite Galois cover, $A_e\cong k_{\t}[x_1,x_2][w]/\bigl(g(w)\bigr)$, and the complement of the open embedding $\mathbb A^2_{k,\s}\rightarrow Y$ consists of $2$ disjoint $\mathbb A^1_k$ curves given by the equations $w=x_2=1$ and $w=x_1=0$ implying that $\rho_{\et}(e)=\rho(e)=2$. For $(m,q)=(1,2)$, $e$ is $e_{t+1,t+1,k}$ of Section \ref{S2}.

Let $\epsilon\in\{1,2\}$ be $1$ iff $m=1$. Similarly, $e$ factors through a finite flat morphism $e_Z:Z\rightarrow\mathbb A^2_{k,\t}$, where $Z:=\Spec\bigl(k_{\t}[x_1,x_2][v]/(h(v))\bigr)$ is a monomial model of $X_e$ with 
$$h(v)=v^{mq}+\Big(\sum_{i=1}^{m-\epsilon} h_{mq-iq}v^{mq-iq}\Big)+h_1v+h_0\in k_{\t}[x_1,x_2][v]$$ obtained from $x_1^{m(q+\epsilon-2)}x_2z^{-mq}f(z)$ via the substitution $v:=x_1^{\epsilon} x_2z^{-1}$; so each $h_j\in k_{\t}[x_1,x_2]$ is uniquely determined by $m$, $q$, the index $j$, and $v\in A_e$. E.g., $h_0(x_1,x_2)=x_1^{m(q+\epsilon-2)}x_2\bigl(1-(-1)^mx_2^{mq-1}\bigr)$ and $h_1(x_1,x_2)=-x_1^{m(q+\epsilon-2)-\epsilon}$; also, if $q\ge 2$, then 
$$h_{mq-iq}(x_1,x_2)=x_1^{mq-m}x_2(-1)^{i}\binom{m}{i}x_1^{m-i}x_2^{mq-1}(vz)^{qi-mq}=(-1)^{i}\binom{m}{i}x_1^{i(q-1)}x_2^{qi}.$$ 
As $h_1(x_1,x_2)=-x_1^{m(q+\epsilon-2)-\epsilon}$ becomes $-1$ when $mq=2$ (so $p=2$), it follows that for $mq=2$ the morphism $e_Z$ is a finite Galois cover and thus $X_e=Z$, and for $mq>2$ the \'etale locus of $e_Z$ is where $x_1$ is invertible and therefore the $k_{\t}[x_1,x_2]_{x_1}$-algebra $A_{e,x_1}\cong k_{\t}[x_1,x_2]_{x_1}[v]/(h(v))$ is \'etale. 

We show that $A_e$ is always regular and that $\psi_e$ is non-\'etale iff $mq>2$. It suffices to show that for each $\beta\in k$, by denoting $R:=k_{\t}[[x_1,x_2-\beta]]$ the completion of $k_{\t}[x_1,x_2]$ at the maximal ideal defining the point $(0,\beta)$, the normalization $A^{\n}_{\beta}$ of the $R$-algebra $A_{\beta}:=R[w]/\bigl(g(w)\bigr)$ is isomorphic to the product $R\times R_{mq-1}$, where $R_{mq-1}:=k[[x_1^{\frac{1}{mq-1}},x_2]]$ is regular. 

To check this, we write $g(w)=w^{mq}-w^{mq-1}+\sum_{i=0}^{\frac{mq}{p}-1} u_ix_1^{l_i}(x_2-\beta)^{m_i}w^{pi}$ with each $u_i\in\{0\}\cup R^{\ast}$ and with $(l_i,m_i)\in\mathbb N^2$; this is possible as $\binom{m}{i}w^{q(m-i)}$ is $0$ when $p$ does not divide $qi$ as in this case $p$ divides both $m$ and $\binom{m}{i}$. One computes $l_0=m$ and if $u_i\neq 0$ then $q$ divides $pi$ and $l_i$ is the solution of the equation $q(m-x)=pi$ in $x$, hence $l_i=m-\frac{pi}{q}>l_0\frac{mq-1-pi}{mq-1}$. If $\beta\bigl(1+(-1)^m\beta^{mq-1}\bigr)\neq 0$ then $m_0=0$ and if $u_i\neq 0$ the inequality $m_i\ge\frac{mq-1-pi}{mq-1}m_0$ holds. If $\bigl(1+(-1)^m\beta^{mq-1}\bigr)=0$ (resp.\ $\beta=0$) then $m_0=mq-1$ and if $u_i\neq 0$ then $m_i$ is the value of $qx-1$ computed at the solution of the equation $q(m-x)=pi$ in $x$ and hence $m_i=mq-1-pi=\frac{mq-1-pi}{mq-1}m_0$ (resp.\ $m_i=mq-1\ge\frac{mq-1-pi}{mq-1}m_0$). Thus $mq-1$ divides $m_0$ and the statement on $A^{\n}_{\beta}$ follows from Lemma \ref{L3}(4) if $mq>2$ and from Lemma \ref{L3}(2) if $mq=2$. 

As for each $(\alpha,\beta)\in k^{\ast}\times k$, $\s(\alpha,\beta)\in\{mq-1,mq\}$ is $mq-1$ iff $\beta\in\mathbb I_{p,m,q}^{\ast}$, $\Spec A_{e,x_1}\setminus\Imm(\imath_e)$ is $mq-1=|\mathbb I_{p,m,q}^{\ast}|$ copies of $\mathbb G_{\m,k}$; thus $\rho_{\et}(e)\ge mq-1$. From this, the fact that $\s(0,\beta)\in\{0,1\}$ is $0$ iff $\beta\in\mathbb I_{p,m,q}^{\ast}$, and the above part on complete local rings it follows that $\rho(e)-\rho_{\et}(e)$ is $1$ if $mq>2$ and is $0$ if $mq=2$. Thus $0<mq-1=\rho_{\et}(e)<mq=\rho(e)$ if $mq>2$ and $\rho_{\et}(e)=\rho(e)=2$ if $mq=2$. 

We say more in the following cases.

{\bf Case 1: $m=p$.} As $\binom{p}{i}=0$ if $i\in \llbracket1,p-1\rrbracket$, we compute 
$$g(w)=w^{pq}-w^{pq-1}+x_1^px_2^{pq-1}(1-x_2^{pq-1})^{pq-1}.$$ 
If $q=1$, then we have $h(v)=v^{p}-x_1^{p-2}v+x_1^{p}x_2(1-x_2^{p-1})$, and if $q\ge 2$, then we have $h(v)=v^{pq}-x_1^{pq-p-1}v+x_1^{pq-p}x_2(1-x_2^{pq-1})$.

{\bf Case 2: $m=1$.} So $p|q$ and $\epsilon=1$. As in the prior paragraph we compute $g(w)=w^q-w^{q-1}+x_1x_2^{q-1}(1-x_2^{q-1})^{q-1}$ and $h(v):=v^q-x_1^{q-2}v+x_1^{q-1}x_2(1-x_2^{q-1})$.\end{example}

\subsection{Proof of Corollary \ref{C2}.}\label{S14.1}
To prove Corollary \ref{C2} we can assume that $n=2$ as the passage from $n=2$ to an arbitrary $n>2$ is achieved by adding $n-2$ indeterminates $x_3,\ldots,x_n$ and equations $x_i=0$ for $i\in\{3,\ldots,n\}$. We apply Example \ref{EX11} with $mq=p$; so $\deg(e)=p$ and $\mathbb I_{p,m,q}=\mathbb F_p$. For $n=2$, Corollary \ref{C2}(1) (resp.\ Corollary \ref{C2}(2)) follows from the following identity $\s(0,0)=1$ (resp.\ $\s(0,1)=0$). Concretely, suppose that we have $m=p$ and $q=1$, then we can take $(g_1,g_2)\in\mathbb F_p[x,y]^2\subset\kappa[x,y]^2$ to be $\bigl(x-(x-1)^py^{p-1},y-(x-1)^py^p\bigr)$ (resp.\ $(1+x-x^py^{p-1},y-x^py^p)$) obtained from the Rule (\ref{EQ3.9}) by replacing $x$ with $x-\epsilon$, where $\epsilon$ is $1$ (resp.\ $0$). Similar examples can be obtained based on Example \ref{EX8} by taking $(m,q,s)=(p,1,1)$.

\section{Two general properties of endomorphisms of affine spaces}\label{S15}

We include practical results on $e\in\End_n(K)$ that are often used in what follows and that either estimate $\deg(e)$ if $e\in\QF_n(K)$ or provide a criteria on when $e$ is finite or, for inductive purposes, construct $d\in\End_{n+1}(k)$ from $e$ that has a few iteration properties that are identical to those of $e$. 

\begin{lemma}\label{F3}
Let $e:=\e(g_1,\ldots,e_n)\in\End_n(K)$ with $(g_1,\ldots,g_n)\in K_{\s}[x_1,\ldots,x_n]^n$ be such that $q:=\prod_{i=1}^n \deg(g_i)>0$. For $i\in \llbracket1,n\rrbracket$, let $h_i\in K_{\s}[x_1,\ldots,x_n]$ be the homogeneous component of $g_i$ of degree $\deg(g_i)$. Then the following properties hold.

\medskip
{\bf (1)} If $e\in\D_n(K)$, then $\deg(e)\le q$. 

\smallskip
{\bf (2)} If the only solution in $K^n$ of the homogeneous system of $n$ equations
$$h_1=\cdots=h_n=0$$
is $(0,\ldots,0)$, then $e$ is finite. Moreover, if $e\in\QFE_n(K)$, then $\deg(e)=q$.\end{lemma}

\begin{proof} For $i\in \llbracket1,n\rrbracket$ we write $g_i=\sum_{j=0}^{\deg(g_i)} h_{i,j}$, with each $h_{i,j}\in K_{\s}[x_1,\ldots,x_n]$ homogeneous of degree $j$; so $h_i=h_{i,\deg(g_i)}$. We consider the closed subscheme $W$ of $\mathbb P^n_{K_{\t}[x_1,\ldots,x_n]}$ defined by the system of $n$ homogeneous equations
$$\sum_{j=0}^{\deg(g_i)} h_{i,j}(y_1,\ldots,y_n)y_0^{\deg(g_i)-j}-x_iy_0^{\deg(g_i)}=0,\;\;\;\;i\in \llbracket1,n\rrbracket.$$
Let $\phi:W\rightarrow\mathbb A^n_K$ be the resulting natural morphism. As $e\in\D_n(K)$, for a generic point $P\in K^n$ we have $|e^{-1}(P)|=\deg(e)$; for such a point, $\phi^{-1}(P)_{\red}$ has $|e^{-1}(P)|$ irreducible components of dimension $0$ in $\mathbb A^n_{K_{\t}[x_1,\ldots,x_n]}$ and possibly other irreducible components in $\mathbb P^{n-1}_{K_{\t}[x_1,\ldots,x_n]}=\mathbb P^n_{K_{\t}[x_1,\ldots,x_n]}\setminus \mathbb A^n_{K_{\t}[x_1,\ldots,x_n]}$ and hence $|e^{-1}(P)|\le |\phi^{-1}(P)_{\red}|$. From B\'ezout's inequality (see \cite{Ful}, Ch.\ 8, Ex.\ 8.4.6) we get that$ |\phi^{-1}(P)_{\red}|\le q$. So $\deg(e)\le q$, i.e., part (1) holds. 

For part (2), due to the hypotheses, $W\cap\mathbb P^{n-1}_{K_{\t}[x_1,\ldots,x_n]}=\emptyset$. So $W$ is $\mathbb A^n_{K,\s}$, i.e., is the spectrum of $K_{\t}[x_1,\ldots,x_s,y_1,\ldots,y_s]\bigl(g_1(y_1,\ldots,y_n)-x_1,\ldots,g_n(y_1,\ldots,y_n)-x_n\bigr)$. Hence $e$ is projective; being also affine, it is finite. If $e\in\QFE_n(K)$, then for a generic $P$, $\phi^{-1}(P)=e^{-1}(P)\cong \Spec k^{\deg(e)}$, so $\deg(e)=q$ by B\'ezout's Theorem (see \cite{Ful}, Ch.\ 8, Prop.\ 8.4; see also \cite{H2}, Ch.\ I, Thm.\ 7.7 applied $n-1$ times).\end{proof}

\begin{definition}\label{D7+}
Let $e\in\End_n(K)$ and $d\in\End_{n+1}(K)$.

\medskip
{\bf (1)} We say that $d$ is an extension\index{extension} of $e$ if $d(P,0)=\bigl(e(P),0\bigr)$ for all $P\in K^n$.

\smallskip
{\bf (2)} We say that $d$ is a fiberwise extension\index{extension!fiberwise} or a deformation of $e$ if we have $d(K^n\times\{\alpha\})\subset  K^n\times\{\alpha\}$ for each $\alpha\in K$. We cal the endomorphism $d_{\alpha}\in\End_n(k)$ defined by the identity $(d_{\alpha}(P),\alpha)=d(P,\alpha)$ for $P\in K^n$ as the fiber\index{extension!fiber of a fiberwise extension} of $d$ at $\alpha\in K$ (so $e=d_0$).

\smallskip
{\bf (3)} We say that $d$ is a fiberwise finite extension\index{extension!fiberwise finite} or a deformation by finite endomorphisms of $e$ if it is a fiberwise extension such that $d_{\alpha}\in\End_n(k)$ is finite for each $\alpha\in K^{\ast}$.
\end{definition}

\begin{lemma}\label{L13}
For $(g_1,\ldots,g_n)\in K[x_1,\ldots,x_n]^n$ let $e=\e(g_1,\ldots,g_n)\in\End_n(K)$ and for $(f_1,\ldots,f_{n+1})\in K[x_1,\ldots,x_{n+1}]^{n+1}$ let $d=\e(f_1,\ldots,f_{n+1})\in\End_{n+1}(K)$. Then the following properties hold.

\medskip
{\bf (1)} The endomorphism $d$ is an extension of $e$ iff $x_{n+1}$ divides $f_{n+1}$ and $g_i-f_i$ for each $i\in \llbracket1,n\rrbracket$.

\smallskip
{\bf (2)} The endomorphism $d$ is a fiberwise extension of $e$ iff $f_{n+1}=x_{n+1}$ and $x_{n+1}$ divides $g_i-f_i$ for each $i\in \llbracket1,n\rrbracket$.

\smallskip
{\bf (3)} Let $(q_1,\ldots,q_n)\in (\mathbb N^{\ast})^n$ be such that $q_i>\deg(g_i)$ for each $i\in \llbracket1,n\rrbracket$. Let $(h_1,\ldots,h_n)\in K[x_1,\ldots,x_n]^n$ be such that $h_i$ is homogeneous of degree $q_i$ for each $i\in \llbracket1,n\rrbracket$. Let $(L_1,\ldots,L_n)\in (K[x_1,\ldots,x_{n+1}]x_{n+1})^n$ be such that the total degree of $L_i$ in $x_1,\ldots,x_n$ is less than $q_i$ for each $i\in \llbracket1,n\rrbracket$. Let $(l_1,\ldots,l_n)\in (\mathbb N^{\ast})^n$. Suppose that $f_i=g_i+h_ix_{n+1}^{l_i}+L_i$ for each $i\in \llbracket1,n\rrbracket$, $f_{n+1}=x_{n+1}$, and the only solution in $K^n$ of the homogeneous system of $n$ equations $h_1=\cdots=h_n=0$ is $(0,\ldots,0)$ (e.g., if $(h_1,\ldots,h_n)=(x_1^{q_1},\ldots,x_n^{q_n})$). Then $d\in D_{n+1}(K)$ and $d$ is a fiberwise finite extension of $e$ with $\deg(d)\le\prod_{i=1}^n q_i$. If moreover $\mathcal K$ is a subfield of $K$ such that $e$ is defined over $\mathcal K$, then $d$ is defined over $\mathcal K$ provided $(h_1,\ldots,h_n)\in \mathcal K[x_1,\ldots,x_n]^n$ and $(L_1,\ldots,L_n)\in (\mathcal K[x_1,\ldots,x_{n+1}]x_{n+1})^n$.

\smallskip
{\bf (4)} If $d$ is a fiberwise finite extension of $e$, then for each $m\in\mathbb N^{\ast}$ the sets $\Gamma_{e^m}$ and $\Gamma_{d^m}$ are compatibly isomorphic under isomorphisms induced by the projection $\mathbb A^{n+1}_{K,\t}\rightarrow \mathbb A^n_{K,\t}$ on the first $n$ coordinates and hence we have $\iota(d)=\iota(e)$. 

\smallskip
{\bf (5)} If there exists $s\in\mathbb N^{\ast}$, $(e_1,\ldots,e_s)\in\End_n(K)^s$, and $(d_1,\ldots,d_s)\in\End_{n+1}(K)^s$ such that $e=e_1e_2\cdots e_s$ and $d=d_1d_2\cdots d_s$ and $d_i$ is an extension (resp.\ a fiberwise extension or fiberwise finite extension) of $e_i$ for each $i\in \llbracket1,s\rrbracket$, then $d$ is an extension (resp.\ a fiberwise extension or fiberwise finite extension) of $d$.

\smallskip
{\bf (6)} Suppose that $\chr(K)=p$. If $e\in\pMor_n(K)$, then there a fiberwise finite extension of $e$ in $\pMor_{n+1}(K)$.

\smallskip
{\bf (7)} Suppose that $\chr(K)=p$. If $e\in\EE_n(K)$ and if in part (3) for each $i\in \llbracket1,n\rrbracket$ we have $q_i\in p\mathbb N^{\ast}$, $h_i\in K[x_1^p,\ldots,x_n^p]$, and $L_i\in K[x_1^p,\ldots,x_{n+1}^p]$, then $d$ is a fiberbise finite extension of such that $d\in\EE_{n+1}(K)$, $\deg(d)=\prod_{i=1}^n q_i$ and $\mathcal F(d)=\mathcal F(e)\cup\{\deg(e)\}$. If moreover, $e\in\BE_n(K)$ and $l_i\in p\mathbb N^{\ast}$ for each $i\in \llbracket1,n\rrbracket$, then $d\in\BE_{n+1}(K)$.
\end{lemma}

\begin{proof} Parts (1), (2), and (5) are clear from definitions. 

For part (3), $d$ is a fiberwise extension of $e$ by part (2) and its fiber at each $\alpha\in K^{\ast}$ is defined by the $n$-tuple 
$$(g_1+\alpha^{l_1}h_1+f_1(x_1,\ldots,x_n,\alpha),\ldots,g_n+\alpha^{l_n}h_1+f_n(x_1,\ldots,x_n,\alpha)$$ 
and hence it is is finite by Lemma \ref{F3}(2) with $\deg(d_{\alpha})\le\prod_{i=1}^n q_i$. Hence $d$ is a fiberwise finite extension of $e$. As $\Imm\bigl(d(K)\bigr)$ contains $K^n\times K^{\ast}$, we have $d\in\D_{n+1}(K)$. As for each $(P,\alpha)\in K^n\times K^{\ast}$ we have $|d^{-1}(P,\alpha)|=|d_{\alpha}^{-1}(P)|\le\deg(d_{\alpha})\le\prod_{i=1}^n q_i$, we have $\deg(d)\le\prod_{i=1}^n q_i$. The last sentence of part (3) on $\mathcal K$ is clear, so part (3) holds.

For part (4), for each $(m,\alpha)\in\mathbb N\times K$ we have we have $(d^m)_{\alpha}=(d_{\alpha})^m$. If $\alpha\in\ K^{\ast}$, then $d_{\alpha}$ is surjective as it is finite. Hence $\Gamma_{d^m}=\Gamma_{e^m}\times\{0\}$ and thus the projection $\mathbb A^{n+1}_{k,\t}\rightarrow \mathbb A^n_{k,\t}$ on the first $n$ coordinates induces bijections $\Gamma_{d^m}\rightarrow \Gamma_{e^m}$ compatibly. For each $m\in\mathbb N^{\ast}$ we have $\Imm(d^m)=\Imm(d^{m+1})\Leftrightarrow \Imm(e^m)=\Imm(e^{m+1})$, hence $\iota(d)=\iota(e)$. So part (4) holds.

Based on part (5) and Definition \ref{D4.3}, to prove part (6) we can assume that either (i) $e\in\GA_n(K)$ or (ii) there exists a polynomial $f\in K[x_1,\ldots,x_n]$ such that $(g_1,\ldots,x_n)=\bigl(x_1,\ldots,x_n+f(x_1,\ldots,x_{n-1},x_n^p)\bigr)$. If (i) holds, i.e., $e\in\GA_n(K)$, then $d=e\times 1_{\mathbb A^1_K}\in\GA_{n+1}(K)\subset \pMor_{n+1}(K)$ is a fiberwise finite deformation of $e$. If (ii) holds, let $l\in\mathbb N^{\ast}$ be such that $pl>\deg(f)$; as in part (3) we argue that 
$$d_f:=\e\bigl((x_1,\ldots,x_n+f(x_1,\ldots,x_{n-1},x_n^p)+x_{n+1}x_n^{p^l},x_{n+1})\bigr)$$ is a fiberwise finite deformation of $e$. If $a:=\e(x_1,\ldots,x_{n-1},x_{n+1},x_n)\in\GA_{n+1}(K)$, then $ad_fa$ is an elementary $p$-morphism, thus $d_f\in\pMor_{n+1}(K)$. So part (6) holds.

For part (7), we have $d\in\EE_{n+1}(K)$ and $d_{\alpha}\in\EE_n(K)$ for reach $\alpha\in K^{\ast}$. If $e\in\BE_n(K)$ and $l_i\in p\mathbb N^{\ast}$ for each $i\in \llbracket1,n\rrbracket$, then $d\in\BE_{n+1}(K)$ and $d_{\alpha}\in\BE_n(K)$ for each $\alpha\in K^{\ast}$. For each $\alpha\in K^{\ast}$ we have $\deg(d_{\alpha})=\prod_{i=1}^n q_i$ by Lemma \ref{F3}(2). As $\prod_{i=1}^n q_i$ is greater that $\prod_{i=1}^n \deg(g_i)$ and hence than $\deg(e)$ by Lemma \ref{F3}(1), it follows that $\deg(d)=\prod_{i=1}^n q_i$ and that $\mathcal F(d)$ is the set of cardinals of fibers of $e$. Thus $\mathcal F(d)=F(e)\cup\{\deg(e)\}$. So part (7) holds.\end{proof} 

\begin{example}\normalfont\label{EX3+}
Let $e:=e\bigl(x_1+x_1^{p^m}-x_2^{p^m(p^m+1)},x_2+x_1^{p^m}-x_2^{p^m(p^m+1)}\bigr)\in\BE_2(k)$ with $m\in\mathbb N$, $a:=\e(x_1+x_2^{p^m+1},x_2)\in\GA_2(k)$, and $d:=ea$. We compute that $\break d=\e(x_1+x_1^{p^m}+x_2^{p^m+1},x_2+x_1^{p^m})$. Therefore $d\in\BE_2(k)$, we have inequalities $\pi(d)\le\min\bigl(\pi_{\i}(d),\pi_{\l}(d),\pi_{\r}(d)\bigr)\le\max\bigl(\pi_{\i}(d),\pi_{\l}(d),\pi_{\r}(d)\bigr)\le p^m+1$, and Lemma \ref{F3}(2) gives that $d$ is finite with $\deg(d)=p^m(p^m+1)$. If $\pi(d)\le p^m$, i.e., there exists $(b,c)\in\GA_2(k)^2$ such that $bdc=\e(g_1,g_2)$ with $(g_1,g_2)\in k_{\s}[x_1,x_2]^2$ satisfying $\pi(g_1,g_2)\le p^m$, then $p^m(p^m+1)=\deg(d)=\deg(bdc)\le\deg(g_1)\deg(g_2)\le (p^m)^2$ by Lemma \ref{F3}(1), a contradiction. So $p^m+1=\pi(d)=\pi_{\i}(d)=\pi_{\l}(d)=\pi_{\r}(d)$. As $\pi(e)=\pi(d)$ and $\pi_{\r}(e)=\pi_{\r}(d)$, it follows that 
$$p^m+1=\pi_{\r}(e)=\pi(e)\le\min\bigl(\pi_{\l}(e),\pi_{\i}(e)\bigr).$$\end{example}

\section{On iterate inner invariants}\label{S16}

We first include examples and applications of Section \ref{S3} using Notation \ref{NOT1}.

\begin{example}\normalfont\label{EX12}
Let $e:X\rightarrow X$ be an endomorphism of a variety over $K$ such that $\varphi(e)\in\mathbb N^{\ast}$. If $\varsigma_e\ge 2$, then $\varphi(e^{l+1}\bigr)<\frac{\varsigma_e\varphi(e)}{\varsigma_e-1}\le 2\varphi(e)$ by Lemma \ref{L1}(7) and therefore $\iota(e)\le\max\bigl(1,\frac{\varphi(e)}{\varsigma_e-1}\bigr)\le \varphi(e)$ and $\varphi\bigl(e^{\iota(e)}\bigr)<\frac{\varsigma_e\varphi(e)}{\varsigma_e-1}\le 2\varphi(e)$. Similarly, if there exists $s\in\mathbb N^{\ast}$ such that $\varsigma_{e^s}\ge 2$, then the inequalities $\iota(e^s)\le\frac{\varphi(e^s)}{\varsigma_{e^s}-1}\le \frac{s\varphi(e)}{\varsigma_{e^s}-1}$ (the second one by Lemma \ref{L1}(2)) imply that $\iota(e)\le s\iota(e^s)\le\max\bigl(s,\frac{s^2\varphi(e)}{\varsigma_{e^s}-1}\bigr)$, thus $\iota(e)\in\mathbb N$.
\end{example}

\begin{corollary}\label{C3}
Let $e:X\rightarrow X$ be an endomorphism of a variety over $k$ such that $\Gamma_e$ is finite. We assume that there exists $q\in\mathbb N^{\ast}$ such that $e$ is the pullback of an endomorphism $d:Y\rightarrow Y$ over $\mathbb F_{p^q}$. Then $\iota(e)\in\mathbb N$. More precisely, if we choose $q$ and $d$ such that $\Gamma_e\subset Y(\mathbb F_{p^q})$, then $\iota(e)\le |Y(\mathbb F_{p^q})\setminus\Gamma_e|$.
\end{corollary}

\begin{proof}
As $\Gamma_{e(k)}^+\subset Y(\mathbb F_{p^q})$ is finite, this follows from Lemma \ref{L1}(6).
\end{proof}

\begin{example}\normalfont\label{EX13}
Let $e\in\QF_2(k)$ be non-surjective (i.e., $0\in\mathcal F(e)$ or $\varphi(e)>0$). If $e$ and $\Gamma_e$ are defined over $\mathbb F_{p^q}$, then $\iota(e)\le p^{2q}-\varphi(e)$ by Corollary \ref{C3} and hence $\varphi\bigl(e^{\iota(e)}\bigr)\le\varphi(e)\iota(e)\le\varphi(e)[p^{2q}-\varphi(e)]$ by Lemma \ref{L1}(2).
\end{example}

For an endomorphism $e:\mathcal B\rightarrow\mathcal B$ of a set $\mathcal B$ and $P\in\mathcal B$ we consider two iteration conditions\index{iteration condition}.

\medskip
\phantomsection{($\natural_{e,P}$)\index{iteration condition!($\natural_{e,P}$)} {\it The sequence $\bigl(e^m(P)\bigr)_{m\ge 0}$ has distinct elements and for each $(i,j)\in\mathbb N^2$, we have $|e^{-i}\bigl(e^j(P)\bigr)|>0$ iff $i\le j$ and iff $e^{-i}\bigl(e^j(P)\bigr)=\{e^{j-i}(P)\}$.}}\label{PH52}

\medskip
($\flat_{e,P}$)\index{iteration condition!($\flat_{e,P}$)} {\it For each $m\in\mathbb N$, $e^m(P)\notin\Imm(e^{m+1})$.}

\begin{notation}\normalfont\label{NOT4}
For a non-empty constructible subset $V$ of a variety $W$ over $K$, let $\overline{V}$ be its Zariski closure in $W$, 
$$\digamma(V)=\digamma(\bar{V}):=|\{Y\in\Irr(\bar{V})|\dim(Y)=\dim(\overline{V})\}|,$$ 
and 
$$\hbar(V)=\hbar(\bar{V}):=\bigl(\dim(\bar{V}),\digamma(\overline{V})\bigr)\in\mathbb N^2.$$ If $W$ is irreducible, let $\eta_W\in W$ be its generic point.
\end{notation}

Endowing $\mathbb N^2$ with the lexicographic order, for a dominant morphism $W_1\rightarrow W_2$ of varieties over $K$ we have the inequality $\hbar(W_1)\ge\hbar(W_2)$.

\begin{lemma}\label{L14}
Let $e:X\rightarrow X$ be an endomorphism of a variety over $K$. For $m\in\mathbb N$, let $e_m:\Imm(e^m)\rightarrow\Imm(e^m)$ be the endomorphism of constructible sets induced by $e$. Then the following statements are equivalent.

\medskip
{\bf (1)} We have $\iota(e)=\infty$.

\smallskip
{\bf (2)} There exists $m\in\mathbb N$ and a point $P\in\Imm(e^m)$ such that ($\natural_{e_m,P}$) holds.

\smallskip
{\bf (3)} There exists $m\in\mathbb N$ and a point $P\in\Imm(e^m)$ such that ($\flat_{e_m,P}$) holds.

\smallskip
{\bf (4)} There exists an algebraically closed field $K_1$ that contains $K$, an $m\in\mathbb N$, and a point $P_1\in\Imm(e^m)(K_1)$ such that ($\flat_{(e_m)_{K_1},P_1}$) holds.

\medskip
If $\Gamma_e$ is finite, then the statements are also equivalent to the following two.

\smallskip
{\bf (5)} There exists $m\in\mathbb N$ and $P\in\Imm(e^m)(K)$ such that ($\natural_{e_m,P}$) holds.

\smallskip
{\bf (6)} There exists $m\in\mathbb N$ and $P\in\Imm(e^m)(K)$ such that ($\flat_{e_m,P}$) holds.
\end{lemma}

\begin{proof}
The implications $(2)\Rightarrow (3)$ and $(5)\Rightarrow (6)$ always hold. Clearly, $(6)\Rightarrow (3)$.

We show that $(1)\Rightarrow (2)$. For $i\in\mathbb N$ the constructible set $Z_i:=\Imm(e^i)\setminus\Imm(e^{i+1})$ (by Chevalley's Theorem, see \cite{Gro2}, Thm.\ (1.8.4)) is non-empty. Let $Y_i$ be the Zariski closure of $Z_i$ in $X$. As $e$ induces surjective morphisms $Z_i\rightarrow Z_{i+1}$, it also induces dominant morphisms $Y_i\rightarrow Y_{i+1}$. As $\bigl(\hbar(Y_n)\bigr)_{n\ge 0}$ is a non-increasing sequence of the well ordered $\mathbb N^2$, there exists $(m,q,l)\in\mathbb N^3$ such that $\hbar(Y_i)=\hbar(Y_m)=(q,l)$ for each integer $i\ge m$. For $i\ge m$ let $\eta_i:=\{\eta_Y|Y\in\Irr(Y_i),\dim(Y)=q\}$; so $|\eta_i|=l>0$ and $\eta_i\subset Z_i$. An induction on $j\in\mathbb N$ gives that $e^j(\eta_m)=\eta_{m+j}$. 

Let $P\in\eta_m$ and $(i,j)\in\mathbb N^2$. As $P\in Z_m$, we have $|e^{-i}\bigl(e^j(P)\bigr)|>0$ iff $i\le j$. If $i\le j$, then the inclusion $e^{-i}\bigl(e^j(P)\bigr)\supset\{e^{j-i}(P)\}$ is clear and the reversed inclusion $e^{-i}\bigl(e^j(P)\bigr)\subset\{e^{j-i}(P)\}$ follows from $\hbar(Y_{m+j-i})=\hbar(Y_m)$. Thus ($\natural_{e_m,P}$) holds. 

If $\Gamma_e$ is finite, then $q=0$ and $\eta_m\subset X(K)$, so $(1)\Rightarrow (5)$ by previous paragraph.

To show that $(3)\Rightarrow (4)$, let $K_1$ be an uncountable algebraically closed field that contains the residue field of $P$. If $\overline{P}$ is the Zariski closure in $X$ of $P$, then $\overline{P}_{K_1}$ is not the union of its proper subvarieties $\overline{P}_{K_1}\cap e_{m,K_1}^{-l}\bigl(\Imm(e_{m,K_1}^{l+1})\bigr)$ indexed by $l\in\mathbb N$ and hence there exists a point $P_1\in\overline{P}(K_1)$ such that ($\flat_{e_{m,K_1},P_1}$) holds. Thus $(3)\Rightarrow (4)$.

We show that $(4)\Rightarrow (1)$. As $\iota(e)=\iota(e_{K_1})$ and $\iota(e)=\infty\Leftrightarrow \iota(e_m)=\infty$, we can assume that $K=K_1$ and $m=0$. As ($\flat_{e,P}$) holds, we have $\Imm(e^{l})\neq\Imm(e^{l+1})$ for all $l\in\mathbb N$, hence $\iota(e)=\infty$.\end{proof}

\begin{remark}\normalfont\label{R10}
Suppose that $\Gamma_e$ is finite and $\iota(e)=\infty$. We use Lemma \ref{L1} to show that statement (5) of Lemma \ref{L14} holds (i.e., to give a second proof that in Lemma \ref{L14} we have $(1)\Rightarrow (5)$). Let $(l,N)\in\mathbb N^2$ be such that $\varphi(e^{n+1})-\varphi(e^n)=l$ for all integers $n\ge N$ by Lemma \ref{L1}(1). By replacing $e$ with $e_N$ we can assume that $N=0$; so $\varphi(e^n)=nl$ for all $n\in\mathbb N$, with $l\in\mathbb N^{\ast}$ as $\iota(e)=\infty$. For $n\in\mathbb N$, as $|\cup_{i=0}^n e^i(\Gamma_e)|\le l(n+1)=\varphi(e^{n+1})=|\Gamma_{e^{n+1}}|$, we get that $\Gamma_{e^{n+1}}=\sqcup_{i=0}^n e^i(\Gamma_e)$ by Lemma \ref{L1}(5); so $(\natural_{e,P}$) holds for each $P\in\Gamma_e\subset X(K)$.\end{remark}

\begin{remark}\normalfont\label{R11}
Let $e:X\rightarrow X$ be an endomorphism of a non-empty variety over $K$ which is dominant (resp.\ which is dominant outside a closed subvariety of the target of dimension $m<\dim(X)$). A decreasing induction on $i\in\llbracket0,\dim(X)\rrbracket$ (resp.\ $i\in\llbracket m+1,\dim(X)\rrbracket$) gives that $e$ permutes the generic points of irreducible components of $X$ of dimension $i$. So there exists $N\in\mathbb N^{\ast}$ such that $e^N$ fixes the generic points of all irreducible components (resp.\ irreducible components of dimension $>m$) of $X$. If $d:X\rightarrow X$ is quasi-finite, then $1\le\digamma\bigl(\Imm(d)\bigr)\le\digamma(X)$; so there exists $i\in \llbracket0,\digamma(X)-1\rrbracket$ such that for $j\in\mathbb N$ we have 
$$\digamma\bigl(\Imm(d^{i+j})\bigr)=\digamma\bigl(\Imm(d^i)\bigr)\in\mathbb N.$$
\end{remark}

Next theorem was first obtained in \cite{PNCG}, Thm.\ 1 under the extra hypotheses that $K=\mathbb C$, $X$ is affine, $e$ is universally open, and $\Gamma_e$ is finite.

\begin{theorem}\label{T4}
Let $e:X\rightarrow X$ be a quasi-finite endomorphism of a variety over $K$. Then the following properties hold.

\medskip
{\bf (1)} We have $\iota(e)\in\mathbb N$.

\smallskip
{\bf (2)} There exists no $P\in X$ for which Condition $(\flat_{e,P})$ (or $(\natural_{e,P})$) holds.\end{theorem}

\begin{proof}
Part (2) follows from part (1) and Lemma \ref{L14}. 

To prove part (1), let $\overline{Z}$ be the Zariski closure in $X$ of a subset $Z$ of $X$. The decreasing sequence $\bigl(\overline{\Imm(e^m)}\bigr)_{\m\in\mathbb N}$ is stationary. Therefore, by replacing $X$ with $W:=\cap_{m\ge 0} \overline{\Imm(e^m)}$, we can assume that $X=\overline{\Imm(e)}$ (i.e., $e$ is dominant).\footnote{Note that if $l\in\mathbb N$ is such that $W=\overline{\Imm(e^l)}$, then $\Imm(e^{m+l})\subset e^m(W)\subset \Imm(e^m)$, thus the decreasing sequence $\bigl(\Imm(e^m)\bigr)_{\m\in\mathbb N}$ is stationary iff the decreasing sequence $\bigl(e^m(W)\bigr)_{\m\in\mathbb N}$ is so.}

By replacing $e$ with $e^N$ for some $N\in\mathbb N^{\ast}$, we can assume that $e$ fixes each generic point of any $Y\in\Irr(X)$ (see Remark \ref{R11}); so $e$ induces an endomorphism $e_Y:Y\rightarrow Y$. If $\iota(e_Y)\in\mathbb N$ for each $Y\in\Irr(X)$, then $\iota(e)\in\mathbb N$ as in fact we have $\iota(e)\le \max\bigl(e_Y|Y\in\Irr(X)\bigr)$. Thus, by replacing $e$ with the normalizations of any $e_Y$ we can assume that $X$ is normal irreducible.

To end the proof it suffices to show that the assumption that $\iota(e)=\infty$ leads to a contradiction. From Lemma \ref{L14} we get that, by replacing $e$ with the endomorphism it induces on the open subvariety $\Imm(e^m)$ for some $m\in\mathbb N$, we can assume that there exists $P\in X$ such that Condition ($\flat_{e,P}$) holds.

Let $R$ be a finitely generated $\mathbb Z$-subalgebra of $K$ such that $e$ is defined over $R$ (see \cite{Gro4}, Prop.\ (8.9.1)(iii)) and the Zariski closure $Z_P$ of $P$ in $X$ is defined over $\mathcal K:=\Frac(R)$; let $e_R:X_R\rightarrow X_R$ be an endomorphism of a $\Spec R$-scheme $X_R$ of finite type such that its pullback to $\Spec K$ is $e$. The fiber $e_0:X_0\rightarrow X_0$ of $e_R$ over $\mathcal K$ is quasi-finite and there exists $P_0\in X_0$ with the property that the pullback of its Zariski closure in $X_0$ to $X$ is $Z_P$. As ($\flat_{e,P}$) holds, ($\flat_{e_0,P_0}$) also holds. By localizing $R$ we can assume that $X_R$ is normal irreducible and (see \cite{Gro4}, Thm.\ (8.10.5)(xi)) that $e_R$ is quasi-finite; so $e_R$ is universally open by \cite{Gro4}, Cor.\ (14.4.3). 

Let $Q$ be a closed point of the Zariski closure of $P_0$ in $X_0$; its residue field $\kappa$ is finite and we have an inclusion $\{e^m(Q)|m\in\mathbb N\}\subset X_0(\kappa)$ of finite sets. Let the pair $(q,s)\in\mathbb N\times\mathbb N^{\ast}$ be such that $e_0^q(Q)=e_0^{q+s}(Q)$. As $e_0^q(Q)$ belongs to the open subscheme $\Imm(e_0^{q+s})$ of $X_0$, $e_0^q(P_0)\in\Imm(e_0^{q+s})\subset\Imm(e_0^{q+1})$, which contradicts ($\flat_{e_0,P_0}$).\end{proof}

\begin{corollary}\label{C4}
Let $f:W\rightarrow Y$ be a morphism of finite presentation between quasi-compact schemes. Let $e:W\rightarrow W$ be an endomorphism of schemes over $Y$ (i.e., $f=f\circ e$) with the property that all its fibers $e_y:W_y\rightarrow W_y$ over geometric points $y$ of $Y$ are quasi-finite. Then there exists $n\in\mathbb N$ such that $\Imm(e^n)=\Imm(e^{n+1})$.
\end{corollary}

\begin{proof}
For $m\in\mathbb N$, as $\Imm(e^m)$ is a constructible subset of $W$, from \cite{Gro4}, Cor.\ (9.5.2) we get that $Y_m:=\{y\in Y|\Imm(e^m)_y=\Imm(e^{m+1})_y\}$ is a constructible subset of $Y$. As $\Imm(e^m)_y=\Imm(e_y^m)$, it follows first that $Y_m\subset Y_{m+1}$ for all $m\in\mathbb N$ and second that $\cup_{m\ge 0} Y_m=Y$ by Theorem \ref{T4}. So there exists $n\in\mathbb N$ such that $Y_n=Y$ by \cite{Gro2}, Cor.\ (1.9.9); hence $\Imm(e^n)=\Imm(e^{n+1})$.\end{proof}

Our examples of endomorphisms $e:X\rightarrow X$ with $\iota(e)=\infty$ are in essence based on the notion `trap' used in \cite{Bori} and formalized here as follows.

\begin{definition}\label{D8} Let $e:X\rightarrow X$ be an endomorphism of a non-empty variety over $K$. By a trap\index{trap} of $e$ we mean a non-empty closed subvariety $Y$ of $X$ such that $\dim(Y)<\dim(X)$, $e(Y)\subset Y$, the endomorphism $e_Y:Y\rightarrow Y$ induced by $e$ is dominant, and we have a strictly increasing sequence 
$$\bigl(|\{Z\in\Irr(e^{-n}(Y))|\dim(Z)>\dim(Y)\;\textup{or}\;\dim(e^n(Z))<\dim(Z)=\dim(Y)\}|\bigr)_{n\in\mathbb N^{\ast}}.$$\end{definition}

\begin{example}\normalfont\label{EX14}
We list several examples of endomorphisms in $\D_n(K)\setminus\QF_n(K)$ with $n\in\{2,3\}$ whose iterate inner invariants are $\infty$ and often mention zero loci that are traps of them.

\medskip
{\bf (1)} The rule $(x,y)\mapsto(xy,xy+y)$ defines $e\in\D_2(K)\setminus\QF_2(K)$ which is birational (i.e., generically \'etale with $\deg(e)=1$); the zero locus $x_1=x_2=0$ is a trap of it of dimension $0$. If $K_1$ is an algebraically closed field that contains $K$ and a transcendental element $\alpha$ over the prime field, then ($\natural_{e_{K_1},(\alpha,\alpha)}$) holds, therefore $\iota(e)=\iota(e_{K_1})=\infty$ by Lemma \ref{L14}. 

\smallskip
{\bf (2)} If $\chr(K)=0$ and $\beta=1$ (resp.\ $\chr(K)=p$ and there exists a transcendental element $\beta\in K$ over $\mathbb F_p$), then the rule $(x,y)\mapsto(xy,\beta y+1)$ defines $d\in\D_2(K)\setminus\QF_2(K)$ which is birational with $\iota(d)=\infty$; the zero locus $x_1=0$ is a trap of it of dimension $1$.

\smallskip
{\bf (3)} Similarly to example (1), either one of the rules $(x,y)\mapsto (x^py^p,y^p+x^py^p)$ and $(x,y)\mapsto\bigl(x(x-y)^{p-1},y^p(x-y)^{p^2-p}\bigr)$ defines an $e\in\D_2(k)\setminus\QF_2(k)$ which is generically radicial (i.e., $\deg(e)=1$) and defined over $\mathbb F_p$ with $\iota(e)=\infty$.

\smallskip
{\bf (4)} Let $(l,m)\in (\mathbb N^{\ast})^2$. Let $C_m\subset K^2$ be the empty set if $m=1$ and be $\{0\}\times K^\ast$ if $m\ge 2$. The rule $(x,y)\mapsto\bigl(x^my,x^{lm}y^l+x^{l(m-1)}y^{l+1}\bigr)$ defines an endomorphism $e_{\a,l,m}\in\D_2(K)\setminus\QF_2(K)$. If either $x=0$ and $m\ge 2$ or $y=0$, then $(x,y)\mapsto (0,0)$. If $m=1$, then $(0,y)\mapsto (0,y^{l+1})$. If $xy\neq 0$, then we have the additive property
\begin{equation}\label{EQ5+}
\frac{x^{lm}y^l+x^{l(m-1)}y^{l+1}}{(x^my)^l}=1+\frac{y}{x^l}.
\end{equation} 
If $\chr(K)\nmid l+m$ and $(\alpha,\beta)\in (K^{\ast})^2$ with $\beta\neq\alpha^l$, then the system of equations 
$$x^my-\alpha=0=x^{lm}y^l+x^{l(m-1)}y^{l+1}-\beta$$ 
gives based on Equation (\ref{EQ5+}) that $\frac{\beta}{\alpha^l}=\frac{y}{x^l}+1$, hence $y=\frac{(\beta-\alpha^l)x^l}{\alpha^l}$; substituting this in the first equation we get that $x^{l+m}=\frac{\alpha^{l+1}}{\beta-\alpha^l}$, hence the system of equations has $l+m$ solutions, therefore $\deg(e_{\a,l,m})=l+m$ and $e_{\a,l,m}$ is generically \'etale. If $\chr(K)=0$, then $\deg(e_{\a,l,m})=l+m\ge 2$, $\iota(e_{\a,l,m})=\infty$ and the zero locus $x_1=x_2=0$ is a trap of $e_{\a}$ of dimension $0$; more precisely, we have $\iota(e_{\a,l,m})=\infty$ as for $n\in\mathbb N^{\ast}$, one computes
$$\Gamma_{e_{\a,l,m}^n}=\{(x,y)\in (K^{\ast})^2|x^{-l}y\in \llbracket1,n\rrbracket\}\cup C_m.$$
The case $(l,m)=(1,2)$ is the ``Additive Trap" of \cite{Bori}, Ex.\ 1.

\smallskip
{\bf (5)} The multiplicative analog of part (4) works in all characteristics. Let $(l,m)$ and $C_m$ be as in part (4). Let $(q,s)\in\mathbb N^2$. Let $h(t)\in K[t]$ be such that $\deg(h)=s$ and the derivative of $t^{l(q+1)+m}h(t)$ is non-zero. Assume there exists $\delta\in K^{\ast}$ of infinite multiplicative order. Let 
$$e:=e_{\m,l,m,q,s,h,\delta}\in\D_2(K)\setminus\QF_2(K)$$ 
be defined by the rule 
$$(x,y)\mapsto \bigl(x^my^q(x^l-y)h(x),\delta x^{l(m-1)}y^{ql+1}(x^l-y)^lh(x)^l\bigr).$$ If either $(x^l-y)h(x)=0$ or $x=0$ and $m\ge 2$ or $y=0$ and $q\ge 1$, then we have $e(x,y)= (0,0)$. If $m=1$, then we have $e(0,y)=\bigl(0,(-1)^l\delta h(0)^ly^{(q+1)l+1}\bigr)$. If $q=0$, then we have $e(x,0)=(x^mh(x),0)$. If $xy(x^l-y)h(x)\neq 0$, then we have the multiplicative property
\begin{equation}\label{EQ5++}
\frac{\delta x^{l(m-1)}y^{lq+1}(x^l-y)^lh(x)^l}{[x^my^q(x^l-y)h(x)]^l}=\delta\frac{y}{x^l}.
\end{equation} 
For a generic $(\alpha,\beta)\in (K^{\ast})^2$ with $\beta\neq\alpha^l$ we consider the system of equations 
$$x^my^q(x^l-y)h(x)-\alpha=0=\delta x^{l(m-1)}y^{ql+1}(x^l-y)^lh(x)^l-\beta.$$ From Equation (\ref{EQ5++}) we get that $\delta\frac{y}{x^l}=\frac{\beta}{\alpha^l}$, hence $y=\frac{\beta x^l}{\alpha^l}$. Substituting this in the first equation we get that $x^{l(q+1)+m}h(x)=\frac{\alpha^{(q+1)l+1}}{(\alpha^l-\beta)\beta^q}$. From this and the assumptions $[t^{l(q+1)+m}h(t)]'\neq 0$ and $(\alpha,\beta)$ generic, we get that the system has $l(q+1)+m+s$ solutions. Thus $\deg(e)=l(q+1)+m+s$ and $e$ is generically \'etale. The zero locus $x_1=x_2=0$ is a trap of $e$ of dimension $0$. We have $\iota(e)=\infty$ as for $n\in\mathbb N^{\ast}$ we have an inclusion
$$\Gamma_{e^n}(K)\supset \bigl\{(x,y)\in (K^{\ast})^2|x^{-l}y\in\{\delta^l|l\in{\llbracket1,n\rrbracket}\}\bigr\}\cup C_m\cup (K^\ast \times \{0\})$$
which in fact is an equality if $h(t)=1$ and $q\ge 1$. The case $(l,m,q,h)=(1,2,1,1)$ is the ``Multiplicative Trap" of \cite{Bori}, Ex.\ 3.

If $k$ is an algebraic closure of $\mathbb F_p$ (i.e., if $k^{\ast}$ has no element of infinite multiplicative order), the rule 
$$(x,y,z)\mapsto\bigl(x^my^q(x^l-y)h(x),x^{l(m-1)}y^{l+1}(x^l-y)^lh(x)^lz,z\bigr)$$ 
defines an endomorphism 
$$d:=d_{\m,l,m,q,s,h}\in\D_3(k)\setminus\QF_3(k).$$ 
For $n\in\mathbb N^{\ast}$, $\Gamma_{d^n}$ contains $\bigl\{(x,x^lz^q,z)|q\in{\llbracket1,n\rrbracket},(x,z)\in (k^{\ast})^2\bigr\}$ and does not contain the generic point of the zero locus $y-x^lz^{n+1}=0$, so $\iota(d)=\infty$. We have $\deg(d)=l(q+1)+m+s$.\end{example}

\begin{corollary}\label{C4.1} 
If $\chr(K)\neq 2$ (resp.\ $\chr(K)=2$) let $N\ge 2$ (resp.\ $N\ge 3$) be an integer. Let $n\in\{2,3\}$ be such that we have $n=2$ iff $K^{\ast}$ has an element of infinite multiplicative order. Then there exists $e\in\D_n(K)\setminus\QF_n(K)$ such that $\iota(e)=\infty$, $\deg(e)=N$, and $e$ is generically \'etale.
\end{corollary}

\begin{proof}
This follows from Example \ref{EX14}(5) as we can choose $(l,q,m,s,h)$ such that $N=l(q+1)+m+s$. E.g., for $l=m=1$, if $\chr(K)\nmid N$ we can take $(q,s,h)=(N-2,0,1)$ and if $\chr(K)\mid N$ we can take $(q,s,h)=(N-3,1,t+1)$.\end{proof}

\begin{proposition}\label{PR5}
If $\chr(K)=p$, we assume that $K$ is not an algebraic closure of $\mathbb F_p$. Let $n\ge 3$ be an integer. Then there exists an endomorphism $e\in\D_n(K)$ such that $\varphi(e^l)=l$ for each $l\in\mathbb N$ (thus $\iota(e)=\infty$). If $n=3$, then for each integer $s\ge 3$ we can choose $e$ such that moreover we have $\deg(e)=s$.\end{proposition}

\begin{proof}
We can assume that $n=3$ by Lemma \ref{L13}(3) and (4). Let $\alpha\in K$ be in the image of $\mathbb N^{\ast}\setminus\{1\}$ in $K$ if $\chr(K)=0$ and be transcendental over $\mathbb F_p$ if $\chr(K)=p$. 

The rule $(x,y)\mapsto \bigl(xy,\alpha(y+1)\bigr)$ defines an element $d\in\D_2(K)\setminus\QF_2(K)$ which is birational (cf.\ Example \ref{EX14}(2)). As ($\natural_{d,(\beta,\alpha)}$) holds for each $\beta\in K^{\ast}$, we have $\iota(d)=\infty$ by Lemma \ref{L14}. Let $Y$ be the plane in $\mathbb A^3_K$ defined by the equation $x_3=0$ and let $C$ be its curve defined by the equations $x_2-\alpha=x_3=0$. We have $\Gamma_d\times\{0\}=C\setminus\{(0,\alpha,0)\}\cong\mathbb A^1_K\setminus\{0\}$.

For a triple $(m,f,g)\in\mathbb N^{\ast}\times (K[x_1,x_2,x_3]\setminus\{0\})^2$ let 
$$e_{m,f,g}:=e\bigl(x_1x_2+gx_3(g-x_1x_2+1),\alpha(gx_3-1)(fx_3-x_2)+\alpha,x_3^m(gx_3-1)\bigr)\in\End_3(K).$$
It extends $d$ as $e_{m,f,g}(x,y,0)=\bigl(xy,\alpha(y+1),0\bigr)$. We have $e_{m,f,g}^{-1}(Y)=Y\sqcup Z$, where $Z$ is the surface defined by the equation $gx_3=1$. Moreover, $e_{m,f,g}(Z)\subset C$. It follows that, if there exists $Q\in C(K)\setminus [\{(0,\alpha,0)\}\cup e_{m,f,g}(Z)]$, then Condition ($\natural_{e_{m,f,g},Q}$) holds and thus $\iota(e_{m,f,g})=\infty$ by Lemma \ref{L14}. 

Let $P:=(1,\alpha,0)$ and $g:=x_2x_3^q$ with $q\in\mathbb N$. For each $(m,f)\in\mathbb N^{\ast}\times K[x_1,x_2,x_3]$ we have
$$e_{m,f,g}\bigl(Z(K)\bigr)=\{e_{m,f,g}(\gamma,\delta^{-q-1},\delta)|(\gamma,\delta)\in K\times K^{\ast}\}$$
$$=\{(\delta^{-1}+1,\alpha,0)|\delta\in K^{\ast}\}=C(K)\setminus\{P\}.$$

We take either $f=(\alpha+1)x_2-\alpha x_2x_3$ with $m\ge 1$ or $f=1$ with $m\ge 2$; for $e:=e_{m,f,x_2x_3^q}$ we show that we have
\begin{equation}\label{EQ6}
C\setminus [\{(0,\alpha,0)\}\cup e(Z)]=\{P\}=\Gamma_e.
\end{equation}

We already know that the first identity of Equation (\ref{EQ6}) holds. To show that the second identity of Equation (\ref{EQ6}) holds, let $(\beta,\alpha\gamma,\delta)\in K^2\times K^{\ast}$ and we consider the system of equations (without the non-zero factor $\alpha$ of the second equation)
$$xy+yz^{q+1}(yz^q-xy+1)-\beta=(yz^{q+1}-1)\bigl(f(y,z)z-y\bigr)+1-\gamma=z^m(yz^{q+1}-1)-\delta=0.$$
The third equation gives $y=\frac{\delta+z^m}{z^{m+q+1}}$. Based on this, the first equation gives 
$$xy=\frac{\beta-yz^{q+1}(yz^q+1)}{1-yz^{q+1}}.$$ Moreover, $z$ must be a root of the polynomial $h(z)=h_{m,f,\gamma,\delta}(z)\in K[z]$
we get by substituting $y=\frac{\delta+z^m}{z^{m+q+1}}$ in $(yz^{q+1}-1)\bigl(f(y,z)z-y\bigr)+1-\gamma=0$ and clearing the denominators. If $f=(\alpha+1)x_2-\alpha x_2x_3$ with $m\ge 1$, then 
$$f(y,z)z-y=(\alpha+1)yz-\alpha yz^2-y=-y(z-1)(\alpha z-1)$$ 
and 
$$h(z)=(1-\gamma)z^{2m+q+1}-\delta(z^m+\delta)(z-1)(\alpha z-1).$$
If $f=1$ with $m\ge 2$ then 
$$h(z)=\frac{(\gamma-1)}{\delta}z^{2m+q+1}-z^{m+q+2}+z^m+\delta.$$

We show that the assumption that all roots of $h(z)=0$ are $m$-th root of $-\delta$ leads to a contradiction. For this we can assume that there exists $z_0\in K$ such that $h(z_0)=0=z_0^m+\delta$. As $\delta\neq 0$, we have $z_0\neq 0$.

If $f=(\alpha+1)x_2-\alpha x_2x_3$ with $m\ge 1$, then $h(z_0)=(\gamma-1)z_0^{2m+q+1}=0$ implies $1-\gamma=0$ and hence $h(z)=\delta(z^m+\delta)(z-1)(\alpha z-1)$. It follows that both $1$ and $\alpha^{-m}$ are equal to $-\delta$, hence $\alpha^m=1$, a contradiction.

If $f=1$ with $m\ge 2$, then we can assume that $\gamma_1:=\frac{(\gamma-1)}{\delta}\neq 0$. It follows that $\gamma_1\delta^2z_0^{q+1}+\delta z_0^{q+2}=0$ and hence, as $z_0\delta\neq 0$, we have $z_0=-\gamma_1\delta$. As $m\ge 2$, we have $2m+q+1>m+q+2$, hence $h(z)$ has degree $2m+q+1$. Thus $h(z)=-\gamma_1(z+\gamma_1\delta)^{2m+q+1}$ as $z_0$ is the only possible root. By identifying the coefficient of $z$, as $m\ge 2$ we get that $0=(2m+q+1)\gamma_1^{2m+q+1}\delta^{2m+q}$, hence $\chr(K)=p$ with $p$ dividing $2m+q+1$. So all positive exponents of powers of $z$ that show up in $h(z)$ with non-zero coefficients are divisible by $p$, i.e., $p$ divides also $m$ and $m+q+2$. So $p$ divides $1=m+(m+q+2)-(2m+q+1)$, a contradiction.

For a root of $h(z)$ which is not an $m$-th root of $-\delta$, we have $y\neq 0$ and thus we can solve for $x$ and hence the system is consistent. From this and the descriptions of $\Gamma_d$ and $e_{m,f,g}\bigl(Z(K)\bigr)=e\bigl(Z(K)\bigr)$ we get that the second identity of Equation (\ref{EQ6}) also holds. Note that it follows that $Z$ is a trap of $e$ of dimension $2$.

If $f=(\alpha+1)x_2-\alpha x_2x_3$ with $m\ge 1$, then for $\gamma=1-\delta$, the resulting polynomial $h(z)=\delta[z^{2m+q+1}-(z^m+\delta)(z-1)(\alpha z-1)]$ is separable of degree $2m+q+1$ for all $\delta$s outside a finite set, hence $\deg(e)=2m+q+1$; for each integer $s\ge 3$, there exists a pair $(m,q)\in\mathbb N^{\ast}\times\mathbb N$ such that $s=2m+q+1$.

Similarly, if $f=1$ with $m\ge 2$, then for $\gamma=\delta+1$ the resulting polynomial $h(z)=z^{2m+q+1}-z^{m+q+2}+z^m+\delta$ is separable of degree $2m+q+1$ for all $\delta$s outside a finite set, hence $\deg(e)=2m+q+1$; for each integer $s\ge 5$, there exists a pair $(m,q)\in (\mathbb N^{\ast}\setminus\{1\})\times\mathbb N$ such that $s=2m+q+1$.

As $\iota(e)=\infty$, the sequence $\bigl(\varphi(e^l)\bigr)_{l\in\mathbb N}$ is strictly increasing. As $\varphi(e)=1$ by Equation (\ref{EQ6}), from Lemma \ref{L1}(2) we get that $\varphi(e^l)\le l$. From the last two sentences we get that $\varphi(e^l)=l$ for each $l\in\mathbb N$.\end{proof}

Proposition \ref{PR5} shows that Theorem \ref{T4} (resp.\ Corollary \ref{C3}) is optimal, i.e., it does not hold in general for endomorphisms that are not quasi-finite (resp.\ not defined over finite fields) in dimensions at least $3$. In order to supplement Corollary \ref{C3} for dimensions at most $2$, we first prove the following general lemma.

\begin{lemma}\label{L15} Let $e:X\rightarrow X$ be an endomorphism of a variety over $K$ with the property that there exists $P\in X(K)$ such that $(\natural_{e,P}$) holds. For $(q,r)\in\mathbb N\times\mathbb N^{\ast}$ with $q<r$, let $Y_{q,r}$ be the Zariski closure in $X$ of the sequence $\bigl(e^{q+ri}(P)\bigr)_{i\ge 0}$ and let $Y^{-}_{q,r}:=\cup_{\{Y\in\Irr(Y_{q,r})|\dim(Y)>0\}} Y$. Then the following properties hold.

\medskip
{\bf (1)} The generic points of the irreducible components of $Y^{-}_{q,r}$ are permuted transitively by $e^r$. In particular, $Y^{-}_{q,r}$ is equidimensional.

\smallskip
{\bf (2)} There exists no $n\in\mathbb N$ such that $e^n(P)$ belongs to two or more irreducible components of $Y^{-}_{q,r}$.

\smallskip
{\bf (3)} The number $|\Irr(Y^{-}_{q,r})|$ does not depend on $q\in\llbracket0,r-1\rrbracket$.

\smallskip
{\bf (4)} The number $\dim(Y^{-}_{q,r})$ does not depend on $(q,r)\in\mathbb N\times\mathbb N^{\ast}$ with $q<r$.

\smallskip
{\bf (5)} The variety $Y^{-}_{q,r}$ is a union of irreducible components of $Y^{-}_{0,1}$. 

\smallskip
{\bf (6)} If $Z\in\Irr(Y^{-}_{0,1})$ contains the sequence $\bigl(e^{q+ri}(P)\bigr)_{i\ge 0}$, then $Z$ is the Zariski closure of this sequence.

\smallskip
{\bf (7)} Let $Z\in\Irr(Y^{-}_{0,1})$ and $s\in\mathbb N^{\ast}$ be such that $e^s(Z)\subset Z$. Then the endomorphism $d:Z\rightarrow Z$ induced by $e^s$ is dominant and generically radicial. Moreover, there exists $m\in\mathbb N$ such that $e^m(P)\in Z$, Condition ($\natural_{d,e^m(P)}$) holds and the orbit $\{d^i\bigl(e^m(P)\bigr)|i\in\mathbb N\}$ is Zariski dense in $Z$.\end{lemma}

\begin{proof}
As $e^r(Y_{q,r})\cup\{e^q(P)\}$ is Zariski dense in $Y_{q,r}$, $e^r$ induces a dominant endomorphism of $Y^{-}_{q,r}$ and thus (see Remark \ref{R11}) it permutes the generic points of varieties in $\Irr(Y^{-}_{q,r})$. Let $Z_1,Z_2\in\Irr(I^{-}_{q,r})$ with $Z_1\neq Z_2$. For $i\in \{1,2\}$ let $m_i\in q+r\mathbb N$ be such that $e^{m_i}(P)\in Z_i$. We can assume that $m_1<m_2$ and that $e^{m_2}(P)$ does not lie in any other irreducible component of $Y^{-}_{q,r}$ different from $Z_2$. So, as $r$ divides $m_2-m_1$, $e^{m_2-m_1}(Z_1)$ must be $Z_2$ and part (1) holds. 

Due to the transitivity property of part (1), $Z_1$ is uniquely determined by $m_1$ and thus part (2) follows by taking $m_1=n$. 

As $Y^{-}_{0,r}$ is equidimensional of some dimension $l\in\mathbb N^{\ast}$, from part (1) it follows that for each $Z\in\Irr(Y^{-}_{0,r})$, the Zariski closure of every $e^i(Z)$ with $i\in\mathbb N$ is of dimension $l$. Thus, as $Y_{q,r}$ is the Zariski closure of $e^q(Y_{0,r})$ and $Y_{0,1}=\cup_{j=0}^{r-1} Y_{j,r}$, parts (3) to (5) hold. 

Part (6) follows from part (4).

For part (7), let $m\in\mathbb N$ be the smallest such that $e^m(P)\in Z$. Then ($\natural_{d,e^m(P)}$) holds as ($\natural_{e,P}$) holds. For $j\in\mathbb N^{\ast}$ such that $js>m$, the sequence $\bigl(e^{m+jsi}(P)\bigr)_{i\ge 0}$ is a subsequence of the orbit $\{d^i\bigl(e^m(P)\bigr)|i\in\mathbb N\}$; thus from part (6) it follows that this orbit is Zariski dense in $Z$. Let $V$ be an open dense subvariety of $Z$ such that $d_V:=d\times_Z 1_V:e^{-1}(V)\rightarrow V$ is finite. By shrinking $V$ we can assume that $d_V$ is finite flat of degree $o\in\mathbb N$ with (see \cite{Gro4}, Cor.\ (9.7.9)) all geometric fibers of the same cardinality $o_{\sep}\in\llbracket1,o\rrbracket$. As $V\cap\{d^i\bigl(e^m(P)\bigr)|i\in\mathbb N^{\ast}\}\neq\emptyset$ by the Zariski density of the mentioned orbit, we have $o_{\sep}=1$ by ($\natural_{d,e^m(P)}$). Thus $d_V$ is finite radicial; so $d$ is dominant and generically radicial. Hence part (7) holds.\end{proof}%mention density in the proof

\begin{proposition}\label{PR6}
Let $e:X\rightarrow X$ be an endomorphism of a variety over $K$ of dimension at most $2$. If $\Gamma_e$ is finite, then $\iota(e)\in\mathbb N$.
\end{proposition}

\begin{proof}
We show that the assumption that $\iota(e)=\infty$ leads to a contradiction. We can assume that there exists $P\in X(K)$ such that ($\natural_{e,P}$) holds by Lemma \ref{L14}. Let $Y^{-}_{0,1}$ be the closed subvariety of $X$ obtained from $P$ as in Lemma \ref{L15}; it is equidimensional of dimension $1$ or $2$ and $e$ permutes transitively its irreducible components by Lemma \ref{L15}(1). Thus, if $\dim(Y^{-}_{0,1})=1$ then the endomorphism of $Y^{-}_{0,1}$ induced by $e$ is dominant hence quasi-finite, so ($\natural_{e,P}$) does not hold by Theorem \ref{T4} and Lemma \ref{L14}, a contradiction.\footnote{For an alternative proof see Example \ref{EX51} of Section \ref{S32}.} Thus $\dim(Y^{-}_{0,1})=2$. 

By replacing $e$ with $e^q$ for some $q\in\mathbb N^{\ast}$, we can assume that there exists an irreducible component $Z\in\Irr(Y^{-}_{0,1})\subset\Irr(X)$ (so $\dim(Z)=2$) such that the endomorphism $e_Z:Z\rightarrow Z$ induced by $e$ is dominant and generically radicial, $(\natural_{e_Z,e^m(P)}$) holds for some $m\in\mathbb N$, and the orbit $\{e_Z^i\bigl(e^m(P)\bigr)|i\in\mathbb N\}$ is Zariski dense in $Z$ (see Lemma \ref{L15}(7)).

Let $W:=\cup_{Z'\in\Irr(X)\setminus\{Z\}} Z\cap Z'$; we have $e(W)\subset W$ by Remark \ref{R11} and hence the mentioned Zariski density gives that $e_Z^i\bigl(e^m(P)\bigr)\notin W$ for all $i\in\mathbb N$. If a point $Q_1\in (Z\setminus W)(K)$ is $e(Q_0)$ for some $Q_0\in X(K)$, then either $Q_0\in (Z\setminus W)(K)$ or $Q_0$ is one of the finitely many isolated points of $X$. Thus there exists a finite set $W_0$ of $K$-valued points of $Z\setminus W$ such that $[Z\setminus (W\cup W_0)]\subset \Imm(e_Z)$. Let $Y\subset [Z\setminus (W\cup W_0)]$ be a dense open subset such that the morphism $e_Z^{-1}(Y)\rightarrow Y$ induced by $e_Z$ (or $e$) is finite radicial. Clearly, $e_Z^{-1}(Y)\subset Z\setminus W$.

By counting points over finite fields (i.e., by considering a model $e_R$ of $e$ over a finitely generated $\mathbb Z$-subalgebra of $K$ and by reducing it modulo a maximal ideal that has finite residue field) we see that the number of $1$-dimensional irreducible components of $(Z\setminus W)\setminus e_Z^{-1}(Y)$ is equal to the number of $1$-dimensional irreducible components of $(Z\setminus W)\setminus Y$. But the generic points of irreducible components of $(Z\setminus W)\setminus Y$ lie in the image of $e_Z$. As $e_Z^{-1}(Y)\subset e_Z^{-1}(Z\setminus W)\subset Z\setminus W$, by counting curves we see that $(Z\setminus W)\setminus e_Z^{-1}(Z\setminus W)$ has no $1$-dimensional irreducible component and no $1$-dimensional irreducible component of $e_Z^{-1}(Z\setminus W)\setminus e_Z^{-1}(Y)$ is contracted. In particular, the morphism $d:e_Z^{-1}(Z\setminus W)\rightarrow Z\setminus W$ induced by $e$ is quasi-finite. 

Let $\mathcal Y:Z^{\n}\rightarrow Z$ be the normalization of $Z$ and let $e_{Z^{\n}}:Z^{\n}\rightarrow Z^{\n}$ be the endomorphism induced by $e_Z$. Let $W_1:=\mathcal Y^{-1}(W)$. Let $W_2$ be the union of $1$-dimensional irreducible components of $W_1$. Let $\Omega:=Z^{\n}\setminus W_2$; the open embedding $\Omega\rightarrow Z^{\n}$ is affine by Nagata's Theorem (see \cite{H1}, Cor.\ 3.3). Thus $e^{-1}_{Z^{\n}}(\Omega)\rightarrow Z^{\n}$ is affine, so the complement $Z^{\n}\setminus e^{-1}_{Z^{\n}}(\Omega)$ is purely of codimension $1$ (see \cite{Gro6}, Exp.\ V, Ex.\ 3.4). But by the above part we get that $e_{Z^{\n}}$ sends $\mathcal Y^{-1}(Z\setminus W)$ to itself except for finitely many $K$-valued points. So $Z^{\n}\setminus e^{-1}_{Z^{\n}}(\Omega)$ does not meet $\mathcal Y^{-1}(Z\setminus W)$, hence does not meet $\Omega$ which is obtained by adding the isolated points of $W_1$. Thus $e_{Z^{\n}}$ restricts to an endomorphism $e_{\Omega}:\Omega\rightarrow\Omega$ which is quasi-finite as, up to finitely many points, it is the normalization of $d$. 

Let $Q\in\Omega(K)$ be such that it lifts $e^m(P)$. As ($\natural_{e_Z,e^m(P)}$) holds, ($\flat_{e_{\Omega},Q}$) holds and hence $\iota(e_{\Omega})=\infty$ by Lemma \ref{L14}, which contradicts Theorem \ref{T4}.\end{proof}

For a normal (possibly not connected) variety $X$ over $K$, let $\rho_X\in\mathbb N$ be the rank of the N\'eron--Severi group $NS(X)$ of $X$ (e.g., see \cite{Kah}, Thm.\ 3 for the classical result that $NS(X)$ is finitely generated). If $l$ is a prime different from $\chr(K)$, then the Kummer sequence in the \'etale topos of $X$ gives an injective homomorphism $NS(X)/l^mNS(X)\rightarrow H^2\bigl(\Reg(X),\mu_{l^m}\bigr)$ for all $m\in\mathbb N$, hence the homomorphism $\imath_Z:NS(X)\otimes_{\mathbb Z} \mathbb Q_l\rightarrow H^2\bigl(\Reg(X),\mathbb Q_l(1)\bigr)=H^2\bigl(\Reg(X),\mathbb Q_l\bigr)$ is injective.\footnote{The Kummer sequence applies as the group $\Pic\bigl(\Reg(X)\bigr)_{\a}$ of divisor classes on $X$ algebraically equivalent to zero is divisible. If $X$ is projective, this group is in bijection to the group of $K$-valued points of the dual of the Albanese variety $\Alb(X)$ of $X$ (see \cite{Lang-S1}, Ch.\ VI, Sect.\ 1, Thm.\ 1; one verifies that the proofs in the projective case give in general a surjection $\Pic\bigl(\Alb(X)\bigr)_{\a}\rightarrow\Pic\bigl(\Reg(X)\bigr)_{\a}$).}

The following supplement to Theorem \ref{T4} also generalizes Lemma \ref{L6}(2).

\begin{theorem}\label{T5} Let $X$ be a non-empty variety over $K$. For $i\in \{0,1\}$ let $r_i(X)$ be the number of irreducible components of $X$ of dimension $\dim(X)-i$. Let $Y$ be the normalization of the union of irreducible components of $X$ of dimension $\dim(X)$. Let $$N(X):=r_0(X)+r_1(X)+\rho(Y)-1.$$ 
Then for each quasi-finite endomorphism $e:X\rightarrow X$ and integer $i\ge N(X)$ we have $\dim\bigl(\Imm(e^i)\setminus\Imm(e^{i+1})\bigr)\le \dim(X)-2$.\end{theorem}

\begin{proof}
There exists $i_0\in\llbracket0,r_0(X)-1\rrbracket$ such that the set of irreducible components of dimension $\dim(X)$ of the Zariski closure of $e^{i_0+j}(X)$ in $X$ is independent on $j\in\mathbb N$ by Remark \ref{R11}. Similarly, there exists $i_1\in\llbracket0,r_1(X)\rrbracket$ such that the set of irreducible components of dimension $\dim(X)-1$ of the Zariski closure of $e^{i_0+i_1+j}(X)$ in $X$ is independent of $j\in\mathbb N$. For $i\in \{0,1\}$ let $Y_i$ be the normalization of the union of irreducible components of each $e^{i_0+i_1+j}(X)$ with $j\in\mathbb N$ of dimension $\dim(X)-i$. It follows that $e$ induces endomorphisms $e_0:Y_0\rightarrow Y_0$ and $e_1:Y_1\rightarrow Y_1$. As $Y_0$ is a clopen of $Y$, we have $\rho_{Y_0}\le\rho_Y$. To end the proof, by replacing $e$ with $e_0$, it suffices to show that if $X$ is normal and equidimensional, then for each quasi-finite endomorphism $e:X\rightarrow X$, the assumption that there exists $r\in\mathbb N$ such that $r>\rho_X$ and $\dim\bigl(\Imm(e^r)\setminus\Imm(e^{r+1})\bigr)\ge\dim(X)-1$ leads to a contradiction. This assumption implies that there exists a quasi-finite morphism 
$$\phi:X\rightarrow X\setminus (\cup_{i=1}^r Y_i),$$ 
where $Y_1,\ldots,Y_r$ are distinct irreducible closed subvarieties of $X$ of dimension $\dim(X)-1$ (say the one induced by $e^r$, with $Y_i\subset\Imm(e^{i-1})\setminus\Imm(e^i)$). Let 
$$\Omega:=X\setminus \bigl[\bigl(\Sing(X)\cup_{i=1}^r \Sing(Y_i)\bigr)\bigcup (\cup_{1\le i<j\le r} Y_i\cap Y_j)\bigr];$$
it is an open subvariety of $\Reg(X)$ whose codimension in $X$ is at least $2$.

Let $V:=\Omega\setminus (\cup_{i=1}^r Y_i)$. We have an exact complex 
$$0=\oplus_{i=1}^r H^1_{Y_i\cap\Omega}(\Omega,\mathbb Q_l)\rightarrow H^1(\Omega,\mathbb Q_l)\rightarrow H^1(V,\mathbb Q_l)$$
$$\rightarrow \oplus_{i=1}^r H^2_{Y_i\cap\Omega}(\Omega,\mathbb Q_l)\xrightarrow{\phi^*} H^2\bigl(\Reg(X),\mathbb Q_l\bigr)$$
in the \'etale cohomology with $\mathbb Q_l$-coefficients; the last exactness follows from the injectivity of $\imath_{\Reg(X)}$, $\imath_V$, and the pullback homomorphism 
$$\phi^*:NS(V)\rightarrow NS\bigl(\phi^{-1}(V)\cap\Reg(X)\bigr)\cong NS\bigl(\Reg(X)\bigr).$$

As $r>\rho_X$, the divisor classes of the $Y_i$s with $i\in\llbracket1,r\rrbracket$ are dependent over $\mathbb Q_l$, hence from the compatibility of cycle classes with algebraic equivalence (see \cite{Del1}, Cycle, Rmk.\ 2.3.10) it follows that $\Ker(\phi^*)\neq 0$. Thus
\begin{equation}\label{EQ5.8}
\dim_{\mathbb Q_l}\bigl(H^1(V,\mathbb Q_l)\bigr)>\dim_{\mathbb Q_l}\bigl(H^1(\Omega,\mathbb Q_l)\bigr)=\dim_{\mathbb Q_l}\bigl(H^1(\Reg(X),\mathbb Q_l)\bigr)
\end{equation} 
with the equality by the purity of branch loci (see \cite{Gro6}, Exp.\ X, Thm.\ 3.4 (i)).

Let $j:W\rightarrow V$ be an open embedding such that $\phi$ is finite over $W$. The functorial homomorphism $\phi^*:H^1(W,\mathbb Q_l)\rightarrow H^1\bigl(\phi^{-1}(W),\mathbb Q_l\bigr)$ is injective as it has a left inverse which is $\deg(\phi)^{-1}$ times the trace map $H^1\bigl(\phi^{-1}(W),\mathbb Q_l\bigr)\rightarrow H^1(W,\mathbb Q_l)$ (see \cite{SGA4-3}, Exp.\ XVII, Thm.\ 6.2.3 and Ex.\ 6.2.6). The restriction homomorphism $H^1(V,\mathbb Q_l)\rightarrow H^1(W,\mathbb Q_l)$ is injective as for finite coefficients $\Lambda$ we have $\Lambda=j_{*}(\Lambda)$. 

From the last two sentences it follows that the composite of the functorial homomorphism $\phi^*:H^1(V,\mathbb Q_l)\rightarrow H^1\bigl(\phi^{-1}(V),\mathbb Q_l\bigr)$ with the restriction homomorphism $H^1\bigl(\phi^{-1}(V),\mathbb Q_l\bigr)\rightarrow H^1\bigl(\phi^{-1}(W),\mathbb Q_l\bigr)$ is injective. Therefore we get that $\phi^*:H^1(V,\mathbb Q_l)\rightarrow H^1\bigl(\phi^{-1}(V),\mathbb Q_l\bigr)$ itself is injective. Also, the restriction 
$$H^1\bigl(\phi^{-1}(V),\mathbb Q_l\bigr)\rightarrow H^1\bigl(\phi^{-1}(V)\cap \Reg(X),\mathbb Q_l\bigr)\cong H^1\bigl(\Reg(X),\mathbb Q_l\bigr)$$ 
is injective. From the last two sentences we get that Inequality (\ref{EQ5.8}) is false, a contradiction.\end{proof}

\begin{remark}\normalfont\label{R12}
Let $X_{\perf}$ be the perfection of a variety $X$ over $k$ and we consider an endomorphism $e:X_{\perf}\rightarrow X_{\perf}$ with finite fibers. Then the proof of Theorem \ref{T4} applies to give that $\iota(e)\in\mathbb N$ and the proof of Theorem \ref{T5} applies to give that for every integer $i\ge N(X)$ we have $\dim\bigl(\Imm(e^i)\setminus\Imm(e^{i+1})\bigr)\le \dim(X)-2$. If $\Frob_k:\Spec k\rightarrow\Spec k$ is the Frobenius automorphism, then there exists $m\in\mathbb Z$ such that $e$ is defined by a quasi-finite morphism $X\rightarrow (\Frob_k^m)^{\ast}(X)$. If $\dim(X)\le 2$ and $d:X_{\perf}\rightarrow X_{\perf}$ is an endomorphism such that $\Gamma_d$ is finite, then the proof of Proposition \ref{PR6} applies to give that $\iota(d)\in\mathbb N$; there exists $l\in\mathbb Z$ such that $d$ is defined by a morphism $d_l:X\rightarrow (\Frob_k^l)^{\ast}(X)$ with $(\Frob_k^l)^{\ast}(X)\setminus\Imm(d_l)$ finite.\end{remark}

\section{On gap and iteration invariants for $n=2$}\label{S17}

In this section we apply the previous one to obtain several properties of gap and iteration invariants for $n=2$.

\begin{example}\normalfont\label{EX15}
Let $e\in\QF_2(K)$ be such that $\varphi(e)>0$ and the following subset $\mathcal E:=\{P\in k^2||e^{-1}(P)|=1\}$ 
of $k^2$ is infinite; so $\varsigma_e=1$. E.g., if $e$ is as in Example \ref{EX2} we have $\mathcal E\supset K^{\ast}\times\{\alpha\in K|f(\alpha)=0\}$, if $K=k$ and $e$ is as in Example \ref{EX8} with $s\ge 1$ we have $\mathcal E=\{\alpha_1,\ldots,\alpha_s\}\times k^{\ast}$, and if $K=k$ and $e$ is as in Example \ref{EX11} we have $\mathcal E\supset\{0\}\times (k\setminus\mathbb I_{p,m,q}^{\ast})$. 

Let $\Gamma:=\Gamma_e$. Let $(l,s)\in\llbracket0,|\Gamma|\rrbracket\times \mathbb N^{\ast}$. Let $\Gamma_0$ be a subset of $\Gamma$ with $|\Gamma_0|=l$. Let $\Gamma_0^{\perp}:=\Gamma\setminus\Gamma_0$; we have $|\Gamma_0^{\perp}|=\varphi(e)-l$. Let $a_{l,s}\in\GA_2(K)$ be such that by defining $d_{l,s}:=ea_{l,s}\in\QF_2(K)$ and $\Gamma_i:=d_{l,s}(\Gamma_{i-1})$ and $\Gamma_i^{\perp}:=d_{l,s}(\Gamma_{i-1}^{\perp})$ for $i\in \llbracket1,s\rrbracket$, for each $i\in\llbracket0,s-1\rrbracket$ we have $a_{l,s}(\Gamma_i)\subset [e^{-1}(\mathcal E)]\setminus\Gamma$ and $a_{l,s}(\Gamma_i^{\perp})\cap e^{-1}(\mathcal E)=\emptyset$; the existence of $a_{l,s}$ follows from \cite{Sr}, Thm.\ 2 applied to $s\varphi(e)$ points. We have disjoint unions $\sqcup_{i=0}^s \Gamma_i$ and $\sqcup_{i=0}^s \Gamma_i^{\perp}$ with $(s+1)l$ and $(s+1)[\varphi(e)-l]$ (respectively) elements. We have $\Gamma_i\subset\Imm(d_{l,s}^i)\setminus\Imm(d_{l,s}^{i+1})$ for each $i\in \llbracket1,s\rrbracket$ as $d_{l,m}^{-1}(\Gamma_j)=\Gamma_{j-1}$ for $j\in\llbracket1,i\rrbracket$ implies that $d_{l,s}^{-i-1}(\Gamma_i)=d_{l,s}^{-1}(\Gamma_0)=\emptyset$. Thus, with Notation \ref{NOT1} applied to $d_{l,s}(K):K^2\rightarrow K^2$ but using only $d_{l,s}$ instead of $d_{l,s}(K)$, we have $\mathcal E_i(d_{l,s})\subset\cup_{j=1}^i d_{l,s}^j(\Gamma)$ by Lemma \ref{L1}(4), so $\mathcal E_i(d_{l,s})\subset\cup_{j=1}^i (\Gamma_j\cup\Gamma_j^{\perp})$. It follows that $\mathcal E^1_i(d_{l,s})=\sqcup_{j=1}^i \Gamma_j$ for each $i\in \llbracket1,s\rrbracket$. Thus $\varrho_{d_{l,s}}\ge ls$, $\xi_{d_{l,s}}\ge s$, and $\iota(d_{l,s})\ge s+1$ but $\varphi(d_{l,s})=\varphi(e)$ and $\deg(d_{l,s})=\deg(e)$ do not depend on either $l$ or $s$.\end{example}

Directly from Example \ref{EX15} we get the following result.

\begin{corollary}\label{C5}
For each $r\in\mathbb N^{\ast}$, the following sets 
$$\{\iota(e)|e\in\EE_2(k),\varphi(e)=r,\psi_e\;\textup{is \'etale}\},$$ 
$$\{\iota(e)|e\in\EE_2(k),\varphi(e)=r,\psi_e\;\textup{is non-\'etale}\},$$ and 
$$\{\iota(e)|e\in\QF_2(K),\varphi(e)=r\}$$ 
are infinite.\end{corollary}

\begin{example}\normalfont\label{EX16}
We refer to $e\in\EE_2(k)$ of Example \ref{EX11}; hence $\varphi(e)=mq-1$, $\Gamma_e=\{(0,\beta)|\beta\in\mathbb I_{p,m,q}^{\ast}\}$, and for $\mathcal B:=\{P\in k^2|e^{-1}\bigl((e(P)\bigr)=\{P\}\}$ we have the identity
$$\mathcal B=\{\bigl(-(1+(-1)^m\beta^{mq-1})^{-q},\beta+(-1)^m\beta^{mq}\bigr)|\beta\in k\setminus \mathbb I_{p,m,q}^{\ast}\}.$$ We apply Example \ref{EX15} with $l\in\llbracket0,mq-1\rrbracket$ and $s\in\mathbb N^{\ast}$; for instance, as we have $\mathcal F(e)=\{0,1,mq-1\}$, if $mq>1$ and $s=1$, then $a_{l,1}\in\GA_2(k)$ is subject to the only constraint that for $\Gamma_0:=\{P\in\Gamma_e|a_{l,1}(P)\in\mathcal B\}$ we have $|\Gamma_0|=l$. If $P\in\mathcal E_1(d_{l,s})\setminus\mathcal E_1^1(d_{l,s})\subset d_{l,s}(\Gamma_0)\setminus\Gamma_1$, then $P\in\Gamma_1^{\perp}$ is such that $e^{-1}(P)$ is contained in $\Gamma\setminus\Gamma_0=\Gamma_0^{\perp}$; as $|\Gamma_0^{\perp}|=mq-1-l$ and $\mathcal F(d_{l,s})=\mathcal F(e)=\{0,1,mq-1\}$, such a point $P$ can exist only when $l=0$ and if it exists it is unique. Thus for $l\in\llbracket1,mq-1\rrbracket$, the set $\mathcal E_1(d_{l,s})=\mathcal E_1^1(d_{l,s})=\Gamma_1$ has $l$ elements. Moreover, 
\begin{equation*}\label{EQ7}
\varphi(d_{l,s}^2)=\varphi(d_{l,s})+l\;\;\;\textup{and}\;\;\; \iota(d_l)\ge 2
\end{equation*}
by Lemma \ref{L1}(3). Concretely, if the elements of the set $\mathbb I_{p,m,q}^{\ast}=\{\delta_1,\ldots,\delta_{mq-1}\}$ are listed, $P_i=(0,\delta_i)$, and $Q_i=(\gamma_i,\delta_i)$ for $i\in\llbracket1,mq-1\rrbracket$ with each $\gamma_i\in k$ such that we have $Q_i\in\mathcal B$ iff $i\in\llbracket1,l\rrbracket$, then we can take $a_{l,1}:=\e\bigl(x_1+f(x_2),x_2\bigr)\in\TGA_2(k)$ where $f\in k[t]$ is the only polynomial of degree at most $mq-2$ such that $f(\delta_i)=\gamma_i$ for all $i\in\llbracket1,mq-1\rrbracket$. Based on Example \ref{EX11}, with this choice of $a_{l,1}$ we have $d_{l,1}=\e(g_1,g_2)$ with $\pi(g_1,g_2)\le 4$ if $mq=2$ and with 
$$\pi(g_1,g_2)\le (mq-2)[mq+\max\bigl((mq-1)q, m\bigl)$$ 
if $mq>2$. We have $\varphi(d_{mq-1,s}^2)=2\varphi(d_{mq-1,s})$. If $mq=p^r$ with $r\in\mathbb N^\ast$, then $\mathbb I_{p,m,q}=\mathbb F_{p^r}$; thus for $l\in\llbracket1,mq-1\rrbracket$, $a_{l,s}$ is not defined over $\mathbb F_{p^r}$.
\end{example}

\begin{theorem}\label{T6}
We have identities $\{\varphi(e)|e\in\EE_2(k),\psi_e\;\textup{is non-\'etale}\}=\mathbb N$ and
$$\{\varphi(e)|e\in\D_2(K)\setminus\QF_2(K)\;\textup{defined over the prime field of}\; K, \iota(e)\in\mathbb N\}=\mathbb N\cup\{\infty\}.$$
\end{theorem}

\begin{proof} Let $r\in\mathbb N$. To show that there exists $e_r\in\EE_2(k)$ such that $\varphi(e_r)=r$, in order to compute the geometric degrees and estimate the algebraic degrees, we consider four disjoint cases as follows. In all cases, the $e_r$s are composites of \'etale endomorphisms among which at least one has a Jacobian surface which is non-\'etale based on Theorem \ref{T1}(3) and hence the $\psi_{e_r}$s are automatically non-\'etale.

{\bf Case 1: $r=0$.} We take $e_0$ to be a composite $cd$ with $d\in\EE_2(k)$ non-\'etale with $\deg(d)=p$, $\varphi(d)=p-1$, and $\pi(d)\le 2p$ if $p>2$ and with $\deg(d)=4$, $\varphi(d)=3$, and $\pi(d)\le 8$ if $p=2$ (see Theorem \ref{T1}(3)) and with $c\in\EE_2(k)$ finite with $\deg(c)=p$ and $\pi(c)\le p$ if $p>2$ and with $\deg(c)=4$ and $\pi(c)\le 4$ if $p=2$. As $\varphi(d)<\deg(c)$, $e_0=cd$ is surjective, i.e., $\varphi(e_0)=0$.

{\bf Case 2: $r=po-1$ with $o\in\mathbb N^{\ast}$ and $(p,o)\neq (2,1)$.} We take $e_r$ to be the $e$ of Example \ref{EX11} with $m=po$ and $q=1$, thus $\deg(e_r)=mq=po$, $\varphi(e_r)=r$, and $\pi(e_r)\le 2mq=2po$.

{\bf Case 3: $po-1<r<p(o+1)-1$ with $o\in\mathbb N^{\ast}$.} We take $e_r$ to be $d^2_{r-po+1,s}$ of Example \ref{EX16} applied to $m=po$ and $q=1$; as $r-po+1\in\llbracket1,mq-1\rrbracket$ we have $\varphi(e_r)=\varphi(d_{r-po+1,s})+r-po+1=r$, $\deg(e_r)=\deg(d_{r-po+1,s})^2=(mq)^2=p^2o^2$, and when $s=1$ we can choose $d_{r-po+1,s}$ such that $\pi(e_r)\le \max\bigr(2po,2po(po-2)\bigr)^2$.

{\bf Case 4: $r\in \llbracket1,p-2\rrbracket$ or $(p,r)=(2,1)$.} Let $e_{pr}\in\EE_2(k)$ be such that we have $\varphi(e_{pr})=pr$, $\deg(e_{pr})=(pr)^2$, and $\pi(e_{pr})\le \max\bigl(2pr,2pr(pr-2)\bigr)^2$ (see Case 3 with $o=r$). Based on \cite{Sr}, Thm.\ 2, let $a=\e(h_1,h_2)\in\GA_2(k)$ be such that the set $\Gamma_{ae_{pr}}$ is $\mathbb F_p\times\mathbb I_r\subset k^2$ with $\mathbb I_r$ a subset of $k$ that has $r$ elements. We check that we can assume that $\pi(h_1,h_2)\le (pr-1)^2$. To check this, using affine automorphisms it suffices to show that if $(\alpha_1,\ldots,\alpha_{pr})$ and $(\delta_1,\ldots,\delta_{pr})$ are $pr$-tuples of distinct elements of $k$, then for each $(\beta_1,\ldots,\beta_{pr},\gamma_1,\ldots,\gamma_{pr})\in k^{2pr}$, there exists $b=\e(g_1,g_2)\in\GA_2(k)$ such that $b(\alpha_i,\beta_i)=(\gamma_i,\delta_i)$ for each $i\in\llbracket1,pr\rrbracket$ with $\pi(g_1,g_2)\le (pr-1)^2$. We can take $b=b_1b_2$, with $b_1=\e\bigl(x_1,x_2+f_1(x_1)\bigr)$ and $b_2=\bigl(x_1+f_2(x_2),x_2\bigr)$ in $\GA_2(k)$, where the pair $(f_1,f_2)\in k[t]^2$ is formed by Lagrange interpolating polynomials of degree at most $pr-1$ given by the identities $\delta_i=\beta_i-f_1(\alpha_i)$ and $\gamma_i=\alpha_i+f_2(\delta_i)$ for each $i\in\llbracket a,r\rrbracket$; clearly, we have inequalities $\pi(g_1,g_2)\le\deg(f_1)\deg(f_2)\le (pr-1)^2$. 

If $e_r$ is the composite of $ae_{pr}:\mathbb A^2_{k,\s}\rightarrow \mathbb A^2_{k,\t}$ with the finite \'etale endomorphism $\mathbb A^2_{k,\t}\rightarrow\mathbb A^2_{k,\t}$ which defines the quotient epimorphism $\mathbb G_{\a,k}^2\rightarrow \mathbb G_{\a,k}^2/(\mathbb F_p\times 0)\cong\mathbb G_{\a,k}^2$ and thus is defined by the pair $(x_1^p-x_1,x_2)$, then we have $\Gamma_{e_r}=\{0\}\times\mathbb I_r$, and hence $\varphi(e_r)=r$; we have identities $\deg(e_r)=p\deg(e_{pr})=p^3r^2$ and an inequality $\pi(e_r)\le p(pr-1)^2\max\bigr(2pr,2pr(pr-2)\bigr)^2$. As $\varsigma_{e_r}-1\ge p-1\ge r$, we have $\iota(e_r)=1$ by Example \ref{EX12}.

All endomorphisms of this paragraph are defined over the prime field of $K$. The rule $(x,y)\mapsto (xy,y)$ (resp.\ $(x,y)\mapsto (xy,y+1)$) defines $e_K\in\D_2(K)\setminus\QF_2(K)$ with $\pi(e_K)=2$, $\varphi(e_K)=\infty$, and $\iota(e_K)=1$ (resp.\ and $\iota(e_K)$ equal to $\infty$ if $\chr(K)=0$ by Example \ref{EX14}(2) and equal to $p$ if $\chr(K)=p$). The rule $(x,y)\mapsto (xy,y^2+y)$ defines $d_K\in\D_2(K)\setminus\QF_2(K)$ with $\varphi(d_K)=\iota(d_K)=0$ and $\pi(d_K)=2$. If $r\in\mathbb N^{\ast}$, let $e_{r,K}\in\QF_2(K)$ be such that $\pi(e_{r,K})\le 2r+1$ and $\varphi(e_{r,K})=r$ (see Example \ref{EX2}); we have $\varphi(e_{r,K}d_K)=r$, $\pi(e_{r,K}d_K)\le 4q+2$, and Corollary \ref{C3} gives $\iota(e_{r,K}d_K)\in\mathbb N^{\ast}$ provided $\chr(K)=p$. If $\chr(K)=0$, then for a finite surjective endomorphism $c\in\End_2(K)$ with $\varsigma_c\ge 2$, we have $e_{r,K}d_Kc\in\D_2(K)\setminus\QF_2(K)$ by Remark \ref{R1}, $\varphi(e_{r,K}d_Kc)=r$, and (see Example \ref{EX12}) $\iota(e_{r,K}d_Kc)\in\mathbb N$.\end{proof} 

We have the following consequence of either Theorem \ref{T1}(2) or Theorem \ref{T6}.

\begin{corollary}\label{C6}
Let $(l,q)\in\mathbb N^2$. Then there exist an affine connected smooth surface $Y$ over $\Spec k$ that has an open subvariety isomorphic to $\mathbb A^2_k$ and that is equipped with an \'etale endomorphism $d:Y\rightarrow Y$ whose complement $Y\setminus d(Y)$ is a disjoint union of $l$ copies of $\mathbb A^1_k$ and of $q$ points.
\end{corollary}

\begin{proof}
We can assume $l\ge 1$ by either Theorem \ref{T1}(2) or Theorem \ref{T6}. Let $e\in\EE_2(k)$ be such that $\varphi(e)=q$ and $\deg(e)\ge p$. Let $m\in\mathbb N^{\ast}$ be such that $pm-1\ge l$. Let $(f,g):=(t^{pm-1}-1,1)\in\Theta^1_{k,m}$. We take $Y$ to be an open subvariety of $\mathbb S^2_{f,g,k}$ which is the union of $\Imm(\imath_{f,g,k})$ and of $l$ connected components of $\mathbb D^1_{f,g,k}$ (see Display (\ref{EQ0c}); $|\Irr(\mathbb D^1_{f,g,k})|=pm-1$). Let $\imath_Y:Y\rightarrow \mathbb S^2_{f,g,k}$ be the open embedding and let $\imath_{f,g,k,Y}:\mathbb A^2_k\rightarrow Y$ be the factorization of $\imath_{f,g,k}$; as $c_{f,g,k}\circ\imath_{f,g,k}$ is surjective, so is $c_{f,g,k}\circ\imath_Y$. Then the endomorphism
$$d:=\imath_{f,g,k,Y}\circ e\circ c_{f,g,k}\circ\imath_Y:Y\rightarrow Y$$ 
is \'etale and $Y\setminus d(Y)=[Y\setminus\Imm(\imath_{f,g,k})]\sqcup\Gamma_e$ is a disjoint union of $l$ copies of $\mathbb A^1_k$ and of $q$ points; we have $\deg(d)=pm\deg(e)\ge p(l+1)$. For $\deg(e)=p$ (see Theorem \ref{T1}(2) for $q\ge 1$) and $m:=\lceil\frac{l+1}{p}\rceil$ we get the smallest value $\deg(d)=p^2\lceil\frac{l+1}{p}\rceil$.\end{proof}

\phantomsection{As the subsets of $K^2$ of cardinality $m\in\mathbb N$ are permuted transitively by $\GA_2(K)$ (e.g., see \cite{Sr}, Thm.\ 2), the isomorphism class of the quasi-affine scheme $\mathbb B^2_{m,K}$ over $\Spec K$ obtained from $\mathbb A^2_K$ by removing $m$ distinct $K$-valued points is well defined.}\label{EXTRA3}

\begin{corollary}\label{C7}
For each $(m,q)\in\mathbb N^2$, the following properties hold.

\medskip
{\bf (1)} There exists a surjective \'etale morphism $d_{m,q}:\mathbb B^2_{m,k}\rightarrow\mathbb B^2_{q,k}$ which extends (or which does not extend) to a finite \'etale cover of $\mathbb B^2_{q,k}$.

\smallskip
{\bf (2)} There exists a surjective quasi-finite morphism $d_{m,q}:\mathbb B^2_{m,K}\rightarrow\mathbb B^2_{q,K}$.
\end{corollary}

\begin{proof}
The existence of $d_{0,q}$ in part (1) follows from Theorem \ref{T1}(2) in case we want it to extend to a finite \'etale cover of $\mathbb B^2_{q,k}$ and from Theorem \ref{T6} in case we do not want it to extend to a finite \'etale cover of $\mathbb B^2_{q,k}$. The existence of $d_{0,q}$ in part (2) follows from Example \ref{EX2}. If $q>0$, then $d_{0,q}$ is generically \'etale and non-birational. For part (1) (resp.\ (2)), we choose $d_{0,0}$ to be \'etale (resp.\ generically \'etale) and finite non-birational. Hence there exists an infinite number of points $P\in \mathbb B^2_{q,k}(k)$ (resp.\ $P\in \mathbb B^2_{q,k}(K)$) such that $|d_{0,q}^{-1}(P)|\ge 2$; by removing from each one of $m$ such distinct fibers one point, $d_{0,q}$ restricts to a surjective morphism $d_{m,q}:\mathbb B^2_{m,k}\rightarrow\mathbb B^2_{q,k}$ (resp.\ $d_{m,q}:\mathbb B^2_{m,K}\rightarrow\mathbb B^2_{q,K}$) which is \'etale (resp.\ quasi-finite).\end{proof}

\section{Proof of Corollary \ref{C1}}\label{S18}

\noindent
{\bf Proof of Corollary \ref{C1}.} Based on Lemma \ref{L13}(7), to prove Corollary \ref{C1} we can assume that $n=l+2$. Let $d=\e(g_1,g_2)\in\EE_2(k)$ be such that $\varphi(d)=q$ by Theorem \ref{T1}(2) or Theorem \ref{T6}. Let $(s_1,s_2)\in (\mathbb N^{\ast})^2$ be such that $ps_1>\deg(g_1)$ and $ps_2>\deg(g_2)$. We write $X=\Spec\bigl(k[x_3,\ldots,x_n]/(f_1,\ldots,f_m)\bigr)$, with $m\in\mathbb N$ and $(f_1,\ldots,f_m)$ a radical ideal of $k[x_3,\ldots,x_n]$. Then the $n$-tuple
$$\bigl(g_1+\sum_{i=1}^m x_1^{p(s_1+i-1)}f_i,g_2+\sum_{i=1}^m x_2^{p(s_2+i-1)}f_i,x_3,\ldots,x_n\bigr)\in k[x_1,\ldots,x_n]^n$$ 
defines an endomorphism $e\in\EE_n(k)$ that restricts to the endomorphism $d\times 1_X$ on $\mathbb A^2_k\times_{\Spec k} X$ and to an \'etale endomorphism $d_P\times 1_P$ on $\mathbb A^2_k\times_{\Spec k} P$ for every $P\in [\Spec(k[x_3,\ldots,x_n])\setminus X](k)$ which is finite (hence surjective) by Lemma \ref{F3}(2). Thus we have $\Gamma_{e^i}=\Gamma_{d^i}\times X$ for every $i\in\mathbb N^{\ast}$. Hence $\iota(e)=\iota(d)$ and thus $\iota(e)\in\mathbb N$ by Theorem \ref{T4}. Corollary \ref{C1} holds by taking $q_i:=\varphi(d^i)$ for each $i\in\mathbb N^{\ast}$ (so $q_1=q$).%\end{proof}

\medskip
Next we use the proof of Corollary \ref{C1} in a way that relates to the traps of Example \ref{EX14} in order to show that there exists $e\in\QF_3(K)$ such that $\iota(e)\in\mathbb N$ and $\Gamma_e^+$ is infinite (cf.\ Lemma \ref{L1}(6)).

\begin{example}\normalfont\label{EX17}
Suppose that $K$ is not an algebraic closure of a finite field. Let $\alpha\in K^{\ast}$ of infinite multiplicative order. If $\chr(K)=0$ (resp.\ $\chr(K)=p$) let $e\in\QF_2(K)$ (resp.\ $e\in\EE_2(K)$) be such that $\varphi(e)>0$. Let $(g_1,g_2)\in k[x_1,x_2]^2$ define $e$ and let $(s_1,s_2)\in (\mathbb N^{\ast})^2$ (resp.\ $(s_1,s_2)\in p(\mathbb N^{\ast})^2$) be such that $s_1>\deg(g_1)$ and $s_2>\deg(g_2)$. Then the triple $\bigl(g_1+(x_3-1)x_1^{s_1},g_2+(x_3-1)x_2^{s_2},\alpha x_3\bigr)$ defines a $d\in\QF_3(K)$ (resp.\ $d\in\EE_3(k)$) such that $\varphi(d)=\varphi(e)$, $\Gamma_d$ is the finite set $\Gamma_e\times\{1\}$ with $\Gamma^+_d=\cup_{i\in\mathbb N} d^i(\Gamma_d)$ infinite, and $\iota(d)\in\mathbb N$ by either Theorem \ref{T4} or Lemma \ref{L1}(10) as it is easy to see that the set $\{P\in\Gamma_d^+||d^{-1}(P)|=1\}$ is finite. If $\varphi(e)=1$, then $d$ induces a bijection $\Gamma_d^+\rightarrow\Gamma_d^+\setminus\Gamma_d$.\footnote{If $\chr(K)=0$, then we can also take the triple $\bigl(g_1+x_3x_1^{s_1},g_2+x_3x_2^{s_2},x_3+1\bigr)$ to define $d$.}
\end{example}

Lemma \ref{L13}(7) increases geometric and algebraic degrees considerably as $n$ grows: e.g., if $p>3$ and $e\in\EE_2(k)$ is as in Example \ref{EX8} with $\deg(e)=pm$ and $\deg(g)=1$ (hence $s=1$), then for $d\in\EE_3(k)$ as in Lemma \ref{L13}(7) applied to $n=2$ we have $\deg(d)\ge 5p^2m$ and for a $d_2\in\EE_4(k)$ obtained from $d$ as in Lemma \ref{L13}(7) applied to $n=3$ we have $\deg(d_2)\ge 2p^3(5m+1)$. We include two examples that start from two cases of Example \ref{EX11}, that pertain to Corollary \ref{C1} with $l=1$ and with either $q=p^2-p$ and $e\in\BE_2(k)$ or $q=p-1$ and $e\in\EE_2(k)$, and that for $n\ge 3$ involve smaller geometric and algebraic degrees than what Lemma \ref{L13}(7) would give when $n$ increases.

\begin{example}\normalfont\label{EX18}
Suppose that $n\ge 3$. 

\medskip
{\bf (1)} Let $e\in\BE_n(k)$ be defined by the $n$-tuple
$$\bigl(x_1-x_1^{p^2}x_2^{p(p^2-1)},x_2-x_1^{p}x_2^{p^2}-x_3^p,x_3-x_3^p-\sum_{i=4}^n x_i^px_3^{p^2+p(i-4)},x_4,x_5,\ldots,x_n\bigr)$$
in $k_{\s}[x_1,\ldots,x_n]^n$ (cf.\ Example \ref{EX11} applied to $m=q=p$ and $\mathbb I_{p,p,p}=\mathbb F_{p^2}$). For $\underline{\alpha}:=(\alpha_1,\ldots,\alpha_n)\in k^n$, the number $\s(\underline{\alpha})$ of solutions of the $F_3$-system
$$x_1-x_1^{p^2}x_2^{p(p^2-1)}-\alpha_1=x_2-x_1^{p}x_2^{p^2}-x_3^p-\alpha_2=x_3-x_3^p-\sum_{i=4}^n \alpha_i^px_3^{p^2+p(i-4)}-\alpha_3=0$$
equals the number of solutions of the $F_n$-system obtained by adding to the $F_3$-system the $n-3$ equations $x_i-\alpha_i=0$ indexed by $i\in\llbracket4,n\rrbracket$.

{\bf Case 1: there exists $i\in\{4,\ldots,n\}$ such that $\alpha_i\neq 0$.} Then the equation 
$$x_3-x_3^p-\sum_{i=4}^n \alpha_i^px_3^{p^2+p(i-4)}=0$$ 
has at least $p^2$ distinct zeros, hence it has a zero $\beta$ with $\beta^p+\alpha_2\notin\mathbb F_{p^2}^{\ast}$. From this and Example \ref{EX11} applied as mentioned we get that $\s(\underline{\alpha})>0$. 

{\bf Case 2: $\alpha_4=\alpha_5=\cdots=\alpha_n=0$.} Let $\beta\in k$ be such that the solution set of the equation $x_3-x_3^p-\alpha_3=0$ is $\{\beta+i|i\in\mathbb F_p\}$. From this and Example \ref{EX11}, as $\beta^p+\alpha_2\in\mathbb F_{p^2}\setminus\mathbb F_p\Leftrightarrow \{\beta^p+i^p+\alpha_2|i\in\mathbb F_p\}\subset\mathbb I_{p,p,p}^{\ast}=\mathbb F_{p^2}^{\ast}$ we get that
$$\s(\underline{\alpha})=0\Leftrightarrow [\alpha_1=-1\;\;\;\textup{and}\;\;\;\beta^p+\alpha_2\in\mathbb F_{p^2}\setminus\mathbb F_p].$$

Based on the two cases we get a disjoint union decomposition 
$$\Gamma_e=\sqcup_{\alpha\in\mathbb F_{p^2}\setminus\mathbb F_p} \{(-1,-t^p+\alpha,-t^p+t,0,0,\ldots,0)|t\in k\}$$ 
into $p^2-p$ curves isomorphic to $\mathbb A^1_k$.

\smallskip
{\bf (2)} If $d\in\EE_n(k)$ is defined by
$$\bigl(x_1-x_1^px_2^{p-1},x_2-x_1^px_2^p-x_3,x_3-\sum_{i=4}^n x_ix_3^{p(i-3)},x_4,x_5,\ldots,x_n\bigr)\in k[x_1,\ldots,x_n]^n$$
(cf.\ Example \ref{EX11} applied to $m=p$, $q=1$, and $\mathbb I_{p,p,1}=\mathbb F_{p}$), then a similar argument gives a disjoint union decomposition 
$$\Gamma_d=\sqcup_{\alpha\in\mathbb F_p^{\ast}} \{(-1,-t+\alpha,t,0,0,\ldots,0)|t\in k\}$$
into $p-1$ curves isomorphic to $\mathbb A^1_k$.
\end{example}

\section{Proof of Theorem \ref{T2}}\label{S19}

For $q\in\mathbb N^{\ast}$, we denote by $F_{p^q}$ all functorial endomorphisms of linear algebraic group over $\Spec\mathbb F_{p^q}$ which for the $\GL_{n,\mathbb F_{p^q}}$s are defined by the rule that maps each matrix into the matrix whose entries are raised to the $p^q$-th power. For Lang torsors we refer to \cite{Bore2}, Ch.\ V, Thm.\ 16.3 and Cor.\ 16.4.

\phantomsection{The Lang torsor for $\SL_{2,k}$ over $\Spec \mathbb F_{p^q}$ is a morphism}\label{EXT7}
$$\mathcal L:\SL_{2,k,\s}=\SL_{2,k}\rightarrow \SL_{2,k,\t}=\SL_{2,k}$$
defined by the following rule
\begin{equation}\label{EQ7.5}
\begin{bmatrix} 
w & x \\
y & z \\ 
\end{bmatrix}\mapsto \begin{bmatrix} 
w & x \\
y & z \\ 
\end{bmatrix}\cdot F_{p^q}\left(\begin{bmatrix} 
w & x \\
y & z \\ 
\end{bmatrix}\right)^{-1}=\begin{bmatrix} 
wz^{p^q}-xy^{p^q} & -wx^{p^q}+xw^{p^q} \\
yz^{p^q}-zy^{p^q} & -yx^{p^q}+zw^{p^q} \\ 
\end{bmatrix}
\end{equation}
on valued points; so $\mathcal L$ is a right $\SL_2(\mathbb F_{p^q})$-torsor under the right translation action of $\SL_2(\mathbb F_{p^q})$ on $\SL_{2,k,\s}$ (cf.\ \cite{Bore2}, Ch.\ V, Cor.\ 16.4). 

Let $\star\in\{+,-\}$. Let $B^{\star}$ be the Borel subgroup scheme of $\SL_2$ of $2\times 2$ matrices of determinant $1$ that are upper triangular if $\star=+$ and are lower triangular if $\star=-$; we have a short exact sequence $1\rightarrow\mathbb U^{\star}\rightarrow B^{\star}\rightarrow\mathbb G_m\rightarrow 1$ 
of group schemes over $\Spec\mathbb Z$, where $U^{\star}\cong\mathbb G_a$ is the unipotent radical of $B^{\star}$. 

\phantomsection{The product morphism}\label{EXT7}
$$\mathcal T:U^3_k:=U_k^-\times_{\Spec k} U_k^+\times_{\Spec k} U_k^-\cong\mathbb A^3_k\rightarrow\SL_{2,k,\t}$$
is defined by the $k$-algebra homomorphism 
$$\mathcal T^{\#}:k_{\t}[w,x,y,z]/(wz-xy-1)\rightarrow k[x_1,x_2,x_3]$$ 
that maps $w$, $x$, $y$, and $z$ to $1+x_2x_3$, $x_2$, $x_1+x_3+x_1x_2x_3$, and $1+x_1x_2$ (respectively). Thus $\mathcal T$ is birational; more precisely, it induces an isomorphism 
\begin{equation}\label{EQ8}
U^3_k\setminus D_k=\Spec(k[x_1,x_2,x_3,x_2^{-1}])\cong\SL_{2,k,\t}\setminus B^-_k=\Spec(k_{\t}[w,x,z,x^{-1}]),
\end{equation}
where $D:=U^-\times_{\Spec\mathbb Z} 1_{U^+}\times_{\Spec\mathbb Z} U^- $ with $1_{U^+}$ as the identity element subgroup scheme of $U^+$. Hence $D_k=\mathcal T^{-1}(B_k^{-})$ and the pullback $\mathcal T^{\ast}$ on divisor groups defined by $\mathcal T$ is bijective. Note that $\mathcal T(D_k)=U_k^{-}$.

We view $\mathcal T$ and $\mathcal L$ as $\mathbb G_{\a,k}$-invariant morphisms, the actions of $\mathbb G_{\a,k}$ on $U_k^3$, $\SL_{2,k,\t}$, and $\SL_{2,k,\s}$ being defined by the identification $\mathbb G_{\a,k}=U_k^{-}$ that maps every element $z\in \mathbb G_{\a,k}(k)=k$ to $\left[ 
\begin{array}{cc}
1 & 0 \\ 
z & 1\end{array}\right]\in U_k^{-}(k)$ and the following respective rules on valued points
$$\bigl(u,(u_1,u_2,u_3)\bigr)\mapsto\bigl(u\cdot u_1,u_2,u_3\cdot F_{p^q}(u)^{-1}\bigr), (u,x)\mapsto u\cdot x\cdot F_{p^q}(u)^{-1},(u,x)\mapsto u\cdot x.$$

\phantomsection{We consider the cartesian diagram of $\mathbb G_{\a,k}$-varieties}\label{EXT11}
\[\xymatrix{
\mathbb X_k^3 \ar[r]^{\mathcal I} \ar[d]^{\mathcal M_3} & \SL_{2,k,\s} \ar[d]^{\mathcal L} \\
U_k^3 \ar[r]^{\mathcal T} & \SL_{2,k,\t}
}\]
which is defined over $\mathbb F_{p^q}$; so $\mathbb X_k^3$ has a canonical model $\mathbb X^3$ which is a $\mathbb G_{\a,\mathbb F_{p^q}}$-variety over $\Spec\mathbb F_{p^q}$. Moreover, $\mathcal M_3:\mathbb X^3_k\rightarrow U_k^3$ is a right $\SL_2(\mathbb F_{p^q})$-torsor. From the rule defining $\mathcal L$ and the description of $\mathcal T^{\#}$ we get that
$$\mathbb X^3_k=\Spec k[x_1,x_2,x_3,w,x,y,z]/(J),$$
where $J$ is the ideal generated by $x_2x_3-wz^{p^q}+xy^{p^q}-1$, $x_2+wx^{p^q}-xw^{p^q}$, $x_1x_2x_3+x_1+x_3-yz^{p^q}+zy^{p^q}$, $x_1x_2+yx^{p^q}-zw^{p^q}-1$, and $wz-yz-1$.

\phantomsection{We have a disjoint union decomposition}\label{EXTRA2} 
$$\mathbb B^2_k:=\mathcal L^{-1}(B^-_k)=\sqcup_{i=1}^{p^q+1} B^-_k\tau_i$$
where $\tau_1,\ldots,\tau_{p^q+1}\in\SL_2(\mathbb F_{p^q})$ are representatives of the quotient set $$B^-(\mathbb F_{p^q})\backslash \SL_2(\mathbb F_{p^q}).$$ 
For each $i\in\llbracket1,p^q+1\rrbracket$ we have a disjoint union decomposition
$$B^-_k\tau_i\cap\mathcal L^{-1}(U_k^{-})=\sqcup_{j=1}^{p^q-1} U^-_k\varpi_j\tau_i$$
where $\varpi_1,\ldots,\varpi_{p^q-1}\in B^-(\mathbb F_{p^q})$ are representatives of the quotient set $$U^-(\mathbb F_{p^q})\backslash B^-(\mathbb F_{p^q}).$$
E.g., if we choose $\{\tau_1,\ldots,\tau_{p^q+1}\}=\left\{\left[ 
\begin{array}{cc}
\alpha & 1 \\ 
-1 & 0\end{array}\right]|\alpha\in\mathbb F_{p^q}\right\}\cup\left\{\left[ 
\begin{array}{cc}
1 & 0 \\ 
0 & 1\end{array}\right]\right\}$ and $\{\varpi_1,\ldots,\varpi_{p^q-1}\}=\left\{\left[ 
\begin{array}{cc}
\beta & 0 \\ 
0 & \beta^{-1}\end{array}\right]|\beta\in\mathbb F^{\ast}_{p^q}\right\}$, it follows that $\mathbb B^2_k$ is the closed subscheme of $\SL_{2,k}$ defined by the equation $w^{p^q}x=wx^{p^q}$\footnote{This also follows from the rule of Display (\ref{EQ7.5}).} and one computes that 
$$\mathbb G_{\a,k}\backslash [\SL_{2,k,\s}\setminus\mathbb B^2_k]\cong\mathbb A^2_k\setminus [\{(0,0)\}\cup\{(x_1,x_2)|(x_1:x_2)\in\mathbb P^1(\mathbb F_{p^q})\}].$$

From the last two disjoint union decompositions we get that for each index $i\in\llbracket1,p^q+1\rrbracket$ we have a disjoint union decomposition into irreducible divisors
$$\mathcal I^{-1}(B^-_k\tau_i)=\sqcup_{j=1}^{p^q-1} \mathcal I^{-1}(U^-_k\varpi_j\tau_i)$$
that are canonically isomorphic to the fiber product
$$D_k\times_{\mathcal T|D_k,U^-_k,\mathcal L|U^-_k} U^-_k\cong\mathbb A^2_k.$$
Here $\mathcal T|D_k:D_k\rightarrow U^-_k$ and $\mathcal L|U^-_k:U^-_k\rightarrow U^-_k$ are restrictions of $\mathcal T$ and $\mathcal L$ (respectively) with codomains replaced by $U_k^-$ and the product morphism $\mathcal T|D_k$ is a trivial line bundle. From these disjoint union decompositions and canonical isomorphisms and the existence of the Isomorphism (\ref{EQ8}), we get that the pullback $\mathcal I^{\ast}$ on divisor groups defined by $\mathcal I$ is injective. From this and the identity $\mathcal O(\SL_{2,k,\s})^{\ast}=k^{\ast}$ implied by the existence of $\mathcal T$, we get that $\mathcal O(\mathbb X^3_k)^{\ast}=k^{\ast}$.

As $\mathcal M_3:\mathbb X^3_k\rightarrow U_k^3$ is a right $\SL_2(\mathbb F_{p^q})$-torsor, by taking quotients through the $p$-Sylow subgroup $H^\star:=U^\star(\mathbb F_{p^q})$ of $\SL_2(\mathbb F_{p^q})$, we get a finite \'etale morphism
\begin{equation}\label{EQ9}
\mathbb Y^{3,\star}_k:=\mathbb X^3_k/H^{\star}\rightarrow \mathbb A^3_k
\end{equation}
of degree $p^{2q}-1$ with $\mathcal O(\mathbb Y_k^{3,\star})^{\ast}=k^{\ast}$ (recall that $U_k^3\cong \mathbb A^3_k$) and a birational morphism $\mathcal I/H^\star:\mathbb Y_k^{3,\star}\rightarrow\SL_{2,k,\s}/H^{\star}$,
both defined over $\mathbb F_{p^q}$. To show that $\SL_{2,k,\s}/H^\star$ is rational we can work only with $\star=-$ and this case follows from the fact that $\mathcal T$ is $H^-$-invariant inducing a birational morphism
$$U^3_k/H^-=U^-_k\times_{\Spec k} U^+_k\times_{\Spec k} (U^-_k/H^-)\rightarrow \SL_{2,k,\s}/H^-;$$
note that $(U^-_k/H^-)\cong U^-_k\cong\mathbb G_{\a,k}$ (see \cite{SGA3-2}, Exp.\ XVII, Props.\ 2.2 and 4.1.1). 

Taking quotients through $G^\star:=B^{\star}(\mathbb F_{p^q})$, we similarly get a finite
\'etale morphism
\begin{equation}\label{EQ10}
\mathbb W^{3,\star}_k:=\mathbb X^3_k/G^{\star}\rightarrow \mathbb A^3_k
\end{equation}
of degree $p^q+1$ with $\mathcal O(\mathbb W_k^{3,\star})^{\ast}=k^{\ast}$ and a birational morphism 
$$\mathcal I/G^\star:\mathbb W_k^{3,\star}\rightarrow\SL_{2,k,\s}/G^{\star},$$ 
both defined over $\mathbb F_{p^q}$. To show that $\mathbb W_k^{3,\star}$ is rational we can work only with $\SL_{2,k,\s}/G^{\star}$ and assume that $\star=+$. As the product morphism 
$$U^-_k\times_{\Spec k} B^+_k\rightarrow\SL_{2,k}$$ 
is $G^+$-invariant and an open embedding, it induces an open embedding
$$U^-_k\times_{\Spec k} (B^+_k/G^+)\rightarrow \SL_{2,k,\s}/G^+.$$ 
The Lang torsor for $B^+_k$ over $\Spec\mathbb F_{p^q}$ induces an isomorphism $B^+_k/G^+\cong B^+_k$ and hence $B^+_k/G^+\cong\mathbb A^1_k\times_{\Spec k} \mathbb G_{\m,k}$. From the last two sentences we get that $\SL_{2,k,\s}/G^+$ is birational to $\mathbb A^2_k\times_{\Spec k} \mathbb G_{\m,k}$ and thus it is rational.

\phantomsection{The morphism $\mathcal M_3$ and Morphisms (\ref{EQ9}) and (\ref{EQ10}) already give Theorem \ref{T2}(1) for $n=3$. To get it for $n=2$, we recall that our varieties come with $\mathbb G_{\a,k}$-actions by which we can quotient and define}\label{EXTRA1}
$$\mathbb E_k:=\mathcal M_3^{-1}(D_k)=\mathcal M_3^{-1}\bigl(\mathcal T^{-1}(B_k^{-})\bigr)=(\mathcal T\circ\mathcal M_3)^{-1}(B_k^{-})=(\mathcal L\circ\mathcal I)^{-1}(B_k^{-})$$
$$=\mathcal I^{-1}\bigl(\mathcal L^{-1}(B_k^{-})\bigr)=\sqcup_{i\in\llbracket1,p^q+1\rrbracket}\mathcal I^{-1}(B^-_k\tau_i)=\sqcup_{(i,j)\in\llbracket1,p^q+1\rrbracket\times\llbracket1,p^q-1\rrbracket}\mathcal I^{-1}(U_k^-\varpi_j\tau_i).$$

The morphism $\mathcal I$ induces an isomorphism
\begin{equation}\label{EQ11}
\mathbb X^3_k\setminus\mathbb E_k\cong\SL_{2,k,\s}\setminus\mathbb B^2_k.
\end{equation}
We refer to $\mathbb E_k$ and $\mathbb B^2_k$ as the exceptional divisors.

Let $$U^2_k:=\mathbb G_{\a,k}\backslash U^3_k\cong U^3_k/V_k^1\cong\mathbb G_{\a,k}^2,$$ 
where $V^1_k$ is the connected smooth subgroup of $U^3_k$ defined by the two equations $F_{p^q}(u_1)u_3=u_2=\begin{bmatrix}1 & 0 \\
0 & 1 \\ 
\end{bmatrix}$; recall that valued points of $U^3_k$ are triples $(u_1,u_2,u_3)$ and $\mathbb G_{\a,k}$ acts on $U_k^3$ by multiplication by elements of $V^1_k$. The quotient morphisms $U^3_k\rightarrow U^2_k$ and $\SL_{2,k,\s}\rightarrow\mathbb G_{\a,k}\backslash\SL_{2,k,\s}\cong [A^2_k\setminus \{(0,0)\}]$ are geometric quotients and principal fiber bundles in the sense of \cite{MFK}, Ch.\ 0, Defs.\ 0.6 and 0.10.

We have a natural dominant morphism $\mathbb G_{\a,k}\backslash [\SL_{2,k,\s}\setminus\mathbb B^2_k]\rightarrow U^2_k$ defined by $\mathcal M_3$ via the Isomorphism (\ref{EQ11}) and we define $\mathbb X^2_k$ to be the normalization of $U_k^2$ in the field of fractions of $\mathbb G_{\a,k}\backslash [\SL_{2,k,\s}\setminus\mathbb B^2_k]$. We have a natural finite morphism $\mathcal M_2:\mathbb X^2_k\rightarrow U^2_k$ and $\mathcal O(\mathbb X^2_k)$ is the $k$-subalgebra of $\mathcal O(\mathbb X^3_k)$ of $\mathbb G_{\a,k}$-invariant functions. Thus the inclusion 
$$\mathcal O(\mathbb G_{\a,k}\backslash\SL_{2,k,\s})=\mathcal O(\SL_{2,k,\s})^{\mathbb G_{\a,k}}\rightarrow\mathcal O(\mathbb X^2_k)=\mathcal O(\mathbb X^3_k)^{\mathbb G_{\a,k}}$$ 
between $k$-subalgebras of $\mathbb G_{\a,k}$-invariant functions induced by $\mathcal I$ defines a birational morphism
\begin{equation}\label{EQ12}
\mathcal M_2:\mathbb X^2_k\rightarrow\Spec\bigl(\mathcal O(\mathbb G_{\a,k}\backslash\SL_{2,k,\s})\bigr)\cong\Spec \mathcal O\bigl(\mathbb A^2_k\setminus \{(0,0)\}\bigr)\cong\mathbb A^2_k;
\end{equation}
more precisely, from Zariski's Main Theorem we get that $\mathbb G_{\a,k}\backslash [\SL_{2,k,\s}\setminus\mathbb B^2_k]$ is an open subvariety of $\mathbb X^2_k$. Thus the natural morphism 
\begin{equation}\label{EQ13}
\mathbb X^3_k\rightarrow U^3_k\times_{U^2_k} \mathbb X^2_k
\end{equation}
between normal integral finite schemes over $U^3_k$ is birational and therefore it is a $\mathbb G_{\a,k}$-invariant isomorphism. This implies that the morphism
$$\mathbb X^3_k\rightarrow\mathbb X^2_k$$ is a $\mathbb G_{\a,k}$-torsor. As the $\mathbb G_{\a,k}$-action on $\mathbb X^3_k$ commutes with the right action of $\SL_2(\mathbb F_{p^q})$, $\SL_2(\mathbb F_{p^q})$ acts on $\mathbb X^2_k$ and the Isomorphism (\ref{EQ13}) is $\SL_2(\mathbb F_{p^q})$-invariant. 
Thus $\mathbb X^2_k\rightarrow U^2_k$ is a right $\SL_2(\mathbb F_{p^q})$-torsor and we get morphisms
\begin{equation}\label{EQ14}
\mathbb Y^{2,\star}_k:=\mathbb X^2_k/H^{\star}\rightarrow U_k^2\;\;\;\textup{and}\;\;\;\mathbb W^{2,\star}_k:=\mathbb X^2_k/G^{\star}\rightarrow U_k^2.
\end{equation}
The Isomorphism (\ref{EQ13}) induces isomorphisms
$$\mathbb Y^{3,\star}_k\rightarrow U^3_k\times_{U^2_k} \mathbb Y^{2,\star}_k\;\;\;\textup{and}\;\;\;\mathbb W^3_k\rightarrow U^3_k\times_{U^2_k} \mathbb W^{2,\star}_k$$
which imply that $\mathbb Y^{2,\star}_k$ and $\mathbb W^{2,\star}_k$ are rational with group of units equal to $k^{\ast}$ and that the morphisms 
$$\mathbb Y^{3,\star}_k\rightarrow\mathbb Y^{2,\star}_k\;\;\;\textup{and}\;\;\;\mathbb W^{3,\star}_k\rightarrow\mathbb W^{2,\star}_k$$ are $\mathbb G_{\a,k}$-torsors. 

Theorem \ref{T2}(1) for $n=2$ follows from the existence of Morphisms (\ref{EQ12}) and (\ref{EQ14}).

Assume $p^q=2$. Then $\SL_2(\mathbb F_2)$ is the symmetric group $S_3$ and it has a normal cyclic subgroup $C_3$ of order 3 generated by 
$$\tau:=\begin{bmatrix} 
0 & 1 \\
1 & 1 \\ 
\end{bmatrix}.$$ We choose the $\tau_i$s such that we have $\{\tau_i|i\in\{1,2,3\}\}=C_3$. Let 
$$\mathbb V^3_k:=\mathbb X^3_k/C_3\;\;\;\textup{and}\;\;\;\mathbb V^2_k:=\mathbb X^2_k/C_3;$$ they have canonical models $\mathbb V^3$ and $\mathbb V^2$ (respectively) over $\Spec\mathbb F_2$. Moreover, $\mathbb V^3$ is a $\mathbb G_{\a,\mathbb F_2}$-variety. The morphisms $\mathcal M_3$ and $\mathcal M_2$ factor through finite Galois covers
\begin{equation}\label{EQ15}
\mathbb V^3_k\rightarrow\mathbb A^3_k=\Spec(k[x_1,x_2,x_3])\;\;\;\textup{and}\;\;\;\mathbb V^2_k\rightarrow\mathbb A^2_k
\end{equation}
(respectively) of degree $2$. Moreover, as above we argue that we have an isomorphism
$\mathbb V^3_k\rightarrow U^3_k\times_{U^2_k} \mathbb V^2_k$ and that $\mathbb V^3_k\rightarrow\mathbb V^2_k$ is the geometric quotient of the $\mathbb G_{\a,k}$-action on $\mathbb V^3_k$. As $C_3$ permutes the irreducible components of the exceptional divisors, we get that the image of $\mathbb E_k$ in $\mathbb V^3_k$ is an irreducible divisor $\mathbb J_k$ given by the equation $x_2=0$; thus $x_2$ is a prime element of $A:=\mathcal O(\mathbb V^3_k)$. As $\mathcal I^{-1}(U_k)=U_k\times_{\Spec k} 1_{U^+_k}\times_{\Spec k} U_k\cong\mathbb A^2_k$ is an irreducible component of $\mathbb E_k$, $A/Ax_2$ is a polynomial $k$-algebra in $2$ indeterminates. Moreover, Isomorphism (\ref{EQ11}) induces an isomorphism
$$\Spec(A_{x_2})=\mathbb V^3_k\setminus\mathbb J_k\cong (\SL_{2,k}\setminus \mathbb B^2_k)/C_3=\Spec(k_{\s}[w,x,z]_v)/C_3$$
with $v:=wx(w+x)$, where the last identity follows from the facts that $B^-$ is defined by the equation $x=0$, $B^-\tau$ is defined by the equation $w=0$ and $B^-\tau^2$ is defined by the equation $w+x=0$. The $k$-automorphism of $k_{\s}[w,x,z]_v$ defined by the right translation action of $\tau^{-1}=\tau^2=\begin{bmatrix} 
1 & 1 \\
1 & 0\\
\end{bmatrix}$ on $\SL_{2,k}$, i.e., by the matrix rule $\begin{bmatrix} 
w & x \\
y & z\\
\end{bmatrix}\begin{bmatrix} 
1 & 1 \\
1 & 0\\
\end{bmatrix}\mapsto \begin{bmatrix} 
w+x & w \\
y+z & y\\
\end{bmatrix}$, maps $w$ to $x+w$, $x$ to $w$, $w+x$ to $x$, $z$ to $y=(wz-1)x^{-1}$, $y$ to $y+z$, and $y+z$ to $z$. Thus $A_{x_2}=k_{\s}[w,x,z]_v^{C_3}$ is the normalization of its $k$-subalgebra $B:=k_{\s}[\omega_1,\omega_2,t,u,v]_v$,
where $\omega_1:=y^2+yz+z^2$, $\omega_2:=yz(y+z)$, $t:=x^3+x^2w+w^3$, and $u:=w^2+wx+x^2$; so $u$ and $v$ are Dickson invariants, see \cite{Di}, Sect.\ 1, Thm. Let $\omega:=w(x+w)z+x\in A_{x_2}$. The identity
\begin{equation}\label{EQ16}
u^3+t^2+tv+v^2=0
\end{equation}
implies that we have a $k$-algebra homomorphism 
$$\sigma_{\omega}:C:=k[t,u,v,\omega]/(u^3+t^2+tv+v^2)\rightarrow A_{x_2}.$$ 
As the left-hand side of Equation (\ref{EQ16}) is the sum of two relatively prime homogeneous polynomials of degrees $3$ and $2$, $\Spec C$ is an integral rational hypersurface in $\mathbb A^4_k$ (see, \cite{Sh}, Sect.\ 3.4, Exc.\ 5 or \cite{Ram}, Probl.\ 11.115); concretely, the rule 
$$(x_1,x_2,x_3)\mapsto \bigl(x_1(x_1^2+x_1x_2+x_2^2),x_2(x_1^2+x_1x_2+x_2^2),x_1^2+x_1x_2+x_2^2,x_3)\bigr)$$
defines a birational morphism $\mathbb A^3_k\rightarrow\Spec C$.\footnote{Let $(m,n)\in [\mathbb N^{\ast}\setminus\{1\}]^2$. Let $f\in R:=\mathbb F_{p^q}[x_1,...,x_n]$ be an irreducible homogeneous polynomial of degree $m$. The hypersurface $W:=\Spec \bigl(R[x_{n+1}]/(x_{n+1}^{m+1}-f)\bigr)$ in $\mathbb A^{n+1}_{\mathbb F_{p^q}}$ is factorial by \cite{Sa1}, Thm.\ or \cite{Sa2}, Ch.\ 1, Thm.\ 8.1 with $\mathcal O(W)^{\ast}=\mathbb F_{p^q}^{\ast}$ and it is regular iff the zero locus $f=0$ is. As the rule $(x_1,\ldots,x_n)\mapsto (x_1f,\ldots,x_nf,f)$ defines a birational morphism $\mathbb A^n_{\mathbb F_{p^q}}\rightarrow W$ which induces an isomorphism between principal opens defined by inverting $f$ in the source and $x_{n+1}$ in the target, $W$ has the same zeta function as $\mathbb A^n_{\mathbb F_{p^q}}$.}
As $C$ is an integral domain, $\sigma_{\omega}$ is injective by reasons of dimension. As we have $[k(w,x,z):\Frac(C)]=[k(x,w):k(t,u,v)]$, $[k(w,x,z):\Frac(A_{x_2})]=3$, $[k(x,w):k(u,v)]=6$ by \cite{Di}, Sect.\ 1, Thm., and $[k(t,u,v):k(u,v)]=2$, we conclude that $\Frac(\sigma_{\omega})$ is an isomorphism. Thus $\mathbb V^3_k$ is rational. Similar arguments give that the inclusion $k[t,u,v]_v\rightarrow \mathcal O(\mathbb V^2_k)_{x_2}$ is birational and that $\mathbb V^2_k$ is rational.%w^3+uw=v=x_2

Theorem \ref{T2}(2) follows from the existence of Morphisms (\ref{EQ15}). 

Assume $p^q=3$. The group $\SL_2(\mathbb F_3)$ has order $24$ and a unique normal subgroup $Q_8$ of order $8$ (it is the quaternion group). Let $\mathbb V^3_k:=\mathbb X^3_k/Q_8$ and $\mathbb V^2_k:=\mathbb X^2/Q_8$; they have canonical models $\mathbb V^3$ and $\mathbb V^2$ (respectively) which are $\mathbb G_{\a,\mathbb F_3}$-varieties over $\Spec\mathbb F_3$. The morphisms $\mathcal M_3$ and $\mathcal M_2$ factor through finite Galois covers
\begin{equation}\label{EQ17}
\mathbb V_k^3\rightarrow\mathbb A^3_k\;\;\;\textup{and}\;\;\;\mathbb V^2_k\rightarrow\mathbb A^2_k
\end{equation}
of degree $3$ and as above we have an isomorphism
$\mathbb V^3_k\rightarrow U^3_k\times_{U^2_k} \mathbb V^2_k$ and $\mathbb V^3_k\rightarrow\mathbb V^2_k$ is the geometric quotient of the $\mathbb G_{\a,k}$-action on $\mathbb V^3_k$. 

L\"uroth's analog theorem for surfaces (see \cite{Za}, Sect.\ 1) applied to the separable composite morphism $B_k^+\rightarrow U^-_k\backslash\SL_{2,k}\rightarrow U^-_k\backslash\SL_{2,k}/Q_8$ gives that $U^-_k\backslash\SL_{2,k}/Q_8$ is rational, hence $\mathbb V_k^2$ is also rational. As $\mathbb V^3_k\cong U^3_k\times_{U^2_k} \mathbb V^2_k$, $\mathbb V^3_k$ is also rational. 

Theorem \ref{T2}(3) follows from the existence of Morphisms (\ref{EQ17}). 

Theorem \ref{T2}(4) follows from constructions. Thus Theorem \ref{T2} holds.

\begin{remark}\normalfont\label{R13}
We can assume that $\tau_1$ is the identity element of $\SL_2(\mathbb F_{p^q})$. Let $\mathbb A^0:=\Spec\mathbb Z$. To simplify the notation, let $\zeta(\surd)$ be the logarithm of the zeta function of a variety $\surd$ over $\mathbb F_{p^q}$. If a variety $W$ over $k$ has a well defined (natural) model $\surd$ over $\mathbb F_{p^q}$, let $\zeta(W):=\zeta(\surd)$. For $l\in\mathbb N$, let 
$$\zeta_{l,p^q}:=\zeta(\mathbb A_{\mathbb F_{p^q}}^l).$$ 
Based on the mentioned canonical isomorphisms and Isomorphism (\ref{EQ11}) and their well defined models over $\mathbb F_{p^q}$, we compute
$$\zeta(\mathbb X^3)=\zeta(\mathbb X^3_k\setminus \mathbb E_k)+\zeta(\mathbb E_k)=\zeta\bigl(\SL_{2,k,\s}\setminus (B_k^{-}\sqcup_{i=2}^{p^q+1} B_k^-\tau_i)\bigr)+\zeta(\mathbb E_k)$$
$$=\zeta(\SL_{2,k,\s}\setminus B_k^{-})-p^q\zeta(B_k^-)+\zeta(\mathbb E_k)=\zeta(U_k^3\setminus D_k)-p^q\zeta(B_k^-)+\zeta(\mathbb E_k)$$
$$=\zeta_{3,p^q}-\zeta_{2,p^q}-p^q(\zeta_{2,p^q}-\zeta_{1,p^q})+(p^q+1)(p^q-1)\zeta_{2,p^q}$$
$$\;\;\;\;\;\;\;=\zeta_{3,p^q}+(p^{2q}-p^q-2)\zeta_{2,p^q}+p^q\zeta_{1,p^q}.$$
Thus
$$\zeta(\mathbb X^2)=\zeta_{2,p^q}+(p^{2q}-p^q-2)\zeta_{1,p^q}+p^q\zeta_{0,p^q}.$$
We compute the zeta function of the canonical model $\mathbb Y^{3,-}$ of $\mathbb Y^{3,-}_k$. 

As $U^-$ is normal in $B^-$, via right translations, we have $B^-H^-=B^-$ and $U^-\varpi_j\tau_1H^-=U^-\varpi_j\tau_1$ for each $j\in\llbracket1,p^q-1\rrbracket$, with 
$$U^{-}\varpi_j\tau_1/H^-\cong U^{-}/\varpi_jH^{-}\varpi_j^{-1}\cong U^{-}/H^{-}\cong U^-\cong\mathbb A^2_k.$$ 
As the normalizer of $U^-_k$ (equivalently, of any non-identity element of $U^-(k)$) in $\SL_{2,k}$ is $B^-_k$, $H^-$ acts transitively on the set $\{B^-\tau_i|i\in\llbracket2,q+1\rrbracket\}$. Thus the image of the exceptional divisor $\mathbb B^2_k$ (resp.\ $\mathbb E_k$) in $\SL_{2,k}/H^-$ (resp.\ in $\mathbb Y^{3,-}_k$) is the disjoint union of $2$ (resp.\ of $2(p^q-1)$) irreducible divisors isomorphic to $B^-_k/H^-\cong\mathbb A^1_k\times_{\Spec k} \mathbb G_{\m,k}$ (resp.\ to $\mathbb A^2_k$) through isomorphisms defined over $\mathbb F_{p^q}$, and we similarly compute
$$\zeta(\mathbb Y^{3,-})=\zeta_{3,p^q}-\zeta_{2,p^q}-(\zeta_{2,p^q}-\zeta_{1,p^q})+2(p^q-1)\zeta_{2,p^q}=\zeta_{3,p^q}+(2p^q-4)\zeta_{2,p^q}+\zeta_{1,p^q},$$
$$\zeta(\mathbb Y^{2,-})=\zeta_{2,p^q}+(2p^q-4)\zeta_{2,p^q}+\zeta_{0,p^q}.$$ 
In particular, if $p^q=2$, then we have
$$\zeta(\mathbb X^2)=\zeta_{2,2}+2\zeta_{0,2}\;\;\;\textup{and}\;\;\;\zeta(\mathbb Y^{2,-})=\zeta_{2,2}+\zeta_{0,2}.$$
Based on analogous computations over finite field extensions of $\mathbb F_{p^q}=\mathbb F_2$, we check that neither $\mathbb X^2_k$ nor $\mathbb Y^{2,-}_k$ has an open subvariety isomorphic to $\mathbb A^2_k$. If $\mathbb A^2_k$ is an open subvariety of either $\mathbb X^2_k$ or $\mathbb Y^{2,-}_k$, then we can assume that it is defined over a finite field, and this would contradict the zeta function computation as the complement is either empty or of pure codimension $1$ by \cite{Gro6}, Exp.\ V, Ex.\ 3.4. 
\end{remark}

\begin{remark}\normalfont\label{R14}
As $B^-_k$ is defined by the equation $x=0$, from Isomorphism (\ref{EQ8}) and \cite{Ma}, Ch.\ 7, Thm.\ 20.2 we get that $\SL_{2,k}$ is factorial. Thus $\mathbb X^3_k$ is factorial iff $\mathcal I^{\ast}$ is a bijection on divisor groups and hence, based on the Isomorphism (\ref{EQ11}) and the formulas for the exceptional divisors, iff $p^q=2$. 

For $p^q=2$, as $H^{\star}$ is a $2$-group and $\mathbb X^3_k$ is factorial, we get that $\mathbb Y_k^{3,\star}$ is also factorial; similarly, as there exist no non-trivial homomorphisms $\mathbb G_{\a,k}\rightarrow\mathbb G_{\m,k}$ (see \cite{SGA3-2}, Exp.\ XVII, Prop.\ 2.4i)), we get that $\mathbb X_k^2$ and $\mathbb Y_k^{2,\star}$ are factorial. From this, end of Remark \ref{R13}, and \cite{Ru}, Thm.\ 2 we get that the Kodaira dimensions of $\mathbb X_k^2$ and $\mathbb Y_k^{2,-}$ are at least $0$ and thus (see \cite{Ru}, Lem.\ 1.8) there exists no separable dominant morphism $\mathbb A^2_k\rightarrow \mathbb X_k^2$ or $\mathbb A^2_k\rightarrow \mathbb Y_k^{2,-}$. Hence also, there exists no separable dominant morphisms $\mathbb A^N_k\rightarrow \mathbb X_k^3$ or $\mathbb A^N_k\rightarrow \mathbb Y_k^{3,-}$ with $N\ge 3$ an integer. Based on \cite{Mi1}, Thm.\ 1, we also get that there exist no non-trivial actions of $\mathbb G_{\a,k}$ on either $\mathbb X_k^2$ or $\mathbb Y_k^{2,-}$.
\end{remark}

\begin{remark}\normalfont\label{R15}
Assume $p^q=2$; so $-1=1$. Let $(t,u,v)\in k_{\s}[w,x,z]_v^3$ be the triple of invariant polynomials that show up in Equation (\ref{EQ16}). Viewing the isomorphism 
$$(\mathcal T^{\#}_v)^{-1}:k[x_1,x_2,x_3]_{x_2}\rightarrow k_{\t}[w,x,y,z]_v/(wz+xy+1)$$ 
and the localization with respect to $v$ of the $k$-algebra homomorphism 
$$\mathcal L^{\#}:k_{\t}[w,x,y,z]/(wx+yz+1)\rightarrow k_{\s}[w,x,y,z]/(wx+yz+1)$$ 
that defines $\mathcal L$ as identifications, then, as $xy=zw+1$, we get that $x_1$, $x_2$, and $x_3$ are equal to $\frac{(wz+1)x+zw^2+1}{v}$, the second Dickson invariant $v$, and $\frac{wz^2+\frac{(wz+1)^2}{x}+1}{v}$ (respectively). Solving for $z$ (resp.\ $z^2$) in terms of $x_1$ (resp.\ $x_3$), $x$ and $w$, we get that $z=\frac{x+x_1x_2+1}{w^2+wx}$ and $z^2=\frac{x(x_2x_3+1+x^{-1})}{w^2+wx}$. From these identities we get that $\frac{(x+x_1x_2+1)^2}{w^2+wx}=x(x_2x_3+1+x^{-1})$, hence $(x+x_1x_2+1)^2=x_2(x_2x_3+1+x^{-1})$ or $x^3+(x_1^2x_2^2+x_2^2x_3+x_2+1)x+x_2=0$. Based on this and $u=w^2+wx+x^2=\frac{x_2}{x}+x^2$, the first Dickson invariant is $u=x_1^2x_2^2+x_2^2x_3+x_2+1$. The two finite \'etale covers $\mathbb Y^{3,+}\rightarrow\mathbb A^3_{\mathbb F_2}$ and $\mathbb Y^{3,-}\rightarrow\mathbb A^3_{\mathbb F_2}$ are defined by the normalization of $\mathbb F_2[x_1,x_2,x_3]$ in the two field extensions $\mathbb F_2(x_1,x_2,x_3)\rightarrow\mathbb F_2(x_1,x_2,x_3)[w]/(w^3+uw+v)$ and $\break\mathbb F_2(x_1,x_2,x_3)\rightarrow\mathbb F_2(x_1,x_2,x_3)[x]/(x^3+ux+v)$ (respectively). Equation (\ref{EQ16}) implies that the finite Galois cover $\mathbb V^3\rightarrow\mathbb A^3_{\mathbb F_2}$ is defined by the normalization of $\mathbb F_2[x_1,x_2,x_3]$ in the field extension 
$$\mathbb F_2(x_1,x_2,x_3)\rightarrow\mathbb F_2(x_1,x_2,x_3)[t]/(t^2+vt+v^2+u^3).$$ 
Hence 
$$\mathbb V^3\cong\Spec\mathbb F_2[x_1,x_2,x_3][t_v]/(t_v^2+t_v+1+u_v)$$ 
with $u_v:=\frac{u^3+v+1}{v^2}\in\mathbb F_2[x_1,x_2,x_3]$ and indeterminate $t_v:=\frac{t+1}{v}$.\end{remark}

\section{On finite \'etale covers of affine spaces and applications}\label{S20}

There exist many finite \'etale morphisms $\Phi:X\rightarrow\mathbb A^2_k$ with $X$ a connected affine smooth surface over $k$ which either is isomorphic to $\mathbb A^2_k$ (e.g., $\Phi:=d\times 1_{A^1_k}$ with $d\in\EE_1(k)$) or does not have an open subvariety isomorphic to $\mathbb A^2_k$ (e.g., $k^{\ast}\subsetneq\mathcal O(X)^{\ast}$ or is the pullback via a projection $\mathbb A^2_k\rightarrow\mathbb A^1_k$ of a finite \'etale cover $C\rightarrow\mathbb A^1_k$ with $C$ a non-rational curve). E.g., we have the following general result.

\begin{lemma}\label{PR8-}
Let $l\in\mathbb N^{\ast}$. If $p\ge 5$ (resp.\ $p\le 3$) we assume that $l\ge p$ (resp.\ $l\ge p+1$). Then there exists a non-Galois finite \'etale cover $c_n:\mathbb H^n\to\mathbb A^n_k$ of degree $l$.
\end{lemma}

\begin{proof}
We can assume that $n=1$ as we can take $c_{1+i}:=c_1\times 1_{\mathbb A^i_k}$ for $i\in\mathbb N^{\ast}$. Let $C\to\mathbb A^1_k$ be a connected finite \'etale cover which is Galois of Galois group $G$. If $G$ has a non-normal subgroup $H$ of index $l$, then the finite \'etale cover $C/H\to\mathbb A^1_k$ is non-Galois of degree $l$. We can take $G$ to be the alternating group $A_l$ (resp.\ symmetric group $S_l$) by \cite{Abh}, Sect.\ 13 First (resp.\ Second) Cor. As the subgroup $A_{l-1}$ of $A_l$ (resp.\ $S_{l-1}$ of $S_l$) is non-normal, by choosing $H$ to be the subgroup $A_{l-1}$ (resp.\ $S_{l-1}$) the lemma holds.\end{proof}

If $m\in\mathbb N^{\ast}$, $f_1(x,y,z)\in k[x,y,z]$ is a monic polynomial in $x$ of degree $m$ and $g_1(y,z)\in k[y,z]$, then for the $k$-algebra $A:=k[x,y,z]/\bigl(f_1(x^p,y,z)[x+g_1(y,z)]-1\bigr)$ we have $k^{\ast}\subsetneq A^{\ast}$ and the projection on the last two coordinates defines a finite \'etale morphism $\Spec A\rightarrow\mathbb A^2_k$ of degree $pm+1$. 

To provide examples of finite \'etale morphisms $\Phi:X\rightarrow\mathbb A^2_k$ with $\mathbb A^2_k\subsetneq X$ we recall that for $(f,g)\in\Theta_k$ and a finite additive subgroup $G$ of $k$ that contains the zeros of $f(t)$, if $m\in\mathbb N$ is such that $|G|=pm$ and $h_G(t):=\prod_{\alpha\in G,\alpha f(\alpha)\neq 0} (t-\alpha)$, then the product $tf(t)h_G(t)$ is a monic additive polynomial (e.g., see \cite{Go}, Ch.\ 1, Thm.\ 1.2.1); being also separable, we have $(fh_G,g)\in\Theta^1_{k,m}$.

\begin{proposition}\label{PR8}
Suppose that $n\ge 2$. Let $(f,g)\in\Theta_k$, $m\in\mathbb N^{\ast}$, and $h(t)\in k[t]$ be such that $\deg(f)=n-1$ and $(fh,g)\in\Theta^1_{k,m}$. Let $q\in\mathbb N^{\ast}$ be such that $q>\lfloor\frac{m}{n}\rfloor$ and let $(\beta_1,\ldots,\beta_q)\in k^{q-1}\times k^{\ast}$. Then the following properties hold.

\medskip
{\bf (1)} The rule on valued points 
$(x,y,z)\mapsto\bigl(y,h(x)z+\sum_{i=1}^q\beta_iz^{pi}\bigr)$ defines a finite \'etale cover 
$$\Phi_{f,g,h,q}:\mathbb S_{f,g,k}^2\rightarrow\mathbb A^2_k$$ 
of degree $pnq$ and hence the ring $k[x,y,z]/\bigl(xf(x)+yg(y)z\bigr)$ has a $p$-basis.

\smallskip
{\bf (2)} The composite of the open embedding $\imath_{f,g,k}$ of Displays (\ref{EQ0a}) and (\ref{EQ0b}) (resp.\ $\jmath_{f,g,k;\delta}$ of Proposition \ref{PR4.5} if $g(0)\neq 0$, with $\delta\in k^{\ast}$ such that $f\bigl(\delta g(0)\bigr)\neq 0$) with the morphism $\Phi_{f,g,h,q}$ is a surjective endomorphism $e_{f,g,h,q}\in\EE_2(k)$ (resp.\ $d_{f,g,h,q}\in\EE_2(k)$). For $\star=e$ (resp.\ $\star=d$) we have an isomorphism $X_{\star}\cong\mathbb S^2_{f,g,k}$ and identities $\deg(\star_{f,g,h,q})=pnq$ and $\rho_{\et}(\star_{f,g,h,q})=\rho(\star_{f,g,h,q})=(n-1)\z(tg)$.\end{proposition}

\begin{proof}
It is easy to see that $\Phi_{f,g,h,q}$ is generically \'etale. Hence it suffices to show that for each $(\alpha,\beta)\in k^2$, the solution set $\mathcal J_{\alpha,\beta}$ in $k^3$ of the system of equations
$$y-\alpha=h(x)z+\left(\sum_{i=1}^q\beta_iz^{pi}\right)-\beta=xf(x)+yg(y)z=0$$
has $pnq$ solutions. As $(fh,g)\in\Theta_{k,m}^1$, there exists $(\alpha_0,\alpha_1,\ldots,\alpha_{m-1})\in k^{\ast}\times k^{m-1}$ such that $tf(t)h(t)=\alpha_0t+[\sum_{i=1}^{m-1} \alpha_it^{pi}]+t^{pm}$. So $[tf(t)h(t)]'=\alpha_0\in k^{\ast}$ and $tf(t)$ and $h(t)$ do not have a common root. If $\alpha g(\alpha)=0$, then $xf(x)=0$, $h(x)\neq 0$, and
$$\mathcal J_{\alpha,\beta}=\Biggl\{(x,\alpha,z)\in k^3\Bigl|xf(x)=0=-\beta+h(x)z+\sum_{i=1}^q\beta_iz^{pi}\Biggr\}$$ 
has $pnq$ elements as all fibers of the surjective function $\mathcal J_{\alpha,\beta}\rightarrow\{x\in k|xf(x)=0\}$ have cardinality $pq$. If $\alpha g(\alpha)\neq 0$, then
$$\mathcal J_{\alpha,\beta}=\Biggl\{\Bigl(x,\alpha,\frac{xf(x)}{-\alpha g(\alpha)}\Bigr)\in k^3\Bigl|\frac{x^{pm}+\alpha_0x+\sum_{i=1}^{m-1} \alpha_ix^{pi}}{-\alpha g(\alpha)}+\sum_{i=1}^q\frac{\beta_ix^{pi}f(x)^{pi}}{[-\alpha g(\alpha)]^{pi}}=\beta\Biggr\}.$$
has $pnq$ elements as $\alpha_0\in k^{\ast}$ and $\deg\bigl(\beta_qx^{pq}f(x)^{pq}\bigr)=pnq>pm$ as $q>\lfloor\frac{m}{n}\rfloor$. So part (1) holds.

Part (2) follows from part (1) and Section \ref{S2} (resp.\ Proposition \ref{PR4.5}); the surjectivity part follows from the fact that for each $(\alpha,\beta)\in k^2$, the solution set $\mathcal J_{\alpha,\beta}$ is not contained in $\mathbb D^1_{f,g,k}=\mathbb S_{f,g,k}^2\setminus\Imm(\imath_{f,g,k})$ (resp.\ $\mathbb S_{f,g,k}^2\setminus\Imm(\jmath_{f,g,k;\delta})$) of Display (\ref{EQ0c}) (resp.\ the proof of Proposition \ref{PR4.5}).
\end{proof}

\begin{example}\normalfont\label{EX19}
Suppose that $p$ does not divide $n-1$. Thus for each monic polynomial $g(t)\in k[t]$ the surface $\mathbb S^2_{t^{n-1}-1,g,k}$ is smooth over $\Spec k$ and there exists $(l,m)\in (\mathbb N^{\ast})^2$ such that $pm=l(n-1)+1$. For $h(t):=\sum_{i=0}^{l-1} t^{(n-1)i}\in k[t]$ we have $(t^n-t)h(t)=t^{pm}-t$. Let $q\in\mathbb N^{\ast}$ be such that $q>\lfloor\frac{m}{n}\rfloor$. The rule on valued points 
$$(x,y,z)\mapsto\bigl(y,h(x)z+z^{pq}\bigr)$$ 
defines a finite \'etale morphism $\Phi_{n,m,q,g,k}:\mathbb S_{t^{n-1}-1,g,k}^2\rightarrow\mathbb A^2_k$ of degree $pnq$ by Proposition \ref{PR8}(1). Let $e=e_{n,m,q,g,k}:=\Phi_{n,m,q,g,k}\circ\imath_{t^{n-1}-1,g,k}\in\EE_2(k)$; therefore $\deg(e)=pnq$, $X_e=\mathbb S^2_{t^{n-1}-1,g,k}$ and $\rho_{\et}(e)=\rho(e)=(n-1)\z(tg)$ by Proposition \ref{PR8}(2). Moreover, $d_{n,m,q,k}:=\Phi_{n,m,q,1,k}\circ\Inv_{\mathbb S^2_{t^{n-1}-1,1,k}}\circ\imath_{t^{n-1}-1,1,k}\in\EE_2(k)$ has same listed invariants as $e_{n,m,q,1,k}$, where $\Inv_{\mathbb S^2_{t^{n-1}-1,1,k}}$ is the standard involution of $\mathbb S^2_{t^{n-1}-1,1,k}$ defined by the rule $(x,y,z)\mapsto (x,z,y)$. 

E.g., if $p=3$, then $d_{2,1,1,k}$ is defined by the rule $(x,y)\mapsto (x-x^2y,y+y^2x+y^3)$.

Let $l\in\mathbb N$. We take $g$ such that $\z(tg)=l$. As $\mathbb D^1_{f,g,k}$ in $\mathbb S^2_{f,g,k}$ is the zero locus $f(x)=yg(y)=0$ (see the line after Display (\ref{EQ0c})), $N_e$ is the zero locus $yg(y)=0$ and hence it is smooth and isomorphic to $l$ copies of $\mathbb A^1_k$. In particular, if $l=1$ (i.e., if $g(t)=t^s$ with $s\in\mathbb N$), then $N_e\cong\mathbb A^1_k$ is smooth and integral and therefore the Nollet--Xavier Conjecture over $\mathbb C$ (see \cite{Je}, Sect.\ 1 and \cite{NX}, Quest.\ 6) becomes false over $k$ even for $n=2$.\end{example}

\begin{remark}\normalfont\label{R16}
We refer to Example \ref{EX11} with $mq>2$. As $\rho(e)-\rho_{\et}(e)=1$, $Y_e:=X_e\setminus X_e^{\et}$ is irreducible, i.e., $Y_e\in\Irr(\Gamma_{\psi_e})$. The union $\mathbb S^2_e:=\mathbb S^2_{Y_e}=\mathbb A^2_{k,\s}\cup Y_e$ is an affine open subvariety $\Spec A_{Y_e}$ of $X_e$ (see Corollary \ref{C2.6}(2)) which is smooth over $\Spec k$ (as $\psi_e$ is regular) and for which $Y_e\cong\mathbb A^1_k$ (see Corollary \ref{C2.6}(1)) and $A_{Y_e}$ does not have a $p$-basis (see Corollary \ref{C2.8}(2)). In particular, there exists no \'etale morphism from $\mathbb S^2_e$ to affine surfaces over $k$ that have $p$-bases such as $\mathbb S^2_{f,g,k}$ with $(f,g)\in\Theta_k$ by Proposition \ref{PR8}(1) if $\deg(f)\ge 1$ and by $\mathbb S^2_{1,g,k}\cong\mathbb A^2_k$ if $\deg(f)=0$ (i.e., if $f=1$).
\end{remark}

\begin{corollary}\label{C8}
For $s\in \mathbb N^{\ast}$ we consider the set 
$$\mathbb O_{p,s}:=\{\rho(e)|e\in\EE_2(k),\rho(e)=\rho_{\et}(e),\deg(e)=ps\}.$$
Then the following properties hold.

\medskip
{\bf (1)} If $s$ is even, then $\mathbb O_{p,s}=\mathbb N$.

\smallskip
{\bf (2)} If $s$ odd, then $(ps-1)\mathbb N\subset\mathbb O_{p,s}$ and for each divisor $j$ of $s$ such that $p\nmid j-1$ we have $(j-1)\mathbb N\subset\mathbb O_{p,s}$ (e.g., $(2s-1)\mathbb N\subset\mathbb O_{2,s}$, and $2\mathbb N\subset\mathbb O_{p,s}$ if $p\ge 3$ and $3|s$).
\end{corollary}

\begin{proof}
As there exists $d\in\EE_2(k)$ which is finite \'etale of degree $ps$, $0\in \mathbb O_{p,s}$. The inclusion $(ps-1)\mathbb N^{\ast}\subset\mathbb O_{p,s}$ follows from the existence of the finite \'etale endomorphisms $c_{t^{ps-1}-1,g,k}\in\EE_2(k)$ (see Section \ref{S2}). 

If $s$ is even, in Example \ref{EX19} we take $(n,l,m,q)=(2,p-1,1,\frac{s}{2})$; for $e:=e_{2,1,\frac{s}{2},g,k}$ we have $\deg(e)=ps$ and $\rho(e)=\rho_{\et}(e)=\z(tg)\in\mathbb N^{\ast}$, hence $\mathbb N^{\ast}\subset\mathbb O_{p,s}$.

We are left to prove that if $s$ and $j$ are as in part (2), then $(j-1)\mathbb N^{\ast}\subset\mathbb O_{p,s}$; we have $j\ge 3$. In Example \ref{EX19} we take $(n,q)=(j,\frac{s}{j})$ and the unique $m\in\llbracket1,j-2\rrbracket$ such that there exists $l\in\mathbb N$ with $pm=l(j-1)+1$. For $e:=e_{j,m,\frac{s}{j},g,k}$ we have $\deg(e)=ps$ and $\rho(e)=\rho_{\et}(e)=(j-1)\z(tg)\in (j-1)\mathbb N^{\ast}$, hence $(j-1)\mathbb N^{\ast}\subset\mathbb O_{p,s}$.\end{proof}%there exists $e_r\in\EE_2(k)$ with $\rho(e_r)=\rho_{\et}(e_r)=r$. For $r=0$ see Remark \ref{R10}. For $r\in\mathbb N^{\ast}\setminus p\mathbb N^{\ast}$ see Example \ref{EX19} applied to $n=r+1$ and $f(y)=y$. For $r\in p\mathbb N^{\ast}$ see Theorem \ref{T1}(2) applied to $(m,s)=(\Frac(R){p},1)$

\begin{example}\normalfont\label{EX20-}
For $i\in \{1,2\}$, let $(n_i,m_i)\in\mathbb N^2$ with $n_i\ge 2$. Let $(f_i,t^{m_i})\in\Theta_K$ with $\deg(f_i)=n_i-1$. The open embedding $\imath_{f_i,t^{m_i},K}:\mathbb A^2_K\rightarrow\mathbb S^2_{f_i,t^{m_i},K}$ is such that the complement $\mathbb S^2_{f_i,t^{m_i},K}\setminus\Imm(\imath_{f_i,t^{m_i},K})$ is a disjoint union of $n_i$ curves isomorphic to $\mathbb A^1_K$ (see Displays (\ref{EQ0a}) and (\ref{EQ0b}) for notation). Thus the rank $\rho_{\mathbb S^2_{f_i,t^{m_i},K}}=n_i-1$ of the N\'eron--Severi group of $\mathbb S^2_{f_i,t^{m_i},K}$ is independent of $m_i$. Considering the factorization in linear factors $f_i(t)=\prod_{j=1}^{n_i-1} (t-\alpha_{i,j})$ with $(\alpha_{i,1},\ldots,\alpha_{i,n_i-1})\in (K^{\ast})^{n_i-1}$, for $j\in\llbracket1,n_i-1\rrbracket$ let $f_{i,j}(t):=t-\alpha_{i,j}$. From Theorem \ref{T3}(2) we get that the tangent bundle $T_{\mathbb S^2_{f_i,t^{m_i},K}}$ is trivial. If $m_1\neq m_2$, then $\mathbb S^2_{t-1,t^{m_1},K}\not\cong\mathbb S^2_{t-1,t^{m_2},K}$ by \cite{GG}, Thm.\ B(i). We have an affine open cover $\mathbb S^2_{f_i,t^{m_i},K}=\cup_{j=1}^{n_i-1} \Imm(\imath_{f_{i,j},t^{m_i};f,K})$ with each $\Imm(\imath_{f_{i,j},t^{m_i};f,K})\cong\mathbb S^2_{t-1,t^{m_i},K}$; it consists of $n_i-1$ spheres (resp.\ almost spheres) if $m_i=0$ (resp.\ $m_i>0$). 

{\bf Case 1: $\chr(K)=0$.} Let $\beta_i:=\sum_{j=1}^{n_i-1} \alpha_{i,j}$. We consider the polynomial $g_i(t):=(t+\frac{\beta_i}{n_i})f_i(t+\frac{\beta_i}{n_i})\in K[t]$ of degree $n_i$; the coefficient of $t^{n_i-1}$ in $g_i$ is $0$. If $(n_1,m_1)\neq (n_2,m_2)$ with $m_1m_2>0$, then $\mathbb S^2_{f_1,t^{m_1},K}\not\cong\mathbb S^2_{f_2,t^{m_2},K}$ by \cite{M-L2}, Thm.\ 2 over $\mathbb C$. If $(n_1,m_1)=(n_2,m_2)$, then $\mathbb S^2_{f_1,t^{m_1},K}\cong\mathbb S^2_{f_2,t^{m_2},K}$ iff either $n_1=1$ or $n_1>1$ and there exists $\gamma\in K^{\ast}$ such that $g_2(t)=\gamma^{-n_1}g_1(\gamma t)$ by \cite{M-L2}, Thm.\ 2 over $\mathbb C$ if $m_1>0$ and by \cite{BD}, Thm.\ 5.4.5(1) if $m_1=0$.

{\bf Case 2: $\chr(K)=p$.} Let $e_i\in\EE_2(k)$ be such that $X_{e_i}\cong\mathbb S^2_{f_i,t^{m_i},K}$ by Proposition \ref{PR8}(2) applied to $\imath_{f_i,t^{m_i},K}$; we have $\rho(e_i)=\rho_{\et}(e_i)=n_i-1$ and $\bigl(\rho_1(e_i),\rho_2(e_i)\bigr)$ is $(n_i,0)$ if $m_i=0$ and is $(0,n_i)$ if $m_i>0$. \end{example}

\begin{example}\normalfont\label{EX20+}
Let $f(t)\in k[t]$ be monic separable of degree $l\in\mathbb N^{\ast}$; so for each $m\in\mathbb N^{\ast}$ we have $(f,t^m)\in\Theta_k$. Let $(m_1,m_2)\in (\mathbb N^{\ast})^2$. For $i\in \{1,2\}$, let $e_i\in\EE_2(k)$ be such that $X_{e_i}\cong \mathbb S^2_{f,t^{m_i},k}$ by Proposition \ref{PR8}(2) applied to $\imath_{f,t^{m_i},k}$; we have $\rho(e_i)=\rho_{\et}(e_i)=\rho_2(e_i)=l$ (the last equality by Example \ref{EX20-}). If $m_1\neq m_2$, then $X_{e_1}\not\cong X_{e_2}$ by \cite{GG}, Thm.\ B(i).\end{example}

\begin{example}\normalfont\label{EX21}
Let $l\in\mathbb N^{\ast}\setminus\{1\}$. Let $\mathcal P_l(K^{\ast})$ be the set of all finite subsets of $K^{\ast}$ of cardinality $l$. For $\mathcal B\in\mathcal P_l(K^{\ast})$ let $f_{\mathcal B}(t):=\prod_{\alpha\in\mathcal B} (t-\alpha)$. For $a\in\GA_1(K)$, let $a(\mathcal B):=\{a(\alpha)|\alpha\in\mathcal B\}$. For $\mathcal B_1,\mathcal B_2\in\mathcal P_l(K^{\ast})$, we have $\mathbb S^2_{f_{\mathcal B_1},1,K}\cong \mathbb S^2_{f_{\mathcal B_2},1,K}$ iff there exists $a\in\GA_1(K)$ such that $\mathcal B_2=a(\mathcal B_1)$ by \cite{BD}, Thm.\ 5.4.5(1). Let $\mathcal C_1,\mathcal C_2\in\mathcal P_l(K^{\ast})$ be such that $\mathcal C_2\neq a(\mathcal C_1)$ for all $a\in\GA_1(K)$. If $\chr(K)=p$, then for $i\in \{1,2\}$, let $e_i\in\EE_2(K)$ be such that $X_{e_i}\cong \mathbb S^2_{f_{\mathcal C_i},1,K}$ by Proposition \ref{PR8}(2) applied to $\imath_{f_{\mathcal C_i},1,K}$; we have $X_{e_1}\not\cong X_{e_2}$ and $\rho(e_i)=\rho_{\et}(e_i)=\rho_1(e_i)=l$ (the last equality by Example \ref{EX20-}).\end{example}

\begin{example}\normalfont\label{EX22}
Let $(f,g)\in\Theta_K$. Let $(g_1,g_2)\in K[t]^2$ be such that there exists $g_3(t)\in K[t]$ with $g_1(t)g_3(t)=tg(t)$ and $g_1\mid g_3$. Let 
$$F(x,y,z):=xf(x)+g_3(y)[g_1(y)z-xg_2(x)]\in K[x,y,z]$$
and $\mathbb S^2_{f,g;g_1,g_2,K}:=\Spec\bigl(K[x,y,z]/(F)\bigr)$. The zero locus $F=F_x=F_y=F_z=0$ in $\mathbb A^3_K$ is contained in the zero locus 
$$g_1(y)g_3(y)=xf(x)-xg_2(x)g_3(y)=f(x)+xf'(x)-g_3(y)[g_2(x)+xg_2'(x)]=0$$ and thus, as $g_1\mid g_3$, in the zero locus $xf(x)=f(x)+xf'(x)=yg(y)=0$ which is empty as $tf(t)$ is separable. Hence $\mathbb S^2_{f,g;g_1,g_2,K}$ is smooth. The rule 
$$(x,y)\mapsto \bigl(xyg(y),y,-xf(xyg(y))+xg_3(y)g_2(xyg(y))\bigr)$$ 
defines an open embedding
$\imath_{f,g;g_1,g_2,K}:\mathbb A^2_K\rightarrow \mathbb S^2_{f,g;g_1,g_2,K}$. Though $F$ is not in general of the form $x\tilde f(x)+y\tilde g(y)z$ with $(\tilde f,\tilde g)\in\Theta_K$ but even then it can happen that $\mathbb S^2_{f,g;g_1,g_2,K}$ is isomorphic to an open subvariety of $\mathbb S^2_{\tilde f,\tilde g,K}$ for some $(\tilde f,\tilde g)\in\Theta_K$.

If $\deg(f)=1$, $g(t)=t^q$ with $q\in\mathbb N$ and $\deg(g_2)=0$, then the ind-group scheme $\Aut(\mathbb S^2_{f,g;g_1,g_2,K})$ is a semidirect product of the ind-unipotent group scheme $K[t]$ and the cyclic group of order $2$ by \cite{Cr}, Cor.\ 4.3; in such a case $\mathbb S^2_{f,g;g_1,g_2,K}$ is not isomorphic to any affine surface with a non-trivial $\mathbb G_{\m,K}$-action.

As a concrete example, for the remaining part of the example we take $f(t)=t+1$, $g(t)=t$, $g_1(t)=g_3(t)=-t$, and $g_2(t)=-2$. Thus $$F(x,y,z)=x^2+x+y(2x-yz)=(x+y+1)(x+y)+y[y(-z-1)-1]$$ and hence 
$$\mathbb S^2_{f,g;g_1,g_2,K}=\mathbb S^2_{t+1,t;-t,-2,K}\cong\Spec\bigl(K[x,y,z]/x^2+x+y(yz-1)\bigr).$$ If $\chr(K)=2$, then $\mathbb S^2_{t+1,t;-t,-2,K}\cong\Spec\bigl(K[x,y,z]/x^2+x+y(yz+1)\bigr)$ is isomorphic to $\mathbb S^2_{t+1,-t,K}\cong\mathbb S^2_{t-1,t,K}$ as one can see using a substitution $x=w+y$ and hence it is not isomorphic to $\mathbb S^2_K$ by \cite{GG}, Thm.\ B(i). If $\chr(K)\neq 2$, then $\mathbb S^2_{t+1,t;-t,2,K}\not\cong\mathbb S^2_{t-1,t^m,K}$ for all $m\in\mathbb N$ by the last sentence of the prior paragraph.

The open embedding $\jmath_{t-1,t+1,K;1}:\mathbb A^2_K\rightarrow\mathbb S^2_{t-1,t+1,K}$ (see Proposition \ref{PR4.5}) is defined by the rule $(x,y)\mapsto \bigl(h(yh+1),yh,-x-h^2\bigr)$, where $h:=1+xy\in K[x,y]$. We have a disjoint union decomposition $\mathbb S^2_{t-1,t+1,K}\setminus\Imm(\jmath_{t-1,t+1,K;1})=Y_1\sqcup Y_2$, where $Y_1$ is the zero locus $x=z=0$ and $Y_2$ is the zero locus $x-1=y+1=0$ (see Cases 1 and 2 of the proof of Proposition \ref{PR4.5}). Let 
$$R:=K[x,y,z]/\bigl(x^2-x+y(y+1)z\bigr).$$ By definition, $\mathbb S^2_{t-1,t+1,K}=\Spec R$. So we can identify $R$ with the $K$-subalgebra $K[h(yh+1),yh,-x-h^2]$ of $K[x,y]$. Under this identification, we write the complement $\mathbb S^2_{Y_i}:=\mathbb S^2_{t-1,t+1,K}\setminus Y_{3-i}=\Spec S_i$ for $i\in\{1,2\}$. We have 
$$S_1=R_{h(yh+1)-1}\cap R_{yh+1}=K[h,yh,x]=K[x,h,yh]\cong K[w',x,y']/(w^2-w+xy')$$ 
under the substitution $(w',y'):=(h,-yh)$. Thus $\mathbb S^2_{Y_1}\cong\mathbb S^2_K$. Similarly, 
$$S_2=R_{x+h^2}\cap R_{h(yh+1)}=K[y,yh,x+h^2]$$
$$=K[y,y^2x,x+2xy+x^2y^2]\cong K[x',y,z]/\bigl(F(x',y,z')\bigr)$$ 
under the substitution $(x',z'):=(y^2x,x+2xy+x^2y^2)$. Thus $\mathbb S^2_{t+1,t;-t,-2,K}\cong \mathbb S^2_{Y_2}$; so $Cl(\mathbb S^2_{t+1,t;-t,-2,K})\cong\mathbb Z$ and $\mathbb S^2_{t+1,t;-t,-2,K}$ is isomorphic to an open subvariety of $\mathbb S^2_{t-1,t+1,K}$; thus, as $\{(\tilde f,\tilde g)\in\Theta_K|Cl(\mathbb S^2_{\tilde f,\tilde g,K})\cong\mathbb Z\}=\{(t-\alpha,t^m)|\alpha\in K^{\ast}, m\in\mathbb N\}$, for $\chr(K)\neq 2$ we have $\mathbb S^2_{t+1,t;-t,-2,K}\not\cong\mathbb S^2_{\tilde f,\tilde g,K}$ for each $(\tilde f,\tilde g)\in\Theta_K$ by the last paragraph and reasons of isomorphisms classes of divisor class groups.

Assume $\chr(K)=p$. Let $e\in\EE_2(K)$ be such that $\imath_e=\imath_{t-1,t+1,K}$ (see Proposition \ref{PR8}(2)). Then $\{\mathbb S^2_{Y_1},\mathbb S^2_{Y_2}\}$ is the canonical open affine cover of $X_e=\mathbb S^2_{t-1,t+1,K}$ defined by $e$ (see Definition \ref{D7}(1)). From the non-isomorphisms of the last two paragraphs we get that $\rho_1(e)=\rho_2(e)=1$ and $\rho(e)=\rho_{\et}(e)=2$; if $p\neq 2$, then $\mathbb S^2_{Y_1}$ is an almost sphere that did not show up priorly in the paper.\end{example}

\begin{example}\normalfont\label{EX22.1}
Let $q\in \mathbb N^{\ast}\setminus\{1\}$. Let $\Sigma:\mathbb S^2_{t-1,t^{q-1},k}\rightarrow\mathbb P^1_k$ be the $\mathbb A^1$-fibration defined by the rule rule $(x,y,z)\mapsto (-x:y^q)$ or $(z:x-1)$ (whichever is defined; if both are defined, then they are equal). For each standard affine pseudo-plane $\mathbb S^2_{t-1,t^{q-1},k;\Sigma,P}$ of $\mathbb S^2_{t-1,t^{q-1},k}$ of type $q$ in the sense of Definition \ref{D6}(2), there exist \'etale morphisms $\mathbb S^2_{t-1,t^{q-1},k;\Sigma,P}\rightarrow\mathbb A^2_k$ that are composites of the open embedding $\mathbb S^2_{t-1,t^{q-1},k;\Sigma,P}\rightarrow\mathbb S^2_{t-1,t^{q-1},k}$ with finite \'etale covers $\mathbb S^2_{t-1,t^{q-1},k}\rightarrow\mathbb A^2_k$ as in Proposition \ref{PR8}(2). 

For the remaining part of the example we take $P=(0:1)\in\mathbb P^1_k(k)$. Then $\Sigma^{-1}(P)$ is the zero locus $x=z=0$ and the action of $\mathbb G_{\m,k}$ on $\mathbb S^2_{t-1,t^{q-1},k}$ defined by the rule $\alpha\cdot(x,y,z)=(x,\alpha y,\alpha^{-q}z)$ for $\alpha\in k^{\ast}$ restricts to an action of $\mathbb G_{\m,k}$ on $\mathbb S^2_{t-1,t^q,k;\Sigma,P}$; so $\mathbb S^2_{t-1,t^q,k;\Sigma,P}$ is of type $(q,q)$ in the sense of \cite{MM1}, Def. \ 2.1 (3) (cf.\ \cite{MM2}, Thm.\ 1.1). The case $q\ge 3$ (resp.\ $q=2$) is in contrast to the situation in characteristic $0$ in which, with $\chr(K)=0$, for an affine pseudo-planes $X_{q,r}$ over $K$ with non-trivial $\mathbb G_m$ action of type $(q,r)$ with $r\in\mathbb N^{\ast}$ there exists no \'etale morphism from $X_{q,r}$ to $\mathbb A^2_K$ if $r\neq 2$ and $d\nmid r-2$ by \cite{Mi5}, Thm.\ 3.6 (resp.\ the classical Jacobian Conjecture for $n=2$ implies that there exists no \'etale morphism from $X$ to $\mathbb A^2_K$ even if either $r=2$ or $d\mid r-2$; cf.\ proof of \cite{Mi5}, Thm.\ 3.3). If there exists $l\in\mathbb N^{\ast}$ such that $r=ql$, then $X_{q,ql}$ is the zero locus $x^q-y(1+y^lz)$ in $\mathbb A^3_K$ by \cite{DP}, Eq.\ (5.1). Using this for $l=1$, if $p\nmid q$ one similarly gets that $\mathbb S^2_{t-1,t^q,k;\Sigma,P}$ is isomorphic to the zero locus $x^q-y(1+yz)$ in $\mathbb A^3_k$. 
\end{example}

\begin{example}\normalfont\label{EX22.2}
Let $(m,q)\in (\mathbb N^{\ast}\setminus\{1\})^2$ be such that $\chr(K)\nmid m$. We define $s:=\min(m,q-1)\in\mathbb N^{\ast}$. Let $\underline{\alpha}:=(\alpha_2,\ldots,\alpha_s)\in K^{s-1}$. Let $\mathbb S^2_{m,q,\underline{\alpha}}$ be the zero locus $x^m+(\sum_{i=2}^s \alpha_ix^{m-i}y^i)-1+y^qz=0$ in $\mathbb A^3_K$; it is a smooth hypersurface in $\mathbb A^3_K$ which in the case $K=\mathbb C$ are generalized Danielewski surfaces studied in \cite{MM1}, Sect.\ 2. Following loc.\ cit., the rule $(x,y,z)\mapsto y$ defines an $\mathbb A^1$-fibration $\Sigma:\mathbb S^2_{m,q,\underline{\alpha}}\rightarrow\mathbb A^1_K$ with all fibers smooth and with only one fiber which is not irreducible. The non-irreducible fiber is $\Sigma^{-1}(0)$ and is a disjoint union of $m$ copies of $\mathbb A^1_K$. The cyclic group $G:=\{\alpha\in K^{\ast}|\alpha^m=1\}$ of order $m$ acts freely on $\mathbb S^2_{m,q,\underline{\alpha}}$ via the rule $\alpha\cdot(x,y,z)=(\alpha x,\alpha y,\alpha^{-q}z)$ in a way compatible with the $\mathbb A^1$-fibration. Thus $\mathbb S^2_{m,q,\underline{\alpha}}/G\rightarrow \mathbb A^1_K/G\cong\mathbb A^1_K$ is a $\mathbb A^1$-fibration with all fibers irreducible and with only one non-reduced fiber of multiplicity $m$. Hence $\mathbb S^2_{m,q,\underline{\alpha}}/G$ is an affine pseudo-plane over $K$ of type $m$. If $\chr(K)=0$, then $\mathbb S^2_{m,q,\underline{\alpha}}$ is simply connected, the $\mathbb A^1$-fibration is unique, and there exist open embeddings $\imath_{\omega}:\mathbb A^2_K\rightarrow\mathbb S^2_{m,q,\underline{\alpha}}$ indexed by $\omega\in G$ whose complements are disjoint unions of $m-1$ copies of $\mathbb A^1_K$ by \cite{MM1}, Lem.\ 2.6(1), (4), and (5) which is over $\mathbb C$; thus $\mathbb S^2_{m,q,\underline{\alpha}}$ is a sphere-like surface iff $m=2$. One can show that \cite{MM1}, Lem.\ 2.6(5) applies in general as $\chr(K)\nmid m$, giving open embeddings $\imath_{\omega}:\mathbb A^2_K\rightarrow\mathbb S^2_{m,q,\underline{\alpha}}$ whose complements are disjoint union of $m-1$ copies of $\mathbb A^1_K$; hence 
$$\bigl(x^m+(\sum_{i=2}^s \alpha_ix^{m-i}y^i)-1,y^q\bigr)\in\Smintr_2(K)\cap\Affintr_2(K).$$ 
If $\chr(K)=p$, $s\ge 2$, and $\underline{\alpha}$ is not the zero vector $(0,\ldots,0)$ of $K^{s-1}$, we do not know when there exists a finite \'etale cover $\mathbb S^2_{m,q,\underline{\alpha}}\rightarrow\mathbb A^2_K$. Let $f\in K[t]$ be such that $tf(t)=(t+1)^m-1$; so $\deg(f)=m-1$. We have $\mathbb S^2_{m,q,(0,\ldots,0)}\cong\mathbb S^2_{f,t^{q-1},K}$, hence the rule $\alpha\cdot(x,y,z)\mapsto (\alpha x+\alpha-1,\alpha y,\alpha^{-1}z)$ is a free action of $G$ on $\mathbb S^2_{f,t^{q-1},K}$. 

In this paragraph we assume that $\chr(K)=p$; so $p\nmid m$. Let $(l,s)\in (\mathbb N^{\ast})^2$ be such that $pl-1=ms$. Let $h(t):=(t+1)\sum_{i=0}^l (t+1)^{mi}\in K[t]$. Thus we have $tf(t)h(t)=(t+1)^{pm}-t-1$. Let $r\in\mathbb N^{\ast}$ be such that $r>\lfloor\frac{l}{m}\rfloor$ and $pr\equiv 1-q\;(\textup{mod}\; m)$. Let $G$ act on $\mathbb A^2_K$ via the rule $\alpha\cdot(y,z)\mapsto (\alpha y,\alpha^{1-q}z)$. The rule $(x,y,z)\mapsto (y,h(x)z+z^{pr})$ defines a finite \'etale $G$-invariant morphism $\Phi_{f,t^{q-1},h,r}:\mathbb S^2_{f,t^{q-1},K}\rightarrow\mathbb A^2_K$ of degree $pmr$ by Proposition \ref{PR8}(1). Thus we have a finite surjective morphism $\Phi_{f,t^{q-1},h,r}/G:\mathbb S^2_{f,t^{q-1},K}/G\rightarrow\mathbb A^2_K/G$ of degree $pmr$ which is \'etale outside $Y/G$, where $Y:=\{(0,0)\}$ if $m\nmid q-1$ and $Y:=\{0\}\times\mathbb A^1_K$ if $m\mid q-1$. If $m\mid q-1$, then $\mathbb A^2_K/G\cong\mathbb A^2_K$.
\end{example}

\section{Differential equations and applications}\label{S21}

In this section we generalize \cite{No}, Ex.\ (6.10) with $(a,b)\in (k^{\ast})^2$ and the recent work and differential equation of Lang (see \cite{Lang-J}, Props.\ 4.5 and 4.8) to all $n\ge 2$ in order to get more examples of endomorphisms $e\in\EE_n(k)$ that have interesting invariants and properties. We begin with a general lemma on determinants.

\begin{lemma}\label{L16} 
Let $R$ be an arbitrary commutative $\mathbb Z$-algebra. For $i\in \llbracket1,n\rrbracket$ let $(\alpha_i,\beta_i,\gamma_i)\in R^3$. Let $\triangle$ be the diagonal $n\times n$ matrix whose $i,i$ entry is $\alpha_i$ for $i\in \llbracket1,n\rrbracket$. Let $\triangle_1$ be the $n\times n$ matrix whose $i,j$ entry is $\beta_i\gamma_j$ for each $(i,j)\in \llbracket1,n\rrbracket^2$. Then the following identity holds 
$$\det(\triangle+\triangle_1)=\bigl(\prod_{i=1}^n \alpha_i\bigr)+\sum_{i=1}^n \bigl(\beta_i\gamma_i\prod_{j\in \llbracket1,n\rrbracket\setminus\{i\}} \alpha_j\bigr)\in R.$$
\end{lemma}

\begin{proof}
Note that $\triangle_1$ has rank at most $1$. For a subset $\mathbb O\subset \llbracket1,n\rrbracket$, let $\triangle_{1,\mathbb O}$ be the square matrix obtained from $\triangle_1$ by removing all $i$-th rows and columns with $i\in\mathbb O$; if $n-2\ge |\mathbb O|$, then $\det(\triangle_{1,\mathbb O})=0$. By convention $\det(\triangle_{1,\llbracket1,n\rrbracket}):=1$.

We compute
$$\det(\triangle+\triangle_1)=\sum_{\mathbb O\subset \llbracket1,n\rrbracket} \bigl(\det(\triangle_{1,\mathbb O})\prod_{i\in\mathbb O} \alpha_i\bigr)$$
$$=\sum_{\mathbb O\subset \llbracket1,n\rrbracket, |\mathbb O|\in\{n-1,n\}} \bigl(\det(\triangle_{1,\mathbb O})\prod_{i\in\mathbb O} \alpha_i\bigr)=\bigl(\prod_{i=1}^n \alpha_i\bigr)+\sum_{i=1}^n \bigl(\beta_i\gamma_i\prod_{j\in \llbracket1,n\rrbracket\setminus\{i\}} \alpha_j\bigr).$$
\end{proof}

\begin{corollary}\label{C8.1}
Let $n\ge 2$. Let $(l_1,\ldots,l_n)\in (\mathbb N^{\ast})^n$, $(f_1,\ldots,f_n)\in k[t]^n$, and $L\in k[x_1^p,\ldots,x_n^p]$. Let $w:=L(x_1,\ldots,x_n)\prod_{i=1}^n x_i^{l_i}\in k[x_1,\ldots,x_n]$. Then the determinant of the Jacobian matrix $J\bigl(x_1f_1(w),\ldots,x_nf_n(w)\bigr)(x_1,\ldots,x_n)$ is
$$\left(\prod_{i=1}^n f_i(w)\right)+\sum_{i=1}^n l_iwf_1(w)\cdots f_{i-1}(w)f_i'(w)f_{i+1}(w)\cdots f_n(w).$$
\end{corollary}

\begin{proof}
Let $(\alpha_i,\beta_i,\gamma_i):=\bigl(f_i(w),x_iwf_i'(w),l_ix_i^{-1}\bigr)\in \bigl(k[x_1,\ldots,x_n][\frac{1}{\prod_{i=1}^n} x_i]\bigr)^3$ for $i\in \llbracket1,n\rrbracket$. Then $\frac{\partial x_if_i(w)}{\partial x_j}$ is $\beta_i\gamma_j$ if $j\neq i$ and is $\alpha_i+\beta_i\gamma_i$ if $j=i$. Thus the corollary follows from Lemma \ref{L16} applied over $R=k[x_1,\ldots,x_n][\frac{1}{\prod_{i=1}^n} x_i]$ to the Jacobian matrix $J\bigl(x_1f_1(w),\ldots,x_nf_n(w)\bigr)(x_1,\ldots,x_n)=\triangle+\triangle_1$.\end{proof}
 
 We have the following generalization of \cite{Lang-J}, Prop.\ 4.8.
 
\begin{proposition}\label{PR8.5}
Let $(l_1,\ldots,l_n)\in\llbracket1,\ldots,p-1\rrbracket^n$, $L\in k_{\s}[x_1^p,\ldots,x_n^p]$, and $\alpha\in k^{\ast}$. Let $(f_1,\ldots,f_n)\in k[t]^n$ be such that the product $\prod_{i=1}^n f_i$ is a non-zero solution of the differential equation
\begin{equation}\label{EQ21-}
\Bigl(t\prod_{i=1}^n f_i^{l_i}\Bigr)'=\alpha\prod_{i=1}^n f_i^{l_i-1}.
\end{equation}
By defining $w:=L\prod_{i=1}^n x_i^{l_i}$, the rule 
$$(x_1,\ldots,x_n)\mapsto\bigl(x_1f_1(w),\ldots,x_nf_n(w)\bigr)$$
defines an \'etale endomorphism $e_{f_1,\ldots,f_n;L}^{l_1,\ldots,l_n}\in\EE_n(k)$ of Jacobian determinant $\alpha$.
\end{proposition}

\begin{proof}
Dividing Equation (\ref{EQ21-}) by $\prod_{i=1}^n f_i^{l_i-1}$, based on Corollary \ref{C8.1} we get that
$\det\bigl(J(x_1f_1(w),\ldots,x_nf_n(w))(x_1,\ldots,x_n)\bigr)$ is equal to 
$$\prod_{i=1}^n f_i(w)+\sum_{i=1}^n l_iwf_1(w)\cdots f_{i-1}(w)f_i'(w)f_{i+1}(w)\cdots f_n(w)=\alpha,$$
from which the proposition follows.
\end{proof}

For each non-zero solution $(f_1,\ldots,f_n)\in k[t]^n$ of Equation (\ref{EQ21-}), the polynomial $t\prod_{i=1}^n f_i(t)\in k[t]$ is separable; but only in the case $l_1=\cdots=l_n$ we are able to get a nice description of the solution set of Equation (\ref{EQ21-}) as follows.
 
\begin{example}\normalfont\label{EX23-}
Let $(l,\alpha)\in \llbracket1,p-1\rrbracket\times k^\ast$. Assume that $l_1=\cdots=l_n=l$. Let $s\in \llbracket1,p-1\rrbracket$ be such that $ls\equiv -1\;(\textup{mod}\; p)$. By defining $f:=\prod_{i=1}^n f_i\in k[t]$, Equation (\ref{EQ21-}) with $f\neq 0$ is equivalent to the differential equation
\begin{equation}\label{EQ21--}
f+ltf'=\alpha.
\end{equation}
As for each $m\in\mathbb N^{\ast}$, we have $x^m+lt(x^m)'=(1+lm)x^m$, the solution set of Equation (\ref{EQ21--}) is $\alpha+t^sk[t^p]$. Thus from Proposition \ref{PR8.5} we get that for each $n$-tuple $(f_1,\ldots,f_n)\in (k_{\s}[t])^n$ with $\prod_{i=1}^n f_i\in \alpha+t^sk[t^p]$ and every $L\in k_{\s}[x_1^p,\ldots, x_n^p]$, by defining $z:=\prod_{i=1}^n x_i$ and $w:=z^lL$, the rule
$$(x_1,\ldots,x_n)\mapsto\bigl(x_1f_1(w),\ldots,x_nf_n(w)\bigr)$$
defines an \'etale endomorphism $e_{f_1,\ldots,f_n;L}^l\in\EE_n(k)$ of Jacobian determinant $\alpha$.\end{example}

\begin{example}\normalfont\label{EX22.9}
In Example \ref{EX23-} we take $p>2$, $(n,l,s,\alpha,L)=(2,1,p-1,1,1)$, and a pair $(f_1,f_2)\in k[t]^2$ of monic polynomials such that $f_1(t)f_2(t)\in 1+t^{p-1}k[t^p]$ and $p\nmid \deg(f_1)\deg(f_2)[\deg(f_1)+1][\deg(f_2)+1]$. E.g., the last condition holds if $\deg(f_1)$ and $\deg(f_2)$ are congruent to $\frac{p-1}{2}$ modulo $p$. We have $w=z=x_1x_2$. As the leading terms of $x_1f_1(z)$ and $x_2f_2(z)$ are $x_1^{\deg(f_1)+1}x_2^{\deg(f_1)}$ and $x_1^{\deg(f_2)}x_2^{\deg(f_2)+1}$ (respectively), $d:=e^1_{f_1,f_2;1}=\e\bigl(x_1f_1(z),x_2f_2(z)\bigr)\in \EE_2(k)$ has $2$ points at infinity modulo $p$ in the sense of Definition \ref{D4.3}(5). Let $a:=\e(x_1-x_2^q,x_2)\in\GA_2(k)$ with $q\in\mathbb N^{\ast}\setminus\{1\}$ such that $p\nmid q+(q+1)\deg(f_1)$. So 
$$e:=da^{-1}=\e\bigl((x_1+x_2^q)f_1(x_1x_2+x_2^{q+1}),x_2f_2(x_1x_2+x_2^{q+1})\bigr)$$ 
has $1$ point at infinity modulo $p$ by Remark \ref{R3}(1) as the leading term of the product $(x_1+x_2^q)f_1(x_1x_2+x_2^{q+1})$ is $x_2^{q+(q+1)\deg(f_1)}$. As a contrast, we note that $ae$ has $1$ point at infinity modulo $p$ by Proposition \ref{PR1+}(1) while $d=ea$ has $2$ points at infinity modulo $p$. 
\end{example}

We generalizes the case $g=1$ of Example \ref{EX8} by elaborating on the case $(l,\alpha,L)=(1,1,1)$ of Example \ref{EX23-}.

%the existence of parameters, Galois groups generated by A_e, non-existence of intermediate fields.
\begin{example}\normalfont\label{EX23}
For $m\in\mathbb N^{\ast}$ and $(\alpha_1,\ldots,\alpha_m)\in k^{m-1}\times k^{\ast}$ we define polynomials $f(t):=1+\sum_{i=1}^m \alpha_it^{pi-1}$ and $F(t):=tf(t)$; we have $F'(t)=1$ and thus $F$ (or $f$) is separable. Assume $n\ge 2$ and let $z:=\prod_{i=1}^n x_i\in k_{\s}[x_1,\ldots,x_n]$. Let the $n$-tuple $(f_1,\ldots,f_n)\in k[t]^n$ be such that $f_1(0)=\cdots=f_n(0)=1$ and $\prod_{i=1}^n f_i(t)=f(t)$. Let 
$$\mathbb K^n_{f_1,\ldots,f_n}:=\Spec\bigl(k[x_0,x_1,\ldots,x_n]/(F(x_0)-z)\bigr).$$ The projection 
$$\pi_0:\mathbb K^n_{f_1,\ldots,f_n}\rightarrow\mathbb A^n_{k,\t}$$ on the last $n$ coordinates is a finite \'etale morphism of degree $pm$.

Let $\imath_{f_1,\ldots,f_n}:\mathbb A^n_{k,\s}\rightarrow\mathbb K^n_{f_1,\ldots,f_n}$ be the morphism defined by the rule 
$$(x_1,\ldots,x_n)\mapsto\bigl(z,x_1f_1(z),\ldots,x_nf_n(z)\bigr).$$ Clearly, $\imath_{f_1,\ldots,f_n}$ is birational; to show that it is in fact an open embedding, it suffices to show that for each $(\beta_0,\ldots,\beta_n)\in \mathbb K^n_{f_1,\ldots,f_n}(k)$, i.e., for each $(\beta_0,\ldots,\beta_n)\in k^{n+1}$ with $F(\beta_0)=\prod_{i=1}^n \beta_i$, the system $\mathcal S$ of equations
$$z-\beta_0=x_1f_1(z)-\beta_1=\cdots=x_nf_n(z)-\beta_n=0$$
in the indeterminates $x_1,\ldots,x_n$ has at most one solution in $k^n$. If $f(\beta_0)\neq 0$, then the only solution of $\mathcal S$ is $\bigl(\frac{\beta_1}{f_1(\beta_0)},\ldots,\frac{\beta_n}{f_n(\beta_0)}\bigr)$. Assume now that $f(\beta_0)=0$; thus $\beta_0\neq 0$. As the roots of $f$ are distinct, there exists a unique $j\in \llbracket1,n\rrbracket$ such that $f_j(\beta)=0$. Thus for $i\in \llbracket1,n\rrbracket\setminus\{j\}$ we have $x_i=\frac{\beta_i}{f_i(\beta_0)}$. If there exists $i\in \llbracket1,n\rrbracket\setminus\{j\}$ such that $\beta_i=0$, then $x_i=0$, hence $z=0$ which contradicts $\beta_0\neq 0$, therefore in such a case $\mathcal S$ inconsistent. If for each $i\in \llbracket1,n\rrbracket\setminus\{j\}$ we have $\beta_i\neq 0$, then $x_j=\frac{z}{\prod_{i\in \llbracket1,n\rrbracket\setminus\{j\}} x_i}=\beta_0\prod_{i\in \llbracket1,n\rrbracket\setminus\{j\}} \frac{f_i(\beta_0)}{\beta_i}$; as $F(\beta_0)=\prod_{i=1}^n \beta_i$ implies that $\beta_j=0$, $\mathcal S$ has indeed a unique solution. We conclude that $\imath_{f_1,\ldots,f_n}$ is an open embedding and that
$$\mathbb E^{n-1}_{f_1,\ldots,f_n}:=\mathbb K^n_{f_1,\ldots,f_n}\setminus\Imm(\imath_{f_1,\ldots,f_n})=\cup_{j\in \llbracket1,n\rrbracket}\cup_{\beta\in k,f_j(\beta)=0} \cup_{i\in \llbracket1,n\rrbracket\setminus\{j\}} \mathbb E^{n-1}_{\beta,i}$$
where $\mathbb E^{n-1}_{\beta,i}$ is the closed subvariety of $\mathbb K^n_{f_1,\ldots,f_n}$ (or of $\mathbb A^{n+1}_k$) defined by the equations $x_0-\beta=x_i=0$; we have $\mathbb E^{n-1}_{\beta,i}\cong\mathbb A^{n-1}_k$. It follows that
\begin{equation}\label{EQ21}
|\Irr(\mathbb E^{n-1}_{f_1,\ldots,f_n})|=\sum_{j=1}^n\deg(f_j)(n-1)=(n-1)\sum_{j=1}^n\deg(f_j)=(n-1)(pm-1).
\end{equation}

Thus
$$e_{f_1,\ldots,f_n}:=\pi_0\circ\imath_{f_1,\ldots,f_n}\in\EE_n(k)$$ 
is defined by the $n$-tuple $\bigl(x_1f_1(z),x_2f_2(z),\ldots,x_nf_n(z)\bigr)\in (k_{\s}[x_1,\ldots,x_n])^n$ and we have $\deg(e_{f_1,\ldots,f_n})=pm$ and $X_{e_{f_1,\ldots,f_n}}\cong\mathbb K^n_{f_1,\ldots,f_n}$. Also, $N_{e_{f_1,\ldots,f_n}}$ is contained in the zero locus $z=0$; therefore $e_{f_1,\ldots,f_n}$ becomes finite \'etale after inverting $z\in k_{\t}[x_1,\ldots,x_n]$. If $\deg(f_i)<pm-1$ for each $i\in\llbracket1,n\rrbracket$, then $N_{e_{f_1,\ldots,f_n}}$ is the zero locus $z=0$.

For each point $P$ of the zero locus $z=0$, we have $(0,P)\in\Imm(\imath_{f_1,\ldots,f_n})$, hence $P\in\Imm(e_{f_1,\ldots,f_n})$. The last two sentences imply that $e_{f_1,\ldots,f_n}$ is surjective. So $\varphi(e_{f_1,\ldots,f_n})=\iota(e_{f_1,\ldots,f_n})=0$. Note that $e_{f_1,\ldots,f_n}$ is $e_{f_1,\ldots,f_n;1}^1$ of Example \ref{EX23-} and has $1$ point at infinity modulo $p$ in the sense of Definition \ref{D4.3}(5) iff $n=2$ and $p\mid\deg(f_1)\deg(f_2)$; if $p$ does not divide $\prod_{i=1}^n \deg(f_i)[\deg(f_i)+1]$, then $e_{f_1,\ldots,f_n}$ has $n$ points at infinity modulo $p$. 

If $n\ge 3$, then, for distinct elements $i_1,i_2,j$ of $\llbracket1,n\rrbracket$ and $\beta\in k$ with $f_j(\beta)=0$, the intersection $\mathbb E^{n-1}_{\beta,i_1}\cap \mathbb E^{n-1}_{\beta,i_2}\cong\mathbb A^{n-2}_k$ is non-empty, in contrast to Corollary \ref{C2.6}(1).

We consider the particular case $n=2$. Let $(g,h):=(f_1,f_2)$. Equation (\ref{EQ21}) gives $\rho_{\et}(e_{g,h})=\rho(e_{g,h})=pm-1$. From the description of $\mathbb E^1_{g,h}$ we get that $e_{g,h}^{-1}(0,0)=\{(0,0)\}$ and that for each $\beta\in k^*$ we have 
$$e_{g,h}^{-1}(0,\beta)=\{\bigl(\frac{z h(z)}{\beta},\frac{\beta}{h(z)}\bigr)|z\in k,\, zg(z)=0\}$$ and 
$$e_{g,h}^{-1}(\beta,0)=\{\bigl(\frac{\beta}{g(z)},\frac{z g(z)}{\beta}\bigr)|z\in k,\, zh(z)=0\}.$$ 
If $1\notin\{g,h\}$, then $\mathcal F(e_{g,h})=\{1,\deg(g)+1,\deg(h)+1\}$. If $1\in\{g,h\}$, then $\mathcal F(e_{g,h})=\{1\}$. Also, from Proposition \ref{PR1}(1) to (3) we get that for $\star\in\{\l,\r\}$ we have
$$2\min\bigl(\deg(g),\deg(h)\bigr)+1\le\pi_{\i}(e_{g,h})\le\pi_{\star}(e_{g,h})=2\max\bigl(\deg(g),\deg(h)\bigr)+1;$$ so $\pi_{\l}(e_{g,h})=\pi_{\r}(e_{g,h})\ge 2\lceil\frac{pm-1}{2}\rceil+1$. E.g., if $p$ is odd and $\deg(g)=\deg(h)=\frac{pm-1}{2}$, then $\pi_{\i}(e_{g,h})=\pi_{\i}(e_{g,h})=\pi_{\l}(e_{g,h})=pm$. 

If $1\notin\{g,h\}$ let $\epsilon:=2$ and if $1\in\{g,h\}$ let $\epsilon:=1$; Proposition \ref{PR1}(3) also gives that $\max\bigl(\deg(g),\deg(h)\bigr)+\epsilon\le\pi_{\i}(e_{g,h})$. If $1\notin\{g,h\}$, then $N_e$ is the zero locus $z=0$ and hence we have $\nu(e)=2$.

To give a concrete general example independent of any product decomposition for $n=2$, in this paragraph we assume that $p$ is odd. Let $o\in\mathbb N^{\ast}$ be such that $p=2o+1$. Let $r\in\mathbb N$. We take $g(t)=1+f(t)$, $h(t)=1-f(t)$, where $f(t):=\sum_{i=0}^r \beta_it^{pi+o}\in k[t]$
with $(\beta_0,\beta_1,\ldots,\beta_r)\in k^r\times k^{\ast}$. Then $\deg(g)=\deg(h)=pr+o$, and 
$$tg(t)h(t)=t-t^p\Bigl(\sum_{i=0}^r \beta_it^{pi}\Bigr)^2=t+\sum_{i=1}^{2r+1} \alpha_it^{pi};$$
here $\alpha_1:=-\beta_0^2,\alpha_2:=-2\beta_0\beta_1,\ldots, \alpha_{2r+1}:=-\beta_r^2$ are elements of $k$ that are quadratic forms in $\beta_0,\ldots,\beta_r$ and $m=2r+1$. Denoting $e_f:=e_{1+f,1-f}$, we get that $\deg(e_f)=pm=p(2r+1)$ and $\mathcal F(e_f)=\{1,pr+o+1\}=\{1,\frac{p+1}{2}+pr\}$. Also, $e_f$ has $2$ points at infinity modulo $p$.
\end{example}%GROUP

\begin{remark}\normalfont\label{R17}
{\bf (1)} Referring to Example \ref{EX23}, $e_{f_1,\ldots,f_n}\in\EE_n(k)\setminus\GA_n(k)$ fixes all points $(\alpha_1,\ldots,\alpha_n)\in k^n$ with $\prod_{i=1}^n\alpha_i=0$, so it maps injectively each hyperplane $x_i=0$ onto a smooth closed hypersurface, which is in contrast to the $\chr(K)=0$ situation: if $d\in\EE_2(K)$ is injective on one line of $\mathbb A^2_{K,\s}$ (resp.\ maps one line to a smooth curve), then $d\in\GA_2(K)$ by \cite{Gw}, Thm.\ 1.1 (resp.\ \cite{YD}, Thm.\ 2.4).

\smallskip
{\bf (2)} Referring to $e_{f,g}$ of Example \ref{EX23}, if $1\notin\{g,h\}$ and $p\nmid m$ then $e_{f,g}$ is not a $p$-morphism by Remark \ref{R3}(6).
\end{remark}

Next we generalize \cite{Lang-J}, Prop.\ 4.5, Case $a\equiv b\equiv 1\;(\textup{mod}\; p)$ to all $n\ge 2$ by elaborating on the case $(l,\alpha,L)=(1,1,\prod_{i=1}^n x_i^{pq_i})$ with $(q_1,\ldots,q_n)\in\mathbb N^n$ of Example \ref{EX23-}.

\begin{example}\normalfont\label{EX24}
Let $n\ge 2$, $m\in\mathbb N^{\ast}$, $(\alpha_1,\ldots,\alpha_m)\in k^{m-1}\times k^{\ast}$, $f(t)$, $F(t)$, $f_1(t),\ldots,f_n(t)$, and $z=\prod_{i=1}^n x_i$ be as in Example \ref{EX23}. Recall that $e_{f_1,\ldots,f_n}$ is defined by the $n$-tuple $(g_1,\ldots,g_n)$, where $g_i(x_1,\ldots,x_n):=x_if_i(z)$ for $i\in \llbracket1,n\rrbracket$. 

Let $(q_1,\ldots,q_n)\in\mathbb N^n$. Let $L:=\prod_{i=1}^n x_i^{pq_i}\in k_{\s}[x_1^p,\ldots,x_n^p]$. For $i\in \llbracket1,n\rrbracket$ let $m_i:=pq_i+1$.  Let $w:=\prod_{i=1}^n x_i^{m_i}$; we have divisibilities $z|w|z^{\max(m_1,\ldots,m_n)}$ and $w=zL$. For $i\in \llbracket1,n\rrbracket$ let $h_i(x_1,\ldots,x_n):=x_if_i(w)$. 

Let $e_{f_1,\ldots,f_n}^{q_1,\ldots,q_n}:=\e(h_1,\ldots,h_n)\in\End_n(k)$. As $e_{f_1,\ldots,f_n}^{q_1,\ldots,q_n}$ is $e_{f_1,\ldots,f_n;L}^1$ of Example \ref{EX23-}, we have $e_{f_1,\ldots,f_n}^{q_1,\ldots,q_n}\in\EE_n(k)$. Note that $e_{f_1,\ldots,f_n}=e_{f_1,\ldots,f_n}^{0,\ldots,0}$.

As in the case of $e_{f_1,\ldots,f_n}$ (see Example \ref{EX23}), it is easy to see that $e_{f_1,\ldots,f_n}^{q_1,\ldots,q_n}$ becomes finite \'etale after inverting $z$ (equivalently, $w$).

Let
$$\mathbb K^{n;q_1,\ldots,q_n}_{f_1,\ldots,f_n}:=\Spec\Bigl(k[x_0,x_1,\ldots,x_n]/(F(x_0)\bigl(\prod_{i=1}^n f_i(x_0)^{pq_i}\bigr)-w)\Bigr).$$ 
The projection 
$$\pi_0:\mathbb K^{n;q_1,\ldots,q_n}_{f_1,\ldots,f_n}\rightarrow\mathbb A^n_{k,\t}$$ on the last $n$ coordinates is a finite flat morphism of degree $pm+\sum_{i=1}^n pq_i\deg(f_i)$. Let $\imath_{f_1,\ldots,f_n}^{q_1,\ldots,q_n}:\mathbb A^n_{k,\s}\rightarrow\mathbb K_{f_1,\ldots,f_n}^{n;q_1,\ldots,q_n}$ be the morphism defined by the rule 
$$(x_1,\ldots,x_n)\mapsto\bigl(w,h_1(x_1,\ldots,x_n),\ldots,h_n(x_1,\ldots,x_n)\bigr)=\bigl(w,x_1f_1(w),\ldots,x_nf_n(w)\bigr).$$ 
Clearly, $\imath_{f_1,\ldots,f_n}^{q_1,\ldots,q_n}$ is birational. As $e_{f_1,\ldots,f_n}^{q_1,\ldots,q_n}=\pi_0\circ\imath_{f_1,\ldots,f_n}^{q_1,\ldots,q_n}$, it follows that $\imath_{f_1,\ldots,f_n}^{q_1,\ldots,q_n}$ is quasi-finite, $\deg(e_{f_1,\ldots,f_n}^{q_1,\ldots,q_n})=pm+\sum_{i=1}^n pq_i\deg(f_i)$ and $\mathbb K^{n;q_1,\ldots,q_n}_{f_1,\ldots,f_n}$ is a monomial model of $X_{e_{f_1,\ldots,f_n}^{q_1,\ldots,q_n}}$. As in Example \ref{EX23} we argue that 
$$\mathbb K^{n;q_1,\ldots,q_n}_{f_1,\ldots,f_n}\setminus\Imm(\imath_{f_1,\ldots,f_n}^{q_1,\ldots,q_n})=\mathbb E^{n-1}_{f_1,\ldots,f_n}$$ 
is independent on $(q_1,\ldots,q_n)$ and therefore $\rho(e)\ge (n-1)(pm-1)$ by Equation (\ref{EQ21}) and Remark \ref{R8}(1). Moreover, $\imath_{f_1,\ldots,f_n}^{q_1,\ldots,q_n}$ is injective iff $q_1=\cdots=q_n=0$.\end{example}

\begin{example}\normalfont\label{EX25}
We take $q_1=\cdots=q_{n-1}=0$ and $q:=q_n\in\mathbb N^{\ast}$ and assume that $\deg(f_n)\ge 1$; let $\alpha_n\in k^{\ast}$ be such that $f_n(\alpha_n)=0$. We have 
$$p(m+q)\le pm+pq\deg(f_n)=\deg(e^{0,\ldots,0,q}_{f_1,\ldots,f_n})\le pm+pq(pm-1).$$ 
Using the substitution $v:=x_0-\alpha_n$, it follows that $$\mathbb K^{n;0,\ldots,0,q}_{f_1,\ldots,f_n}=\Spec\Bigl(k[v,x_1,\ldots,x_n]/\bigl(v^{pq+1}g(v)-x_n^{pq+1}\prod_{i=1}^{n-1} x_i\bigr)\Bigr),$$
with $g(v)\in k[v]$ such that $g(0)\neq 0$ and $\deg(g)=p[m+q\deg(f_n)-q]-1$. 

Let $\mathcal R_0$ be the normalization of the local ring of the hypersurface $\mathbb K^{n;0,\ldots,0,q}_{f_1,\ldots,f_n}$ at the point $(0,\ldots,0)\in\mathbb K^{n;0,\ldots,0,q}_{f_1,\ldots,f_n}(k)$; it is a localization of the $k_{\t}[x_1,\ldots,x_n]$-algebra $\mathcal O(X_{e^{0,\ldots,0,q}_{f_1,\ldots,f_n}})$. Let 
$$R_{n,q,k}:=k[[w,x_1,\ldots,x_n]]/\bigl(w^{pq+1}-\prod_{i=1}^{n-1} x_i\bigr).$$
As $g(0)\neq 0$ it is easy to see that we have a homomorphism $\mathcal R_0\rightarrow R_{n,q,k}$ of $k[x_1,\ldots,x_n]$-algebras that maps $v+\bigl(v^{pq+1}g(v)-x_n^{pq+1}\prod_{i=1}^{n-1} x_i\bigr)$ into a unit times $wx_n+(w^{pq+1}-\prod_{i=1}^{n-1} x_i)$ and this implies that the completion of $\mathcal R_0$ at a suitable maximal ideal of it is isomorphic as a $k[[x_1,\ldots,x_n]]$-algebra to $R_{n,q,k}$. Thus for $n\ge 3$ we get that $\mathcal R_0$ is a non-regular ring and hence $\psi_{e^{0,\ldots,0,q}_{f_1,\ldots,f_n}}$ is a non-regular Jacobian variety. If $n=2$, we only get that $\psi_{e^{0,q}_{f_1,f_2}}$ is non-\'etale and one can easily check that it is in fact regular; moreover, if $\deg(f_2)=1$, i.e., $f_2(t)=t-\alpha_2$, then $\mathbb K^{2;0,q}_{f_1,t-\alpha_2}$ is isomorphic to $\mathbb S^2_{pq+1,pq+1,(t+\alpha_2)f_1(t+\alpha_2),1,k}$ of Example \ref{EX10}.\end{example}

\begin{remark}\normalfont\label{R18}
The additive version of Example \ref{EX23} uses sums instead of products. If $n\ge 1$ and $f_1(t),\ldots,f_n(t)\in k[t]$ are such that $\sum_{i=1}^n f_i(t)=\sum_{i=0}^m \alpha_it^{pi}$ with $m\in\mathbb N^\ast$ and $(\alpha_1,\ldots,\alpha_m)\in k^{m-1}\times k^{\ast}$, then, by defining $z:=\sum_{i=1}^n x_i$, the $n$-tuple $\bigl(x_1+f_1(z),\ldots,x_n+f_n(z)\bigr)$ defines $e^{\a}_{f_1,\ldots,f_n}\in\EE_n(k)$ of Jacobian matrix 
$$\begin{bmatrix} 
1+f_1'(z) & f_1'(z) & \cdots & f_1'(z)\\
\cdots & \cdots & \cdots & \cdots\\
f_n'(z) & f_n'(z) & \cdots &1+f_n'(z) \\ 
\end{bmatrix}$$
that has determinant $1$ as one can easily check based on the fact that for each column, its entries add up to $1$. For $(\beta_1,\ldots,\beta_n)\in k^n$, the fiber $(e^{\a}_{f_1,\ldots,f_n})^{-1}(\beta_1,\ldots,\beta_n)$ has $pm$ elements as it equals
$$\{\bigl(\beta_1-f_1(\delta),\ldots,\beta_n-f_n(\delta)\bigr)|\delta\in k, \delta+\sum_{i=0}^m \alpha_i\delta^{pi}=\sum_{i=1}^n \beta_i\}.$$ Thus $\deg(e^{\a}_{f_1,\ldots,f_n})=pm$ and $e^{\a}_{f_1,\ldots,f_n}$ is finite \'etale; hence $\psi_{e^{\a}_{f_1,\ldots,f_n}}=e^{\a}_{f_1,\ldots,f_n}$ and $\rho_{\et}(e^{\a}_{f_1,\ldots,f_n})=\rho(e^{\a}_{f_1,\ldots,f_n})=\varphi(e^{\a}_{f_1,\ldots,f_n})=\iota(e^{\a}_{f_1,\ldots,f_n})=0$.
\end{remark}%GROUP3

\begin{remark}\normalfont\label{R19}
In Example \ref{EX23} we take $n=2$ and $\deg(f_1)\le\deg(f_2)$; the pair $(g_1,g_2):=\bigl(x_1f_1(x_1x_2),x_2f_2(x_1x_2)\bigr)$ defines $e_{f_1,f_2}\in\EE_2(k)$ with $\deg(e_{f_1,f_2})=pm$ and  inequalities
$$\deg(g_1)=2\deg(f_1)+1\le pm\le 2\deg(f_2)+1=\deg(g_2)=2pm-\deg(g_1);$$
note that $m\in\mathbb N^{\ast}$ and $\deg(g_1)\in\llbracket1,pm\rrbracket\cap (2\mathbb N^{\ast}-1)$ are arbitrary. 

Similarly, in Example \ref{EX24} we take $n=2$ with $\deg(f_1)\le\deg(f_2)$. The pair $(h_1,h_2)=\bigl(x_1f_1(x_1^{pq_1+1}x_2^{pq_2+1}),x_2f_2(x_1^{pq_1+1}x_2^{pq_2+1})\bigr)$ defines $e_{f_1,f_2}^{q_1,q_2}\in\EE_2(k)$ with $\deg(e_{f_1,f_2}^{q_1,q_2})=pm+pq_1\deg(f_1)+pq_2\deg(f_2)$ and we have inequalities
$$\deg(h_1)=(pq_1+pq_2+2)\deg(f_1)+1\le\deg(e_{f_1,f_2}^{q_1,q_2})\le (pq_1+pq_2+2)\deg(f_2)+1$$
$$=\deg(h_2)=2\deg(e_{f_1,f_2}^{q_1,q_2})-\deg(h_1)+p(q_1-q_2)\bigl(\deg(f_2)-\deg(f_1)\bigr).$$
If $p$, $q_1$, and $q_2$ are odd, then $\deg(h_1)$ and $\deg(h_1)$ are even.

Let $(q,s)\in (\mathbb N^{\ast}\setminus\{1\})\times \mathbb N^{\ast}$ and $l\in\llbracket ps,2ps-1\rrbracket$. The pair $(x_1+x_2^q,x_2+x_1^{ps})$ defines $d\in\EE_2(k)$ which is finite with $\deg(d)=pqs$ by Lemma \ref{F3}(2). If $a\in\GA_2(K)$ is defined by the pair $(x_1,x_2+x_1^l)$, then $e:=ad\in\EE_2(k)$ is finite with $\deg(e)=pqs$. Moreover, $ae$ is defined by the pair $(g_3,g_4):=\bigl(x_1+x_2^q,x_2+x_1^{ps}+(x_1+x_2^q)^l\bigr)$ for which we have inequalities
$$\deg(g_3)=q< pqs=\deg(e)\le ql=\deg(g_4)\le 2pqs-q;$$
if $q\in 2\mathbb N^{\ast}$, $l=2ps-1$, and $m=qs$ we get an even analog of the particular case $\deg(g_1)|m$ of the first paragraph.

The examples of the last three paragraphs are in contrast to the $\chr(K)=0$ situation: if $c=\e(g,h)\in\EE_2(K)$, then $\deg(c)\le\min\bigl(\deg(g),\deg(h)\bigr)$ by \cite{Zh}, Introd., Thm.\ or \cite{Kat}, Thm.\ 1.\end{remark}

\section{Types of endomorphisms}\label{S22}
 
We begin with an example that justifies the introduction of versions of types. 

\begin{example}\normalfont\label{EX26}
Let $m\ge 2$ and $q\ge 1$ be integers with $\chr(K)\nmid q$. The rule $(x,y)\mapsto x+x^my^q$ defines a smooth morphism $e:\mathbb A^2_K\rightarrow\mathbb A^1_K$ which is not a projection as $|\Irr\bigl(e^{-1}(0)\bigr)|\ge 2$; if $q=1$, then for each $\alpha\in K^{\ast}$, $e^{-1}(\alpha)\cong\mathbb A^1_K\setminus\{0\}$. 
\end{example}

Example \ref{EX26} implies that for each $(n,m)\in (\mathbb N^{\ast})^2$ there exist smooth morphisms $\mathbb A^{m+n}_{K,\s}\rightarrow\mathbb A^n_{K,\t}$ which are not projections.

The following types are used in connection to Conjecture \ref{CJ3} in Section \ref{S27}.

\begin{definition}\label{D11} Let $n\in\mathbb N^{\ast}$, $e\in\D_n(K)$, and $m\in \llbracket0,n\rrbracket$. We say that $e$ has type (resp.\ strong type or weak type) $\le m$ if there exists a smooth morphism (resp.\ a projection or a morphism) $\pi:\mathbb A^n_{K,\s}\rightarrow\mathbb A^m_K$ such that the product morphism $\pi\times e:\mathbb A^n_{K,\s}\rightarrow\mathbb A^m_{K,\t}\times_{\Spec K} \mathbb A^n_K$ is quasi-finite, is birational to the Zariski closure $\overline{\Imm(\pi\times e)}$ of its image, and has the property that for each irreducible component $Y$ of the closed subvariety $(\pi\times e)^{-1}\bigl(\Sing(\overline{\Imm(\pi\times e)}\bigr)$ of $\mathbb A^n_{K,\s}$ of dimension $n-1$ we have $e^{-1}\bigl(e(Y)\bigr)=Y$. We say that $e$ has type (resp.\ strong type or weak type) $m$ if it has type (resp.\ strong type weak type) $\le m$ but, if $m\ge 1$, does not have type (resp.\ strong type or weak type) $\le m-1$.\end{definition}

\begin{example}\normalfont\label{EX27} Let $e\in\EE_n(K)$. 

\medskip
{\bf (1)} Then $e$ has (weak or strong) type $0$ iff $e\in\GA_n(K)$ by Lemma \ref{L6}(4). 

\smallskip
{\bf (2)} If via $e^{\#}$, the $K_{\t}[x_1,\ldots,x_n]$-algebra $K_{\s}[x_1,\ldots,x_n]$ is generated by $1$ element $z$, then by taking $\pi:\mathbb A^n_{K,\s}\rightarrow\mathbb A^1_K$ such that $\pi^{\#}(x_1)=z$ we get that $\pi\times e$ is a closed embedding and hence $e$ has weak type $\le 1$.
\end{example}

Referring to Definition \ref{D11}, the next example shows that in general $\pi\times e$ is not an open embedding even it is \'etale over the generic point of each $Y$. 

\begin{example}\normalfont\label{EX28.1}
Let $q\in\mathbb N^{\ast}\setminus\{1\}$ be such that $p\nmid q$. Let $X$ and $W$ be the hypersurfaces in $\mathbb A^3_k$ that are the zero loci $x^q=y^q(y+z)$ and $x^q=y+z$ (respectively). The non-regular locus $Y$ of $X$ is the zero locus $x=y=0$; so $Y\cong\mathbb A^1_k$. If $\pi:X\rightarrow\mathbb A^2_k$ is the finite surjective morphism defined by the projection on the last two coordinates, then its normalization $\pi^{\n}:X^{\n}\rightarrow\mathbb A^2_k$ is flat, \'etale above the generic point of $\pi(Y)$, has the property that $\bigl(\pi^{\n})^{-1}(\pi(Y)\bigr)$ is irreducible, and is isomorphic to the projection $W\rightarrow\mathbb A^2_k$ induced by $\pi$; so $X^{\n}\cong W\cong\mathbb A^2_k$.
\end{example}

\begin{lemma}\label{L17}
Each $e\in\EE_n(K)$ has weak type $\le 2$.
\end{lemma}

\begin{proof}
Let $z\in A_e$ be such that the approximation $Z_e=\Spec B_e$ of $X_e$ defined by the $K_{\t}[x_1,\ldots,x_n]$-subalgebra $B_e$ of $A_e$ generated by $z$ is a model of $X_e$. Let $Y$ be the smallest closed subvariety of $\mathbb A^n_{K,\t}$ such that the birational morphism $X_e\rightarrow Z_e$ is an isomorphism above $Z_e\setminus (Z_e\times_{\mathbb A^n_{K,\t}} Y)$. Let $l\in\mathbb N$ be such that $Y$ has $l$ irreducible components of dimension $n-1$. Let $O_1,\ldots,O_l$ be the local rings of $\mathbb A^n_{K,\t}$ at the generic points of these $l$ irreducible components. For $i\in\llbracket1,l\rrbracket$, the normalization $O_{i,\s}$ of $O_i$ in $K_{\s}(x_1,\ldots,x_n)$ is a free $O_i$-module. Denoting by $\widehat{O}_i$ the completion of $O_i$, the tensor product $O_{i,\s}\otimes_{O_i} \widehat{O}_i$ is a product $\widehat{O}_{i,1}\times \widehat{O}_{i,2}$ such that the factors of $\widehat{O}_{i,1}$ (resp.\ $\widehat{O}_{i,2}$) are completions of local rings of $X_e$ at points in $\Imm(\imath_e)$ (resp.\ at points in $X_e\setminus\Imm(\imath_e)$). Thus $\widehat{O}_{i,1}$ is an \'etale $\widehat{O}_i$-algebra. As the residue field $\mathcal K_i$ of $\widehat{O}_i$ is infinite, each \'etale $\mathcal K_i$-algebra is generated by $1$ element. Hence there exists $y_i\in\widehat{O}_{i,1}$ that generates the $\widehat{O}_i$-algebra $\widehat{O}_{i,1}$. From the theory of non-equivalent valuations of $\Frac(K_{\s}[x_1,\ldots,x_n])$ we get that there exists $y\in \Frac(K_{\s}[x_1,\ldots,x_n])$ such that $y-y_i$ belongs to Jacobson radical of $\widehat{O}_{i,1}$ for each $i\in\llbracket1,l\rrbracket$. As $K_{\s}[x_1,\ldots,x_n]$ is a unique factorization domain, we can assume that $y\in K_{\s}[x_1,\ldots,x_n]$. If $C_e$ is the $B_e$-subalgebra of $K_{\s}[x_1,\ldots,x_n]$ generated by $y$, it follows that the morphism $\mathbb A^n_{K,\s}\rightarrow C_e$ induces an open embedding outside a closed subvariety $W$ of $\mathbb A^n_{K,\s}$ of codimension at least $2$. As for the morphism $\pi:\mathbb A^n_{K,\s}\rightarrow\mathbb A^2_K$ defined by the rule $(x_1,\ldots,x_n)\mapsto (y,z)$, $\pi\times e:\mathbb A^n_{K,\s}\rightarrow\mathbb A^2_K\times_{\Spec K} \mathbb A^n_{K,\t}$ is a locally closed embedding outside $W$ it follows that $e$ has weak type $\le 2$.\end{proof}

\section{Additional examples}\label{S23}

The first example generalizes the \'etale morphism $\Sigma$ of Section \ref{S2}.

\begin{example}\normalfont\label{EX28}
For $(m,l)\in (\mathbb N^\ast)^2$ with $l$ prime to $p$, we consider the morphism 
$$\Sigma_{m,l}:\mathbb G_{\m,k}=\Spec(k[t,t^{-1}])\rightarrow\mathbb A^1_k$$ 
defined by the $k$-algebra homomorphism $k[x_1]\rightarrow k[t,t^{-1}]$ that maps $x_1$ to $t^{-l}-t^{pm}$; it is defined over $\mathbb F_p$. As $d(t^{-l}-t^{pm})=-lt^{-l-1}dt\in k[t,t^{-1}]^{\ast}dt$, we get that $\Sigma_{m,l}$ is \'etale. As for $\alpha\in k$, the equation $\beta^{-l}-\alpha-\beta^{pm}=0$ with $\beta\in k^{\ast}$ is equivalent to the equation $\beta^{pm+l}+\alpha \beta^{l}-1=0$ with $\beta\in k$ and as the polynomial 
$$f(t):=t^{pm+l}+\alpha t^l-1=t^l(t^{pm}+\alpha)-1$$ 
is separable because for its derivative $f'(t)=lt^{l-1}(t^{pm}+\alpha)$ we have an identity $f(t)-l^{-1}tf'(t)=-1$, it follows that $\Sigma_{m,l}$ is finite \'etale of degree $pm+l$. Also, $\Sigma_{m,l}$ is induced by the endomorphism of $\mathbb P^1_k$ defined on valued points by the rule 
$$(x_0:x_1)\mapsto (-x_0^{pm+l}+x_1^{pm+l}:x_0^lx_1^{pm}).$$ 
Therefore, for $l=1$, we denote simply $\Sigma_m:=\Sigma_{m,1}$, and we conclude that $\Sigma_m$ extends to a surjective non-finite \'etale morphism
$$\Sigma_m^+:\mathbb A^1_k\rightarrow\mathbb P^1_k$$
of degree $pm+1$ which is finite over the `finite part' open subvariety $\mathbb A^1_k$ of $\mathbb P^1_k$. We have $(\Sigma_m^+)^{-1}(1:0)=0$ and $\Sigma_1$ is the morphism $\Sigma$ of Section \ref{S2}. Moreover, $\Sigma_m^+$ extends to a finite flat endomorphism $\Sigma_m^{\proj}:\mathbb P^1_k\rightarrow\mathbb P^1_k$.
\end{example}

Next two examples study for $p=2o+1>2$ perturbations of the particular case $\e(x+x^{o+1}y^o,y-x^oy^{o+1})\in\EE_2(k)$ of the endomorphisms $e^1_{f_1,f_2;1}$ of Example \ref{EX22.9}.

\begin{example}\normalfont\label{EX29}
If $p=2o+1$ is odd, then the pair $(x+x^{o+1}y^o+x^p,y-x^oy^{o+1}-y^p)$ defines an endomorphism $e\in\EE_2(k)$ with $\pi(e)\le p$ and Jacobian matrix 
$$\begin{bmatrix} 
1+(o+1)x^oy^o & ox^{o+1}y^{o-1} \\
-ox^{o-1}y^{o+1} & 1-(o+1)x^oy^o \\ 
\end{bmatrix}$$ 
of determinant $1-[o^2-(o+1)^2]x^{2o}y^{2o}=1$.
Using the indeterminate $z:=\dfrac{y}{x}$, so $y=xz$, the system 
$$x+x^{o+1}(x^o+y^o)-x_1=0=y-y^{o+1}(x^o+y^o)-x_2$$
gives $-z^{o+1}(x-x_1)=y-x_2=zx-x_2$, hence $x=\frac{x_1z^{o+1}+x_2}{z^{o+1}+z}$, $y=\frac{x_1z^{o+1}+x_2}{z^o+1}$, and $x-x_1=\frac{-x_1z+x_2}{z^{o+1}+z}$. Thus, by substituting these expressions in the equation $x-x_1=-x^{o+1}(x^o+y^o)=-x^p(z^o+1)$, we get that 
$$\frac{-x_1z+x_2}{z^{o+1}+z}=-(z^o+1)\frac{(x_1z^{o+1}+x_2)^p}{z^p(z^o+1)^p}=-\frac{(x_1z^{o+1}+x_2)^p}{z^p(z^o+1)^{p-1}}.$$
So the $k_{\t}(x_1,x_2)$-algebra $k_{\s}(x_1,x_2)$ is isomorphic to $k_{\t}(x_1,x_2)[z]/\bigl(h(z)\bigr)$, where 
$$h(z)=h(x_1,x_2,z):=x_1^pz^{p(o+1)}+z^{p-1}(z^o+1)^{p-2}(-x_1z+x_2)+x_2^p\in k_{\t}[x_1,x_2][z].$$ Thus $z\in (A_e)_{x_1}$, $z^{-1}\in (A_e)_{x_2}$, and $\deg(e)=p(o+1)$. The group 
$$G:=\{\alpha\in k|\alpha^{p-1}=1\}\times\{\alpha\in k|\alpha^o=1\}$$ 
of order $o(p-1)=2o^2$ acts on the morphism $X_e\rightarrow\mathbb A^2_{k,\t}$ defined by $\psi_e$ via the following rule at the level of $k$-algebras: $(\alpha,\beta)\cdot (x_1,x_2,z)=(\alpha x_1,\alpha\beta x_2,\beta z)$. The group $G$ is cyclic iff $o=1$.

It is convenient to also use $x$ and $x^{-1}$ (or $y$ and $y^{-1}$) as indeterminates generating the field extension $k_{\t}(x_1,x_2)\rightarrow k_{\s}(x,y)$. The system of equations 
$$x+x^{o+1}y^o+x^p-x_1=0=y-x^oy^{o+1}-y^p-x_2$$ 
gives $y^o=\frac{x_1-x-x^p}{x^{o+1}}=\frac{x_1}{x^{o+1}}-\frac{1}{x^o}-x^o$ from the first equation. From this and the equation $y(1-x^oy^o-y^{2o})=x_2$ we get that 
$$y\Bigl(1-\frac{x_1}{x}+1+x^{2o}-\frac{x_1^2}{x^{2o+2}}-\frac{1}{x^{2o}}-x^{2o}-2+\frac{2x_1}{x}+2\frac{x_1}{x^{2o+1}}\Bigr)=x_2$$
and hence $y\bigl(\frac{x_1}{x}-\frac{1}{x^{2o}}+2\frac{x_1}{x^{2o+1}}-\frac{x_1^2}{x^{2o+2}}\bigr)=y\frac{-x_1^2+2x_1x-x^2+x_1x^p}{x^{2o+2}}=x_2$. So
$$y=\frac{x_2x^{2o+2}}{-x_1^2+2x_1x-x^2+x_1x^p}.$$ The $o$-th power of last expression equalized to the prior expression $y^o=\frac{x_1-x-x^p}{x^{o+1}}$ gives $x_2^ox^{p(o+1)}-(x_1-x-x^p)(-x_1^2+2x_1x-x^2+x_1x^p)^o=0$ in $k_{\s}(x,y)$. By defining $$g(x_1,x_2,x):=x_2^ox^{p(o+1)}-(x_1-x-x^p)(-x_1^2+2x_1x-x^2+x_1x^p)^o\in k_{\t}[x_1,x_2][x],$$
we get that the $k_{\t}(x,y)$-algebra $k_{\s}(x,y)$ is isomorphic to $k_{\t}(x_1,x_2)[x]/\bigl(g(x_1,x_2,x)\bigr)$ as well. As $2o=-1$ in $k$, we compute that $g_x(x_1,x_2,x)$ is
$$(-x_1^2+2x_1x-x^2+x_1x^p)^{o-1}[-x_1^2+2x_1x-x^2+x_1x^p-o(x_1-x-x^p)(2x_1-2x)]$$
$$=(-x_1^2+2x_1x-x^2+x_1x^p)^{o-1}[-x_1^2+2x_1x-x^2+x_1x^p+(x_1-x)^2-x_1x^p+x^{p+1}]$$
$$=x^{p+1}(-x_1^2+2x_1x-x^2+x_1x^p)^{o-1}.$$

Working with $y$ instead of $x$, a similar computation shows that for the polynomial $f(x_1,x_2,y):=x_1^oy^{p(o+1)}-(-x_2+y-y^p)(x_2^2-2x_2y+y^2+x_2y^p)^o\in k_{\t}[x_1,x_2][y]$ we have $f(x_1,x_2,y)=0$ in $k_{\s}(x,y)$ and $f_y(x_1,x_2,y)=-y^{p+1}(x_2^2-2x_2y+y^2+x_2y^p)^o$.

We show that 
$$e^{-1}(0,0)=\{(\alpha,0)|\alpha\in k, \alpha=-\alpha^p\}\cup\{(0,\beta)|\beta\in k, \beta=\beta^p\}.$$ 
It suffices to show that the assumption that the following system of equations $1+x^o(y^o+x^o)=0=1-y^o(x^o+y^o)$ has a solution in $k^2$ leads to a contradiction. As $1+x^o(y^o+x^o)=1-y^o(x^o+y^o)$, we get that $(x^o+y^o)^2=0$, so $x^o+y^o=0$; by substituting in the first equation we get that $1=0$, a contradiction. Thus $|e^{-1}(0,0)|=2p-1$.

We show that for $(\alpha,\beta)\in k^2$ with $\alpha\beta(\alpha^o+\beta^o)\neq 0$ we have $|e^{-1}(\alpha,\beta)|=p(o+1)$. The coefficient of $x^{p(o+1)}$ in $g(x_1,x_2,x)$ is $x_1^o+x_2^o$, hence $g(\alpha,\beta,x)\in k[x]$ has degree $p(o+1)$ and we have $g(\alpha,\beta,0)=(-1)^{o+1}\alpha^p\neq 0$. From these and the expressions of $g_x(\alpha,\beta,0)$ it follows that $g(\alpha,\beta,x)$ is separable. Hence 
$$e^{-1}(\alpha,\beta)=\Big\{\Big(x,\frac{\beta x^{2o+2}}{-\alpha^2+2\alpha x-x^2+\alpha x^p}\Big)\Bigl|g(\alpha,\beta,x)=0\Big\}$$
has $p(o+1)$ elements. 

Let $(\alpha,\beta)\in k^2$ be such that $\alpha^o+\beta^o=0\neq\alpha\beta$. In this case $g(\alpha,\beta,x)$ has degree $po+2$ and a similar argument shows that it is separable and $|e^{-1}(\alpha,\beta)|=po+2$.

Let $(\alpha,\beta)\in k^2\setminus\{(0,0)\}$ be such that $\alpha\beta=0$. By the similarity of $g$ and $f$ and their partial derivatives, we can assume that $\alpha\neq 0$ and $\beta=0$. We know that $e^{-1}(\alpha,0)$ is the solution set in $k^{\ast}\times k$ of the system of equations 
$$y^o-\frac{\alpha-x-x^p}{x^{o+1}}=y(-\alpha^2+2\alpha x-x^2+\alpha x^p)=0.$$ 
As $-\alpha^2+2\alpha x-x^2+\alpha x^p+\alpha(\alpha-x-x^p)=x(\alpha-x)$ and $\alpha\neq 0$, the polynomials $-\alpha^2+2\alpha x-x^2+\alpha x^p$ and $\alpha-x-x^p$ are separable of degree $p$ and do not have a common root. Therefore the fiber
$$e^{-1}(\alpha,0)=\Big\{(\delta,0)|\delta\in k,\delta+\delta^p=\alpha\Big\}\cup\Big\{\Bigl(\delta,\sqrt[o]{\frac{\alpha-\delta-\delta^p}{\delta^{o+1}}}\Bigr)\Bigl|\alpha^2-2\alpha\delta+\delta^2=\alpha\delta^p\Big\}$$
has $p+po=p(o+1)$ elements.

We conclude that $\mathcal F(e)=\{2p-1,po+2\}$ and that $N_e$ is the zero locus $x_1^o+x_2^o=0$, hence $\nu(e)=o$.

If $o=1$, i.e., $p=3$, then $g(x_1,x_2,x)=(x_1+x_2)x^6-x^5+(x_1-1)x^3+x_1^3=0$, $\nu(e)=1$, $\mathcal F(e)=\{5\}$, and $\deg(e)=6$; hence $\rho_{\et}(e)=\rho(e)=1$ and $\psi_e$ is regular, being in fact finite \'etale of degree $6$ by either Example \ref{EX5} or by remarking that the rules $(x,y)\mapsto (x-y,y)$, $(x,y)\mapsto (x+y,y)$ and $(x,y)\mapsto (y,x)$ define automorphisms $a$, $b$, and $c$ (respectively) in $\GA_2(k)$ with the property that $bea=cd_{1,1,k}c$ (see Example \ref{EX19} for $d_{1,1,k}$); thus $X_e\cong\mathbb S^2_k$.

From now on we assume that $o\ge 2$ (i.e., $p\ge 5$). Let $$\mathbb I:=\{\alpha\in k|\alpha^o=-1\}\subset\mathbb F_p.$$ 
For $\gamma\in\mathbb I$, let $\ell_{\gamma}\subset\mathbb A^2_k$ be the line $x_1-\gamma x_2=0$ and let $\widehat{O}_{\gamma}$ be the completion of the local ring $O_{\gamma}$ of $\mathbb A^2_{k,\t}$ at the generic point of $\ell_{\gamma}$. As these lines are permuted transitively by $G$, the $O_{\gamma}$s are also permuted transitively by $G$. 

Each irreducible component of $X_e\setminus\mathbb A^2_{k,\s}$ maps onto one line $\ell_{\gamma}$ with $\gamma\in\mathbb I$. So to show that $\rho(e)=o$ it suffices to show that there exists a unique irreducible component $Y_{\gamma}$ of $X_e\setminus\mathbb A^2_{k,\s}$ that maps onto $\ell_{\gamma}$ and it maps onto $\ell_{\gamma}$ isomorphically.

As $(x_1x_2)^{-1}\in \widehat{O}_{\gamma}$, $A_e\times_{k_t[x_1,x_2]} \widehat{O}_{\gamma}$ is the normalization of $\widehat{O}_{\gamma}[z]/\bigl(h(z)\bigr)$ in $\Frac(\widehat{O}_{\gamma})[z]/\bigl(h(z)\bigr)$ and we show that it is a product $\widehat{O}_{\gamma,1}\times \widehat{O}_{\gamma,2}$ with $\widehat{O}_{\gamma}\rightarrow \widehat{O}_{\gamma,1}$ finite \'etale of degree $po+2$ and $\widehat{O}_{\gamma}\rightarrow \widehat{O}_{\gamma,2}$ a finite flat totally ramified homomorphism of degree $p-2$. We compute 
$$h_z(x_1,x_2,z)=z^{p-2}(z^o+1)^{p-3}[-(z^o+1)(-x_1z+x_2)+z^o(-x_1z+x_2)-x_1z(z^o+1)]$$
$$=z^{p-2}(z^o+1)^{p-3}(-x_2-x_1z^{o+1}).$$
Therefore, after inverting $x_1x_2$, the multiple roots of $h(z)$ occur precisely when $z^o+1=x_1z^{o+1}+x_2=0$; so above the open subscheme of $\ell_{\gamma}$ defined by $x_1x_2\neq 0$, the only multiple root of $h(z)$ is $z=\gamma^{-1}$ with multiplicity $p-1$.
For the uniformizer $\pi_{\gamma}:=x_2-\gamma^{-1} x_1\in k_{\t}[x_1,x_2]$ of $\widehat{O}_{\gamma}$ and the substitution $t:=z-\gamma^{-1}$, we compute that $h(x_1,x_2,t):=h(x_1,x_2,t+\gamma^{-1})$ is equal to 
$$\pi_{\gamma}^{p}+x_1^p\Bigl[\sum_{i=1}^{o+1} \binom{o+1}{i}t^{pi}\gamma^{-p(o+1-i)}\Bigr]+(t+\gamma^{-1})^{p-1}\Bigl[\sum_{i=1}^o \binom{o}{i}t^i\gamma^{-(o-i)}\Bigr]^{p-2}(-x_1t+\pi_{\gamma}).$$
Under the additional substitution $v:=\pi_{\gamma}^{-1}t$, we can write
$$h(x_1,x_2,v):=h(x_1,x_2,\pi_{\gamma}v)=\pi_{\gamma}^{p-1}(\pi_{\gamma}+u_1v^{p-2}+u_2v^{p-1}),$$
with $(u_1,u_2)\in \bigl(\widehat{O}_{\gamma}[v]\cap \widehat{O}_{\gamma}[[v]]^{\ast}\bigr)^2$. As the normalization of $\widehat{O}_{\gamma}$ in 
$$\Frac(\widehat{O}_{\gamma})[[v]]/(\pi_{\gamma}+u_1v^{p-2}+u_2v^{p-1})$$ 
has a factor which is a finite flat totally ramified $\widehat{O}_{\gamma}$-algebra $\widehat{O}_{\gamma,2}$ of degree $p-2$ and has a factor which is a an integral \'etale $\widehat{O}_{\gamma}$-algebra of degree at least $po+2$ (recall that $|e^{-1}(\alpha,\beta)|=po+2$ if $(\alpha,\beta)\in\ell_{\gamma}(k)\setminus\{(0,0)\}$), we conclude that $A_e\times_{k[x_1,x_2]} \widehat{O}_{\gamma}$ is a product $\widehat{O}_{\gamma,1}\times \widehat{O}_{\gamma,2}$ of the type mentioned. Hence $\rho(e)=o$ and $\rho_{\et}(e)=0$. As the finite morphism $Y_{\gamma}\rightarrow\ell_{\gamma}$ is birational, it is an isomorphism.

As $\rho(e)-\rho_{\et}(e)=o$, from Corollary \ref{C2.8}(6) it follows that the singular locus of $X_e$ has at most $o$ points. If $X_e$ has a non-regular $k$-valued point $P$ that maps to a point in the set $\{(\alpha,\beta)\in (k^{\ast})^2|\alpha^o+\beta^o=0\}$ on which $G$ acts freely, then the $G$-orbit of $P$ has $o(p-1)$ points and thus $X_e$ has at least $o(p-1)$ non-regular points, a contradiction. From Corollary \ref{C2.8}(7) we get that $X_e$ is actually regular.\end{example}%HGROUP3

\begin{example}\normalfont\label{EX30}
We have the following variant of Example \ref{EX29} with smaller geometric degrees. If $p=2o+1$ is odd, then the pair
$$\bigl(x+y^o(x^{o+1}+y^{o+1}),y-x^o(x^{o+1}+y^{o+1})\bigr)$$
defines an endomorphism $e\in\EE_2(k)$ with $\pi(e)\le p$. Using the same substitution $z:=\frac{y}{x}$ we have the analogous expressions $x=\frac{x_2z^o+x_1}{z^{o+1}+1}$ and $y=\frac{x_2z^{o+1}+x_1z}{z^{o+1}+1}$ and polynomial 
$$h(z)=h(x_1,x_2,z):=x_2^pz^{po}+(-x_1z+x_2)(z^{o+1}+1)^{p-2}+x_1^p\in k_{\t}[x_1,x_2][z]$$
whose degree in $z$ is $po$, the coefficient of $z^{po}$ being $-x_1+x_2^p$; thus $\deg(e)=po$.

As we have $po+o=2o(o+1)=(p-1)(o+1)=\frac{p^2-1}{2}$ it follows that the cyclic group $G:=\{\alpha\in k|\alpha^{po+o}=1\}$ acts on the morphism $X_e\rightarrow\mathbb A^2_{k,\t}$ defined by $\psi_e$ via the following rule at the level of $k$-algebras: $\alpha\cdot (x_1,x_2,z):=(\alpha^{\frac{1}{p}}x_1,\alpha x_2,\alpha^{\frac{p-1}{p}}z)$.

Arguments similar to the ones of Example \ref{EX29} give that $e$ becomes a finite \'etale morphism after inverting $x_1^{o+1}+x_2^{o+1}$, that $|e^{-1}(0,0)|=1$, that for $(\alpha,\beta)\in (k^{\ast})^2$ with $\alpha^{o+1}+\beta^{o+1}=0$ we have $|e^{-1}(\alpha,\beta)|=po-p+2$, that $\mathcal F(e)=\{1,po-p+2\}$, and that for $o\ge 2$ we have $\rho_{\et}(e)=0$, $\rho(e)=o+1$, and $X_e$ is regular.

If $o=1$, so $p=3$, then, using the substitution $t:=x^2+y^2$ and noting that $[x+y(x^2+y^2)]^2+[y-x(x^2+y^2)]^2=t+t^3$, we get in an easier way that $\mathcal F(e)=\{1,2\}$. More precisely, for $(\alpha,\beta)\in k^2$ the solution set of the system of two polynomial equations
$$x+y(x^2+y^2)-\alpha=y-x(x^2+y^2)-\beta=0$$
is $\{(0,0)\}$ if $\alpha=\beta=0$, has $2$ elements if $\alpha^2+\beta^2=0\neq\alpha\beta$, and for $\alpha^2+\beta^2\neq 0$ is
$$\left\{\Bigl(\frac{\alpha-\beta t}{t^2+1},\frac{\beta+\alpha t}{t^2+1}\Bigr)\Big|t\in k,\; t^3+t=\alpha^2+\beta^2\right\}.$$
The morphism $X_e\rightarrow\mathbb A^2_{k,\t}$ is the projection on the first two coordinates
$$\Spec\bigl(k[x_1,x_2,t]/(t^3+t-x_1^2-x_2^2)\bigr)\rightarrow\Spec (k[x_1,x_2])$$
and thus it is a finite Galois cover isomorphic to $c_{t^2+1,t,k}$ of Section \ref{S2}. For each $i\in k$ with $i^2=-1$ there exists a unique irreducible component (affine line) of $X_e\setminus \mathbb A^2_{k,\s}$ that maps isomorphically onto the line $x_1+ix_2=0$ of $\mathbb A^2_{k,\t}$; so $\rho_{\et}(e)=\rho(e)=2$. \end{example}%HGROUP3

\begin{example}\normalfont\label{EX31}
For $\star\in\{\s,\t\}$ let $I_{\star}$ be the ideal $(x_1,x_2)$ of $k_{\star}[x_1,x_2]$. We consider $e\in\EE_2(k)$ defined by a pair $(x_1+g_1,x_2+g_2)\in k_{\s}[x_1,x_2]^2$ with $(g_1,g_2)\in I_{\s}^2$. We assume that there exist $\alpha\in k$ and $h(t)\in t^2k[t]$ such that $t+g_1\bigl(t,h(t)\bigr)$ and $h(t)+g_2\bigl(t,h(t)\bigr)$ have $\alpha$ as a common root with the same multiplicity $m\in\mathbb N^{\ast}$; thus $\alpha\neq 0$. If $C$ is the irreducible curve of $\mathbb A^2_{k,\s}$ which is the zero locus $x_2-h(x_1)=0$, then the Zariski closure $\overline{e(C)}$ of $e(C)$ in $\mathbb A^2_{k,\t}$ is the zero locus of an irreducible polynomial $f(x_1,x_2)\in I_{\t}$ such that $f\bigl(t+g_1(t,h(t)),h(t)+g_2(t,h(t))\bigr)=0$ for all $t\in k$. Let $(\beta_1,\beta_2)\in k^2$ be such that $f(x_1,x_2)-\beta_1x_1-\beta_2x_2\in I_{\t}^2$. We have $\beta_1=0$ as $\beta_1$ is the coefficient of $t$ in $f\bigl(t+g_1(t,h(t)),h(t)+g_2(t,h(t))\bigr)=0$. We have $\beta_2=0$ as there exists $\gamma\in k^{\ast}$ such that 
$$f\bigl(t+g_1(t,h(t)),h(t)+g_2(t,h(t))\bigr)-\gamma\beta_2(t-\alpha)^m\in k[t](t-\alpha)^{m+1}.$$ 
Therefore $(0,0)\in\Sing(\overline{e(C)})$. Moreover, if $d\in\EE_2(k)$ is such that $C\in\Irr(N_d)$ with $N_d\subset\mathbb A^2_{k,\t}=\mathbb A^2_{k,\s}$, then $\overline{e(C)}\in\Irr(N_{ed})$ is a singular irreducible component.

As a concrete example, for $p=2o+1$ odd we take $g_1(x_1,x_2)=x_1^{o+1}x_2^o-2x_1^{p(o+1)}$, $g_2(x_1,x_2)=-x_1^ox_2^{o+1}$, $h(t)=t^{p+1}$, and $\alpha\in\{\delta\in k|\delta^{o(p+2)}=1\}$. Then we compute $t+g_1\bigl(t,h(t)\bigr)=t\bigl(1-t^{o(p+2)}\bigr)$ and $h(t)+g_2\bigl(t,h(t)\bigr)=t^{p+1}\bigl(1-t^{o(p+2)}\bigr)$; so $m=1$ is the same multiplicity for $\alpha$.\end{example}

\begin{proposition}\label{PR9}
Let $(m, q)\in (\mathbb N^\ast)^2$. Let $\mathbb U\subset\{1,p,2p,3p,\ldots,(m-1)p\}$ be a subset with $0<|\mathbb U|\le p^q$. Then there exists $e\in\BE_2(k)$ with $\deg(e)=pm$ and $\mathcal F(e)=\mathbb U$ (so $e$ is surjective) which is defined over $\mathbb F_{p^q}$.
\end{proposition}

\begin{proof}
We write $\mathbb U=\{r_1,\ldots,r_j\}$ with $r_1<\cdots<r_j$; so $j\in \llbracket1,m\rrbracket$. As $p^q\ge j$, we consider $j$ distinct elements $\alpha_1,\ldots,\alpha_j$ of $\mathbb F_{p^q}$. Let $f_1(t),\ldots,f_m(t)\in \mathbb F_{p^q}[t]\setminus\{0\}$ be polynomials subject to the following condition.

\medskip
{\bf ($\sharp$)} If $\mathbb I_i:=\{\alpha\in k|f_i(\alpha)=0\}$ for $i\in \llbracket1,m\rrbracket$, then we have a chain of inclusions
$$\mathbb I_1\subset \mathbb I_2\subset\cdots\subset \mathbb I_m=\{\alpha_i|i\in\llbracket1,j\rrbracket\}$$ 
 defined recursively as follows. We define $\mathbb I_1:=\emptyset$ if $1\notin\mathbb U$ and $\mathbb I_1:=\{\alpha_1\}$ if $1\in\mathbb U$. For $s\in\llbracket1,m-1\rrbracket$, $\mathbb I_{s+1}\setminus \mathbb I_s$ is empty if $ps\notin\mathbb U$ and is $\{\alpha_o\}$ if $ps=r_o\in\mathbb U$.

\medskip
The elementary $p$-morphism $e\in\pMor_2(k)$ defined by the rule
$$(x,y)\mapsto \bigl(x,y+\sum_{i=1}^m f_i(x)y^{pi}\bigr)$$
is defined over $\mathbb F_{p^q}$. Condition ($\sharp$) implies that for each $(\alpha,\beta)\in k^2$ the fiber $e^{-1}(\alpha,\beta)$ has $pm$ points if $f_m(\alpha)\neq 0$ and has $r_o$ points if $\alpha=\alpha_o$, $r_o$ being the degree of the polynomial $y+\sum_{i=1}^m f_i(\alpha_l)y^{pi}\in\mathbb F_{p^q}[y]$. Thus we have $\mathcal F(e)=\mathbb U$. If we chose $f_1(t),\ldots,f_m(t)\in \mathbb F_{p^q}[t^p]\setminus\{0\}$, then $e\in\BE_2(k)$.
\end{proof}

\begin{example}\normalfont\label{EX32} As in this example we use composites of endomorphisms, we use only $k$ without the lower indices ${}_{\s}$ and  ${}_{\t}$. Generalizing the endomorphisms in $\EE_2(k)$ that are defined by the rules $(x,y)\mapsto \bigl(x+g(x)^p,y+y^pf(x)^p\bigr)$ (see Example \ref{EX6}) and $(x,y)\mapsto \bigl(x,y+\sum_{i=1}^m f_i(x)y^{pi}\bigr)$ (see Proposition \ref{PR9}) or that are as in Example \ref{EX24} for $(n,q_1,f_1)=(2,1,1)$, we consider $e\in\EE_2(k)$ be defined by a pair 
$$(g_1,g_2)\in k[x_1]\times k[x_1,x_2].$$ We have an identity $e=d_1d_2$, with $d_2:=\e\bigl(x_1,g_2(x_1,x_2)\bigr)\in\EE_2(k)$ surjective and $d_1:=\e\bigl(g_1(x_1),x_2\bigr)\in\EE_2(k)$ finite, and there exist pairs $(\alpha_1,\alpha_2)\in (k^{\ast})^2$ and $(g,h)\in k[x_1]\times k[x_1,x_2]$ such that $(g_1,g_2)=\bigl(\alpha_1 x_1+g(x_1^p),-\alpha_2 x_2-h(x_1,x_2^p)\bigr)$. Both $d_1$ and $d_2$ have strong type $\le 1$ in the sense of Definition \ref{D11} and moreover $(e,d_1,d_2)\in\pMor_2(k)^3$. As $\deg(d_1)=\max\bigl(1,p\deg(g)\bigr)$ and $\deg(d_2)=\max(1,pm)$, where $m$ is the degree of $h$ in $x_2$, Conjecture \ref{CJ3} holds for $e$. As $A_e$ is the normalization of $k[x_1,x_2]$ in $\Frac(A_e)=\Frac\bigl(k[x_1,x_2][y]/(\alpha_2y+h(x_1,y^p)+x_2)\bigr)$, to study it, using the substitution $x_2^{\prime}:=x_2+h(x_1,0)$, we can assume that $h(x_1,0)=0$; so $m\neq 0$. If $m=-\infty$, then $d_2$ is an isomorphism; so we can assume $m\ge 1$. We write 
$$h(x_1,x_2)=\sum_{i=1}^m f_i(x_1)x_2^i$$ 
with $(f_1,\ldots,f_m)\in k[x_1]^m$, $f_m\neq 0$. Substituting $z:=f_m(x_1)y$, $A_e$ is the normalization of the $k[x_1,x_2]$-algebra $B_e:=k[x_1,x_2][z]/\bigl(h(z)\bigr)$, where
$$h(z):=z^{pm}+\sum_{i=1}^{m-1}f_i(x_1)f_m(x_1)^{pm-1-pi}z^{pi}+\alpha_2f_m(x_1)^{pm-2}z+f_m(x_1)^{pm-1}x_2.$$
Multiplying $h(z)$ by $f_m(x_1)x_2^{pm-1}$ and substituting $w:=f_m(x_1)x_2z^{-1}$, $A_e$ is the normalization of the $k[x_1,x_2]$-algebra
$$B'_e:=k[x_1,x_2][w]/\bigl(w^{pm}+\alpha_2w^{pm-1}+\sum_{i=1}^{m} f_i(x_1)x_2^{pi-1}w^{p(m-i)}\bigr).$$

If $m$ is a power of $p$ and $h(x_1,x_2)$ is additive in $x_2$ (i.e., and $f_i(x_1)=0$ if $i\in \llbracket1,m\rrbracket$ is not a power of $p$), then $\mathbb G_{\a,k}$ is a subgroup scheme of the group ind-scheme of automorphisms of the morphism $X_e\rightarrow\mathbb A^2_{k,\t}$ defined by $\psi_e$: at the level of $k$-algebras, they are parametrized by $\beta\in k$ via the rule 
$$\beta\cdot(x_1,x_2,y)=\bigl(x_1,x_2-\alpha_2\beta-\sum_{i=1}^m f_i(x_1)\beta^{pi},y+\beta\bigr),$$
where we identify $y$ with $y+(\alpha_2y+h(x_1,y^p)+x_2)\in\Frac(A_e)$.

Assume that $m=1$. For $p=2$ we have $A_e\cong B_e\cong B'_e$ and $e$ and $d_2$ are \'etale with $0<\rho_{\et}(e)=\rho(e)=\rho_{\et}(d_2)=\rho(d_2)=\z(f_m)$ while for $p>2$ we have $\rho(e)\ge \z(f_m)$ and $A_e$ is regular by Lemma \ref{L3}(4).

If $m=2$, then $A_e$ is regular by Lemma \ref{L4}.

If $m\ge 3$, then $A_e$ is regular in many cases including the following ones: (i) $f_m$ is a separable polynomial by Lemma \ref{L5}; (ii) $m$ is a power of $p$ and $h(x_1,x_2)$ is additive in $x_2$ by Lemma \ref{L4} or as a direct consequence of the existence of the $\mathbb G_{\a,k}$ group scheme of automorphisms; (iii) Corollary \ref{C2.8}(7) or Lemma \ref{L3} applies.

As $e$ is finite \'etale above the complement of the union of lines of $\mathbb A^2_{k,\t}$ given by equation $x_1=\gamma$ with $\gamma\in k$ such that $f_m(\gamma)=0$, we have $\nu(e)=\z(f_m)$.\end{example}

\begin{example}\normalfont\label{EX33}
We refer to Example \ref{EX32} with $p=2$, $m=2$, and $h(x_1,x_2)=x_1x_2^2$; so $f_1(x_1)=0$, $f_2(x_1)=x_1$, and $h(z)=z^4+x_1^2z+x_1^3x_2$. Denoting also by $z$ the element $z+(h(z))\in k_{\t}[x_1,x_2][z]/\bigl(h(z)\bigr)\subset A_e$ and introducing new indeterminates and elements $w:=\frac{z^2}{x_1}$, $v:=\frac{w^3}{x_1}+1$, and $t:=v+x_2w$, we get that $A_e$ is also the normalization of either one of the following three $k_{\t}[x_1,x_2]$-algebras $k_{\t}[x_1,x_2][w]/(w^4+wx_1+x_1^2x_2^2)$, $k_{\t}[x_1,x_2][v]/(v^4+v^3+x_1^2x_2^6)$, and 
$$k_{\t}[x_1,x_2][t]/(t^4+t^3+x_1x_2^3),$$ but with only one generator we cannot get better than this (cf.\ Corollary \ref{C11} below). We have $(z,w,v,t)\in A_e^4$ satisfying $z=w^2+x_1x_2$ and $v=t^2$. So we compute $w^3+x_1x_2w+x_1t+x_1=t^2+t+x_2w=wt+w^2x_2+x_1x_2^2=0$. Thus we have a natural $k_{\t}[x_1,x_2]$-algebra homomorphism 
$$\tilde\imath_e:k_{\t}[x_1,x_2][w,t]/(w^3+x_1x_2w+x_1t+x_1,t^2+t+x_2w,wt+w^2x_2+x_1x_2^2)\rightarrow A_e.$$ 
The domain $B_e$ of $\tilde\imath_e$ is a $k_{\t}[x_1,x_2]$-module generated by $\{1,w,w^2,t\}$ and by tensoring with $k_{\t}(x_1,x_2)$ it follows that it is free of rank $4$. By computing the ranks of the $3\times 4$ Jacobian matrices at points of $\Spec B_e$ we get that $B_e$ is a smooth $k$-algebra. It follows that $\tilde\imath_e$ is an isomorphism and we reobtain the fact that $A_e$ is regular. Either from this last isomorphism or from the proof of Lemma \ref{L4} applied to $p=m=2$ we get that $\rho_{\et}(e)=0$ and $\rho(e)=1$. Hence $Y:=X_e\setminus X_e^{\et}=X_e\setminus\Imm(\imath_e)$ is irreducible and its multiplicity in $\psi_e^{-1}\bigl(\psi_e(Y)\bigr)$ is $3$ and (see Corollary \ref{C2.8}(1) and (2) and its notation) $A_Y=A_e$ does not have a $2$-basis.
\end{example}

\begin{example}\normalfont\label{EX34}
Let $\bigl((f_1,g_1),(f_2,g_2)\bigr)\in\Theta_K^2$. We have a complete intersection
$$\mathbb S^2_{f_1,g_1;f_2,g_2,K}:=\Spec \bigl(K[w,x,y,z]/(xf_1(x)+yg_1(y)z,wf_2(w)+xg_2(y))\bigr)$$
in $\mathbb A^4_K$. The projection $\mathbb A^4_K\rightarrow\mathbb A^3_K$ on the last three coordinates induces a finite flat morphism $\phi:\mathbb S^2_{f_1,g_1;f_2,g_2,K}\rightarrow\mathbb S^2_{f_1,g_1,K}$. Let $g:=g_1g_2\in K[t]$. The cartesian diagram
\[\xymatrixcolsep{6pc}\xymatrix{
\mathbb S^2_{f_2,g,K}\ar[r]^{\jmath_{f_1,g_1;f_2,g_2,K}} \ar[d]^{\pi} & \mathbb S^2_{f_1,g_1;f_2,g_2,K} \ar[d]^{\phi} \\
\mathbb A^2_K \ar[r]^{\imath_{f_1,g_1,K}} & \mathbb S^2_{f_1,g_1,K}
}\]
is defined by the identification $\mathbb S^2_{f_2,g,K}=\Spec \bigl(K[w,x,y]/(wf_2(w)+xyg(y))\bigr)$ with $\pi$ induced by the projection on the last two variables and $$\jmath_{f_1,g_1;f_2,g_2,K}^{\#}\bigl(w+(xf_1(x)+yg_1(y)z,wf_2(w)+xg_2(y))\bigr)=w+\bigl(wf_2(w)+xyg(y))\bigr),$$ and gives an open embedding 
$$\imath_{f_1,g_1;f_2,g_2,K}:=\jmath_{f_1,g_1;f_2,g_2,K}\circ\imath_{f_2,g}:\mathbb A^2_K\rightarrow\mathbb S^2_{f_1,g_1;f_2,g_2,K}$$ 
defined by the rule $(x,y)\mapsto \bigl(xyg(y),Fyg_1(y),y,-Ff_1(Fyg_1(y)\bigr)$ for the polynomial $F(x,y):=-xf_2\bigl(xyg(y)\bigr)\in K[x,y]$. The flatness of $\phi$ and the integrality of $\mathbb S^2_{f_2,g,K}$ imply that $\mathbb S^2_{f_1,g_1;f_2,g_2,K}$ is integral.

Note that if for $i\in\{1,2\}$ we have $\deg(f_i)\ge 1$ and there exists $m_i\in \mathbb N^{\ast}$ such that $g_i(t)=t^{m_i}$, then the resulting surfaces $\mathbb S^2_{f_1,g_1;f_2,g_2,K}$ form a subclass of the class of double Danielewski surfaces introduced and studied in \cite{GS}, Sect.\ 3 (they are integral by \cite{GS}, Lem.\ 3.3 and \cite{GS}, Thm.\ 3.11 lists their isomorphism classes).

Assume that $\chr(K)=p$. Let $m\in\mathbb N^{\ast}$. We assume that $(f_2,g_2)\in\Theta^1_{k,m}$. Thus $\phi$ is a finite \'etale cover of degree $pm$. Let $\Phi:=\Phi_{f_1,g_1,h_1,q_1}:\mathbb S^2_{f_1,g_1,K}\rightarrow\mathbb A^2_K$ be a finite \'etale cover as in Proposition \ref{PR8}(1). The composite 
$$e:=\Phi\circ\phi\circ\imath_{f_1,g_1;f_2,g_2,K}\in\EE_2(K)$$ 
is such that $X_e=\mathbb S^2_{f_1,g_1;f_2,g_2,K}$ and $p^2$ divides $\deg(e)=pm\deg(\Phi)$; we have $\edim(X_e)\le 4$ and it is expected that the equality holds in many cases with $\deg(f_1)\deg(f_2)\deg(g_2)>0$.\end{example}

\begin{example}\normalfont\label{EX35}
Let $s\in\mathbb N^{\ast}\setminus\{1\}$. Let $\bigl((f_1,g_1),\ldots,(f_s,g_s)\bigr)\in\Theta_K^s$. Let
$$\mathbb S^{s+1}_{(f_1,g_1;\ldots;f_s,g_s),K}:=\Spec \bigl(K[x_1,\ldots,x_s,y,z_1,\ldots,z_s]/(x_if_i(x_i)+yg_i(y)z_i|i\in \llbracket1,s\rrbracket)\bigr);$$
it is a complete intersection in $\mathbb A^{2s+1}_K$ of dimension $s+1$ equipped with a finite flat morphism $\mathbb S^{s+1}_{(f_1,g_1;\ldots;f_s,g_s),K}\rightarrow\mathbb A^{s+1}_K$ defined by the projection $\mathbb A^{2s+1}_K\rightarrow\mathbb A^{s+1}_K$ on the last $s+1$ coordinates. One can use the open embeddings $\imath_{f_i,g_i,K}$ and, when $\deg (f_i)\ge 1$ and $g_i(0)\neq 0$, $\jmath_{f_i,g_i,K;\delta_i}$ with $\delta_i\in K^{\ast}$ such that $f_i\bigl(\delta_ig_i(0)\bigr)=0$ to get open embeddings $\mathbb A^{s+1}_K\rightarrow\mathbb S^{s+1}_{(f_1,g_1;\ldots;f_s,g_s),K}$. E.g., using only $\imath_{f_i,g_i,K}$ for each $i\in \llbracket1,s\rrbracket$, we have an open embedding
$$\imath_{(f_1,g_1;\ldots;f_s,g_s),K}:\mathbb A^{s+1}_K\rightarrow\mathbb S^{s+1}_{(f_1,g_1;\ldots;f_s,g_s),K}$$ 
defined by the rule 
$$(x_1,\ldots,x_s,y)\mapsto \bigl(x_1yg_1(y),\ldots,x_syg_s(y),y,-x_1f_1(x_1yg_1(y)),\ldots,-x_sf_s(x_1yg_s(y))\bigr).$$ 

If $\chr(K)=p$ and $n_i:=\deg(f_i)-1\ge 1$ for each $i\in \llbracket1,s\rrbracket$, then the proof of Proposition \ref{PR8}(1) applies to show that $\mathbb S^{s+1}_{(f_1,g_1;\ldots;f_s,g_s),K}$ is a finite \'etale cover of $\mathbb A^{s+1}_n$ of degree $p^s\prod_{i=1}^s q_in_i$, where for each $i\in \llbracket1,s\rrbracket$ the integer $q_i\in\mathbb N^{\ast}$ is associated to the pair $(f_i,g_i)$ as in Proposition \ref{PR8}. 
\end{example}

\begin{proposition}\label{PR10}
For each $l\in\llbracket1,\lfloor\frac{pm-1}{2}\rfloor\rrbracket$ there exists $e\in\EE_2(k)$ with $\deg(e)=pm$, $\varphi(e)=0$, and $\mathcal F(e)=\{1,l+1,pm-l\}$. Moreover, we can chose $e$ to be defined over each finite field $\mathbb F_{p^q}$ that contains $l$ zeros of the equation $z^{pm-1}-1=0$ (e.g., we can take $q\in\mathbb N^\ast$ such that $pm-1$ divides $p^q-1$), or over each $\mathbb F_{p^q}$ with $q\in\mathbb N^\ast$ such that $p^q\ge l+1$ if $l\in\llbracket1,m-1\rrbracket$, or over $\mathbb F_p$ if $l$ divides $pm-1$. 
\end{proposition}

\begin{proof} 
We apply Example \ref{EX23} for $n=2$ as follows. If $g\in k[t]$ is such that $\deg(g)=l$, then $\deg(h)=pm-l-1$, $\deg(e_{g,h})=pm$, and $\mathcal F(e_{g,h})=\{1,l+1,pm-l\}$. If $tg(t)h(t)=t-t^{pm}$, then we can choose $g$ and $h$ and hence also $e_{g,h}$ to be defined over each finite field $\mathbb F_{p^q}$ which contains $l$ zeros of the equation $z^{pm-1}-1=0$ (e.g., we can take $q\in\mathbb N^\ast$ such that $pm-1$ divides $p^q-1$). 

If $l$ divides $pm-1$, then we can choose $g(t)=1-t^l$ and $h(t)=\frac{1-t^{pm-1}}{1-t^l}$ in $\mathbb F_p[t]$.

If $l\in\llbracket1,m-1\rrbracket$, let $\beta_1,\ldots,\beta_l\in k^{\ast}$ be distinct and we take 
$$g(t):=\prod_{i=1}^l (1+\beta_it)\in k[t];$$ 
clearly $\deg(g)=l$. As the determinant of the generalized Vandermonde $l\times l$ matrix whose columns indexed by $s\in\llbracket1,l\rrbracket$ are $(\beta_s^{-p},\beta_s^{-2p},\ldots,\beta_s^{-pl})$ is in $k^{\ast}$ (being $\prod_{i=1}^l \beta_i^{-p}$ times a non-zero Vandermonde determinant), there exists a unique $l$-tuple $(\alpha_1,\ldots,\alpha_l)\in k^l$ such that $tg(t)$ divides $F(t)=t+(\sum_{i=1}^{l} \alpha_it^{pi})+\sum_{i=l+1}^m t^{pi}$. If $p^q\ge l+1$, then by taking $\beta_1,\ldots,\beta_l\in\mathbb F_{p^q}^{\ast}$ we get first that $g(t)\in \mathbb F_{p^q}[t]$ and second that $h(t)=\frac{F(t)}{tg(t)}\in\mathbb F_{p^q}[t]$.
\end{proof}

\section{A surjectivity property}\label{S24}

\begin{proposition}\label{PR11}
Let $\mathcal S$ be an $F_2$-system over $k$ such that $\pi(\mathcal S)\le 2p$. Then $e_{\mathcal S}$ is surjective, i.e., $\varphi(e_{\mathcal S})=0$.%, and $\deg(e_{\mathcal S})\in\{1,p,2p,p^2,p^2+p,2p^2,4p^2\}$.
\end{proposition}

\begin{proof}
We can assume that $e:=e_{\mathcal S}\in\EE_2(k)$ is defined by a pair 
$$(\ell_1+\ell_3^p+g_1^p,\ell_2+\ell_4^p+g_2^p)\in k_s[x_1,x_2]^2,$$ 
where $\ell_1$ to $\ell_4$ are linear forms and $g_1$ and $g_2$ are quadratic forms. As $e$ is \'etale, $\ell_1$ and $\ell_2$ are linearly independent over $k$. If the system of homogeneous equations $g_1=0=g_2$ has no solution in $\mathbb P^1_k$, then $e$ is finite surjective and $\deg(e)=4p^2$ by Lemma \ref{F3}(2). Hence we can assume that there exist linear forms $\ell_0$, $\ell_5$, and $\ell_6$ in $k_{\s}[x_1,x_2]$ such that $\ell_0\neq 0$, $g_1=\ell_0\ell_5$, and $g_2=\ell_0\ell_6$. 

Let $LCC_{\s}$ (resp.\ $LCC_{\t}$) denote the operation of a linear change of coordinates in the source (target) space. 

We consider two disjoint cases as follows, and we also compute $\deg(e)$ in all subcases except one.

{\bf Case 1: $\ell_5$ and $\ell_6$ are linearly dependent.} By performing $LLC_{\t}$ we can arrange that $\ell_5=0$; so $g_1=0$ and in this case the roles of $\ell_0$ and $\ell_6$ are symmetric. 

If $\ell_3=0$, then by performing $LLC_{\s}$ we can assume that $\ell_1=x_1$ and $\ell_2=x_2$; so for each $\alpha\in k$, the equation $x_2+\ell_4^p(\alpha,x_2)^p+g_2(\alpha,x_2)^p=\beta$ in $x_2$ has solutions in $k$ and $\deg(e)$ is the degree in $x_2$ of $x_2+\ell_4^p(\alpha,x_2)^p+g_2(\alpha,x_2)^p$, and therefore $\deg(e)\in\{1,p,2p\}$. Thus we can assume that $\ell_3\neq 0$. 

If $\ell_6\neq 0$ and $\ell_3\notin k^{\ast}\ell_0\cup k^{\ast}\ell_6$, then $e$ is finite surjective and $\deg(e)=2p^2$ by Lemma \ref{F3}(2). Thus, by interchanging $\ell_0$ and $\ell_6$ if needed, we can assume that $\ell_3=\ell_0$. 

By performing $LLC_{\s}$ we can also assume that $\ell_1=x_1$ and $\ell_3\in\{x_1,-x_2\}$. If $\ell_3=x_1$, then the system of equations $x_1+x_1^p-\alpha=0=\ell_2+\ell_4^p+x_1^p\ell_6^p-\beta=0$ has solutions in $k^2$, so $e$ is surjective, and $\deg(e)$ is $p$ times the degree in $x_2$ of $\ell_2+\ell_4^p+x_1^p\ell_6^p$; so $\deg(e)\in\{p,p^2,p^3\}$. If $\ell_3=x_2$, then the system of equations $x_1-x_2^p-\alpha=0=\ell_2+\ell_4^p-x_2^p\ell_6^p-\beta=0$ over $k$ is equivalent to the equation $F_2(x_2)=\beta$ with $F_2:=\ell_3(x_2^p+\alpha,x_2)+\ell_4(x_2^p+\alpha,x_2)^p-x_2^p\ell_6(x_2^p+\alpha,x_2)^p\in k_{\s}[x_2]$, hence $e$ is finite with $\deg(e)=\deg(F_2)\in\{1,p,2p,p^2,p^2+p\}$.

{\bf Case 2: $\ell_5$ and $\ell_6$ are linearly independent.} Let $\pi_1:\mathbb A^2_{k,\t}\rightarrow\mathbb A^1_{k,\t}$ be the projection on the first coordinate. We have $\varphi(e_{\mathcal S})=0$ iff for each $\alpha\in\mathbb A^1_{k,\t}(k)$, the \'etale morphism 
$$e_{\alpha}:(\pi_1e)^{-1}(\alpha)\rightarrow \pi_1^{-1}(\alpha)=\{\alpha\}\times_{\Spec k} \mathbb A^1_{k,\t}\cong\mathbb A^1_k$$
is surjective. By performing $LLC_{\s}$ we can assume that $\ell_0=x_1$, $g_1=x_1^2$, and $g_2=x_1^px_2^p$. The plane curve $(\pi_1e)^{-1}(\alpha)$ is the zero locus $\ell_1+\ell_3^p+x_1^{2p}-\alpha=0$. Let $(\gamma_1,\gamma_2,\delta_1,\delta_2)\in k^4$ be such that $\ell_1=\gamma_1x_1+\gamma_2x_2$ and $\ell_3=\delta_1x_1+\delta_2x_2$. 

If $\gamma_2=\delta_2=0$, then $\gamma_1\neq 0$ and $(\pi_1e)^{-1}(\alpha)$ is isomorphic to $2p$ copies of $\mathbb A^1_k$ indexed by the roots of the equation $x_1^{2p}+\delta_1x_1^p+\gamma_1x_1-\alpha=0$ and for each such $\mathbb A^1_k$ copy the resulting \'etale morphism $\mathbb A^1_k\rightarrow \pi_1^{-1}(\alpha)$ is finite surjective by Lemma \ref{L10}(3) of degree $p$ for each $\alpha\in k$ outside a finite set; thus $\deg(e)=2p^2$. 

If $\gamma_2\delta_2=0\neq\gamma_2+\delta_2$, then $(\pi_1e)^{-1}(\alpha)\cong \mathbb A^1_k$. The resulting \'etale morphism $(\pi_1e)^{-1}(\alpha)\rightarrow \pi_1^{-1}(\alpha)$ is finite surjective by Lemma \ref{L10}(3) and by solving for $x_2$ in terms on $x_1$ using the equation $\ell_1+\ell_3^p+x_1^{2p}-\alpha=0$ we get that $\deg(e)=2p^2+p$ if $\delta_2=0$ and $\deg(e)=3p^2$ if $\gamma_2=0$.

If $\gamma_2\delta_2\neq 0$, then $(\pi_1e)^{-1}(\alpha)$ is either isomorphic to $p$ copies of $\mathbb A^1_k$ or is a connected finite Galois cover of $\mathbb A^1_k$ of Galois group $\mathbb Z/p\mathbb Z$ and thus, as $\mathbb P^1_k$ is simply connected, it has $1$ point at infinity; as $(\pi_1e)^{-1}(\alpha)$ is a disjoint union of smooth affine curves with $1$ point at infinity, $e_{\alpha}$ is surjective by Lemma \ref{L10}(3).\end{proof}

\section{\'Etale endomorphisms of $\mathbb A^2_k$ via line bundles}\label{S25}

We have the following moduli interpretation of $\mathbb S^2$. With $B^+$ and $B^-$ as in Section \ref{S19}, the intersection
$$T:=B^+\cap B^-$$
is the maximal torus of $\SL_2$ formed by diagonal $2\times 2$ matrices of determinant $1$. The centralizer of $T$ in $\GL_2$ is the maximal torus $GT$ of $\GL_2$ formed by invertible diagonal $2\times 2$ matrices. The moduli scheme $\mathbb M$ over $\Spec\mathbb Z$ that parametrizes $2\times 2$ matrices that are projectors of rank $1$, equivalently, have trace $1$ and determinant $0$ and hence can be put in the form $\begin{bmatrix} 
x & y \\
z & 1-x \\ 
\end{bmatrix}$ with $x^2-x+yz=0$, is naturally identified with $\mathbb S^2$. The action by conjugation 
$$\GL_2\times_{\Spec(\mathbb Z)} \mathbb M\rightarrow\mathbb M$$ 
is transitive. As the stabilizer (centralizer) of $\begin{bmatrix} 
1 & 0 \\
0 & 0 \\ 
\end{bmatrix}$ is $GT$, we identify
$$\SL_2/T=\GL_2/GT\cong\mathbb M=\mathbb S^2.$$
If $wz-xy=1$, then we compute
$$\begin{bmatrix} 
w & x \\
y & z \\ 
\end{bmatrix}\begin{bmatrix} 
1 & 0 \\
0 & 0 \\ 
\end{bmatrix}\begin{bmatrix} 
w & x \\
y & z \\ 
\end{bmatrix}^{-1}=\begin{bmatrix} 
w & x \\
y & z \\ 
\end{bmatrix}\begin{bmatrix} 
1 & 0 \\
0 & 0 \\ 
\end{bmatrix}\begin{bmatrix} 
z & -x \\
-y & w \\ 
\end{bmatrix}=\begin{bmatrix} 
wz & -wx \\
yz & 1-wz \\ 
\end{bmatrix}.$$ Thus the composite morphism $\SL_2\rightarrow\SL_2/T\rightarrow\mathbb S^2$ is given by the following rule on valued points:
$\begin{bmatrix} 
w & x \\
y & z \\ 
\end{bmatrix}\mapsto (wz,-wx,yz).$

\phantomsection{Let $\mathcal Q:\mathbb S^2=\SL_2/T\rightarrow \SL_2/B^+\cong\mathbb P^1$ be the natural quotient morphism that defines a line bundle over $\mathbb P^1$ and that allows us to identify $\mathbb D^1=\mathcal Q^{-1}(1:0)$. So, on valued points, $\mathcal Q$ is defined by the rule: $(x,y,z)$ maps to either $(-x:y)$ or $(z:x-1)$ (whichever one is defined; if both are defined, then they are equal).}\label{EXT10}

For $m\in\mathbb N^{\ast}$ we consider the cartesian diagram
\[\xymatrix{
\mathbb A^2_k \ar[r]^{\nat} \ar[d]^{\Sigma_{m,2}} & \mathbb A^1_k \ar[d]^{\Sigma_m^+} \\
\mathbb S^2_k=\SL_{2,k}/T_k \ar[r]^{\mathcal Q_k} & \SL_{2,k}/B^+_k\cong\mathbb P^1_k,
}\]
where $\Sigma_m^+$ is as in Example \ref{EX28} and we are using the fact that each line bundle over $\mathbb A^1_k$ is trivial; so $\nat:\mathbb A^2_k\rightarrow\mathbb A^1_k$ is the projection on the first coordinate and the diagram is defined over $\mathbb F_p$. Thus $\Sigma_{m,2}$ is surjective \'etale of degree $pm+1$ as so is $\Sigma_m^+$. We get that
$$\mathbb A^1_k\cong\nat^{-1}(0)=\nat^{-1}\bigl((\Sigma_m^+)^{-1}(1:0)\bigr)=\Sigma_{m,2}^{-1}\bigl(\mathcal Q_k^{-1}(1:0)\bigr)=\Sigma_{m,2}^{-1}(\mathbb D^1_k)$$ maps via $\Sigma_{m,2}$ isomorphically onto $\mathbb D^1_k=\mathcal Q_k^{-1}(1:0)$.

\begin{remark}\normalfont\label{R20}
Using the definition of $\imath_{t-1,1,\mathbb Z}$ (see Displays (\ref{EQ0a}) and (\ref{EQ0b})), the morphism 
$$\Sigma_{m,2}:\mathbb A^2_k\rightarrow\mathbb S^2_k$$ is, in suitable coordinates, given by the following rule on valued points
$$(x,y)\mapsto\bigl((x^{pm+1}+1)(-xy+1),x-x^2y,(x^{2pm}+2x^{pm-1})(x^2y-x)+x^{pm}+y\bigr).$$
We have $\Sigma_{m,2}(\{0\}\times_{\Spec k}\mathbb A^1_k)=\mathbb D^1_k$.
\end{remark}

\begin{proposition}\label{PR12}
Let $m\in\mathbb N$ and $(f,g)\in\Theta_k$ with $\deg(f)\ge 1$. Let $(\alpha,\beta)\in k^{\ast}\times k$ be such that $f(\alpha)=\beta g(\beta)=0\neq g(\beta)+\beta g'(\beta)$. Let 
$$Y=Y_{\alpha,\beta}\in\Irr\bigl(\mathbb S^2_{f,g,k}\setminus\Imm(\imath_{f,g,k})\bigr)$$ 
be the zero locus $x-\alpha=y-\beta=0$ (see the line after Display (\ref{EQ0c})). For the open subvariety $\mathbb S^2_{f,g,k,Y}:=\Imm(\imath_{f,g,k})\cup Y$ of $\mathbb S^2_{f,g,k}$ there exists a non-finite surjective \'etale morphism 
$$\Sigma_{f,g,m,Y}:\mathbb A^2_k\rightarrow\mathbb S^2_{f,g,k,Y}$$ 
of degree $pm+1$ with the following properties.

\medskip
{\bf (1)} It is defined over each subfield of $k$ that contains all roots of $fg$.

\smallskip
{\bf (2)} It is finite over $\Imm(\imath_{f,g,k})$ and $\Sigma_{f,g,m,Y}^{-1}(Y)\cong \mathbb A^1_k$. 

\smallskip
{\bf (3)} It is the composite of an open embedding $\imath_{f,g,m,Y}:\mathbb A^2_k\rightarrow \mathbb T^2_{f,g,m,Y}$ and a finite flat morphism $\mathbb T^2_{f,g,m,Y}\rightarrow\mathbb S^2_{f,g,k}$ of degree $pm+1$, with $\mathbb T^2_{f,g,m,Y}$ a smooth surface over $k$ which is a line bundle over $\mathbb P^1_k$ and for which $(\mathbb T^2_{f,g,m,Y}\setminus\imath_{f,g,m,Y})\cong\mathbb A^1_k$.
\end{proposition}

\begin{proof} The proofs of parts (1) and (2) are as in the case of $\mathbb S^2_k$, we only have to define the morphism 
$$\mathcal Q_{f,g,k}:\mathbb S^2_{f,g,k}\rightarrow\mathbb P^1_k$$ 
that generalizes $\mathcal Q_k$. Let $\mathcal Q_{f,g,k}$ be defined by the rule: $(x,y,z)$ maps to either $\bigl(-x:yg(y)\bigr)$ or $\bigl(z:f(x)\bigr)$ (whichever one is defined; if both are defined, then they are equal). As $g(\beta)+\beta g'(\beta)\neq 0$, $Y$ is a connected component of the inverse image $\mathcal Q_{f,g,k}^{-1}(1:0)=\Spec \bigl(k[x,y,z,x^{-1}]/(xf(x),yg(y))\bigr)$. Let 
$$\mathcal Q_{f,g,k,Y}:=\mathcal Q_{f,g,k}|\mathbb S^2_{f,g,k,Y}:\mathbb S^2_{f,g,k,Y}\rightarrow\mathbb P^1_k;$$ so $\mathcal Q_{f,g,k,Y}^{-1}(0:1)=Y$ and hence $\mathcal Q_{f,g,k,Y}$ is a line bundle by Theorem \ref{T10}(2) applied to $q=1$. So parts (1) and (2) hold by taking $\Sigma_{f,g,m,Y}$ to be the pullback of $\Sigma_m^+$ via $\mathcal Q_{f,g,k,Y}$. 

Part (3) holds by taking $\mathbb T^2_{f,g,m,Y}\rightarrow\mathbb S^2_{f,g,k}$ to be the pullback of the finite flat morphism $\Sigma_m^{\proj}:\mathbb P^1_k\rightarrow\mathbb P^1_k$ of Example \ref{EX28} via $\mathcal Q_{f,g,k,Y}$.
\end{proof}

\begin{remark}\normalfont\label{R21}
Let $\mathbb T^2_{f,g,m}\rightarrow\mathbb S^2_{f,g,k}$ be the pullback of $\Sigma_m^{\proj}$ via the flat morphism $\mathcal Q_{f,g,k}:\mathbb S^2_{f,g,k}\rightarrow\mathbb P^1_k$. Then $\mathbb T^2_{f,g,m,Y}$ is an open subvariety of $\mathbb T^2_{f,g,m}$ and the resulting morphism $\mathbb T^2_{f,g,m}\setminus\mathbb T^2_{f,g,m,Y}\rightarrow\mathbb S^2_{f,g,k}\setminus\mathbb S^2_{f,g,k,Y}$ is an isomorphism between reduced affine varieties isomorphic to $\deg(f)\z(tg)-1$ copies of $\mathbb A^1_k$.
\end{remark}

\begin{remark}\normalfont\label{R22}
Let $n:=\deg(f)+1\ge 2$. The composite of $\Sigma_{f,g,m,Y}:\mathbb A^2_k\rightarrow\mathbb S^2_{f,g,k,Y}$ of Proposition \ref{PR12} with the inclusion $\mathbb S^2_{f,g,k,Y}\rightarrow \mathbb S^2_{f,g,k}$ and with a finite \'etale morphism $\Phi_{f,g,h,q}:\mathbb S^2_{f,g,k}\rightarrow\mathbb A^2_k$ as in Proposition \ref{PR8}(1), is a surjective non-finite endomorphism $e=e_{f,g,m,h,q}\in\EE_2(k)$ of degree $pnq(pm+1)$ such that $X_e=\mathbb T^{2,\n}_{f,g,m}$ is the normalization of $\mathbb T^2_{f,g,m}$ and $\psi_e$ factors through the composite of $\mathbb T^{2,\n}_{f,g,m}\rightarrow\mathbb S^2_{f,g,k}$ with $\Phi_{f,g,h,q}$; so $\mathbb T^2_{f,g,m}$ is a flat model of $X_e$. If $tg(t)$ is separable, then $\mathcal Q_{f,g,k}$ is smooth and thus $X_e=\mathbb T^2_{f,g,m}$ is smooth over $\mathbb P^1_k$ and hence regular. As $\mathbb T^2_{f,g,m}\setminus\Imm(\imath_e)=(\mathbb T^2_{f,g,m}\setminus\mathbb T_{f,g,m,Y})\cup \bigl(\mathbb T^2_{f,g,m,Y}\setminus\Imm(\imath_e)\bigr)$, we have $\rho(e)\ge \deg(f)\z(tg)$ by Remark \ref{R8}(1) and the equality holds if $tg(t)$ is separable. As the \'etale locus of $\Sigma_m^{\proj}$ is $\mathbb A^1_k$, if $tg(t)$ is separable we get that $X_e^{\et}=\mathbb A^2_{k,s}$ and thus $\rho_{\et}(e)=0$.\end{remark}%GROUP4

\begin{remark}\normalfont\label{R23}
Assume $\chr(K)=0$. For $l\in\mathbb N$, we do not know if there exist \'etale morphisms $\mathbb A^2_K\rightarrow\mathbb S^2_{t-1,t^l,K}$ with finite complements (cf.\ \cite{Mi5}, Lem.\ 1.9(2) which proves that generically Galois \'etale morphisms $\mathbb A^2_K\rightarrow\mathbb S^2_{t-1,t^l,K}$ with finite complements do not exist).\end{remark}

\begin{remark}\normalfont\label{R24}
Torsors under $\mathcal O_{\mathbb P^1_K}(-q)$ are parametrized by $H^1\bigl(\mathbb P^1_K,\mathcal O_{\mathbb P^1_K}(-q)\bigr)$. As $H^1\bigl(\mathbb P^1_K,\mathcal O_{\mathbb P^1_K}(-2)\bigr)\cong K$, the group of automorphisms of $\mathcal O_{\mathbb P^1_K}(-2)$ acts transitively on the non-trivial torsors under $\mathcal O_{\mathbb P^1_K}(-2)$ and thus, as schemes over $\Spec K$, they are all isomorphic to $\mathbb S^2_K$. By Serre duality for $\mathbb P^1_K$ (see \cite{H2}, Ch.\ III, Thm.\ 7.1(b)), for each non-zero class in $H^1\bigl(\mathbb P^1_K,\mathcal O_{\mathbb P^1_K}(-q)\bigr)$ there exists a homomorphism $\mathcal O_{\mathbb P^1_K}(-q)\rightarrow\mathcal O_{\mathbb P^1_K}(-2)$ that maps it to a non-zero class in $H^1\bigl(\mathbb P^1_K,\mathcal O_{\mathbb P^1_K}(-2)\bigr)$. Thus the non-trivial torsors under $\mathcal O_{\mathbb P^1_K}(-q)$ are affine as they are affine schemes over $\mathbb S^2_K$.

For each integer $q\ge 3$, the surfaces that are non-trivial torsors under the line bundle $\mathcal O_{\mathbb P^1_K}(-q)$ over $\mathbb P^1_K$ are spectra of $K$-algebras described in \cite{Wr2}, Thm.\ 3.1; strictly speaking, loc. cit. works over $\mathbb C$ but its arguments work over any $K$. 

It would be of interest to classify the isomorphism classes defined by $\mathbb S^2_{f,g,k,Y}$ or $\mathbb T^2_{f,g,m,Y}$ of Proposition \ref{PR12} in terms of their $q$ value and spectra.\end{remark}

\section{Affine moduli schemes}\label{S26}

Let $q\in\mathbb N^\ast$. Recall the affine scheme $End_{n,q}=\Spec(A_{n,q})$ over $\Spec\mathbb Z$ and its subsheaves $EE_{n,q}$ and $GA_{n,q}$ of sets defined by the rules on $R$-valued points: 
$$End_{n,q}(R):=\{(f_1,\ldots,f_n)\in R[x_1,\ldots,x_n]^n|\pi(f_1,\ldots,f_n)\le q\},$$
$$EE_{n,q}(R):=\{(f_1,\ldots,f_n)\in End_{n,q}(R)|\e(f_1,\ldots,f_n)\in\EE_n(R)\},$$
$$GA_{n,q}(R):=\{(f_1,\ldots,f_n)\in End_{n,q}(R)|\e(f_1,\ldots,f_n)\in\GA_n(R)\}.$$

For $(l_1,\ldots,l_n)\in\mathbb N^n$ let $\deg(l_1,\ldots,l_n):=\sum_{i=1}^n l_i$. We endow the set
$$\mathbb L_{n,q}:=\{(l_1,\ldots,l_n)\in\mathbb N^n|\deg(l_1,\ldots,l_n)\le q\}$$ 
with the lexicographic order: its smallest and greatest elements are the $n$-tuples $(0,\ldots,0)$ and $(q,0,\ldots,0)$ (respectively). Let 
$$N(n,q):=|\mathbb L_{n,q}|\in\mathbb N^{\ast}.$$ 
As the set $\mathbb L_{n,q}$ is in bijection to the set $\{(l_0,l_1,\ldots,l_n)\in\mathbb N^{n+1}|\sum_{i=0}^n l_i=q\}$ and (e.g., see \cite{BV}, Subsect. 4.2.6) it is well-known that the last set has $\binom{n+q}{n}$ elements, we get that $N(n,q)=\binom{n+q}{n}$.

For $(f_1,\ldots,f_n)\in End_{n,q}(R)$ and $i\in \llbracket1,n\rrbracket$ we write 
$$f_i(x_1,\ldots,x_n)=\sum_{\lambda=(l_1,\ldots,l_n)\in\mathbb L_{n,q}} \alpha^{(i)}_{\lambda}\prod_{j=1}^n x_j^{l_j},$$
where each $\alpha_{\lambda}^{(i)}\in R$, and the rule 
$$(f_1,\ldots,f_n)\mapsto \bigl(\alpha^{(1)}_{(0,\ldots,0)},\ldots,\alpha^{(1)}_{(q,0,\ldots,0)},\ldots,\alpha^{(n)}_{(0,\ldots,0)},\ldots,\alpha^{(n)}_{(q,0,\ldots,0)}\bigr)$$
allows us, by viewing $\alpha^{(i)}_{\lambda}$ as indeterminates, to identify 
$$A_{n,q}=\mathbb Z[\alpha^{(1)}_{(0,\ldots,0)},\ldots,\alpha^{(n)}_{(q,\ldots,0)}]\;\;\;\textup{and}\;\;\;End_{n,q}=\mathbb A^{nN(n,q)}.$$

For $(f_1,\ldots,f_n)\in End_{n,q}(R)$ we have $(f_1,\ldots,f_n)\in EE_{n,q}(R)$ iff its Jacobian determinant belongs to $R[x_1,\ldots,x_n]^{\ast}$.

There exists an $\mathbb N$-bigradation\index{$\mathbb N$-bigradation} on each $A_{n,q}\otimes_{\mathbb Z} R$ that comes from the natural $\mathbb G_m$-actions on $\mathbb Z_{\t}[x_1,\ldots,x_n]$ and $\mathbb Z_{\s}[x_1,\ldots,x_n]$ defined by the rule $(t,x_i)\mapsto tx_i$ for each $i\in \llbracket1,n\rrbracket$. So for $\lambda\in\mathbb L_{n,q}$, each indeterminate $\alpha^{(i)}_{\lambda}$ has bidegree $\bigl(\deg(\lambda),1\bigr)$.

\phantomsection{We have truncation closed embeddings}\label{MS1}
$$\Xi_{n,q}:End_{n,q}\rightarrow End_{n,q+1}$$
that forget all coefficients of degree $q+1$, i.e., are defined by surjective ring homomorphisms 
$$\Xi_{n,q}^{\#}:A_{n,q+1}\rightarrow A_{n,q}$$ that map each $\alpha^{(i)}_{\lambda}$ to itself if $\deg(\lambda)\le q$ and to $0$ if $\deg(\lambda)=q+1$. This allows us to define the projective limit
$A_n$ of the projective system $A_{n,q+1}\rightarrow A_{n,q}$ indexed by $q\in\mathbb N^\ast$ (in each bidegree the inverse limit stabilizes). 

\phantomsection{Let $SEE_{n,q}$ be the closed subscheme of $End_{n,q}$ defined by the following rule on $R$-valued points:}\label{MS3}
$$SEE_{n,q}(R):=\{(f_1,\ldots,f_n)\in End_{n,q}(R)|\det\bigl(J(f_1,\ldots,f_n)\bigr)=1\}.$$

\begin{proposition}\label{PR13} 
Let $l\in\mathbb N$. Let $SA_{n,q,l}$ be the closed subscheme of the product $SEE_{n,q}\times_{\Spec\mathbb Z} SEE_{n,q^{n-1}+l}$ whose $R$-valued points are pairs
$$\bigl((f_1,\ldots,f_n),(g_1,\ldots,g_n)\bigr)\in SEE_{n,q}(R)\times SEE_{n,q^{n-1}+l}(R)$$
with the property that we have $\e(f_1,\ldots,f_n)\e(g_1,\ldots,g_n)=\e(x_1,\ldots,x_n)$. If $R$ is reduced, then the first projection induces an identity
$$SA_{n,q,l}(R)=GA_n(R)\cap SEE_{n,q}(R)$$
and thus the reduced scheme $SA_{n,q}:=(SA_{n,q,l})_{\red}$ does not depend on $l\in\mathbb N$.
\end{proposition}

\begin{proof} 
To show that such endomorphisms $\e(f_1,\ldots,f_n)$ and $\e(g_1,\ldots,g_n)$ are automorphisms we can assume first that $R$ is a field by \cite{Gro5}, Cor.\ (17.9.5) and second that $R$ is an algebraically closed field. Thus $\deg\bigl(\e(g_1,\ldots,g_n)\bigr)=1$ and hence $\e\bigl(g_1,\ldots,g_n\bigr)\in\GA_n(R)$ by Lemma \ref{L6}(4). So $\e(f_1,\ldots,f_n)=\e(g_1,\ldots,g_n)^{-1}$ is also an automorphism. 

\phantomsection{Thus the first projections $SEE_{n,q}(R)\times SEE_{n,q^{n-1}+l}(R)\rightarrow SEE_{n,q}(R)$ induce functorial maps}\label{MS4}
$$\Upsilon_{n,q,l}(R): SA_{n,q,l}(R)\rightarrow\GA_n(R)\cap SEE_{n,q}(R)$$ 
which are injective for every $R$, and to end the proof it suffices to prove that $\Upsilon_{n,q,l}(R)$ are in fact surjective if $R$ is reduced. To show this, we can assume that $R$ is a finitely generated $\mathbb Z$-algebra, let $(f_1,\ldots,f_n)\in\GA_n(R)\cap SEE_{n,q}(R)$ and let $(g_1,\ldots,g_n)\in\GA_n(R)$ be such that $\e(g_1,\ldots,g_n)=\e(f_1,\ldots,f_n)^{-1}$. We are left to show that actually we have $(g_1,\ldots,g_n)\in SEE_{n,q^{n-1}}(R)$. As $R$ is reduced, by replacing it with fields that are direct factors of the total ring of fractions of $R$, to show that $(g_1,\ldots,g_n)\in SEE_{n,q^{n-1}}(R)$, i.e., $\pi(g_1,\ldots,g_n)\le q^{n-1}$, we can assume $R$ is a field and in this case the inequality follows from the proof of \cite{BCW2}, Cor.\ (1.6).\end{proof} 

\begin{remark}\normalfont\label{R26}
The scheme $SA_{n,q,K}$ is not a priori reduced but its reduction is canonically identified with $(SA_{n,q,l,K})_{\red}$ for all $l\in\mathbb N$.
\end{remark}

\phantomsection{The injective maps $\Upsilon_{n,q,l}(R)$ define natural morphisms $SA_{n,q,l}\rightarrow SEE_{n,q}$ and thus also a natural morphism $\Upsilon_{n,q}:SA_{n,q}\rightarrow SEE_{n,q}$.}\label{F5+}

\begin{lemma}\label{F5}
The morphism $\Upsilon_{n,q}:SA_{n,q}\rightarrow SEE_{n,q}$ is a closed embedding. 
\end{lemma}

\begin{proof}
The image $S:=\Imm\bigl(\mathcal O(SEE_{n,q})\rightarrow\mathcal O(SA_{n,q})\bigr)$ is reduced. As $\Upsilon_{n,q,l}(S)$ are isomorphisms, from \cite{vdE}, Part I, Ch.\ 1, Prop.\ 1.1.7 applied to the inclusion $S\rightarrow \mathcal O(SA_{n,q})$ we get that this inclusion is an isomorphism.\end{proof}

We view the injective maps $\Xi_{n,q}(R):End_{n,q}(R)\rightarrow End_{n,q+1}(R)$ as inclusions and we also speak about $SA_{n,q}(R)\subset SA_{n,q+1}(R)$, $SEE_{n,q}(R)\subset SEE_{n,q+1}(R)$, etc. E.g., we identify $SA_{n,q,l}$ with a closed subscheme of $SA_{n,q,l+1}$ for all $l\in\mathbb N$.

\begin{remark}\normalfont\label{R25}
The $R$-valued points rule 
$$EE^{\ast}_{n,q}(R):=\{(f_1,\ldots,f_n)\in EE_{n,q}(R)|\det\bigl(J(f_1,\ldots,f_n)\bigr)\in R^{\ast}\}$$ 
defines a locally closed subscheme $EE^{\ast}_{n,q}$ of $End_{n,q}$. As the group monomorphism $R^{\ast}\subset R[x_1,\ldots,x_n]^{\ast}$ is an isomorphism iff $R$ is reduced, the functorial inclusion $EE^{\ast}_{n,q}(R)\subset EE_{n,q}(R)$ is an identity iff $R$ is reduced. Note that 
$$EE^{\ast}_{n,q}\cong\mathbb G_m\times_{\Spec \mathbb Z} SEE_{n,q}$$ 
(e.g., by multiplying $x_1$ in the source or the target by invertible scalars). If in Proposition \ref{PR13} we replace $SEE$ by $EE^{\ast}$ we obtain affine schemes $GA^{\ast}_{n,q,l}$ over $\Spec\mathbb Z$ whose reductions do not depend on $l$ and define $GA^{\ast}_{n,q}$; we have isomorphisms 
$$GA^{\ast}_{n,q,l}\cong\mathbb G_m\times_{\Spec \mathbb Z} SA_{n,q,l}\;\;\;\textup{and}\;\;\;GA^{\ast}_{n,q}\cong\mathbb G_m\times_{\Spec \mathbb Z} SA_{n,q}.$$ 
From Lemma \ref{F5} we get that the morphism $GA^{\ast}_{n,q}\rightarrow EE^{\ast}_{n,q}$ is a closed embedding.
\end{remark}
%following two rules $SA^{\ast}_{n,q}(R)=EE^{\ast}_{n,q}(R)\cap GA_{n,q}(R)$ and closed subscheme $GA^{\ast}_{n,q}$ of $EE^{\ast}_{n,q}$ and $GA^{\ast}_{n,q}(R)=GA_{n,q}(R)$ and 
%The scheme $GA_{n,q}$ is of finite presentation over $\Spec\mathbb Z$ and from \cite{vdE}, Ch.\ 1, Lem.\ 1.1.9 we get that it is formally smooth over $\Spec\mathbb Z$ in the sense of \cite{Gro}, Def. (17.3.1). Thus $GA_{n,q}$ is smooth over $\Spec\mathbb Z$. It is easy to see that this implies that $SA_{n,q}$ is also smooth over $\Spec\mathbb Z$.

\phantomsection{We write $SEE_{n,q}=\Spec(B_{n,q})$, where $B_{n,q}=A_{n,q}/I_{n,q}$ is a quotient of $A_{n,q}$ compatible with the $\mathbb N$-bigradings. As $\Xi_{n,q}^{\#}(I_{n,q+1})\subset I_{n,q}$, the inverse limit}\label{MS2} 
$$B_n=A_n/I_n$$ 
of the projective system $B_{n,q+1}\rightarrow B_{n,q}$ indexed by $q\in\mathbb N^\ast$ is well-defined. 

Let $\overline{\Spec(B_{n,q,\mathbb Q})}$ be the schematic closure of $\Spec(B_{n,q,\mathbb Q})$ in $\Spec\bigl(B_{n,q,\mathbb Z_{(p)}}\bigr)$ and let $\overline{\Spec(B_{n,\mathbb Q})}$ be the schematic closure of $\Spec(B_{n,\mathbb Q})$ in $\Spec\bigl(B_{n,\mathbb Z_{(p)}}\bigr)$. The natural closed embedding
$$\Spec(B_{n,q})\rightarrow\Spec(B_n)$$ induces a closed embedding
$$\mathfrak C_{n,q}:\overline{\Spec(B_{n,q,\mathbb Q})}\rightarrow \overline{\Spec(B_{n,\mathbb Q})}\cap \Spec\bigl(B_{n,q,\mathbb Z_{(p)}}\bigr)$$
which is not always an isomorphism, as the next example shows.

\begin{example}\normalfont\label{EX39}
The pair $(x_1-x_1^p,x_2)\in\EE_{2,p}(\mathbb F_p)$ can be lifted to a formal pair
$\bigl(x_1-x_1^p,x_2(1-px_1^{p-1})^{-1}\bigr)$ over $\mathbb Z_p$ whose Jacobian matrix has determinant $1$, hence
$$(x_1-x_1^p,x_2)\in [\overline{\Spec(B_{2,\mathbb Q})}\cap \Spec\bigl(B_{2,p,\mathbb Z_{(p)}}\bigr)](\mathbb F_p).$$
But if $p$ is small, say $p=2$ by \cite{Wa2}, Thm.\ 62 (based on \cite{Moh1}, Sect.\ 6 and \cite{Wa1} one could take $p\le 97$ and based on \cite{Ng}, Thm.\ 3.6 one could take $p\le 103$), then 
$$(x_1-x_1^p,x_2)\notin\overline{\Spec(B_{2,p,\mathbb Q})}(\mathbb F_p)$$
which implies that $\mathfrak C_{2,p}$ is not an isomorphism. This disproves the equivalence used in \cite{MR}, Def.\ 3.3 within parentheses and implicitly disproves \cite{MR}, Conj.\ 3.4 if the first part (not within parentheses) of \cite{MR}, Def.\ 3.3 is used.\footnote{In \cite{MR}, Conj.\ 3.4, if the second part (within parentheses) of \cite{MR}, Def.\ 3.3 is used, is false in characteristic $p$ and is only a reformulation of the Jacobian Conjecture in characteristic $0$.}
\end{example}

The following general lemma helps in keeping track of $q$ in mixed characteristic $(0,p)$ in the case when $q\in p\mathbb N^\ast$.

\begin{lemma}\label{L18}
For $i\in \llbracket1,n\rrbracket$ let $r_i\in\mathbb N$ and let $g_i(x_1,\ldots,x_n)\in K[x_1,\ldots,x_n]$ be homogeneous of degree $r_i$. We consider the set $\mathcal X$ of $n$-tuples $(\alpha_1,\ldots,\alpha_n)\in (K^{\ast})^n$ with the property that the homogeneous system of $n$ equations
\begin{equation}\label{EQ28}
g_1(x_1,\ldots,x_n)+\alpha_1x_1^{r_1}=\cdots=g_n(x_1,\ldots,x_n)+\alpha_nx_n^{r_n}=0
\end{equation}
has no solution in $\mathbb P^{n-1}(K)$. Then $\mathcal X$ is Zariski dense in $\mathbb A^n_K$.
\end{lemma}

\begin{proof} If there exists $i\in \llbracket1,n\rrbracket$ such that $r_i=0$, so $g_i\in K$, by taking $\alpha_i\in K\setminus\{-g_i\}$ the lemma holds. Thus we can assume that we have $r_i\in\mathbb N^\ast$ for all $i\in \llbracket1,n\rrbracket$. 

For each $i\in \llbracket1,n\rrbracket$ let $t_i$ be an indeterminate. Let $\mathcal H\in K[t_1,\ldots,t_n]$ be the normalized resultant of the $n$ homogeneous polynomials $x_i^{r_i}+t_ig_i(x_1,\ldots,x_n)$ with $i\in\llbracket1,n\rrbracket$ (see \cite{vdW}, Ch.\ XI, Sect.\ 82, p.\ 14 or \cite{Jo2}, Sect.\ 2, Prop.\ 2.3 (i) or \cite{Dem}, Sect.\ 4, Def.\ 3); we have $\mathcal H(0,\ldots,0)=1$. Thus, as $K$ is infinite, there exists $(\beta_1,\ldots,\beta_n)\in (K^{\ast})^n$ such that $\mathcal H(\beta_1,\ldots,\beta_n)\neq 0$. Properties of resultants show that the system of homogeneous equations $x_i^{r_i}+\beta_ig_i(x_1,\ldots,x_n)=0$ indexed by $i\in \llbracket1,n\rrbracket$ has no solution in $\mathbb P^{n-1}(K)$ (see \cite{vdW}, Ch.\ XI, Sect.\ 82, p.\ 15 or \cite{Jo2}, Sect.\ 1 or \cite{Dem}, Sect.\ 4, Sch.). Hence, if $\alpha_i:=\beta_i^{-1}$ for all $i\in \llbracket1,n\rrbracket$, then the System (\ref{EQ28}) has no solution in $\mathbb P^{n-1}(K)$, and thus we have an inclusion of sets of $k$-valued points 
$$\{(\beta_1^{-1},\ldots,\beta_n^{-1})\in (K^{\ast})^n|\mathcal H(\beta_1,\ldots,\beta_n)\neq 0\}\subset\mathcal X$$ 
from which the lemma follows.\end{proof} 

Let $q_+:=p\lceil \frac{q}{p}\rceil\in p\mathbb N^\ast$; so $q\le q_+$, and we have $q=q_+$ iff $q\in p\mathbb N^\ast$. 

\begin{proposition}\label{PR18} 
For $q\in\mathbb N^\ast$ the following properties hold.

\medskip
{\bf (1)} Each $e\in EE_{n,q}(K)$ (resp.\ $e\in SEE_{n,q}(K)$) can be deformed to an element $a\in GA_{n,1}(K)\subset GA_{n,q}(K)$ (resp.\ $a\in SA_{n,1}(K)\subset SA_{n,q}(K)$) using polynomials in $K_{\s}[x_1,\ldots,x_n,t]$ of degree at most $2q-1$ and of total degree in $x_1,\ldots,x_n$ at most $q$. Moreover, if $e\in GA_{n,q}(K)$ (resp.\ $e\in SA_{n,q}(K)$), then we can assume that the deformation, by evaluating $t$ at arbitrary elements of $K$, defines elements of $GA_{n,q}(K)$ (resp.\ $SA_{n,q}(K)$).

\smallskip
{\bf (2)} Let $d\in EE_{n,q}(k)$ (resp.\ $d\in SEE_{n,q}(k)$) and let $q_+$ be as above. Then $d\in\EE_{n,q_+}(k)$ (resp.\ $d\in SEE_{n,q_+}(k)$) is the specialization of a finite endomorphism in $EE_{n,q_+}\bigl(k(t)\bigr)$ (resp.\ in $SEE_{n,q_+}\bigl(k(t)\bigr)$) of geometric degree $q_+^n$.
\end{proposition}

\begin{proof} If $e$ is defined by an $n$-tuple $(f_1,\ldots,f_n)\in K_{\s}[x_1,\ldots,x_n]^n$, then let 
$$h_{i,t}=h_i(x_1,\ldots,x_n,t):=\frac{f_i(tx_1,\ldots,tx_n)+(t-1)f_i(0,\ldots,0)}{t}\in K_{\s}[x_1,\ldots,x_n,t];$$ 
its degree is at most $2q-1$ and its total degree in $x_1,\ldots,x_n$ is at most $q$. For each $t_0\in K$, let $h_{i,t_0}:=h_i(x_1,\ldots,x_n,t_0)\in K_{\s}[x_1,\ldots,x_n]$. We get an \'etale endomorphism $e_{t_0}:=\e(h_{1,t_0},\ldots,h_{n,t_0})\in EE_{n,q}(K)$, with $e_1=e$ and, as $\pi(e_0)=1$, with $e_0\in GA_{n,1}(K)\subset GA_{n,q}(K)$ having the same Jacobian determinant as $e$. If $e\in GA_{n,q}(K)$ (resp.\ $e\in SA_{n,q}(K)$), then for each $t_0\in K$ we have $e_{t_0}\in GA_{n,q}(K)$ (resp.\ $e_{t_0}\in SA_{n,q}(K)$). So part (1) holds.

For part (2), let $d$ be defined by an $n$-tuple $(f_1,\ldots,f_n)\in k_{\s}[x_1,\ldots,x_n]^n$ and for each $i\in \llbracket1,n\rrbracket$ let $g_i\in k_{\s}[x_1,\ldots,x_n]$ be $0$ if $\deg(f_i)<q_+$ (e.g., if $q<q_+$) and be the homogeneous part of $f_i$ of degree $q_+$ if $\deg(f_i)=q=q_+$. From Lemma \ref{F3}(2) and from Lemma \ref{L18} applied to $r_1=\cdots=r_n=q_+$ and $(g_1,\ldots,g_n)$ we get that there exists $(\alpha_1,\ldots,\alpha_n)\in (k^{\ast})^n$ such that 
$$d_1:=\e(f_1+\alpha_1x_1^{q_+},\ldots,f_1+\alpha_nx_n^{q_+})\in\EE_{n,q_+}(k)$$ 
is finite with $\deg(d_1)=q_+^n$. Then 
$$E:=\e(f_1+t\alpha_1x_1^{q_+},\ldots,f_n+t\alpha_nx_n^{q_+})\in\EE_{n,q^+}(k[t])$$ specializes for $t=0$ to $d$ and for $t=1$ to $d_1$. The last specialization implies based on the proof of Lemma \ref{L18} that the generic fiber $E_{k(t)}\in\EE_{n,q_+}\bigl(k(t)\bigr)$ of $E$ is such that Lemma \ref{F3}(2) applies to it giving that $\deg(E_{k(t)})=q_+^n$.\footnote{We also have $\deg(E_{k(t)})\ge q_+^n$ by \cite{Gro5}, Prop.\ (18.2.8) and $\deg(E_{k(t)})\le q_+^n$ by Lemma \ref{F3}(1).} Also, $d\in SEE_{n,q}(k)\Rightarrow E\in SEE_{n,q^+}(k[t])$. As $E_{k(t)}$ specializes to $d$, part (2) holds.\end{proof} 

\begin{remark}\normalfont\label{R27} 
The existence of $E$ of the proof of Proposition \ref{PR18}(2) disproves \cite{MR}, Conj.\ 5.3(2) for each integral $k$-algebra that contains $k[t]$ and has $k$ as a quotient ring.\end{remark}

\begin{corollary}\label{C9}
Let $q\in\mathbb N^\ast$. Then the affine schemes $SA_{n,q,K}$ and $SEE_{n,q,K}$ (hence also $EE^{\ast}_{n,q,K}$) over $\Spec K$ are connected. 
\end{corollary}

\begin{proof} As $SA_{n,1,K}$ is connected (being the smooth group of affine transformations of $K^n$ of determinant $1$), this follows from Proposition \ref{PR18}(1).\end{proof}

\begin{lemma}\label{L18.1}
For $(n,q)\in (\mathbb N^{\ast})^2$ and $m\in\llbracket1,q^n\rrbracket$ the following properties hold.

\medskip
{\bf (1)} There exists a reduced closed subscheme $SEE_{n,q}^{\le m}$ of $SEE_{n,q}$ whose geometric points are those for which we get \'etale endomorphisms of the affine spaces of dimension $n$ of geometric degree at most $m$.

\smallskip
{\bf (2)} We have $SEE_{n,q}^{\le q^n}=(SEE_{n,q})_{\red}$.

\smallskip
{\bf (3)} We have $SEE_{n,q}^{\le 1}=SA_{n,q}$.

\smallskip
{\bf (4)} Parts (1) to (3) hold with the pair $(SEE,SA)$ replaced by $(EE^{\ast},GA^{\ast})$.
\end{lemma}

\begin{proof}
We consider the universal \'etale endomorphism $e_{n,q}^{\univ}:\mathbb A^n_{SEE_{n,q},\s}\rightarrow \mathbb A^n_{SEE_{n,q},\t}$. From either \cite{Gro4}, Prop.\ (15.5.1)(i) or \cite{Gro5}, Prop.\ (18.2.8) we get that there exists an open subscheme $U_{n,q}$ of $\mathbb A^n_{SEE_{n,q},\t}$ whose geometric points are those $P$ for which  $|(e_{n,q}^{\univ})^{-1}(P)|>m$. As the projection $A^n_{SEE_{n,q},\t}\rightarrow SEE_{n,q}$ is flat, the image of $U_{n,q}$ under it is an open non-empty subscheme $SEE^{>m}_{n,q}$ of $SEE_{n,q}$ by \cite{Gro3}, Thm.\ (2.4.6) and we take $SEE_{n,q}^{\le m}$ to be its complement in $SEE_{n,q}$ endowed with the reduced structure. So part (1) holds. 

Part (2) follows from Lemma \ref{F3}(1).

From part (1) applied to $m=1$ and Lemma \ref{L6}(4) we get that the two reduced closed subschemes $SEE_{n,q}^{\le 1}$ and $SA_{n,q}$ of $SEE_{n,q}$ have the same geometric points. Thus part (3) holds.

Part (4) follows from parts (1) to (3) and Remark \ref{R25} which mentions isomorphisms $EE^{\ast}_{n,q}\cong\mathbb G_m\times_{\Spec \mathbb Z} SEE_{n,q}$ and $GA^{\ast}_{n,q}\cong\mathbb G_m\times_{\Spec \mathbb Z} SA_{n,q}$ and implicitly their analogs $EE_{n,q}^{\ast,\le m}\cong\mathbb G_m\times_{\Spec \mathbb Z} SEE^{\le m}_{n,q}$ with $m\in\llbracket1,q^n\rrbracket$.\end{proof}

\begin{proposition}\label{C9+} Let $q\in\mathbb N^{\ast}$. Let $\deg_{p,q,n}\in\mathbb N^{\ast}$ (resp.\ $\deg^{\f}_{p,q,n}\in\mathbb N^{\ast}$) be the largest such that there exists an (resp.\ a finite) endomorphism $e\in SEE_{n,q}(k)$ such that $\deg(e)=\deg_{p,q,n}$ (resp.\ $\deg(e)=\deg^{\f}_{p,q,n}$). Then the following properties hold.

\medskip
{\bf (1)} We have inequalities $\deg^{\f}_{p,q,n}\le \deg_{p,q,n}\le q^n$.

\smallskip
{\bf (2)} If $p\mid q$, then $\deg^{\f}_{p,q,n}=\deg_{p,q,n}=q^n$.

\smallskip
{\bf (3)} If $p\nmid q$ and $q>p$, then $pq^{n-1}\lfloor\frac{q}{p}\rfloor\le\deg^{\f}_{p,q,n}$.

\smallskip
{\bf (4)} There exists an open non-empty subscheme $SEE^{\gen}_{n,q,k}$ of $SEE_{n,q,k}$ such that $SEE^{\gen}_{n,q,k}(k)=\{e\in SEE_{n,q,k}(k)=SEE_{n,q}(k)|\deg(e)=\deg_{p,q,n}\}$.

\smallskip
{\bf (5)} There exists a constructible subset $SEE^{\f}_{n,q,k}$ of $SEE_{n,q,k}$ formed by points $P$ over whose residue fields $\kappa(P)$ we get finite \'etale endomorphisms of $\mathbb A^n_{\kappa(P)}$.

\smallskip
{\bf (6)} There exists a constructible subset $SEE^{\f-\gen}_{n,q,k}$ of $SEE_{n,q,k}$ contained in $SEE^{\f}_{n,q,k}$ and formed by points $P$ over whose residue fields $\kappa(P)$ we get finite \'etale endomorphisms of $\mathbb A^n_{\kappa(P)}$ of degrees $\deg^{\f}_{p,q,n}$.

\smallskip
{\bf (7)} If $p\mid q$, then $SEE^{\f-\gen}_{n,q,k}$ is a constructible subset of $SEE^{\gen}_{n,q,k}$ which is Zariski dense in $(SEE_{n,q,k})_{\red}$.
\end{proposition}

\begin{proof} 
The inequality $\deg^{\f}_{p,q,n}\le\deg_{p,q,n}$ is clear and the inequality $\deg_{p,q,n}\le q^n$ follows from Lemma \ref{F3}(1); so part (1) holds. 

If $p\mid q$, then $e:=\e(x_1-x_1^q,\ldots,x_n-x_n^q)\in SEE_{n,q}(k)$ is finite by Lemma \ref{F3}(2) with $\deg(e)=q^n$. From this and part (1) we get that part (2) holds. 

If $p\nmid q$ with $q>p$, then $e:=\e(x_1-x_2^q,\ldots,x_{n-1}-x_n^q, x_n-x_1^{p\lfloor\frac{q}{p}\rfloor})\in SEE_{n,q}(k)$ is finite by Lemma \ref{F3}(2) with $\deg(e)=pq^{n-1}\lfloor\frac{q}{p}\rfloor$, so part (3) holds.

Part (4) follows from Lemma \ref{L18.1}(1) applied to $m=\deg_{p,q,n}-1$, where by convention $SEE_{n,q,k}^{\le 0}$ is the empty scheme: the complement of $SEE^{\gen}_{n,q,k}$ in $SEE_{n,q,k}$ is $(SEE_{n,q,k}^{\le \deg_{p,q,n}-1})_{\red}$ and $SEE^{\gen}_{n,q,k}$ it is non-empty by the definition of $\deg_{p,q,n}$.

Part (5) is a particular case of \cite{Gro4}, Prop.\ (9.6.1)(vi).

For part (6), if $\deg^{\f}_{p,q,n}=\deg_{p,q,n}$, then $SEE^{\f,\gen}_{n,q,k}=SEE^{\gen}_{n,q,k}\cap SEE^{\f}_{n,q,k}$. If $\deg^{\f}_{p,q,n}<\deg_{p,q,n}$ then we can take 
$$SEE^{\f,\gen}_{n,q,k}=\bigl(SEE^{\f,\gen}_{n,q,k}\cap SEE^{\le\deg^{\f}_{p,q,n}}_{n,q,k}\bigr)_{\red}\setminus \bigl(SEE^{\le\deg^{\f}_{p,q,n}-1}_{n,q,k}\bigr)_{\red}.$$ So part (6) holds.

For part (7), the inclusion $SEE^{\f,\gen}_{n,q,k}\subset SEE^{\gen}_{n,q,k}$ follows from part (2) and the Zariski density follows from Proposition \ref{PR18}(2).\end{proof}

\begin{notation}\normalfont\label{NOT5}
For $n\in\mathbb N^{\ast}\setminus\{1\}$, let $\mathfrak D_n$ be the set of all pairs $(m,q)\in (\mathbb N^{\ast}\setminus\{1\})^2$ such that $m\le q^n$ (cf.\ Lemma \ref{F3}(1)) and $\{e\in\EE_n(\mathbb C)|\deg(e)=m,\pi(e)\le q\}=\emptyset$.
\end{notation}

\begin{remark}\normalfont\label{R19.4}
The classical Jacobian Conjecture in dimension $2$ is equivalent to the identity $\mathfrak D_2=\{(m,q)\in (\mathbb N^{\ast}\setminus\{1\})^2|m\le q^2\}$. From  \cite{\.Z}, Thm.\ 6.12 we get that $\{2,3,4,5\}\times (\mathbb N^{\ast}(\setminus\{1,2\})\subset\mathfrak D_2$. If $q\in\mathbb N\setminus\{1\}$ is small, say $q=2$ by \cite{Wa2}, Thm.\ 62 (based on \cite{Moh1}, Sect.\ 6 and \cite{Wa1} one could take $q\le 100$ and based on \cite{Ng}, Thm.\ 3.6 one could take $q\le 104$), then $\llbracket2,q^2\rrbracket\times\{q\}\subset\mathfrak D_2$.
\end{remark}

\begin{corollary}\label{C9.1}
For $n\in\mathbb N^{\ast}\setminus\{1\}$ and $(m,q)\in\mathfrak D_n$ the following properties hold.

\medskip
{\bf (1)} There exists a finite set of primes $P_{n,m,q}$ such that for each prime $p$, we have $p\in P_{n,m,q}$ iff for one (equivalently, for each) algebraically closed field $k$ of characteristic $p$ the set $\{e\in\EE_n(k)|\deg(e)=m,\pi(e)\le q\}$ is non-empty.

\smallskip
{\bf (2)} There exists a smallest prime $p_{n,m,q}$ such that for each algebraically closed field $\kappa$ with $\chr(\kappa)\ge p_{n,m,q}$, the set $\{e\in\EE_n(\kappa)|\deg(e)=m,\pi(e)\le q\}$ is empty.
\end{corollary}

\begin{proof}
We consider the open subscheme $SEE_{n,q}^{=m}:=SEE^{\le m}_{n,q}\setminus SEE^{\le m-1}_{n,q}$ of $SEE^{\le m}_{n,q}$ by Lemma \ref{L18.1}(1). As $(m,q)\in\mathfrak D_n$, $SEE_{n,q}^{=m}$ has no characteristic $0$ points. From this and the fact that $SEE_{n,q}^{=m}$ is a quasi-affine finitely generated scheme over $\Spec \mathbb Z$, it follows that there exists a finite set of primes $P_{n,m,q}$ such that for a prime $p$ we have $p\in P_{n,m,q}$ iff $SEE_{n,q}^{=m}$ has a point of characteristic $p$ (equivalently, has a $k$-valued point or a point of finite residue field). So part (1) holds.

Part (2) holds for the smallest prime $p_{n,m,q}$ with $p_{n,m,q}>\max(P_{n,m,q})$.
\end{proof}

\begin{example}\normalfont\label{EX40}
Suppose that $p=n=q=2$. Each pair $(f_1,f_2)\in SEE_{2,2,k}(R)$ is of the form 
$$f_i(x_1,x_2)=\alpha^{(i)}_{(0,0)}+\alpha^{(i)}_{(1,0)}x_1+\alpha^{(i)}_{(0,1)}x_2+\alpha^{(i)}_{(2,0)}x_1^2+\alpha^{(i)}_{(1,1)}x_1x_2+\alpha^{(i)}_{(0,2)}x_2^2$$
and, as $1=-1$, one computes that $J(f_1,f_2)(x_1,x_2)$ is equal to
$$\alpha^{(1)}_{(1,0)}\alpha^{(2)}_{(0,1)}+\alpha^{(1)}_{(0,1)}\alpha^{(2)}_{(1,0)}+\bigl(\alpha^{(1)}_{(1,0)}\alpha^{(2)}_{(1,1)}+\alpha^{(1)}_{(1,1)}\alpha^{(2)}_{(1,0)}\bigr)x_1+\bigl(\alpha^{(1)}_{(1,1)}\alpha^{(2)}_{(0,1)}+\alpha^{(1)}_{(0,1)}\alpha^{(2)}_{(1,1)}\bigr)x_2.$$
Thus 
$$A_{2,2}\otimes_{\mathbb Z} k=k[\alpha^{(i)}_{(0,0)},\alpha^{(i)}_{(1,0)},\alpha^{(i)}_{(0,1)},\alpha^{(i)}_{(2,0)},\alpha^{(i)}_{(1,1)},\alpha^{(i)}_{(0,2)}|i\in \{1,2\}]$$
and $\Imm(I_{2,2}\rightarrow A_{2,2}\otimes_{\mathbb Z} k)$ is the ideal
$$\bigl(\alpha^{(1)}_{(1,0)}\alpha^{(2)}_{(0,1)}+\alpha^{(1)}_{(0,1)}\alpha^{(2)}_{(1,0)}+1,\alpha^{(1)}_{(1,0)}\alpha^{(2)}_{(1,1)}+\alpha^{(1)}_{(1,1)}\alpha^{(2)}_{(1,0)},\alpha^{(1)}_{(1,1)}\alpha^{(2)}_{(0,1)}+\alpha^{(1)}_{(0,1)}\alpha^{(2)}_{(1,1)}\bigr).$$
Therefore $B_{2,2}\otimes_{\mathbb Z} k=(A_{2,2}/I_{2,2})\otimes_{\mathbb Z} k$ is isomorphic to a polynomial algebra in $6$ indeterminates (they are $\alpha_{0,0}^{(i)}$, $\alpha_{2,0}^{(i)}$, and $\alpha_{0,2}^{(i)}$ with $i\in \{1,2\}$) over the $k$-algebra
$$B:=k[x,y,z,x',y',z']/(xx'+yy'+1,xz'+zy',zx'+yz')$$
which is isomorphic to $k[x,y,x',y']/(xx'+yy'+1)$ as $z'=z=0$, i.e., we have $\alpha_{1,1}^{(1)}=\alpha^{(2)}_{1,1}=0$. This is so as the three relations defining $B$ can be rewritten as a matrix identity $\begin{bmatrix} 
x & y' \\
y & x' \\ 
\end{bmatrix}\begin{bmatrix} 
z' \\
z \\ 
\end{bmatrix}=\begin{bmatrix} 
0 \\
0 \\ 
\end{bmatrix}$ in which the $2\times 2$ matrix has determinant $1$. Thus $\Spec B$ is an integral affine smooth scheme over $\Spec k$ isomorphic to $\SL_{2,k}$. Hence $SEE_{2,2,k}$ is an affine smooth irreducible variety over $k$ isomorphic to $\SL_{2,k}\times_{\Spec k} \mathbb A^6_k$ and hence of dimension $9$. 

We show that $SA_{2,2,k}$ has dimension $6$ and $3$ irreducible components. As $SA_{2,1,k}$ has dimension $5$ and is irreducible and as $SA_{2,1,k}$ acts on $SA_{2,2,k}$ via left or right translations, it suffices to show that the affine scheme $W$ over $k$ that involves the eight indeterminates $y_1$ to $y_4$ and $z_1$ to $z_4$ and is defined by the composite equation
$$\e(x_1+y_1x_1^2+y_2x_2^2,x_2+y_3x_1^2+y_4x_2^2)\e(x_1+z_1x_1^2+z_2x_2^2,x_2+z_3x_1^2+z_4x_2^2)=\e(x_1,x_2)$$ 
has dimension $1$ and $3$ irreducible components. This composite equation implies, as $1=-1$, that we have $z_i=y_i$ for all $i\in\{1,2,3,4\}$ and that 
$$W\cong\Spec\Bigl(k[y_1,y_2,y_3,y_4]/\bigl(y_1^3+y_2y_3^2,y_4^3+y_1^2y_3,y_1(y_2+y_4)^2,y_4(y_1+y_3)^2\bigr)\Bigr)$$
is non-reduced. The irreducible components of $W$ are the zero loci $y_1=y_4=y_2=0$, $y_1=y_4=y_3=0$, and $y_1=y_2=y_3=y_4$. 

As $SA_{2,2,0,k}\cong SA_{2,1,k}\times_{\Spec k} W$, $SA_{2,2,0,k}$ is also non-reduced.

Based on Remark \ref{R25} and the above identities $\alpha_{1,1}^{(1)}=\alpha^{(2)}_{1,1}=0$ we get that if $e\in\EE_2(k)$ is such that $\pi(e)=2$, then $e\in\BE_2(k)$, and from the proof of \cite{V}, Thm.\ 2.4.1(b) we get that $\mathcal F(e)=\emptyset$, i.e., $e$ is finite, and $\deg(e)\in\{2,4\}$. Thus we have identities $SEE^{\f}_{2,2,k}=SEE_{2,2,k}$ and $SEE^{\f-\gen}_{2,2,k}=SEE^{\gen}_{2,2,k}$ of constructible sets.

Moreover, with the notation of Corollary \ref{C9.1}, we have $2\in P_{2,2,2}$. From this and the last sentence of Section \ref{S2} we get that $P_{2,2,2}=\{2\}$. Hence $p_{2,2,2}=3$. 
\end{example}

\begin{example}\normalfont\label{EX41}
Assume $p=2$. The rule $(x,y)\mapsto (x,y+xy^2)$ defines $d\in\EE_2(k)$ with $\deg(d)=2$ and $\mathcal F(d)=\{1\}$. From Proposition \ref{PR1}(1) to (3) we get relations $\pi(d)\le\pi_{\r}(d)=\pi_{\l}(d)=\pi_{\i}(d)=3$. From the previous paragraph we get that $\pi(d)\neq 2$. As $d\notin\GA_2(k)$, we have $\pi(d)\neq 1$. We conclude that $\pi(d)=3$.
\end{example}

For an $n\times n$ matrix $\Box$ with entries in $R$ and $i\in \llbracket1,n\rrbracket$, let $\tr_i(\Box)\in R$ be the trace of the $\binom{n}{i}\times\binom{n}{i}$ matrix which is the $i$-th exterior power $\bigwedge^{\raisebox{-0.0ex}{\scriptsize $i$}}(\Box)$ of $\Box$. If $1_n$ is the identity matrix of size $n\times n$, then we have the following formula
 \begin{equation}\label{EQ29}
 \det(1_n+\Box)=1+\sum_{i=1}^n \tr_i(\Box)
 \end{equation}
 which follows from the fact that $\Box$ is similar to an upper triangular matrix.

\begin{theorem}\label{T7}
Let $q\in\mathbb N^\ast$. Then $SEE_{n,q,K}$ (equivalently, $EE^{\ast}_{n,q,K}$) is smooth over $\Spec K$ iff either $n=q=\chr(K)=2$ or $\min(n,q)=1$, in which case it is also irreducible. 
\end{theorem}

\begin{proof} If $\chr(K)=2$, then $SEE_{2,2,K}$ is smooth and irreducible (see Example \ref{EX40}). Clearly, $SEE_{n,1}=SA_{n,1}\cong\SL_n\times_{\Spec\mathbb Z} \mathbb A^n$. If $\chr(K)=0$, then 
$$SEE_{1,q,K}\cong SEE_{1,1,K}=SA_{1,1,K}=\Spec\bigl(K[\alpha^{(1)}_{(0)}]\bigr)\cong\mathbb A^1_K.$$ 
If $\chr(K)=p$, then
$$SEE_{1,q,K}=SEE_{1,p\lfloor \frac{q}{p}\rfloor,K}=\Spec\Bigl(K\bigl[\alpha^{(1)}_{(pi)}\bigl|i\in\bigl\llbracket0,\bigr\lfloor \frac{q}{p}\bigr\rfloor\bigr\rrbracket\bigr]\Bigr)\cong \mathbb A^{\lfloor \frac{q}{p}\rfloor+1}_K.$$ 

We are left to prove that if $\min(n,q)\ge 2$ and for $n=q=2$ we have $\chr(K)\neq 2$, then $SEE_{n,q,K}$ is not smooth. It suffices to show that there exists $m\in\{2,3\}$ and a $K[t]/(t^m)$-valued point $P_m$ of $SEE_{n,q,K}$ which does not lift to a $K[t]/(t^{m+1})$-valued point of $EE^{\ast}_{n,q,K}$ (hence also of $SEE_{n,q,K}$). 

We consider an $n$-tuple $(f_1,\ldots,f_n)\in\EE_{n,q,K}\bigl(K[t]/(t^{m+1})\bigr)$ whose reduction modulo $\bigl(t)/(t^{m+1}\bigr)$ is $(x_1,\ldots,x_n)$ and whose reduction modulo $\bigl(t^2)/(t^{m+1}\bigr)$ is a point $P_2\in SEE_{n,q,K}\bigl(K[t]/(t^2)\bigr)$.
For $i\in \llbracket1,n\rrbracket$ we write 
$$f_i(x_1,\ldots,x_n)=x_i+\sum_{j=1}^m t^j g_{i,j}(x_1,\ldots,x_n)$$ 
such that for every $j\in\llbracket1,m\rrbracket$ we have $g_{i,j}\in K_{\s}[x_1,\ldots,x_n]$ with $\deg(g_{i,j})\le q$. Let $(h_1,\ldots,h_m)\in K_{\s}[x_1,\ldots,x_s]^m$ be such that
$$\det\bigl(J(f_1,\ldots,f_n)\bigr)=1+\sum_{j=1}^m t^j h_j(x_1,\ldots,x_n);$$ 
as $P_2\in SEE_{n,q,K}(K[t]/(t^2))$, we have $h_1=0$. From Equation (\ref{EQ29}) applied to $\Box=\sum_{j=1}^m t^jJ(g_{1j},\ldots,g_{nj})$ we get identities
\begin{equation}\label{EQ30}
h_1(x_1,\ldots,x_n)=\tr_1\bigl(J(g_{1,1},\ldots,g_{n,1})\bigr)=\sum_{i=1}^n\frac{\partial g_{i,1}}{\partial x_i}(x_1,\ldots,x_n)
\end{equation}
and
\begin{equation}\label{EQ31}
h_2(x_1,\ldots,x_n)=\tr_2\bigl(J(g_{1,1},\ldots,g_{n,1})\bigr)+\sum_{i=1}^n\frac{\partial g_{i,2}}{\partial x_i}(x_1,\ldots,x_n).
 \end{equation}

We show that we can choose $m\in\{2,3\}$ and $g_{i,j}\in K_{\s}[x_1,\ldots,x_n]$ for each pair $(i,j)\in\llbracket1,n\rrbracket\times \llbracket1,m-1\rrbracket$ such that $h_1=\cdots=h_{m-1}=0$ but, no matter what the $n$-tuple $(g_{1,m},\ldots,g_{n,m})\in K_{\s}[x_1,\ldots,x_n]^n$ is, we have $h_m\notin K$: we can take $P_m\in SEE_{n,q,K}\bigl(K[t]/(t^m)\bigr)$ to be the reduction of $(f_1,\ldots,f_n)$ modulo $(t^m)/(t^{m+1})$. 

We first consider three non-disjoint cases with $m=2$ and $g_{i,1}=0$ if $i\in\llbracket3,n\rrbracket$.

{\bf Case 1: $\chr(K)\nmid q$.} We take $g_{1,1}=-q^{-1}x_1^q$ and $g_{2,1}=x_1^{q-1}x_2$. Equations (\ref{EQ30}) and (\ref{EQ31}) give that indeed $h_1=0$ but 
$h_2=-x_1^{2q-2}+\sum_{i=1}^n\frac{\partial g_{i,2}}{\partial x_i}(x_1,\ldots,x_n)\notin K$
(each $\frac{\partial g_{i,2}}{\partial x_i}$ has degree at most $q-1$ and thus less than $2q-2$ as $q\ge 2$); with $m=2$, $P_2$ has the desired property.

{\bf Case 2: $n\ge 3$.} By taking $g_{1,1}=-x_1x_3^{q-1}$ and $g_{2,1}=x_2x_3^{q-1}$, this case is similar to the previous one: $h_2=-x_3^{2q-2}+\sum_{i=1}^n\frac{\partial g_{i,2}}{\partial x_i}(x_1,\ldots,x_n)\notin K$.

{\bf Case 3: $\chr(K)\nmid q-1$ and $q\ge 4$.} By taking $g_{1,1}=-(q-1)^{-1}x_1^{q-1}$ and $g_{2,1}=x_1^{q-2}x_2$, we have relations $h_2=-x_3^{2q-3}+\sum_{i=1}^n\frac{\partial g_{i,2}}{\partial x_i}(x_1,\ldots,x_n)$. Thus, as for $q\ge 4$ we have $2q-3>q-1$, we get that $h_2\notin K$.

The only case not covered by the first three ones is the following one. 

{\bf Case 4: $n=2$ and $q=\chr(K)=3$.} We take $m=3$, $g_{1,1}=x_1^2$, $g_{1,2}=0$, $g_{2,1}=x_1x_2$, and $g_{2,2}=x_1^2x_2$. As
$$\det\bigl(J(x_1+t x_1^2,x_2+t x_1x_2+t^2x_1^2x_2)\bigr)=1+3t x_1+3t^2 x_1^2+2t^3 x_1^3=1+2t^3 x_1^3,$$ $h_1=h_2=0$ but $h_3=2x_1^3+\frac{\partial g_{1,3}}{\partial x_1}+\frac{\partial g_{2,3}}{\partial x_2}\notin K$ as $\deg(g_{i,3})\le 3$ for $i\in \{1,2\}$.

From the four cases above we get that there exists $m\in\{2,3\}$ and a $K[t]/(t^m)$-valued point $P_m$ with the desired properties. Hence the theorem holds.\end{proof} 

\section{On cases when Conjecture \ref{CJ3} holds}\label{S27}

Artin--Schreier equations define connected finite Galois covers of $\mathbb A^n_k$ of degree $p$. We show that \cite{Mi2}, Thms.\ 1.1 and 1.2 hold for affine schemes over $\Spec\mathbb F_p$.

\begin{lemma}\label{F6}
Let $W$ be a connected $\mathbb F_p$-scheme such that $H^1(W,\mathcal O_W)=0$ (e.g., $W$ is affine). Then we have a functorial isomorphism of groups
\begin{equation}\label{EQ37.1}
\textup{cob}_W:\mathcal O(W)/\bigl(x^p-x|x\in\mathcal O(W)\bigr)\cong H^1_{\et}(W,\mathbb Z/p\mathbb Z)
\end{equation} 
which for $f\in\mathcal O(W)$ maps $f+\bigl(x^p-x|x\in\mathcal O(W)\bigr)$ to the $\mathbb Z/p\mathbb Z$-torsor $V$ over $W$ defined by the homomorphism $\mathbb Z/p\mathbb Z\rightarrow \Aut_{\mathcal O_W-\textup{algebra}}\bigl(\mathcal O_W[x]/(x^p-x-f)\bigr)$ that maps $1+p\mathbb Z$ to the $\mathcal O_W$-algebra automorphism of $\mathcal O_W[x]/(x^p-x-f)$ that maps $x+(x^p-x-f)$ to $x+1+(x^p-x-f)$.
\end{lemma} 

\begin{proof}
The first cohomology coboundary of the short exact sequence 
$$0\rightarrow (\mathbb Z/p\mathbb Z)_W\rightarrow \mathbb G_{\a,W}\xrightarrow{x^p-x} \mathbb G_{\a,W}\rightarrow 0$$ 
in the \'etale topos of $W$ is the functorial isomorphism of Equation (\ref{EQ37.1}). The identification of $H^1_{\et}(W,\mathbb Z/p\mathbb Z)$ with the non-abelian $H^1(W,\mathbb Z/p\mathbb Z)$, defined as the set of isomorphism classes of $(\mathbb Z/p\mathbb Z)_W$-torsors as in \cite{Gi}, Ch.\ III, Rmk.\ 3.5.4, gives that $V$ is the pullback 
\[\xymatrix{
V \ar[r]\ar[d] & W \ar[d]^{\s_f} \\
\mathbb G_{\a,W} \ar[r]^{x^p-x} & \mathbb G_{\a,W}
}\]
via the section $\s_f$ that defines $f$, from which the lemma follows.\end{proof}

We also include a proof of the following likely known result. 

\begin{proposition}\label{PR15} 
There exists no connected finite \'etale cover $X$ of $\mathbb A^n_k$ which is either of degree $m\in\llbracket2,p-1\rrbracket$ or Galois of degree at least $2$ and prime to $p$.
\end{proposition}

\begin{proof} We show that the assumption that there exists an $n\in\mathbb N^\ast$ such that $\mathbb A^n_k$ has a connected finite \'etale cover $X$ of the type mentioned leads to a contradiction. We choose the smallest such $n$. If $n=1$, the connected finite Galois cover $X^+$ of $\mathbb A^1_k$ generated by $X$ has a degree $m^+>1$ which either divides $m!$ and hence $(p-1)!$ or is prime to $p$; so $X^+\rightarrow\mathbb A^1_k$ is a non-trivial tamely ramified finite Galois cover of the curve $\mathbb A^1_k$ and this contradicts \cite{SGA1}, Exp.\ XIII, Cor.\ 2.12. If $n\ge 2$, from Bertini Theorem (see \cite{Jo1}, Part I, Thm.\ 6.3(4)) we get that there exists a hyperplane $\mathbb H\subset \mathbb A^n_k$ such that $X\times_{\mathbb A^n_k} \mathbb H$ is a connected finite \'etale cover of $\mathbb H\cong\mathbb A^{n-1}_k$ of the type mentioned, and this contradicts the minimality assumption.\end{proof}

\begin{corollary}\label{C10}
For $e\in\EE_n(k)$ the following properties hold.

\medskip
{\bf (1)} If $\deg(e)<p$ and $e$ is finite, then $e\in\GA_n(k)$.

\smallskip
{\bf (2)} If the field extension $k_{\t}(x_1,\ldots,x_n)\rightarrow k_{\s}(x_1,\ldots,x_n)$ induced by $e^{\#}$ is Galois, then the finite morphism $X_e\rightarrow\mathbb A^n_{K,\t}$ is Galois with $\deg(e)\in\{1\}\cup p\mathbb N^{\ast}$ and hence $\psi_e$ is \'etale. In particular, if $p\nmid\deg(e)$, then $\deg(e)=1$ and $e\in\GA_n(k)$. 
\end{corollary}

\begin{proof} 
For part (1) we have $\deg(e)=1$ by Proposition \ref{PR15}. For both parts, if $\deg(e)=1$, then $e\in\GA_n(K)$ by Lemma \ref{L6}(4). So part (1) holds.

For part (2), if $O_{\t}$ is a local ring of $\mathbb A^n_{k,\t}$ which is a discrete valuation ring, from Lemma \ref{L6}(2) we get that there exists a local ring $O_{\s}$ of $\mathbb A^n_{k,\s}$ which dominates $O_{\t}$ and hence it is a discrete valuation ring. As $\mathbb A^n_{k,\s}$ is an open subvariety of $X_e$ by Zariski's Main Theorem, $O_{\s}$ is a local ring of $X_e$ which is an \'etale $O_{\t}$-algebra. From this, as the finite field extension $k_{\t}(x_1,\ldots,x_n)\rightarrow k_{\s}(x_1,\ldots,x_n)$ is Galois, we get that the morphism $X_e\rightarrow\mathbb A^n_{k,\t}$ is \'etale above all points of $\mathbb A^n_{k,\t}$ of dimension at most $1$ and hence, based on the purity of the branch locus (see \cite{Gro6}, Exp.\ X, Thm.\ 3.4 (i)), it is \'etale. So $X_e\rightarrow\mathbb A^n_{k,\t}$ is Galois and from Proposition \ref{PR15} applied to $X_e$ we get that $\deg(e)\in\{1\}\cup p\mathbb N^{\ast}$. So part (2) holds.\end{proof} 

It is a standard piece of algebraic geometry to show that the classical results of Wright (see \cite{Wr1}, Thms.\ 3.2, 3.3 and 3.7) and Razar (see \cite{Raz}, Thm.\ 2) in characteristic $0$ (for $K=\mathbb C$, see also \cite{C}, Thm.) follow also directly from Corollary \ref{C10}(1) and (2) (respectively). Even more directly, the proof of Corollary \ref{C10} applies over $K$ as well leading to new proofs of these classical results as well.

Directly from Corollaries \ref{C9+}(1) and \ref{C10}(1) we get the following result.

\begin{corollary}\label{C10+}
If $q\in\mathbb N^{\ast}$ is such that $q^n<p$, then $\deg_{p,n,q}^{\f}=1$.
\end{corollary}

Next we generalize the part of Example \ref{EX32} that pertains to Conjecture \ref{CJ3}.

\begin{proposition}\label{PR16}
Let $e\in\EE_n(K)$ be of weak type $1$ in the sense of Definition \ref{D11}. Then the following properties hold.

\medskip
{\bf (1)} If $\chr(K)=0$, then $e\in\GA_n(K)$. 

\smallskip
{\bf (2)} If $\chr(K)=p$, then Conjecture \ref{CJ3} holds for $e$.
\end{proposition}

\begin{proof}
Let $\pi:\mathbb A^n_{K,\s}\rightarrow\mathbb A^1_K$ be a morphism such that the condition for weak type $\le 1$ of Definition \ref{D11} holds for $\pi\times e:\mathbb A^n_{K,\s}\rightarrow\mathbb A^1_K\times_{\Spec K} \mathbb A^n_{K,\t}$. 

The Zariski closure $Z$ of $\Imm(\pi\times e)$ in $\mathbb A^1_K\times_{\Spec K} \mathbb A^n_{K,\t}$ is an irreducible hypersurface in $\mathbb A^1_{K,\s}\times_{\Spec K} \mathbb A^n_K=\mathbb A^{n+1}_K$ and hence the zero locus of an irreducible polynomial $f(x_0,\ldots,x_n)\in K[x_0,\ldots,x_n]$ whose degree in $x_0$ is $\deg(e)$. 

Let $(g_1,\ldots,g_n)\in K_{\s}[x_1,\ldots,x_n]^n$ define $e$ and let $g_0\in K_{\s}[x_1,\ldots,x_n]$ define $\pi$. For all $(\alpha_1,\ldots,\alpha_n)\in K^n$ we have $f\bigl(g_0(\alpha_1,\ldots,\alpha_n),\ldots,g_n(\alpha_1,\ldots,\alpha_n)\bigr)=0$. Taking partial derivatives we get that for each $i\in \llbracket1,n\rrbracket$ we have
$$\sum_{j=0}^n f_{x_j}\bigl(g_0(\alpha_1,\ldots,\alpha_n),\ldots,g_n(\alpha_1,\ldots,\alpha_n)\bigr)g_{j,x_i}(\alpha_1,\ldots,\alpha_n)=0.$$

As the determinant of the Jacobian matrix $J\bigl(g_1(x_1,\ldots,x_n),\ldots,g_n(x_1,\ldots,x_n)\bigr)$ is in $K^{\ast}$, both systems of equations $f=f_{x_0}=\cdots=f_{x_n}=0$ and $f=f_{x_0}=0$ along $\Imm(\pi\times e)$ define $\Sing(Z)\cap\Imm(\pi\times e)$.\footnote{If $\pi$ is smooth, then the system of equations $f=f_{x_1}=\cdots=f_{x_n}=0$ along $\Imm(\pi\times e)$ defines $\Sing(Z)\cap\Imm(\pi\times e)$ as well.} Thus the zero locus of the image of $f_{x_0}\bigl(g_0(\alpha_1,\ldots,\alpha_n),\ldots,g_n(\alpha_1,\ldots,\alpha_n)\bigr)$ in $K_{\s}[x_1,\ldots,x_n]$ via the $K$-algebra homomorphism $K[x_0,\ldots,x_n]\rightarrow K[x_1,\ldots,x_n]$ that maps $x_i$ to $g_i$ for each $i\in \llbracket0,n\rrbracket$ is an effective divisor $\sum_{i=1}^s q_iY_i$ of $\mathbb A^n_{K,\s}$, where $s\in\mathbb N$, $(q_1,\ldots,q_s)\in (\mathbb N^{\ast})^s$, and $Y_1,\ldots,Y_s$ are distinct irreducible divisors, such that $e^{-1}\bigl(e(Y_i)\bigr)=Y_i$ for each $i\in \llbracket1,s\rrbracket$ by Definition \ref{D1}. Hence for each $i\in \llbracket1,s\rrbracket$ there exists an irreducible polynomial $h_i\in K_{\t}[x_1,\ldots,x_n]$ such that the zero locus $e^{\#}(h_i)=0$ is $Y_i$. It follows that there exists $\beta\in K^{\ast}$ such that $f_{x_0}-\beta\prod_{i=1}^s h_i^{q_i}$ is divisible by $f$ and by reasons of degrees in $x_0$ it follows that $f_{x_0}=\beta\prod_{i=1}^s h_i^{q_i}$ has degree $0$ in $x_0$. 

Thus $\deg(e)=1$ if $\chr(K)=0$ and $\deg(e)\in\{1\}\cup p\mathbb N^{\ast}$ if $\chr(K)=p$. From this and Lemma \ref{L6}(4) we get that the proposition holds.\end{proof}

\begin{corollary}\label{C11}
Let $e\in\EE_n(k)$. Then the following properties hold.

\medskip
{\bf (1)} If the $k_{\t}[x_1,\ldots,x_n]$-algebra $A_e$ is generated by $1$ element, then $\psi_e$ is \'etale and Conjecture \ref{CJ3} holds for $e$.

\smallskip
{\bf (2)} If $A_e$ does not have a $p$-basis, then the $k_{\t}[x_1,\ldots,x_n]$-algebra $A_e$ is not generated by $1$ element.

\smallskip
{\bf (3)} If $A_e$ is regular and does not have a $p$-basis, then $\ciedim(X_e)=\infty$.
\end{corollary}

\begin{proof}
To prove part (1), let $\pi:\mathbb A^n_{k,\s}\rightarrow\mathbb A^1_k$ be a morphism defined by an element in $A_e\subset k_{\s}[x_1,\ldots,x_n]$ that generates the $k_{\t}[x_1,\ldots,x_n]$-algebra $A_e$. The Zariski closure of $\Imm(\pi\times e)$ in $\mathbb A^1_k\times_{\Spec k} \mathbb A^n_{k,\t}$ is $X_e$ itself; so $e$ is of weak type $\le 1$ and Conjecture \ref{CJ3} holds for $e$ by Proposition \ref{PR16}. With its notation for $k$ instead of $K$, the mentioned proof gives that $f_{x_0}=\beta\in k^{\ast}$, hence $\psi_e$ is \'etale.

If $A_e$ does not have a $p$-basis, then $\psi_e$ is non-\'etale and hence part (2) follows from part (1) and part (3) follows from Theorem \ref{T3+}(6).\end{proof}

\section{Counterexamples to Conjecture \ref{CJ3}}\label{S28}

To provide counterexamples to Conjecture \ref{CJ3} in all dimensions $n\ge 2$ or $n\ge 2$, it suffices to provide counterexamples to it in dimension $3$ or $2$. Thus in this section we take $n\in\{2,3\}$. 

As in this section the main positive integers used as parameters are divisible by $p$, to emphasize this the factor $p$ is inserted in the notation.

\begin{theorem}\label{T6.3}
Let $(l,m,s,q)\in (\mathbb N^{\ast})^4$, $f(t):=\frac{1-(1-t)^{pl-1}}{t}\in k[t]$, and $\delta\in k$. The rule 
$$(x,y,z)\mapsto \bigl(x(1-xy),yf(xy)+x^s(1-xy)^sz,z+\delta x^{pq}z^{pm}\bigr)$$ on valued points defines an endomorphism $e\in\End_3(k)$. Then the following properties hold.

\medskip
{\bf (1)} We have $\e\in\EE_3(k)$.

\smallskip
{\bf (2)} Suppose that $q=l(pm-1)$, $s=pl-1$, and $\delta=-1$. Then the following properties hold.

\medskip\noindent
{\bf (2.a)} We have $\deg(e)=p^2lm-pl+1\equiv 1 \pmod{p}$. 

\smallskip\noindent
{\bf (2.b)} The endomorphism $e$ is surjective.
\end{theorem}

\begin{proof}
We have $t(1-t)f(t)=(1-t)-(1-t)^{pl}\in -t+t^pk[t]$. Hence the endomorphism $d_0:=e^{1,1}_{1-t,f(t);1}=\e\bigl(x(1-xy),yf(xy)\bigr)\in\EE_2(k)$ of Example \ref{EX23-} applied to the sextuple $(n,f_1,f_2,L,l,\alpha)=\bigl(2,1-t,f(t),1,1,-1\bigr)$ has Jacobian determinant $-1$. Considering $d_1:=d_0\times 1_{\mathbb A^1_k}\in\EE_3(k)$ and the special tame automorphism $a:=\e(x_1,x_2+x_1^sx_3,x_3)$, the composite $d_2:=ad_1\in\EE_3(k)$ is defined on valued points by the rule $(x,y,z)\mapsto \bigl(x(1-xy),yf(xy)+x^s(1-xy)^sz,z\bigr)$ and has Jacobian determinant $-1$. As $e$ is a perturbation of $d_2$, its Jacobian determinant is also $-1$. So part (1) holds.

For part (2), to solve the system of equations
$$x(1-xy)-x_1=yf(xy)+x^s(1-xy)^sz-x_2=z+\delta x^{pq}z^{pm}-x_3=0,$$
we first note that from the equation $x-x^2y-x_1=0$ we get that $y=\frac{x-x_1}{x^2}$, that $xy=\frac{x-x_1}{x}=1-\frac{x_1}{x}$, and that the equation $yf(xy)+x^s(1-xy)^sz-x_2=0$ can be rewritten as $yf(xy)+x_1^sz=x_2$. Hence 
$$z=\frac{x_2-yf(xy)}{x_1^s}=\frac{x_2-\frac{x-x_1}{x^2}f(1-\frac{x_1}{x})}{x_1^s}=\frac{x_2x-\bigl(1-\frac{x_1}{x}\bigr)f(1-\frac{x_1}{x})}{x_1^sx}.$$
As $-(1-t)f(1-t)=t^{pl-1}-1$, we get that
$$z=\frac{x_2x+\bigl(\frac{x_1}{x}\bigr)^{pl-1}-1}{x_1^sx}=\frac{x_2x^{pl}-x^{pl-1}+x_1^{pl-1}}{x_1^sx^{pl}}.$$ 

Substituting this expression of $z$ in the equation $z+\delta x^{pq}z^{pm}-x_3=0$ we get that
$$\frac{x_2x^{pl}-x^{pl-1}+x_1^{pl-1}}{x_1^sx^{pl}}+\delta x^{pq}\frac{(x_2x^{pl}-x^{pl-1}+x_1^{pl-1})^{pm}}{x_1^{psm}x^{p^2lm}}-x_3=0.$$
Therefore, as $q=l(pm-1)$, we get that
$$\frac{x_2x^{pl}-x^{pl-1}+x_1^{pl-1}}{x_1^s}+\delta\frac{(x_2x^{pl}-x^{pl-1}+x_1^{pl-1})^{pm}}{x_1^{psm}} -x_3x^{pl}=0.$$
As $s=pl-1$ and by substituting $\delta=-1$, the rational function
$$L(x_1,x_2,x_3,w):=\frac{\frac{(x_2w^{pl}-w^{pl-1}+x_1^{pl-1})^{pm}}{x_1^{psm}}-\frac{x_2w^{pl}-w^{pl-1}+x_1^{pl-1}}{x_1^s}+x_3w^{pl}}{w^{pl-1}}$$ 
is a polynomial in $k_{\t}[x_1,x_2,x_3]\Bigl[\frac{1}{x_1}\Bigr][w]$, so in $w$ with coefficients in $k_{\t}[x_1,x_2,x_3]\Bigl[\frac{1}{x_1}\Bigr]$. Thus $L$ has degree $p^2lm-pl+1$ and leading coefficient $\bigl(\frac{x_2}{x_1^s}\bigr)^{pm}$. Moreover, we have $L(x_1,x_2,x_3,0)=\frac{1}{x_1^s}$.

So the field extension $k_{\s}(x,y,z)$ of $k_{\t}(x_1,x_2,x_3)$ is generated by $x$ and we have $L(x_1,x_2,x_3,x)=0$. Thus $\deg(e)\le p^2lm-1$. 

If $\mathcal B$ is the solution set in $k$ of the equation $w^{pm}=w$, then $1\in\mathcal B$, and we have identities $|\mathcal B|=pm$ and
$$e^{-1}(1,1,0)=\Bigr\{\Bigl(\alpha,\frac{\alpha-1}{\alpha^2},\frac{\alpha^{pl}-\alpha^{pl-1}+1}{\alpha^{pl}}\Bigr)\Bigl|\alpha\in k^{\ast},\alpha^{pl}-\alpha^{pl-1}+1\in\mathcal B\Bigr\}.$$ 
For $\beta\in\mathcal B$, the equation $w^{pl}-w^{pl-1}+1=\beta$ in $w$ has $pl$ distinct solutions in $k^{\ast}$ if $\beta\neq 1$ and has $1$ solution in $k^{\ast}$ if $\beta=1$. Thus $|e^{-1}(1,1,0)|=p^2lm-pl+1$.

Thus $\deg(e)\ge p^2lm-pl+1$. Hence $\deg(e)=p^2lm-pl+1$. So part (2.a) holds.

For part (2.b), let $(\gamma_1,\gamma_2,\gamma_3)\in k^3$. If $\gamma_1=0$, then $\bigl(0,-\gamma_2,\gamma_3\bigr)\in e^{-1}(\gamma_1,\gamma_2,\gamma_3)$ as $qs>0$ and $f(0)=-1$. 

If $\gamma_1\neq 0$, then $L(\gamma_1,\gamma_2,\gamma_3,w)\in k[w]$ is a polynomial of degree $p^2lm-pl+1$ if $\gamma_2\neq 0$ and of degree $(pl-1)(pm-1)$ if $\gamma_2=0$. Hence, as $L(\gamma_1,\gamma_2,\gamma_3,0)=\frac{1}{\gamma_1^s}\neq 0$, it has a solution $\alpha$ in $k^{\ast}$ and therefore $\Bigl(\alpha,\frac{\alpha-\gamma_1}{\alpha^2},\frac{\gamma_2\alpha^{pl}-\alpha^{pl-1}+\gamma_1^{pl-1}}{\gamma_1^s\alpha^{pl}}\Bigr)\in e^{-1}(\gamma_1,\gamma_2,\gamma_3)$. 

From the last two paragraphs we get that part (2.b) holds.\end{proof}

The next two theorems are `$-1\equiv \pmod{p}$' versions of Theorem \ref{T6.3} that provide both non-surjective and surjective \'etale endomorphisms in $\EE_3(k)$ but the geometric degrees are congruent to $-1$ modulo $p^2$.

\begin{theorem}\label{T6.4}
Let $(l,m,q)\in (\mathbb N^{\ast})^3$. Let the set $\mathbb I_{p,1,pl}^{\ast}=\{\beta\in k^{\ast}|\beta^{p^l}=\beta\}$ be as in Example \ref{EX11}. The rule on valued points
$$(x,y,z)\mapsto \bigl(x-x^{pl}y,y+(1-x^{pl-1}y)^{pl}+z\bigl[x-x^{pl}y-(x-x^{pl}y)^{pl}\bigr],z+\delta x^{pl}z^{pm}\bigr)$$ defines an endomorphism $e\in\End_3(k)$. Then the following properties hold.

\medskip
{\bf (1)} We have $\e\in\EE_3(k)$ and $e$ is non-surjective with $\mathbb I_{p,1,pl}^{\ast}\times\mathbb A^1_k\subset \Gamma_e$.

\smallskip
{\bf (2)} If $\delta=-1$ and $q=plm-l$, then $\deg(e)=p^2lm-1$ and $\Gamma_e=\mathbb I_{p,1,pl}^{\ast}\times\mathbb A^1_k$ is a complete intersection curve in $\mathbb A^3_k$ with $pl-1$ connected components isomorphic to $\mathbb A^1_k$.
\end{theorem}

\begin{proof}
From Example \ref{EX11} applied to $(x,y,m,q)=(y,x,1,pl)$, we get that the rule $(x,y)\mapsto \bigl(x-x^{pq}y,y+(1-x^{pl-1}y)^{pl})$ on valued points defines a non-surjective $d_0\in\EE_2(k)$. For $d_1:=d_0\times 1_{\mathbb A^1_k}\in\EE_3(k)$ and the special tame automorphism $a:=\e(x_1,x_2+(x_1-x_1^{pl})x_3,x_3)\in\GA_3(k)$, the composite $d_2:=ad_1\in\EE_3(k)$ is defined on valued points by the rule 
$$(x,y,z)\mapsto \bigl(x-x^{pl}y,y+(1-x^{pl-1}y)^{pl}+z\bigl[x-x^{pl}y-(x-x^{pl}y)^{pl}\bigr],z\bigr).$$ As $e$ is a perturbation of $d_0$, the Jacobian matrices of $e$ and $d_2$ are equal. Hence $e\in\EE_3(k)$. 

From Equation (\ref{EQ5.1}) with the roles of $x$ and $y$ interchanged, we get the identity $\Gamma_{d_0}=\{(\beta,0)|\beta\in\mathbb I_{p,1,pl}^{\ast}\}$. From this, as the term $z[x-x^{pl}y-(x-x^{pl}y)^{pl}]$ is $0$ when $x-x^{pl}y\in\mathbb I_{p,1,pl}^{\ast}$, we get that $\mathbb I_{p,1,pl}^{\ast}\times\mathbb A^1_k\subset \Gamma_e$. So part (1) holds.

For part (2), to solve the system of equations
$$x-x^{pl}y-x_1=y+(1-x^{pl-1}y)^{pl}+z[x-x^{pl}y-(x-x^{pl}y)^{pl}]-x_2=z+\delta x^{pl}z^{pm}-x_3=0,$$
we note that the equation $x-x^{pl}y-x_1=0$ gives that $y=\frac{x-x_1}{x^{pl}}$, that 
$$1-x^{pl-1}y=1-\frac{x-x_1}{x}=\frac{x_1}{x},$$ and that the equation $y+(1-x^{pl-1}y)^{pl}+z[x-x^{pl}y-(x-x^{pl}y)^{pl}]-x_2=0$ can be rewritten as $y+(1-x^{pl-1}y)^{pl}+(x_1-x_1^{pl})z=x_2$. Hence 
$$z=\frac{x_2-y-(1-x^{pl-1}y)^{pl}}{x_1-x_1^{pl}}=\frac{x_2-\frac{x-x_1}{x^{pl}}-\frac{x_1^{pl}}{x^{pl}}}{x_1-x_1^{pl}}=\frac{x_2x^{pl}-x+x_1-x_1^{p^l}}{(x_1-x_1^{pl})x^{pl}}.$$
Substituting this expression of $z$ in the equation $z+\delta x^{pq}z^{pm}-x_3=0$ we get that
$$\frac{x_2x^{pl}-x+x_1-x_1^{p^l}}{(x_1-x_1^{pl})x^{pl}}+\delta x^{pq}\frac{(x_2x^{pl}-x+x_1-x_1^{p^l})^{pm}}{(x_1-x_1^{pl})^{pm}x^{p^2lm}}-x_3=0.$$

Therefore, as $q=plm-l$, we get that
$$\frac{x_2x^{pl}-x+x_1-x_1^{p^l}}{x_1-x_1^{pl}}+\delta\frac{(x_2x^{pl}-x+x_1-x_1^{p^l})^{pm}}{(x_1-x_1^{pl})^{pm}} -x_3x^{pl}=0.$$
Substituting $\delta=-1$ and multiplying by $-1$, the rational function
$$L(x_1,x_2,x_3,w):=\frac{\frac{(x_2w^{pl}-w+x_1-x_1^{p^l})^{pm}}{(x_1-x_1^{pl})^{pm}}-\frac{x_2w^{pl}-w+x_1-x_1^{p^l}}{x_1-x_1^{pl}}+x_3w^{pl}}{w}$$ 
is a polynomial in $k_{\t}[x_1,x_2,x_3]\Bigl[\frac{1}{x_1}\Bigr][w]$, so in $w$ with coefficients in $k_{\t}[x_1,x_2,x_3]\Bigl[\frac{1}{x_1}\Bigr]$. Thus $L$ has degree $p^2lm-1$ and leading coefficient $\bigl(\frac{x_2}{x_1-x_1^{pl}}\bigr)^{pm}$. Moreover, we have $L(x_1,x_2,x_3,0)=\frac{1}{x_1-x_1^{pl}}$.

So the field extension $k_{\s}(x,y,z)$ of $k_{\t}(x_1,x_2,x_3)$ is generated by $x$ and we have $L(x_1,x_2,x_3,x)=0$. Thus $\deg(e)\le p^2lm-1$. 

Let $\mathbb I_{p,1,pl}=\mathbb I_{p,1,pl}^{\ast}\cup\{0\}$ be as in Example \ref{EX11}. To each $\gamma_1\in k\setminus\mathbb I_{p,1,pl}$, we associate $\gamma_4:=\gamma_1-\gamma_1^{pl}\in k^{\ast}$.

Let $\mathcal B$ be the solution set in $k$ of the equation $w^{pm}=w$; so $1\in\mathcal B$ and $|\mathcal B|=pm$. 

If $\gamma_1\in k\setminus\mathbb I_{p,1,pl}$, then
$$e^{-1}(\gamma_1,1,0)=\Bigr\{\Bigl(\alpha,\frac{\alpha-\gamma_1}{\alpha^{pl}},\frac{\alpha^{pl}-\alpha+\gamma_4}{\gamma_4\alpha^{pl}}\Bigr)\Bigl|\alpha\in k^{\ast},\frac{\alpha^{pl}-\alpha+\gamma_4}{\gamma_4}\in\mathcal B\Bigr\};$$ 
so $|e^{-1}(\gamma_1,1,0)|=p^2lm-1$ as for $\beta\in\mathcal B$, the equation $\frac{w^{pl}-w+\gamma_4}{\gamma_4}=\beta$ in $w$ has $pl$ distinct solutions in $k^{\ast}$ if $\beta\neq 1$ and has $pl-1$ distinct solutions in $k^{\ast}$ if $\beta=1$.

Thus $\deg(e)\ge p^2lm-1$. Hence $\deg(e)=p^2lm-1$. 

Based on part (1), to check that $\Gamma_e=\mathbb I_{p,1,pl}^{\ast}\times\mathbb A^1_k$ it suffices to show that for each $(\gamma_1,\gamma_2,\gamma_3)\in (k\setminus\mathbb I_{p,1,pl}^{\ast})\times k^2$, we have $e^{-1}(\gamma_1,\gamma_2,\gamma_3)\neq\emptyset$. For $\gamma_1=0$ we have $(0,\gamma_2-1,\gamma_3)\in e^{-1}(0,\gamma_2,\gamma_3)$. 

Assume now that $\gamma_1\in k\setminus \mathbb I_{p,1,pl}$. The polynomial $L(\gamma_1,\gamma_2,\gamma_3,w)\in k[w]$ has degree $p^2lm-1$ if $\gamma_2\neq 0$ and degree $(pl-1)pm-1$ if $\gamma_2=0$. From this and $L(\gamma_1,\gamma_2,\gamma_3,0)=\frac{1}{\gamma_4}$ we get that $L(\gamma_1,\gamma_2,\gamma_3,w)$ has a root $\alpha\in k^{\ast}$ and therefore $\Bigl(\alpha,\frac{\alpha-\gamma_1}{\alpha^{pl}},\frac{\alpha^{pl}-\alpha+\gamma_4}{\gamma_4\alpha^{pl}}\Bigr)\in e^{-1}(\gamma_1,\gamma_2,\gamma_3)$. 

We conclude that $\Gamma_e=\mathbb I_{p,1,pl}^{\ast}\times\mathbb A^1_k$ and that part (2) holds.\end{proof}

\begin{corollary}\label{C8.2}
Suppose that $p=2$ and $n\in\mathbb N^{\ast}\setminus\{1,2\}$. Then the following properties hold. 

\medskip
{\bf (1)} For each $j\in \mathbb N^{\ast}$ there exists a surjective $e\in\EE_n(k)$ with $\deg(e)=j$.

\smallskip
{\bf (2)} Let $j\in\mathbb N^{\ast}\setminus\{1\}$ be such that it is not divisible by a prime number $l$ congruent to $1$ modulo $4$. Then there exists a non-surjective $e\in\EE_n(k)$ with $\deg(e)=j$.
\end{corollary}

\begin{proof}
If $e\in\EE_3(k)$ is surjective (resp.\ non-surjective) with $\deg(e)=j$, then the product $d:=e\times 1_{\mathbb A^{n-3}_k}\in\EE_n(k)$ is surjective (resp.\ non-surjective) with $\deg(d)=j$. So we can assume that $n=3$. Using composites, it suffices to check this when either $j$ is a prime or, for part (1), $j=1$.

For part (1), we can take $e=1_{\mathbb A^3_k}$ if $j=1$ and $e=\e(x_1+x_1^2,x_2,x_3)$ if $j=2$. If $j$ is an odd prime, then the existence of $e$ follows from Theorem \ref{T6.3}(2.a) and (2.b) applied to $p=2$, $m=1$, and $l=\lfloor\frac{j}{2}\rfloor$. So part (1) holds.

For part (2), if $j=2$ then we can take $e=d\times 1_{\mathbb A^1_k}$ where the non-surjective $d\in\EE_2(k)$ with $\deg(d)=2$ is as in Theorem \ref{T1}(3) applied to $p=2$ and $mq=2$. If $j$ is an odd prime, then there exists $l\in\mathbb N^{\ast}$ such that $j=4l-1$ by our hypotheses; so the existence of $e$ follows from Theorem \ref{T6.4}(2) applied to $p=2$, $m=1$, $q=l$, and $\delta=-1$. So part (2) holds.
\end{proof}

\begin{theorem}\label{T6.5}
Let $(q,m)\in (\mathbb N^{\ast})^2$. The rule 
$$(x,y,z)\mapsto\bigl(x+x^{pm}y,y+xz+x^{pm}yz,z+(-1)^q x^{p^2mq-pm}z^{pq}\bigr)$$ on valued points defines an endomorphism $e\in\End_3(k)$. Then the following properties hold.

\medskip
{\bf (1)} We have $\e\in\EE_3(k)$.

\smallskip
{\bf (2)} We have $\deg(e)=p^2mq-1$.

\smallskip
{\bf (3)} The endomorphism $e$ is non-surjective iff there exists $l\in\mathbb N$ such that we have $q=m=p^l$. In such a case, $\Gamma_e$ is the zero locus $x_2=x_1^{p^{l+1}}x_3+(-1)^{q+1}=0$; in particular, $\Gamma_e$ is a complete intersection curve in $\mathbb A^3_k$ isomorphic to $\mathbb A^1_k\setminus\{0\}$.
\end{theorem}

\begin{proof}
The Jacobian determinant of $e$ is $1$. Thus $e\in\EE_3(k)$ and part (1) holds.

To solve the system of equations
$$x+x^{pm}y-x_1=y+xz+x^{pm}yz-x_2=z+(-1)^q x^{p^2mq-pm}z^{pq}-x_3=0,$$
we first note that from the equation $x+x^{pm}y-x_1=0$ we get that $y=\frac{x_1-x}{x^{pm}}$ and that the equation $y+xz+x^{pm}yz-x_2=0$ can be rewritten as $y+x_1z=x_2$. Hence $z=\frac{x_2-y}{x_1}$. Substituting $y=\frac{x_1-x}{x^{pm}}$ in $z=\frac{x_2-y}{x_1}$, we get that 
$$z=\frac{x_2-\frac{x_1-x}{x^{pm}}}{x_1}=\frac{x_2x^{pm}+x-x_1}{x_1x^{pm}}.$$ 
Substituting this expression of $z$ in the equation $z+(-1)^q x^{p^2mq-pm}z^{pq}-x_3=0$ we get that
$$\frac{x_2x^{pm}+x-x_1}{x_1x^{pm}}+(-1)^q x^{p^2mq-pm}\frac{(x_2x^{pm}+x-x_1)^{pq}}{(x_1x^{pm})^{pq}}-x_3=0.$$
Therefore, by multiplying with $-x^{pm}$ we get that
$$(-1)^{q+1}\frac{(x_2x^{pm}+x-x_1)^{pq}}{x_1^{pq}}+x_3x^{pm}-\frac{x_2x^{pm}+x-x_1}{x_1}=0.$$
As $(-1)^{q+1}(-1)^{pq}+1=0$, we can simplify by $x$ and we get that the field extension $k_{\s}(x,y,z)$ of $k_{\t}(x_1,x_2,x_3)$ is generated by $x$ and $x$ is a root of a polynomial 
$$L(x_1,x_2,x_3,w):=\frac{(-1)^{q+1}\frac{(x_2w^{pm}+w-x_1)^{pq}}{x_1^{pq}}+x_3w^{pm}-\frac{x_2w^{pm}+w-x_1}{x_1}}{w}$$ 
in $k_{\t}[x_1,x_2,x_3]\Bigl[\frac{1}{x_1}\Bigr][w]$, so in $w$ with coefficients in $k_{\t}[x_1,x_2,x_3]\Bigl[\frac{1}{x_1}\Bigr]$. Note that $L$ has degree $p^2mq-1$ and leading coefficient $(-1)^{q+1}\frac{x_2^{pq}}{x_1^{pq}}$ and satisfies the identity $L(x_1,x_2,x_3,0)=-\frac{1}{x_1}$. Thus $\deg(e)\le p^2mq-1$.

If $\mathcal B$ is the solution set in $k$ of the equation $(-1)^{q+1} w^{pq}=w$, then $-1\in\mathcal B$ and we have identities $|\mathcal B|=pq$ and
$$e^{-1}(1,1,0)=\Bigr\{\Bigl(\alpha,\frac{1-\alpha}{\alpha^{pm}},\frac{\alpha^{pm}+\alpha-1}{\alpha^{pm}}\Bigr)\Bigl|\alpha\in k^{\ast},\alpha^{pm}+\alpha-1\in\mathcal B\Bigr\}.$$ 
For $\beta\in\mathcal B$, the equation $(-1)^{q+1}\alpha^{pm}+\alpha-1=\beta$ has $pm$ distinct solutions in $k^{\ast}$ if $\beta\neq -1$ and has $pm-1$ distinct solutions in $k^{\ast}$ if $\beta=-1$. It follows that $|e^{-1}(1,1,0)|=p^2mq-1$.

Therefore $\deg(e)\ge p^2mq-1$. Hence $\deg(e)=p^2mq-1$. So part (2) holds.

For part (3), let $(\gamma_1,\gamma_2,\gamma_3)\in k^3$ and we want to determine when $e^{-1}(\gamma_1,\gamma_2,\gamma_3)$ is the empty set. If $\gamma_1=0$, then $(0,\gamma_2,\gamma_3)\in e^{-1}(\gamma_1,\gamma_2,\gamma_3)$. 

Assume now that $\gamma_1\neq 0$. Let $f=f_{\gamma_1,\gamma_2,\gamma_3}:=L(\gamma_1,\gamma_2,\gamma_3,w)\in k[w]$. 

If $\gamma_2\neq 0$, then $f$ has degree $p^2mq-1$ and satisfies $f(0)\neq 0$; hence $f$ has a solution $\alpha\in k^{\ast}$ and we get that $\bigl(\alpha,\frac{1-\alpha}{\alpha^{pm}},\frac{\alpha^{pm}+\alpha-1}{\alpha^{pm}}\bigr)\in e^{-1}(\gamma_1,\gamma_2,\gamma_3)$.

Assume now that $\gamma_2=0$. We have
$$f(w)=\frac{\frac{(-1)^{q+1} (w-\gamma_1)^{pq}}{\gamma_1^{pq}}+\gamma_3w^{pm}-\frac{w}{\gamma_1}+1}{w}\in k[w].$$

If either $q\neq m$ or $q=m$ is a not a power of $p$, then $f(w)$ is not a monomial and hence the equation $f(w)=0$ has a non-zero solution $\alpha\in k^{\ast}$ and again we have $\bigl(\alpha,\frac{1-\alpha}{\alpha^{pm}},\frac{\alpha^{pm}+\alpha-1}{\alpha^{pm}}\bigr)\in e^{-1}(\gamma_1,\gamma_2,\gamma_3)$; so $e$ is surjective in this case.

Assume now that $q=m=p^l$ with $l\in\mathbb N$, then $f(w)=w^{p^{l+1}-1}\bigl(\frac{(-1)^{q+1}}{\gamma_1^{pq}}+\gamma_3\bigr)-1$ and thus it has a root $\alpha$ in $k^{\ast}$ iff $\gamma_3\gamma_1^{pq}+(-1)^{q+1}\neq 0$. So for $\gamma_3\gamma_1^{pq}+(-1)^{q+1}\neq 0$ we have $\bigl(\alpha,\frac{1-\alpha}{\alpha^{pm}},\frac{\alpha^{pm}+\alpha-1}{\alpha^{pm}}\bigr)\in e^{-1}(\gamma_1,\gamma_2,\gamma_3)$ and if $\gamma_3\gamma_1^{pq}+(-1)^{q+1}=0$, then we have $e^{-1}(\gamma_1,\gamma_2,\gamma_3)=\emptyset$. So $\Gamma_e$ is the zero locus $x_2=x_1^{pq}x_3+(-1)^{q+1}=0$. Therefore part (3) holds.\end{proof}

\begin{remark}\normalfont\label{R19.5}
We refer to Theorem \ref{T6.5}.

\medskip
{\bf (1)} Suppose that $p=2$ and $m=q=1$. So $\deg(e)=3$ and $e$ is not surjective. Moreover, $e$ is the perturbation of $d_2:=\e(x_1+x_1^2x_2,x_2+x_1x_3+x_1^2x_2x_3,x_3)$ in $\EE_3(k)$ and we have $d_2=ad_1$, with $d_1=d_0\times 1_{\mathbb A^1_k}$ for $d_0:=e^{1,1}_{1+t,1;1}$ as in Proposition \ref{PR8.5} and $a=\e(x_1,x_2+x_1x_3,x_3)\in\GA_3(k)$ a special tame automorphism. The weaker form of Theorem \ref{T6.5} for this special case that involves only parts (1) and (2) was first obtained using computer programming in \cite{H-K}, Thm.\ 1.2.

\smallskip
{\bf (2)} Let $l\in\mathbb N^{\ast}$. Taking $q=1$ and $m=l$, from Theorem \ref{T6.5}(2) and (3) we get that there exists a surjective $e\in\EE_3(k)$ with $\deg(e)=p^2l-1$.

\smallskip
{\bf (3)} Let $l\in\mathbb N^{\ast}$. Taking $q=m=p^{l-1}$, from Theorem \ref{T6.5}(2) and (3) we get that there exists a non-surjective $e\in\EE_3(k)$ with $\deg(e)=p^{2l}-1$.
\end{remark}

\begin{theorem}\label{T6.6}
Let $(l,q,m,r,s)\in (\mathbb N^{\ast})^5$ and $j\in\llbracket0,p-1\rrbracket$. We consider polynomials $h(t,w)\in k[t,w]$, $f(t,w):=h(t^{pl},w)\in k[t^{pl},w]$, and $g(t)\in k[t]$ with $\deg(g)=s$. Let $e\in\End_2(k)$ be defined on valued points by the rule
$$(x,y)\mapsto\bigl(x+f(x,y)+x^{pl}[x+f(x,y)]^{pm},y+g(x^{pl})[x+f(x,y)]^{pq+j}\bigr).$$
Then the following properties hold.

\medskip
{\bf (1)} We have $\e\in\EE_2(k)$.

\smallskip
{\bf (2)} Suppose that $q\ge ms$, $h(x,y)=xy^r$, and $g(t)=-\delta (-t)^s$ with $\delta\in k^{\ast}$. Let $q_1:=q-ms\in\mathbb N$. Then the following properties hold.

\medskip\noindent
{\bf (2.a)} We have $\pi(e)\le\max\bigl(pl+pm(pl+r),pls+(pl+r)(pq+j)\bigr)$.

\smallskip\noindent
{\bf (2.b)} The field extension $k_{\t}(x_1,x_2)\rightarrow k_{\s}(x_1,x_2)$ induced by $e$ at the level of fields of fractions is isomorphic to the field extension 
$$k_{\t}(x_1,x_2)\rightarrow k_{\t}(x_1,x_2)[w]/\bigl(L(w)\bigr),$$ 
where $L(w)=L(x_1,x_2,w)$ in $k_{\t}[x_1,x_2,w]$ is the irreducible polynomial
$$\Bigl[w^{pm+1}+(w-x_1)\bigl(x_2+\delta(w-x_1)^sw^{pq_1+j}\bigr)^r\Bigr]^{pl}+(w-x_1)w^{p^2ml-pm}.$$
In particular, we have 
$$p^2ml-pm+1\le\deg(e)=\deg(L)\le pl\max\bigl(pm+1,r(pq_1+s+j)+1\bigr),$$
where by $\deg(L)$ we mean the degree of $L$ in $w$ alone.

\smallskip\noindent
{\bf (2.c)} The endomorphism $e$ is non-surjective iff $pm=r(pq_1+s+j)$, $m=lrs$, $\delta^r=-1$, and $rs+1$ is a power of $p$. Moreover, if $e$ is non-surjective, then $\varphi(e)=pl(rs+1)-1$ and $\Gamma_e=\{t\in k|t^{pl(rs+1)-1}=1\}\times\{0\}$.

\smallskip\noindent
{\bf (2.d)} Suppose that $pm=r(pq_1+s+j)$, $p\mid rs+1$, and $\delta^r=-1$. Let $i_0\in\mathbb N^{\ast}$ be the largest such that $p^{i_0}\mid rs+1$. If $m>l$, then we also assume that there exists $i\in\llbracket1,i_0\rrbracket$ such that $pq_1+s+j\ge p^i$ and $m\le (p^i-1)l$. Then $\deg(e)=p^2ml-pm+1$. Moreover, the leading coefficient of $L$ is $1$ and $\Spec\Bigl(k[x_1,x_2][w]/\bigl(L(x_1,x_2,w)\bigr)\Bigr)$ is a monomial model of $X_e$.

\smallskip\noindent
{\bf (2.e)} Suppose that $q_1=0$, $j=r=1$, $s=pm-1$, and $\delta=-1$. Let $i_0\in\mathbb N^{\ast}$ be the largest such that $p^{i_0}\mid pm=s+1$. We also assume that $m\le (p^{i_0}-1)l$. Then $\deg(e)=p^2ml-pm+1$. Moreover, $e$ is non-surjective iff $p=2$ and $l=m=s=1$.
\end{theorem}

%\smallskip
%{\bf (5)} Suppose that $p=2$, $m=1$, $q=2$, $r=1$, and $s=2m-1$. Then $X_e\cong \Spec k[x_1,x_2][t]/\bigl(t^3+t(x_1^4+x_1+x_2^2)t+x_1x_2\bigr)$ has one singular point.

\begin{proof}
As $\frac{\partial f}{\partial t}=0$, the rule $(x,y)\mapsto \bigl(x+f(x,y),y\bigr)$ on valued points defines a surjective $d_0\in\EE_2(k)$ with Jacobian determinant $1$.

For the special tame automorphism $a:=\e(x_1,x_2+x_1^{pq+j})\in\GA_2(k)$, the composite $d_1:=ad_0\in\EE_2(k)$ is defined on valued points by the rule 
$$(x,y)\mapsto\bigl(x+f(x,y),y+[x+f(x,y)]^{pq+j}\bigr)$$ 
and has Jacobian determinant $1$. Thus the rule 
$$(x,y)\mapsto\bigl(x+f(x,y),y+g(x^{pl})[x+f(x,y)]^{pq+j}\bigr)$$
on valued points defines $d_2\in\EE_2(k)$ which has Jacobian determinant $1$ and which is a multiplicative perturbation of $d_0$.
As $e$ is a perturbation of $d_2$, the Jacobian matrices of $e$ and $d_2$ are equal. So the Jacobian determinant of $e$ is $1$, hence $e\in\EE_2(k)$. So part (1) holds.

To solve the system of equations
$$x+f(x,y)+x^{pl}[x+f(x,y)]^{pm}-x_1=0=y+g(x^{pl})[x+f(x,y)]^{pq+j}-x_2,$$
we use the indeterminate $z:=x+f(x,y)\in x+k[x^{pl},y]$. The system can be rewritten as
$$z+x^{pl}z^{pm}-x_1=y+g(x^{pl})z^{pq+j}-x_2=0.$$ 
Hence $x^{pl}=\frac{x_1-z}{z^{pm}}$ and $y=x_2-g(x^{pl})z^{pq+j}$. Thus
$$y=x_2-g\bigl(\frac{x_1-z}{z^{pm}}\bigr)z^{pq+j}.$$
As $f(x,y)=h(x^{pl},y)$, we compute
$$x=z-h(x^{pl},y)=z-h\Bigl(\frac{x_1-z}{z^{pm}},x_2-g\bigl(\frac{x_1-z}{z^{pm}}\bigl)z^{pq+j}\Bigr).$$
Plugging this expression into the equation $x^{pl}=\frac{x_1-z}{z^{pm}}$ we get that
\begin{equation}\label{EQ46}
\Bigl[z-h\Bigl(\frac{x_1-z}{z^{pm}},x_2-g\bigl(\frac{x_1-z}{z^{pm}}\bigl)z^{pq+j}\Bigr)\Bigr]^{pl}=\frac{x_1-z}{z^{pm}}.
\end{equation}

For part (2), part (2.a) is clear and Equation (\ref{EQ46}) becomes
$$\Bigl[z-\frac{x_1-z}{z^{pm}}\bigl(x_2+\delta z^{pq_1+j}(z-x_1)^s\bigr)^r\Bigr]^{pl}=\frac{x_1-z}{z^{pm}}.$$
Therefore
$$\Bigl[z^{pm+1}+(z-x_1)\bigl(x_2+\delta z^{pq_1+j}(z-x_1)^s\bigr)^r\Bigr]^{pl}+(z-x_1)z^{p^2ml-pm}=0.$$

As $q_1\ge 0$ we have $L(x_1,x_2,w)\in k_{\t}[x_1,x_2,w]$. The field extension $k_{\s}(x,y)$ of $k_{\t}(x_1,x_2)$ is generated by $z$ and we have $L(x_1,x_2,z)=0$. Hence $\deg(e)\le\deg(L)$. Clearly, $\deg(L)\le pl\max\bigl(pm+1,r(pq_1+s+j)+1\bigr)$.

In what follows, all derivatives, degrees, and leading coefficients are with respect to $w$ alone. So let $h_0(x_1,x_2)\in k_{\t}[x_1,x_2]$ be the leading coefficient of $L$. 

As $L'(w)=w^{p^2lm-pm}$ and $L(0)=(-1)^{pl}x_1^{pl}x_2^{plr}$, we have $p^2lm-pm+1\le\deg(L)$. For each $\underline{\gamma}:=(\gamma_1,\gamma_2)\in k^2$ with $\gamma_1\gamma_2h_0(\gamma_1,\gamma_2)\neq 0$, for the solution set $\mathcal J_{\underline{\gamma}}$ in $k$ of the equation $L(\gamma_1,\gamma_2,w)=0$ in $w$ we have $|\mathcal J_{\underline{\gamma}}|\le\deg(L)$ and $\mathcal J_{\underline{\gamma}}\subset k^{\ast}$; this implies that the fiber $e^{-1}(\gamma_1,\gamma_2)$ is equal to
\begin{equation}\label{EQ47}
\Bigl\{\Bigl(\frac{\alpha^{pm+1}+(\alpha-\gamma_1)\bigl(\gamma_2+\delta \alpha^{pq_1+j}(\alpha-\gamma_1)^s\bigr)^r}{\alpha^{pm}},\gamma_2+\delta \alpha^{pq_1+j}(\alpha-\gamma_1)^s\Bigr)\Bigl|\alpha\in\mathcal J_{\underline{\gamma}}\Bigr\}.
\end{equation}
Thus $\deg(e)\ge\deg(L)$. Therefore $\deg(e)=\deg(L)$ and part (2.b) holds.

For part (2.c), let $\underline{\gamma}:=(\gamma_1,\gamma_2)\in k^2$. As $t$ divides both $h(t,w)=t^{pl}w$ and $g(t)$, for $\gamma_1=0$ we have $(0,\gamma_2)\in e^{-1}(0,\gamma_2)$.

Assume now that $\gamma_1\neq 0$. Let $f_{\underline{\gamma}}:=L(\gamma_1,\gamma_2,w)\in k[w]$. If $f_{\underline{\gamma}}$ is not a monomial, then it has a root $\alpha\in k^{\ast}$ and thus 
$$\Bigl(\frac{\alpha^{pm+1}+(\alpha-\gamma_1)\bigl(\gamma_2+\delta \alpha^{pq_1+j}(\alpha-\gamma_1)^s\bigr)^r}{\alpha^{pm}},\gamma_2+\delta \alpha^{pq_1+j}(\alpha-\gamma_1)^s\Bigr)$$ is in $e^{-1}(\gamma_1,\gamma_2)$. If $\gamma_2\neq 0$, then $f_{\underline{\gamma}}(0)\neq 0$. From this and $f_{\underline{\gamma}}'(w)=w^{p^2lm-pm}$ we get that $f_{\underline{\gamma}}$ is not a monomial.

From the last two paragraphs we get that $e$ is non-surjective iff there exists $\gamma_1\in k^{\ast}$ such that for $\underline{\gamma_1}:=(\gamma_1,0)$, $f_{\underline{\gamma_1}}(w)$ is a monomial. 
As 
$$f_{\underline{\gamma_1}}(w)=[w^{pm+1}+\delta^r(w-\gamma_1)^{rs+1}w^{r(pq_1+j)}]^{pl}+(w-\gamma_1)w^{p^2ml-pm}$$ 
with $rs+1\ge 2$, such a $\gamma_1$ exists iff the following conditions hold: (i) we have $pm+1=r(pq_1+j)+rs+1$; (ii) the coefficient $\delta^r$ of $w^{pm+1}$ in $\delta^r(w-\gamma_1)^{rs+1}w^{r(pq_1+j)}$ is $-1$; (iii) the positive integer $rs+1$ is a power of $p$; (iv) we have an identity $pl[pm+1-rs-1]=p^2lm-pm$, equivalently, we have $m=lrs$.

If these conditions hold, then $\gamma_1\in k^{\ast} $ is an arbitrary non-zero root of the polynomial $\delta^{plr}(-t)^{pl(rs+1)}-t=t^{pl(rs+1)}-t\in k[t]$. So we have $\varphi(e)=pl(rs+1)-1$ and $\Gamma_e=\{t\in k|t^{pl(rs+1)-1}=1\}\times\{0\}$. Thus part (2.c) holds.

For part (2.d), as $pm=r(pq_1+s+j)$, the coefficient of $w^{pm+1}$ in 
$$L_1(x_1,x_2,w):=w^{pm+1}+(w-x_1)\bigl(x_2+\delta w^{pq_1+j}(w-x_1)^s\bigr)^r\in k_{\t}[x_1,x_2,w]$$ 
is $1+\delta^r=0$. 

To show that $\deg(L_1^{pl})<p^2lm-pm+1$ we consider two disjoint cases as follows.

{\bf Case 1: $m\le l$.} As $\deg\bigl(w^{pq_1+j}(w-x_1)^s\bigr)\ge 2$, the coefficient of $w^{pm}$ in $L_1(x_1,x_2,w)$ comes from $(w-x_1)[\delta w^{pq_1+j}(w-x_1)^s]^r$, i.e., from the following product $\delta^r w^{r(pq_1+j)}(w-x_1)^{rs+1}$, and hence, as $pm+1=r(pq_1+s+j)+1$, it is $0$ iff $p\mid rs+1$. From this and our hypotheses we get that $\deg(L_1)\le pm-1$. Thus 
$$\deg(L_1^{pl})\le pl(pm-1)=p^2lm-pl\le p^2ml-pm+1\le p^2lm-pm+1.$$ 

{\bf Case 2: $\frac{m}{l}\in (1,p^i-1]$ and $pq_1+s+j\ge p^i$.} As in the previous case, the inequality $\deg\bigl(w^{pq_1+j}(w-x_1)^s\bigr)\ge p^i$ implies that the coefficients of $w^{pm-i_1}$ with $i_1\in\llbracket0,p^i-1\rrbracket$ in $L_1(x_1,x_2,w)$ come from the product $\delta^r(w-x_1)^{rs+1}w^{r(pq_1+j)}$, and hence they are all $0$ as $p^i\mid p^{i_0}\mid rs+1$ and $pm=r(pq_1+s+j)$. From this and our hypotheses we get that $\deg(L_1)\le pm-p^i+1$. So we have relations 
$$\deg(L_1^{pl})\le pl(pm-p^i+1)=p^2lm-pl(p^i-1)<p^2ml-pl(p^i-1)+1\le p^2lm-pm+1.$$ 

From the two cases we get that $\deg(L_1^{p^l})<\deg\bigl((w-x_1)w^{p^2ml-pm}\bigr)$. From this and the fact that $L(x_1,x_2,w)$ is the sum of $L_1(x_1,x_2,w)^{pl}$ and $(w-x_1)w^{p^2ml-pm}$ it follows that $\deg(L)=p^2lm-pm+1$ and that $h_0(x_1,x_2)=1$; so part (2.d) follows from part (2.b) and the definition of a monomial model (see Definition \ref{D6-}(3)).

For part (2.e), the first statement holds by part (2.d). If $e$ is non-surjective, then from part (2.c) we get that $pm=s+1$ is a power of $p$ and $m=lrs$; so both $s$ and $s+1$ are powers of $p$ and this implies that $p=2$ and $m=s=1$. Conversely, if $p=2$ and $m=s=1$, then $e$ is non-surjective by part (2.c). So part (2.e) holds.
\end{proof}

\begin{corollary}\label{C8.3}
Suppose that $p=2$. Then for each $j\in \mathbb N^{\ast}$ there exists $e\in\EE_2(k)$ with $\deg(e)=j$.
\end{corollary}

\begin{proof}
Using composites and an automorphism, it suffices to check this when $j$ is a prime. If $j=2$, then we can take $\e=(x_1+x_1^2,x_2)$. Thus we can assume that $j\ge 3$ is an odd prime. With $m:=\lfloor\frac{j}{2}\rfloor$, the corollary follows from the following two disjoint cases that apply Theorem \ref{T6.6}. 

{\bf Case 1: $m$ is a power of $2$.} Let $i\in\mathbb N$ be such that $m=2^i$. We take $q_1=0$, $j=r=1$, $s=2^{i+1}-1$, and $\delta=1$. The largest power of $2$ that divides $2m$ is $i+1$. As $2^i\le 2^{i+1}-1$, Theorem \ref{T6.6}(2.e) applied to $l=1$ gives that there exists $e\in\EE_2(k)$ with $\deg(e)=2^2m-2m+1=2m+1=j$.

{\bf Case 2: $m$ is not a power of $2$.} Let $i\in\mathbb N^{\ast}$ be such that $2^i<m<2^{i+1}$. We take $q_1=m-2^i\in\llbracket 1,2^i-1\rrbracket$, $j=r=1$, $s=2^{i+1}-1$, and $\delta=1$. The largest power of $2$ that divides $s+1$ is $i+1$. As $m<2^{i+1}-1$, Theorem \ref{T6.6}(2.d) applied to $l=1$ gives that there exists $e\in\EE_2(k)$ with $\deg(e)=2^2m-2m+1=2m+1=j$.
\end{proof}

\begin{example}\normalfont\label{EX41.5}
We apply Theorem \ref{T6.6} with $p=2$, $l=m=q=r=s=j=1$, $h(t,w)=tw$, and $g(t)=t$ (so $\delta=1$). Thus 
$$e=\e(x+x^2y+x^4+x^6y^2,y+x^5+x^6y+x^6y^2+x^8y^3)\in\EE_2(k)$$ 
is such that $\deg(e)=3$ and $\varphi(e)=3$ by Theorem \ref{T6.6}(2.a) to (2.c). As Conditions ($\triangleleft$) and ($\triangleright$) of Proposition \ref{PR1} hold, we have 
$$9\le\pi_{\i}(e)\le \pi_{\l}(e)=\pi_{\r}(e)=11\ge \pi(e)$$ 
by Proposition \ref{PR1}(1) to (3). The weaker form of this special case that does not prove the non-surjectivity part and does not compute the invariants $\varphi(e)$, $\pi_{\r}(e)$, and $\pi_{\i}(e)$ was first obtained using computer programming in \cite{Mon}, Thm.\ 1.2. 

We compute that 
$$L=\bigl[w^3+(w+x_1)\bigl(w(w+x_1)+x_2)\bigr]^2+(w+x_1)w^2=w^3+w^2(x_1^4+x_1+x_2^2)+x_1^2x_2^2.$$
The field extension $k_{\t}(x_1,x_2)\rightarrow k_{\s}(x,y)=k_{\t}(x_1,x_2)(z)\cong k_{\t}(x_1,x_2)[w]/(L)$ is non-Galois by Corollary \ref{C10}(2).

As $z$ is a root of the monic polynomial $L(x_1,x_2,w)\in k_{\t}[x_1,x_2,w]$, i.e., we have $L(x_1,x_2,z)=0$, it follows that $z\in A_e$. Let $v:=\frac{z}{x_1x_2}+x_1^2+x_2\in k_{\t}(x_1,x_2)(z)$. So $v+x_1^2+x_2=\frac{z}{x_1x_2}$ is a root of the polynomial 
\begin{equation*}
\begin{aligned}
L_2=L_2(x_1,x_2,w)&:=w^3+(x_1^4+x_1+x_2^2)w+x_1x_2\\
&=w(w+x_1^2+x_2)^2+x_1(w+x_2)\in k_{\t}[x_1,x_2,w]
\end{aligned}
\end{equation*}
and thus $v$ is a root of the 
polynomial 
$$L_3=L_2(x_1,x_2,w):=w^3+(x_1^2+x_2)w^2+x_1w+x_1^3\in k_{\t}[x_1,x_2,w].$$ 
Hence $v\in A_e$. If $B_e$ is the $k_{\t}[x_1,x_2]$-subalgebra of $A_e$ generated by either $v$ or $v+x_1^2+x_2$, then we identify $B_e:=k[x_1,x_2][w]/(L_2)$ but for some computations it is convenient to use that $B_e\cong k[x_1,x_2][w]/(L_3)$. So $Y_e:=\Spec B_e$ is a monomial model of $X_e$. 

The finite flat morphism $Y_e\rightarrow\mathbb A^2_{k,\t}$ becomes \'etale after inverting $x_1x_2$ and for $i\in\{1,2\}$, there exist precisely two points $\eta_{i,1}$ and $\eta_{i,2}$ of $Y_e$ that map to the generic point $\eta_i$ of the curve in $\mathbb A^2_{k,\t}$ which is the zero locus $x_i=0$. From Zariski's Jacobian Criterion (see \cite{Gro2}, Prop.\ 22.6.7(iii)) we get that $\Sing(Y_e)$ is the zero locus $x_1=w+x_2=0$. Reindexing and using that $\Reg(Y_e)$ is an open subvariety of $X_e$, we can assume that $\eta_{1,1}\in\Reg(Y_e)$, that $\eta_{1,2}$ is the generic point of $\Sing(Y_e)$, that $\eta_{2,1}\in\mathbb A^2_{K,\s}\cap\Reg(Y_e)$, and that $\eta_{2,2}\in [X_e\setminus\Imm(\imath_e)]\cap\Reg(Y_e)$. The residue fields of $\eta_{2,1}$ and $\eta_2$ are $k(x_1)$ and the residue field of $\eta_{2,2}$ is $k(x_1)(\sqrt{x_1})$ and hence is a purely inseparable extension of $k(x_1)$ of degree $2$. So $\psi_e$ is non-\'etale. 

As $\eta_{2,2}\in X_e$, we can speak about its Zariski closure $Y_{2,2}$ in $X_e$. We have $Y_{2,2}\in \Irr(X_e\setminus X_e^{\et})\setminus \Irr_{\textup{t}}(X_e\setminus X_e^{\et})$.  

%We consider the discrete valuation ring $D:=k(x_2)[[x_1]]$ with uniformizer $x_1$ and the $D$-algebra $B:=D[w]/(L_3)$. From Hensel's Lemma we get that for $i\in\{1,2,3\}$ there exists a $D$-algebra retraction $\sigma_i:B\rightarrow D$ such that $\sigma_1\bigl(w+(L_3)\bigr)\in (x_1^2)$, $\sigma_2\bigl(w+(L_3)\bigr)\equiv x_2\pmod{(x_1)}$, and $\sigma_3\bigl(w+(L_3)\bigr)\equiv\frac{x_1}{x_2}\pmod{(x_1^2)}$. Hence there exist two distinct points $\eta_{1,2,1}$ and $\eta_{1,2,2}$ of $X_e$ that map to $\eta_{1,2}\in Y_e$. As $X_e$ has three points $\eta_{1,1}$, $\eta_{1,2,1}$, and $\eta_{1,2,2}$ that map to $\eta_1$, we get that these points belong to $X_e^{\et}$. We conclude that $Y_{2,2}=X_3\setminus X_3^{\et}$ and therefore $\rho(e)=1+\rho_{\et}(e)$. Hence $|\Sing(X_e)|\le 1$ by Corollary \ref{C2.8}(6).

The endomorphism $e:\mathbb A^2_{k,\s}\rightarrow\mathbb A^2_{k,\t}$ becomes finite \'etale after inverting $x_1x_2$ by Equation (\ref{EQ47}); so $\Spec (A_e)_{x_1x_2}=\Spec (B_e)_{x_1x_2}$ is an open subvariety of $\Imm(\imath_e)$. 

As for each $\gamma_2\in k^{\ast}$, the polynomial $L(0,\gamma_2,w)$ has a unique solution in $k^{\ast}$ and as for $z=0$ we get two additional points $(0,\gamma_2)$ and $(\gamma_2^{-1},\gamma_2)$ in the fiber $e^{-1}(0,\gamma_2)$, we get that $|e^{-1}(0,\gamma_2)|=3$. So there exists three points $\eta_{1,1}$, $\eta_{1,2,1}$, and $\eta_{1,2,2}$ of $\Imm(\imath_e)$ that map to either $\eta_{1,1}$ or $\eta_{1,2}$. 

From the last two paragraphs we get that $\rho_{\et}(e)=0$ and hence $\rho(e)=1$.\end{example}%As $k_{\t}[[x_1,x_2]]$ is a unique factorization domain, it is easy to see that the $k_{\t}[[x_1,x_2]]$-algebra $k_{\t}[[x_1,x_2]][w]/(L_2)$ has no retraction. This implies that there exists a unique point $P\in X_e$ that maps to $(0,0)\in\mathbb A^2_{k,\t}$. Therefore $Y_{1,1}\cap Y_{2,1}\cap Y_{2,2}=\{P\}$; in particular, $Y_{1,1}\subset X_e^{\et}$ by Corollary \ref{C2.8}(3) and the curves $Y_{1,1}\cap\mathbb A^2_{k,\s}$ and $Y_{2,1}\cap\mathbb A^2_{k,\s}$ have 2 or more points at infinity. [TO BE FINALIZED: We want to compute $\rho(e)-\rho_{\et}(e)$ and to show that $P\in\Sing(X_e)$ and hence that $X_e$ is non-regular.]

\section{On extending vector bundles}\label{S29}

For $n\in\mathbb N^{\ast}\setminus\{1\}$, the subset $\mathfrak D_n$ of $(\mathbb N^{\ast}\setminus\{1\})^2$ was introduced in Notation \ref{NOT5} and to each pair $(m,q)\in\mathfrak D_n$ a finite set of primes $P_{n,m,q}$ was introduced in Corollary \ref{C9.1}(1). For $m=2$ we have $P_{n,m,q}=\{2\}$ by Corollary \ref{C10}(2) and the fact that for $p=2$ and $d:=\e(x_1+x_1^2,x_2)\in\EE_2(k)$ we have $\deg(e)=2=\pi(e)$.

The case $m\ge 3$ is a lot more complex and in this section we aim to develop tools that could be used to compute the sets $P_{n,m,q}$ with $(3,q)\in\mathfrak D_n$ and that are often used in the next section for $m=3$. 

We begin with a definition and then we present two general results on vector bundles on principal open subvarieties of affine spaces $\mathbb A^n_K$ with $n\in\mathbb N^{\ast}\setminus\{1\}$ and an application to the embedding dimensions of the $X_e$s with $e\in\EE_2(K)$ such that $\psi_e$ is non-regular.

\begin{definition}\label{D11.9}
Let $n\in\mathbb N^{\ast}\setminus\{1\}$ and $h\in R:=K[x_1,\ldots,x_n]$. 

\medskip
{\bf (1)} We say that $h$ is affinely univariate\index{polynomial!affinely univariate polynomial} if there exists $(y_1,\ldots,y_n)\in\mathcal A_n(K)$ such that each irreducible factor of $h$ is a univariate polynomial in an indeterminate in $K\oplus_{i=1}^n Ky_i$.

\smallskip
{\bf (2)} We say that $h$ has the freeness property\index{polynomial!has the freeness property} if each finitely generate projective $R_h$-module is free.

\smallskip
{\bf (3)} For $e\in\EE_n(K)$ we introduce the following \'etale endomorphism conditions.\index{endomorphism condition}

\medskip
{\bf ($\ast^e_h$)} For each irreducible component $Z$ of the zero locus $h=0$, the inverse image $e^{-1}(Z)$ is irreducible (equivalently, the group monomorphism $R_h^{\ast}\rightarrow (A_e)_h^{\ast}$ is an isomorphism).\index{endomorphism condition!($\ast^e_h$)} 

\smallskip
{\bf ($\dag_h^e$)} Either $n\ge 3$ and the $R_h$-module $(A_e)_h/R_h$ is free or $n=2$.\index{endomorphism condition!($\dag^e_h$)}
\end{definition}

\begin{theorem}\label{T8}
Let $n\in\mathbb N^{\ast}\setminus\{1\}$ and $h\in R:=K[x_1,\ldots,x_n]$. Then the following properties hold.

\medskip
{\bf (1) (Quillen--Suslin--Roitman)} If $h$ has the freeness property, then it has the stably freeness property, i.e., for each $m\in\mathbb N$, every finitely generated projective $R_h[x_{n+1},\ldots,x_{n+m}]$-module is free.

\smallskip
{\bf (2) (Gabber)} If either $n=2$ or $h$ is affinely univariate, then $h$ has the freeness property.  

\smallskip
{\bf (3) (Gabber)} Suppose that $n=3$. Let $\{Y_1,\ldots,Y_m\}$ be the set of irreducible components of the zero locus $h=0$. Assume that there exist smooth integral curves $C_1,\ldots, C_m$ over $\Spec K$ such that either they are affine and $Y_i\cong C_i\times_{\Spec K} \mathbb A^1_K$ for each $i\in\llbracket1,m\rrbracket$ or they are projective and $Y_i$ is smooth and birational to the projectivization of a vector bundle over $C_i$ of rank $2$ for each $i\in\llbracket1,m\rrbracket$. Then $h$ has the freeness property.
\end{theorem}

\begin{proof}
For part (1), every finitely projective $R[x_{n+1},\ldots,x_{n+m}]$ is free by Quillen--Suslin Theorem (e.g., see \cite{Lang-S2}, Ch.\ XXI, Sect.\ 3, Thm.\ 3.7) and thus it is isomorphic to the scalar extension of a  finitely projective $R$-module. From this and \cite{Roi}, Prop.\ 2 we get that each finitely generated projective $R_h[x_{n+1},\ldots,x_{n+m}]$-module is isomorphic to the scalar extension of a  finitely projective $R_h$-module and hence it is free as $h$ has the freeness property. So part (1) holds.

For part (2), the case $n\ge 3$ is proved in \cite{Ga2}, Thm.\ 2.1. If $n=2$, then each vector bundle over $\Spec R_h$ extends to a vector bundle over $\Spec R$ punctured in a finite number of points and thus also to a vector bundle over $\Spec R=\mathbb A^2_{K,\t}$ by \cite{H3}, Cor.\ 1.4, which is trivial by \cite{Ses}, Thm. So part (2) holds.

For part (3) see \cite{Ga2}, Cor.\ 2.5.
\end{proof}

\begin{theorem}\label{T9}
Let $(n,m)\in (\mathbb N^{\ast}\setminus\{1\})^2$ and $R:=K[x_1,\ldots,x_n]$. Bet $B$ be the normalization of $R$ in a finite field extension $K_0$ of $K(x_1,\ldots,x_n)$ of degree $m$. Then the following properties hold.

\medskip
{\bf (1) (Hochster)} The short exact sequence $0\rightarrow R\rightarrow B\rightarrow B/R\rightarrow 0$ of $R$-modules splits.

\smallskip
{\bf (2)} If $n=2$, then the $R$-module $B/R$ is a free of rank $m-1$ and we have $\edim(\Spec B)\le m+1$. 

\smallskip
{\bf (3)} If $m=2$, then $B/R$ is a free $R$-module of rank $1$ and $\edim(\Spec B)\le n+1$.

\smallskip
{\bf (4) (Hartshorne)} If $n\ge 3$ and $m\ge 3$, then there exists a closed subset $Y$ of $\Spec R=\mathbb A^n_K$ of codimension at least $3$ and such that the restriction to $\Spec R\setminus Y$ of the coherent $\mathcal O_{\Spec R}$-module defined by $B/R$ is locally free of rank $m-1$.
\end{theorem}

\begin{proof}
We have a direct sum decomposition $B_e=R\oplus B_1$ of $R$-modules by \cite{Ho}, Thm.\ 2. As the $R$-modules $B_1$ and $B/R$ are isomorphic, part (1) holds.

For parts (2) and (3), the $R$-module $B_1$ is torsionfree and the coherent $\mathcal O_{\Spec R}$-module $\mathfrak F_1$ it defines is normal (as $B$ is normal) in the sense of Barth and hence it is reflexive by \cite{H3}, Prop.\ 1.6 and Def.\ before it. From this and \cite{H3}, Cor.\ 1.4 we get that part (4) holds and that, if $n=2$, $\mathfrak F_1$ is locally free of rank 1.

From the last property for $n=2$ and \cite{Ses}, Thm.\ we get that $B/R$ is a free $R$-module of rank $m-1$. If $(y_1,\ldots,y_{m-1})\in B^{m-1}$ maps to an $R$-basis of $B/R$, then the $R$-algebra homomorphism $R[x_3,\ldots,x_{m+1}]=K_{\t}[x_1,\ldots,x_{m+1}]\rightarrow B$ that maps $x_{2+i}$ to $y_i$ for each $i\in\llbracket1,m-1\rrbracket$ is surjective and hence $\edim(\Spec B)\le m+1$. So part (2) holds. 

For part (3), as the coherent $\mathcal O_{\Spec R}$-module $\mathfrak F_1$ is reflexive of rank $1$, it is locally invertible by \cite{H3}, Prop.\ 1.9. From this and the fact that $R$ is factorial, we conclude that there exists an isomorphism $B_1\cong R$ of $R$-modules. So the $R$-algebra $B$ is a quotient of $R[x_{n+1}]$ and hence $\edim(\Spec B)\le n+1$. So part (3) holds.
\end{proof}

In Definition \ref{D11.9}(3), the $R_h$-module $(A_e)_h/R_h$ is free of rank $\deg(e)-1$ by Theorem \ref{T9}(2) applied to $B=A_e$ if $n=2$.

\begin{corollary}\label{C25}
Let $e\in\EE_2(K)$ be such that $\Psi_e$ is non-\'etale. Then the following properties hold.

\medskip
{\bf (1)} We have $\edim(X_e)\in\llbracket 4,\deg(e)+1\rrbracket$. In particular, $\deg(e)\ge 3$ and, if $\deg(e)=3$, then $\edim(X_e)=4$.

\smallskip
{\bf (2)} If $\psi_e$ is regular, then $\edim(X_e)\in\{4,5\}$.
\end{corollary}

\begin{proof}
We have $\edim(X_e)\le \deg(e)+1$ by Theorem \ref{T9}(2). As $\psi_e$ is non-\'etale, we have $\Imm(\imath_e)\subsetneq X_e$ and hence $X_e\not\cong\mathbb A^2_K$; so $\edim(X_e)\ge 3$. If $\edim(X_e)=3$, then $X_e$ is a hypersurface in $\mathbb A^2_K$, hence a complete intersection, and this contradicts Corollary \ref{C11}(3). Thus $\edim(X_e)\ge 4$ and part (1) holds.

Part (2) follows from part (1) and Theorem \ref{T3}(1).
\end{proof}

\section{On geometric degrees $3$, part I: general theory}\label{S30}

All our basic results that are specific to geometric degree $3$ are grouped together as follows.

\begin{theorem}\label{T10}
Let $e\in\EE_n(K)$ be such that $\deg(e)=3$. Let $R:=K_{\t}[x_1,\ldots,x_n]$ and $S:=K_{\s}[x_1,\ldots,x_n]$. Let $(h_0,h_1)\in R^2$ be such that $N_e^{\textup{n}-\et}=\Spec R/Rh_0$ and $N_e=\Spec R/Rh_1$ (so $h_0\in K^{\ast}$ iff $\psi_e$ is \'etale and $h_1\in K^{\ast}$ iff $\imath_e$ is an isomorphism). Let $h\in R$ be such that $h_0\mid h\mid h_1$. We consider the following degree $3$ Galois condition.

\medskip
{\bf ($S_3$)} The field extension $\Frac(e^{\#}):\Frac(R)\rightarrow\Frac(S)$ is non-Galois.\index{endomorphism condition!($S_3$)}

\medskip
Then the following properties hold.

\medskip
{\bf (1)} If Condition ($S_3$) does not hold, then $\chr(K)=3$ and the finite morphism $X_e\rightarrow\mathbb A^2_{K,\t}$ is an Artin--Schreier cover.

\smallskip
{\bf (2)} Suppose that Condition ($S_3$) holds. Then the following properties hold.

\medskip\noindent
{\bf (2.a)} The field extension $\Frac(e^{\#})$ generates a Galois extension $\Frac(R)\rightarrow K_0$ of Galois group $G\cong S_3$. 

\smallskip\noindent
{\bf (2.b)} There exists a finite morphism 
$$\phi:Y_e=\Spec B_e\rightarrow A^n_{K,\t}=\Spec R$$ 
between integral normal varieties which is generically \'etale of degree $2$ and such that the normalization $W_e=\Spec C_e$ of $X_e$ in $K_0$ is the normalization of $X_e\times_{\mathbb A^n_{K,\t}} Y_e$ (so $K_0=\Frac(C_e)$, $W_e$ is also the normalization of $Y_e$ in $K_0$, and $f_e:W_e\rightarrow X_e$ is as in Example \ref{EX5+}). In particular, we have a diagram of finite surjective morphisms 

\[\xymatrix{
W_e=\Spec C_e \ar[r]^{f_e} \ar[d]^{} & X_e=\Spec A_e \ar[d]^{} \\
Y_e=\Spec B_e\ar[r]^{\phi} & \mathbb A^n_k=\Spec R.
}\]

\noindent
{\bf (2.c)} We have an inequality $\dim\bigl(\Sing(Y_e)\bigr)\le n-2$ and moreover the finite morphism $W_e\times_{Y_e} \Reg(Y_e)\rightarrow\Reg(Y_e)$ is Galois. 

\smallskip\noindent
{\bf (2.d)} The non-\'etale locus of $\phi$ is the zero locus $h_0=0$ in $Y_e$.

\smallskip\noindent
{\bf (2.e)} The $R$-module $B_e$ is free of rank $2$.

\smallskip\noindent
{\bf (2.f)} If $\chr(K)\neq 2$, then $B_e\cong R[z]/(z^2-h_0)$. In particular, $h_0\notin K^{\ast}$ and the diagram of part (2.b) is not a cartesian diagram. 

\smallskip\noindent
{\bf (2.g)} If $\chr(K)=2$, then $(B_e)_{h_0}\cong R_{h_0}[z]/(z^2+z+f_0)$ for a suitable $f_0\in R_{h_0}$.

\smallskip\noindent
{\bf (2.h)} If $\chr(K)=2$, then there exists a triple $(g_0,g_1,l)\in R^2\times\mathbb N$ such that $g_1\mid h_0^l$ and $B_e\cong R[z]/(z^2+h_0g_1z+g_0)$.

\smallskip\noindent
{\bf (2.i)} We have $\ciedim(Y_e)\le n+1$. 

\smallskip\noindent
{\bf (2.j)} If $\dim\bigl(\Sing(Y_e)\bigr)\le n-3$, then the finite morphism $W_e\rightarrow Y_e$ is Galois (and hence \'etale).

\smallskip\noindent
{\bf (2.k)} If $\psi_e$ is \'etale, then $\chr(K)=2$.

\smallskip\noindent
{\bf (2.l)} Suppose that $\psi_e$ (i.e., $X_e$) is regular. Then the $A_e$-module $C_e/A_e$ is projective of rank $1$. In particular, $C_e$ is a projective $A_e$-module of rank $2$ and thus the morphism $W_e\rightarrow X_e$ is flat.

\medskip
{\bf (3)} Suppose that Condition ($S_3$) holds and $n=2$. Then the following properties hold.

\medskip\noindent
{\bf (3.a)} If the morphism $W_e\rightarrow X_e$ (resp.\ $W_e\rightarrow Y_e$) is flat, then we have an inequality $|\Sing(X_e)|\le|\Sing(W_e)|$ (resp.\ $|\Sing(Y_e)|\le|\Sing(W_e)|$).

\smallskip\noindent
{\bf (3.b)} Let $Z$ be an irreducible component of $N_e$. If $Z$ is not an irreducible component of $N_e^{\textup{n}-\et}$ and there exists a unique irreducible component $Y$ of $X_e^{\et}\setminus\Imm(\imath_e)$ that maps onto $Z$, then we assume that the finite morphism $Y\rightarrow Z$ is generically \'etale of degree $2$. Then $Z$ is isomorphic to $\mathbb A^1_K$.

\smallskip\noindent
{\bf (3.c)} If $\chr(K)\neq 2$, then $\psi_e$ is regular.

\smallskip\noindent
{\bf (3.d)} The element $h_0\in R$ has the freeness property.

\smallskip\noindent
{\bf (3.e)} If $\chr(K)\neq 2$, then $\psi_e(X_e\setminus X_e^{\et})$ has $\rho(e)-\rho_{\et}(e)$ connected components which are irreducible and isomorphic to $\mathbb A^1_K$. In particular, $\nu_{\textup{n}-\et}(e)=\rho(e)-\rho_{\et}(e)$.

\smallskip\noindent
{\bf (3.f)} If $\chr(K)\neq 2$, then $W_e$ and $Y_e$ are regular and the morphism $W_e\to Y_e$ is finite Galois.

\smallskip\noindent
{\bf (3.g)} Suppose the finite morphism $W_e\rightarrow Y_e$ is Galois (e.g., this holds if $\chr(K)\neq 2$ or $Y_e$ is regular). Let $Z_1$ be an irreducible component of $N_e^{\textup{n}-\et}$. Then $\psi_e^{-1}(Z_1)=C_1\sqcup Y_1$ is a disjoint union of two irreducible closed curves of $X_e$ isomorphic to $\mathbb A^1_K$ and such that $C_1\cap\Imm(\imath_e)\neq\emptyset$ and $Y_1\subset X_e\setminus \Imm(\imath_e)$.%, $Z_1\subset\Imm(e)$, $Z_1\cap Z_2=\emptyset$ if $Z_1\neq Z_1$, and $Z_2\subset\Imm(e)$ if $Z_2\cong\mathbb A^1_K$. In particular, $\Gamma_e\subset \cup_{Z\in\Irr(N_e)\setminus \\Irr(N_e^{\textup{n}-\et})} \Sing(Z)$.

\medskip
{\bf (4)} If $h$ has the freeness property (e.g., $n=2$ or $h$ is affinely univariate), then Condition ($\dag^e_h$) holds.

\smallskip
{\bf (5)} Suppose that Condition ($\dag^e_h$) holds. Then there exists a unimodular quadruple $(\alpha,\beta,\gamma,\delta)\in R_h^4$ such that by denoting $L(x,y):=\alpha x^3+\beta x^2y+\gamma xy^2+\delta y^3\in R_h[x,y]$, for each $R_h$-algebra $A$, the $A$-algebra $(A_e)_h\otimes_{R_h} A$ is generated by $1$ element iff there exists a pair $(u,v)\in A^2$ such that $L(u,v)\in A^{\ast}$. %Moreover, there exists $(u,v)\in 
%Then the $R_h$-algebra $(A_e)_h$ is generated by $1$ element and $\chr(K)=3$.

\smallskip
{\bf (6)} Suppose Condition ($\dag^e_h$) holds and $\chr(K)\neq 2$. Then there exists a pair $(v_1,v_2)\in (A_e)_h^2$ such that $L(v_1,v_2)\in R_h^{\ast}$.%If $\chr(K)=2$, then we also assume that condition ($S_3$) does not hold. 

\smallskip
{\bf (7)} Let $h_2\in R$ be such that $h=h_0h_2$. If Condition ($\ast^e_{h_2}$) holds and the $R_h$-algebra $(A_e)_h$ is generated by $1$ element, then $\chr(K)=3$.

\smallskip
{\bf (8)} Suppose that Condition ($S_3$) holds. Then the following properties hold.

\medskip\noindent
{\bf (8.a)} If $\chr(K)\neq 2$, then $\ciedim(\Spec R_h)=\ciedim\bigl(\Spec (B_e)_h\bigr)=n+1$, $\ciedim\bigl(\Spec (A_e)_h\bigr)\le 2n+2$, and $\ciedim\bigl(\Spec (C_e)_h\bigr)\le 2n+2$. 

\smallskip\noindent
{\bf (8.b)} If $\chr(K)=2$, then $\ciedim(\Spec R_h)\le n+1$, $\ciedim\bigl(\Spec (B_e)_h\bigr)\le n+2$, $\ciedim\bigl(\Spec (A_e)_h\bigr)\le 2n+2$, and $\ciedim\bigl(\Spec (C_e)_h\bigr)\le 2n+2$. 

\medskip
{\bf (9)} Suppose that Condition ($\dag^e_h$) holds. We consider an $R_h$-algebra presentation $(A_e)_h=R_h[y_1,y_2]/I$ such that $(1+h,y_1+I,y_2+I)$ is an $R_h$-basis of $(A_e)_h$. Then the ideal $I$ is generated by $3$ polynomials of degree $2$ in $y_1$ and $y_2$. Moreover, if $h$ has the freeness property, then $I$ is generated by $2$ elements iff the $R_h[y_1,y_2]$-module $\Ext^1_{R_h[y_1,y_2]}(I,R_h[y_1,y_2])$ is cyclic, and in such a case, if Condition ($S_3$) also holds, we have inequalities $\ciedim\bigl(\Spec (A_e)_h\bigr)\le n+3$ and $\ciedim\bigl(\Spec (C_e)_h\bigr)\le n+3$.

\smallskip
{\bf (10)} Assume that $\chr(K)\neq 3$ and Condition ($S_3$) holds. Let 
$$\langle e_h\rangle\in H^1_{\et}\bigl(\Spec (B_e)_h,(\mathbb Z/3\mathbb Z)_K\bigr)$$ 
be the class that defines the Galois cover $\Spec (C_e)_h\rightarrow\Spec (B_e)_h$. Let $M_h$ be a projective $(B_e)_h$-module of rank $1$ whose class 
$$\langle M_h\rangle\in H^1_{\et}\bigl(\Spec (B_e)_h,\mathbb G_{\m,K}\bigr)$$ 
is the image of $\langle e_h\rangle$ via the monomorphism $(\mathbb Z/3\mathbb Z)_K\rightarrow\mathbb G_{\m,K}$. Then the following properties hold.

\medskip\noindent
{\bf (10.a)} Assume that $M_h\cong (B_e)_h$. Then there exists a unit $u\in (B_e)_h^{\ast}$ such that 
$$(C_e)_h\cong (B_e)_h[t]/(t^3-u),$$ 
and we have inequalities $\ciedim(\Spec (C_e)_h)\le n+2$ and $\ciedim(\Spec (A_e)_h)\le n+3$ if $\chr(K)\neq 2$ and $\ciedim\bigl(\Spec (C_e)_h\bigr)\le n+3$ if $\chr(K)=2$.

\smallskip\noindent
{\bf (10.b)} There exists a projective $B_e$-module $M$ of rank $1$ such that the notation matches, i.e., $M_h=M\otimes_{B_e} (B_e)_h$, iff there exists a closed subvariety $Z$ of the zero locus $h=0$ of $Y_e$ and a line bundle over $Y_e\setminus Z$ such that $\dim(Z)\le n-4$ and the line bundle restricts to a line bundle over $\Spec (B_e)_h\setminus Z$ that defines $\langle M_h\rangle$.

\smallskip\noindent
{\bf (10.c)} If we have $\dim\bigl(\Sing(Y_e)\bigr)\le n-3$ (e.g., this holds if $n=2$ and $\chr(K)\neq 2$ by part (3.f)), then the class $\langle e_h\rangle$ is the image of a class $\langle e\rangle\in H^1_{\et}\bigl(\Spec B_e,(\mathbb Z/3\mathbb Z)_K\bigr)$ and hence there exists a projective $B_e$-module $M$ of rank $1$ such that the notation matches, i.e., $M_h=M\otimes_{B_e} (B_e)_h$.

%\smallskip\noindent
%{\bf (10.d)} If $\chr(K)\neq 2$, then the class $\langle M_{h_0}\rangle$ is non-trivial. 

\medskip
{\bf (11)} Assume that Condition ($S_3$) holds, $\chr(K)\neq 3$, and $h$ has the freeness property. Then the following properties hold.

\medskip\noindent
{\bf (11.a)} The $R_h$-module $M_h$ is free of rank $2$.

\smallskip\noindent
{\bf (11.b)} The $(B_e)_h$-module $M_h$ is generated by $2$ elements and we have a $(B_e)_h$-linear isomorphism $M_h\oplus M_h^{\otimes 2}\cong (B_e)_h^2$.

\smallskip\noindent
{\bf (11.c)} Suppose $\chr(K)\neq 2$. Then there exists a triple $(\beta_z,\gamma_z,\delta_z)\in R_h^3$ such that the identity $\delta_z^2+\beta_z\gamma_z=h$ holds and by defining $$L_2(x,y):=\gamma_zx^2+2\delta_z xy-\beta_zy^2\in (A_e)_h[x,y],$$ we have $M_h\cong (B_e)_h$ iff there exists $(u,v)\in (B_e)_h^2$ such that $L_2(u,v)\in (B_e)_h^{\ast}$. Moreover, if the pair $(\beta_z,\gamma_z)\in R_h^2$ is unimodular, then the cocycle that defines the class $\langle M_h\rangle$ is the unit $z+\delta_z\in (B_e)_{h\beta_z\gamma_z}^{\ast}$. 

\smallskip\noindent
{\bf (11.d)} Suppose $\chr(K)=2$. Then, with $f_0\in R_{h_0}\subset R_h$ as in part (2.g), there exists a triple $(\beta_z,\gamma_z,\delta_z)\in R_h^3$ such that the identity $\delta_z^2+\delta_z+\beta_z\gamma_z=f_0$ holds and by defining $$L_2(x,y):=\gamma_zx^2+xy+\beta_zy^2\in (A_e)_h[x,y],$$ we have $M_h\cong (B_e)_h$ iff there exists $(u,v)\in (B_e)_h^2$ such that $L_2(u,v)\in (B_e)_h^{\ast}$. 

\medskip
{\bf (12)} Suppose that $\chr(K)\neq 3$, Condition ($S_3$) holds, and there exists a projective $B_e$-module $M$ of rank $1$ such that $M_h=M\otimes_{B_e} (B_e)_h$. Then the following properties hold.

\medskip\noindent
{\bf (12.a)} The $R$-module $M$ is free of rank $2$.

\smallskip\noindent
{\bf (12.b)}  The $B_e$-module $M$ is generated by $2$ elements. If $\dim\bigl(\Sing(Y_e)\bigr)\le n-3$, then we can choose $M$ such that $M\oplus M^{\otimes 2}\cong B_e^2$.

\smallskip\noindent
{\bf (12.c)} Assume $\chr(K)\neq 2$. Then there exists a triple $(\beta_z,\gamma_z,\delta_z)\in R^3$ such that the identity $\delta_z^2+\beta_z\gamma_z=h$ holds and for $L_2(x,y)$ as in part (11.c) we have $M\cong B_e$ iff there exists $(u,v)\in B_e^2$ such that $L_2(u,v)\in B_e^{\ast}$. Moreover, if the pair $(\beta_z,\gamma_z)\in R^2$ is unimodular, then the cocycle that defines the class $\langle M\rangle$ is the unit $z+\delta_z\in (B_e)_{\beta_z\gamma_z}^{\ast}$.

\smallskip\noindent
{\bf (12.d)} Suppose $\chr(K)=2$. Then there exists a triple $(\beta_z,\gamma_z,\delta_z)\in R^3$ such that the identity $\delta_z^2+\delta_z+\beta_z\gamma_z=f_0$ holds and for $L_2(x,y)$ as in part (11.d)  we have $M\cong B_e$ iff there exists $(u,v)\in B_e^2$ such that $L_2(u,v)\in B_e^{\ast}$. 
\end{theorem}

\begin{proof}
For part (1) we have $\chr(K)\neq 0$ by \cite{Raz}, Thm.\ 2 and $\chr(K)\neq p$ with $p$ a prime different from $3$ by Corollary \ref{C10}(2). Thus $\chr(K)=3$ and the finite morphism $X_e\rightarrow\mathbb A^2_{K,\t}$ is Galois by Corollary \ref{C10}(2) and an Artin--Schreier cover by Lemma \ref{F6}. So part (1) holds.

Parts (2.a) and (2.b) are standard Galois theory. 

To show that the morphism $W_e\times_{Y_e} \Reg(Y_e)\rightarrow\Reg(Y_e)$ is Galois (equivalently, is \'etale) we can assume that $X_e\neq X_e^{\et}$; thus $h_0\notin K$. Let $Y\in\Irr(N_e^{\textup{n}-\et})$ and let $\eta$ be its generic point; so $\bigl(\psi_e^{-1}(\eta)\bigr)_{\red}$ consists of $2$ points $\eta_1\in\Imm(\imath_e)$ and $\eta_2\in X_e\setminus X_e^{\et}\subset X_e\setminus\Imm(\imath_e)$. Thus $e^{-1}(Y)$ is irreducible and $\eta_2$ is ramified over $\eta$. Also $\bigl(\phi^{-1}(\eta)\bigr)_{\red}$ is one point $\eta_0$ of $Y_e$ which is ramified over $\eta$. The inverse image of $\eta_0$ in $W_e$ has at least $2$ points permuted transitively by a cyclic group of order $3$ and hence it has $3$ points; thus $W_e\rightarrow Y_e$ is \'etale at these $3$ points. As $\deg(e)=3$ implies that the generic points of $\psi_e^{-1}\bigl(\psi_e(X_e^{\et}\setminus\Imm(\imath_e))\bigr)$ are contained in $X_e^{\et}$, it follows that $W_e\rightarrow Y_e$ is \'etale in codimension $1$ and hence it is \'etale outside a codimension at least $2$ closed subset of $W_e$. Thus, as the non-\'etale locus of $W_e\times_{Y_e} \Reg(Y_e)\rightarrow\Reg(Y_e)$ is either empty or of pure codimension $1$ by the affineness of \'etale loci (see \cite{Stacks1} or \cite{GZ}, Thm.\ 1), we get that it is empty. So $W_e\times_{Y_e} \Reg(Y_e)\rightarrow\Reg(Y_e)$ is a finite \'etale cover. As $Y_e$ is normal, $\Sing(Y_e)$ has codimension at least $2$ in $Y_e$. So part (2.c) holds. 

The images in $\mathbb A^n_{K,\t}$ of the non-\'etale loci of $\phi$ and $W_e\rightarrow\mathbb A^n_{K,\t}$ are $\Spec R/Rh_0$. From this and the ramification of $\phi$ at $\eta_0$s we get that part (2.d) holds. 

Part (2.e) holds by Theorem \ref{T9}(3).

For parts (2.f) to (2.h), let $\tilde z\in B_e$ such that $(1,\tilde z)$ is an $R$-basis of $B_e$. Let $(f_1,g_0)\in R$ be such that $\tilde z^2=f_1\tilde z+g_0$. 

For part (2.f), as $\chr(K)\neq 2$, we can assume that $f_1=0$. So we can identify $B_e=R[z]/(z^2-g_0)$ via the isomorphism that maps $\tilde z$ to $z+(z^2-g_0)$. As $B_e$ is a normal integral domain, $g_0\in R\setminus K$ is a product of non-associated irreducible polynomials. The non-\'etale locus of $\phi$ is defined by $g_0$. From this and part (2.d) we get that the zero loci $h_0=0$ and $g_0=0$ in $\mathbb A^n_{K,\t}$ coincide. Thus we can assume that $g_0=h_0$; so $h_0\notin K^{\ast}$. With the notation of the proof of part (2.b), the inverse images of $\eta_1$ in $W_e$ is one point ramified over $\eta_1$ and the inverse images of $\eta_2$ in $W_e$ is two points unramified over $\eta_1$. This implies that the diagram of part (2.b) is not cartesian. So part (2.f) holds. 

Part (2.g) follows from part (2.d) and Lemma \ref{F6}. 

For part (2.h), we have $f_1\in R_{h_0}^{\ast}$ by part (2.g). From part (2.d) we get that $h_0\mid f_1$ in $R$, i.e., there exists $g_1\in R$ such that $f_1=h_0g_1$. As $g_1=f_1h_0^{-1}\in R_{h_0}^{\ast}$, the existence of $l\in\mathbb N$ follows. Thus part (2.h) holds. 

Part (2.i) follows from the definitions and parts (2.f) and (2.h). 

Part (2.l) is proved in the same way as Theorem \ref{T9}(3), the only difference being that we cannot conclude that the projective $A_e$-module $C_e/A_e$ is free.

As $Y_e$ is a complete intersection by part (2.i), part (2.j) follows from parts (2.c) and \cite{Gro6}, Exp.\ X, Thm.\ 3.4(ii).

If $\chr(K)\neq 2$, then $\phi$ is non-\'etale by part (2.f) and hence $\psi_e$ is non-\'etale. So part (2.k) holds.
 
Part (3.a) holds as $\Reg(W_e)$ maps to $\Reg(X_e)$ (resp.\ $\Reg(Y_e)$) by \cite{Ma}, Ch.\ 8, Thm.\ 23.7(i). 

For part (3.b), if $Z$ is not an irreducible component of $N_e^{\textup{n}-\et}$, let $i\in\{1,2\}$ be such that we have a disjoint union decomposition $\psi_e^{-1}(Z)\cap \bigl(X_e^{\et}\setminus\Imm(\imath_e)\bigr)=\sqcup_{j=1}^i Y_j$; if $i=1$, then the finite morphism $Y_1\rightarrow Z$ is not birational by our hypotheses.  

Let $Z^{\n}$ be the normalization of $Z$. As $Z$ is rational and has only one point at infinity it suffices to show that $Z$ is regular. If $P\in\Sing(Z)$, then the inverse image $\psi_e^{-1}(P)$ is the spectrum of a $K$-algebra $A_P$ which as a $K$-vector space has dimension $3$. 

Let $Y$ be an irreducible component of the reduced inverse image $\bigl(\psi_e^{-1}(Z)\bigr)_{\red}$. As $Y\cong\mathbb A^1_K$ by Corollary \ref{C2.6}(1), one of the following three distinct conditions hold: (i) either the morphism $Y\rightarrow Z$ factors through a flat morphism $Y\rightarrow Z^{\n}$ of degree $2$ or (ii) the inverse image $\psi_e^{-1}(Z)$ has multiplicity $2$ at $Y$ or (iii) $Z$ is not an irreducible component of $N_e^{\textup{n}-\et}$ and we have $i=2$ with the finite morphism $Y_j\rightarrow Z$ inducing an isomorphism $Y_j\rightarrow Z^{\n}$ for each $j\in\{1,2\}$. If the condition (ii) holds and $I_Y$ is the ideal of $A_e$ that defines $Y$, then $A_e/I_Y\cong K[t]$ is a principal ideal domain and the $A_e/I_Y$-module $Y_Y/I_Y^2$ has rank $1$; thus there exists a closed subscheme $Y^+$ of $X_e$ which contains $Y$ and is isomorphic to $Y\times_{\Spec K} \Spec \bigl(K[t]/(t^2)\bigr)$. If the conditions (i) holds, then let $Y^+:=Y$. If the conditions (iii) holds, then let $Y^+:=Y_1\sqcup Y_2$. Let $B_P$ be the quotient of $A_P$ such that $\Spec B_P$ defines the inverse image of $P$ via the morphism $Y^+\rightarrow Z$ induced by $\psi_e$.

The inverse image of $P$ via the morphism $Z^{\n}\rightarrow Z$ is the spectrum of a $K$-algebra which as a $K$-vector space has dimension at least $2$. It follows that 
$$\dim_K(B_P)\ge 2\cdot 2=4>3=\dim_K(A_P)\ge\dim_K(B_P)$$ regardless of which one of the three conditions (i) to (iii) holds, a contradiction. Thus $Z$ is regular and part (3.b) holds.

Part (3.c) holds by Criterion \ref{CRI1}(1) as its hypotheses hold based on part (3.b) and Corollary \ref{C2.6}(1).

Part (3.d) follows from part (3.b) and Theorem \ref{T8}(3) applied to $h_0$ viewed as an element of $K_{\t}[x_1,x_2,x_3]$.

For parts (3.e) and (3.f), we write $\Irr(X_e\setminus X_e^{\et})=\{Y_i|i\in\llbracket1,\rho(e)-\rho_{\et}(e)\rrbracket\}$. For $i\in\llbracket1,\rho(e)-\rho_{\et}(e)\rrbracket$ let $Z_i:=\psi_e(Y_i)$; it is a closed curve of $\mathbb A^2_{K,\t}$ isomorphic to $\mathbb A^1_K$ by part (3.b) and the reduced inverse image $\bigl(\psi_e^{-1}(Z_i)\bigr)_{\red}$ is the union of $Y_i$ and a closed curve $C_i$ of $X_e$ whose generic point belongs to $\Imm(\imath_e)$. For each $i\in \llbracket1,\rho(e)-\rho_{\et}(e)\rrbracket$ the finite surjective morphism $C_i\rightarrow Z_i$ is birational and hence, as $Z_i\cong\mathbb A^1_K$, it is an isomorphism; so $C_i\cong\mathbb A^1_K$. 

We check that for each pair $(i,j)\in\llbracket1,\rho(e)-\rho_{\et}(e)\rrbracket^2$ with $i<j$ the assumption that there exists a point $Q_{ij}\in (Z_i\cap Z_j)(K)$ leads to a contradiction. The inverse image $\psi_e^{-1}(P_{ij})$ is the spectrum of a $K$-algebra $A_{ij}$ which as a $K$-vector space has dimension $3$. As $Y_i\cap Y_j=\emptyset$ by Corollary \ref{C2.6}(1), the $K$-algebra $A_{ij}$ has a quotient $B_i\times B_j$, where $\Spec B_i=Y_i\cap \Spec A_{ij}$ and $\Spec B_j=Y_j\cap \Spec A_{ij}$. As $Y_i\sqcup Y_j\subset X_e\setminus X_e^{\et}$, both $K$-vector spaces $B_i$ and $B_j$ have dimension at least $2$. Thus $\dim_K(A_{ij})\ge 2+2>3=\dim_K(A_{ij})$, a contradiction.

Hence $\psi_e(X_e\setminus X_e^{\et})=\sqcup_{i=1}^{\rho(e)-\rho_{\et}(e)} Z_i$ has $\rho(e)-\rho_{\et}(e)$ connected components which are irreducible. So part (3.e) holds. 

For part (3.f), based on part (2.j) it suffices to show that $Y_e$ is regular. But this follows from part (2.f) and the fact that the zero locus $h_0=0$ is a disjoint union of regular curves isomorphic to $\mathbb A^1_K$ by part (3.b). So part (3.f) holds.

For part (3.g), the proof of part (2.b) gives a union $\psi_e^{-1}(Z_1)=C_1\cup Y_1$ with the generic point of $C_1$ belonging to $\Imm(\imath_e)$ and with $Y_1\subset X_e\setminus X_e^{et}$. Based on this and the proof of part (3.e), we only have to prove that for each $Q\in C_1(K)$ we have $Q\notin Y_1$. As $Q\in C_1(K)$, there exists a unique point $P\in W_e(K)$ that maps to $Q$. As the morphism $W_e\rightarrow Y_e$ is \'etale at $P$ by hypothesis, it follows that the morphism $X_e\rightarrow\mathbb A^2_{K,\t}$ is \'etale at $Q$ and hence $Q\notin Y_1(K)$. So part (3.g) holds.
 
For part (4), as $h_0\mid h$, the morphism $R_h\to (A_e)_h$ is integral and \'etale. Hence $(A_e)_h$ is a projective $R_h$-module. From this and Theorem \ref{T9}(1) we get that $(A_e)_h/R_h$ is also a projective $R_h$-module and therefore, as $h$ has the freeness property, it is free. So part (4) holds.

To prove part (5), we consider an $R_h$-basis $(1,w_1,w_2)$ of $(A_e)_h$. The $A$-algebra $(A_e)_h\otimes_{R_h} A$ is generated by $1$ element iff there exists $(u,v)\in A^2$ such that $A$ is generated by $w_1\otimes u+w_2\otimes v$ and iff, by writing 
$$(w_1\otimes u+w_2\otimes v)^2=1\otimes t_2+w_1\otimes u_2+w_2\otimes v_2$$ with $(t_2,u_2,v_2)\in A^3$, the matrix $\triangle_{u,v}:=\begin{bmatrix} 
u & v \\
u_2 & v_2 \\ 
\end{bmatrix}$ is invertible, i.e., is in $\GL_2(A)$. As $u_2$ and $v_2$ are homogeneous polynomials in $u$ and $v$ of degree $2$ with coefficients in $R_h$ that do not depend on $A$, there exists $(\alpha,\beta,\gamma,\delta)\in R_h^4$ such that by defining $L(x,y):=\alpha x^3+\beta x^2y+\gamma xy^2+\delta y^3\in R_h[x,y]$ we have $\det(\triangle_{u,v})=L(u,v)$. As $(A_e)_h$ modulo each maximal ideal of $R_h$ is generated by $1$ element, $(\alpha,\beta,\gamma,\delta)\in R_h^4$ is unimodular. We conclude that part (5) holds.

For part (6), the $(A_e)_h$-algebra $(A_e)_h\otimes_{R_h} (A_e)_h$ is a product $(A_e)_h\times A$ with $A$ as an $(A_e)_h$-algebra which as an $(A_e)_h$-module is free of rank $2$. As the short exact sequence $0\rightarrow (A_e)_h\rightarrow A\rightarrow A/(A_e)_h\rightarrow 0$ of $(A_e)_h$-modules splits by \cite{Ho}, Thm.\ 2, by taking determinants we get that the $(A_e)_h$-module $A/(A_e)_h$ is free of rank $1$. So the \'etale $(A_e)_h$-algebra $(A_e)_h\otimes_{R_h} (A_e)_h$ is generated by each $(0,w_0)\in (A_e)_h\times A$ with $A=(A_e)_hw_0\oplus (A_e)_hw_0^2$; such an element $w_0\in A$ exists as $\chr(K)\neq 2$ and the \'etale $(A_e)_h$-algebra $A$ is isomorphic to $(A_e)_h[z]/(z^2-\tilde h)$ for some unit $\tilde h\in (A_e)_h^{\ast}$. Thus there exists $(v_1,v_2)\in (A_e)_h^2$ such that $L(v_1,v_2)\in (A_e)_h^{\ast}=R_h^{\ast}$ by part (5). So part (6) holds. 

To prove part (7), we note that, as $\deg(e)=3$, Condition ($\ast^e_{h}$) holds iff Condition ($\ast^e_{h_2}$) holds. This is so as $e^{-1}(Y)$ is irreducible for each irreducible component $Y$ of the zero locus $h_0=0$ (see the proof of part (2.c)). So, if Condition ($\ast^e_{h_2}$) holds, then the homomorphism $R_h\rightarrow (A_e)_h$ induces an isomorphism $R_h^{\ast}\rightarrow (A_e)_h^{\ast}$. 

Let $x_h\in A_e$ be such that its image in $(A_e)_{h}$ generates the $R_h$-algebra $(A_e)_h$. Let $\pi_h:\mathbb A^n_{K,\s}\rightarrow\mathbb A^1_K$ be defined by $x_h$. From the prior paragraph we get that $e$ is of weak type $1$ by Definition \ref{D11} applied to $\pi_h\times e$. From this and Proposition \ref{PR16} we get that $\chr(K)$ divides $3$ and hence $\chr(K)=3$. Thus part (7) holds.

For parts (8.a) and (8.b), as $R_h\cong K[x_1,\ldots,x_{n+1}]/(x_{n+1}h-1)$ it follows that $\ciedim(\Spec R_h)\le n+1$. For $\diamond\in\{A,C\}$ the morphism $\Spec (\diamond_e)_h\rightarrow\Spec R_h$ is \'etale. Hence $T_{\Spec (\diamond_e)_h}$ is trivial, and thus $\ciedim\bigl(\Spec (\diamond_e)_h\bigr)\le 2n+2$ by Theorem \ref{T3}(2) and (2.f). 

Based on this, for $\chr(K)\neq 2$ (resp.\ $\chr(K)=2$), the fact that the inequality $\ciedim\bigl(\Spec (B_e)_h\bigr)\le n+1$ (resp.\ $\ciedim\bigl(\Spec (B_e)_h\bigr)\le n+2$) holds follows from part (2.f) and the $K$-algebra isomorphism $(B_e)_h\cong K[x_1,\ldots,x_{n+1}]/(x_{n+1}^2h-1)$ (resp.\ from part (2.g)). In particular, part (8.b) holds.

If $\chr(K)\neq 2$, then from $h_0\notin K^{\ast}$ by part (2.f) and $h_0\mid h$ we get that $h\notin K^{\ast}$ and thus $K^{\ast}\subsetneq R_h^*\subset (B_e)_h$. So $\edim(\Spec R_h)>n$ and $\ciedim\bigl(\Spec (B_e)_h\bigr)>n$; thus $\edim(\Spec R_h)=\ciedim\bigl(\Spec (B_e)_h\bigr)=n+1$. So part (8.a) also holds.

For part (9), note that the ideal 
$$I=\bigl(y_1^2-\ell_1(y_1,y_2),y_1y_2-\ell_{12}(y_1,y_2),y_2^2-\ell_2(y_1,y_2)\bigr)$$ 
of $R_h[y_1,y_2]$ is defined by a triple $(\ell_1,\ell_{12},\ell_2)\in R_h[y_1,y_2]^3$ of polynomials of degree at most $1$. As $A_e$ is an integral domain and $\{y_1,y_2\}\cap K^{\ast}=\emptyset$, it follows that $\deg(\ell_1)=\deg(\ell_2)=1$ and $\deg(\ell_{12})\in\{0,1\}$. As locally in the Zariski topology of $\Spec R_h$ the $R_h$-algebra $(A_e)_h$ is generated by $1$ element (e.g., this follows from part (5)), locally in the Zariski topology of $\Spec R_h$ the ideal $I$ is generated by $2$ elements and hence, as $R_h[y_1,y_2]$ and $(A_e)_h$ are regular, the homological dimension of $I$ is at most $1$. Based on this, if $h$ has the freeness property, then from \cite{Ser1}, Part II, Sect.\ 4, Cor.\ we get that the ideal $I$ is generated by $2$ elements iff the $R_h[y_1,y_2]$-module $\Ext^1_{R_h[y_1,y_2]}(I,R_h[y_1,y_2])$ is cyclic. 

If the ideal $I$ is generated by $2$ elements $f$ and $g$ and Condition ($S_3$) holds, then we have $K$-algebra isomorphisms $(A_e)_h\cong K[x_1,\ldots,x_{n+3}]/(x_{n+1}h-1,f,g)$ and $(C_e)_h\cong K[x_1,\ldots,x_{n+3}]/(x_{n+1}^2h-1,f,g)$ and hence we have inequalities $\ciedim\bigl(\Spec (A_e)_h\bigr)\le n+3$ and $\ciedim\bigl(\Spec (C_e)_h\bigr)\le n+3$. So part (9) holds.

For part (10), as $\chr(K)\neq 3$, there exists $w\in K_0$ that generates $K_0$ over $\Frac(B_e)$ and such that $u:=w^3\in\Frac(B_e)$ by Kummer theory; we can assume that $w\in C_e$ and that $u\in B_e$ is cubic free, i.e., for each $\star\in B_e\setminus B_e^{\ast}$ we have $u\notin \star^3B_e$. As $\Spec (C_e)_h\rightarrow\Spec (B_e)_h$ is a finite Galois cover between regular varieties of degree $3$, for each point $P\in\Spec (B_e)_h(K)$, there exists an affine open subvariety $U_P$ of $\Spec (B_e)_h$ which contains $P$ and an element $u_P\in\Frac(B_e)$ such that $uu_P^3$ is an invertible global function on $U_P$. The set $\{(U_P,uu_P^3)|P\in\Spec (B_e)_h(K)\}$ is a Cartier divisor that defines $\langle M_h\rangle$. 

For part (10.a), as $M_h\cong (B_e)_h$, we can assume that $u\in (B_e)_h^{\ast}$, and hence $(C_e)_h\cong (B_e)_h[t]/(t^3-u)$. If $\chr(K)=2$, then from this and part (2.g) we get that $\ciedim\bigl(\Spec (C_e)_h\bigr)\le n+3$. 

In this and the next two paragraphs we assume that $\chr(K)\neq 2$. We have $\ciedim\bigl(\Spec (C_e)_h\bigr)\le n+2$ as $(C_e)_h\cong K[x_1,\ldots,x_{n+2}]/(x_{n+1}^2h-1,x_{n+2}^3-u)$. 

As $\chr(K)\neq 2$, we view $(B_e)_h\cong R_h[z]/(z^2-h_0)$ (see part (2.f)) as an identification and we denote the element $z+(z^2-h_0)\in (B_e)_h$ simply by $z$. We write $u=f+zg$ with $(f,g)\in R_h^2$ by part (2.f). Let $\tau\in G$ be the unique element of order $2$ that fixes $A_e$. Then $\tau(z)=-z$, $x:=w\tau(w)\in (A_e)_h$ and $v:=u\tau(u)=f^2-hg^2\in R_h^{\ast}$. As $x^3=v$ and $\Frac(R)\rightarrow\Frac(A_e)$ is a non-Galois field extension, it follows that $x\in\Frac(R)\cap (A_e)_h=R_h$. As $x^3=v\in R_h^{\ast}$, we have $x\in R_h^{\ast}$. So there exists $q\in\mathbb N$ such that $x|h^q$. 

We have $\tau(w)=\frac{x}{w}$ and $\tau(u)=\frac{v}{u}$. As $u\in (B_e)_h^{\ast}$, $(1,z,w,wz,w^2,w^2z)$ is an $R_h$-basis of $(C_e)_h$. As $\tau$ interchanges $R_hw\oplus R_hwz$ and $R_hw^2\oplus R_hw^2z$ and we have $\tau(R_h\oplus R_hz)=R_h\oplus R_hz$ and $(A_e)_h:=\{x\in (C_e)_h|\tau(x)=x\}$, it follows that $\bigl(1,w+\tau(w),wz+\tau(wz)\bigr)$ is an $R_h$-basis of $(A_e)_h$ and Condition ($\dag^e_h$) holds. We apply parts (5) and (9) to $(w_1,w_2):=\bigl(w+\tau(w),wz+\tau(wz)\bigr)$. As $w:=\frac{w_1}{2}+\frac{w_2z}{2h}$ generates the $(B_e)_h$-algebra $(C_e)_h$, we have $L(\frac{1}{2},\frac{z}{2h})\in (B_e)_h^{\ast}$. As $(C_e)_h\cong (B_e)_h[t]/(t^3-u)$ under the isomorphism that maps $w$ to $t+(t^3-u)$, the $(B_e)_h[y_1,y_2]$-module 
$$\Ext^1_{(B_e)_h[y_1,y_2]}(I\otimes_{R_h} \otimes (B_e)_h,(B_e)_h[y_1,y_2]\bigr)=\Ext^1_{R_h[y_1,y_2]}(I,R_h[y_1,y_2])\otimes_{R_h} (B_e)_h$$
$$=\Ext^1_{R_h[y_1,y_2]}(I,R_h[y_1,y_2])\otimes_{R_h[y_1,y_2]} (B_e)_h[y_1,y_2]$$ 
is cyclic isomorphic to $(B_e)_h[y_1,y_2]/I_z$ for some ideal $I_z$ of $(B_e)_h[y_1,y_2]$. As there exists a $(B_e)_h[y_1,y_2]$-linear isomorphism $(B_e)_h[y_1,y_2]/I_z\cong (B_e)_h[y_1,y_2]/\tau(I_z)$, we have $\tau(I_z)=I_z$ and Galois descend implies that the following $R_h[y_1,y_2]$-module $\Ext^1_{R_h[y_1,y_2]}(I,R_h[y_1,y_2])\cong R_h[y_1,y_2]/\{y\in I_z|\tau(y)=y\}$ is cyclic. Therefore $\ciedim\bigl(\Spec (A_e)_h\bigr)\le n+3$ by part (9), and thus part (10.a) holds. 

Part (10.b) follows from part (2.i) and \cite{Gro6}, Exp.\ XI, Thm.\ 3.13(ii). 

Part (10.c) follows from part (2.j): we can take $\langle e\rangle\in H^1_{\et}\bigl(\Spec B_e,(\mathbb Z/3\mathbb Z)_K\bigr)$ to define the Galois cover $\Spec C_e\rightarrow\Spec B_e$ and $\langle M\rangle\in H^1_{\et}\bigl(\Spec B_e,\mathbb G_{\m,K}\bigr)$ to be the image of $\langle e\rangle$ via the monomorphism $(\mathbb Z/3\mathbb Z)_K\rightarrow\mathbb G_{\m,K}$. 

%For part (10.d) it suffices to show that the assumption that the class $\langle M_{h_0}\rangle$ is trivial leads to a contradiction. So part (10.d) holds. 

Part (11.a) is a particular case of part (4). 

The first part of (11.b) follows from part (11.a). As $M_h$ is generated by two elements by part (11.a), there exists a projective $(B_e)_h$-module $M_h^{\vee}$ such that we have a $(B_e)_h$-linear isomorphism $M_h\oplus M_h^{\vee}\cong (B_e)_h^2$. Taking determinants, we get that $M_h^{\vee}$ is the dual of $M_h$. As $\langle M_h\rangle$ is a class of order 1 or $3$, $M_h^{\otimes 3}\cong (B_e)_h$, hence $M_h^{\vee}\cong M_h^{\otimes 2}$ as $(B_e)_h$-modules and part (11.b) holds. 

For part (11.c), we identify $M_h=R_h^2$ as $R_h$-modules and assume $\chr(K)\neq 2$. As $(B_e)_h=R_h[z]/(z^2-h)$, we have $M_h\cong (B_e)_h$ as $(B_e)_h$-modules iff there $w_h\in M_h$ such that $(w_h,zw_h)$ is an $R_h$-basis of $M_h$. Let $(w_1,w_2)$ be the standard $R_h$-basis of $R_h^2$. Let the matrix $\triangle_z=\begin{bmatrix} 
\alpha_z & \beta_z \\
\gamma_z & \delta_z \\ 
\end{bmatrix}\in\GL_2(R_h)$ define the action of $z$ on $M_h$. We have $\triangle_z^2=\begin{bmatrix} 
h & 0 \\
0 & h \\ 
\end{bmatrix}.$ As $h$ is not a square in $R_h$, a simple computation gives that $\alpha_z=-\delta_z$ and $\delta_z^2+\beta_z\gamma_z=h$. For $(u,v)\in R_h^2$, by taking $w_h:=uw_1+vw_2$, $(w_h,zw_h)$ is an $R_h$-basis of $M_h$ iff $\begin{bmatrix} 
u & v \\
-\delta_zu+\beta_zv & \gamma_zu+\delta_zv \\ 
\end{bmatrix}\in\GL_2(R_h)$ and iff $L_2(u,v)\in R_h^{\ast}$. If the pair $(\beta_z,\gamma_z)\in R_h^2$ is unimodular, then the cocycle part follows from the identities $(B_e)_{h\gamma_z}w_1=(M_h)_{\gamma_z}$, $(B_e)_{h\beta_z}w_2=(M_h)_{\beta_z}$ and $w_2=\frac{z+\delta_z}{\gamma_z}w_1$. So part (11.c) holds.

For part (11.d), we identify $M_h=R_h^2$ as $R_h$-modules and assume $\chr(K)=2$. As $(B_e)_h=R_h[z]/(z^2+z+f_0)$, we have $M_h\cong (B_e)_h$ as $(B_e)_h$-modules iff there $w_h\in M_h$ such that $(w_h,zw_h)$ is an $R_h$-basis of $M_h$. Let $(w_1,w_2)$ be the standard $R_h$-basis of $R_h^2$. Let the matrix $\triangle_z=\begin{bmatrix} 
\alpha_z & \beta_z \\
\gamma_z & \delta_z \\ 
\end{bmatrix}$ with entries in $R_h$ define the action of $z$ on $M_h$. We have $\triangle_z^2+\triangle_z=\begin{bmatrix} 
f_0 & 0 \\
0 & f_0 \\ 
\end{bmatrix}$. As there exists no $f_1\in R_h$ such that $f_0=f_1^2+f_1$, a simple computation gives that $\alpha_z=\delta_z+1$ and $\delta_z^2+\delta_z+\beta_z\gamma_z=f_0$. For $(u,v)\in R_h^2$, by taking $w_h:=uw_1+vw_2$, $(w_h,zw_h)$ is an $R_h$-basis of $M_h$ iff $\begin{bmatrix} 
u & v \\
\delta_zu+u+\beta_zv & \gamma_zu+\delta_zv \\ 
\end{bmatrix}\in\GL_2(R_h)$ and iff $L_2(u,v)\in R_h^{\ast}$.

The proof of part (12) is entirely the same as of part (11) once we note that the $R$-module $M$ is free by Quillen--Suslin Theorem.\end{proof}

\begin{remark}\normalfont\label{R27.9}
We refer to Theorem \ref{T10} when $\chr(K)\neq 3$ and Condition ($S_3$) holds. Let $g\in\Frac(R)$ be such that $A_e$ is the normalization of $R$ in the field extension $\Frac(R)[w]/(w^3-3w-g)$ of $\Frac(R)$ by \cite{MW}, Thm.\ 2.1(c). So $\Frac(B_e)$ is the normalization of $R$ in the field extension $\Frac(R)[w]/(w^2+gw+1)$ of $\Frac(R)$ by \cite{MW}, Thm.\ 3.1. We consider the norm maps $\N:\Frac(B_e)\rightarrow\Frac(R)$ and $\N_1:\Frac(C_e)\rightarrow\Frac(B_e)$.

\medskip
{\bf (1)} Suppose that $\chr(K)\neq 2$. The field extensions $\Frac(R)[w]/(w^2+1-\frac{g^2}{4})$ and $\Frac(B_e)=\Frac(R)[z]/(z^2-h_0)$ of $\Frac(R)$ are isomorphic by part (2.f). So there exists $h_1\in \Frac(R)$ such that $\frac{g^2}{4}-1=h_0h_1^2$. We write $h_1=\frac{f_1}{g_1}$ with $(h_1,g_1)\in R^2$, $g_1\neq 0$, and $h_1$ and $g_1$ relatively prime. Denoting $g_2:=gg_1\in\Frac(R)$, we get that $$g_2^2=4g_1^2+4h_0f_1^2.$$ 
As $g_2^2\in R$ and $R$ is normal, we have $g_2\in R$. 
Denoting also by $z$ the element $z+(z^2-h_0)\in B_e=R[z]/(z^2-h_0)$, $u_0:=-\frac{g}{2}+h_1z\in\Frac(B_e)$ is a root of the polynomial $w^2+gw+1\in R[w]$ and we have $\N(u_0)=1$. From \cite{MW}, Thm.\ 3.1 we get that $C_e$ is the normalization of $B_e$ in the field extension $\Frac(B_e)[w]/(w^3-u_0)$ of $\Frac(B_e)$. Let $v_0\in\Frac(B_e)$ be such that $v_0^3=u_0$. Let $v:=g_1v_0\in\Frac(C_e)$. As $v^3=g_1^2u\in B_e$, we have $v\in C_e$. So $C_e$ is also the normalization of the $B_e$-algebra $C:=B_e[w]/(w^3-v)$ and $A_e$ is the normalization of the $R$-algebra $$A:=R[w]/(w^3-3g_1^2w-g_1^2g_2).$$%Let $u:=g_1u_0=-\frac{g_2}{2}+f_1z\in B_e$. We have $N(u)=g_1^2$. 

\smallskip
{\bf (2)} Suppose that $\chr(K)=2$. The field extensions $\Frac(R)[w]/(w^2+w+\frac{1}{g^2})$ and $\Frac(B_e)=\Frac(R)[z]/(z^2+z+f_0)$ of $\Frac(R)$ are isomorphic by part (2.g). So there exists $h_1\in \Frac(R)$ such that $\frac{1}{g^2}=f_0+h_1+h_1^2$. Denoting also by $z$ the element $z+(z^2+z+f_0)\in B_e=R[z]/(z^2+z+f_0)$, $u_0:=gh_1+gz\in\Frac(B_e)$ is a root of the polynomial $w^2+gw+1\in R[w]$ and we have $\N(u_0)=1$. From \cite{MW}, Thm.\ 3.1 we get that $C_e$ is the normalization of $B_e$ in the field extension $\Frac(B_e)[w]/(w^3-u_0)$. Let $v_0\in\Frac(B_e)$ be such that $v_0^3=u_0$.

\smallskip
{\bf (3)} Let $\alpha\in K^{\ast}$ of order $3$. Referring to both parts (1) and (2), the conjugates of $v_0$ over $\Frac(B_e)$ are $v_0$, $v_0\alpha$, and $v_0\alpha^2$. So $\N_1(v_0)=v_0^3=u_0$. From this and Theorem \ref{T10} (2.c) we get that the multiplicities of $z_0$ at all points of $W_e$ that map to the same point of $Y_e$ of codimension $1$ are equal. In particular, if $\tilde h\in R$ is such that $u_0\in (B_e)_{\tilde h}$, then from $\N(u_0)=1$ we get that $u_0\in (B_e)_{\tilde h}^{\ast}$ and from this and $\N_1(v_0)=u_0$ we get that $v_0\in (C_e)_{\tilde h}^{\ast}$. For instance, if $\chr(K)\neq 2$, then we can take $\tilde h:=g_1h_0$ by part (1).
\end{remark}

For the first general applications of Theorem \ref{T10} when $n\ge 3$ see \url{https://people.math.binghamton.edu/adrian/AI.pdf}. For $n=2$ we include the following general application of Theorem \ref{T10}.

\begin{corollary}\label{C26}
Suppose that $\chr(K)\notin\{2,3\}$. Let $e\in\EE_2(K)$ be such that $\deg(e)=3$. Then the following properties hold for each irreducible component $Y$ of $N_e^{\textup{n}-\et}$.

\medskip
{\bf (1)} There is no $a\in\GA_2(K)$ such that $a(Y)$ is the zero locus $x_1=0$ in $\mathbb A^2_K$.\footnote{Such a $Y$ is non-rectifiable in the terminology of \cite{Le} and is a wild line in the terminology of \cite{Da}, Sect.\ 2.}

\smallskip
{\bf (2)} We have $\chr(K)\neq 0$, i.e., $\chr(K)\ge 5$.

\smallskip
{\bf (3)} There exists no $l\in\mathbb N^{\ast}$ such that $Y$ is the pullback of the zero locus $x_1=0$ in $\mathbb A^2_K$ via an endomorphism $d_l\in\End_2(K)$ which is a composite of $l$ elementary purely inseparable endomorphisms of $\mathbb A^2_K$ of degree $\chr(K)$, i.e., of endomorphisms of the form $a_1\e(x_1^{\chr(K)},x_2)a_2$ with $(a_1,a_2)\in\GA_2(K)^2$.\footnote{The Conjecture on the Classification of Lines, see the equivalent conjectures \cite{Da}, Conjs.\ 1 and 2 where the first one is a correction of  \cite{Moh2}, Conj.\ 1, is also equivalent to the statement that for each curve in $\mathbb A^2_K$ isomorphic to $\mathbb A^1_K$ such an $l$ exists.}
\end{corollary}

\begin{proof}
From Theorem \ref{T10}(1) we get that Condition ($S_3$) holds; so we can use the notation of Theorem \ref{T10} for $n=2$. So $R:=K_{\t}[x_1,x_2]$, $(h_0,h_1)\in R^2$, and $Y_e=\Spec B_e$ and $W_e=\Spec C_e$ are as in Theorem \ref{T10}(2.b). We have $Y_e=\Spec R[z]/(z^2-h_0
)$ and the finite morphism $W_e\rightarrow Y_e$ is a Galois cover by Theorem \ref{T10}(2.f) and (2.j).

Let $j:=\nu_{\textup{n}-\et}(e)$. We have $j=\rho(e)-\rho_{et}(e)$ and the irreducible components $Y=Y_1,Y_2,\ldots,Y_j$ of $N_e^{\textup{n}-\et}=\Spec R/Rh_0$ are also connected components of $N_e^{\textup{n}-\et}$ and are isomorphic to $\mathbb A^1_K$ by Theorem \ref{T10}(3.e). 

For $i\in \llbracket1,j\rrbracket$ let $f_i\in R$ be an irreducible element such that $Y_i$ is the zero locus $f_i=0$; it is uniquely determined up to a scalar multiple by an element of $K^{\ast}$. 

For $i\in\llbracket2,j\rrbracket$, as $f_1+(f_i)\in K[x_1,x_2]/(f_i)\cong K[t]$ is a unit, it follows that there exists a unique $\alpha_i\in K^{\ast}$ such that $f_1+(f_i)=\alpha_i+(f_i)\in K[x_1,x_2]/(f_i)$. Thus $f_i\mid f_1-\alpha_i$ and hence $\deg(f_i)\mid \deg(f_1)$. By symmetry, $\deg(f_1)\mid\deg(f_i)$. So $\deg(f_i)=\deg(f_1)$ and up to scalars in $K^{\ast}$ we can assume that $f_i=f_1-\alpha_i$. 

Let $\alpha_1:=0$. The elements $\alpha_1$ to $\alpha_j$ in $K$ are distinct and we have an identity $h_0=\prod_{i=0}^j (f_1-\alpha_i)$.

For part (1) we can assume that $f_1=x_1$. Let $C=\Spec B_{e,0}$ be the smooth affine curve over $\Spec K$ which is the zero locus $z^2-h_0=0$ in $\Spec K_{\t}[z,x_1]$. We have $Y_e=C\times_{\Spec K} \mathbb A^1_K$. 

If $j=1$, then $C\cong\mathbb A^1_K$ and $Y_e\cong\mathbb A^2_K$; thus, as $W_e\rightarrow Y_e$ is Galois, we reached a contradiction to Proposition \ref{PR15}. So we can assume that $j\ge 2$. 

The first projection $\pi_1:Y_e=C\times_{\Spec K} \mathbb A^1_K\rightarrow C$ induces an isomorphism $\pi_1^*:\Pic(C)\rightarrow\Pic(Y_e)$ between Picard groups (e.g., see \cite{GH}, Thm.\ 1.6) and we have $B_{e,0}^{\ast}=B_{e,0}[t]^{\ast}=B_e^{\ast}$. This implies that $\pi_1^*:H^1_{\et}\bigl(C,(\mathbb Z/3\mathbb Z)_K\bigr)\rightarrow H^1_{\et}\bigl(Y_e,(\mathbb Z/3\mathbb Z)_K\bigr)$ is an isomorphism. Thus there exists a finite Galois cover $\phi_2:C_2\rightarrow C$ of degree $3$ that defines a class $\langle e_C\rangle\in H^1_{\et}\bigl(C,(\mathbb Z/3\mathbb Z)_K\bigr)$ that maps to the class $\langle e\rangle\in H^1_{\et}\bigl(Y_e,(\mathbb Z/3\mathbb Z)_K\bigr)$ of Theorem \ref{T10}(10.c). 

Hence we can identify $W_e=C_2\times_{\Spec K} \mathbb A^1_K$ in such away that the morphism $W_e\rightarrow Y_e$ is $\phi_2\times 1_{\mathbb A^1_K}$. It follows that there exists a finite flat morphism $\C_2\rightarrow C_1$ of degree $2$ such that $X_e=C_1\times_{\Spec K} \mathbb A^1_K$. As we have a dominant morphism $\mathbb A^2_K\rightarrow X_e$ (e.g., the open embedding $\imath_e:\mathbb A^2_K\rightarrow X_e$ is dominant), $C_1$ is a connected smooth affine unirational curve over $\Spec K$ whose group of multiplicative units is $K^{\ast}$; from L\"uroth Theorem we get that $C_1$ is rational. So $C_1\cong\mathbb A^1_K$ and thus $X_e\cong\mathbb A^2_K$. This contradicts the facts that $X_e\setminus X_e^{\et}$ has $j$ irreducible components with $\mathbb A^2_K\cong\Imm(\imath)\subset X_e^{\et}$. 

So our assumption is false and part (1) holds.

For part (2), from Abhyankar--Moh Theorem (see \cite{AM}, Epimorphism Thm.\ (1.2)) we get that $\chr(K)\neq 0$. Thus $\chr(K)\ge 5$ and part (2) holds.\footnote{Part (2) also follows from \cite{O}, Thm.\ 1.1 which is over $\mathbb C$.}

For part (3), we show that the assumption that there exists $l\in\mathbb N^{\ast}$ and a composite $d_l\in\End_2(K)$ of $l$ elementary purely inseparable endomorphisms of $\mathbb A^2_K$ such that $Y=d_l^{-1}(\{0\}\times\mathbb A^1_K)$ leads to a contradiction. As $f_i=f_1+\alpha_i$ for each $i\in\llbracket2,j\rrbracket$, we have an identity $N_e^{\textup{n}-\et}=d_l^{-1}\bigl(d_l(N_e^{\textup{n}-\et})\bigr)$. Based on this and the fact that $d_l$ induces an isomorphisms between fundamental groups, the diagram of Theorem \ref{T10}(2.b) is the pullback via $d_l$ of a similar diagram. So as in the proof of part (1) we reach a contradiction. For instance, if $j\ge 2$, then the morphism $X_e\rightarrow\mathbb A^2_K$ induced by $\psi_e$ is the pullback via $d_l$ of an endomorphisms of $\mathbb A^2_K$ because as in the proof of part (1) we argue that its source is $\mathbb A^2_K$; so there exists a purely inseparable finite flat morphism $\mathbb A^2_K\rightarrow X_e$ and this contradicts the facts that $X_e\setminus X_e^{\et}$ has $j$ irreducible components with $\mathbb A^2_K\cong\Imm(\imath)\subset X_e^{\et}$. So part (3) holds. 
\end{proof}

\section{On geometric degrees $3$, part II: complements for $p\in\{2,3\}$}\label{S31}

This section gathers several special features of \'etale endomorphisms $e\in\EE_2(k)$ with $\deg(e)=3$ for $p\in\{2,3\}$. They show that many of the hypotheses of Theorem \ref{T10} are indeed required.

\begin{example}\normalfont\label{EX49}
Suppose that $p=3$. Let the pair $(d,e)\in\EE_2(k)^2$ be such that $\deg(d)=\deg(e)=3$, $\psi_d$ is \'etale, and $\psi_e$ is non-\'etale (for concrete examples and extra properties see Theorem \ref{T1}(2) and (3) and Example \ref{EX11}). 

\medskip
{\bf (1)} Condition ($S_3$) does not hold for $d$ by Theorem \ref{T10}(2.k). So the morphism $X_d\rightarrow\mathbb A^2_{k,\t}$ is Galois and an Artin--Schreier cover by Lemma \ref{F6}. Thus $\edim(X_e)$ is $2$ if $\imath_e$ is an isomorphism and is $3$ if $\imath_e$ is not an isomorphism. The invariant $\varphi(d)\in\mathbb N$ can be an arbitrary natural number by Equation (\ref{EQ3}). The set of possible invariants $\rho(d)=\rho_{\et}(d)$ contains 2$\mathbb N$ by Corollary \ref{C8}(2).

\smallskip
{\bf (2)} If Condition ($S_3$) does not hold for $e$, then the finite morphism $X_e\rightarrow\mathbb A^2_{k,\t}$ is \'etale, a contradiction to $\psi_e$ being non-\'etale. So Condition ($S_3$) holds for $e$. Hence $\psi_e$ is regular by Theorem \ref{T10}(3.c), we have $\nu_{\textup{n}-\et}(e)=\rho(e)-\rho_{\et}(e)$ by Theorem \ref{T10}(3.e), and, with the notation of Theorem \ref{T10}(2.b), the morphism $W_e\rightarrow Y_e$ is a Galois cover between regular affine varieties by Theorem \ref{T10}(3.f). We have $\edim(X_e)=4$ by Corollary \ref{C25}.
\end{example}

The next example shows that the hypotheses of Theorem \ref{T10}(2.j), (3.a), (3.f), and (3.g) are necessary and that in Theorem \ref{T10}(2.b) for $p=2$, $W_e$ and $Y_e$ can have non-rational singularities and $Y_e$ can have tame local fundamental groups equipped with epimorphisms onto the cyclic group of order $3$. 

\begin{example}\normalfont\label{EX50}
Suppose that $p=n=2$. Let $e\in\EE_2(k)$ be such that  
$$Y_e=\Spec k[x_1,x_2][w]/(L_2)$$ 
is a monomial model of $X_e$, with $L_2$ is as in Example \ref{EX41.5}; we have $\deg(e)=3$. We use the notation of Theorem \ref{T10} but with $K$ replaced by $k$. 

Denoting, $f(x_1,x_2):=x_1^4+x_1+x_2^2$ we have $L_2=w^3+f(x_1,x_2)w+x_1x_2$. 

The field extension $k_{\t}(x_1,x_2)\rightarrow\Frac(B_e)$ is isomorphic to the field extension $k_{\t}(x_1,x_2)\rightarrow k_{\t}(x_1,x_2)[z]/\b(z^2+\frac{x_1^4x_2^4}{f^6}z+1\bigr)$ by \cite{MW}, Thms.\ 2.1(c) and 3.1 applied to $(e,f,g)=(0,f,x_1x_2)$. Hence $B_e$ is the normalization of the $k$-algebra $k_{\t}[x_1,x_2][z]/(z^2+x_1^4x_2^4z+f^{12})$. Replacing $z$ first by $z+f^6+x_1^2x_2^2f^3$ and second by $\frac{z+f^6+x_1^2x_2^2f^3}{x_1^3x_2^3}$, we get that $B_e$ is also the normalizations of the $k$-algebras $k_{\t}[x_1,x_2][z]/(z^2+x_1^4x_2^4z+x_1^6x_2^6f^3)$ and $k_{\t}[x_1,x_2][z]/(z^2+x_1x_2z+f^3)$. We compute 
$$f(x_1,x_2)^3=x_1^3(x_1^3+1)^3+x_1^2(x_1^3+1)^2x_2^2+x_1^4x_2^4+x_1x_2^4+x_2^6.$$ 
Replacing $z$ by $z+x_2^3$, we get that $B_e$ is also the normalization of the $k$-algebra 
$$k_{\t}[x_1,x_2][z]/\bigl(z^2+x_1x_2z+x_1^3(x_1^3+1)^3+x_1^2(x_1^3+1)^2x_2^2+x_1^4x_2^4\bigr).$$
Thus
$$B_e=k_{\t}[x_1,x_2][z]/\bigl(z^2+x_2z+x_1(x_1^3+1)^3+(x_1^3+1)^2x_2^2+x_1^2x_2^4\bigr)$$
as Zariski's Jacobian Criterion (see \cite{Gro2}, Prop.\ 22.6.7(iii)) gives that $\Spec B_e$ has precisely three singular points which form the zero locus $z=x_2=x_1^3+1=0$. Thus $|\Sing(Y_e)|=3$. 

Referring to Theorem \ref{T10}(2.g) and (2.h) we have identities $h_0=x_2$, $l=0$, $g_0=x_1(x_1^3+1)^3+(x_1^3+1)^2x_2^2+x_1^2x_2^4$, and $f_0=\frac{g_0}{h_0^2}\in R_{h_0}$.

Let $\alpha_1:=1$. Let $\alpha_2\in k^{\ast}$ and $\alpha_3\in k^{\ast}$ be the primitive third roots of unity. For each $i\in\{1,2,3\}$ we have $\alpha_i^3=1$ and we let $Q_i:=(\alpha_i,0)\in k^2=\mathbb A^2_{k,\t}(k)$. There exists a unique point $Q_{i,0}\in Y_e(k)$ that maps to $Q_i$. With $y_i:=x_1+\alpha_i$, the $k_{\t}[[y_i,x_2]]$-algebra $k_{\t}[[y_1,x_2]]/\bigl(w^3+f(y_i,x_2)w+(y_i+\alpha)x_2\bigr)$ has no retraction, i.e., there exists no $v$ in the maximal ideal $(x_1,x_2)$ of $R$ such that we have an identity $v^3+(y_i^4+y_i+x_2^2)v+(y_i+\alpha)x_2=0$; this implies that there exists a unique point $Q_{i,1}\in X_e(k)$ that maps to $Q_i$. As $2$ and $3$ are relatively prime, it follows that there exists a unique point $Q_{i,2}\in W_e(k)$ that maps to $Q_i$ (or $Q_{i,1}$ or $Q_{i,0}$). In particular, we have an inclusion $\{Q_{i,1}|i\in\{1,2,3\}\}\subset \Gamma_e$ which is in conformity to Theorem \ref{T6.6}(2.c). From Theorem \ref{T10}(2.j) we get that $\Sing(W_e)\subset\{Q_{1,2},Q_{2,2},Q_{3,2}\}$.

The Galois group $\Gal\bigl(\Frac(C_e)/\Frac(B_e)\bigr)\cong \mathbb Z/3\mathbb Z$ acts on $W_e$ by Theorem \ref{T10}(2.c) and it fixes each $Q_{i,2}$ with $i\in\{1,2,3\}$; moreover, it acts on $\Spec \mathcal O_{W_e,Q_{i,2}}$ in such a way that the action restricted to the punctured spectrum is free by Theorem \ref{T10}(2.j).

The multiplicative subgroup $G:=\{\alpha_1,\alpha_2,\alpha_3\}$ of $k^{\ast}$ acts on the diagram of Theorem \ref{T10}(2b) and to describe this action it suffices to describe its action on the finite flat morphism $X_e\rightarrow\mathbb A^2_{k,\t}$. On this morphism it acts via the following rule on elements of $k_{\t}[x_1,x_2]$ and the indeterminate $w$: for $\gamma\in G$, 
$$\bigl(\gamma,(x_1,x_2,w)\bigr)\mapsto (\gamma x_1,\gamma^{-1}x_2,\gamma^{-1}w).$$ 
The restriction of this action to the curve $Y:=X_e\setminus X_e^{\et}$ of the source of $\psi_e$ that maps bijectively onto the line $\{0\}\times\mathbb A^1_k$ of $\mathbb A^2_{k,t}$ has only one orbit of $k$-valued points with one element and has all other orbits of $k$-valued points with three elements; in particular, it acts transitively on the set $\{Q_{1,1},Q_{2,1},Q_{3,1}\}$ of points. From this and $|\Sing(X_e)|\le 1$ by Example \ref{EX41.5}, it follows that $\Sing(X_e)\subset \{Q_{0,2}\}$, where $(\psi_e^{-1})_{\red}=\{Q_{0,1},Q_{0,2}\}$ with $Q_{0,1}\in\Imm(\imath_e)$; so $Q_{0,2}=\bigl(Y\cap \psi_e^{-1}(0,0)\bigr)_{\red}$. 

As $\{Q_{1,2},Q_{2,2},Q_{3,2}\}$ is an orbit of the action of $G$ on $W_e$, the non-\'etale locus of the morphism $W_e\rightarrow Y_e$ is either empty or equal to $\{Q_{1,2},Q_{2,2},Q_{3,2}\}$; if it is empty, then $W_e\rightarrow Y_e$ is Galois and from Theorem \ref{T10}(3.g) we get for $i\in\{1,2,3\}$ that $Q_{i,1}$ is not the only point of $X_e$ that maps to $Q_i$, a contradiction. Thus the non-\'etale locus of the morphism $W_e\rightarrow Y_e$ is $\{Q_{1,2},Q_{2,2},Q_{3,2}\}$ and hence this morphism is non-Galois.

Thus the tame part of the local fundamental group of $\Spec \mathcal O_{Y_e,Q_{i,0}}$, i.e., of the fundamental group of the punctured spectrum of $\Spec \mathcal O_{Y_e,Q_{i,0}}$, is non-trivial as it has a quotient of order $3$. 

An orbit argument as in the second paragraph above gives that $\Sing(W_e)$ is either empty or equal to $\{Q_{1,2},Q_{2,2},Q_{3,2}\}$. We check that it is non-empty, i.e., $W_e$ is singular as follows. Picking up an index, say $1$, it suffices to show that the assumption that $Q_{1,2}\in\Reg(W_e)$ leads to a contradiction. To check this, as $L_2\in\mathbb F_2[x_1,x_2,w]$ we can assume that $e$ and $Q_{1,2}$ are defined over $\mathbb F_2$. So this assumption implies that the cotangent space $V_{1,2}$ of $W_e$ at $Q_{1,2}$ has dimension $2$ and that the action of the Galois group $S_3$ on $V_{1,2}$ is also defined over $\mathbb F_2$. As $\GL_2(\mathbb F_2)=\SL_2(\mathbb F_2)$ is isomorphic to $S_3$ and as an endomorphism of $S_3$ is either an automorphism or such that its elements of order $3$ are contained in the kernel, it follows that the element $\tau\in S_3$ of order $3$ that fixes $\Frac(B_e)$ is such that the $k$-vector subspace $V_{1,2}^{\tau}$ of $V_{1,2}$ fixed by $\tau$ has dimension $0$ or $2$. But $V_{1,2}^{\tau}$ contains the image of the cotangent space of $Y_e$ at $Q_{1,0}$ and hence it is non-empty as from the proof of Theorem \ref{T10}(2.b) we get that $Z:=\bigl(\phi^{-1}(\mathbb A^1_k\times\{0\})\bigr)_{\red}$ is such that its inverse image in $W_e$ consists of three lines that pass through $Q_{1,2}$ and that map birationally onto $Z$ and hence isomorphically as $Z\cong \mathbb A^1_K$. So $V_{1,2}^{\tau}=V_{1,2}$ and this implies that $\tau$ acts trivially on the local ring $\mathcal O_{W_e,Q_{1,2}}$, a contradiction. 

For $j\in\{1,2\}$ let $\mathcal O_j$ be the completion of the local ring $\mathcal O_{X_e,Q_{0,j}}$. The normalization of the $k_{\t}[[x_1,x_2]]$-algebra $k_{\t}[[x_1,x_2]]/(L_2)$ is $\mathcal O_1\times\mathcal O_2$, we have $\mathcal O_1\cong k_{\t}[[x_1,x_2]]$, and the integral flat morphism $k_{\t}[[x_1,x_2]]\rightarrow\mathcal O_2$ has degree $2$. Using Hensel's Lemma first and Nagata's Jacobian Criterion (see \cite{Gro2}, Prop.\ 22.7.2) in this context one gets that $\mathcal O_2$ is regular; so $\psi_e$ is regular and $X_e$ is a quasi-sphere by the identities $\rho(e)=1$ and $\rho_{\eta}(e)=0$ of Example \ref{EX41.5}.

For each $i\in\{1,2,3\}$ the local homomorphism $\mathcal O_{X_e,Q_{i,1}}\rightarrow\mathcal O_{W_e,Q_{i,2}}$ has a regular source and thus it is flat by \cite{Gro3}, Prop.\ (6.1.5) or \cite{Ma}, Ch.\ 8, Thm.\ 23.1.

As $|\Sing(W_e)|=3$, for each $i\in\{1,2,3\}$ the singularity of $W_e$ at $Q_{i,2}$ is rational iff the singularity of $Y_e$ at $Q_{i,0}$ is rational by Proposition \ref{PR25}(2). The last singularity has a local fundamental group whose tame part is non-trivial and of order divisible by $3$ and has a quadratic form which is a product of two non-proportional linear forms: in the indeterminates $y_1,x_2,z$, the quadratic form of $\mathcal O_{Y_e,Q_{1,2}}$ is $z(z+x_2)$. So from the classification of such types of singularities in positive characteristic $2$ obtained in \cite{Ar}, Sects.\ 3 and 4 which only involve the $A_l$ types with $l\in -1+3\mathbb N^{\ast}$, by considerations of homogeneous components of degree $3$ (resp.\ degrees $3$ to $5$) one rule out the case $l=2$ (resp.\ $l\ge 5$). We conclude that the singularities of $Y_e$ are non-rational. So the singularities of $W_e$ are also non-rational and in particular we obtain a second proof of the fact that $W_e$ is non-regular.\end{example}%From this, the fact that $\mathcal O_{Y_e,P}=\mathcal O_{W_e,P_2}^G$ is non-regular, it follows that $\mathcal O_{W_e,P_2}$ is non-regular. From this and Theorem \ref{T10}(2.c) we get that $|\Sing(W_e)|=3$.

The next lemma justifies the extra hypothesis of Theorem \ref{T10}(6) in the case when $h\in K^{\ast}$.

\begin{lemma}\label{L27} 
Suppose that $p=2$. Let $A$ be a $k$-algebra such that $k^{\ast}=A^{\ast}$. Let $f\in A$ and $B:=A[w]/(w^2+w+f)$; so $\Spec B\rightarrow \Spec A$ is an Artin--Schreier cover. Then the $A$-algebra $A\times B$ is generated by $1$ element iff $B\cong A^2$, equivalently, iff there exists $g\in A$ such that $g^2+g=f$.
\end{lemma}

\begin{proof}
If $B\cong A^2$, then $A\times B\cong A^3$ and the $A$-algebra $A^3$ is generated by each element $z:=(1,\alpha,\alpha^2)\in A^3$ with $\alpha\in k\setminus\{0,1\}$ as $\bigl((1,1,1),z,z^2\bigr)$ is an $A$-basis of the $A$-module $A^3$. So the `if' part holds. 

For the `only if' part, we define $v_1:=(1,0)$, $v_2:=\bigl(0,1+(w^2+w+f)\bigr)$, and $v_3:=\bigl(0,w+(w^2+w+f)\bigr)$; so $(v_1,v_2,v_3)$ is an $A$-basis of $A\times B$. As the $A$-algebra $A\times B$ is generated by $1$ element, there exists $(\alpha_1,\alpha_2,\beta)\in A^3$ such that $\bigl(v_1+v_2,\alpha_1 v_1+\alpha_2 v_2+\beta v_3,\alpha_1^2 v_1+(\alpha_2^2+f\beta^2) v_2+\beta^2 v_3\bigr)$ is an $A$-basis of $A\times B$. The matrix $\begin{bmatrix} 
1 & 1 & 0 \\
\alpha_1 & \alpha_2 & \beta \\ 
\alpha_1^2 & \alpha_2^2+f\beta^2 & \beta^2 \\
\end{bmatrix}$
has determinant $\beta\alpha_4\in A^{\ast}$, with $\alpha_3:=\alpha_1+\alpha_2\in A$ and $\alpha_4:=\alpha_3^2+\beta\alpha_3+f\beta^2$. Thus $(\beta,\alpha_4)\in (A^{\ast})^2$. Hence for $\alpha_5:=\alpha_3\beta^{-1}\in A$, we have $\alpha_6:=\alpha_5^2+\alpha_5+f=\frac{\alpha_4}{\beta^2}\in A^{\ast}=k^{\ast}$. As $B\cong A[w]/(w^2+w+\alpha_6)\cong A^2$ we get that the `only if' part holds. 
\end{proof}

\begin{lemma}\label{L28} 
Suppose that $p=2$. Let $R:=k[x_1,x_2]$. Let $\Spec B\rightarrow\Spec R=\mathbb A^2_k$ be a finite \'etale morphism of degree $3$ with $B$ an integral domain. Let $\Tr:B\to R$ be the $R$-linear trace map from $B$ to $R$. Then the $R$-algebra $B$ is generated by $1$ element iff there exists $u\in B^{\ast}$ with $\Tr(z)=0$ (so $u\notin k^{\ast}$).
\end{lemma}

\begin{proof}
For the `only if' part let $z\in B$ be such that the $R$-algebra $B$ is generated by $z$, i.e., $B=R[z]$. So $(1,z,z^2)$ is an $R$-basis of $B$. 

Let $g(w)=w^3+f_2w^2+f_1w+f_0\in R[w]$ be the minimal polynomial of $z$ over $R$. Substituting $z$ by $z-\frac{f_2}{3}$ we can assume that $f_2=0$; so $\Tr(z)=0$ and the discriminant of $g$ is $f_0^2$ . Hence, as $B$ is an \'etale $R$-algebra, we get that $f_0^2\in R^{\ast}$. So $z^{-1}=f_0^{-1}(z^2+f_1)\in B$ exists. Thus $z\in B^{\ast}$. So the `only if' part holds. 

For the `if' part, let $u\in B^{\ast}\setminus k^{\ast}$ with $\Tr(u)=0$. So $u\notin R$ and the minimal polynomial $h(w)\in R[w]$ of $u$ over $R$ has degree $3$. As $u$ is a unit, $h(0)$ is the determinant of the $R$-linear map of the multiplication by $u$ on either $B$ or on its $R$-subalgebra $R[u]$ generated by $u$. So, as $u\in B^{\ast}$, we have $h(0)\in B^{\ast}\cap R=R^{\ast}$. As $\Tr(u)=0$, the discriminant $h(0)^2$ of $h$ is a unit of $R$. Thus the $R$-algebra $R[u]$ is \'etale. So the birational finite morphism $\Spec B\to\Spec R[u]$ is \'etale and hence an isomorphism. As $B=R[u]$, the `if' part holds.
\end{proof}

We have the following consequence of Theorem \ref{T10} and Lemma \ref{L28} or Corollary \ref{C11}(1).

\begin{corollary}\label{C27}
Suppose that $p=2$ and that there exists $e\in\EE_2(k)$ such that $\deg(e)=3$. Then the following properties hold.

\medskip
{\bf (1)} The $k_{\t}[x_1,x_2]$-algebra $A_e$ is not generated by $1$ element.

\smallskip
{\bf (2)} Condition ($S_3$) holds for $e$.

{\bf (3)} Let $\Tr:C_e\to A_e$ be the $C_e$-linear trace map from $C_e$ to $A_e$. With the notation of Theorem \ref{T10}(2.b), there exists no unit $u\in C_e^{\ast}$ with $\Tr(u)=1$.
\end{corollary}

\begin{proof}
As $A_e$ is a $k$-subalgebra of $k_{\s}[x_1,x_2]$, we have $A_e^{\ast}=k^{\ast}$ and hence part (1) follows from Lemma \ref{L28}. As $2\nmid 3$, part (1) also follows from Corollary \ref{C11}(1).

Part (2) follows from Theorem \ref{T10}(1).

For part (3), we show that the assumption that there exists a unit $u\in C_e^{\ast}$ with $\Tr(u)=1$ leads to a contradiction. The discriminant of the minimal polynomial of $u$ over $A_e$ is $\Tr(u)^2$, and hence as in the proof of Lemma \ref{L28} we argue that the $A_e$-algebra $C_e$ is \'etale and generated by $u$. As $\bigl((1,1),(0,u),(0,u^2)\bigr)$ is an $A_e$-basis of the  $A_e$-module $A_e\times C_e$ and $\Spec C_e\rightarrow X_e$ is an Artin-Schreier cover by Lemma \ref{F6}, from Lemma \ref{L27} we get that $C_e\cong A_e^2$, a contradiction to part (2) by the very definition of $W_e$ as an integral affine scheme. So part (3) holds. 
\end{proof}

\begin{remark}\normalfont\label{R28}
{\bf (1)} Suppose that $p=2$ and $n\in\mathbb N^{\ast}\setminus\{1\}$. Let $R:=k[x_1,\ldots,x_n]$. Conjecture \ref{CJ3} predicts that there exists no $e\in\EE_n(k)$ such that $\deg(e)=3$, $\psi_e$ is \'etale, and Condition ($S_3$) holds for $K=k$. To prove (resp.\ disprove) this, based on \cite{MW}, Cor.\ 2.2($b_1$) and Thm.\ 3.1 it would suffice to show that there exists no (resp.\ there exists an) $F\in\Frac(R)$ such that the following \'etale characteristic $2$ conditions hold for it.\index{characteristic $2$ condition}

\medskip
{\bf ($\bigtriangleup_F$)} The polynomials $t^3+t+F$ and $t^2+Ft+1$ in $\Frac(R)[t]$ are irreducible and the normalizations of $R$ in $\Frac(R)[t]/(t^3+t+F)$ and $\Frac(R)[t]/(t^2+Ft+1)$ are \'etale $R$-algebras $A$ and $B$ (respectively).\index{characteristic $2$ condition!($\bigtriangleup_F$)}

\smallskip
{\bf ($\bigtriangledown_F$)} The affine variety $\Spec A$ contains an open subvariety isomorphic to $\mathbb A^n_k$.\index{characteristic $2$ condition!($\bigtriangledown_F$)}

\medskip
{\bf (2)} For $n\in\{2,3\}$ there exists a polynomial $F\in\Frac(R)$ such that Condition ($\bigtriangleup_F$) holds, $A^{\ast}=k^{\ast}$, $\Spec (A)$ is rational and factorial, and Condition ($\bigtriangledown_F$) does not hold. This follows from Theorem \ref{T2}(2) and (4) and Remarks \ref{R13} and \ref{R14} based on \cite{MW}, Cor.\ 2.2($b_1$) and Thm.\ 3.1. For $n=3$, such $F$s can be computed based on Remark \ref{R15}.

{\bf (3)} For $n=2$ there exists $F\in\Frac(R)$ such that Condition ($\bigtriangledown_F$) holds, the polynomials $t^3+t+F$ and $t^2+Ft+1$ in $\Frac(R)[t]$ are irreducible, but Condition ($\bigtriangleup_F$) does not hold; for instance, this holds if $t^3+t+F$ is irreducible and the field extensions $\Frac(R)\rightarrow\Frac(R)[t]/(t^3+t+F)$ and $\Frac(R)\rightarrow \Frac(R)[w]/(L_2)$ are isomorphic, with $L_2\in R[w]$ as in Example \ref{EX41.5}.

{\bf (4)} Suppose Condition ($\bigtriangleup_F$) holds. Let $f\in R\setminus\{\star^2+\star|\star\in R\}$ be such that $\Spec A=R[z]/(z^2+z+f)$ by Lemma \ref{F6} (cf.\ Theorem \ref{T10}(2.g)). As the field extensions $\Frac(R)[z]/(z^2+z+f)$ and $\Frac(R)[z]/(z^2+Fz+1)$ of $\Frac(R)$ are isomorphic, there exists $f_1\in\Frac(R)$ such that $\frac{1}{F^2}=f+f_1+f_1^2$; many such $f_1$'s (e.g., $f_1=f$) are excluded by the fact that $B$ is an \'etale $R$-algebra.

Using the substitution $f_2:=f+f_1\in\Frac(R)$, we get first that there exists $f_3\in\Frac(R)$ such that $f_2=f_3^2$ and second that $F=\frac{1}{f+f_3+f_3^2}$. Writing $f_3=\frac{g}{h}$, with $g$ and $h$ relatively prime polynomials in $R$, we get that $F=\frac{h^2}{g^2+gh+h^2f}$. It follows that $A$ is the normalization of the $R$-algebra 
$$A_1:=R[w]/\bigl(w^3+(g^2+gh+h^2f)^2w+h^2(g^2+gh+h^2f)^2\bigr).$$
As $h$ and $g^2+gh+h^2f$ are relatively prime, by unramified considerations at the generic points of irreducible components of the zero locus $g^2+gh+h^2f=0$ we get that there exists $h_1\in R$ such that $g^2+gh+h^2f=h_1^3$.

Replacing $w$ by $\frac{w}{h_1^2}$ we get that $A$ is also the normalization of the $R$-algebra 
$A_2:=R[w]/\bigl(w^3+h_1^2w+h^2\bigr)$. This forces $w+\bigl(w^3+h_1^2w+h^2\bigr)\in A_2$ to be a square and hence $A$ is also the normalization of the $R$-algebra 
$$A_3:=R[w]/\bigl(w^3+h_1w+h\bigr).$$ 
As $A^{\ast}=k^{\ast}$, it follows that $h\notin k^{\ast}$. 

So the pair $(h,h_1)\in (R\setminus k^{\ast})\times R$ is such that the following two conditions hold: (i) there exists $g\in R$ such that $h_1^3+g^2+gh\equiv 0\pmod{h^2}$, (ii) the normalization $A$ of the $R$-algebra $A_3$ is \'etale above each generic point of the zero locus $h=0$. 

Condition ($\bigtriangledown_F$) holds iff $\mathbb A^n_k$ is isomorphic to an open subvariety of $\Spec A$.\end{remark}

\section{Appendix A: weak separations}\label{S32}

For simplicity, we introduce the notion ``weakly separated" only for normal schemes of finite type over fields. 

\begin{definition}\label{D12}
Let $X$ be a normal scheme of finite type over a field $\mathcal K$. Let $\overline{\Delta}_X$ be the closure in $X\times_{\Spec \mathcal K} X$ of the image $\Delta_X$ of the diagonal immersion $X\rightarrow X\times_{\Spec \mathcal K} X$, and let $U\subset\overline{\Delta}_X$ be the (open) locus where both projections to $X$ are quasi-finite. We say that $X$ is weakly separated\index{weakly separated} if $U=\Delta_X$.
\end{definition}

Note that $U$ is open in $\overline{\Delta}_X$ by Chevalley's Semi-continuity Theorem (see \cite{Gro4}, Thm.\ (13.1.3)).

\begin{proposition}\label{PR22}
With the notation of Definition \ref{D12}, the following properties hold.

\medskip
{\bf (1)} Both projections $\pi_1:U\rightarrow X$ and $\pi_2:U\rightarrow X$ are local isomorphisms.

\smallskip
{\bf (2)} The morphism $U\rightarrow X\times_{\Spec \mathcal K} X$ is an equivalence relation.
\end{proposition}

Before proving this proposition, we first make a second definition and prove a lemma required in the proof.

\begin{definition}\label{D13}
Let $X$ be a normal scheme of finite type over a field $\mathcal K$. By the weak separation\index{weak separation} of $X$ we mean the quotient $X_{\w}:=X/U$ of $X$ by the equivalence relation $U$ (quotient by an equivalence relation of big Zariski sheaves or quotient in the category of ringed spaces). 
\end{definition}

\begin{lemma}\label{L29} 
Assume we have a commutative diagram
\begin{tikzcd}[row sep=tiny]
X \arrow[rd, "f"'] \arrow[rr, "h"] & & Z \arrow[ld, "g"] \\
& Y
\end{tikzcd}
%\begin{equation*}\end{equation*}
of integral schemes of finite type over a field $\mathcal K$ with $Y$ normal, $f$ dominant quasi-finite, and $g$ birational. Then $h(X)$ is contained in an open subscheme $\Omega$ of $Z$ such that the morphism $\Omega\rightarrow Y$ induced by $g$ is a local isomorphism.
\end{lemma}

\begin{proof}
It suffices to show that for each $x\in X$, the morphism $g$ induces a homomorphism $\mathcal O_{Y,f(x)}\rightarrow\mathcal O_{Z,h(x)}$ of local rings which is an isomorphism. To show this we can assume that $X$, $Y$, and $Z$ are affine and we can replace $X$ and $Z$ by their normalizations. Then $X$ is open in the normalization $W$ of $Z$ in the field of fraction of $X$ by Zariski's Main Theorem. By the going down property $h(X)$ is open in $Z$. So $h(X)\rightarrow Y$ is quasi-finite birational and thus an open embedding by Zariski's Main Theorem. Hence the homomorphism $\mathcal O_{Y,f(x)}\rightarrow\mathcal O_{Z,h(x)}$ is an isomorphism.
\end{proof}

\noindent
{\bf Proof of Proposition \ref{PR22}.} Part (1) follows from Lemma \ref{L29} applied to $h=1_U$. The relation defined by the monomorphism $U\rightarrow X\times_{\Spec \mathcal K} X$ is clearly reflexive and symmetric by definition. To show that it is transitive we can assume $X$ is integral and we consider the factorization $h:U\times_{\pi_2,X,\pi_1} U\rightarrow\overline{\Delta}_X$ of the morphism $U\times_{\pi_2,X,\pi_1} U\rightarrow X\times_{\Spec \mathcal K} X$ defined by the rule $(u,v)\mapsto\bigr(\pi_1(u),\pi_2(v)\bigl)$ on valued points. From Lemma \ref{L29} applied to $g$ equal to one of the two projections of $\overline{\Delta}_X$ on $X$ it follows that $h$ factors through an open subscheme $\Omega$ of $\overline{\Delta}_X$ with the property that both projections from it to $X$ are local isomorphisms. Hence $\Omega$ is contained in $U$. Thus the relation is also transitive, so part (2) also holds. 

\medskip
Lemma \ref{L29} implies that the association $X\mapsto U$ is functorial with respect to dominant quasi-finite morphisms between integral normal schemes of finite type over $\Spec\mathcal K$ and hence these morphisms induce morphisms on weak separations. One verifies that $\overline{\Delta}_X$ is the inverse image under the morphism 
$$X\times_{\Spec \mathcal K} X\rightarrow X_{\w}\times_{\Spec \mathcal K} X_{\w}$$ 
of the diagonal closure $\overline{\Delta}_{X_{\w}}$ in $X_{\w}
\times_{\Spec \mathcal K} X_{\w}$ of the image of the diagonal immersion $X_{\w}\rightarrow X_{\w}\times_{\Spec \mathcal K} X_{\w}$. It follows that weak separations are weakly separated.

\begin{example}\normalfont\label{EX51}
Let $d:C\rightarrow C$ be an endomorphism of a curve over $K$ that permutes the generic points of its irreducible components. Then $d^{|\Irr(C)|!}$ fixes these generic points and hence induces dominant endomorphisms of its irreducible components. So to prove that $\iota(d)\in\mathbb N$ we can assume that $C$ is irreducible. The weak separation $C_{\w}$ of $C$ is the integral separated curve obtained from $C$ by identifying the points whose local rings coincide inside the field of fractions of $C$. The functorial endomorphism $d_{\w}:C_{\w}\rightarrow C_{\w}$ induced by $d$ is finite surjective by Lemma \ref{L10}(3). By considering an open non-empty subset $U$ of $C_{\w}$ isomorphic to its inverse image in $C$, we have $U\subset\Imm(d^i)$ for all $i\in\mathbb N$. Thus $\iota(d)\le |C\setminus U|$.\end{example}

For an endomorphism $e:X\rightarrow X$ of a variety over $K$ and $P\in X(K)$ we consider the following iteration condition.

\medskip
\phantomsection{($\w\natural_{e,P}$)\index{iteration condition!($\w\natural_{e,P}$)} {\it Condition ($\flat_{e,P}$) holds and for each $n\in\mathbb N^{\ast}$ the Zariski closure of the set $\{e^m(P)|m\in\mathbb N, m\ge n, e^{-n}\bigl(e^m(P)\bigr)=\{e^{m-n}(P)\}\}$ contains a non-empty open subvariety $\Omega$ of an irreducible component of positive dimension of the Zariski closure $Y_{0,1}$ of the orbit $\{e^i(P)|i\in\mathbb N\}$.}}\label{PH53}

\medskip
Note that if ($\w\natural_{e,P}$) holds and $m\in\mathbb N^{\ast}$ is such that $e^m(P)\in\Omega$, then ($\w\natural_{d,e^m(P)}$) holds for each endomorphism $d:Y\rightarrow Y$ induced by $e$ with $Y$ a locally closed subvariety of $X$ that contains $e^m(P)$. Clearly, ($\natural_{e,P}$)$\Rightarrow$ ($\w\natural_{e,P}$)$\Rightarrow$ ($\flat_{e,P}$) and hence the following proposition is a weaker form of Theorem \ref{T4}(2).

\begin{proposition}\label{PR23}
Let $e:X\rightarrow X$ be a quasi-finite endomorphism of a variety over $K$. Then there exists no $P\in X(K)$ for which Condition $(\w\natural_{e,P})$ holds.
\end{proposition}

\begin{proof} We show that the assumption that there exists $P\in X(K)$ for which Condition $(\w\natural_{e,P})$ holds leads to a contradiction. We can assume that $Y_{0,1}=X$. Let $s\in\mathbb N^{\ast}$ be such that $e^s$ fixes the generic points of irreducible components of the equidimensional variety $X^{-}:=\cup_{Y\in\Irr(X),\dim(Y)>0} Y$ by Lemma \ref{L15}(1). Let $Z\in\Irr(X^-)$ and let $d:Z\rightarrow Z$ be induced by $e^s$. Let $l\in\mathbb N$ be the smallest such that $e^l(P)\in Z$. As $(\w\natural_{e,P})$ holds, so do ($\w\natural_{d,e^l(P)}$) and ($\flat_{d,e^l(P)}$). 

We choose $Z$ such that each point $e^m(P)$ in some non-empty open subvariety of $Z$ with $m\in\mathbb N\setminus\llbracket 0,s-1\rrbracket$ satisfies $d^{-1}\bigl(e^m(P)\bigr)=\{e^{m-s}(P)\}$. In particular, $d:Z\rightarrow Z$ is dominant (being quasi-finite with $Z$ irreducible) and generically radicial, and the same holds for its normalization $d^{\n}:Z^{\n}\rightarrow Z^{\n}$. From Zariski's Main Theorem we get that the functorial weak separation $d^{\n}_{\w}:Z^{\n}_{\w}\rightarrow Z^{\n}_{\w}$ is a bijection. By considering an open non-empty subset $\Omega$ of $Z$ naturally isomorphic to $(d^{\n})^{-1}(U)$ and its image in $Z^{\n}_{\w}$, we have $U\subset\Imm(d^i)$ for all $i\in\mathbb N$, which contradicts ($\flat_{d,e^l(P)}$).\end{proof}

\begin{corollary}\label{C28} Let $f:X\rightarrow Y$ be a birational morphism between irreducible varieties over $K$. Assume there exist endomorphisms $e:X\rightarrow X$ and $d:Y\rightarrow Y$ such that $d\circ f=f\circ e$. Let $Q\in Y(K)$ be such that its orbit $\{d^i(Q)|i\in\mathbb N\}$ is Zariski dense in $Y$. Let $i\in\mathbb N$ be the smallest such that there exists $P\in X(k)$ such that $f(P)=d^i(Q)$. If ($\w\natural_{d,Q}$) holds, then $e$ is not quasi-finite.\end{corollary}

\begin{proof}
As ($\w\natural_{d,Q}$) holds, so does ($\w\natural_{d,P}$); so corollary holds by Proposition \ref{PR23}.
\end{proof}

\begin{remark}\normalfont\label{R29} 
In the proof of Proposition \ref{PR6}, instead of its last paragraph, from Corollary \ref{C28} applied to the birational composite morphism $\Omega\rightarrow Z^{\n}\rightarrow Z$, the endomorphisms $e_{\Omega}:\Omega\rightarrow \Omega$ and $e_Z:Z\rightarrow Z$, and the point $e^m(P)\in Z(K)$, we get that $e_{\Omega}$ is not quasi-fine, a contradiction.\end{remark}

\section{Appendix B: on pullbacks of Lang torsors}\label{S33}

Let $\kappa$ be a subfield of $k$ such that the field extension $\kappa\rightarrow k$ is algebraic. Let $X=\Spec R$ be an integral affine scheme of finite type over $\Spec\kappa$ and of dimension at least $1$. To construct Galois covers over $X$ that are pullbacks of Lang torsors, we first prove the following general proposition.

\begin{proposition}\label{PR24}
Let $\mathcal K$ be a subfield of $K$ such that the field extension $\mathcal K\rightarrow K$ is algebraic. Let $\mathcal G$ be a smooth connected linear algebraic group over $\Spec\mathcal K$. Let $R_u(\mathcal G_K)$ be the unipotent radical of $\mathcal G_K$. Then the following properties hold.

\medskip
{\bf (1)} If there exists $N\in\mathbb N^{\ast}$ and a dominant morphism $\Sigma:\mathbb A^N_K\rightarrow\mathcal G_K$ whose general fibers are geometrically integral, then $\mathcal G_K/R_u(\mathcal G_K)$ is a simply connected semisimple group.

\smallskip
{\bf (2)} Assume that $R_u(\mathcal G_K)$ is defined over $\mathcal K$ (i.e., the unipotent radical $R_u(\mathcal G)$ of $\mathcal G$ is defined) and $\mathcal G_K/R_u(\mathcal G_K)$ is a simply connected semisimple group. We also assume that $R_u(\mathcal G)$ is $\mathcal K$-split (e.g., this holds if $\mathcal K$ is perfect) and $\mathcal G/R_u(\mathcal G)$ is quasi-split (e.g., this holds if $\mathcal K$ is finite). Then there exists $N\in\mathbb N^{\ast}$ and a surjective morphism $\Sigma:\mathbb A^N_K\rightarrow\mathcal G_K$ whose general fibers are geometrically integral.
\end{proposition}

\begin{proof}
To prove part (1), we can assume that $\mathcal K=K$. By composing $\Sigma$ with the quotient morphism $\mathcal G\rightarrow\mathcal G/R_u(\mathcal G)$, we can assume that $R_u(\mathcal G)$ is trivial, i.e., $\mathcal G$ is a reductive group. By reasons of units, each morphism $\mathbb A^N_K\to\mathbb G_{\m,K}$ is constant. So $\mathbb G_{\m,K}$ is not a quotient of $\mathcal G$ and thus $\mathcal G$ is semisimple. If $\mathcal G^{\sc}$ is the simply connected semisimple group cover of $\mathcal G$, then its center is of multiplicative type (see \cite{SGA3-3}, Exp.\ XXII, Cor.\ 4.1.7), hence the kernel $\Ker$ of the central isogeny $\mathcal G^{\sc}\rightarrow \mathcal G$ is finite of multiplicative type. As $H^1_{\et}(\mathbb A^N_K,\pmb{\mu}_{m,K})=0$ for each $m\in\mathbb N^{\ast}$ and $\Ker$ is a finite product of such finite group schemes $\pmb{\mu}_{m,K}$ over $\Spec K$, the $\Ker$-torsors on $\mathbb A^N_K$ are trivial. So $\Sigma$ factor through $\mathcal G^{\sc}$ and thus, as its general fibers are geometrically integral, it follows that $\Ker$ is trivial, i.e., $\mathcal G=\mathcal G^{\sc}$; so part (1) holds.

To prove part (2), let $\mathcal G\cong R_u(\mathcal G)\times\mathcal G/R_u(\mathcal G)$ be an isomorphism of schemes defined over $\Spec\mathcal K$ by \cite{Ros}, Thm.\ 1. By the very definition of a $\mathcal K$-split unipotent group, $R_u(\mathcal G)$ has a composition series whose factors are $\mathbb G_{\a,\mathcal K}$, and hence there exists $m\in\mathbb N$ such that the varieties $R_u(\mathcal G)$ and $\mathbb A^m_{\mathcal K}$ over $\Spec\mathcal K$ are isomorphic. Based on the last two sentences, we can assume that $R_u(\mathcal G)$ is trivial. As $\mathcal G$ is a simply connected semisimple group which is quasi-split, it has a Borel subgroup and hence it has a pair of opposite Borel subgroups $B^+_{\mathcal G}$ and $B^-_{\mathcal G}$. For $*\in\{+,-\}$, let $U^*_{\mathcal G}$ be the unipotent radical of $B^*_{\mathcal G}$. From \cite{CGP}, Part I, Ch.\ II, Props.\ 2.2.9 and 2.1.10 and Cor.\ 2.2.5 we get that $U^*_{\mathcal G}$ is $\mathcal K$-split. Hence each product over $\Spec\mathcal K$ of a finite number of copies of either $U^+_{\mathcal G}$ or $U^-_{\mathcal G}$ is isomorphic to an affine space over $\mathcal K$. Let $\mathcal G'$ be the subgroup of $\mathcal G$ generated by $U^+_{\mathcal G}$ and $U^-_{\mathcal G}$ by \cite{Bore2}, Ch.\ I, Prop.\ 2.2; there exists a smallest $q_0\in\mathbb N$ such that for each integer $q\ge q_0$ we have a surjective product morphism $\Sigma_q:V_q\rightarrow\mathcal G^{\prime}$ with $V_q:=\prod_{i=1}^q U_i$ where $U_i$ is $U^-_{\mathcal G}$ if $i$ is odd and is $U^+_{\mathcal G}$ if $i$ is even. E.g., if $\mathcal G=\SL_{2,\mathcal K}$ then $q_0=4$ as one can check easily based on the analog of Isomorphism (\ref{EQ8}) over $\Spec\mathcal K$. 

To prove that $\mathcal G^{\prime}=\mathcal G$, in this paragraph we assume that $\mathcal K=K$. As $\mathcal G$ is a product of simply connected semisimple groups that have simple adjoints by the classification theorem \cite{SGA3-3}, Exp.\ XXV, Thm.\ 1.1, we can assume that the adjoint group of $\mathcal G$ is simple. As $\mathcal G'$ is normalized by $U^+_{\mathcal G}$, $U^-_{\mathcal G}$, and the maximal torus $T:=B^+_{\mathcal G}\cap B^-_{\mathcal G}$ of $\mathcal G$, it is a normal subgroup of $\mathcal G$ and hence $\mathcal G'=\mathcal G$.

Let $\mathcal K_1$ be a finite field extension of $\mathcal K$ such that $\mathcal G_{\mathcal K_1}$ is split. We show that there exists $q_1\in\mathbb N^{\ast}$ such that for each integer $q\ge q_1$, $(\Sigma_q)_{\mathcal K_1}$ has a rational section (often called a cross-section). To prove this, we can assume that $\mathcal G$ is non-trivial and, in order to ease the notation, that $\mathcal K_1=\mathcal K$; so $\mathcal G$ is split. Let $\Phi(\mathcal G,T)$ be the root system of $\mathcal G$ with respect to $T$, let $\Phi^+(\mathcal G,T)$ be the set of positive roots of $\Phi(\mathcal G,T)$ with respect to $B^+_{\mathcal G}$, and let $\Delta\subset \Phi^+(\mathcal G,T)$ be the basis of $\Phi^+(\mathcal G,T)$. For each $\alpha\in\Phi^+(\mathcal G,T)$ there exists a unique $\mathbb G_{\a,\mathcal K}$ closed subgroup scheme $U_{\alpha}$ of $B^+_{\mathcal G}$ (resp.\ $U_{-\alpha}$ of $B^-_{\mathcal G}$) which is normalized by $T$ and such that the action of $T$ on the Lie algebra of $U_{\alpha}$ (resp.\ $U_{-\alpha}$) is via the character of $T$ that corresponds to $\alpha$ (resp.\ $-\alpha)$. Let $T\mathcal G(\alpha)$ be the subgroup of $\mathcal G$ generated by $T$, $U_{\alpha}$, and $U_{-\alpha}$ (see \cite{SGA3-3}, Exp.\ XXII, Sect.\ 5, Prop.\ 5.10.1);  as $\mathcal G$ is simply connected, its derived group $\mathcal G(\alpha)$ is an $\SL_{2,\mathcal K}$ group normalized by $T$. We recall from \cite{SGA3-3}, Exp.\ XXII, Sect.\ 4, Prop.\ 4.1.2 that the product morphism
$$\imath_{\mathcal G}:\left(\prod_{\alpha\in -\Phi^+(\mathcal G,T)} U_{\alpha}\right)\times_{\Spec\mathcal K} T\times_{\Spec\mathcal K} \left(\prod_{\alpha\in \Phi^+(\mathcal G,T)} U_{\alpha}\right)\rightarrow\mathcal G$$
is an open embedding whose (`big cell') image $W:=\Imm(\imath_{\mathcal G})$ does not depend on the orderings in the first and third factor of the source of $\imath_{\mathcal G}$, being in fact equal to the image $B^-_{\mathcal G}B^+_{\mathcal G}$ of the product morphism $B^-_{\mathcal G}\times_{\mathcal K} B^+_{\mathcal G}\rightarrow\mathcal G$. 

If $\Phi(\mathcal G,T)$ is a direct sum of irreducible root systems $\Phi_i(\mathcal G,T)$ indexed by $i$ in the set $\mathcal B$ of simple factors of the adjoint group of $\mathcal G$ and if for $i\in\mathcal B$ we define $\Delta_i:=\Delta\cap\Phi_i(\mathcal G,T)$, let 
$$q_1:=2+4\max(|\Delta_i||i\in\llbracket1,|\mathcal B|\rrbracket).$$ 

The product morphisms $\left(\prod_{\alpha\in -\Phi^+(\mathcal G,T)} U_{\alpha}\right)\rightarrow U^-_{\mathcal G}$ and $\left(\prod_{\alpha\in \Phi^+(\mathcal G,T)} U_{\alpha}\right)\rightarrow U^+_{\mathcal G}$ are isomorphisms and moreover, as $\mathcal G$ is simply connected, the product homomorphism $\prod_{\alpha\in\Delta} T\cap \mathcal G(\alpha)\rightarrow T$ is an isomorphism. Based on this, to show for each integer $q\ge q_1$ the existence of a section for the morphism $\Sigma_q^{-1}(W)\rightarrow W$ induced by $\Sigma_q$, we can assume that $q=q_1$ and it suffices to show that for each $\alpha\in\Delta$ there exists a morphism $\sigma_{\alpha}:T\cap \mathcal G(\alpha)\rightarrow U_{\alpha}\times_{\Spec\mathcal K} U_{-\alpha}\times_{\Spec\mathcal K} U_{\alpha}\times_{\Spec\mathcal K} U_{-\alpha}$ whose composite with the product morphism $U_{\alpha}\times_{\Spec\mathcal K} U_{-\alpha}\times_{\Spec\mathcal K} U_{\alpha}\times_{\Spec\mathcal K} U_{-\alpha}\rightarrow \mathcal G(\alpha)$ is the inclusion $T\cap \mathcal G(\alpha)\rightarrow \mathcal G(\alpha)$. As $\mathcal G(\alpha)\cong\SL_{2,\mathcal K}$, the existence of $\sigma_{\alpha}$ follows from the following matrix identity
$$\begin{bmatrix} 
t & 0 \\
0 & t^{-1} \\ 
\end{bmatrix}=\begin{bmatrix} 
1 & -1 \\
0 & 1 \\ 
\end{bmatrix}\begin{bmatrix} 
1 & 0 \\
1-t & 1 \\ 
\end{bmatrix}\begin{bmatrix} 
1 & t^{-1} \\
0 & 1 \\ 
\end{bmatrix}\begin{bmatrix} 
1 & 0 \\
t(t-1) & 1 \\ 
\end{bmatrix}.$$

As the source and target of $(\Sigma_q)_K$ are smooth varieties over $K$, the smooth locus of $(\Sigma_q)_K$ contains the image of each section of the pullback of $\Sigma_q^{-1}(W)\rightarrow W$ to $\Spec K$ and hence it is non-empty. If $q\ge\max(q_0,q_1)$, then the fact that the surjection $(\Sigma_q)_K$ has a rational section implies that the $K$-algebra homomorphism $(\Sigma_q)_K^{\#}:\mathcal O(\mathcal G_K)\rightarrow\mathcal O((V_q)_K)$ induces an algebraically closed field extension $\Frac((\Sigma_q)_K^{\#})$. This field extension is also separable as the smooth locus of $(\Sigma_q)_K$ is non-empty. So the fiber of $\Sigma_q$ over the generic point of $\mathcal G$ is geometrically integral. Thus the general fibers of $\Sigma_q$ are geometrically integral by \cite{Gro4}, Thm.\ (9.7.7). So part (2) holds by taking $\Sigma=\Sigma_q$ and $N=q\dim(U^-_{\mathcal G})$ if $\mathcal G$ is non-trivial (i.e., $\dim(U^-_{\mathcal G})\ge 1$) or $N=1$ if $\mathcal G$ is trivial.\end{proof}

Let $\kappa'$ be the algebraic closure of $\kappa$ in $\Frac(R)$; let $X'$ be a Zariski dense open subscheme of $X$ which is geometrically irreducible over $\Spec\kappa'$. We denote by $G^{\ab}$ the abelianization of an abstract group $G$. The following result refines \cite{Kam}, Thm.\ and Item 2 of p.\ 640 from several points of view (e.g., the case $p^q=2$ is not excluded and, when $\kappa$ infinite, we do not requite either $\mathcal G$ to have trivial unipotent radical or $X$ to be geometrically integral over $\Spec\kappa$).

\begin{theorem}\label{T13}
Let $q\in\mathbb N^{\ast}$. Let $\mathcal G$ be a smooth connected linear algebraic group over $\mathbb F_{p^q}$ such that $\mathcal G/R_u(\mathcal G)$ is simply connected semisimple. Then there exists a Galois cover $Y=\Spec S\rightarrow X$ of Galois group $G:=\mathcal G(\mathbb F_{p^q})$ such that $S$ is an integral domain and $\kappa'$ is algebraically closed in $\Frac(S)$ (hence, if $Y':=Y\times_X X'$, the Galois cover $Y'\rightarrow X'$ is geometrically irreducible over $\kappa'$).
\end{theorem}

\begin{proof}
Let $F_{p^q}$ be as in the beginning of Section \ref{S19}. For the Weil restriction of scalars $\mathcal G':=\Res_{\mathbb F_{p^q}/\mathbb F_p} \mathcal G$ (see \cite{CGP}, Part IV, App.\ A.5) we have the following identities $\mathcal G'(\mathbb F_p)=\mathcal G(\mathbb F_{p^q})=G$ and $\mathcal G'_k=\mathcal G_k^q$. Thus by replacing $\mathcal G$ with $\mathcal G'$ we can assume that $q=1$. Let $\mathcal L:\mathcal G\rightarrow\mathcal G$ be the Lang torsor; it is a $G$-torsor defined on valued points by the rule $\Box\mapsto\Box\cdot F_p(\Box)^{-1}$. To prove the theorem we can assume that $|G|\ge 2$, i.e., that $\dim(\mathcal G)\ge 1$. Recall that $\mathcal G/R_u(\mathcal G)$ is quasi-split (e.g., see \cite{Bore2}, Ch.\ V, Prop.\ 16.6) and $R_u(\mathcal G)$ is $\mathbb F_p$-split (e.g., see \cite{Bore2}, Ch.\ V, Prop.\ 15.5 (ii)). 

Let $U:=\mathbb A_{\mathbb F_p}$. Based on Proposition \ref{PR24}(2), there exists $N\in\mathbb N^{\ast}$ for which there exists a surjective morphism $\Sigma:U^N\rightarrow\mathcal G$ whose general fibers are geometrically integral. Thus the pullback $\mathcal M_N:\mathbb X_N\rightarrow U^N$ of $\mathcal L$ via $\Sigma$ is a $G$-torsor with $\mathbb X_N$ geometrically irreducible over $\Spec\mathbb F_p$.

Following \cite{Kam}, we consider two disjoint cases as follows.

{\bf Case 1: $\kappa$ is a finite field.} For each irreducible $G$-torsor $Y'\rightarrow X'$, let $\kappa''$ be the algebraic closure of $\kappa$ in the field of fractions of $Y'$. The inclusion $\kappa'\subset\kappa''$ is a finite Galois extension and we have a surjective homomorphism $\mathcal G(\mathbb F_p)^{\ab}\rightarrow\Gal(\kappa''/\kappa')$. In particular, if $m$ is the exponent of $\mathcal G(\mathbb F_p)^{\ab}$ and $\kappa'\rightarrow\kappa'_m$ is the field extension of degree $m$, for $Y'$ to be geometrically irreducible over $\kappa'$ it suffices that $Y'_{\kappa'_m}$ is irreducible. 

Let $\kappa_m'\subset \kappa_+$ be a finite field extension. By a version of Chebotarev density theorem (cf.\ \cite{Del2}, Thm.\ 3.5.3 in the case of finite monodromy), we can assume that $\kappa_+$ is large enough so that the following two properties hold: (i) for each conjugacy class $\textup{conj}$ of $G$ there exist closed points of $\mathbb A^N_{\mathbb F_p}$ with residue field $\kappa_+$ and Frobenius conjugacy class (defined by $\mathcal M_N$) $\textup{conj}$; (ii) the number of closed points of $X'$ with residue field $\kappa_+$ is at least equal to the number of conjugacy classes of $G$. It follows that there exists a morphism $\Sigma_X:X\rightarrow U^N$ with prescribed values on certain $\kappa_+$-valued closed points, and using Jordan's Theorem (e.g., see \cite{Ser2}, Thm.\ 4' or 4 and \cite{Kam}, Lem.\ 1) as in \cite{Kam}, Sect.\ 2, the $G$-torsor $Y:=\mathbb X_N\times_{\mathcal M_N, U^N,\Sigma_X} X\rightarrow X$ is such that the $G$-torsor $Y'_{\kappa_+}\rightarrow X'_{\kappa_+}$ is irreducible. Thus $Y'_{\kappa'_m}$ is also irreducible and the theorem holds in this case.

{\bf Case 2: $\kappa$ is an infinite field.} Let $m\in\mathbb N^{\ast}$ be such that there exists a closed embedding $c:X\rightarrow\mathbb A^m_{\kappa}$; so $c_N:=c\times 1_{U^N_{\kappa}}:X\times_{\Spec\kappa} U^N_{\kappa}\rightarrow \mathbb A^m_{\kappa}\times_{\Spec\kappa} U^N_{\kappa}=\mathbb A^{m+N}_{\kappa}$ is a closed embedding. We show that there exists a morphism $\Sigma_X:X\rightarrow U_{\kappa}^N$ such that the $G$-torsor $Y:=\mathbb X_{N,\kappa}\times_{\mathcal M_{N,\kappa}, U^N_{\kappa},\Sigma_X} X\rightarrow X$ sits in the cartesian diagram
\[\xymatrix{
Y \ar[r] \ar[d] & \mathbb X_{N,\kappa} \ar[d]^{\mathcal M_{N,\kappa}} \\
X \ar[r]^{\Sigma_X} & U_{\kappa}^N
}\]
and has the property that the induced torsor $Y'\rightarrow X'$ is geometrically irreducible over $\kappa'$. This is equivalent to cutting with a general subspace of codimension $N$ of the affine space $\mathbb A^m_{\kappa}\times_{\Spec\kappa} U^N_{\kappa}=\mathbb A^{m+N}_{\kappa}$. 

As $\dim(X'\times_{\Spec\kappa'} \mathbb X_{N,\kappa'})\ge N+1$, the existence of $\Sigma_X$ follows from \cite{Jo1}, Part I, Cor.\ 6.7(3) applied (or from \cite{Jo1}, Part I, Thm.\ 6.3(4) applied $N$ times) to $c_N\circ (1_X\times\mathcal M_{N,\kappa}):X'\times_{\Spec\kappa'}\mathbb X_{N,\kappa}\rightarrow \mathbb A^m_{\kappa'}\times_{\Spec\kappa'} U^N_{\kappa'}=\mathbb A^{m+N}_{\kappa'}$.\end{proof}

\section{Appendix C: open problems, questions, and conjectures}\label{S34}

\noindent
{\bf P1.} For all integers $n\ge 2$ and $l\ge 0$ compute the smallest $\pi_{p,n,l}\in\{1\}\cup p\mathbb N^\ast$
such that there exists an $F_n$-system $\mathcal S$ over $k$ with $\pi(\mathcal S)=\pi_{p,n,l}$ and $\s(\mathcal S)=l$. Clearly, $\pi_{p,n,1}=1$ and the sequences $(\pi_{p,n,l})_{n\ge 2}$ are non-increasing.

\begin{remark}\normalfont\label{R30}
Each $F_n$-system $\mathcal S$ over $k$ is equivalent to an $F_m$-system over $k$ of the form $x_i+f_i^p=0$ for $i\in \llbracket1,m\rrbracket$, where $m\in\mathbb N^\ast$ depends only on $n$ and $\pi(\mathcal S)$ and the partial degree of $f_i\in k[x_1,\ldots,x_m]$ in $x_j$ is at most $1$ for all $(i,j)\in \llbracket1,m\rrbracket^2$.
\end{remark}

\begin{example}\normalfont\label{EX52} {\bf (1)} Clearly, for $l\in\mathbb N$ we have $\pi_{p,n,l}=1$ iff $l=1$.

{\bf (2)} We check that $\pi_{p,n,l}=p$ iff $l\in\{p^m|m\in \llbracket1,n\rrbracket\}$. Based on part (1), it suffices to show that $\pi_{p,n,l}\le p$ iff $l\in\{p^m|m\in \llbracket0,n\rrbracket\}$. For the `only if' part see \cite{V}, Thm.\ 2.4.1(b). The `if part' follows from the fact that for each $m\in \llbracket0,n\rrbracket$ the $F_n$-system $x_1-x_1^p=\cdots x_m-x_m^p=x_{m+1}=\cdots=x_n=0$ has $p^m$ solutions in $k^n$.

{\bf (3)} We have $\pi_{p,n,0}\ge 2p$ by parts (1) and (2). 

{\bf (4)} The system of equations $x-x^py^{p(p-1)}+1=y-xy^p-1=0$ over $k$ has no solutions by Example \ref{EX11} applied to $(m,q)=(1,p)$ and it is equivalent to the $F_2$-system $
x-x^py^{p(p-1)}+1=y+y^p-x^py^{p^2}-1=0$ over $k$, hence $\pi_{p,2,0}\le p^2+p$, and also to the $F_3$-system 
$x-x^py^{p(p-1)}+1=y+y^p-x^pz^p-1=0=z-y^p$ 
over $k$, hence $\pi_{p,3,0}\le p^2$.

{\bf (5)} We have $\pi_{p,2,0}\ge 3p$ by Proposition \ref{PR11}.

{\bf (6)} Parts (1), (3), (4), and (5) give for $n=2$ that $\pi_{2,2,0}=6$ and $\pi_{3,2,0}\in\{9,12\}$ and for $n\ge 3$ that $\pi_{2,n,0}=4$ and $\pi_{3,n,0}\in\{6,9\}$.\end{example}

\noindent
{\bf P2.} Identify properties of the finite sets $\mathcal F(e)\in\mathcal P_f(\mathbb N)$ with $e\in\BE_2(k)$. 

\begin{remark}\normalfont\label{R31} If $\mathcal F(e)\neq\emptyset$ for $e\in\BE_2(k)$, then the union $\mathcal F(e)\cup\{\deg(e)\}$ has a `midpoint' $[\deg(e)+\min\bigl(\mathcal F(e)\bigr)]/2$ but in general it is not symmetric with respect to it (e.g., this is so if in Proposition \ref{PR9} we have $1\in\mathbb U$ and $|\mathbb U|\ge 2$).
\end{remark}

\noindent
{\bf P3.} Describe the inclusions $\{\varphi\bigl(e^{\iota(e)}\bigr)|e\in\BE_2(k)\}\subset \{\varphi\bigl(e^{\iota(e)}\bigr)|e\in\EE_2(k)\}$
and $\{\iota(e)|e\in\BE_2(k)\}\subset \{\iota(e)|e\in\EE_2(k)\}$ and the set $\{\deg(e)|e\in\BE_2(k),\varphi(e)>0\}$. In particular, does there exists $e\in\BE_2(k)$ with $\deg(e)=p$ and $\varphi(e)>0$? 

\smallskip\noindent
{\bf P4.} 
For $n\in\mathbb N^{\ast}$ compute the smallest $\vartheta_p(n)\in\mathbb N^{\ast}$ for which there exists $e\in\BE_2(k)$ with $\varphi(e)=n$ and $\deg(e)=\vartheta_p(n)$. 

\begin{example}\normalfont\label{EX53}
We have $\vartheta_2(4n-1)\le 4n$ by Theorem \ref{T1}(3).
\end{example}

\noindent
{\bf P5.} Classify the sequences $\bigl((\varphi(e^n)\bigr)_{n\in\mathbb N}$ with $e\in\QF_2(K)$. 

\smallskip\noindent
{\bf P6.} Let $d\in\EE_n(K)$. For $n\ge 2$, classify the sequences $\bigl(\rho_{\et}(d^n)\bigr)_{n\in\mathbb N}$ and $\bigl(\rho(d^n)\bigr)_{n\in\mathbb N}$ and decide if $\rho_{\et}(X_d)$ is or is not equal to $\rho_{\et}(d)$ (cf.\ Theorem \ref{T3+}(5)). For $n=2$, how does the pair $\bigl(\rho_1(d),\rho_2(d)\bigr)$ vary with $d$ when the isomorphism class $X_d$ is fixed?

\begin{example}\normalfont\label{EX54}
If $e\in\EE_2(k)$, then for each pair $(n,m)\in (\mathbb N^{\ast})^2$, the identity $e^{n+m}=e^ne^m$ induces a finite flat morphism $X_{e^{n+m}}\rightarrow X_{e^n}$ such that the inverse image of the open subvariety $\mathbb A^2_{k,\s}$ of $X_{e^n}$ is $X_{e^m}$, hence $\rho(e^{n+m})\ge \rho(e^n)+\rho(e^m)$.
\end{example}

\noindent
{\bf P7.} Given $e\in\EE_n(K)$ and an open embedding $\imath:\mathbb A^n_K\rightarrow X_e$, is it true that $\Imm(\imath)\subset X_e^{\et}$? If it is not true, then does there exist $d\in\EE_n(K)$ such that $\imath_d$ is isomorphic to $\imath$? 

Let $X$ be an affine connected smooth surface over $k$ that has an open subvariety isomorphic to $\mathbb A^2_k$. Does there exist $e\in\EE_2(k)$ such that $X$ is isomorphic to an open subvariety of $X_e$? If yes, can we actually choose $e\in\EE_2(k)$ such that $X\cong X_e$? If not, is it true that there exists a non-surjective \'etale endomorphism $d:X\rightarrow X$?

\smallskip\noindent
{\bf P8.} Study other invariants of $e\in\EE_n(K)$ such as the locally closed subvariety $N_l(e)$ with $l\in\mathcal F(e)\cup\{\deg(e)\}$ of Remark \ref{R6} or the isomorphism class of the Galois group of the finite Galois extension of $K_{\t}(x_1,\ldots,x_n)$ generated by $\Frac(e^{\#})$, or the number $\nu(e)$ of irreducible components of $N_e$. 

\smallskip\noindent
{\bf P9.} For $e\in\EE_2(k)$, describe relations such as inequalities or values of polynomial expressions between invariants of its iterates. 

\begin{example}\normalfont\label{EX55}
The set $\{\frac{\rho(e)}{\deg(e)}|e\in\EE_2(k)\}$ is unbounded (see Theorem \ref{T1}(2)).
\end{example}

\noindent
{\bf P10.} \phantomsection{Use special polynomials as in \cite{Wr3}, Defs.\ 1.3 to 1.7 to provide more examples of moduli spaces of endomorphisms $e\in\EE_n(R)$ that have relevant properties. E.g., for each pair $(n,q)\in (\mathbb N^{\ast}\setminus\{1\})^2$ we have a closed subscheme}\label{P10}
$HEE_{n,q}^{\ast}$ of $EE_{n,q}^{\ast}$ whose $R$-valued points are those $n$-tuples of the form
$$(f_1,\ldots,f_n)=(x_1+l_1^q,\ldots,x_n+l_n^q)\in R[x_1,\ldots,x_n]^n$$ 
with $l_1,\ldots,l_n$ linear forms. Of particular interest is to describe the set
$$\nabla_{n,q,\chr(K)}:=\{\deg\bigl(\e(f_1,\ldots,f_n)\bigr)|(f_1,\ldots,f_n)\in HEE_{n,q}^{\ast}(K)\}.$$

\begin{example}\normalfont\label{EX56} {\bf (1)} We have $\nabla_{2,q,0}=\{1\}$, i.e., $(HEE_{2,q,K})_{\red}\rightarrow SA^{\ast}_{2,q,K}$ is a closed embedding if $\chr(K)=0$ (see \cite{vdE}, Part I, Ch.\ 2, Sect.\ 2.1, Exc.\ 4 due to Furter). 

{\bf (2)} We have $\nabla_{n,p,p}=\{1,p,p^2,\ldots,p^n\}$ (the part `$\subset$' follows from \cite{V}, Thm.\ 2.4.1(b); the part `$\supset$' can be seen by taking $l_i\in\{0,x_i\}$ for all $i\in \llbracket1,n\rrbracket$).
 
{\bf (3)} For $l\in\mathbb N^\ast$ we have $\nabla_{2,p^l+1,p}=\nabla_{3,p^l+1,p}=\{1\}$, i.e., we have closed embeddings $(HEE_{2,1+p^l,k})_{\red}\rightarrow SA^{\ast}_{2,p^l+1,k}$ and $(HEE_{3,p^l+1,k})_{\red}\rightarrow SA^{\ast}_{3,p^l+1,k}$.
\end{example}

\noindent
{\bf P11.} Classify all $SEE_{n,q,K}$ that are irreducible. Compute the dimensions of $SEE_{n,q,K}$ and $SA^{\ast}_{n,q,K}$. 

\smallskip\noindent
{\bf P12.} Prove an analog of Theorem \ref{T7} with $SEE$ replaced by $SA$.

\smallskip\noindent
{\bf P13.} Is it possible to define moduli schemes $BE_{n,q,\mathbb F_p}$ and $SBE_{n,q,\mathbb F_p}$?

\smallskip\noindent
{\bf P14.} Identify more smooth affine surfaces over $K$ that have proper open subvarieties isomorphic to $\mathbb A^2_K$ and compute their embedding dimensions (in the set $\{3,4,5\}$, see Theorem \ref{T3}(1)) and decide if they are or are not complete intersections. For $(n_1,n_2,n_3)\in\mathbb N^3$ with $n_1+n_2+n_3>0$ and $q\in \{3,4,5\}$ classify all $e\in\EE_2(k)$ with $\bigl(\rho_1(e),\rho_2(e),\rho(e)-\rho_{\et}(e)\bigr)=(n_1,n_2,n_3)$ and $\edim(X_e)=q$.

\smallskip\noindent
{\bf P15.} Classify all Galois field extensions $\mathcal Z$ of $k(x,y)$ which are of degree $p$ and unirational in a way that one can decide if $\mathcal Z$ is rational or if $\mathcal Z$ has the property that the spectrum $X$ of the normalizations of $k[x,y]$ in $\mathcal Z$ has an affine open subvariety $Y$ which is \'etale over $\mathbb A^2_k$ and: (i) isomorphic to $\mathbb A^2_k$; (ii) finite over $\mathbb A^2_k$; (iii) such that $\mathcal O(Y)=k^{\ast}$; (iv) factorial.

\smallskip\noindent
{\bf P16.} If $k$ is an algebraic closure of $\mathbb F_p$, does there exist $d_2\in\D_2(k)$ which is generically \'etale with $\deg(d_2)>1$ and $\iota(d_2)=\infty$?

\smallskip\noindent
{\bf P17.} Do there exist affine regular varieties $X$ over $k$ with trivial tangent bundles but without $p$-bases? In particular, do there exist \'etale morphisms $\SL_{2,k}\rightarrow\mathbb A^3_k$? 

\smallskip\noindent
{\bf P18.} If $e\in\EE_n(K)$, is it true that $\Sing(X_e)\subset \cup_{(Y,Z)\in\Irr(X_e\setminus X_e^{\et})^2,Y\neq Z} Y\cap Z$? 

\begin{remark}\normalfont\label{R32}
{\bf (1)} For $n=2$, if this is true then $\Sing(X_e)=\emptyset$ by Corollary \ref{C2.6}(1), so $X_e$ is regular. 

{\bf (2)} If $n=2$, then for $m\in\mathbb N$ with $p^m>\varphi(e)$ and a $d\in\EE_2(k)$ which is a finite \'etale cover of degree $p^m$, $de\in\EE_2(k)$ is surjective with $X_{ed}=X_e$. So there exists $e_1\in\EE_2(k)$ with $X_{e_1}$ singular iff there exists a surjective $e_2\in\EE_2(k)$ with $X_{e_2}$ singular.\end{remark}

\noindent
{\bf P19.} Classify all $m\in\mathbb N^{\ast}$ for which there exists $n\in\mathbb N^{\ast}$ and $e\in\EE_n(k)$ such that, with the notation of Criterion \ref{CRI1}, the set difference $\Irr(X_e\setminus X_e^{\et})\setminus\Irr_{\textup{t}}(X_e\setminus X_e^{\et})$ is non-empty and $\deg(e)=m$.

\smallskip\noindent
{\bf P20.} Describe the submonoid $\MBE_n(k)$ of $\GA_n(k)$ generated by $\BE_n(k)$. For $M\in\{\EE_2(k),\MBE_2(k)\}$, describe the set 
$$\infty^M_p:=\{\bigl(\infty_p(f),\infty_p(g)\bigr)|(f,g)\in k_{\s}[x_1,x_2],\e(f,g)\in M,\infty_p(f)\le\infty_p(g)\}\subset  (\mathbb N^{\ast})^2.$$

\begin{example}\normalfont\label{EX57} 

{\bf (1)} By Remark \ref{R3}(1) we have an inclusion 
$$\infty^M_p\subset \{(1,1)\}\cup\{(l_1,l_2)\in\mathbb (N^{\ast})^2,2\le l_1\le l_2\}.$$ 

{\bf (2)} Suppose that $p=2o+1$ is odd. As $e=\e(x_1+x_1^{o+1}x_2^o,x_2-x_1^ox_2^{o+1})$ is in $\EE_2(k)$ (see $e^1_{t^o+1,-t^o+1;1}$ of Example \ref{EX23-}) and for 
$$(a,b):=\bigl(\e(x_1,-x_1+x_2),\e(x_1+x_2,-x_1+x_2)\bigr)\in\GA_2(k)^2$$ 
we have identities $ae=\e\bigl(x_1+x_1^{o+1}x_2^o,-x_1+x_2-x_1^ox_2^o(x_1+x_2)\bigr)\in\EE_2(k)$ and $be=\e\bigl((x_1+x_2+x_1^ox_2^o(x_1-x_2),-x_1+x_2-x_1^ox_2^o(x_1+x_2)\bigr)\in\EE_2(k)$, it follows that $\{(2,2),(2,3),(3,3)\}\subset\infty^{\EE_2(k)}_p$ (cf.\ \cite{Lang-J}, Rmk.\ 4.9).\end{example}

\begin{remark}\normalfont\label{R33}
We do not know if the set $\infty^M_2\setminus\{(1,1)\}$ is non-empty.
\end{remark} 

\noindent
{\bf P21.} If $(n,m)\in (\mathbb N^{\ast}\setminus\{1\})\times\mathbb N^{\ast}$, is the set $\{e\in\pMor_n(k)\setminus\BE_n(k)|\deg(e)=pm\}$ non-empty? Note that if $e\in\BE_n(k)$ is non-surjective, then $e\notin\pMor_n(k)$.

\smallskip\noindent
{\bf P22.} Classify all triples $(p,n,q)$ with $\deg^f_{p,n,q}\neq\deg_{p,n,q}$ (or with $SEE^{\f-\gen}_{n,q,k}$ not Zariski dense in $SEE_{n,q,k}$).

\smallskip\noindent
{\bf P23.} For each $(n,m)\in (\mathbb N^{\ast}\setminus\{1\})\times\mathbb N^{\ast}$ compute $q_{2,n,m}\in\mathbb N^{\ast}$ of Theorem \ref{T2.1}.

\smallskip\noindent
{\bf P24.} Let $n\in\mathbb N^{\ast}\setminus\{1\}$. Determine the set $\mathfrak D_n$ of Notation \ref{NOT5}. For each pair $(m,q)\in\mathfrak D_n$, compute the finite set $P_{n,m,q}$ of primes introduced in Corollary \ref{C9.1}(1).

\medskip
The following conjecture is our expectation of what \cite{Zh}, Introd., Thm.\ for $n=2$ would become for $n\ge 3$.

\begin{conjecture}\label{CJ4}
Suppose that $\chr(K)=0$ and $n\ge 3$. Then for each endomorphism $\e(g_1,\ldots,g_n)\in\EE_n(K)$ we have $\deg\bigl(\e(g_1,\ldots,g_n)\bigr)\le\prod_{i=1}^{n-1} \deg(g_i)$.
\end{conjecture}

For $n\ge 3$, note that Conjecture \ref{CJ4} is not implied by \cite{Fo2}, Thm.\ 2. It is a weakening of The $n$-Degree Conjecture of Moskowicz (see \cite{Mos1}, Conj.\ 3.1) which is known to be false.

The next conjecture is a natural extrapolation of Corollaries \ref{C10}(1) and \ref{C10+}, \cite{Wa2}, Thm.\ 62,  Section \ref{S21}, and moduli considerations.

\begin{conjecture}\label{CJ5}
If $n\ge 2$, $q\ge 3$, and $\sqrt[n]{p}< q<p$, then $\deg_{p,n,q}^{\f}=1$.
\end{conjecture}

The last conjecture extrapolates Lemma \ref{F1.1} and Proposition \ref{L8.1}.

\begin{conjecture}\label{CJ6}
Let $e\in\EE_n(k)$. If $e$ is a quasi-automorphism in the sense of Definition \ref{D4.2}, then $e\in\BE_n(k)$.
\end{conjecture}

Conjecture \ref{CJ6} is true iff statements (1) to (4) of Proposition \ref{L8.1} are equivalent.

\medskip\noindent
{\bf Acknowledgement.} The first author is indebted to David Wright and Edward Formanek who got him interested in the Jacobian Conjecture. The third author would like to thank SUNY Binghamton and IH\'ES for good working conditions.

\medskip\noindent
{\bf Interests and AI statement.} Competing interests: authors declare none. No AI was used in this work, and no results are computer-based or inspired by computer-generated data.

\medskip\noindent
{\bf Supplemental AI statement of the third author.} Very unfortunately and without me being aware, the in-progress version of this work dated January 14, 2025 was shared publicly at the link \url{https://www.ihes.fr/~gabber/BGV160.pdf} on January 16, 2025. My understanding is that after 18 months it became a feed of the AI models. On July 23, a graduate student sent me the 37 pages pdf file at the link \url{https://people.math.binghamton.edu/adrian/AI.pdf} produced by AI when prompted on how the 7 pages AI paper \url{https://www.ulam.ai/research/jacobian.pdf} relates to BGV160.pdf; a typo in the statement of Theorem \ref{T10}(2.j) was fixed based on what the file mentioned (its application to Theorem \ref{T10}(10.c) did not require any change). 

To my best knowledge:

\medskip
{\bf (1)} there has been no prior work before BGV160.pdf aiming to classify all \'etale endomorphisms of the affine spaces that have geometric degrees 3, see Section 26 of BGV160.pdf; 

\smallskip
{\bf (2)} there has been no prior work before BGV160.pdf that used systematically $xy$, $1+xy$, and $\delta+xy$ in the construction of open embeddings, morphisms, and \'etale endomorphisms, including the morphism $\mathbb A^2_K\rightarrow\mathbb P^1_K$ defined by the rule $(x:1+xy)$; 

\smallskip
{\bf (3)} between BGV160.pdf and the 7 pages AI paper there exists overlap of notation (such as $\Gamma_{}$ for complements) and of intent (such as the variation of the geometric degrees, counting fibers, non-proper locus, families, the usage of $1+xy$, the computation of the normalization $X_e$ which actually plays no role in terms of the counterexamples to the Jacobian Conjecture, etc.) which are also very pertinently and systematically presented in the 37 pages pdf file; 

\smallskip
{\bf (4)} as of today the prompt history, the model version, and the workflow behind the 7 pages AI paper have not been disclosed publicly.

\medskip
My Simons Foundation Travel Award application submitted in January 2026 mentioned, ``The research monograph on automorphisms of affine spaces, joint work with Alexander Borisov and Ofer Gabber, currently 121 pages, is aimed to be made public by July", and, ``On endomorphisms of affine spaces and the Jacobian Conjecture [Problem], joint work with Alexander Borisov and Ofer Gabber, currently 138 pages, is aimed to be made public by the end of 2026 or early 2027." With the AI intervening, I felt compelled to propose to upload them both only 4 days apart. One could ponder if the AI intervention has been positive or negative. The current version is BGV356.pdf dated September 4, 2026.

%Propri\'et\'es topologiques des polynômes de deux variables complexes, et automorphismes alg\'ebriques de l'espace  C2 
%Suzuki, Masakazu
%J. Math. Soc. Japan 26 (1974), 241--257.

%ł

\hbox{}
\hbox{Alexander Borisov,\;\;Email: aborisov@binghamton.edu}
\hbox{Address: Department of Mathematics and Statistics, Binghamton University,}
\hbox{Binghamton, P.\ O.\ Box 6000, New York 13902-6000, U.S.A.}

\hbox{}
\hbox{Ofer Gabber,\;\;\;E-mail: gabber@ihes.fr}
\hbox{Address: IH\'ES, Le Bois-Marie, 35 Route de Chartres,}
\hbox{F-91440 Bures-sur-Yvette, France.}

\hbox{}
\hbox{Adrian Vasiu,\;\;\;Email: avasiu@binghamton.edu}
\hbox{Address: Department of Mathematics and Statistics, Binghamton University,}
\hbox{Binghamton, P.\ O.\ Box 6000, New York 13902-6000, U.S.A.}

\newpage\section*{Notation}\label{NOT6}

\subsection*{Small letters}

$_{}$

\noindent
$\bullet$ \hyperref[D1]{$a$} and \hyperref[D1]{$b$} automorphisms of $\mathbb A^n_R$

\noindent
$\bullet$ \hyperref[PH22]{$c$} a finite flat morphism between varieties, often with a lower right index

\noindent
$\bullet$ \hyperref[PH13a]{$d$} an endomorphism of varieties or schemes or sets

\noindent
$\bullet$ \hyperref[D1]{$e$} an endomorphism of varieties or schemes or sets, often with a lower right index

\noindent
$\bullet$ \hyperref[NOT1]{$e^l$} the $l$-th iterate of an endomorphism $e$, where $l\in\mathbb N$

\noindent
$\bullet$ \hyperref[NOT2]{$e^{\#}$} the endomorphism of $K[x_1,\ldots,x_n]$ such that $e=\Spec (e^{\#})$

\noindent
$\bullet$ \hyperref[PH0.5]{$f$} a polynomial or a Stein factor, often with a lower right index

\noindent
$\bullet$ \hyperref[PH16b]{$g$} and \hyperref[L2]{$h$} polynomials, often with a lower right index

\noindent
$\bullet$ \hyperref[PH0]{$i$} and \hyperref[PH23a]{$j$} indices in $\mathbb N$

\noindent
$\bullet$ \hyperref[PH1]{$k$} an algebraically closed field of characteristic $p$

\noindent
$\bullet$ \hyperref[PH1]{$l$}, \hyperref[PH1]{$m$}, and \hyperref[PH2]{$n$} natural numbers

\noindent
$\bullet$ \hyperref[PH0.5]{$p$} a prime in $\mathbb N$

\noindent
$\bullet$ \hyperref[C9.1]{$p_{n,m,q}$} the smallest prime with $p_{n,m,q}>\max(P_{n,m,q})$

\noindent
$\bullet$ \hyperref[PH0]{$q$} and \hyperref[T1]{$s$} natural numbers

\noindent
$\bullet$ \hyperref[PH5]{$t$} and \hyperref[PH5]{$w$} indeterminates

\noindent
$\bullet$ \hyperref[PH5]{$x$}, \hyperref[PH5]{$y$}, and \hyperref[PH5]{$z$} indeterminates, often with a lower right indices

\subsection*{Small letters symbols}

$_{}$

\noindent
$\bullet$ \hyperref[D5]{$\ciedim(X)$} the complete intersection embedding dimension of an affine regular variety $X$ over $K$

\noindent
$\bullet$ \hyperref[PH11]{$\deg=\deg_{n,K}$} the geometric degree invariant on $\D_n(K)$

\noindent
$\bullet$ \hyperref[D4.3]{$\deg_p(f)$} the modulo $p$-th powers degree of $f$

\noindent
$\bullet$ \hyperref[PH12a]{$\dim(X)$} the dimension of a variety $X$ over $K$

\noindent
$\bullet$ \hyperref[PH12a]{$\edim(X)$} the embedding dimension of an affine variety $X$ over $K$

\noindent
$\bullet$ \hyperref[D4.3]{$\pMor_n(k)$} the submonoid of $\EE_n(k)$ of $p$-morphisms

\noindent
$\bullet$ \hyperref[D2]{$\s(\mathcal S)$} the number of solutions in $k^n$ of an $F_n$-system $\mathcal S$ defined over $k$

\noindent
$\bullet$ \hyperref[D5]{$\sciedim(X)$} the strong complete intersection embedding dimension of an affine regular variety $X$ over $K$

\noindent
$\bullet$ \hyperref[D5]{$\wciedim(X)$} the weak complete intersection embedding dimension of an affine regular variety $X$ over $K$

\noindent
$\bullet$ \hyperref[PH0.5]{$\z(f)$} the number of distinct zeros of a non-zero polynomial in $K[t]$

\subsection*{Capital letters}

$_{}$

\noindent
$\bullet$ \hyperref[L3]{$A$} and \hyperref[PH23b]{$B$} algebras

\noindent
$\bullet$ \hyperref[NOT3]{$A_e$} the integral $k_{\t}[x_1,\ldots,x_n]$-algebra such that $X_e=\Spec A_e$, when $e\in\EE_n(k)$

\noindent
$\bullet$ \hyperref[S26]{$A_{n,q}$} the moduli $\mathbb Z$-algebras of endomorphisms of the $n$-dimensional affine spaces defined by polynomials of degrees at most $q$, when $(n,q)\in (\mathbb N^{\ast})^2$

\noindent
$\bullet$ \hyperref[MS2]{$B_{n,q}$} the moduli $\mathbb Z$-algebras of \'etale endomorphisms of the $n$-dimensional affine spaces defined by polynomials of degrees at most $q$ and with Jacobian determinants equal to $1$, when $(n,q)\in (\mathbb N^{\ast})^2$

\noindent
$\bullet$ \hyperref[PH16]{$C$} a curve

\noindent
$\bullet$ \hyperref[FOOT4]{$Cl(X)$} the divisor class group of a normal variety $X$ over $K$

\noindent
$\bullet$ \hyperref[S8]{$D$} a discrete valuation ring

\noindent
$\bullet$ \hyperref[S8]{$D^{\h}$} the henselization of $D$ 

\noindent
$\bullet$ \hyperref[PH25]{$D_y$} the derivation of $K[x_1,\ldots,x_n]$ with respect to an affine coordinate $y$

\noindent
$\bullet$ \hyperref[S26]{$EE_{n,q}$} the subsheaf of $End_{n,q}$ formed by $n$-tuples that define \'etale endomorphisms

\noindent
$\bullet$ \hyperref[R25]{$EE_{n,q}^{\ast}$} the locally closed subscheme of $End_{n,q}$ formed by $n$-tuples that define \'etale endomorphisms whose Jacobian determinants are units of the bases and of the polynomial things over the bases

\noindent
$\bullet$ \hyperref[S26]{$End_{n,q}$} the scheme over $\Spec\mathbb Z$ of endomorphisms of $\mathbb A^n$ of degrees at most $q$

\noindent
$\bullet$ \hyperref[D2]{$F_n$-system} or \hyperref[D2]{Frobenius system in $n$ indeterminates} is a system of polynomial equations of the form $x_i+f_i(x_1^p,\ldots,x_n)=0$ for $i\in \llbracket1,n\rrbracket$ over an $\mathbb F_p$-algebra

\noindent
$\bullet$ \hyperref[R25]{$GA^{\ast}_{n,q,l}$} a closed subscheme of the product $EE^{\ast}_{n,q}\times_{\Spec Z} EE^{\ast}_{n,q^{n-1}+l}$, when $(n,q,l)\in (\mathbb N^{\ast})^2\times\mathbb N$

\noindent
$\bullet$ \hyperref[S26]{$GA^{\ast}_{n,q}$} the reduced closed subscheme of $End_{n,q}$ whose geometric points define automorphisms, when $(n,q)\in (\mathbb N^{\ast})^2$

\noindent
$\bullet$ \hyperref[P10]{$HEE_{n,q}^{\ast}$} a closed subscheme of $EE_{n,q}^{\ast}$

\noindent
$\bullet$ \hyperref[D5]{$I$} an ideal

\noindent
$\bullet$ \hyperref[MS2]{$I_{n,q}$} the ideal of $A_{n,q}$ such that $B_{n,q}=A_{n,q}/I_{n,q}$, when $(n,q)\in (\mathbb N^{\ast})^2$

\noindent
$\bullet$ \hyperref[PH0.5]{$K$} an algebraically closed field

\noindent
$\bullet$ \hyperref[PG19a]{$L$} a polynomial in $K[x,y]$ or $K[x_1,\ldots,x_n,w]$

\noindent
$\bullet$ \hyperref[D1]{$M$} a subset of $\End_n(R)$ stable under left and right translations by elements of $\GA_n(R)$

\noindent
$\bullet$ \hyperref[NOT3]{$N_e$} the affine variety that is the non-proper locus of $\e\in\EE_n(K)$

\noindent
$\bullet$ \hyperref[NOT3]{$N_e^{\textup{n}-\et}$} a closed subvariety of $N_e$

\noindent
$\bullet$ \hyperref[R6]{$N_l(e)$} a reduced locally closed subvariety of the codomain of $e$

\noindent
$\bullet$ \hyperref[PH24]{$O$} a local ring, often with a lower right index

\noindent
$\bullet$ \hyperref[PH12]{$P$} and \hyperref[PH23d]{$Q$} points

\noindent
$\bullet$ \hyperref[C9.1]{$P_{n,m,q}$} the finite set of primes associated to a pair $(n,m,q)\in\mathfrak D_n$ 

\noindent
$\bullet$ \hyperref[PH2]{$R$} a commutative $\mathbb Z$-algebra

\noindent
$\bullet$ \hyperref[PR24]{$R_u(\mathcal G)$} the unipotent radival of $\mathcal G$

\noindent
$\bullet$ \hyperref[PH2]{$R^{\ast}$} the multiplicative group of units of $R$ 

\noindent
$\bullet$ \hyperref[D2]{$S$} a commutative $\mathbb Z$-algebra related to $R$

\noindent
$\bullet$ \hyperref[PR13]{$SA_{n,q,l}$} a closed subscheme of the product $SEE_{n,q}\times_{\Spec Z} SEE_{n,q^{n-1}+l}$, when $(n,q,l)\in (\mathbb N^{\ast})^2\times\mathbb N$

\noindent
$\bullet$ \hyperref[PR13]{$SA_{n,q}$} the reduction of anyone of the $SA_{n,q,l}$ schemes over $\Spec\mathbb Z$ with $l\in\mathbb N$, when $(n,q)\in (\mathbb N^{\ast})^2$

\noindent
$\bullet$ \hyperref[MS3]{$SEE_{n,q}$} the closed subscheme of $End_{n,q}$ that defines \'etale endomorphisms of Jacobian determinant $1$, when $(n,q)\in (\mathbb N^{\ast})^2$

\noindent
$\bullet$ \hyperref[PH12a]{$T_X$} the tangent bundle of a smooth variety $X$ over $K$

\noindent
$\bullet$ \hyperref[PH11a]{$W$} a scheme, such as a variety over $\Spec K$

\noindent
$\bullet$ \hyperref[EX5+]{$W_e$} the affine variety over $\Spec K$ equipped with a finite morphism $f_e:W_e\rightarrow X_e$ which is generically \'etale of degree $\deg(e)-1$, when $e\in\EE_n(K)\setminus\GA_n(K)$

\noindent
$\bullet$ \hyperref[PH11a]{$W_{\red}$} the reduced scheme associated to $W$

\noindent
$\bullet$ \hyperref[PH12a]{$X$} a variety over $K$

\noindent
$\bullet$ \hyperref[D3]{$X^{\et}$} the \'etale locus of a Jacobian variety $\psi:X\rightarrow X$

\noindent
$\bullet$ \hyperref[PH15]{$X_e$} the affine variety over $K$ associated to $e\in\EE_n(K)$

\noindent
$\bullet$ \hyperref[T3+]{$Y$} an affine variety over $K$ or scheme

\noindent
$\bullet$ \hyperref[PH101]{$Z$} a scheme or an algebraic space 

\subsection*{Capital letters symbols}

$_{}$

\noindent
$\bullet$ \hyperref[EXT13]{$\AffMonomial_n(K)$} a main subset of $\Monomial_n(K)$

\noindent
$\bullet$ \hyperref[PH18]{$\Affintr_n(K)$} the subset of $\Intr_n(K)$ that defines hypersurfaces in $\mathbb A^{n+1}_K$ that have open subvarieties isomorphic to $\mathbb A^n_K$

\noindent
$\bullet$ \hyperref[PH8]{$\BE_n(R)$} the set of basic endomorphisms of an $\mathbb F_p$-algebra $R$

\noindent
$\bullet$ \hyperref[SETS]{$\D_n(R)$} the submonoid of $\End_n(R)$ formed by endomorphisms of $\mathbb A^n_R$ that are fiberwise dominant

\noindent
$\bullet$ \hyperref[PH7.9]{$\EE_n(R)$} the submonoid of $\End_n(R)$ formed by \'etale endomorphisms of $\mathbb A^n_R$

\noindent
$\bullet$ \hyperref[PH6]{$\End_n(R)$} the multiplicative monoid of endomorphisms of $\mathbb A^n_R$

\noindent
$\bullet$ \hyperref[PH6]{$\GA_n(R)$} the group of automorphisms of $\mathbb A^n_R$

\noindent
$\bullet$ \hyperref[PH17]{$\Intr_n(K)$} the subset of $K[x_1,\ldots,x_n]^2$ that defines integral hypersurfaces in $\mathbb A^{n+1}_K$

\noindent
$\bullet$ \hyperref[PH12a]{$\Irr(X)$} the set of irreducible components of a variety $X$ over $K$

\noindent
$\bullet$ \hyperref[CRI1]{$\Irr_{\textup{t}}(X_e\setminus X_e^{\et})$} the set of tame irreducible components of $X_e\setminus X_e^{\et})$ for $e\in\EE_n(k)$

\noindent
$\bullet$ \hyperref[EXT12]{$\Monomial_n(K)$} the subset of $K[x_1,\ldots,x_n]$ formed by irreducible polynomials that are monic in $x_1$

\noindent
$\bullet$ \hyperref[PH18]{$\Nointr_n(K)$} the subset of $\Intr_n(K)$ that defines normal hypersurfaces in $\mathbb A^{n+1}_K$

\noindent
$\bullet$ \hyperref[SETS]{$\QFE_n(R)$} the submonoid of $\End_n(R)$ formed by quasi-finite endomorphisms of $\mathbb A^n_R$ that are fiberwise generically \'etale

\noindent
$\bullet$ \hyperref[SETS]{$\QF_n(R)$} the submonoid of $\End_n(R)$ formed by quasi-finite endomorphisms of $\mathbb A^n_R$

\noindent
$\bullet$ \hyperref[PH12a]{$\Reg(X)$} the regular locus of a variety $X$ over $K$

\noindent
$\bullet$ \hyperref[PH12a]{$\Sing(X)$} the singular locus of a variety $X$ over $K$

\noindent
$\bullet$ \hyperref[PH18]{$\Smintr_n(K)$} the subset of $\Intr_n(K)$ that defines smooth hypersurfaces in $\mathbb A^{n+1}_K$

\noindent
$\bullet$ \hyperref[PH100]{$\Triv(W)$} the set of affine open regular subvarieties of a variety $W$ over $K$ that have trivial tangent bundles

\subsection*{Capital calligraphic letters}

$_{}$

\noindent
$\bullet$ \hyperref[NOT2]{$\mathcal A_n(K)$} the set of affine coordinates of $K[x_1,\ldots,x_n]$

\noindent
$\bullet$ \hyperref[PH13]{$\mathcal B$} a set

\noindent
$\bullet$ \hyperref[PH14c]{$\mathcal C$} a subset

\noindent
$\bullet$ \hyperref[NOT1]{$\mathcal E_l(e)$}, \hyperref[NOT1]{$\mathcal E_l^1(e)$}, and \hyperref[NOT1]{$\mathcal E^1(e)$} subsets of $\mathcal B$ associated an en endomorphism $e:\mathcal B\rightarrow\mathcal B$, where $l\in\mathbb N$

\noindent
$\bullet$ \hyperref[PH12]{$\mathcal F=\mathcal F_{n,K}$} the non-generic fibers set invariant on $\QF_n(K)$

\noindent
$\bullet$ \hyperref[PR24]{$\mathcal G$} a smooth connected linear algebraic group over $\Spec\mathcal K$

\noindent
$\bullet$ \hyperref[EXT11]{$\mathcal I$} a birational morphism between affine varieties over $k$

\noindent
$\bullet$ \hyperref[EXT7]{$\mathcal L$} Lang torsors

\noindent
$\bullet$ \hyperref[S8]{$\mathcal K$} a subfield of $K$, often the field of fractions of a dicrete valuation ring $D$ 

\noindent
$\bullet$ \hyperref[T2]{$\mathcal M$} a morphism between affine varieties over $k$

\noindent
$\bullet$ \hyperref[C2.4]{$\mathcal P^0_f\bigl(\Div(\mathcal K(t)\Spec D\bigr)$} the set of finite non-empty subsets of divisorial valuations of $\mathcal K(t)/\Spec D$ that dominate $D$ and have residue fields transcendental over the residue field of $D$, where $D$ is a discrete valuation ring of field of fractions $\mathcal K$ 

\noindent
$\bullet$ \hyperref[PH7.9]{$\mathcal P_f(\mathbb N)$} the set of finite subsets of $\mathbb N$

\noindent
$\bullet$ \hyperref[EXT10]{$\mathcal Q$} a morphism from an affine surface to $\mathbb P^1_K$, iften with a lower right index

\noindent
$\bullet$ \hyperref[D2]{$\mathcal S$} an $F_n$-system of equations

\noindent
$\bullet$ \hyperref[EXT7]{$\mathcal T$} a birational morphism between affine varieties over $k$

\noindent
$\bullet$ \hyperref[T1]{$\mathcal U$} an $F_2$-system

\subsection*{Small Greek letters}

$_{}$

\noindent
$\bullet$ \hyperref[NOT4]{$\digamma(V)$} the number of irreducible components of the Zariski closure $\overline{V}$ of $V$ in $W$ that have dimension equalto $\dim(\overline{V}$, where $V$ is a constructible subset of a variety $W$ over $\Spec K$

\noindent
$\bullet$ \hyperref[L11]{$\chi(\pi)$} the morphism from $\mathbb A^2_K$ to $\mathbb P^1_K$ associated to a projection $\pi:\mathbb A^2_K\rightarrow\mathbb A^1_K$

\noindent
$\bullet$ \hyperref[PH13]{$\iota$} the iterate inner invariant on $\End_n(K)$ or on endomorphisms of a set

\noindent
$\bullet$ \hyperref[PH16a]{$\kappa$} a field, often of characteristic $p$

\noindent
$\bullet$ \hyperref[PH51]{$\lambda$} an element of some $\mathbb N^l$ 

\noindent
$\bullet$ \hyperref[NOT3]{$\nu(e)$} the number of irreducible components of $N_e$

\noindent
$\bullet$ \hyperref[NOT3]{$\nu_{\textup{n}-\et}(e)$} the number of irreducible components of $N_e^{\textup{n}-\et}$

\noindent
$\bullet$ \hyperref[PH50]{$\phi$} a morphism between varieties over $K$

\noindent
$\bullet$ \hyperref[L11]{$\pi$} a projection between affine spaces

\noindent
$\bullet$ \hyperref[PH7]{$\pi(\mathcal V)$} the degree of a system of equations $\mathcal V$

\noindent
$\bullet$ \hyperref[PH7]{$\pi(g_1,\ldots,g_m)$} the degree of an $m$-tuple $(g_1,\ldots,g_m)$ of polynomials

\noindent
$\bullet$ \hyperref[PH10]{$\pi=\pi_{n,K}$} the algebraic degree invariant on $\End_n(K)$

\noindent
$\bullet$ \hyperref[PH10]{$\pi_{\i}=\pi_{\textup{inner},n,K}$} the inner algebraic degree invariant on $\End_n(K)$

\noindent
$\bullet$ \hyperref[PH10]{$\pi_{\l}=\pi_{\textup{left},n,K}$} the left algebraic degree invariant on $\End_n(K)$

\noindent
$\bullet$ \hyperref[PH10]{$\pi_{\r}=\pi_{\textup{right},n,K}$} the right algebraic degree invariant on $\End_n(K)$

\noindent
$\bullet$ \hyperref[D3]{$\psi$} the endomorphism of a variety over $K$ that defines a Jacobian variety

\noindent
$\bullet$ \hyperref[PH15]{$\psi_e$} the Jacobian variety defined by $e\in\EE_n(K)$

\noindent
$\bullet$ \hyperref[D3]{$\rho$} the class rank

\noindent
$\bullet$ \hyperref[D7]{$\rho_1$} the sphere class rank invariant on $\EE_2(K)$

\noindent
$\bullet$ \hyperref[D7]{$\rho_2$} the almost sphere class rank invariant on $\EE_2(K)$

\noindent
$\bullet$ \hyperref[D7]{$\rho_3$} the quasi-sphere class rank invariant on $\EE_2(K)$

\noindent
$\bullet$ \hyperref[D7]{$\rho_4$} the pseudo-sphere class rank invariant on $\EE_2(K)$

\noindent
$\bullet$ \hyperref[D3]{$\rho_{\et}$} the \'etale class rank

\noindent
$\bullet$ \hyperref[PH14]{$\varphi=\varphi_{n,K}$} the gap (or missed points) invariant on $QF_n(K)$

\noindent
$\bullet$ \hyperref[NOT1]{$\varsigma_e$} the smallest cardinality of the fibers of an endomorphism of a set $e$

\noindent
$\bullet$ \hyperref[NOT1]{$\varrho_e$} the cardinality of $\mathcal E^1_e$, where $e$ is an endomorphism of a set

\noindent
$\bullet$ \hyperref[NOT1]{$\xi_e$} the smallest natural number such that $\mathcal E^1_e=\mathcal E^1_{\xi_e}$, when $e$ is an endomorphism of a set an $\varrho_e\in\mathbb N$

\subsection*{Capital Greek letters}

$_{}$

\noindent
$\bullet$ \hyperref[NOT1]{$\Gamma^+_e$} the union $\cup_{i\in\mathbb N} e^i(\gamma_e)$, where $e$ is an in endomorphism of a set

\noindent
$\bullet$ \hyperref[PH13a]{$\Gamma_e$} the complement of an endomorphism $e$

\noindent
$\bullet$ \hyperref[S20]{$\Phi$} a finite \'etale morphism or a root system, often with a lower right index

\noindent
$\bullet$ \hyperref[PH15]{$\Psi=\Psi_{n,K}$} the Jacobian variety inner invariant (up to isomorphisms)

\noindent
$\bullet$ \hyperref[PH19]{$\Theta_R$} a suitable subset of $R[t]^2$

\noindent
$\bullet$ \hyperref[PH15a]{$\Psi^{\q}=\Psi^{\q}_{n,K}$} the Jacobian variety inner invariant (up to quasi-isomorphisms)

\noindent
$\bullet$ \hyperref[PH16]{$\Sigma$} a morphism between varieties over $k$ 

\noindent
$\bullet$ \hyperref[PH23]{$\Theta^1_{R,m}$} a suitable subset of $\Theta_R$ that depends on $m\in\mathbb N^{\ast}$

\noindent
$\bullet$ \hyperref[EX10]{$\Theta_{\mathcal K}^{\gen}$} a suitable subset of $R[t]^2$ which is the general analog of $\Theta_R$ over a field $\mathcal K$

\noindent
$\bullet$ \hyperref[EX10]{$\Theta_{\mathcal K}^{\n}$} a suitable subset of $R[t]^2$ which is the normal analog of $\Theta_R$ over a field $\mathcal K$

\noindent
$\bullet$ \hyperref[PH19]{$\Theta_R$} a suitable subset of $R[t]^2$

\noindent
$\bullet$ \hyperref[MS4]{$\Upsilon_{n,q,l}$} the functorial maps $\Upsilon_{n,q,l})(R):SA_{n,q,l}(R)\rightarrow GA_n(R)\cap SEE_{n,q}(R)$ induced by the first projections

\noindent
$\bullet$ \hyperref[F5+]{$\Upsilon_{n,q}$} the closed embeddings $SA_{n,q}\rightarrow SEE_{n,q}$, when $(n,q)\in (\mathbb N^{\ast})^2$

\noindent
$\bullet$ \hyperref[MS1]{$\Xi_{n,q}$} the truncation closed embeddings $End_{n,q}\rightarrow End_{n,q+1}$, when $(n,q)\in (\mathbb N^{\ast})^2$

\subsection*{Blackboard bold letters}

$_{}$

\noindent
$\bullet$ \hyperref[PH3]{$\mathbb A^n$} and \hyperref[PH4]{$\mathbb A^n_R$} the affine spaces of dimension $n$ over $\Spec\mathbb Z$ and $\Spec R$

\noindent
$\bullet$ \hyperref[EXTRA2]{$\mathbb B$} an exceptional divisor of an affine variety over $k$, often with a lower or/and an upper right index

\noindent
$\bullet$ \hyperref[EXTRA3]{$\mathbb B^2_{m,K}$} the quasi-affine variety over $\Spec K$ obtained from $\mathbb A^2_K$ by removing any $m$ distinct points

\noindent
$\bullet$ \hyperref[PH21]{$\mathbb D^1_{f,g,R}$} a suitable divisor of $\mathbb S^2_{f,g,R}$

\noindent
$\bullet$ \hyperref[EXTRA1]{$\mathbb E$} an exceptional divisor of an affine variety over $k$, often with a lower or an upper right index

\noindent
$\bullet$ \hyperref[PH15c]{$\mathbb F_{p^q}$} a finite field with $p^q$ elements

\noindent
$\bullet$ \hyperref[PH3]{$\mathbb G_{\a}$} and \hyperref[PH4]{$\mathbb G_{\a,R}$} the split connected unipotent groups schemes over $\Spec\mathbb Z$ and $\Spec R$ of relative dimension $1$

\noindent
$\bullet$ \hyperref[PH3]{$\mathbb G_{\m}$} and \hyperref[PH4]{$\mathbb G_{\m,\mathbb R}$} the split multiplicative group schemes over $\Spec\mathbb Z$ and $\Spec R$ of relative dimension $1$

\noindent
$\bullet$ \hyperref[PH18]{$\mathbb H$} a hypersurface in $\mathbb A^{n+1}_K$, often with a lower or/and an upper right index

\noindent
$\bullet$ \hyperref[EXTRA5]{$\mathbb I$} a finite subset of $k$, often with a lower right index

\noindent
$\bullet$ \hyperref[PH15a]{$\mathbb J^{\q}_n(K)$} the set of quasi-isomorphism classes of Jacobian varieties over $K$ of dimension $n$

\noindent
$\bullet$ \hyperref[EXTRA6]{$\mathbb J_n(K)$} the set of isomorphism classes of Jacobian varieties over $K$ of dimension $n$

\noindent
$\bullet$ \hyperref[EX23]{$\mathbb K$} hypersurfaces over $\Spec k$ with lower and upper right indices

\noindent
$\bullet$ \hyperref[PH3]{$\mathbb P^n$} and \hyperref[PH4]{$\mathbb P^n_R$} the projective space of dimension $n$ over $\Spec\mathbb Z$ and $\Spec R$

\noindent
$\bullet$ \hyperref[C2.6]{$\mathbb S^2_Y$} the sphere-like surface associated to an open embedding $\imath:\mathbb A^2_K\rightarrow X$ and an irreducible component $Y$ of the complement $X\setminus\Imm(\imath)$

\noindent
$\bullet$ \hyperref[PH19]{$\mathbb S^2_{f,g,R}$} a suitable hypersurface in $\mathbb A^3_R$, where $(f,g)\in\Theta_R$

\noindent
$\bullet$ \hyperref[T2]{$\mathbb V$}, \hyperref[T2]{$\mathbb W$}, \hyperref[T2]{$\mathbb X$}, and \hyperref[T2]{$\mathbb Y$} affine smooth varieties over $k$

\subsection*{Special characters}

$_{}$

\noindent
$\bullet$ \hyperref[PH12a]{$1_X$} the identity automorphism of $X$

\noindent
$\bullet$ \hyperref[EXT9]{$\ell$} a linear form

\noindent
$\bullet$ \hyperref[D4.1]{$\Re$} an invariant on $\D_n(k)$

\noindent
$\bullet$ \hyperref[PH0.5]{$\chr(K)$} the characteristic of the field $K$

\noindent
$\bullet$ \hyperref[NOT4]{$ \hbar(V)$} the pair $\bigl(\dim(\bar{V}),\digamma(\overline{V})\bigr)\in\mathbb N^2$

\noindent
$\bullet$ \hyperref[PH20]{$\imath$} an open embedding, open with a lower right index

\noindent
$\bullet$ \hyperref[PH15]{$\imath_e$} the open embedding of $\mathbb A_{K,\s}^n$ in $X_e$ associated to $e\in\EE_n(K)$

\noindent
$\bullet$ \hyperref[D4.3]{$\infty_p(f)$} the number of points at infinity modulo $p$ for a suitable $f\in k[x_1,\ldots,x_n]$ which is not a $p$-th power

\noindent
$\bullet$ \hyperref[PH1]{$\llbracket l,m\rrbracket$} the set of integers $i$ that satisfy $l\le i\le m$

\noindent
$\bullet$ \hyperref[D4.1]{$\mathfrak B$} the basic number invariant of endomorphisms of affine spaces over $k$

\noindent
$\bullet$ \hyperref[NOT5]{$\mathfrak D_n$} the set formed by all pairs $(m,q)\in (\mathbb N^{\ast}\setminus\{1\})\times\mathbb N^{\ast}$ for which the following set $\{e\in\EE_n(\mathbb C)|\deg(e)=m,\pi(e)\le q\}$ is empty

\noindent
$\bullet$ \hyperref[D4.1]{$\mathfrak V_n(k)$} the set of $k$-vector subspaces of $k_{\s}[x_1,\ldots,x_n]$ containing $k_{\s}[x_1^p,\ldots,x_n^p]$

\noindent
$\bullet$ \hyperref[EX5-]{$\mathfrak w$} a discrete valuation

\subsection*{Conditions}

$_{}$

\noindent
$\bullet$ \hyperref[T10]{($S_3$)} a degree $3$ condition

\noindent
$\bullet$ \hyperref[R28]{($\bigtriangledown_F$)} the second \'etale characteristic $2$ condition

\noindent
$\bullet$ \hyperref[R28]{($\bigtriangleup_F)$} the first \'etale characteristic $2$ condition

\noindent
$\bullet$ \hyperref[T10]{($\dag^e_h$)} a degree $3$ condition

\noindent
$\bullet$ \hyperref[PH52]{($\flat_{e,P}$)} an iteration condition on an endomorphism $e:\mathcal B\rightarrow\mathcal B$ of a set $\mathcal B$ and $P\in\mathcal B$

\noindent
$\bullet$ \hyperref[PH52]{($\natural_{e,P}$)} an iteration condition on an endomorphism $e:\mathcal B\rightarrow\mathcal B$ of a set $\mathcal B$ and $P\in\mathcal B$

\noindent
$\bullet$ \hyperref[L9]{($\partial_1$)} the first partial derivatives condition

\noindent
$\bullet$ \hyperref[L9]{($\partial_2$)} the second partial derivatives condition

\noindent
$\bullet$ \hyperref[T10]{($\ast^e_h$)} a multiplicity condition for a suitable element $h\in K[x_1,\ldots,x_n]$ with respect to $e\in\EE_n(K)$

\noindent
$\bullet$ \hyperref[PR1]{($\triangleleft$)} the left monomial condition

\noindent
$\bullet$ \hyperref[PR1]{($\triangleright$)} the right monomial condition

\noindent
$\bullet$ \hyperref[PH53]{($\w\natural_{e,P}$)} an iteration condition on an endomorphism $e:X\rightarrow X$ of a variety $X$ over $\Spec K$ and a point $P\in X(K)$

\printindex

\end{document}